%% file: pq.tex
\documentclass{amsbook}

\include{cubical-pre}

\begin{document}

\begin{titlepage}
\enlargethispage{3\baselineskip}
\thispagestyle{empty}
\AddToShipoutPictureBG*{\includegraphics[width=\paperwidth,height=\paperheight]{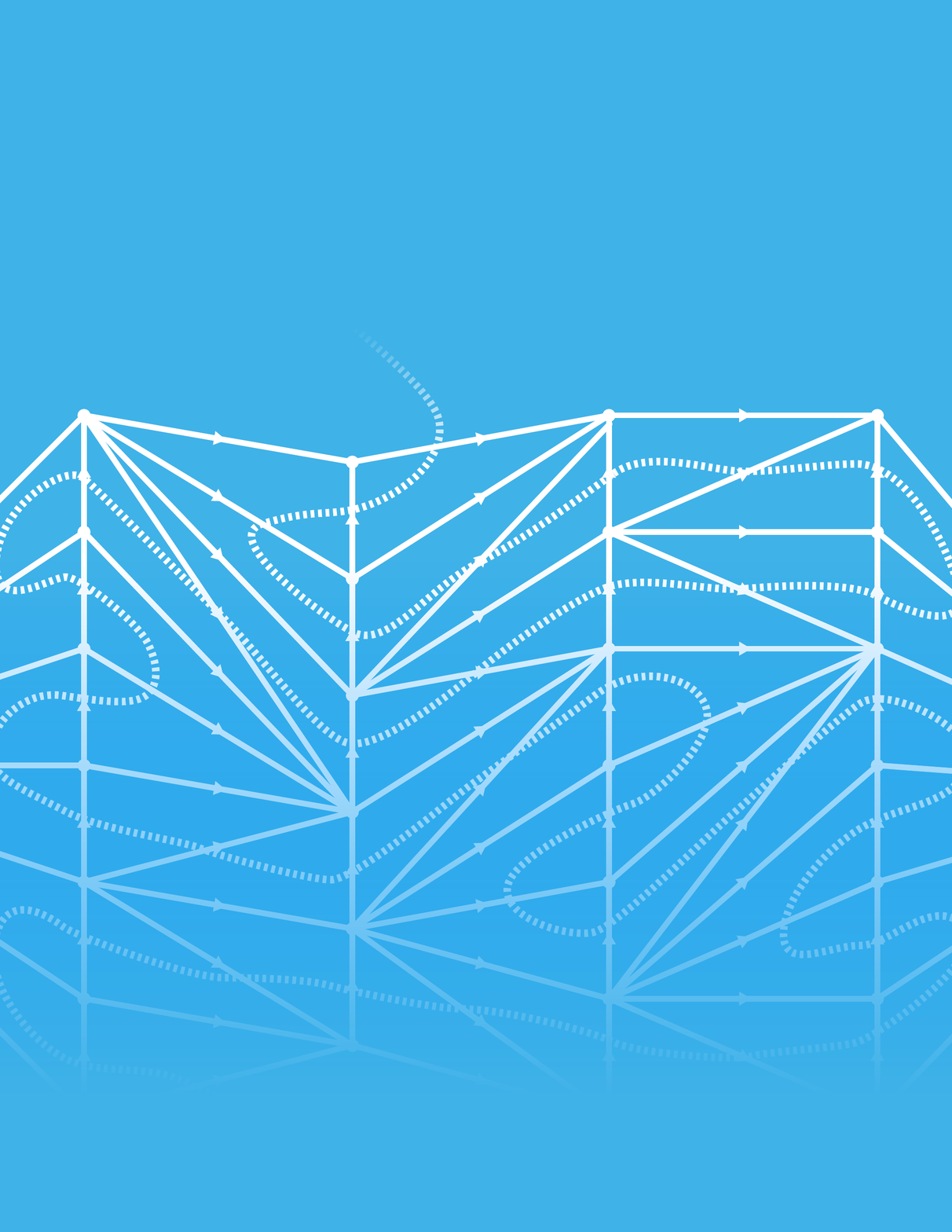}}
\begingroup
\begin{center}
\footnoterulefalse 
\renewcommand{\thefootnote}{} 
    {\color{white}\Huge\bfseries Framed combinatorial topology\par}
    \vspace{12pt}
    {\color{white}\LARGE\mdseries
    Christoph Dorn
    \& Christopher L. Douglas{\footnote[1]{\color{white}\textsc{Mathematical Institute, University of Oxford, Oxford OX2 6GG, United Kingdom}\\ Email addresses:  \href{mailto:dorn@maths.ox.ac.uk}{\color{white}\texttt{dorn@maths.ox.ac.uk}} \&  \href{mailto:cdouglas@maths.ox.ac.uk}{\color{white}\texttt{cdouglas@maths.ox.ac.uk}}}}\par}
    \vspace{14pt}
    {\color{white} December 2021}
    \vfil
\end{center}
\endgroup
\end{titlepage}

\frontmatter

%
%

\ignore{
Framed combinatorial topology is a novel theory describing combinatorial phenomena arising at the intersection of stratified topology, singularity theory, and higher algebra.  The theory synthesizes elements of classical combinatorial topology with a new combinatorial approach to framings.  The resulting notion of framed combinatorial spaces has unexpectedly good behavior when compared to classical, nonframed combinatorial notions of space.  In discussing this behavior and its contrast with that of classical structures, we emphasize two broad themes, computability in combinatorial topology and combinatorializability of topological phenomena.  The first theme of computability concerns whether certain combinatorial structures can be algorithmically recognized and classified.  The second theme of combinatorializability concerns whether certain topological structures can be faithfully represented by a discrete structure.  Combining these themes, we will find that in the context of framed combinatorial topology we can overcome a set of fundamental classical obstructions to the computable combinatorial representation of topological phenomena.
}

%
%




\makeatletter
\patchcmd{\@makeschapterhead}{7.5pc}{6pc}{}{}
\makeatother

\setcounter{tocdepth}{2}
\tableofcontents

\renewcommand{\thesection}{I.{\arabic{section}}}
\renewcommand{\thesubsection}{\thesection.{\arabic{subsection}}}
\renewcommand\thefigure{I.\arabic{figure}}
\footnoteruletrue
\renewcommand*{\thefootnote}{\arabic{footnote}}
\setcounter{footnote}{0}

\mainmatter

\chapter*{Introduction}

Framed combinatorial topology is a novel theory describing combinatorial phenomena arising at the intersection of stratified topology, singularity theory, and higher algebra.  The theory synthesizes elements of classical combinatorial topology with a new combinatorial approach to framings.  The resulting notion of framed combinatorial spaces has unexpectedly good behavior when compared to classical, nonframed combinatorial notions of space.  In discussing this behavior and its contrast with that of classical structures, we emphasize two broad themes, `computability in combinatorial topology' and `combinatorializability of topological phenomena'.  The first theme of computability concerns whether certain combinatorial structures (such as simplicial complexes homeomorphic to spheres) can be algorithmically recognized and classified.  The second theme of combinatorializability concerns whether certain topological structures (such as manifolds) can be faithfully represented by a discrete structure.  Combining these themes, we will find that in the context of framed combinatorial topology we can overcome a set of fundamental classical obstructions to the computable combinatorial representation of topological phenomena.

\vspace{10pt}
We begin this introduction by elaborating the themes of computability and of combinatorializability in, respectively, \autoref{intro:computability} and \autoref{intro:combinatorializability}.  We then give a more formal overview of our results in \autoref{intro:summary}, a chapter-by-chapter outline in \autoref{intro:overview}, and an outlook on the larger program and aims of the subject in \autoref{intro:outlook}.

\section{Computability in combinatorial topology} \label{intro:computability}

Computability is the ability to solve a `general problem' by a `general method', that is, the ability to write a step-by-step procedure which for each specific instance of a problem computes a solution.  Combinatorial topology provides, in a sense, a computation-oriented foundation for the study of spaces, by encoding space in discrete structures  \cite{rourke2012introduction} \cite{bryant2002piecewise}.  However, many fundamental problems in combinatorial topology turn out to be computably intractable; such problems include the following \cite{markov1958insolubility} \cite{volodin1974problem} \cite{nabutovsky1995einstein} \cite{nabutovsky2003fractal} \cite{weinberger2004computers} \cite{poonen2014undecidable}.

\begin{enumerate}
\item The statement `The simplicial complex $K$ is homeomorphic to the $n$-disk' cannot be computably verified for general finite complexes $K$.  Similarly, the statement `The simplicial complex $K$ is homeomorphic to a manifold' cannot be computably verified.
\item These uncomputability issues remain in the piecewise linear setting: the statement `The simplicial complex $K$ piecewise linearly subdivides the $n$-simplex' cannot be verified in general, and neither can the statement `The simplicial complex $K$ is a piecewise linear manifold'.  In particular, one cannot classify all subdivisions of the $n$-simplex with a given number of simplices, nor classify all piecewise linear $n$-manifolds with a given number of simplices.
\item More generally, it is impossible to algorithmically decide whether two simplicial complexes $K$ and $L$ have homeomorphic, or piecewise linearly homeomorphic, geometric realizations.
\item Similarly, given two embedded, or piecewise linearly embedded, simplicial complexes $K \into \lR^n$ and $L \into \lR^n$, one cannot in general determine whether the embeddings are ambient homeomorphic, respectively ambient piecewise linearly homeomorphic.
\end{enumerate}

\nid One could view these failures of computability as unavoidable imperfections of mathematics as we know it, or one can see them as failures of the classical simplicial method of combinatorializing topological structures.  Adopting the latter viewpoint, one may hope for a form of combinatorialization with better computability properties, for instance in which one can recognize combinatorial disks and classify combinatorial manifolds.

The first central theme of this book is that, though typical simplicial methods do not provide a computable foundations for combinatorial topology, there is a different approach, using \emph{framed} combinatorial spaces, that may provide a more suitable basis for computable combinatorial topology.  Our theory of `framed combinatorial topology' differs in two fundamental respects from classical piecewise linear topology: first, we endow simplices and simplicial complexes with a combinatorial framing structure, and second, we generalize the resulting class of `framed simplicial complexes' to a broader class of `framed regular cell complexes'.  Though classical regular cells are much less tractable even than simplices---indeed even the list of cell shapes is uncomputable---it will turn out that framed regular cells arise as iterated constructible combinatorial bundles and therefore both these cells and their complexes are, remarkably, algorithmically classifiable.

\pause

A classical frame of an $m$-dimensional vector space is an ordered choice of $m$ linearly independent vectors.  We will define a combinatorial frame of an $m$-simplex to be an ordered choice of $m$ vectors in the spine of the simplex.  To make sense of a frame on a simplicial complex, we need a notion of the compatibility of frames along faces shared between simplices.  The restriction of a frame of a simplex to a face gives not only information about a frame of the face but also about how that restricted frame embeds in the ambient frame of the simplex; we will be primarily concerned with the resulting notion of embedded framed simplex, and a framed simplicial complex will be a simplicial complex with compatible embedded frames on all its simplices.

Regular cell complexes, that is those complexes whose attaching maps are injective, generalize simplicial complexes by allowing cells of `polytopic' shapes instead of merely `triangular' shapes.  Regular cells can be identified with the geometric realizations of their face posets~\cite{bjorner1984posets}, and via that identification they obtain piecewise linear simplicial subdivisions.  We use that simplicial structure, together with our notion of framed simplicial complexes, to define framings of regular cells and identify a tractable class of such cells, namely those that are `flat' in that they admit a framed embedding into euclidean space.   A framed regular cell complex, finally, will be a regular cell complex with compatible choices of flat framings on each of its cells.

A space with a homeomorphism to a regular cell complex is `cellulated', as a space with a homeomorphism to a simplicial complex is `triangulated'.  The fact that cellulated spaces have played a less prominent role than triangulated spaces in classical combinatorial topology is partially due to the aforementioned fundamental computability obstruction: it is impossible to classify all the possible shapes of regular cells, in the sense that one cannot produce a list of all regular cells with a given number of faces, in general; said another way, there is no general algorithm for deciding whether a given poset is the face poset of a regular cell, even though there are only finitely many posets of a given size.  Endowing regular cells with a framing overcomes this fundamental issue: framed regular cells, in contrast to their nonframed counterparts, are classifiable.  Specifically, given a poset together with a framing of its geometric realization, we can algorithmically recognize whether the poset is the face poset of a framed regular cell.  This is possible, at root, because we will discover that flat framed regular cell complexes are the geometric realizations of a novel combinatorial structure, called `trusses', which are iterated constructible bundles of oriented fence posets.

Framed regular cells strike an unlikely and delicate balance, being simultaneously a class of shapes that is tractable (in that they are algorithmically recognizable, unlike ordinary regular cells) and also a class of shapes that is quite general (unlike ordinary simplices).  The generality of the shapes of framed regular cells provides unique combinatorial representatives in a way that is unthinkable with simplicial structures and unknown with any other class of shapes: a flat framed regular cell complex has a computable unique minimal cell structure.  Having a computably unique representation of these complexes makes algorithmically decidable almost any question about them; for instance, it follows that framed homeomorphism of these complexes is decidable, in stark contrast to the classical (nonframed simplicial) situation.  In this and other related respects, working with framed regular cells and their complexes provides, finally, a computable framework for combinatorial models of spaces.

\section{Combinatorializability of topological phenomena} \label{intro:combinatorializability}

Combinatorics is primarily concerned with discrete, and often finite, structures whose constituents can be counted.  Topology, by contrast, is primarily concerned with the continuous structure of spaces.  The `combinatorializability' of topological phenomena refers to the ability to faithfully encode continuous objects (spaces, manifolds, continuous maps, bordisms, et cetera) in discrete or finite data structures.  This faithful encoding depends both on having a combinatorial representation of the object in question, and on knowing that representation is unique up to some specified combinatorial equivalence relation.

There are by now various known instances of topological phenomena that cannot be faithfully combinatorialized, or even combinatorialized at all, giving an impression of a mysterious and insurmountable divide between topological spaces and any discrete representations of those spaces.  A headline instance of this divide is the disproven `Hauptvermutung' \cite{ranicki1996hauptvermutung}, a conjecture that, roughly speaking, claimed that `topological isomorphism' (meaning homeomorphism) coincides with `combinatorial isomorphism' (meaning piecewise linear homeomorphism).  This conjecture would in particular imply that combinatorial spaces (that is, geometric realizations of simplicial complexes) that are homeomorphic are also piecewise linear homeomorphic.  This intuitive, presumptive claim was eventually disproven \cite{milnor1961two} by the explicit construction of homeomorphic finite simplicial complexes that are not piecewise linear homeomorphic.

A flurry of results followed in subsequent decades \cite{kirby1969triangulation} \cite{akin1969manifold} \cite{hirsch1974smoothings} \cite{kirby1977foundational} \cite{freedman1982topology} \cite{donaldson1983application}, quantifying the divide not only between the `continuous' and the `combinatorial', but also between the `combinatorial' and the `smooth' conceptions of space.  Recently, a disproof of the triangulation conjecture \cite{manolescu2016pin} established an especially stark gap, that in every dimension greater than 4 there exist compact topological manifolds that do not even admit a triangulation.  (It will be pertinent later that most instances of the classical topological--combinatorial gap rely on certain infinitary or `wild' topological constructions.)  By contrast, smooth manifolds always admit triangulations and all triangulations of a smooth manifold are `combinatorially isomorphic'.  However, smooth manifolds that are not smoothly isomorphic may nevertheless be combinatorially isomorphic \cite{milnor1956manifolds}, and combinatorial manifolds need not admit any smooth structure  \cite{kervaire1960manifold}.

One might dream of a topological foundations or combinatorial framework in which the mismatch between the continuous, combinatorial, and smooth conceptions of space would, at least to some extent, be lessened.  One could imagine, for instance, a discrete, perhaps infinitary, combinatorial theory that faithfully represents a delineated class of relevant continuous phenomena, or a discrete, perhaps finitary, combinatorial theory that suitably encodes smooth behavior.  Each of these two comparative visions has been pursued, to some but not complete satisfaction: for instance, an `o-minimal' approach to tame topology provides a method for excluding certain wild topological structures \cite{grothendieck1997esquisse} \cite{shiota2014minimal}, while a `matroid' perspective aims for a direct combinatorial description of smooth structures \cite{macpherson1991combinatorial}.

The second central theme of this book is that, in contrast to the classical gap between topological and combinatorial phenomena, in our framed combinatorial setting there is a faithful comparison between framed topological and framed combinatorial phenomena.  Furthermore, we expect framed combinatorial structures also faithfully encode all framed smooth phenomena, and therefore will provide an unexpected unification of the continuous, combinatorial, and smooth perspectives on space.  The `framed topological' side of these comparisons will be a class of tame topological structures called `flat framed stratifications'.  These stratifications are `flat' and `framed' by an embedding in standard euclidean space, and they are tame in that we insist the stratification admit a refining `mesh', which is a cellulation by framed regular cells; this cellulation requirement is analogous to working with triangulable spaces and therefore excluding, a priori, certain wild behavior.  These mesh cellulations are iterated constructible bundles of stratified 1-manifolds, and will be a precise topological counterpart of the iterated constructible combinatorial structure of trusses mentioned earlier.

The chain of associations, from a flat framed stratification to its mesh cellulation to the corresponding combinatorial truss, does not by itself necessarily ensure a faithful combinatorialization of flat framed stratified topology.  As a space can have various inequivalent triangulations, a flat framed stratification could in theory have various inequivalent meshes (and therefore corresponding trusses)---however, we will prove, crucially, that such a stratification always has a unique coarsest compatible mesh.  This uniqueness is an unexpected and stark counterpoint to the classical situation: given two triangulations of a space, traditionally one aims (and fails) to construct a mutual \emph{refinement} and thereby verify their combinatorial equivalence; now instead, given two mesh cellulations of a stratification, we construct a canonical mutual \emph{coarsening} and thus establish the desired combinatorial equivalence.  The proof of this canonical coarsening relies, of course, on the generality of the shapes of framed regular cells, by contrast with the constrained shapes of classical simplices.  This canonical coarsening provides the desired faithful combinatorialization of topological phenomena in flat framed euclidean space; indeed, we will establish the `flat framed Hauptvermutung', that for flat framed stratifications, framed homeomorphism classes coincide with framed piecewise linear homeomorphism classes.

Regarding the combinatorialization of smooth phenomena, we will conjecture that any smooth manifold can be represented as a flat framed stratification (via a generic embedding in euclidean space) and that the resulting combinatorial representation as a truss faithfully encodes the smooth structure.  We will revisit the context and plausibility of this smooth combinatorialization conjecture in the outlook, \autoref{intro:outlook} below.

\section{Overview} \label{intro:summary}

We collect and summarize our main theorems, and along the way further describe and illustrate our core definitions.  Recall that a framed simplex is an ordinary simplex together with frame, that is a choice of order of its spine vectors.  More generally, a framed regular cell is an ordinary regular cell together with a suitably compatible choice of frames on each simplex in the cell's face poset.  Though it is impossible to classify regular cells, by contrast framed regular cells are classifiable.  The classifying combinatorial structure will be a special case of the notion of `trusses', which are iterated constructible poset bundles defined as follows.

\begin{introdef}[Trusses] \label{intro:defn-truss}
A `1-truss' is a fence poset equipped with a total `frame' order on its elements.  An `$n$-truss' is a length-$n$ tower of constructible bundles of 1-trusses.
\end{introdef}

\begin{figure}[ht]
    \centering
    \def\svgwidth{1\columnwidth}
    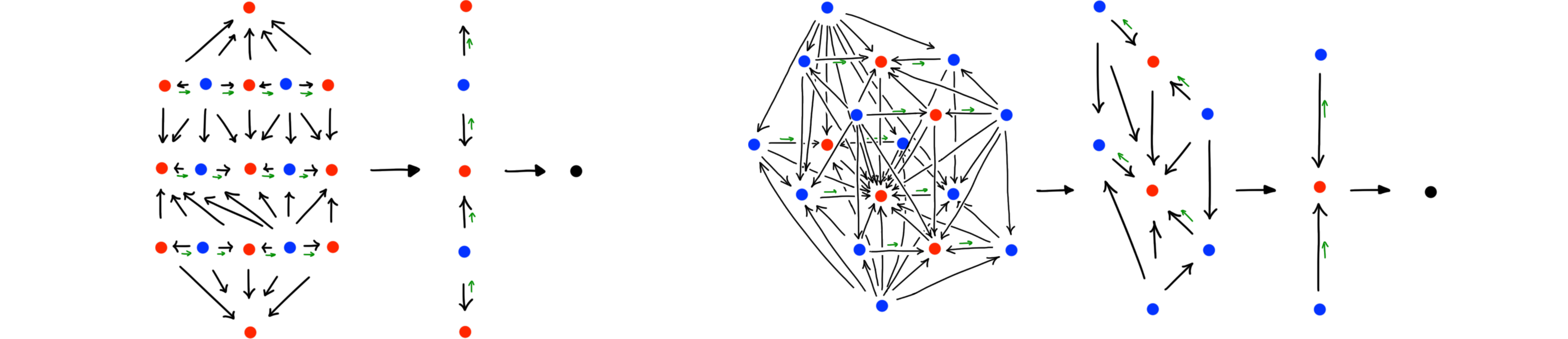

    \caption{A closed 2-truss and an open 3-truss.}
    \label{fig:a-closed-2-truss-and-an-open-3-truss}
\end{figure}

\nid The notions of `1-truss', their `constructible bundles', and `$n$-truss' are given more precisely in, respectively, \autoref{defn:general-1-trusses}, \autoref{defn:1-truss-bundle}, and \autoref{defn:n-trusses}.  Elements of a 1-truss that are targets of poset arrows are called `singular' or `dimension 0', while elements that are sources of poset arrows are called `regular' or `dimension 1'.  A 1-truss is `closed' if both its endpoints are singular, and `open' if both are regular; an $n$-truss is `closed' or `open' if all its fiber 1-trusses are closed or open respectively.  In \autoref{fig:a-closed-2-truss-and-an-open-3-truss} we illustrate a closed 2-truss and an open 3-truss; singular elements are red, regular elements are blue, and the frame orders are indicated by green arrows.

A closed truss is called a `truss block' if its total poset has an initial element; these truss blocks provide the combinatorial correlate of framed regular cells.

\begin{introthm}[Classification of framed regular cells] \label{intro:thm-cell-to-block}
The category of framed regular cells is equivalent to the category of truss blocks.
\end{introthm}

\begin{figure}[ht]
    \centering
    \def\svgwidth{1\columnwidth}
    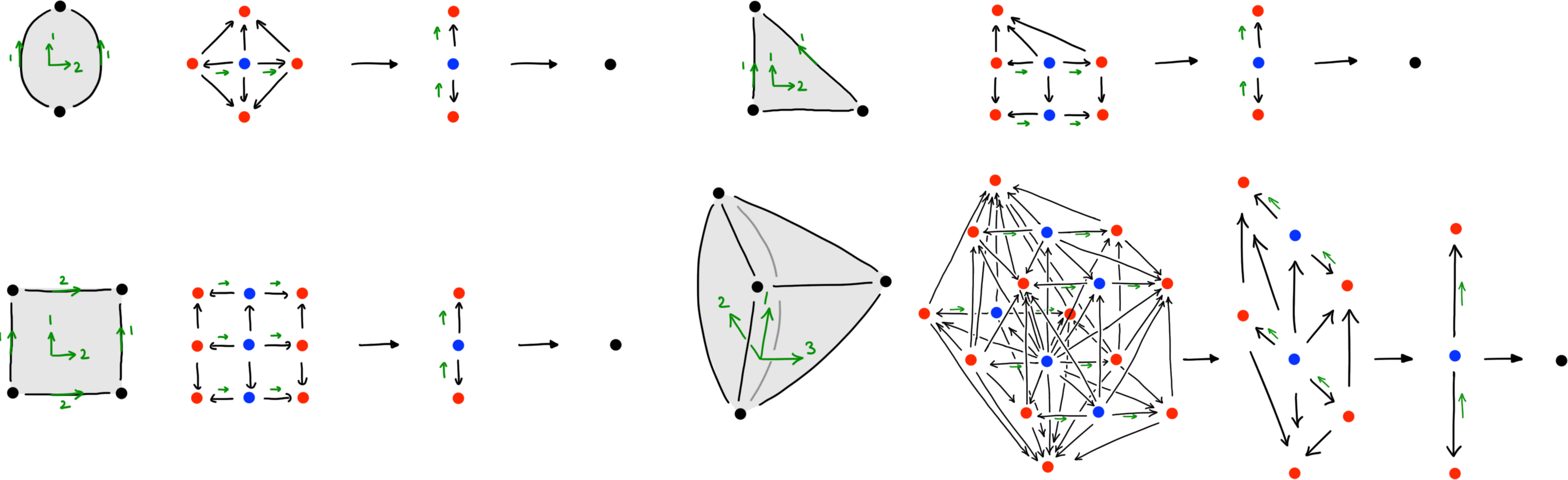

    \caption{Framed regular cells and their classifying truss blocks.}
    \label{fig:classifying-framed-regular-cells-intro}
\end{figure}

\nid This result will appear as \autoref{thm:classification-of-cells}.  Given a framed regular cell, the total poset of the classifying truss block is the face poset of the cell; by sequentially projecting out frame vectors, this poset determines a tower of 1-truss bundles.  In \autoref{fig:classifying-framed-regular-cells-intro} we illustrate a few framed regular cells and their corresponding truss blocks.

The classification of framed regular cells by truss blocks implies a corresponding classification for framed regular cell complexes.  As a simplicial set is a presheaf on the category of simplices, similarly a `truss block set' is a presheaf on the category of truss blocks; a simplicial set is regular if all its simplices embed into its realization, and similarly a truss block set is regular if all its blocks embed into its realization.

\begin{introthm}[Classification of framed regular cell complexes] \label{intro:thm-cell-to-block-set}
The category of framed regular cell complexes is equivalent to the category of regular truss block sets.
\end{introthm}

\nid This result appears in the main text as \autoref{thm:classification-of-framed-cell-cplx}.

\pause

The face poset of a framed regular cell is a truss block; the geometric realization of a truss block is a framed regular cell.  More general trusses also have geometric realizations as `meshes', which are iterated constructible bundles of stratified lines, analogous to \autoref{intro:defn-truss}, as follows.

\begin{introdef}[Meshes] \label{intro:defn-meshes}
A `1-mesh' is a contractible 1- or 0-manifold, stratified by open intervals and points, and equipped with a framing.  An `$n$-mesh' is a length-$n$ tower of constructible bundles of 1-meshes.
\end{introdef}

\begin{figure}[ht]
    \centering
    \def\svgwidth{1\columnwidth}
    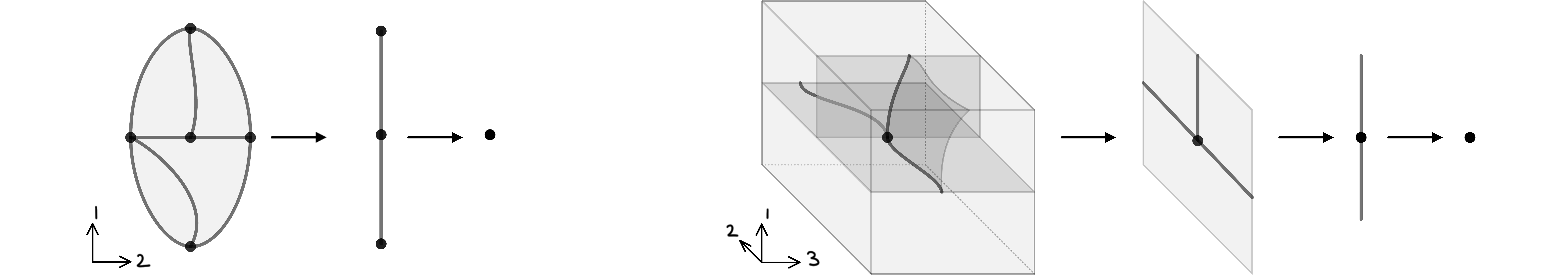

    \caption{A closed 2-mesh and an open 3-mesh.}
    \label{fig:a-closed-2-mesh-and-an-open-3-mesh}
\end{figure}

\nid The notions of `1-mesh', their `constructible bundles', and `$n$-mesh' are given more precisely in, respectively, in \autoref{defn:1-mesh}, \autoref{defn:1-mesh-bundles}, and \autoref{defn:n-meshes}.  An $n$-mesh is `closed' if its total space is compact, and is `open' if its total space is an open disc.  In \autoref{fig:a-closed-2-mesh-and-an-open-3-mesh} we illustrate a closed 2-mesh and an open 3-mesh.

The correspondence of meshes and trusses is at the core of our combinatorialization of topological phenomena.

\begin{introthm}[Equivalence of meshes and trusses] \label{intro:thm:truss-mesh-eqv}
The topological category of closed, respectively open, meshes is weakly equivalent to the discrete category of closed, respectively open, trusses.
\end{introthm}

\nid This result appears in a more precise form as \autoref{thm:diagram_classification}.  Recall the entrance path poset of a stratified space has an element for each stratum and an arrow indicating when a stratum intersects the closure of another stratum.  The above equivalence takes a mesh, a tower of stratified spaces, to the truss given by the tower of corresponding entrance path posets.  As an illustration, note that the meshes in \autoref{fig:a-closed-2-mesh-and-an-open-3-mesh} yield, on application of entrance path posets, the trusses in \autoref{fig:a-closed-2-truss-and-an-open-3-truss}.

A closed mesh will be called a mesh cell if it is the closure of a single stratum.  Naturally, \autoref{intro:thm:truss-mesh-eqv} restricts to an equivalence of mesh cells and truss blocks.  Combining this equivalence with the earlier \autoref{intro:thm-cell-to-block}, relating truss blocks and framed regular cells, yields a correspondence of mesh cells and framed regular cells.

\begin{introcor}[Equivalence of mesh cells and framed regular cells] \label{intro:cor-meshcells-regularcells}
The topological category of mesh cells is weakly equivalent to the discrete category of framed regular cells.
\end{introcor}

\nid This result appears as \autoref{cor:eqv-mesh-cells-and-framed-reg-cells}.  Note well that regular cells are at root combinatorial objects and they come with a canonical piecewise linear structure, whereas mesh cells are purely topological objects and a priori have no piecewise linear structure; thus this seemingly innocuous result provides a fundamental bridge between the topological and piecewise linear contexts.

\pause

Recall that a basic unsolvable problem of classical combinatorial topology is to classify subdivisions of the $n$-simplex.  By contrast, leveraging the above connection between framed regular cells and mesh cells, we can classify framed subdivisions of framed regular cells.

\begin{introthm}[Classification of subdivisions of framed regular cells] \label{intro:thm:subdiv}
A framed regular cell complex framed subdivides a framed regular cell if and only if, as a framed stratified space, the closed mesh corresponding to the regular cell complex refines the closed mesh cell corresponding to the regular cell.
\end{introthm}

\begin{figure}[ht]
    \centering
    \def\svgwidth{1\columnwidth}
    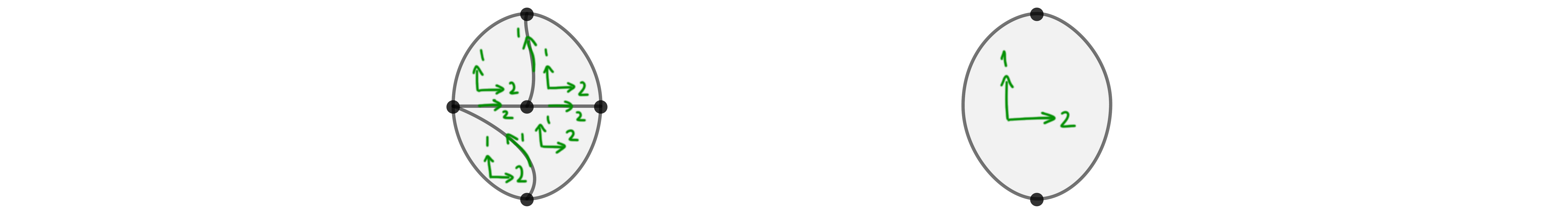

    \caption{The subdivision of a framed regular cell.}
    \label{fig:the-subdivision-of-a-framed-regular-cell}
\end{figure}

This result appears in more precise form as \autoref{cor:fr-subdiv-are-mesh-crs}.  The theorem may appear to translate a piecewise linear classification problem (concerning regular complexes) into a topological classification problem (concerning meshes), but because of the equivalence of meshes and trusses established in \autoref{intro:thm:truss-mesh-eqv}, classifying refinements of mesh cells is algorithmically decidable via the combinatorics of trusses.  \autoref{fig:the-subdivision-of-a-framed-regular-cell} illustrates a framed subdivision of a framed regular cell, of the sort classified by the previous result.



The combinatorial category of closed trusses is isomorphic, via dualizing the truss posets, to the category of open trusses.  For instance, under this isomorphism the closed 3-truss in \autoref{fig:classifying-framed-regular-cells-intro} corresponds to the open 3-truss in \autoref{fig:a-closed-2-truss-and-an-open-3-truss}.  By the equivalence of trusses and meshes established in \autoref{intro:thm:truss-mesh-eqv}, the self-duality of trusses translates to a self-duality of meshes, as follows.

\begin{introcor}[Dualization of meshes] \label{intro:thm:mesh-dual}
The topological category of closed meshes is weakly equivalent to the topological category of open meshes.
\end{introcor}

\begin{figure}[ht]
    \centering
    \def\svgwidth{1\columnwidth}
    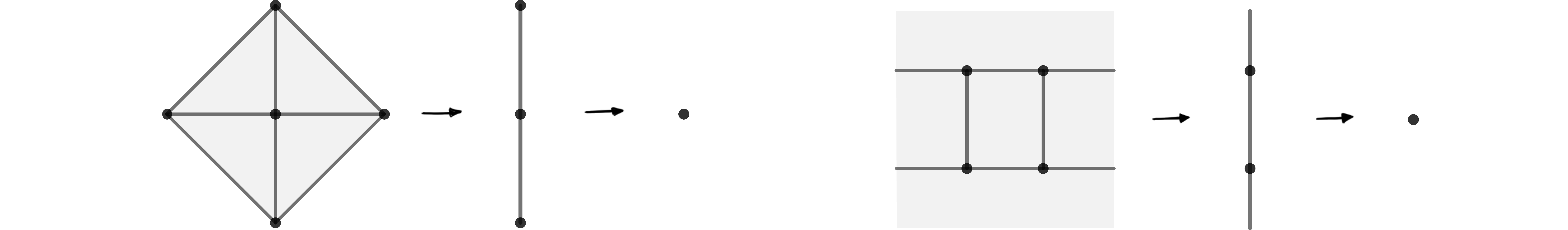

    \caption{Dualization of meshes.}
    \label{fig:dualizability-of-meshes}
\end{figure}

\nid This self-duality appears later as \autoref{thm:dualization_functors}, and is a crucial advantage meshes have over previously known shape categories.  Note that the entrance path poset of the dual $n$-mesh is dual to the entrance path poset of the original mesh, and so the dimensions of all strata in the mesh are dualized---the collection of geometric shapes and their incidences is completely reorganized.  \autoref{fig:dualizability-of-meshes} illustrates a pair of dual 2-meshes.

\pauseae

Meshes, built up from mesh cells and their duals, are a flexible, computationally tractable class of highly structured stratifications; they furthermore provide access to a much broader, almost completely general class of stratifications by considering those stratifications that admit a refinement by a mesh, as follows.

\begin{introdef}[Flat framed stratifications]
A `flat $n$-framed stratification' is a stratification of a subspace of $\lR^n$ that admits a refinement by an $n$-mesh.
\end{introdef}

\begin{figure}[ht]
    \centering
    \def\svgwidth{1\columnwidth}
    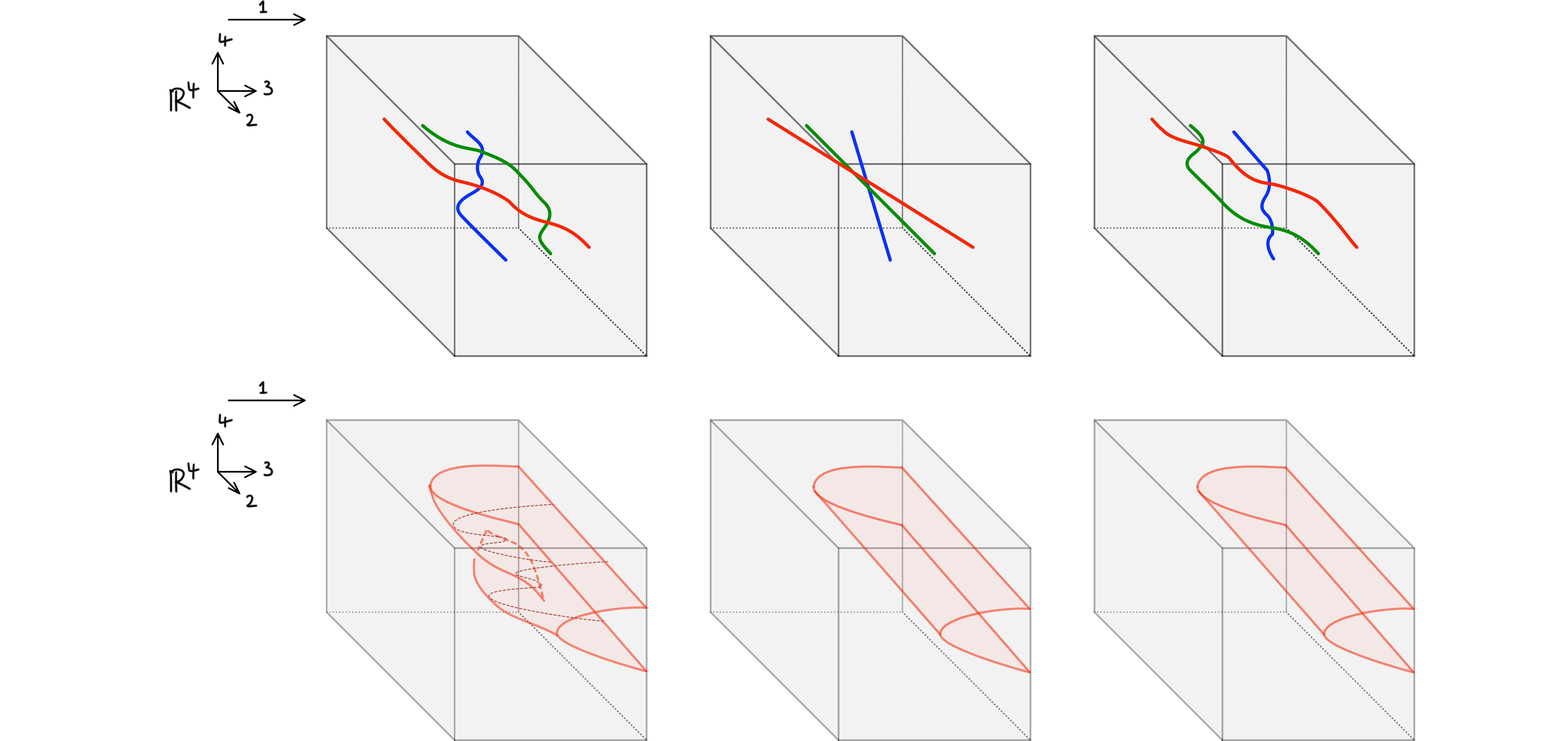

    \caption{Two flat framed stratifications of the 4-cube.}
    \label{fig:yang-baxterator-and-the-swallotail}
\end{figure}

\begin{figure}[ht]
    \centering
    \def\svgwidth{1\columnwidth}
    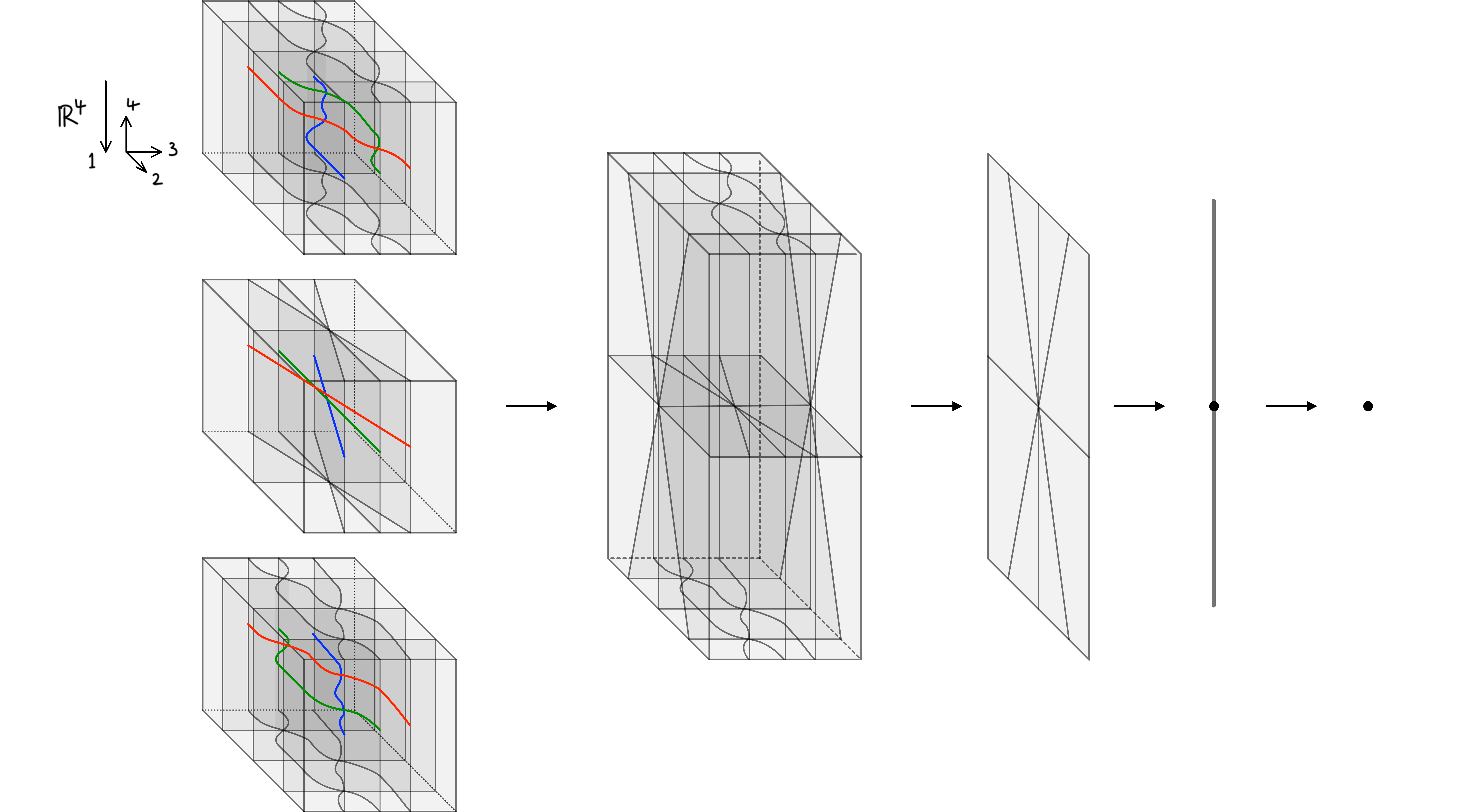

    \caption{The coarsest mesh refining the third Reidemeister move.}
    \label{fig:the-coarsest-mesh-refining-the-reidemeister-iii-move}
\end{figure}

\begin{figure}[ht]
    \centering
    \def\svgwidth{1\columnwidth}
    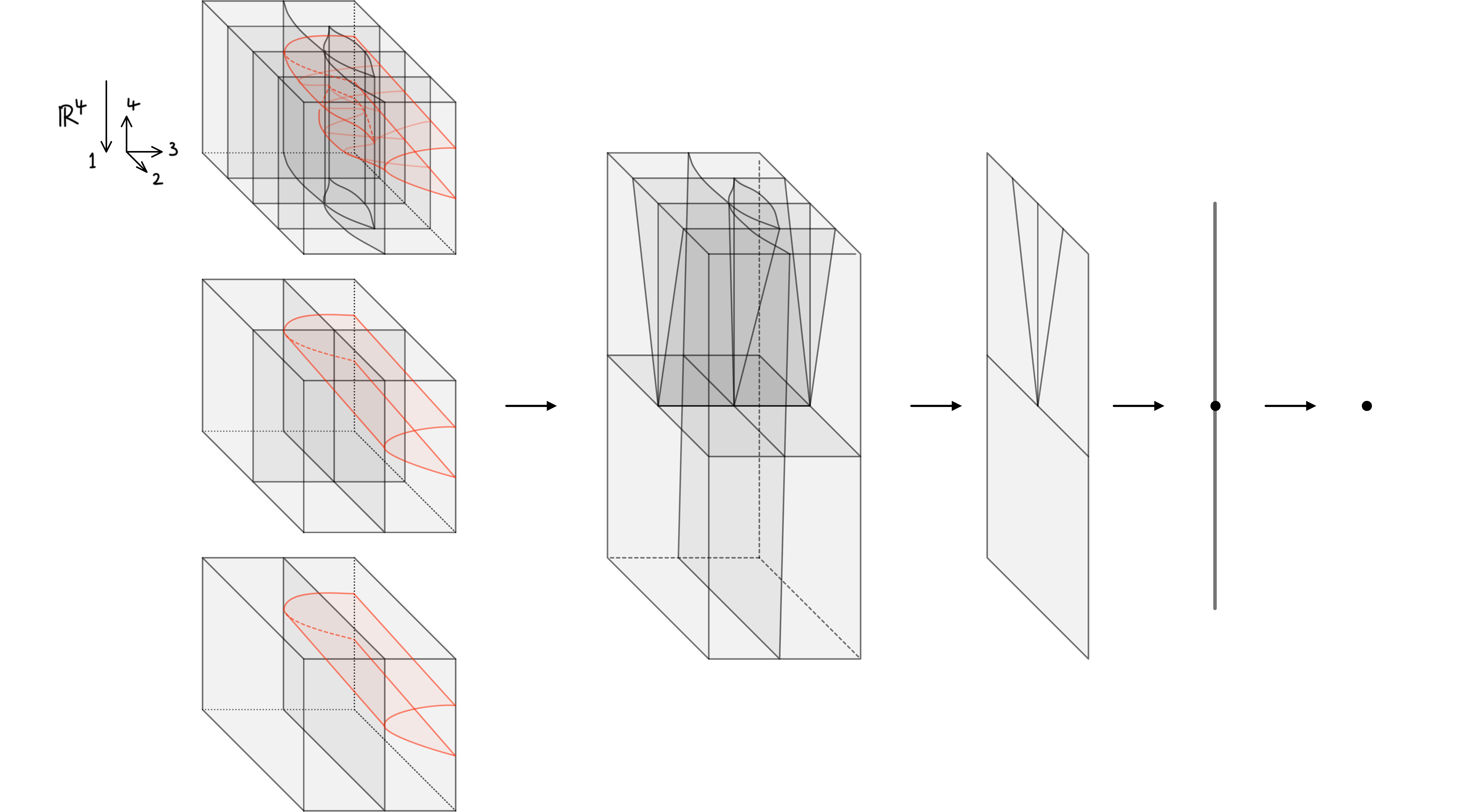

    \caption{The coarsest mesh refining the swallowtail singularity.}
    \label{fig:the-coarsest-mesh-refining-the-swallowtail-singularity}
\end{figure}

\nid This definition will appear in a more precise form in \autoref{defn:flat-framed-strat}.  In \autoref{fig:yang-baxterator-and-the-swallotail} we illustrate two flat framed stratifications of an open 4-cube, by depicting three pertinent slices.  The first stratification is the classical third Reidemeister move, and the second is the classical swallowtail singularity; that these indeed admit mesh refinements is illustrated in a moment in \autoref{fig:the-coarsest-mesh-refining-the-reidemeister-iii-move} and \autoref{fig:the-coarsest-mesh-refining-the-swallowtail-singularity}.

At the heart of the computability of our combinatorial model of stratified framed topology is the completely unexpected fact that among the set of all refining meshes of a flat framed stratification, there is always a canonical \emph{coarsest} choice.  Heuristically, the coarsest refining mesh of a stratification is built up just from the indispensable critical loci of certain projections of the strata.  Needless to say, this situation is in stark contrast to any simplicial model of (stratified) topology, in which two triangulations almost never have a mutual coarsening and typically do not even admit a mutual refinement, preventing any canonical or computable comparison.

\begin{introthm}[Coarsest meshes of flat framed stratifications] \label{intr:thm:coarsest-mesh}
Any flat framed stratification has a unique refining mesh that is coarser than any other refining mesh.
\end{introthm}

\nid This will be established as \autoref{thm:minimal-meshes}.  In \autoref{fig:the-coarsest-mesh-refining-the-reidemeister-iii-move} and \autoref{fig:the-coarsest-mesh-refining-the-swallowtail-singularity} we depict the coarsest mesh for the third Reidemeister move and for the swallowtail singularity.

\pauseae

A flat framed stratification is refined by a mesh, and a mesh is combinatorialized by a truss; to complete the combinatorialization of flat framed stratifications, we translate the initial stratification into a stratified structure on the truss.  A `stratified poset' is a poset together with a `stratification map' to another poset (encoding the set of strata and the combinatorial entrance paths between them); a `stratified truss' is a truss together with a stratification of its total poset.  Furthermore, a stratified truss is `normalized' if it cannot be simplified while preserving the stratification; this property of being normalized corresponds to a mesh being maximally coarsened while still refining a given stratification.

\begin{introthm}[Classification of flat framed stratifications] \label{intro:thm:class-fl-fr-str}
Framed stratified homeomorphism classes of flat framed stratifications are in bijective correspondence with isomorphism classes of normalized stratified trusses.
\end{introthm}

\begin{figure}[ht]
    \centering
    \def\svgwidth{1\columnwidth}
    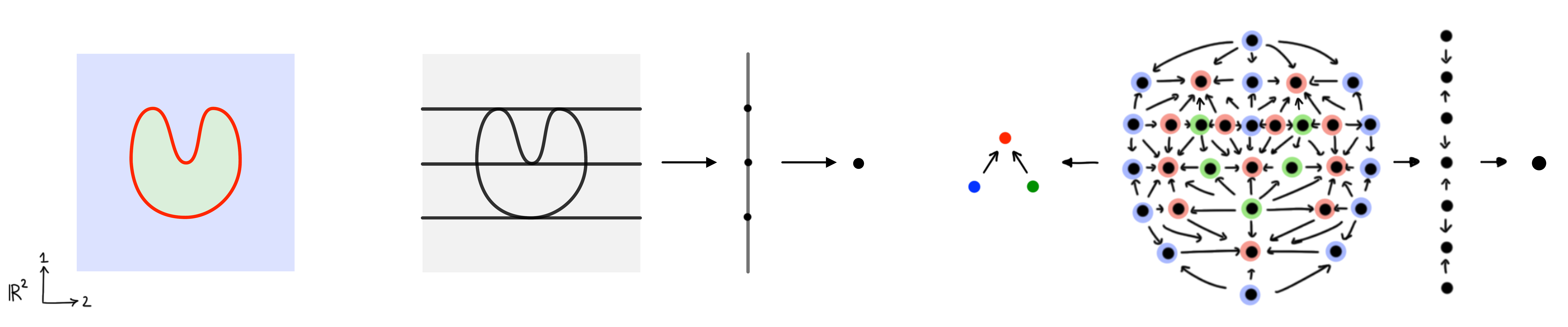

    \caption[A flat framed stratification with its normalized stratified truss.]{A flat framed stratification, its coarsest mesh, and its normalized stratified truss.}
    \label{fig:a-flat-framed-stratification-and-its-corresponding-normalized-stratified-trusses}
\end{figure}

This will be established as \autoref{thm:flat-fr-strat-are-norm-strat-trusses}.  \autoref{fig:a-flat-framed-stratification-and-its-corresponding-normalized-stratified-trusses} illustrates a flat framed stratification, its coarsest refining mesh, and the corresponding normalized stratified truss; the stratification on the truss records which strata of the mesh assemble into each stratum of the initial flat framed stratification.

Flat framed stratifications are intrinsically topological structures, considered up to homeomorphism, while stratified trusses are intrinsically combinatorial or piecewise linear structures, considered up to combinatorial or piecewise linear equivalence---the classification of framed stratifications by stratified trusses thus provides a faithful bridge between the topological and piecewise linear contexts.  Recall the classical, false Hauptvermutung, that homeomorphic simplicial complexes are piecewise linear homeomorphic.  The failure of correspondence between the topological and the piecewise linear remains even for subspaces (or substratifications) of euclidean space: given two piecewise linear embedded triangulated spaces in euclidean space that are ambient homeomorphic, they need not be ambient piecewise linear homeomorphic.  By contrast, in the framed setting we will have a tight correspondence between the topological and piecewise linear, as follows.

\begin{introthm}[Flat framed Hauptvermutung] \label{intro:thm:hauptvermutung}
Flat framed piecewise linear stratifications that are framed stratified homeomorphic are also piecewise linear framed stratified homeomorphic.
\end{introthm}

\nid This result will be established later as \autoref{cor:flat-framed-hpt-vmtg}.  Having a combinatorial or piecewise linear handle on framed stratifications, via stratified trusses, furthermore resolves the fundamental decidability problem for framed homeomorphism.

\begin{introthm}[Decidability of framed stratified homeomorphism] \label{intro:thm:decidable-iso}
Given two flat framed stratifications, one can algorithmically decide whether they are framed stratified homeomorphic.
\end{introthm}
\nid This is recorded later as \autoref{cor:decidability-of-iso}.

One would hope to extend the above computable combinatorialization of framed topological phenomena to framed \emph{smooth} phenomena; indeed we expect this is possible, as follows.

\begin{introconj}[Framed homeomorphism implies diffeomorphism] \label{intro:conj:smooth}
Given two smooth compact manifolds smoothly embedded in euclidean space, and defining flat framed stratifications there, if they are flat framed homeomorphic then they are diffeomorphic.
\end{introconj}

\begin{introconj}[Framed stratifications are dense in smooth embeddings] \label{intro:conj:smooth-2}
Any smooth embedding of a smooth compact manifold into euclidean space has an arbitrarily small perturbation that is a flat framed stratification.
\end{introconj}

\nid These conjectures reappear later as \autoref{conj:smooth-struct} and \autoref{conj:smooth-struct-2}.  Because flat framed stratifications can be faithfully combinatorialized as stratified trusses, these conjectures together imply that there is a sound and complete combinatorial representation of smooth structures on manifolds.



\section{Outline} \label{intro:overview}

\autoref{ch:framed-combinatorial-structures} introduces framed combinatorial structures.  The first such structure, `framed simplices' is a combinatorial analog of classical framed vector spaces.  We then introduce a complementary notion of `proframed simplices' as an analog of a tower of vector space projections.  A `framed simplicial complex' will then be a collection of compatibly framed simplices, while a `proframed simplicial complex' will directly be a tower of projections of simplicial complexes.  We finally generalized framed and proframed simplicial complexes to `framed and proframed regular cell complexes'.

In \autoref{ch:trusses}, we develop our fundamental combinatorial notion of `trusses', as certain iterated constructible bundles of posets.  This development begins with `1-trusses', which are framed fence posets, morphisms between them called `1-truss bordisms', and families of them called `1-truss bundles'.  1-truss bundles over simplices turn out to have an unexpected total order on the top-dimensional simplices in their total posets, and this leads to a crucial method of `truss induction'.  Finally we describe `$n$-trusses', as iterated 1-truss bundles, their corresponding `$n$-truss bordisms' and `$n$-truss bundles', and their elementary constituents `$n$-truss blocks'.

\autoref{ch:classification-of-framed-cells} proves the equivalence of the category of truss blocks and the category of framed regular cells, and more generally the equivalence of the category of regular presheaves on truss blocks and the category of framed regular cell complexes, as stated in \autoref{intro:thm-cell-to-block} and \autoref{intro:thm-cell-to-block-set} above.  Truss blocks are translated into regular cells by an appropriate geometric realization of the total posets of the blocks; the converse translation entails the more technical construction of a tower of 1-truss bundles from the framing information on the regular cell.

In \autoref{ch:meshes}, we introduce our fundamental stratified topological notion of `meshes', as certain iterated constructible bundles of stratified manifolds.  From the outset, meshes appear as a topological analog of the combinatorial structure of trusses, and the notions of `1-mesh', `1-mesh bundle', and `$n$-mesh' parallel the corresponding truss notions.  Indeed, we prove, as claimed in \autoref{intro:thm:truss-mesh-eqv}, that the topological category of meshes is weakly equivalent to the discrete category of trusses; one direction constructs an entrance path truss combinatorially encoding a mesh, and the other direction produces a classifying mesh geometrically realizing a truss.  This fundamental translation allows us to establish the equivalence of mesh cells and framed regular cells (\autoref{intro:cor-meshcells-regularcells}), the classification of subdivisions of framed regular cells (\autoref{intro:thm:subdiv}), and the self-duality of meshes (\autoref{intro:thm:mesh-dual}).

Finally, \autoref{ch:hauptvermutung} radically broadens the class of stratifications under investigation, introducing `flat framed stratifications' as those stratifications that admit a mesh refinement.  The core work of the chapter is the proof that every flat framed stratification has a coarsest refining mesh, as claimed in \autoref{intr:thm:coarsest-mesh}.  We leverage that result to establish the combinatorial classification of flat framed stratification in terms of normalized stratified trusses (\autoref{intro:thm:class-fl-fr-str}).  We then bridge the topological to piecewise linear chasm, proving the flat framed Hauptvermutung (\autoref{intro:thm:hauptvermutung}) that homeomorphic flat framed stratifications are piecewise linear homeomorphic.  As a final application we establish the decidability of framed stratified homeomorphism (\autoref{intro:thm:decidable-iso}).  In the final portion of the chapter, we will describe a future outlook for framed combinatorial topology, including theories of combinatorial tangles and combinatorial singularities, and other aspects of smooth stratified topology.

\autoref{app:frames} provides a detailed discussion of classical linear frames, corresponding notions of indframes and proframes, the generalizations to partial and embedded frames, indframes, and proframes, and the affine analogs of these structures.  \autoref{app:stratifications} reviews and elaborates various elementary notions from stratified topology.

\section{Outlook} \label{intro:outlook}

We briefly summarize our immediate outlook for framed combinatorial topology beyond the present book; a more detailed discussion of these future directions appears in the final \autoref{sec:looking-ahead}.

Flat framed stratifications are already highly, if implicitly, structured by their canonical mesh refinements and the relation of that refinement to the ambient frame on euclidean space.  However, the stratification itself need not be in any sort of generic or transverse relation to the ambient frame; for certain theoretical and computational purposes, it is essential to restrict attention to appropriately transverse stratifications.  We can define and detect a transverse stratification in purely topological, as opposed to smooth, terms, by insisting that every stratum project by a local homeomorphism to the corresponding stage of the ambient euclidean proframe.  We dub these transverse stratifications `manifold diagrams'; indeed our definition provides a solution to the long-standing search for a formal generalization of `string diagrams' to all higher dimensions.  This notion of manifold diagrams is as powerful as one could hope: first, because of the combinatorializability of flat framed stratifications, manifold diagrams are also completely combinatorializable, into an appropriate notion of transverse trusses; second, using the dualization of meshes, manifold diagrams naturally dualize to a notion of higher cell pasting diagrams that formalizes arbitrary composability structures in higher categories.

Leveraging the theory of manifold diagrams, we can identify a combinatorially tractable class of tame tangles, namely those embeddings of manifolds in euclidean space that admit refinements to manifold diagrams.  Of course we expect any embedding of a manifold has an arbitrarily small deformation to a tame tangle, and so nothing is lost by excluding more wild behaviors.  Our combinatorial encoding of (tame) tangles immediately provides a novel computational toolkit: we can stratify the space of tangles by algorithmically computable local or global complexity measures, and formalize computable notions of tangle perturbation, simplification, and stability.  Having a robust algorithmic approach to tangles is already novel in dimension 4, but indeed our definitions and tools apply in all dimensions and all codimensions.  The divergence of our theory from the classical view of tangles becomes especially stark in higher dimensions.  A sufficiently small open neighborhood in a tangle is called a `tangle singularity', or just a `singularity' for brief.  The traditional view has been that singularity classification becomes profoundly unmanageable as the dimension increases: first arise uncountable continuous moduli of distinct singularity types, then the moduli space of singularities itself becomes infinite dimensional, and generally demons abound.  By contrast, we see a natural equivalence relation on singularities (not just tame ones, because we expect we can account for wild ones as well via perturbation methods) for which there is a countable, algorithmically computable classification in all dimensions.

That context of manifold diagrams and tangle singularities considered, there arise various open problems and directions for investigation.  For instance: classify perturbation-stable singularities.  A `perturbation-stable singularity' is one that cannot be simplified by small deformations, and is therefore in a sense an `elementary singularity'; this stability condition is straightforward to formalize using the combinatorial complexity measures at our disposal.  We expect that the set of isomorphism classes of perturbation-stable singularities in any fixed dimension is \emph{finite}, but the structure of this classification remains mysterious as the dimension grows.  Complementary to singularities, which are the most local sort of tangles, are `tangle homotopies', or `homotopies' for brief, which are disconnected tangles encoding the ways manifolds can pass by one another at a distance in euclidean space.  As there are distinguished elementary singularities, namely those that are perturbation-stable, similarly there are `elementary homotopies', namely those that cannot be deformed into a composite of simpler homotopies.  Naturally we may then pose the problem: classify elementary homotopies.  Again, we expect that the set of isomorphism classes of elementary homotopies in any fixed dimension is \emph{finite}, but a precise classification remains unknown even in relatively low dimensions.

Given a sufficiently generic $k$-dimensional tangle $M^k$ in euclidean space $\lR^n$, the composite map $M^k \into \lR^n \to \lR^m$ (where the last map is the standard projection in the canonical proframe of euclidean space) should be a prototypical `$m$-Morse function', in the sense that all its local singularities and global homotopies would be, in an appropriate sense, elementary.  Less precisely than the previous problems, we may ask for the development of a direct definition of $m$-Morse functions (without reference to tangle embeddings), which retains the combinatorial and computational flavor of our tangles and manifold diagrams, and therefore admits a tractable classification and attendant application to smooth manifold topology.  We expect not only that such a combinatorial higher Morse theory exists, but that the resulting combinatorial invariants detect, for instance, all smooth structures on manifolds.  The realization of such an expectation depends, most likely, on the validity of our aforementioned conjectures about framed homeomorphism and framed stratifications---indeed they would imply that every combinatorially tame tangle has a canonical smooth structure and that every smooth tangle has such a combinatorially tame representation.

\addtocontents{toc}{\SkipTocEntry}
\section*{Acknowledgments}

We thank Andr\'e Henriques, Vidit Nanda, Emily Riehl, Mike Shulman, Hiro Tanaka, Kevin Walker, and Mahmoud Zeinalian for helpful and clarifying conversations; Jan Steinebrunner, Filippos Sytilidis, and Christoph Weis for questions and feedback during a preliminary reading group on this material; and David Ayala and Dominic Verity for their encouragement and guidance.  We are grateful to David Reutter and Jamie Vicary for extensive and inspiring discussions concerning combinatorial approaches to homotopy theory.  We thank Ciprian Manolescu and Mark Powell for technical input regarding classical combinatorial topology.

Throughout the research and writing of this book, we were both partially supported by the EPSRC grant EP/S018883/1, ``Higher algebra and quantum protocols", which was conceived and prepared in collaboration with Jamie Vicary and David Reutter.  With thank the Stanford Institute for Theoretical Physics for hosting us for an extended research visit during which some preliminary aspects of this work were envisioned.  CLD was partially supported by a research professorship in the MSRI programs ``Quantum symmetries" and ``Higher categories and categorification" in the spring of 2020, and we both benefited from the hospitality of MSRI that semester, during which core elements of this work took shape.


\renewcommand{\thesection}{\thechapter.{\arabic{section}}}
\renewcommand{\thesubsection}{\thesection.{\arabic{subsection}}}
\renewcommand\thefigure{\thechapter.\arabic{figure}}


\chapter{Framed combinatorial structures} \label{ch:framed-combinatorial-structures}

In this chapter, we introduce notions of framings on classical combinatorial structures.  We start by defining frames on simplices in \autoref{sec:combinatorial-frames}, and introduce core definitions for the `affine combinatorics' of framed simplices. We then generalize our discussion to a notion of framings on simplicial complexes in \autoref{sec:framed-simplicial-complexes}, in which framings of complexes will be pieced together from the frames of each individual simplex of the complex. In a yet further generalization, we will then introduce framings on regular cell complexes in \autoref{sec:framed-reg-complexes}. As it turns out, the generalization from `framed simplices' to `framed regular cells' is less by choice than it is a crucial step in the theory of framed combinatorial topology, on which many later results will depend.

\section{Framed simplices} \label{sec:combinatorial-frames}  \label{par:motivation-for-frames}


The notion of a frame on a simplex, and later on other combinatorial objects, is of course inspired by and modeled on the classical notion of frames. Classically, a `trivialization' of an $m$-dimensional vector space $V$ is specified by a linear isomorphism $V \toiso \lR^m$. Preimages of standard unit vectors $e_i \in \lR^m$ under this trivialization define an ordered list of $(v_1, v_2, ..., v_m)$ in $V$ called a `frame' of $V$.

The guiding intuition in the translation of frames in linear algebra into the combinatorics of simplices, is that directed edges of simplices play the role of vectors. However, $m$-simplices are combinatorially specified by sets of vertices and thus they do not have a distinguished origin. Moreover, their vectors (i.e.\ directed edges) are `affine' in that different vectors may start at different points in the simplex. This observation has important implications for the translation of classical intuition of frames in linear algebra into the combinatorics of simplices, and we will highlight this by speaking of the `affine combinatorics' of framed simplices.

The basic analogy of vectors in a vector space with `vectors' in a simplex will lead to a notion of `frames' on a simplex as follows. We say two vectors in a simplex are `composable' if the endpoint of the first vector is the starting point of the second vector; in this case, their `composite' is the unique directed edge starting at the starting point of the first vector and ending in the endpoint of the second vector. For instance, given a $4$-simplex $S = \{a,b,c,d\}$, the directed edges $d \to b$ and $b \to c$ compose to the directed edge $d \to c$. A `basis' of an $m$-simplex is a set of $m$ vectors such that all other vectors (up to reversing their direction) can be written as composites of vectors in the basis. A `frame' of an $m$-simplex is an ordered basis.

To further explore the analogy to classical linear frames and trivializations, we can rephrase the notion of simplicial frames as follows. Observe that the elements of any basis of an $m$-simplex $S$ must be the elements of a chain of $m$ composable vectors in the simplex; we call such a chain a `spine' of $S$. A choice of basis therefore determines an identification $S \iso [m]$ of $S$ with the \emph{ordered standard simplex} $[m] = (0 < 1 < ... < m)$ (given by mapping the spine of $S$ to the standard spine of $[m]$). Conversely, any such identification $S \iso [m]$ determines a basis in this way. A frame of an $m$-simplex $S$ is then an identification $S \iso [m]$ \emph{together} with a choice of order on the set of standard spine vectors $\spine[m]$ of $[m]$. The standard simplex $[m]$ has of course the canonical identity identification with itself, and so we refer to the simplex $[m]$ with an order $\cF$ of its spine vector set $\spine[m]$ as a `framed standard simplex' $([m],\cF)$.

\begin{figure}[ht]
    \centering
    \def\svgwidth{1\columnwidth}
    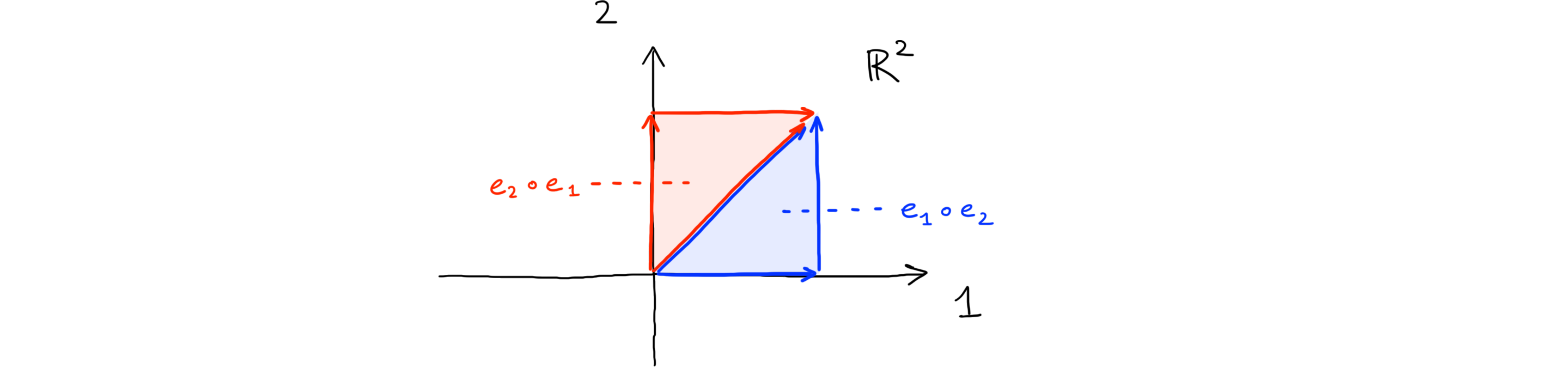

    \caption{Standard vectors spanning distinct standard simplices.}
    \label{fig:two-standard-vectors-can-span-two-different-standard-simplices}
\end{figure}

Framed standard simplices $([m],\cF)$ \emph{with any choice of frame} $\cF$ will play the role of euclidean space $\lR^m$ \emph{with its standard frame} $\{e_1,e_2,...,e_n\}$. The fact that there is not only one `framed standard euclidean space $\lR^m$' but $S_m$-many `framed standard simplices $([m],\cF)$' reflects the fact that several affine constellations of standard basis vectors can arise: for instance, for the standard basis vectors $e_1$ and $e_2$ in $\lR^2$, chaining $e_1$ with $e_2$ forms the spine $e_1 \circ e_2$ of a `standard simplex', whereas chaining $e_2$ with $e_1$ forms the spine $e_2 \circ e_1$ of a different `standard simplex'---see \autoref{fig:two-standard-vectors-can-span-two-different-standard-simplices}. The affine combinatorics of framed simplices accounts for both of these configurations.

We briefly outline this section. In \autoref{ssec:combinatorial-frames} we begin by describing the combinatorial counterparts of classical linear algebraic notions, and then define `frames' on simplices together with several generalizations to `embedded' and `partial' frames. In \autoref{ssec:combinatorial-proframes} we recast these definitions  in terms of so-called `proframes' of simplices. The fact that `frames' and `proframes' are indeed equivalent structures on simplices is discussed in \autoref{ssec:gradients-and-integration}, where we compare the two notions via so-called `gradient' and `integration' functors. Note that the combinatorial notions of `embedded' and `partial' frames (and `proframes') introduced here, may as well be understood in purely classical linear algebraic terms as we explain in \autoref{app:frames}.

\subsection{Frames on simplices} \label{ssec:combinatorial-frames}


Before giving the definition of frames on simplices, we specify some basic terminology and notation.

\begin{term}[Combinatorial simplices]
An `ordered $m$-simplex' is a totally ordered poset with $m+1$ elements. An `$m$-simplex' (for clarity also referred to as an `unordered $m$-simplex') is a set with $m+1$ elements.
\end{term}

\begin{notn}[Category of ordered simplices]
We will denote the category of ordered simplices by $\Delta$; its objects are the ordered simplices, and its morphisms are the order-preserving maps.
\end{notn}

\begin{notn}[Category of unordered simplices] \label{term:standard-unordered-simplices}
We will denote the category of simplices by $\UnSimp$; its objects are the unordered simplices $S$ and its morphisms are all functions.
\end{notn}

\begin{term}[Face maps and degeneracy maps]
An injective map of (ordered or unordered) simplices is called a `face map', and a surjective map of simplices is called a `degeneracy map'.
\end{term}

\begin{term}[Unordering ordered simplices] \label{term:unord-simp-functor} The `unordering' functor $\Unord {(-)} : \Delta \to \UnSimp$ forgets the order of ordered simplices.
\end{term}

\begin{notn}[Maps between ordered and unordered simplices] Maps $S \to T$ or $T \to S$ between an unordered simplex $S$ and an ordered simplex $T$ will be parsed as maps $S \to \Unord T$ resp.\ $\Unord T \to S$.
\end{notn}

\begin{term}[Standard simplex]
The `ordered standard $m$-simplex' $[m]$ is the poset $(0 < 1 < \cdots < m)$; considering posets as categories, this simplex is, equivalently, the category $(0 \ra 1 \ra \cdots \ra m)$.
\end{term}

\nid Since it will be used so frequently, we usually refer to ordered standard $m$-simplex simply as the `$m$-simplex $[m]$'. Every ordered $m$-simplex $S$ is canonically isomorphic to the standard $m$-simplex $[m]$. We may therefore work with a skeleton of $\Delta$ as follows.

\begin{notn}[The skeleton of standard simplices] Abusing notation, we denote the skeleton of $\Delta$ containing only standard simplices $[m]$ for $m \in \lN$ again by $\Delta$.
\end{notn}

\begin{term}[Unordered standard simplices] The unordering $\Unord {[m]}$ of the standard simplex $[m]$ is the `unordered standard $m$-simplex' $\{0,1,...,m\}$.
\end{term}

\begin{term}[Sets of numerals]
The `set of numerals' or `numeral set' $\bnum m$ is the ordered set $\{1 < 2 < \cdots < m\}$.
\end{term}

\begin{notn}[Category of numeral sets]
We will denote the category of numeral sets by $\nabla$; its objects are the numeral sets $\bnum m$ for $m \in \lN$, and its morphisms are the order-preserving functions.
\end{notn}

\pauseae

We next introduce combinatorial notions which mirror ordinary notions from linear and affine algebra. In the most elementary case, this is the correspondence of (nonzero) vectors in a vector space $V$ to (non-degenerate) 1-simplices in a simplex $S$, as follows.

\begin{term}[Simplicial vectors]
    A `vector' $v$ in an (ordered or unordered) simplex $S$ is a map $v : [1] \to S$. We assume by default that $v$ is `non-degenerate' (or `nonzero'), meaning $v$ is injective, and otherwise say that $v$ is a `degenerate' (or `zero') vector. It will be useful to also refer to the set of all zero vectors of $S$, which we call the `affine zero vector' of $S$ and denote by $\simpzero$.
\end{term}

\begin{term}[Spine vectors of ordered simplices] A `spine vector' in an ordered simplex is a non-degenerate vector that cannot written as the composite of more than one non-degenerate vector. The set of spine vectors of the standard simplex is denoted $\spine[m]$; this set $\spine[m] = \{(0 \to 1), (1 \to 2), \ldots, (m-1 \to m)\}$ may be canonically identified with the numeral set $\bnum m$, by mapping the spine vector $(i-1 \to i)$ to the numeral element $i$.
\end{term}

\nid Note that the realization $\abs{S}$ of any simplex $S$ is naturally carries affine structure (see Appendix \ref{appsec:linear-aff} for a brief recollection of affine spaces). Passing to the associated vector spaces of affine spaces, allows to translate `simplicial vectors' to ordinary `linear vectors' as follows.

\begin{notn}[Realizations of vectors in a simplex] \label{notn:realize-vectors} Given a vector $v$ in $S$, we can realize $v$ to an affine vector $\abs{v} = \Delta^1 \to \abs{S}$ in the realized simplex $\abs{S}$. Forgetting the vector's base point, we obtain an ordinary vector $\vec v$ in the associated vector space $\vec V(S)$ of $S$ (that is, $\vec v = \mathsf{unbase}\abs{v}$ using \autoref{term:affine-vectors}).
\end{notn}

The role of `linear projections' of vector spaces will be played by `degeneracies' of simplices. (Indeed, note that any degeneracy $F : S \to T$ of simplices induces a linear projection $\vec V (S) \to \vec V (T)$ of associated vector spaces, defined to map $\vec v$ to $\vec w$ if $w = F(v)$.) Degeneracies admit the following `affine' notion of kernels.

\begin{term}[Affine kernels of simplicial degeneracy maps] \label{term:affine-ker-unord} For a degeneracy $T \to S$, the `affine kernel' $U = \keraff(S \to T)$ is the subset of vectors in $S$ that are mapped to zero vectors in $T$ by the degeneracy.
\end{term}

\nid We write $S \epi S\slash U$ for the simplicial degeneracy of $S$ whose affine kernel is $U$, and $U \intoaff S$ to indicate that $U$ is an affine kernel. Importantly, note that affine kernels $U$ of $S$ cannot, in general, be canonically expressed as simplicial face maps into $S$; that is, the degeneracy $S \epi S/U$ does not have a canonically `splitting' by a simplicial face $S/U \into S$ (i.e.\ such that $(S/U \into S \epi S/U) = \id$).

Conversely, requiring a face $S \into T$ to be split by a degeneracy $T \epi T\slash U$ does not uniquely determine the affine kernel $U$---however, analogous to embeddings into $\lR^n$ being canonically split, this changes when the target simplex $T \iso [n]$ is the standard $n$-simplex as we will now see. Namely, in the presence of order (and thus of spine vectors) it makes sense to adapt the notion of affine kernels as follows.

\begin{term}[Affine kernels for ordered simplices] Given a degeneracy $f : [m] \to [k]$, the `affine kernel' $\keraff(f)$ is the subset of spine vectors in $\spine[m]$ that are mapped to zero vectors by $f$.
\end{term}

\nid Note that, in fact, \emph{any} subset $U \subset \spine[m]$ of spine vectors in $[m]$ determines and is determined by a degeneracy $f: [m] \epi [k]$ with affine kernel $U = \keraff f$.

\begin{term}[Affine cokernels of maps into ordered simplices] Let $S$ be an (unordered) $m$-simplex, and consider an inclusion $S \into [n]$. This induces a unique order $S \iso [m]$ such that $S \iso [m] \into [n]$: note that the (ordered) simplicial face $[m] \into [n]$ is split by a unique degeneracy $[n] \epi [m]$. Thus $S \into [n]$ induces a unique affine kernel $U = \keraff([n] \epi [m])$, called the `affine cokernel' of $S \into [n]$, and written $U = \cokeraff(S \into [n])$.
\end{term}

\begin{term}[Affine faces of simplices] We write $S \intoaff [n]$ for an equivalence class of faces $S \into [n]$ with same affine cokernel. We call $S \intoaff [n]$ an `affine face' (or a `simplicial subspace') and denote the mutual affine cokernel of its representatives by $\cokeraff(S \intoaff [n])$.
\end{term}

\begin{term}[Affine images of maps to ordered simplices] Given an affine face $S \intoaff [n]$ its `affine image' $\imaff(S \intoaff [n])$ is the complement of $\cokeraff(S \intoaff [n])$ in $\spine[n]$.
\end{term}

\nid Note that, in fact, \emph{any} subset $I \subset \spine[m]$ determines and is determined by an affine face $S \intoaff [n]$ with image $I = \imaff(S \intoaff [n])$.

\begin{term}[Splittings of simplicial degeneracies] There is a correspondence of degeneracies $[n] \epi [m]$ and affine faces $[m] \intoaff [n]$ determined by equating their kernels resp.\ cokernels: $\ker([n] \epi [m]) = \cokeraff(S \intoaff [n])$. We say a degeneracy `splits' and `is split by' its corresponding affine face.
\end{term}

\nid The situation reflects the correspondence of linear subspaces $U \into V$ and linear projections $V \epi W$ of vector spaces given by setting $(U \into V) = \ker(V \to W)$ and $(V \epi W) = \coker(U \into V)$. The fact that degeneracies do not correspond to faces but to `affine faces' (i.e.\ equivalence classes of faces) reflects that in the setting of affine spaces canonical choices of kernels do not exist (see \autoref{obs:affine-asymmetry}).

Finally, the notions of `simplicial projections', i.e.\ degeneracies, and of `simplicial subspaces', i.e.\ affine faces, combine to a notion of affine simplicial map as follows.

\begin{term}[Affine simplicial maps] An `affine map' $S \to [n]$ of an $m$-simplex to the standard simplex $[n]$ is a sequence $S \epi T \intoaff [n]$ consisting of a degeneracy $S \epi T$ and an affine face $T \intoaff [n]$. The `affine image' $\im(S \epi T \intoaff [n])$ of such an affine map is the subset of the spine vector set $\spine[m]$ given by the affine image $\imaff(T \intoaff [n])$.
\end{term}

\nid Any simplicial isomorphism $S \iso [m]$ represents the same affine face, and can in particular be considered as an affine map $S = S \iso [m]$. Note also any zero vector $v : [1] \to [m]$ represents the same affine map $[1] \epi [0] \intoaff [m]$, which identify with the affine zero vector $\simpzero$.

\pause

The above combinatorial counterparts to classical linear algebraic notions will allow us to mirror many classical definitions, including that of linear frames, in purely combinatorial terms. Recall a linear frame in an $m$-dimensional vector space $V$ can be specified by a linear trivialization $V \toiso \lR^m$. We will, in fact, also be interested in the following generalizations of linear trivializations, which allow for more general types of maps from $V$ to an euclidean space.

\begin{term}[Partial, embedded, and embedded partial trivializations] A linear projection $V \epi \lR^k$ will be called a `$k$-partial trivialization' (for its relation to `partial' frames of $V$, see Appendix \ref{appsec:linear-fr}). A linear subspace $V \into \lR^n$ will be called an `$n$-embedded trivialization' (for its relation to `embedded' frames of $V$). Yet more generally, a general linear map $V \to \lR^n$ with $k$-dimensional image will be called an `$n$-embedded $k$-partial trivialization' (for its relation to `embedded partial' frames of $V$).
\end{term}

\nid The role of these notions of `generalized trivializations' will be fundamental later on: when defining framed simplicial complexes in \autoref{sec:framed-simplicial-complexes} by patching together global framings from local frames, $n$-embeddings of frames will allow to compare frames of simplices of any dimensions (serving a similar role to local trivializations of $n$-dimensional tangential structure of manifolds).

For the translation of linear (partial, embedded and embedded partial) trivializations into combinatorial definition, there now remains only one central difference: namely, instead of just one `framed standard' model $\lR^n$, we will find an $S_n$-worth of `framed standard' simplices---this is the set of standard $n$-simplices with spine order (that is, a total order $\cF$ of the spine vector set $\spine[n]$). Taking this difference into account, a summary of the translation of linear trivializations (and their generalizations) into affine combinatorial structures is given in \autoref{table:correspondence-of-comb-and-class-frames}. The subsequent sections will give an in-depth discussion of these structures.

\renewcommand{\arraystretch}{1.5}
\begin{figure}[h!]
\centering\small
\begin{tabular}{l | l}
Linear algebra & Affine combinatorics \\
\hline
$m$-dimensional vector space $V$ & unordered $m$-simplex $S$ \\
nonzero vectors $v$ in $V$ & directed edges $v$ in $S$ \\
the zero vector $0$ in $V$ & vertices $x$ in $S$ \\
vector space projections $V \epi V'$ & simplicial degeneracies $S \epi S'$ \\
subspaces $(W \into V)$ & affine faces $T \intoaff S$ \\
framed standard $\lR^m$ & standard simplex $[m]$ with spine order \vspace{+6pt} \\
linear trivialization $V \toiso \lR^m$ & isomorphism $S \iso [m]$ with affine image order \\
partial triv. $V \epi W \toiso \lR^k$ & degeneracy $S \epi T \iso [k]$ with aff. image order \vspace{+6pt}  \\
$n$-embedded triv. $V \into \lR^n$  & affine face $S \intoaff [n]$ with aff. image order \\
$n$-embd.\ partial triv. $V \epi W \into \lR^n$  & affine map $S \epi T \intoaff [n]$ with aff. image order \\
\end{tabular}
\vspace{16pt}
\caption[Affine combinatorics dictionary.]{The analogy between notions in linear algebra and notions in affine combinatorics.}
\label{table:correspondence-of-comb-and-class-frames}
\end{figure}

\subsubsecunnum{The definition of frames}

We introduce frames on simplices, a combinatorial analog of orthonormal frames of euclidean vector spaces, and of linear frames of general vector spaces up to orthoequivalence.  The role of frame vectors will be played by the spine vectors of a simplex.  As a linear frame is an ordered list of its frame vectors, a simplicial frame will be an ordered collection of spine vectors.


\begin{defn}[Frame on the standard simplex]
    A \textbf{frame $\cF$ of the standard $m$-simplex $[m]$} is a bijection $\cF: \spine[m] \rightarrow \bnum{m}$ from the set $\spine[m]$ of spine vectors of the simplex to the set of numerals $\bnum{m} = \{1, 2, \ldots, m\}$.
\end{defn}

\nid We may of course equivalently think of a frame $\cF$ in terms of the inverse function $\cF\inv : \bnum{m} \rightarrow \spine[m]$ from the set of numerals to the spine, or more concretely as an ordered list $(v_1,v_2,...,v_m)$ of spine vectors $v_i = \cF\inv(i)$ of $[m]$. In particular, $\cF$ does provide an order on spine vector in $[m]$ as recorded in \autoref{table:correspondence-of-comb-and-class-frames}. Following the table further, frames of unordered simplices can now be defined as follows.

\begin{defn}[Frame on a simplex] A \textbf{frame} of an $m$-simplex $S$ is an isomorphism $S \iso [m]$ together with a frame $\cF$ on $[m]$.
\end{defn}

\nid We usually denote framed simplices $S$ by tuples $(S \iso [m],\cF)$. We may also keep the isomorphism $S \iso [m]$ implicit, and simply say that $\cF$ is a frame on $S$.

\begin{eg}[Frames on simplices] \label{eg:frames}
    In \autoref{fig:framedsimplices} we illustrate a few framed $m$-simplices $(S \iso [m],\cF)$.  The frame $\cF : \spine[m] \rightarrow \bnum{m}$ is indicated in two ways, as follows.  Firstly, the spine vector $v \in \spine[m]$ is labeled by its numeral value $\cF(v) \in \bnum{m}$.  Secondly, the labeled spine vectors, thought of as vectors in the linear space spanned by the picture of the simplex, are translated so their sources are coincident, and the resulting labeled `coordinate frame' is drawn in or near the simplex.  
\begin{figure}[h!]
    \centering
    \def\svgwidth{1\columnwidth}
    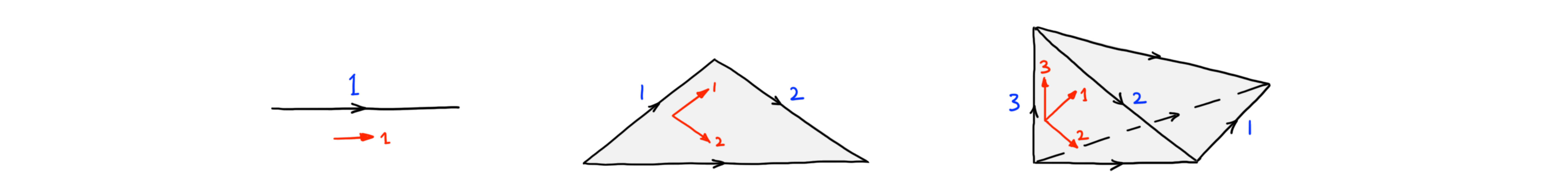

    \caption{Framed simplices.}
    \label{fig:framedsimplices}
\end{figure}
\end{eg}


A linear embedding of a framed simplex into euclidean space may preserve the frame structure in the following sense.

\begin{term}[The standard oriented components of euclidean space] For all $i < n$, the image of the linear subspaces $0^{n-i} \times \lR^i \into 0^{n-i-1} \times \lR^{i+1}$ (with both domain and codomain being subspaces of $\lR^n$) has two components $\eps^-_i$ and $\eps^+_i$ given by $0^{n-i-1} \times \lR_{<0} \times \lR^i$ resp.\ $0^{n-i-1} \times \lR_{>0} \times \lR^i$. We call $\eps^-_i$ and $\eps^+_i$ the `$i$th negative' resp.\ `$i$th positive standard component' of $\lR^n$.
\end{term}

\begin{defn}[Framed realization of a framed simplex] \label{defn:real-fr-simp} A \textbf{framed realization} of a framed $m$-simplex $(S \iso [m],\cF)$ (with frame vectors $v_i = \cF\inv(i)$) is a linear embedding $r_\cF: \abs{S} \into \lR^m$ of the geometric $m$-simplex $\abs{S}$ into $\lR^m$ such that $\vec r_\cF (\vec v_i)$ lies in the $i$th positive standard component $\eps^+_i \subset \lR^m$, for all $i \in \bnum m$.\footnote{Technically, by `linear embedding' $\abs{S} \into \lR^m$ we mean an `affine embedding' in the sense of \autoref{rmk:simplices-as-affine-sp}.}
\end{defn}

\begin{eg}[Framed realization of a framed simplex]
In \autoref{fig:framedrealization} we illustrate framed realizations of two framed 2-simplices. 
\begin{figure}[h!]
    \centering
    \def\svgwidth{1\columnwidth}
    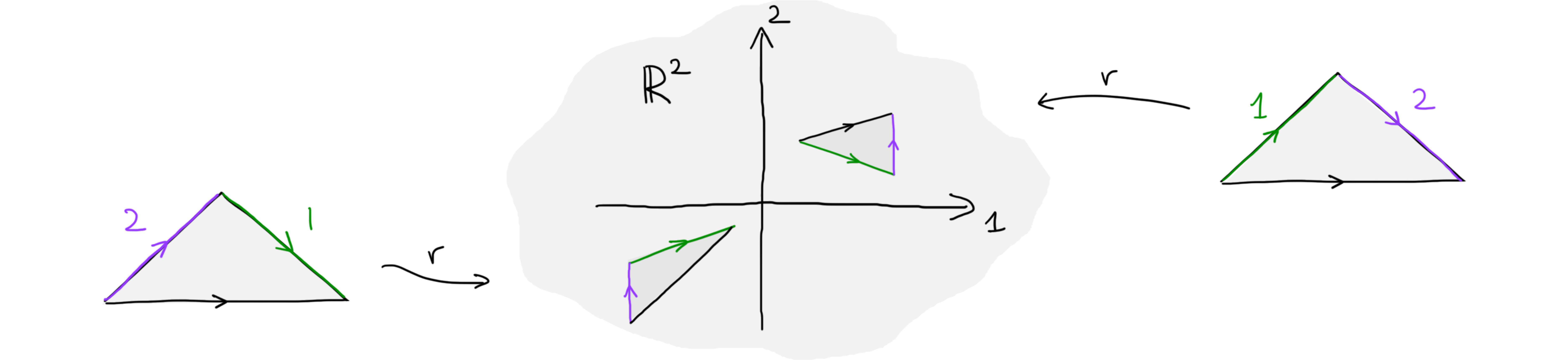

    \caption{Framed realization of framed simplices.}
    \label{fig:framedrealization}
\end{figure}
\end{eg}

\begin{rmk}[Understanding frames in classical linear algebraic terms] The set of all framed realizations of a given framed simplex $(S \iso [m],\cF)$ describes a certain equivalence class of linear trivializations $\abs{S} \into \lR^m$ of the underlying nonframed simplex $S$ (this is formalized in \autoref{rmk:emb-fr-real-orthoeq}).
\end{rmk}



\pauseae

We next mention the following generalization of frames to `partial' frames. Instead of a frame on the `complete' simplex $S$, we may consider a frame on $S$ that is defined everywhere but on a subspace $U \intoaff S$.

\begin{defn}[Partial frame on a simplex] \label{defn:combinatorial-frame} A \textbf{$k$-partial frame} on an $m$-simplex $S$ is a degeneracy $S \epi T$ together with a frame $(T \iso [k],\cF)$.
\end{defn}

\nid Note in particular, that the frame $\cF$ provides an ordering of the affine image of the affine map $S \epi T \iso [k]$ (namely, the image is $\spine[k]$), which matches the structure recorded in \autoref{table:correspondence-of-comb-and-class-frames}. We usually denote partially framed simplices $S$ by tuples $(S \epi T \iso [k],\cF)$ or, more simply, by tuples $(S \epi [k],\cF)$ (note, up to canonical ismorphism, $T$ is redundant). Note that in an $m$-partial frame of an $m$-simplex $(S \epi [m],\cF)$ the degeneracy $S \epi [m]$ must be an isomorphism, and thus $m$-partial frames of $m$-simplices are frames on $m$-simplices.

\begin{term}[Unframed subspace] The `unframed subspace' of a $k$-partially framed simplex $(S \epi [k],\cF)$ is the affine kernel $U = \keraff(S \epi [k])$.
\end{term}


\begin{eg}[Partial frames on simplices] \label{eg:partial-frames}
    In \autoref{fig:combinatorial-frames-and-their-numeric-notation} we illustrate a few $k$-partially framed $m$-simplices $(S \epi [k],\cF)$: we depict degeneracies $S \epi [k]$ by highlighting their unframed subspace (in orange) and illustrate the framed simplices $([k],\cF)$ as in \autoref{eg:frames}. Note that the degeneracy $S \epi [k]$ can equivalently be expressed as a partial order on $S$ (indicated in green), and the frame $\cF$ may then be indicated labeling vectors $v$ in $S$ with $i \in \bnum k$ whenever $w = (S \epi T)(v) \in \spine[k]$ and $\cF(w) = i$ (indicated in purple).
\begin{figure}[h!]
    \centering
    \def\svgwidth{1\columnwidth}
    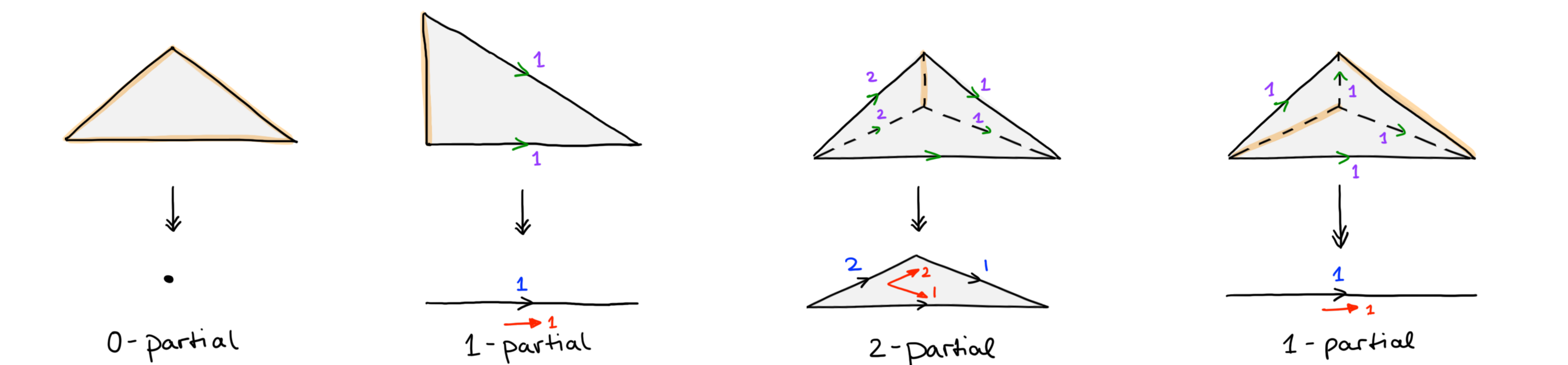

    \caption{Partially framed simplices.}
    \label{fig:combinatorial-frames-and-their-numeric-notation}
\end{figure}
\end{eg}

A linear embedding of a partially framed simplex into euclidean space may preserve the frame structure in the following sense.

\begin{defn}[Framed realization of a partially framed simplex] \label{defn:real-part-fr-simp} Consider a $k$-partially framed $m$-simplex $(S \epi [k],\cF)$ with unframed subspace $U = \keraff(S \epi T)$. A \textbf{framed realization} of $([m],\cF)$ is a linear map $r_\cF: \Delta^m \to \lR^k$ of the standard geometric $m$-simplex $\Delta^m$ to $\lR^k$ such that $\vec r_\cF (\vec v) = 0 \in \lR^k$, for $v \in U$, and such that $\vec r_\cF (\vec v) \in \eps^+_i \subset \lR^k$ whenever $w = (S \epi T)(v) \in \spine[k]$ and $\cF(w) = i$.
\end{defn}

\begin{eg}[Framed realization of a partially framed simplex]
In \autoref{fig:partialframedrealization} we illustrate the framed realizations of a 1-partially framed 2-simplex and of a 2-partially framed 3-simplex.
\begin{figure}[h!]
    \centering
    \def\svgwidth{1\columnwidth}
    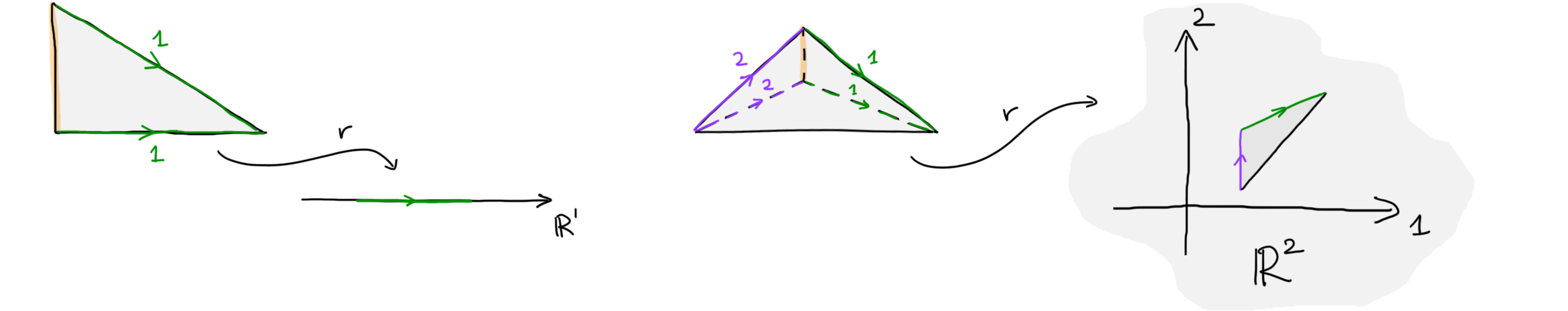

    \caption{Framed realization of partially framed simplices.}
    \label{fig:partialframedrealization}
\end{figure}
\end{eg}

\begin{rmk}[Understanding partial frames in classical affine algebraic terms] The set of all framed realizations of a given $k$-partial framed simplex $(S \epi [k],\cF)$ describes a certain equivalence class of linear trivializations $\abs{S} \to \lR^k$ of the underlying nonframed simplex $S$ (this is formalized in \autoref{rmk:emb-part-fr-real-orthoeq}).
\end{rmk}

\subsubsecunnum{Embedded frames}

We introduce a combinatorial analog of the notion of linear embedded trivializations. Recall, linear embedded trivializations are linear vector space inclusions $V \into \lR^n$ which, if $V$ is euclidean, up to orthoequivalence are represented by orthonormal $n$-embedded frames on $V$, that is, by ordered lists of vectors $(v_1,v_2,...,v_n)$ of which exactly $\dim(V)$ are nonzero and form an orthonormal basis of $V$.


\begin{defn}[Embedded frame on the standard simplex] \label{defn:combinatorial-frame-embedded}
    An \textbf{$n$-embedded frame $\cF$ of the standard $m$-simplex $[m]$} is an injective function $\cF: \spine[m] \into \bnum{n}$ from the spine of the simplex into the set of numerals $\{1,2,\ldots,n\}$.
\end{defn}

\nid We may of course equivalently think of an $n$-embedded frame $\cF$ of $[m]$ in terms of the partial inverse function $\cF\inv : \bnum n \to \spine[m]$ from the set of numerals to the spine, or more concretely as an ordered list $(v_1,v_2,...,v_n)$ where $v_i = \simpzero$ if $\cF\inv$ is undefined on $i$, and $v_i = \cF\inv(i)$ otherwise. Note that an $m$-embedded framed $m$-simplex is the same as a framed $m$-simplex as previously defined. For general simplices, we introduce the following.

\begin{defn}[Embedded frame on a simplex] An \textbf{$n$-embedded frame} of an $m$-simplex $S$ is an isomorphism $S \iso [m]$ together with an $n$-embedded frame $\cF$ on $[m]$.
\end{defn}

\nid We usually denote $n$-embedded framed simplices $S$ by tuples $(S \iso [n],\cF)$. In line with \autoref{table:correspondence-of-comb-and-class-frames}, we remark that there is the following equivalent definition of embedded frames.

\begin{rmk}[Embedded frames via simplicial subspaces] \label{rmk:emb-fr-via-simp-subspaces} An $n$-embedded frame $\cF$ of an $m$-simplex $S$ is equivalently given by an affine face $S \intoaff [n]$ together with an ordering of the affine image $I \subset \spine[n]$ of $S \intoaff [n]$. The two definitions define equivalent structures on $S$, since $I \subset \spine[n]$ also inherits a total order from $\spine[n]$ which allows us to identify $I \iso \spine[m]$ canonically.
\end{rmk}

\begin{eg}[Embedded frames on simplices] \label{eg:combinatorial-frames-emb}
    In \autoref{fig:embedded-combinatorial-frames-and-their-embedded-coordinate-notation} we illustrate a few $n$-embedded framed $m$-simplices $(S \iso [m],\cF)$. As before, the frame $\cF: \spine[m] \into \bnum{n}$ of $[m]$ is indicated in two ways: the spine vector $v \in \spine[m]$ is labeled by its numeral value $\cF(v) \in \bnum{n}$, and the $m$ labeled spine vectors are translated into a labeled coordinate frame.
\begin{figure}[h!]
    \centering
    \def\svgwidth{1\columnwidth}
    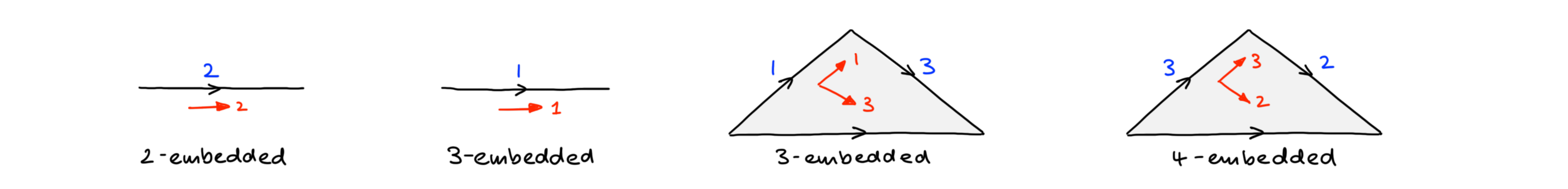

    \caption{Embedded framed simplices.}
    \label{fig:embedded-combinatorial-frames-and-their-embedded-coordinate-notation}
\end{figure}
\end{eg}

A linear embedding of a framed simplex into euclidean space may preserve the frame structure in the following sense.

\begin{defn}[Framed realization of an embedded framed simplex] \label{defn:framed-real-emb}  A \textbf{framed realization} of an $n$-embedded framed $m$-simplex $(S \iso [m],\cF)$ (with nonzero frame vectors $v_i = \cF\inv(i)$, $i \in \im(\cF)$) is a linear embedding $r_\cF: \abs{S} \into \lR^n$ of the geometric $m$-simplex $\abs{S}$ into $\lR^n$ such that $\vec r_\cF (\vec v_i) \in \eps^+_i \subset \lR^m$, for all $i \in \im(\cF)$.
\end{defn}

\begin{eg}[Framed realization of embedded framed simplices]
    In \autoref{fig:embframedrealization} we illustrate framed realizations of two $3$-embedded framed simplices; note that the framed realizations of the $3$-embedded framed 1-simplex (shown on the right) must yield affine vectors that, after translating their base points to the origin, lie in the component $\eps^+_i$ (which consists of vectors in the $\avg{e_2,e_3}$-plane with strictly positive $2$-component). 
\begin{figure}[h!]
    \centering
    \def\svgwidth{1\columnwidth}
    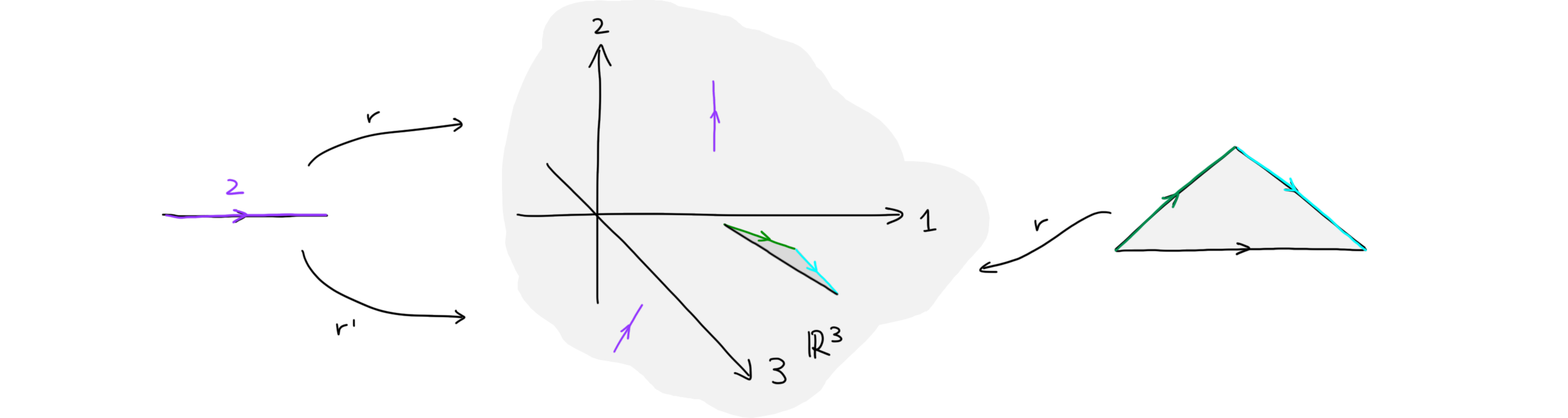

    \caption{Framed realization of embedded framed simplices.}
    \label{fig:embframedrealization}
\end{figure}
\end{eg}

\begin{rmk}[Understanding embedded frames in classical affine algebraic terms] The set of all framed realizations of a given $n$-embedded framed simplex $(S \epi [m],\cF)$ describes a certain equivalence class of linear trivializations $\abs{S} \into \lR^n$ of the underlying nonframed simplex $S$ (this is formalized in \autoref{rmk:emb-fr-real-orthoeq}).
\end{rmk}

\pauseae

While we will be interested ultimately only in non-partial embedded frames, we also mention the following generalization of embedded frames to the partial case. Instead of an embedded frame on the `complete' simplex $S$, we may consider an embedded frame on $S$ that is defined everywhere but on a subspace $U \intoaff S$. Following \autoref{table:correspondence-of-comb-and-class-frames} embedded partial frames of unordered simplices can be defined as follows.

\begin{defn}[Embedded partial frame on a simplex] \label{defn:combinatorial-frame-emb} An \textbf{$n$-embedded $k$-partial frame} on an $m$-simplex $S$ is a degeneracy $S \epi T$ together with an $n$-embedded frame $(T \iso [k],\cF)$.
\end{defn}

\nid We usually omit the (canonically determined) simplex $T$, and denote $n$-embedded $k$-partially framed simplices $S$ by tuples $(S \epi [k],\cF)$ (where $\cF$ is an $n$-embedded frame of the standard simplex $[k]$). Note, in line with \autoref{table:correspondence-of-comb-and-class-frames}, we may equivalently think of the structure of an embedded partially framed simplex as an affine map $S \epi T \intoaff [k]$ together with an ordering of its affine image (see \autoref{rmk:emb-fr-via-simp-subspaces}).

\begin{term}[Unframed subspace] The `unframed subspace' of an $n$-embedded partially framed simplex $(S \epi [k],\cF)$ is the affine kernel $U = \keraff(S \to T)$.
\end{term}

\begin{eg}[Embedded partial frames on simplices]
    In \autoref{fig:embedded-partial-combinatorial-frames} we illustrate a few $n$-embedded $k$-partially framed $m$-simplices $(S \epi [k],\cF)$: we depict degeneracies $S \epi [k]$ by highlighting their unframed subspace (in orange) and illustrate the framed simplices $([k],\cF)$ as in \autoref{eg:combinatorial-frames-emb}. Note that the degeneracy $S \epi [k]$ can equivalently be expressed as a partial order on $S$ (indicated in green), and the frame $\cF$ may then be indicated labeling vectors $v$ in $S$ with $i \in \bnum n$ whenever $w = (S \epi T)(v) \in \spine[k]$ and $\cF(w) = i$ (indicated in purple).
\begin{figure}[ht]
    \centering
    \def\svgwidth{1\columnwidth}
    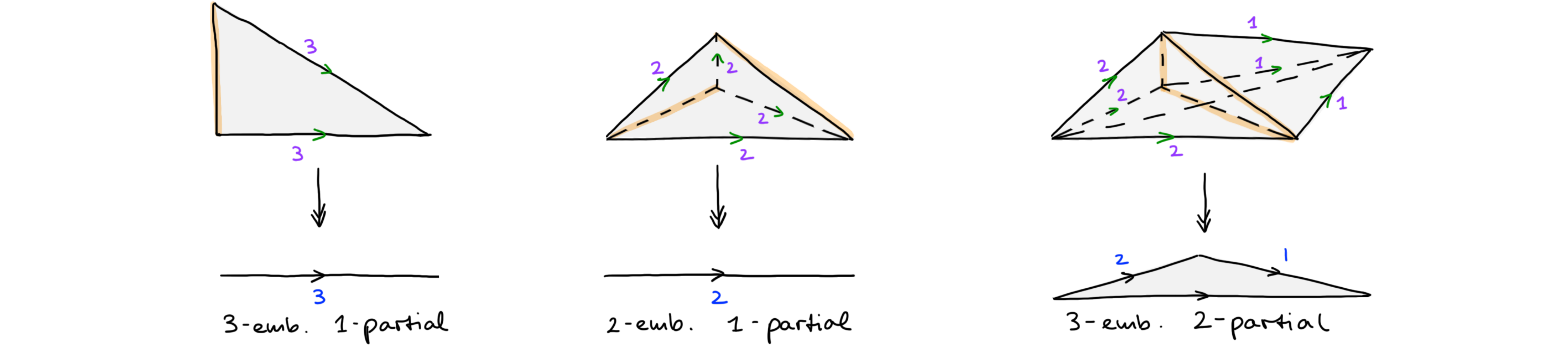

    \caption{Embedded partially framed simplices.}
    \label{fig:embedded-partial-combinatorial-frames}
\end{figure}
\end{eg}

A linear embedding of a embedded partially framed simplex into euclidean space may preserve the frame structure in the following sense.

\begin{defn}[Framed realization of a embedded partially framed simplex] \label{defn:framed-real-emb-part}
Consider a embedded partially framed $m$-simplex $(S \epi [k],\cF)$ with unframed subspace $U = \keraff(S \epi T)$. A \textbf{framed realization} of $([m],\cF)$ is a linear embedding $r_\cF: \Delta^m \into \lR^m$ of the standard geometric $m$-simplex $\Delta^m$ into $\lR^m$ such that $\vec r_\cF (\vec v) = 0 \in \lR^k$, for $v \in U$, and such that $\vec r_\cF (\vec v) \in \eps^+_i \subset \lR^k$ if $w = (S \epi T)(v) \in \spine[k]$ and $\cF(w) = i$.
\end{defn}

\begin{eg}[Framed realization of an embedded partially framed simplex]
    In \autoref{fig:embeddedpartiallyframedrealization} we illustrate framed realizations of a $2$-embedded $1$-partial frame of the $3$-simplex and of a $3$-embedded $2$-partial frame of the $4$-simplex (for the latter example, we think of the central $\lR^2$ plane as being embedded in $\lR^3$ as indicated by the grey coordinate axis).  
\begin{figure}[h!]
    \centering
    \def\svgwidth{1\columnwidth}
    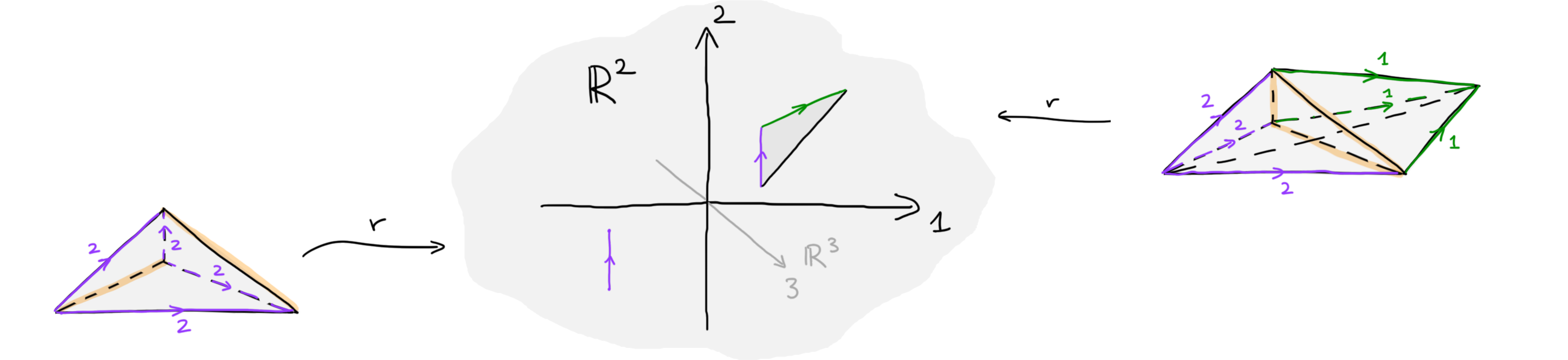

    \caption{Framed realization of an embedded partially framed simplex.}
    \label{fig:embeddedpartiallyframedrealization}
\end{figure}
\end{eg}

\begin{rmk}[Understanding embedded partial frames in classical affine algebraic terms] The set of all framed realizations of a given $n$-embedded $k$-partial framed simplex $(S \epi [k],\cF)$ describes a certain equivalence class of linear trivializations $\abs{S} \to \lR^n$ of the underlying nonframed simplex $S$ (this is formalized in \autoref{rmk:emb-part-fr-real-orthoeq}).
\end{rmk}

\subsubsecunnum{Restricting frames}

A linear trivialization $V \toiso \lR^n$ restricts on any linear subspace $W \into V$ to a linear $n$-embedded trivialization $W \into \lR^n$. We discuss the combinatorial analog of this process, that is, how frames and embedded frames of simplices restrict to embedded frames on simplicial faces. In geometric terms, our definition can be expressed as follows.

\begin{rmk}[Restricting frames in linear algebraic terms] Consider an $n$-embedded framed $m$-simplex $(S \iso [m],\cF)$, and a simplicial $j$-face $f : T \into S$. Pick any framed realization $r : \abs{S} \into \lR^n$. The `restriction' of the $n$-embedded frame of $S$ to $T$ is the unique $n$-embedded frame of $T$ which is framed realized by the linear embedding $r \circ \abs{f} : \abs{T} \into \abs{S} \into \lR^n$.
\end{rmk}

Describing the process of frame restriction in purely combinatorial terms, without reference to the affine framed structure of $\lR^n$, is more subtle.  To properly account for the combinatorial situation, we introduce a combinatorial notion of `non-orthogonality' of simplicial vectors that we refer to as `kinship'.


\begin{defn}[Akin simplicial vectors]
The vectors $v = (a \ra b)$ and $w = (c \ra d)$ in the simplex $[m] \equiv (0 \ra 1 \ra \cdots \ra m)$ are \textbf{akin}, denoted $v \akin w$, if there is a vector $u$ that is a factor of both, i.e.\ such that $v = \tilde{v} \circ u \circ \vardbtilde{v}$ and $w = \tilde{w} \circ u \circ \vardbtilde{w}$ for some possibly degenerate vectors $\tilde{v}, \vardbtilde{v}, \tilde{w}, \vardbtilde{w}$.
\end{defn}

\nid Note that, like the relation of non-orthogonality of linear vectors in a euclidean space, the `kinship' relation between simplicial vectors is reflexive and symmetric but not transitive.

\begin{eg}[Kinship of vectors]
In \autoref{fig:the-join-of-two-1-simplices} we illustrate the kinship of vectors in the 3-simplex.  To emphasize the informal conceptual relationship of this notion with the geometry of non-orthogonality, the simplex is drawn with its three spine vectors (highlighted in red, green, and blue) being orthogonal in the ambient euclidean 3-space.  The three red vectors are akin, the four green vectors are akin, and the three blue vectors are akin; no other vectors are akin.  Note that indeed these kinship relations are precisely the non-orthogonal vectors of this geometric 3-simplex.
\begin{figure}[ht]
    \centering
    \def\svgwidth{1\columnwidth}
    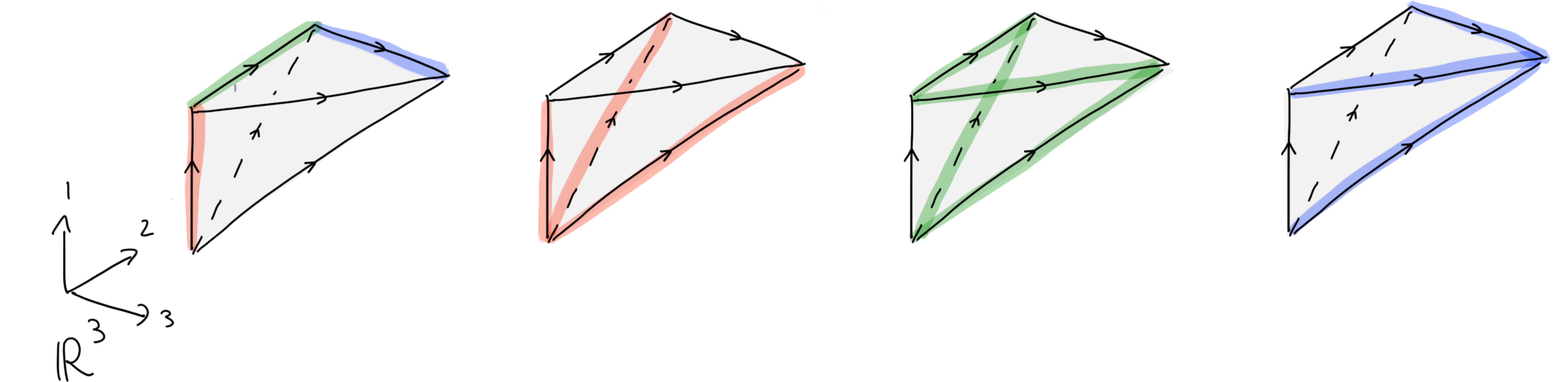

    \caption{Kinship of simplicial vectors.}
    \label{fig:the-join-of-two-1-simplices}
\end{figure}
\end{eg}

Using the notion of kinship, we can describe the restriction of a frame to any vector of the standard simplex, and subsequently to any face of the standard simplex.  The frame order on the spine of the simplex plays a paramount role in this process: the frame label of a general vector will be the lowest numeral among the frame labels of the spine vectors akin to the given vector.

\begin{constr}[Frame restriction to simplicial vectors] \label{constr:frame-inheritance-on-1-faces}
Given an $n$-embedded frame $\cF: \spine[m] \hookrightarrow \bnum{n}$ of the simplex $[m]$, and a vector $v: [1] \to [m]$ of that simplex, the `restriction' $\rest \cF v: \spine[1] \hookrightarrow \bnum{n}$ of the frame to the vector is the $n$-embedded frame of the simplex $[1]$ whose single label is the minimal frame label of the spine vectors akin to the vector $v$, i.e.\ $\rest \cF v (0 \to 1) = \min\{\cF(w) \,|\, w \akin v\}$.
\end{constr}

\nid This restriction procedure produces a plethora of combinatorial arrangements quite distinct from any permutation of its application to the standard frame on the simplex.

\begin{eg}[Frame restriction to simplicial vectors of the standard simplex]
In \autoref{fig:framed-1-faces} we illustrate various embedded framed 3-simplices along with the corresponding embedded framed restrictions to their 1-faces.
\begin{figure}[ht]
    \centering
    \def\svgwidth{1\columnwidth}
    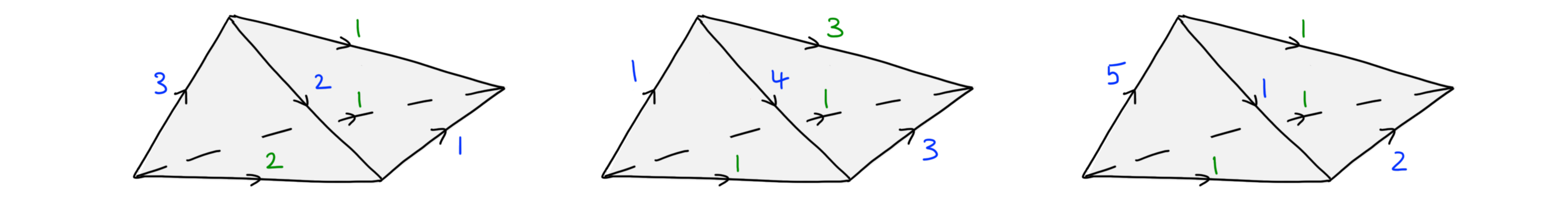

    \caption{Restriction of embedded frames to vectors of a simplex.}
    \label{fig:framed-1-faces}
\end{figure}
\end{eg}

The restriction of an embedded frame to any $j$-face of the simplex is determined directly by the restrictions to the 1-faces, as follows.

\begin{constr}[Frame restriction to simplicial faces of the standard simplex] \label{constr:inherited-frames-on-faces}
Let $\cF: \spine[m] \hookrightarrow \bnum{n}$ be an $n$-embedded frame of the simplex $[m]$, and let $f: [j] \to [m]$ be a $j$-face of that simplex.  The `restriction' $\rest \cF f: \spine[j] \hookrightarrow \bnum{n}$ of the frame to the $j$-face is the $n$-embedded frame whose label value on the spine vector $v: [1] \to [j]$ is the numeral $\rest \cF {f \circ v} \in \bnum{n}$.
\end{constr}

\nid The restriction of frames on standard simplices canonically carries over to the case of general simplices.

\begin{notn}[Restricting trivializations to faces] \label{notn:restricting-triv-to-faces} Given a degeneracy $S \epi [k]$ from an $m$-simplex $S$ to a standard simplex, and a $j$-face $f : T \into S$, then we denote by $\rest {(S \epi [k])} f : T \iso [j]$ the unique degeneracy as well as (abusing notation) by $f : [j] \into [k]$ the unique face, which together factor the composite $(S \epi [k]) \circ f$ as $f \circ \rest {(S \epi [k])} f$. (We call $\rest {(S \epi [k])} f$ the `restricted degeneracy', or the `restricted isomorphism' if $S \epi [k]$ is an isomorphism, and $f : [j] \into [k]$ the `induced standard face').
\end{notn}

\begin{defn}[Frame restrictions to simplicial faces of simplices] \label{defn:frame-restriction-gen} The \textbf{frame restriction} of an $n$-embedded framed $m$-simplex $(S \iso [m],\cF)$ to a simplicial $j$-face $f : T \into S$ is the $n$-embedded frame $(T \iso [j],\rest \cF f)$ of $T$ where $T \iso [j]$ equals the restricted isomorphism $\rest {(S \iso [k])} f$ and the $n$-embedded frame $\rest \cF f$ is obtained by restricting $\cF$ to the induced standard face $f : [j] \into [m]$.
\end{defn}

\begin{eg}[Frame restriction to simplicial faces] \label{eg:inherited-faces-of-framed-simplices}
In \autoref{fig:inherited-frames-on-general-faces} we depict an embedded framed 3-simplex along with the restriction of its frame to various faces.
\begin{figure}[ht]
    \centering
    \def\svgwidth{1\columnwidth}
    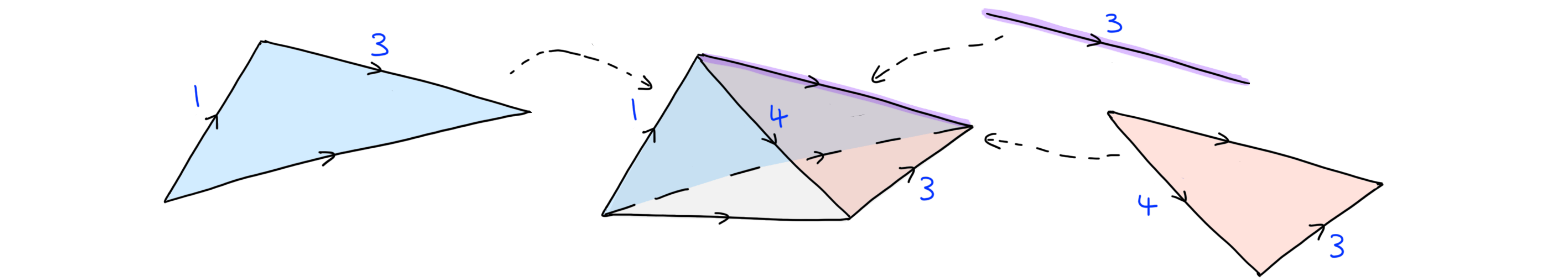

    \caption{Restriction of an embedded frame to faces of a simplex.}
    \label{fig:inherited-frames-on-general-faces}
\end{figure}
\end{eg}

\begin{obs}[Frame restrictions are well-defined on simplicial subspaces] Given $n$-embedded framed $m$-simplex $(S \iso [m],\cF)$, then for any two faces $f,f' : T \into S$ representing the same affine face $T \intoaff S \iso [m]$ we have $\rest \cF f = \rest \cF {f'}$. In other words, frame restriction is well-defined on affine faces, and we may write $\rest \cF {T \intoaff S}$ for the restriction of $(S \iso [m],\cF)$ to an affine face.
\end{obs}

One may of course similarly define frame restrictions of embedded \emph{partial} frames of simplices.

\begin{defn}[Frame restrictions to simplicial faces of simplices] The \textbf{frame restriction} of an $n$-embedded $k$-partially framed $m$-simplex $(S \epi [k],\cF)$ to a $j$-face $f : T \into S$ is the $n$-embedded partial frame $(T \epi [l],\rest \cF f)$ of $T$ where $T \epi [l]$ equals the restricted degeneracy $\rest {(S \iso [k])} f$ and the $n$-embedded frame $\rest \cF f$ is obtained by restricting $\cF$ to the induced standard face $f : [l] \into [k]$.
\end{defn}

The defining property of such restrictions is, in geometric terms, again the following: pick any framed realization $r : \abs{S} \into \lR^n$. The frame restriction of the given $n$-embedded partial frame of $S$ to the face $f : T \into S$ is the unique $n$-embedded partial frame of $T$ induced by the simplicial embedded partial trivialization $r \circ \abs{f} : \abs{T} \into \abs{S} \into \lR^n$.

\subsubsecunnum{Framed maps}

As previously described subspaces $W \into V$ of embedded trivialized vector spaces $V \into \lR^n$ inherit themselves an embedded trivialization $W \into \lR^n$ by restricting the trivialization of $V$. The inclusion $(W \into \lR^n) \into (V \into \lR^n)$ is a prototypical `trivialization inclusion'. Conversely, one can define a `trivialization projection' $(V \into \lR^n) \epi (W \into \lR^n)$ to be a projection $V \epi W$ that splits a framed trivialized inclusion. Combinatorially, this can be translated as follows.

\begin{defn}[Framed faces] Given $n$-embedded framed simplices $(S \iso [l],\cF)$ and $(T \iso [m],\cG)$, \textbf{framed face} $F : (S \iso [l],\cF) \into (T \iso [m],\cG)$ is a face $F : [l] \into [m]$ such that $\cF = \rest \cG F$.
\end{defn}

\begin{defn}[Framed degeneracies] Given $n$-embedded framed simplices $(S \iso [l],\cF)$ and $(T \iso [m],\cG)$, a \textbf{framed degeneracy} $F : (S \iso [l],\cF) \epi (T \iso [m],\cG)$ is a degeneracy $F : [l] \epi [m]$ such that $\cG = \rest \cF {T \intoaff S}$ (where $T \intoaff S$ splits $F$).
\end{defn}

\nid Framed faces and framed degeneracies are the monomorphisms resp.\ epimorphisms in a category of `framed simplices and framed maps'; such framed maps may be uniformly described as follows.

\begin{defn}[Framed map of framed simplices] \label{defn:framed-map-of-framed-simplices}
    Given $n$-embedded framed simplices $(S \iso [l],\cF)$ and $(T \iso [m],\cG)$, a \textbf{framed map} $F : (S \iso [l],\cF) \to (T \iso [m],\cG)$ is a simplicial map $F: [l] \to [m]$ such that for every vector $v: [1] \to [l]$ in the simplex $[l]$, either its frame label is preserved, i.e.\ $\rest \cF v = \rest \cG {F \circ v}$, or the vector is degenerated, i.e.\ $F \circ v: [1] \to [m]$ is constant.
\end{defn}

\nid Note, injective framed maps are exactly framed faces, and a surjective framed maps are exactly framed degeneracies.

\begin{notn}[Categories of embedded framed framed simplices]
Denote by $\FrDelta n$ the category of $n$-embedded framed simplices and their framed maps. (Note that the objects of this category are simplices of dimension necessarily at most $n$.)
\end{notn}

\begin{obs}[(Epi,mono)-factorization of framed maps] \label{obs:framed-epi-mono-fact} As a map of simplices factors as a degeneracy map (epimorphism) followed by a face map (monomorphism), similarly any framed map of framed simplices factors as a framed degeneracy map followed by a framed face map.
\end{obs}

\begin{eg}[Framed and non-framed maps] \label{eg:framed-and-lax-framed-maps-of-simplices}
In \autoref{fig:framed-and-lax-framed-maps} we illustrate four maps between 2-embedded framed simplices; in each case the map is a face map (with highlighted image), an identity, or a standard simplicial degeneracy (with highlighted affine kernel).  The first map is framed, as the one frame label is preserved.  The second map is not framed, as it preserves the frame label of neither spine vector.  The third map is framed, as the first spine vector is degenerated and the other two vectors have their frame labels preserved.  The fourth map is not framed, even though every spine vector is degenerated or has its frame label preserved, because the frame label of the non-spine vector is not preserved.
\begin{figure}[ht]
    \centering
    \def\svgwidth{1\columnwidth}
    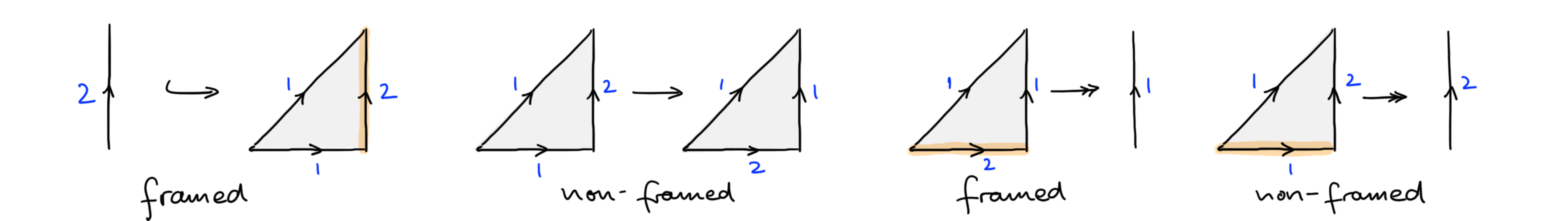

    \caption{Framed and non-framed maps of framed simplices.}
    \label{fig:framed-and-lax-framed-maps}
\end{figure}
\end{eg}

\nid Framed maps of $n$-embedded simplices can also be understood in geometric terms.

\begin{rmk}[Framed maps in linear algebraic terms] \label{rmk:framed-map-geo} Consider $n$-embedded framed simplices $(S \iso [l],\cF)$ and $(T \iso [m],\cG)$ and a simplicial map $F : S \to T$. For framed realizations $r : \abs{S} \into \lR^n$ and $q : \abs{T} \into \lR^n$, the map $F$ is framed if for any vector $v$ in $S$ with $\vec r (\vec v) \in \eps^+_i$ we have $\vec q \circ \vec F (\vec v) \in \eps^+_i \cup \{0\}$.
\end{rmk}

\begin{eg}[Framed maps in linear algebraic terms] In \autoref{fig:geometric-interpretation-of-framed-maps}, for the four framed maps of $2$-embedded framed simplices $F : (S \iso [l],\cF) \to (T \iso [m],\cG)$ that were given in \autoref{fig:framed-and-lax-framed-maps}, we depict framed realizations $r : \abs{S} \into \lR^2$ and $q : \abs{T} \into \lR^2$ of their domain and codomains. Underneath, we depict the linear action $\vec F : \vec V (S) \to \vec V (T)$ on vectors $v$ in $S$ (we depict vectors by their images under the embeddings $\vec V(S) \into \lR^2$ resp.\ $\vec V(T) \into \lR^2$).
\begin{figure}[ht]
\centering
    \def\svgwidth{1\columnwidth}
    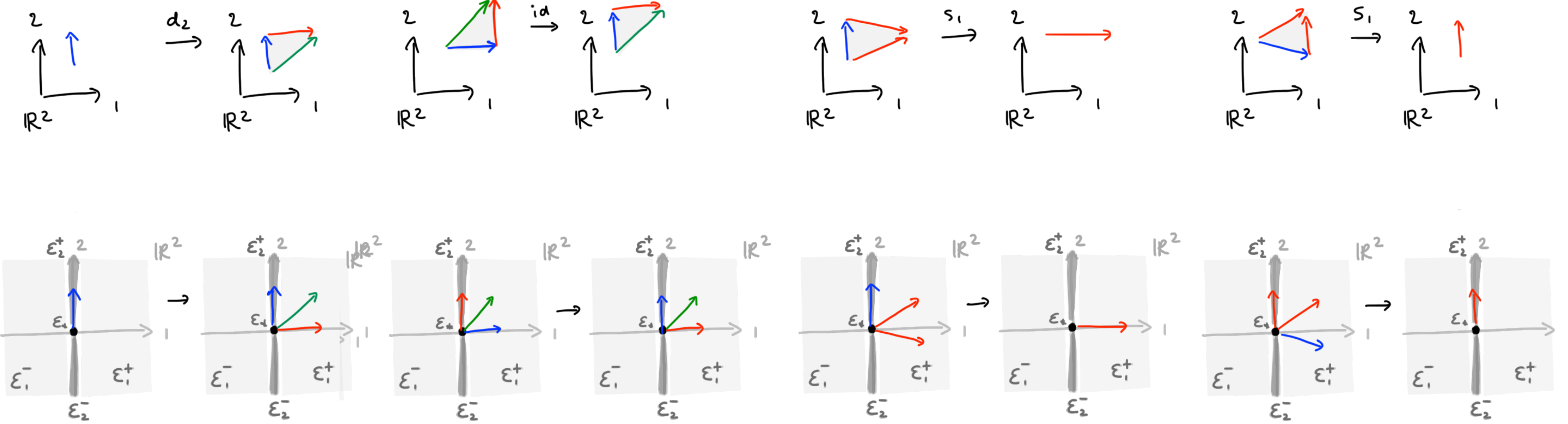

\caption{Framed and non-framed maps via framed realization.}
\label{fig:geometric-interpretation-of-framed-maps}
\end{figure}
\end{eg}

Note that the geometric description of framed maps of embedded framed simplices immediately generalizes to the case of embedded partially framed simplices. In combinatorial terms, this may be phrased as follows.

\begin{defn}[Framed maps of embedded partially framed simplices] \label{rmk:partial-framed-maps} Given $n$-embedded partially framed simplices $(S \epi [j],\cF)$ and $(T \epi [k],\cG)$, a \textbf{framed map} $F : (S \epi [j],\cF) \to (T \epi [k],\cG)$ is a simplicial map $F: S \to T$ that `descends' to a framed map of $n$-embedded simplices $F_n : ([j],\cF) \to ([k],\cG)$, that is,  $F_n : [k] \to [j]$ commutes with $F$ and the degeneracies $S \epi [k]$ and $T \epi [j]$.
\end{defn}

\nid Note that any framed map of $n$-embedded partially framed simplices $F : (S \epi [j],\cF) \to (T \epi [k],\cG)$ descends to a \emph{unique} framed map of $n$-embedded framed simplices $([j],\cF) \to ([k],\cG)$ in this way, and we always denote the latter map by $F_n$.

\begin{notn}[Category of embedded partially framed simplices] Denote by $\PaFrDelta n$ the category of $n$-embedded partially framed simplices and their framed maps.
\end{notn}

\pauseae

Framed maps either preserve the frame label of a vector or degenerate that vector to zero; there is a more general notion of `subframed map' in which vectors may degenerate not just to the zero vector but to any vector with more specialized frame label.  We first describe this geometric viewpoint, and then give a purely combinatorial definition of subframed maps.

\begin{rmk}[Subframed maps in linear algebraic terms] \label{rmk:subframed}
    Consider $n$-embedded framed simplices $(S \iso [l],\cF)$ and $(T \iso [m],\cG)$ and a simplicial map $F : S \to T$. For framed realizations $r : \abs{S} \into \lR^n$ and $q : \abs{T} \into \lR^n$, we say $F$ is `subframed' if for any vector $v$ in $S$ with $\vec r (\vec v) \in \eps^+_i$ we have $\vec q \circ \vec F (\vec v) \in \overline {\eps}^+_i$ (where $\overline {\eps}^+_i$ is the closure of the $i$th positive component $\eps^+_i$ of $\lR^n$; note that this closure in particular contains $0 \in \lR^n$, see \autoref{rmk:framed-map-geo}).
\end{rmk}

\nid Note that subframed maps may in particular send vectors from a positive component $\eps^+_i$ into a negative component $\eps^-_j \subset \overline {\eps}^+_i$ (where $j > i$).

\begin{defn}[Subframed map] \label{def:subframed}
    Given $n$-embedded framed simplices $(S \iso [l],\cF)$ and $(T \iso [m],\cG)$, a \textbf{subframed map} $F : (S \iso [l],\cF) \to (T \iso [m],\cG)$ is a simplicial map $F: S \to T$ such that for every vector $v: [1] \into S$ in the simplex $S$, either its frame label is preserved, i.e.\ $\rest \cF v = \rest \cG {F \circ v}$, specialized, i.e.\ $\rest \cF v < \rest \cG {F \circ v}$, or the vector is degenerated, i.e.\ $F \circ v: [1] \to [m]$ is constant.
\end{defn}

\nid The definition of subframed maps extends, as with that of framed maps, to the embedded partially framed case, by insisting that a vector without a frame label is either mapped to zero or to again a vector without a frame label.

\begin{eg}[Subframed maps]
In the first column of \autoref{fig:sub-framed-maps-and-their-geometric-interpretation} we illustrate two subframed maps of $2$-embedded framed 1-simplices $F : (S \iso [1],\cF) \to (T \iso [1],\cF)$, labeled `$\id$' and `$-\id$' respectively. In the second column, we pick framed realizations $r : \abs{S} \into \lR^2$ and $q : \abs{T} \into \lR^2$, and illustrated in green the image $r(v)$ and in red the image $q(w)$ of the unique nonzero frame vectors $v$ and $w$ in $S$ resp.\ $T$ (as based vectors in $\lR^2$); for both maps, we also depict the image $q \circ F(v)$ (in the same picture as $q(w)$). In the third column, we similarly depict the unbased action of $F$ as a linear map $\vec F : \vec S \to \vec T$ on the corresponding linear vectors $\vec v$ and $\vec w$ (shown embedded in $\lR^2$ via $\vec r : \vec V(S) \into \lR^2$ resp.\ $\vec q : \vec V (T) \into \lR^2$). In particular, note that for the second map the vector $\vec v \in \eps^+_1$ to a vector in $\eps^-_2$.
\begin{figure}[ht]
    \centering
    \def\svgwidth{1\columnwidth}
    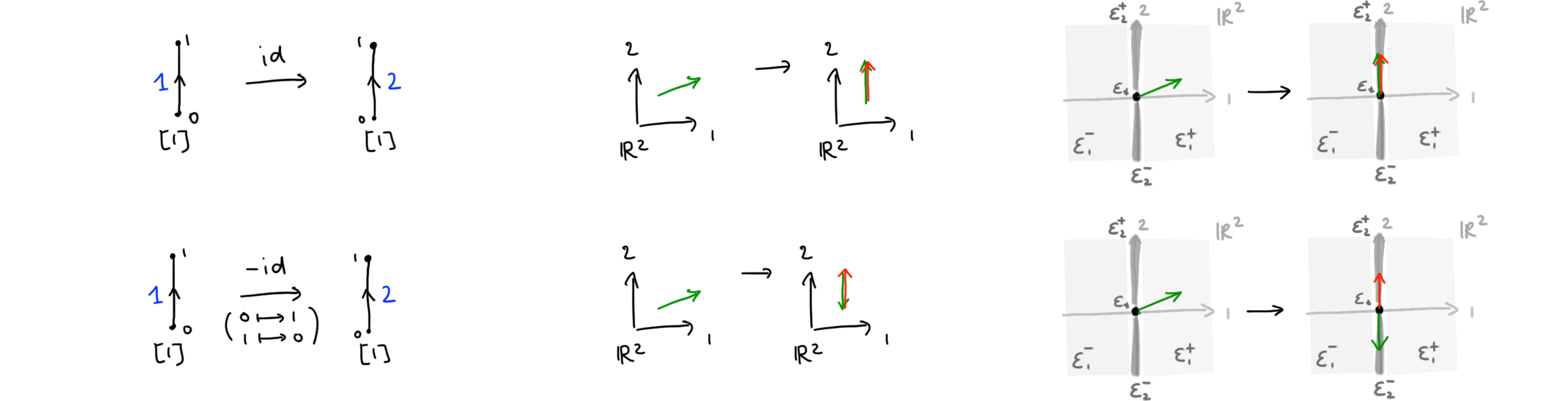

    \caption{Subframed maps and their geometric interpretation.}
    \label{fig:sub-framed-maps-and-their-geometric-interpretation}
\end{figure}
\end{eg}

\subsection{Proframes on simplices} \label{ssec:combinatorial-proframes}


We now recast the definition of frames of simplices in the form of so called `proframes'; these are towers of simplicial projections with 1-dimensional affine kernels. The usefulness of proframes stems from the `non-local' perspective they provide: while we think of frames as choices of `infinitesimal' vectors in tangent spaces, we can think of proframes as choices of `extensional' projections defined on the entire simplex (this perspective also underlies our later discussion of gradients and integrals, which formally relates proframes back to frames). Just as frames on simplices have a classical analog in terms of linear frames of vector spaces, the notion of proframes also finds an analog in classical linear algebraic terms yielding a notion of `linear proframes' (discussed in more detail in Appendix \ref{appsec:linear-fr}); the prototypical model of a linear proframe is the standard euclidean proframe of $\lR^n$ (from which all other proframes can then be obtained via trivializations $V \toiso \lR^n$).

\begin{term}[The standard euclidean proframe] \label{defn:standard-proframe} The `standard euclidean proframe' of $\lR^n$ is the sequence of projections
\[
    \lR^n \xepi {\pi_n} \lR^{n-1} \xepi {\pi_{n-1}} \lR^{n-2} \xepi {\pi_{n-2}} \cdots \xepi {\pi_2} \lR^{1} \xepi {\pi_1} \lR^0
\]
where $\pi_i : \lR^i \to \lR^{i-1}$ forgets the last component of vectors in $\lR^i$.
\end{term}

\subsubsecunnum{The definition of proframes}

As in the case of frames, we first introduce a combinatorial analog of (`standard' and `general') proframes, in which the role of a codimension-1 projection is played by codimension-1 degeneracies.

\begin{defn}[Proframe on a standard simplex]
    A \textbf{proframe $\cP$ of the standard $m$-simplex $[m]$} is a sequence $(p_m, p_{m-1}, \ldots, p_1)$ of surjective simplicial maps of the form
\[
[m] \xto {p_m} [m-1] \xto {p_{m-1}} [m-2] \xto{} \cdots \xto{} [1] \xto {p_1} [0]
.\qedhere
\]
\end{defn}

\begin{defn}[Proframe on a simplex] A \textbf{proframe} of an $m$-simplex $S$ is an isomorphism $S \iso [m]$ together with a proframe $\cP$ on $[m]$.
\end{defn}

\nid We usually denote proframes on $S$ by tuples $(S \iso [m],\cP)$, or keep the isomorphism implicit and simply say $\cP$ is a proframe on $S$.

\begin{eg}[Proframes on simplices]
In \autoref{fig:proframes-on-simplices} we illustrate four proframed $m$-simplices $(S \iso [m],\cP)$. Each degeneracy is indicated by highlighting its affine kernel.
\begin{figure}[ht]
    \centering
    \def\svgwidth{1\columnwidth}
    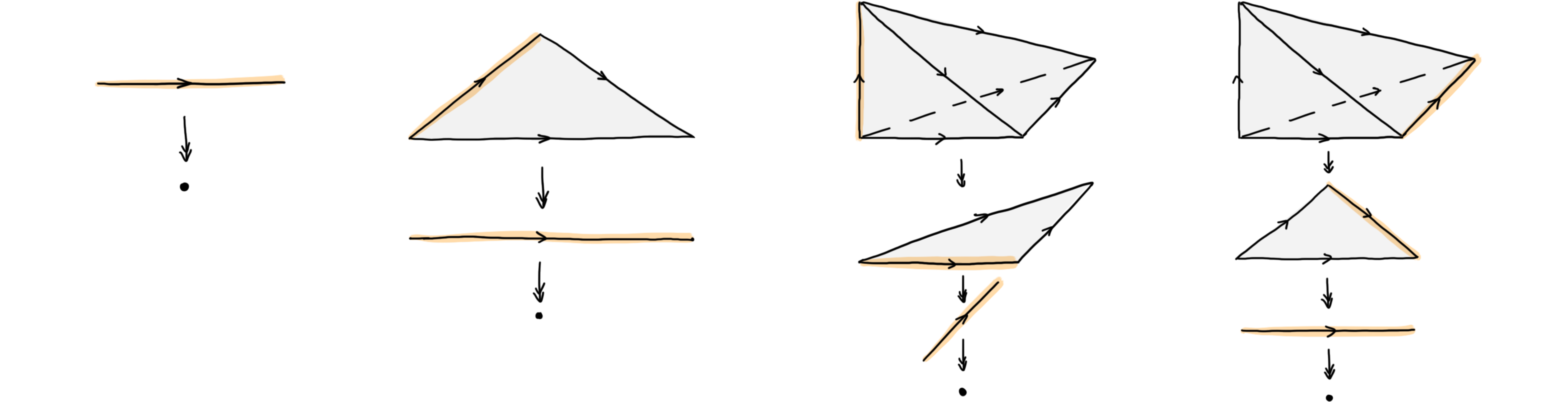

    \caption{Proframes on simplices.}
    \label{fig:proframes-on-simplices}
\end{figure}
\end{eg}

Recall that a linear map from a framed simplex $([m],\cF)$ to euclidean space is a framed realization when it takes each frame vector $v_i$ into the $i$th positive component $\eps^+_i$ of $\lR^m$.  There is an analogous notion of realization for proframes, as mappings into the standard proframed euclidean space, as follows.

\begin{defn}[Proframed realization of a proframed simplex]
    A \textbf{proframed realization} of a proframed simplex $(S \iso [m],\cP=(p_m,p_{m-1},\ldots,p_1))$ is a sequence of linear embeddings $r^\cP_i: \Delta^{i} \into \lR^i$, giving a commutative diagram,
\begin{equation}
\begin{tikzcd}
    {\abs{S} \iso \Delta^{m}} \arrow[r, "p_{m}", two heads] \arrow[d, "r^{\cP}_m", hook]
& {\Delta^{m-1}} \arrow[r, "p_{{m-1}}", two heads] \arrow[d, "r^{\cP}_{m-1}", hook]
& ... \arrow[r, "p_{2}", two heads]
& {\Delta^{1}} \arrow[r, "p_{1}", two heads] \arrow[d, "r^{\cP}_1", hook]
& {\Delta^{0}} \arrow[d, "r^{\cP}_0", hook]
\\
{\lR^m} \arrow[r, "\pi_{m}", two heads]
& {\lR^{m-1}} \arrow[r, "\pi_{m-1}", two heads]
& ... \arrow[r, "\pi_{2}", two heads]
& {\lR^{1}} \arrow[r, "\pi_{1}", two heads]
& {\lR^{0}}
\end{tikzcd}
\end{equation}
such that, for all $i$, $r^\cP_i$ embeds the unique kernel vectors $v = \keraff(p_i)$ orientation-preservingly into the $\lR$-fiber of $\pi_i$ over the point $r^\cP_{i-1} \circ p_i(v)$.
\end{defn}

\begin{eg}[Proframed realization of a proframed simplex]
In \autoref{fig:proframedrealization} we illustrate a proframed realization of a proframed 2-simplex.
\begin{figure}[h!]
    \centering
    \def\svgwidth{1\columnwidth}
    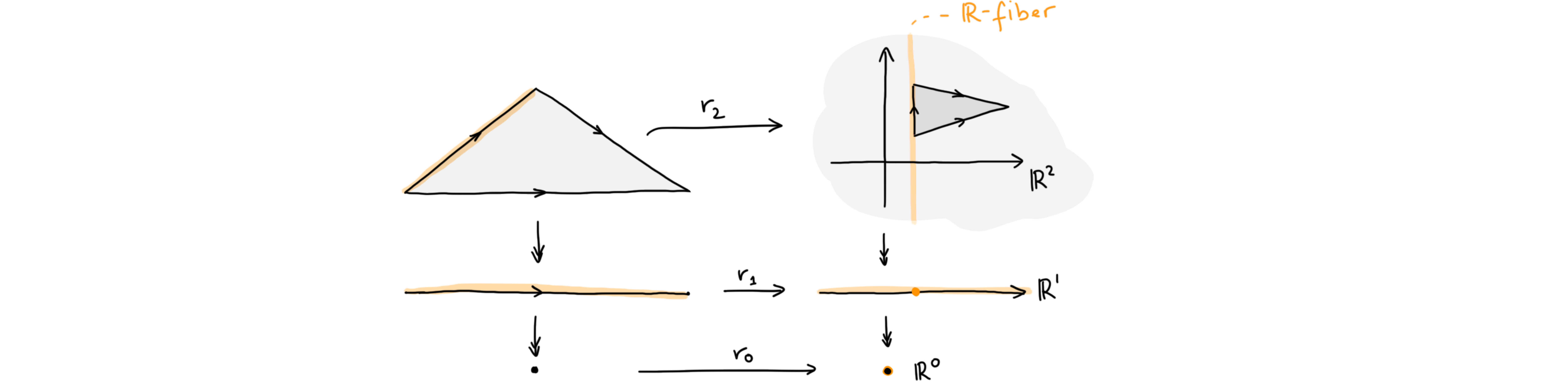

    \caption{Proframed realization of a proframed simplex.}
    \label{fig:proframedrealization}
\end{figure}
\end{eg}

\begin{rmk}[Relations of proframed and framed realizations] \label{rmk:proframed-realization-role} In \autoref{ssec:gradients-and-integration} we will show that proframes $\cP$ and frames $\cF$ on the simplex $[m]$ are in correspondence (i.e.\ they define the same structure on $[m]$); for corresponding proframes $\cP$ and frames $\cF$ on $[m]$, a proframed realization $\{r^\cP_i\}$ of $\cP$ determines is determined by a framed realization $r_\cF$ of $\cF$ by setting $r^\cP_m = r_\cF$.
\end{rmk}

\pauseae

A partial frame of a simplex was a frame of a quotient of the simplex by a simplicial subspace (the `unframed subspace'); similarly, a partial proframe of a simplex will be a proframe of a `quotient' of the simplex, as follows.

\begin{defn}[Partial proframe on a simplex]
A \textbf{$k$-partial proframe} on an $m$-simplex $S$ is a degeneracy $S \epi [k]$ together with a proframe $\cP = (p_k, p_{k-1}, \ldots, p_1)$ of $[k]$.
\end{defn}

\nid We usually denote $k$-partial proframe $\cP$ on a simplex $S$ by tuples $(S \epi [k],\cP)$. Analogous to the case of partial frames, we refer to the affine kernel  $U = \keraff(S \epi [k])$ as the `unframed subspace' of $(S \epi [k],\cP)$. Note that in an $m$-partial proframe of an $m$-simplex $(S \epi [m],\cP)$ the degeneracy $S \epi [m]$ must be an isomorphism, and thus $m$-partial proframes of $m$-simplices are proframes on $m$-simplices.

\begin{eg}[Partial proframes on simplices]
In \autoref{fig:partial-proframes-on-simplices} we depict a few partially proframed simplices.  As before, each degeneracy is given by highlighting its affine kernel.
\begin{figure}[ht]
    \centering
    \def\svgwidth{1\columnwidth}
    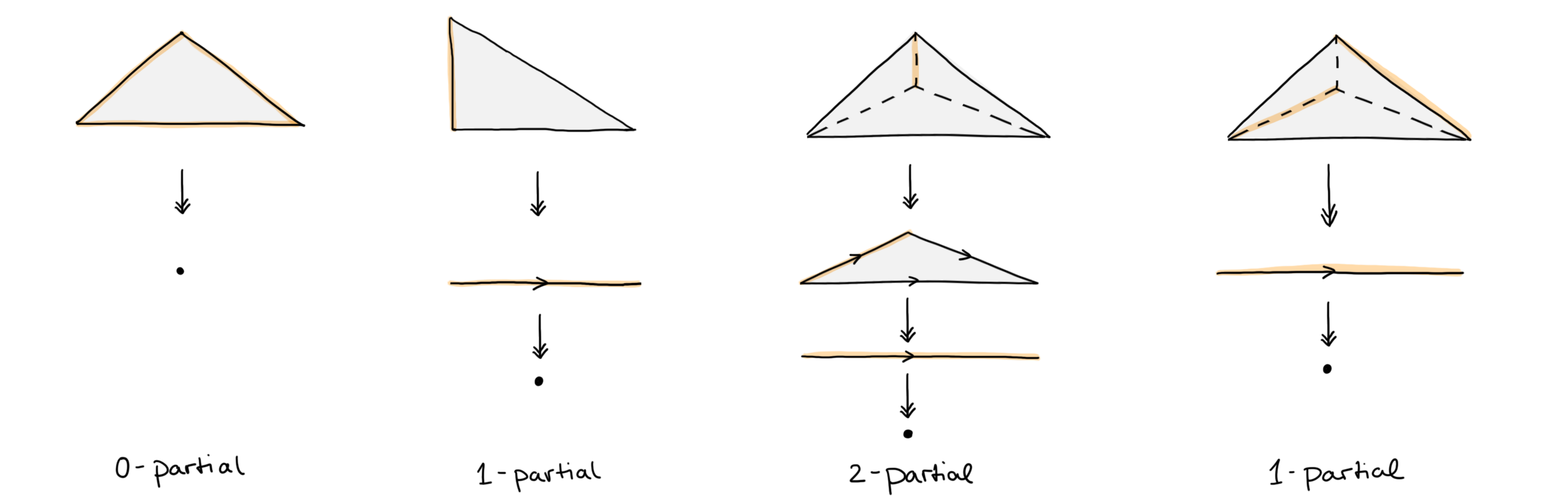

    \caption{Partial proframes on simplices.}
    \label{fig:partial-proframes-on-simplices}
\end{figure}
\end{eg}


The notion of proframed realization adapts to the partial case as follows: a  `proframed realization' of a $k$-partially proframed $m$-simplex $(S \epi [k],\cP)$ is simply a proframed realization $\{r^\cP_i\}$ of $([k],\cP)$ which now further determines a composite $r^{\cP} = r^\cP_k \circ \abs{S \epi [k]}$ mapping $\abs{S}$ to $\lR^k$. \autoref{rmk:proframed-realization-role} will also apply with evident adjustments: a proframed realization $\{r^\cP_i\}$ of a $k$-partially proframed simplex $(S \epi [k],\cP)$ determines and is determined by a framed realization $r_\cF$ of the corresponding $k$-partially framed simplex $(S \epi [m],\cF)$ by setting $r^\cP = r_\cF$.



\subsubsecunnum{Embedded proframes}


We introduce `$n$-embedded proframes' as a combinatorial analog of our earlier notion of `linear $n$-embedded proframes'. As in the linear case, while a proframed $m$-simplex is a sequence of codimension-1 degeneracies, an $n$-embedded proframed $m$-simplex will be a sequence of $n$ degeneracies each either of codimension-1 or codimension-0. Again, we introduce both `standard' and `general' notions of embedded proframes.

\begin{defn}[Embedded proframe on a standard simplex] \label{defn:embedded-proframe}
    An \textbf{$n$-embedded proframe $\cP$ of the standard $m$-simplex $[m]$} is a sequence $(p_n, p_{n-1},\ldots,p_1)$ of surjective simplicial maps
\[
[m] = [m_n] \xto{p_n} [m_{n-1}] \xto{p_{n-1}} [m_{n-2}] \xto{p_{n-2}} \cdots \xto{p_2} [m_1] \xto{p_1} [m_0] = [0]
\]
where for all $i$, either $m_{i-1} = m_i$ or $m_{i-1} = m_i - 1$.
\end{defn}

\begin{defn}[Embedded proframe on simplex] An \textbf{$n$-embedded proframe} of an $m$-simplex $S$ is an isomorphism $S \iso [m]$ together with an $n$-embedded proframe $\cP$ on $[m]$.
\end{defn}

\nid We usually denote $n$-embedded proframed simplices by tuples $(S \iso [m],\cP)$ (where $\cP$ is an $n$-embedded proframe of $[m]$). Note that an $m$-embedded proframe on the $m$-simplex is the same as a proframe on the $m$-simplex as previously defined.

\begin{eg}[Embedded proframes]
In \autoref{fig:embedded-proframes-on-simplices} we illustrate a few $n$-embedded proframed $m$-simplices $(S \iso [m],\cP)$. As before, each degeneracy is given by highlighting its affine kernel.
\begin{figure}[ht]
    \centering
    \def\svgwidth{1\columnwidth}
    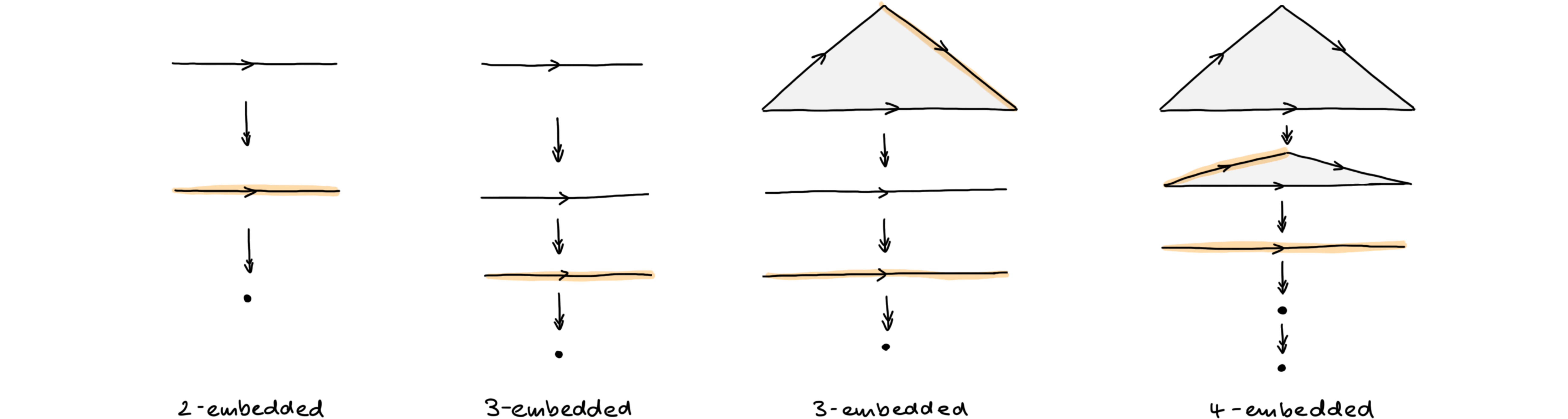

    \caption{Embedded proframes on simplices.}
    \label{fig:embedded-proframes-on-simplices}
\end{figure}
\end{eg}

As a proframe of an $m$-simplex has a notion of proframed realizations given by embeddings into the standard proframe of euclidean $m$-space, similarly an $n$-embedded proframe on an $m$-simplex has proframed realizations given by embeddings into the standard proframe of euclidean $n$-space, as follows.

\begin{defn}[Proframed realization of an embedded proframed simplex] \label{defn:proframed-real-emb}
A \textbf{proframed realization} of an $n$-embedded proframed $m$-simplex $(S \iso [m],\cP)$ with $\cP = (p_n,p_{n-1},\ldots,p_1)$ is a sequence of linear embeddings $r^\cP_i: \Delta^{m_i} \into \lR^i$, giving a commutative diagram,
\begin{equation}
\begin{tikzcd}
    {\abs{S} \iso \Delta^{m_n}} \arrow[r, "p_{n}", two heads] \arrow[d, "r^{\cP}_n", hook]
& {\Delta^{m_{n-1}}} \arrow[r, "p_{{n-1}}", two heads] \arrow[d, "r^{\cP}_{n-1}", hook]
& ... \arrow[r, "p_{2}", two heads]
& {\Delta^{m_1}} \arrow[r, "p_{1}", two heads] \arrow[d, "r^{\cP}_1", hook]
& {\Delta^{m_0}} \arrow[d, "r^{\cP}_0", hook]
\\
{\lR^n} \arrow[r, "\pi_{n}", two heads]
& {\lR^{n-1}} \arrow[r, "\pi_{n-1}", two heads]
& ... \arrow[r, "\pi_{2}", two heads]
& {\lR^{1}} \arrow[r, "\pi_{1}", two heads]
& {\lR^{0}}
\end{tikzcd}
\end{equation}
such that $r^\cP_i$ orientation-preservingly embeds kernel vectors $v = \keraff(p_i)$ (for all $i$ for which $v$ exists) into the $\lR$-fiber of $\pi_i$ over the point $r^\cP_{i-1} \circ p_i(v)$.
\end{defn}

\begin{eg}[Proframed realization of an embedded proframed simplex]
In \autoref{fig:emb-proframedrealization} we illustrate a proframed realization of a 3-embedded proframed 2-simplex.
\begin{figure}[h!]
    \centering
    \def\svgwidth{1\columnwidth}
    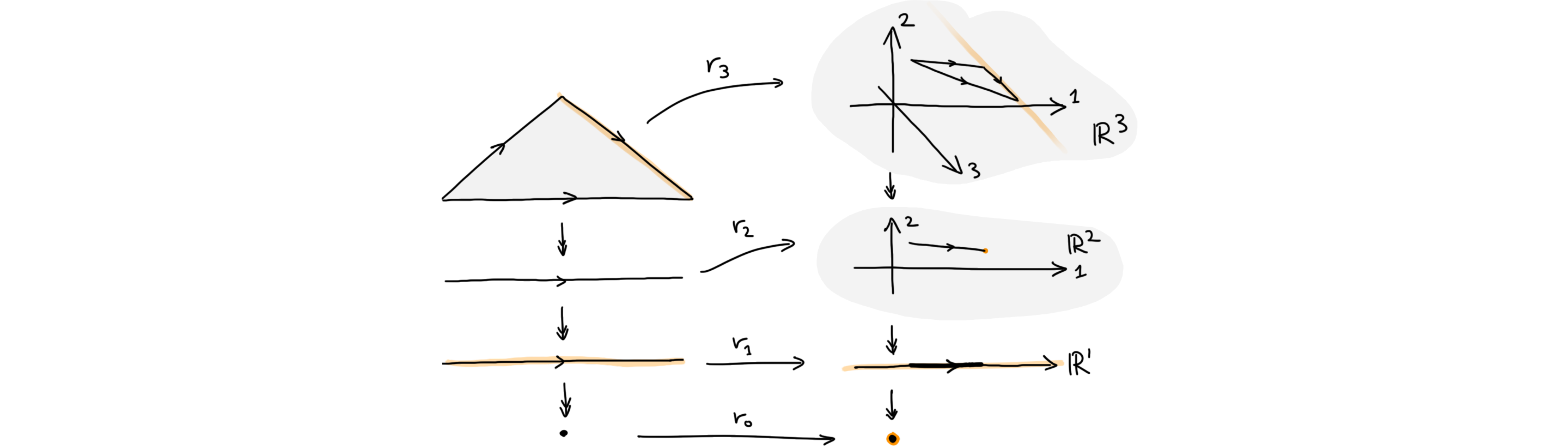

    \caption{Proframed realization of an embedded proframed simplex.}
    \label{fig:emb-proframedrealization}
\end{figure}
\end{eg}

\nid The relation of proframed and framed realizations as described in \autoref{rmk:proframed-realization-role} generalizes to the embedded case: for corresponding $n$-embedded proframes $\cP$ and $n$-embedded frames $\cF$ on $[m]$, a proframed realization $\{r^\cP_i\}$ of $\cP$ determines is determined by a framed realization $r_\cF$ of $\cF$ by setting $r^\cP_n = r_\cF$. We will revisit this observation in more formal terms in \autoref{rmk:corr-of-framed-and-proframed-real}.

\pause

Recall, a $k$-partial proframe of an $m$-simplex is a degeneracy to the $k$-simplex followed by a proframe of that simplex.  Similarly, an $n$-embedded $k$-partial proframe of an $m$-simplex is a degeneracy to the $k$-simplex followed by an $n$-embedded proframe of that simplex.

\begin{defn}[Embedded partial proframe on a simplex]
An \textbf{$n$-embedded $k$-partial proframe} on an $m$-simplex $S$ is a degeneracy $S \epi [k]$ together with an $n$-embedded proframe $\cP = (p_n, p_{n-1}, \ldots, p_1)$ of $[k]$.
\end{defn}

\nid We usually denote $n$-embedded $k$-partial proframe $\cP$ on a simplex $S$ by tuples $(S \epi [k],\cP)$, and refer to the affine kernel $U = \keraff(S \epi [k])$ as the `unframed subspace' of $(S \epi [k],\cP)$. Of course, an $n$-embedded $m$-proframe of an $m$-simplex is the same structure as an $n$-embedded proframe of an $m$-simplex as previously defined

\begin{eg}[Embedded partial proframes]
In \autoref{fig:embedded-partial-proframes-on-simplices} we illustrate three embedded partial proframes. As before, we depict degeneracies by highlighting their affine kernels.
\begin{figure}[ht]
    \centering
    \def\svgwidth{1\columnwidth}
    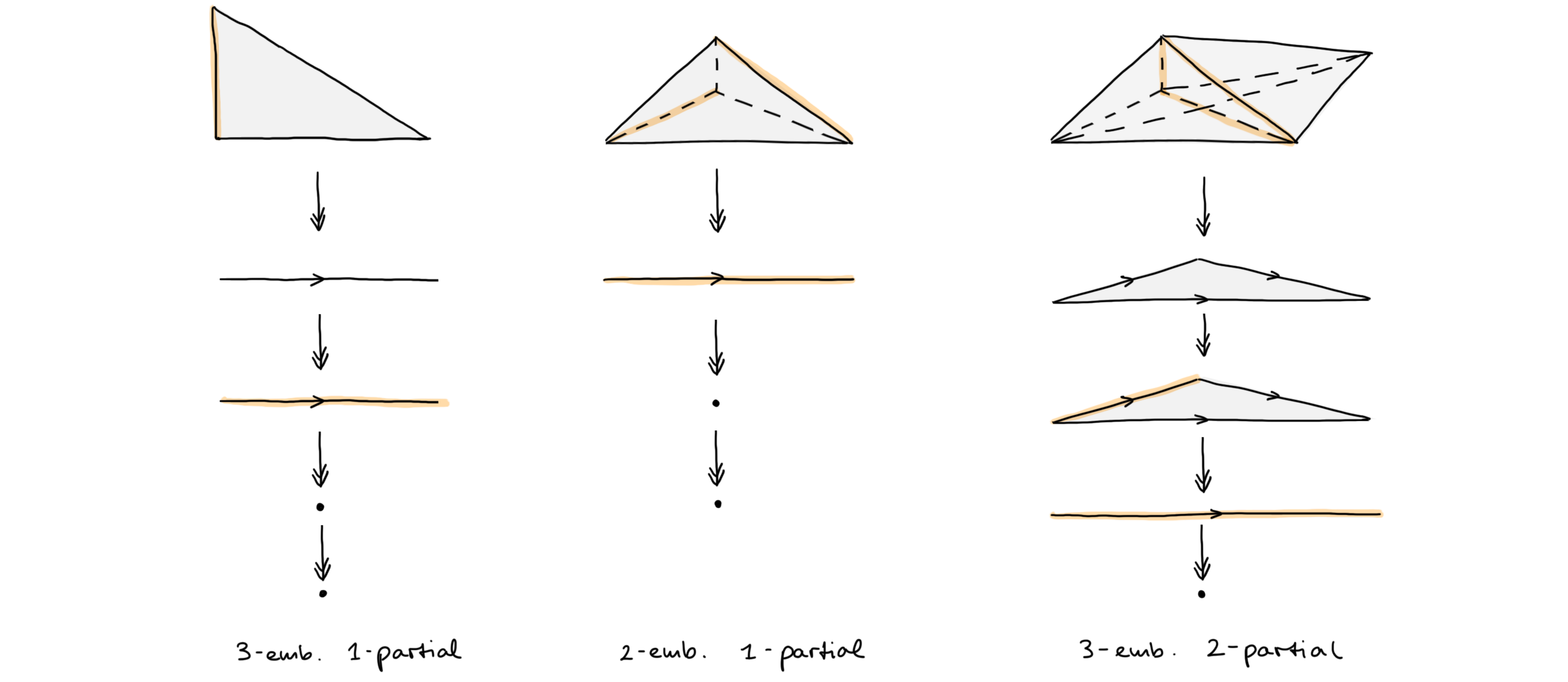

    \caption{Embedded partial proframes on simplices.}
    \label{fig:embedded-partial-proframes-on-simplices}
\end{figure}
\end{eg}


The notion of proframed realization of embedded proframes adapts to the partial case as follows: a proframed realization of an $n$-embedded $k$-partially proframed $m$-simplex $(S \epi [k],\cP)$ is simply a proframed realization $\{r^\cP_i\}$ of the $n$-embedded proframe $([k],\cP)$ which further determines a composite $r^{\cP} = r^\cP_k \circ \abs{S \epi [k]}$ mapping $\abs{S}$ to $\lR^n$. \autoref{rmk:proframed-realization-role} will apply with evident adjustments: a proframed realization $\{r^\cP_i\}$ of a $k$-partially proframed simplex $(S \epi [k],\cP)$ determines and is determined by a framed realization $r_\cF$ of the corresponding $k$-partially framed simplex $(S \epi [m],\cF)$ by setting $r^\cP = r_\cF$.


\pause

There is a notion `indframes' which is dual to that of proframes. The duality can also be phrased in classical linear algebraic terms (see \autoref{obs:linear-eqv-proframe-indframe}). Just as a proframe is given by a sequence of simplicial projections, an indframe is given by a sequence simplicial subspaces (i.e.\ of affine faces). Unlike proframes, however, indframes cannot be canonically expressed in terms of ordinary (i.e.\ non-affine) simplicial maps (which reflects the classical case of `projections' and `subspaces' of affine spaces, see \autoref{obs:affine-asymmetry}). This observation leads to proframes on simplices being an inherently more natural structure to work with combinatorially. Nonetheless, for completeness, we record notions of indframes on simplices as follows.

\begin{defn}[(Embedded) indframe on standard simplex] An \textbf{$n$-embedded indframe $\cI$ of the standard $m$-simplex $[m]$} is a sequence
\[
    [0] = [m_n] \intoaff [m_{n-1}] \intoaff  [m_{n-2}] \intoaff \cdots \intoaff [m_1] \intoaff [m_0] = [m]
\]
of simplicial subspaces $s_i : [m_i] \intoaff [m_{i-1}]$ where for all $i$, either $m_{i-1} = m_i$ or $m_{i-1} = m_i + 1$. If $m = n$ then we call $\cI$ simply an `indframe'.
\end{defn}

\begin{defn}[(Embedded) indframed on simplex] An \textbf{$n$-embedded indframe} $\cI$ of an $m$-simplex $S$ is an isomorphism $S \iso [m]$ together with an $n$-embedded indframe $\cI$ on $[m]$. We call $\cI$ simply an `indframe' if $m = n$.
\end{defn}

\nid Replacing the isomorphism `$S \iso [m]$' with a degeneracy `$S \epi [k]$' one similarly defines (embedded) \emph{partial} indframes.

\begin{obs}[Equivalence of indframes and proframes] As in the case of linear indframes and proframes, a (`partial', `$n$-embedded', or `$n$-embedded partial') indframe $(S \epi [k],\cI)$ defines the same structure on $S$ as a proframe $(S \epi [k],\cP)$: namely, with the same degeneracy $S \epi [k]$, the indframe $\cI = \{[l_{i+1}] \intoaff [l_{i}]\}$ of $[k]$ determines and is determined by the proframe $\cP = \{[j_{i+1}] \epi [j_{i}]\}$ of $[k]$ by defining $\imaff([l_i] \intoaff [k]) = \keraff([k] \epi [j_i]) \subset \spine[k]$.
\end{obs}

\nid For reasons explained above, we henceforth solely focus on the case of proframes.

\subsubsecunnum{Restricting proframes}


Recall that the frame labels of the spine vectors of a framed simplex propagate via a kinship relation to frame labels on all vectors, and thereby determine the restriction of a frame to any face of the simplex.  A proframe on a simplex similarly restricts to any face; as the proframe is encoded directly in the simplicial sequence, without any explicit frame labels, the restriction procedure is especially direct, as follows.

\begin{constr}[Proframe restriction in standard simplices] \label{constr:proframe-restriction}
    Let $\cP = ([m_n] \xto{p_n} [m_{n-1}] \xto{p_{n-1}} \cdots \xto{p_1} [m_0])$ be an $n$-embedded proframe of the standard $m$-simplex $[m]$, and let $f: [j] \to [m]$ be a $j$-face of $[m]$.  The `restriction' $\rest \cP f$ of the proframe to the $j$-face is the sequence given by the upper row in the diagram
\[
\begin{tikzcd}
    {[j] = [j_n]} \arrow[r, "f_n^*p_{n}", two heads, dashed] \arrow[d, "f = f_n"', hook] & {[j_{n-1}]} \arrow[r, "f_{n-1}^*p_{{n-1}}", two heads] \arrow[d, "f_{n-1}"', hook] & \cdots \arrow[r, "f_2^*p_{2}", two heads] & {[j_1]} \arrow[r, "f_1^*p_{1}", two heads] \arrow[d, "f_{1}", hook] & {[j_0] = [0]} \arrow[d, "f_{0}",hook] \\
    {[m]=[m_n]} \arrow[r, "p_{n}", two heads]                           & {[m_{n-1}]} \arrow[r, "p_{{n-1}}", two heads]                     & \cdots \arrow[r, "p_{2}", two heads]    & {[m_1]} \arrow[r, "p_{1}", two heads]                     & {[m_0] = [0]}
\end{tikzcd}
\]
where, inductively, $f_{i-1} \circ f_i^*p_i$ is defined as the image factorization of the map $p_i \circ f_i$.
\end{constr}

\nid In other words, the restricted proframe $\rest \cP f$ is the `simplicial restriction sequence' of the proframe $\cP$ by the map $f$, analogous to our earlier notion of `restriction sequences' of sequences of vector space projections (see \autoref{term:fac-seq}).

\begin{defn}[Proframe restrictions to simplicial faces] \label{defn:proframe-restriction-gen}  The \textbf{proframe restriction} an $n$-embedded proframed $m$-simplex $(S \iso [m], \cP)$ to a $j$-face $f : T \into S$ is the $n$-embedded proframed $l$-simplex $(T \iso [j],\rest \cP f)$ where $T \iso [j]$ is the restricted isomorphism $\rest {(S \iso [m]} f$, and $\rest \cP f$ is the restriction of $\cP$ to the induced standard face $f : [j] \into [m]$ (see \autoref{notn:restricting-triv-to-faces}).
\end{defn}

\begin{eg}[Proframe restriction to simplicial faces]
 In \autoref{fig:inherited-proframes-on-faces} we depict a 4-embedded proframed 3-simplex, along with three of its faces and the proframes they inherit by restriction of the given proframe sequence.
\begin{figure}[ht]
    \centering
    \def\svgwidth{1\columnwidth}
    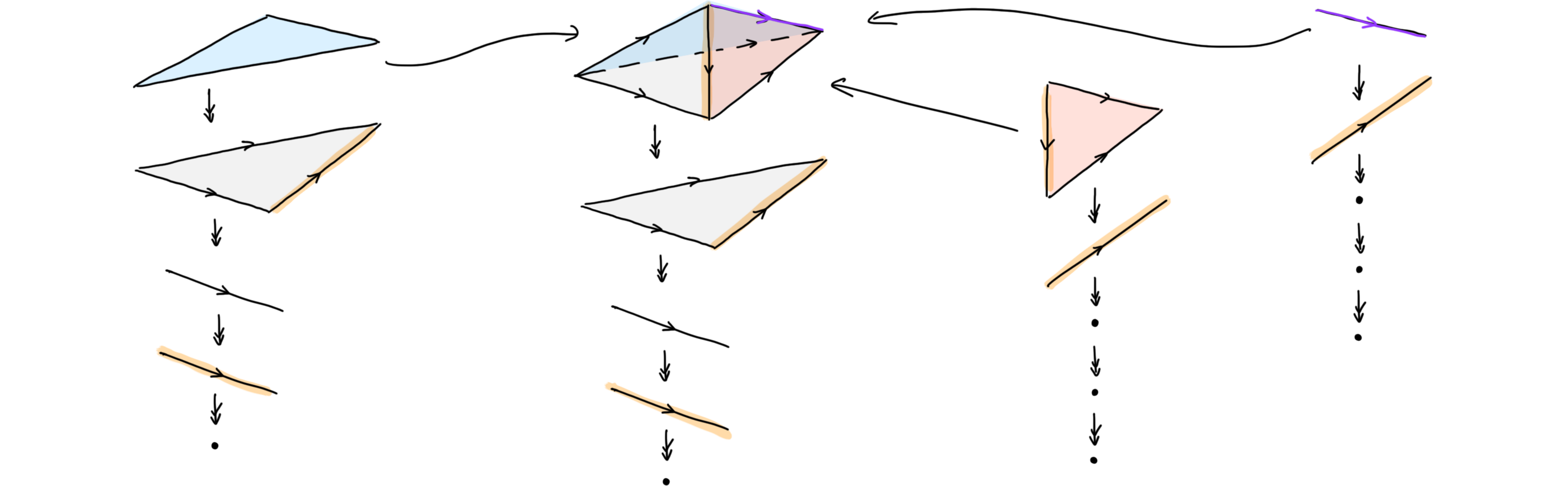

    \caption{Restricting an embedded proframe to the faces of a simplex.}
    \label{fig:inherited-proframes-on-faces}
\end{figure}
\end{eg}

\begin{rmk}[Proframe restriction of partial proframes]
    \autoref{defn:proframe-restriction-gen} applies with the evident adjustments to the case of embedded partially proframed simplices. Namely, the `proframe restriction' an $n$-embedded partially proframed $m$-simplex $(S \epi [k], \cP)$ to a $j$-face $f : T \into S$ is the $n$-embedded partially proframed $l$-simplex $(T \epi [l],\rest \cP f)$ where $T \epi [j]$ is the restricted degeneracy $\rest {(S \epi [m]} f$, and $\rest \cP f$ is the restriction of $\cP$ to the induced standard face $f : [l] \into [k]$.
\end{rmk}

\begin{rmk}[Proframe restriction via proframed realization]
    For an embedded proframed simplex $(S \iso [m],\cP)$, the proframe restriction to any simplicial $j$-face $f: T \into S$ is determined geometrically by a proframe realization $\{r^\cP_i : \Delta^{m_i} \into \lR^i\}$ of the simplex, as follows.  The restricted embedded proframe $\cQ \equiv \rest \cP f$ is the unique embedded proframe $\cQ$ on $T \iso [j]$ that has a proframed realization $\{r^{\cQ}_i : \Delta^{j_i} \into \lR^i\}$ whose top map factors as the face $f$ followed by the ambient proframed realization $r^\cP_n$, i.e.\ $r^{\cQ}_n = r^\cP_n \circ \abs{f}$. A similar remark applies to embedded partial proframes.
\end{rmk}

\subsubsecunnum{Proframed maps}

Recall that a framed map is a map of simplices that, for each vector of the source, either preserves the frame label of the vector or degenerates the vector.  Proframed maps similarly are maps, now of simplicial sequences, that for each vector of the source, either preserve the proframe of the vector or degenerate that vector, as follows.

\begin{defn}[Proframed maps]  \label{defn:proframed-maps}
Given $n$-embedded proframed simplices $(S \iso [l], \cP = (p_n, \ldots, p_1))$ and $(T \iso [m],\cQ = (q_n, \ldots, q_1))$, a \textbf{proframed map} $F: (S \iso [l],\cP) \to (T \iso [m],\cQ)$ is a map of sequences
\begin{equation}
    \begin{tikzcd}[baseline=(W.base)]
{[l]=[l_n]} \arrow[d, "F_n"'] \arrow[r, "p_n"] & {[l_{n-1}]} \arrow[d, "F_{n-1}"'] \arrow[r, "p_{n-1}"] & \cdots \arrow[d, "\cdots", phantom] \arrow[r, "p_2"] & {[l_1]} \arrow[d, "F_1"] \arrow[r, "p_1"] & {[l_0]=[0]} \arrow[d, "F_0"] \\
{[m] = [m_n]} \arrow[r, "{q_n}"'] & {[m_{n-1}]} \arrow[r, "{q_{n-1}}"'] & \cdots \arrow[r, "{q_2}"'] & {[m_1]} \arrow[r, "{q_1}"'] & |[alias=W]| {[m_0] = [0]}
    \end{tikzcd}
\end{equation}
such that for every vector $v: [1] \to [l]$, either its proframe is preserved, i.e.\ $\rest \cP {v} = \rest \cQ {F_n \circ v}$, or the vector is degenerated, i.e.\ $F_n \circ v : [1] \to [m]$ is constant.
\end{defn}

\begin{term}[Proframed face and proframed degeneracy maps]
An injective proframed map is called a `proframed face map', and a surjective proframed map is called a `proframed degeneracy map'.
\end{term}

\begin{obs}[Inclusions of proframe restrictions are proframed faces]
Given a proframed simplex $([m],\cQ)$ and a face $f: [l] \into [m]$ of the simplex, the inclusion of the restriction $([l], \rest \cQ f) \into ([m],\cQ)$ is a proframed face, and every proframed face is of this form.
\end{obs}

\begin{obs}[(Epi,mono)-factorization of proframed maps]
Every proframed map factors as a proframed degeneracy followed by a proframed face.
\end{obs}

\begin{eg}[Proframed maps] \label{eg:proframed-maps-of-simplices}
    In \autoref{fig:proframed-and-lax-proframed-maps} we illustrate proframed and non-proframed maps of proframed simplices. Degeneracies in proframes are indicated as usual by highlighting their affine kernels. Top-level map (which, in a maps of sequences of surjective maps, determine all other maps) coincide with the maps in \autoref{fig:geometric-interpretation-of-framed-maps}. The first map of sequences is proframed and injective, thus is a proframed face. The second map is not a map of sequences (the upper square does not commute as indicated). The third map of sequences is proframed and surjective, thus is a proframed degeneracy.  The last map, despite being a map of sequences, that is forming a commutative diagram as depicted, is not a proframed map: the map does not preserve the proframe of the diagonal vector of the 2-simplex of the source, i.e.\ the map restricted to that vector is not an isomorphism onto its image.v
  \begin{figure}[ht]
    \centering
    \def\svgwidth{1\columnwidth}
    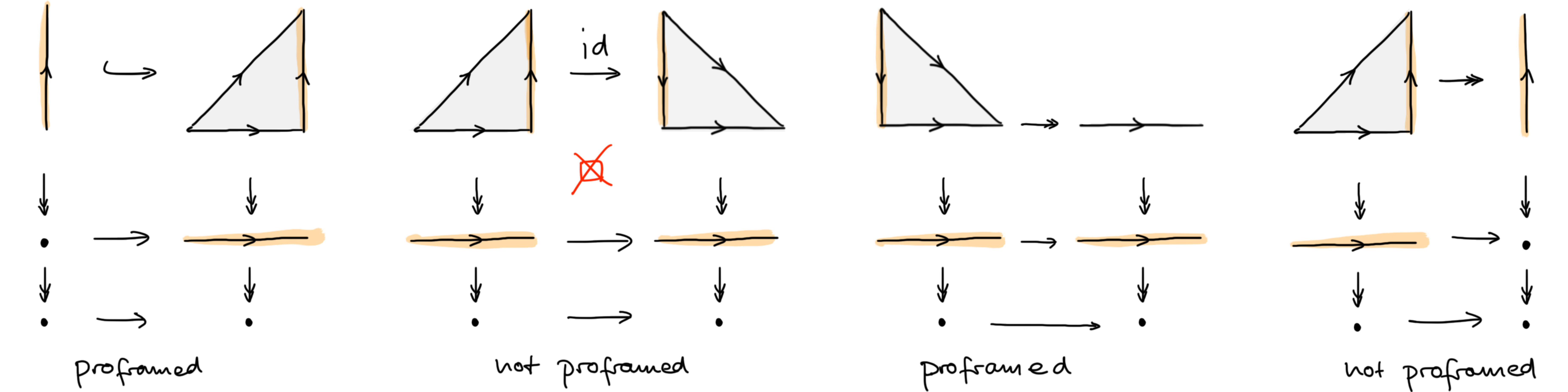

    \caption{Proframed and non-proframed maps.}
    \label{fig:proframed-and-lax-proframed-maps}
\end{figure}
\end{eg}

\begin{notn}[Category of proframed simplices]
Denote by $\PFrDelta n$ the category of $n$-embedded proframed simplices and their proframed maps.
\end{notn}

The definition of proframed maps generalizes naturally from the embedded to the embedded partial case.

\begin{defn}[Proframed maps of partially proframed simplices] \label{rmk:partial-proframed-maps}
    Given an $n$-embedded partially proframed simplices $(S \epi [k],\cP)$ and $(T \epi [j],\cQ)$ a \textbf{proframed map} $F : (S \epi [k],\cP) \to (T \epi [j],\cQ)$ is a simplicial map $F : S \to T$ that descends to a proframed map $(F_n,F_{n-1},..., F_0) : ([k],\cP) \to ([j],\cQ)$ whose top component $F_n : [k] \to [j]$ commutes with $F$ and the degeneracies $S \epi [k]$ and $T \epi [j]$.
\end{defn}

\nid Note that any proframed map $F : (S \epi [k],\cP) \to (T \epi [j],\cQ)$ must descend to a unique proframed map $([k],\cP) \to ([j],\cQ)$ whose components will always be denoted by $(F_n,F_{n-1},..., F_0)$.

\begin{notn}[Category of partially proframed simplices] Denote by $\PaPFrDelta n$ the category of $n$-embedded partially proframed simplices and their proframed maps.
\end{notn}

\nid Of course, we may regard proframed maps of proframed simplices as proframed maps of partially proframed simplices; thus $\PFrDelta n$ is a full and faithful subcategory of $\PaPFrDelta n$.


\begin{rmk}[Subproframed maps]
    Recall from \autoref{rmk:subframed} and \autoref{def:subframed} the notion of subframed map of framed simplices, in which framed vectors are allowed to specialize to have frame labels in lower strata of the standard stratification.  The corresponding notion of `subproframed map' of (partially) proframed simplices is remarkably simple: a subproframed map $F: (S \epi [k],\cP = (p_n,...,p_1)) \to (T \epi [j],\cQ = (q_n,...,q_1))$ is a collection of maps $(F,F_n,...,F_0)$ of unordered simplices $F : S \to T$ and $F_i : \Unord {[k_i]} \to \Unord {[j_i]}$ forming a commutative diagram (i.e.\ $F_n \circ (S \epi [k]) = (T \epi [j]) \circ F$ and $F_{i-1} \circ \Unord p_i = \Unord q_i \circ F_{i+1}$). This does not require any further conditions on conditions on vectors: instead, the structure of the sequence itself controls the specialization of frame vectors without explicit mention of an order on frame labels or the standard stratification of frame vectors.
\end{rmk}


\subsection{Gradient frames and integral proframes} \label{ssec:gradients-and-integration}

We conceive of frames as an infinitesimal notion, concerning `tangential' vectors, and of proframes as extensional notion, concerning quotients of space. (This contrast becomes visible, for instance, in the classical case of frames and proframes on affine spaces, see \autoref{obs:affine-frames} resp.\ \autoref{obs:affine-proframes}.) We will thus refer to the passage from frames to proframes as `integration', and to a converse operation as taking a `gradient'.



\subsubsecunnum{Gradients and integrals for simplices}

We introduce the gradient frame of a proframe, and conversely the integral proframe of a frame.

\begin{notn}[Composite degeneracies in proframes]
    For an $n$-embedded proframe $\cP=(p_n,\ldots,p_1)$ of the $m$-simplex $[m]$, we abbreviate the composite $p_i p_{i+1} \cdots p_n : [m] \to [k_{i-1}]$ by $p_{\geq i}$ and similarly the composite $p_{i+1} \cdots p_n : [m] \to [k_{i}]$ by $p_{>i}$.
\end{notn}

\begin{term}[Gradient frame of a proframed standard simplex] \label{defn:gradient-of-proframed-simplices}
Given an $n$-embedded proframed $m$-simplex $([m],\cP=(p_n,\ldots,p_1))$, its `gradient' $\Gradfr \cP$ is the $n$-embedded framed $m$-simplex $([m],\Gradfr \cP: \spine[m] \xslashedhookrightarrow{} \bnum n)$ with frame label $i$ on the spine vector $\keraff(p_{\geq i}) \backslash \keraff(p_{> i})$ (when that complement is nonempty).
\end{term}

\nid In other words, if the spine vector $v \in \spine[m]$ projects to a spine vector $p_{>i} v \in \spine[k_i]$ and the degeneracy $p_i: [k_i] \to [k_{i-1}]$ collapses that vector $p_{>i} v$, then the spine vector $v \in \spine[m]$ is given the frame label $i$.

\begin{term}[Integral proframe of a framed standard simplex]
Given an $n$-embedded framed $m$-simplex $([m],\cF: \spine[m] \into \bnum n)$, an `integral' $\Intfr \cF$ is an $n$-embedded proframed $m$-simplex $([m],\Intfr \cF)$ whose gradient is $([m],\cF)$.
\end{term}

\begin{constr}[Integral proframes of framed standard simplices exist]
Given an $n$-embedded framed $m$-simplex $([m],\cF: \spine[m] \into \bnum n)$, inductively set $p_i: [k_i] \to [k_{i-1}]$ to be the unique simplicial map collapsing the spine vector $p_{>i}(\cF^{-1}(i))$ i.e.\ the spine vector that has frame label $i$; if there is no spine vector with frame label $i$, then $p_i$ is taken to be the identity.  The sequence of degeneracies $\cP = (p_n,p_{n-1},\ldots,p_1)$ is an $n$-embedded $k$-partial proframe of the $m$-simplex, and has the frame $\cF$ as its gradient.
\end{constr}

\begin{obs}[Gradient and integral are inverse]
By definition, or equivalently by the construction of the integral in the previous observation, the integral is right-inverse to the gradient; similarly, taking the integral of the gradient of a proframe evidently reconstructs that proframe, and so the integral is also left-inverse to the gradient.  That is, for any $n$-embedded proframe $\cP$ and any $n$-embedded frame $\cF$, we have
\[
    \Gradfr \Intfr \cF = \cF \quad \text{and} \quad \Intfr \Gradfr \cP = \cP. \qedhere
\]
\end{obs}

The gradients and integrals now carry over to the case of general simplices (and the case of partial frames and proframes) as follows.

\begin{defn}[Gradient of a proframed simplex] Given an $n$-embedded partially proframed simplex $(S \epi [k],\cP)$ its \textbf{gradient} is the $n$-embedded partially framed simplex $(S \epi [k],\Gradfr \cP)$.
\end{defn}

\begin{defn}[Integral of a framed simplex] Given an $n$-embedded partially framed simplex $(S \epi [k],\cF)$ its \textbf{integral} is the $n$-embedded partially proframed simplex $(S \epi [k],\Intfr \cP)$.
\end{defn}


\begin{eg}[Gradient frame and integral proframe for a simplex]
In \autoref{fig:gradientintegralsimplex} we illustrate a 4-embedded 3-proframed 4-simplex and its corresponding gradient 3-embedded 2-framed 4-simplex; equivalently, that framed simplex integrates to that proframed simplex.  As further examples, note that the four embedded framed simplices in \autoref{fig:embedded-combinatorial-frames-and-their-embedded-coordinate-notation} are, in order, the gradients of the four embedded proframed simplices in \autoref{fig:embedded-proframes-on-simplices}.  Similarly, the three embedded partially framed simplices in \autoref{fig:embedded-partial-combinatorial-frames} are, in order, the gradients of the three embedded proframed simplices in \autoref{fig:embedded-partial-proframes-on-simplices}.
\begin{figure}[ht]
    \centering
    \def\svgwidth{1\columnwidth}
    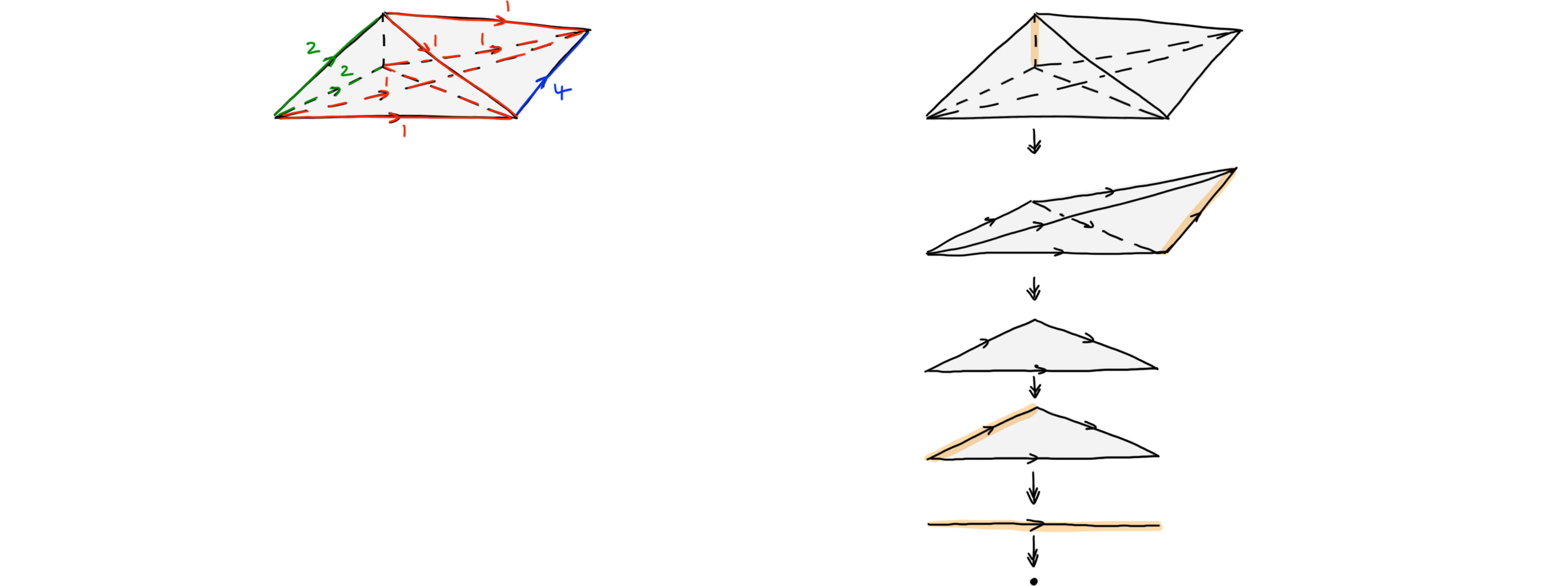

    \caption{The gradient frame and the integral proframe of a simplex.}
    \label{fig:gradientintegralsimplex}
\end{figure}
\end{eg}

\subsubsecunnum{Gradients and integrals for maps}

The gradient--integral relation between frames and proframes applies also to maps.

\begin{obs}[Gradients and integrals respect restriction] \label{rmk:framed-faces-are-compatible-with-gradients-and-integration}
    Given an $n$-embedded partially proframed simplex $(S \epi [k],\cP)$  and an $n$-embedded partially framed simplex $(S \epi [k],\cF)$, and a face $f: [j] \to [m]$, the gradient of the restriction to the face is the restriction of the gradient, and similarly for the integrals:
\[
    \Gradfr (\rest \cP f) = \rest {(\Gradfr \cP)} f\quad \text{and} \quad  \Intfr (\rest \cF f) = \rest {(\Intfr \cF)} f. \qedhere
\]
\end{obs}

\begin{eg}[Gradients and integrals respect restriction]
The framed simplex in \autoref{fig:inherited-frames-on-general-faces} is the gradient of the proframed simplex in \autoref{fig:inherited-proframes-on-faces}, and the gradients of the restricted proframes in the latter figure are the restricted frames in the former figure.  Conversely, the integrals of the simplex and its restrictions in the former figure are the simplex and restrictions in the latter figure.
\end{eg}

\begin{term}[Gradient framed map of a proframed map]
    Given an $n$-embedded partially proframed simplex $(S \epi [j],\cP)$, an $n$-embedded partially proframed simplex $(T \epi [k],\cQ)$, and a proframed map $F : (S \epi [j],\cP) \to (T \epi [k],\cQ)$, the `gradient' $\Gradfr F$ is simply the framed map $(S \epi [j],\Gradfr \cP) \to (T \epi [k],\Gradfr \cQ)$ determined by the simplicial map $F : S \to T$.
\end{term}

\begin{term}[Integral proframed map of a framed map]
Given an $n$-embedded partially framed simplex $(S \epi [j],\cF)$, an $n$-embedded partially framed simplex $(T \epi [k],\cG)$, and a framed map $F : (S \epi [j],\cF) \to (T \epi [k],\cG)$, an `integral' $\Intfr F$ is a proframed map $(S \epi [j],\Intfr \cF) \to (T \epi [k],\Intfr \cG)$ whose gradient is the framed map $F$.
\end{term}

\begin{eg}[Gradient framed maps and integral proframed maps]
The two framed maps in \autoref{eg:framed-and-lax-framed-maps-of-simplices} are the gradients of the two proframed maps in \autoref{eg:proframed-maps-of-simplices}, or equivalently the latter maps are the integrals of the former maps.
\end{eg}

As there exists a unique integral proframe of any framed simplex, similarly there exists a unique integral proframed map of any framed map.  Altogether, we have an equivalence of categories, as follows.

\begin{prop}[Correspondence of partial frames and partial proframes] \label{prop:integrals-and-gradients-inverse}
Gradient and integration are inverse isomorphic functors between the category of $n$-embedded partially proframed simplices with proframed maps and the category of $n$-embedded partially framed simplices with framed maps:
\[
    \pushQED{\qed}
        \Gradfr : \PaPFrDelta n \iso \PaFrDelta n : \Intfr .\qedhere
    \popQED
\]
\end{prop}

\nid In practice we will be focused on a special case of this correspondence, namely between $n$-embedded frames and $n$-embedded proframes.

\begin{cor}[Correspondence of frames and proframes] \label{cor:iso-of-frame-vs-proj-frame-simplex}
Gradient and integration are inverse isomorphic functors between the categories of $n$-embedded framed respectively proframed simplices and framed respectively proframed maps:
\[
    \pushQED{\qed}
        \Gradfr : \PFrDelta n \iso \FrDelta n : \Intfr .\qedhere
    \popQED
\]
\end{cor}

\begin{rmk}[Correspondence of framed and proframed realizations] \label{rmk:corr-of-framed-and-proframed-real} Given an $n$-embedded proframed simplex $(S \iso [m],\cP)$ with gradient framed simplex $(S \iso [m],\cF = \Gradfr \cP)$, then any proframed realization $\{r^\cP_i : \Delta^{m_i} \into \lR^i\}$ determines is and determined by a framed realization $r_\cF : \abs{S} \iso \Delta^{m_n} \into \lR^n$ by setting $r^\cP_n = r_\cF$.
\end{rmk}

\section{Framed simplicial complexes} \label{sec:framed-simplicial-complexes}

Our goal in this section will be to introduce a combinatorial notion of `framings' on simplicial complexes. This is a globalization of the notion of frames on simplices introduced in the previous section: just as manifolds are spaces that are locally modeled on euclidean space, simplicial complexes are modeled on simplices; and just as framed manifolds are locally modeled on framed euclidean space (i.e.\ endowed with a continuous choice of frames in each tangent spaces), framed simplicial complexes will be modeled on framed simplices. The notion of framed simplicial complexes has two important features which distinguish the combinatorial approach to framings from the classical geometric approach to framings of manifolds.

Firstly, framings of simplicial complexes are not `local linear' structures but `piecewise affine' structures in the following sense. As previously discussed, simplices are not infinitesimal geometric objects but are extended affine spaces; frames of simplices are, correspondingly, `affine frames', i.e.\ frames not based at any specific point of the simplex but defined `up to translation'. We will define framings on simplicial complexes by piecing together affine frames of each of their simplices, and in this sense framings will be `piecewise affine'. This stands in contrast to the classical `tangential' notion of framings which defines framings locally, i.e.\ by picking a frame in the tangent space of each point.

Secondly, simplicial complexes are naturally `singular spaces' and generally not manifolds. As a consequence, framed simplicial complex will in fact provide a combinatorial model of `framed singular spaces' and not just of classical `framed manifolds'. Classically, `singular spaces' are gluings of `manifold strata', and singular spaces themselves need not be manifolds. The question of framing singular spaces is subtle since the usual machinery of tangent spaces relies on local euclidean trivializations; these need not exist everywhere in singular spaces. We will not attempt to geometrically define `framed singular spaces', but instead focus on leveraging the tools of affine combinatorics: namely, in combinatorial terms, individual open simplices will play to role of (`pieced') manifold strata, and the question of how framings can transition between strata of different dimension will find an answer using the notion of `$n$-embedded' frames developed in the preceding section.

\begin{figure}[ht]
    \centering
    \def\svgwidth{1\columnwidth}
    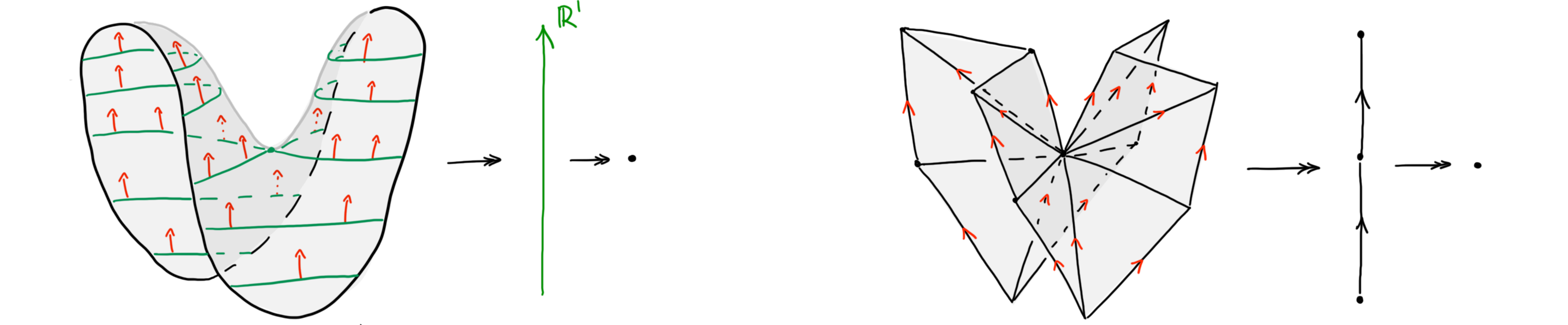

    \caption[Geometric proframings with singularities.]{A `1-framed singular space' modelled by a partial 1-framed simplicial complex.}
    \label{fig:the-saddle-singularity}
\end{figure}

Note `framed singular spaces' arise naturally in familiar situations: on the left in \autoref{fig:the-saddle-singularity}, we depict a `Morse saddle point' given by a function from the 2-disk to $\lR^1$. Up to choosing a Riemannian structure, this defines a gradient vector field (indicated by red arrows) everywhere \emph{except} at the central critical point; if we regard this point as a $0$-dimensional stratum with $0$-dimensional frame (and its complement as a $2$-dimensional stratum with $1$-dimensional  frame given by the gradient vector field), then this is an example of a `framed singular space'.
Combinatorially, this can be modelled by a sequence of surjective maps of simplicial complexes as shown on the right. Note that restricting this sequence to any simplex, yields a 1-embedded partially proframed 2-simplex as previously defined, and thus, by passing to gradient frames, a 1-embedded partially framed simplex (indicated by directing edges in red)\footnote{Following previous notation, each directed edge should also carry the frame label `1', see \autoref{eg:partial-frames}.}. The resulting simplicial complex, now with each simplex endowed with the structure of a 1-embedded partial frame, yields our first example of a (partial) 1-framed simplicial complex, and it models the `1-framed singular space' to its left.


We outline this section. We introduce a combinatorial notion of framings in \autoref{ssec:combinatorial-framings} following by a notion of proframings in \autoref{ssec:combinatorial-proframings}; both notions will be straight-forward generalizations of our earlier notions of framed resp.\ proframed simplices. As in the case of simplices, the relation of the framed and proframed simplicial complexes can be described in terms of via `gradient framings' and `integral proframings', which we will discuss in \autoref{ssec:gradients-integrability-pro-framings}; importantly, the global setting of framed simplicial complexes, integral proframings generally need not exists or not be unique if they exist. In \autoref{ssec:combinatorial-flatness} we will then introduce a combinatorial notion of (locally) flat framings, which provide an important class of `uniquely integrable' framings; local flatness yields a condition for combinatorial framings to be `locally trivializable' and thus `free of singularities', which is therefore analogous to classical tangential framings of manifolds as discussed above.

\subsection{Framings and proframings on simplicial complexes} \label{ssec:simplicial-complexes-recollection} \label{ssec:combinatorial-framings}

\subsubsecunnum{Simplicial complexes}

Before giving the definition of framings on simplicial complexes, we specify some basic terminology and notation.

\begin{term}[Simplicial sets] A `simplicial set $X$' is a presheaf $X : \Delta\op \to \SetCat$ on the simplex category $\Delta$. The `category of simplicial sets' is the category of presheaves on $\Delta$, and is usually denoted by $\SSet$.
\end{term}

\nid We will tacitly Yoneda embed $\Delta \to \SSet$ and, abusing notation, use the simplex $[k] \in \Delta$ to also denote the representable simplicial set that it defines under this embedding (this is sometimes denoted by $\Delta[k]$ in the literature). The set of maps $x : [k] \to X$ provides the set of $k$-simplices in a simplicial set $X$.

\begin{term}[Non-degenerate simplices] A simplex $x : [k] \to X$ in a simplicial set $X$ is `non-degenerate' if there is no non-identity degeneracy map $d : [k] \to [j]$ through which $x$ factors; otherwise, $x$ is `degenerate'.
\end{term}

\nid In piecewise linear topology one commonly uses the following weaker notion in place of simplicial sets.

\begin{term}[Simplicial complexes] A `simplicial complex' consists of a set of vertices $K[0]$ together with a list of subsets $K[1], K[2] ...$ of the vertex powerset $\powerset K[0]$ with the property that elements in $K[i]$ are sets of cardinality $(i+1)$ and any subset of an element in $K[i]$ is an element in some $K[j]$, $j \leq i$. Elements of $K[i]$ are also called $i$-simplices of $K$. Maps of simplicial complexes $K \to L$, called `simplicial maps', are maps of vertex sets $f : K[0] \to L[0]$ whose image on each $i$-simplex $x \in K[i]$ yields a $j$-simplex $f(x) \in L[j]$ (for some $j$). The category of simplicial complexes and their simplicial maps will be denoted by $\SCplx$.
\end{term}

\nid Just as simplicial sets, simplicial complexes may be characterized as (certain) presheaves as follows.

\begin{rmk}[Simplicial complexes as presheaves] Recall the category $\UnSimp$ of unordered simplices (see \autoref{term:standard-unordered-simplices}). Every simplicial complex $K$ gives rise to a presheaf $K : \UnSimp\op \to \SetCat$ by defining $K(S)$ to be the set of functions $x : S \to K[0]$ whose image $\im(x)$ lies in some $K[j]$, and defining $K(f : S' \to S)$ to act by precomposition with $f$. This construction gives rise to a full and faithful embedding of $\SCplx$ into the category $\PSh(\UnSimp)$ of presheaves on $\UnSimp$.
\end{rmk}

\nid Again, we will tacitly Yoneda embedded simplices in $\UnSimp$ into presheaves $\PSh(\UnSimp)$. Note that the Yoneda embedding $\UnSimp \to \PSh(\UnSimp)$ lands in the subcategory $\SCplx \into \PSh(\UnSimp)$, making simplices in particular simplicial complexes.

Recall the `unordering' functor $\Unord {(-)} : \Delta \to \UnSimp$ which forgets orders (see \autoref{term:unord-simp-functor}). This extends to simplicial sets as follows.

\begin{term}[Unordering functor] \label{defn:unordering-functor} The `unordering functor' $\Unord {(-)} : \SSet \to \PSh(\UnSimp)$ forgets the order of vertices in each simplex: formally, this can be defined as the left adjoint to precomposing presheaves with $\Unord {(-)} : \Delta \to \UnSimp$.
\end{term}

\nid Rourke and Sanderson introduce `ordered simplicial complexes' for simplicial complex with an order on their objects \cite[\S 1]{rourke1971delta}; we will use the same term for the following more general notion, which requires each simplex in a simplicial complex to be (consistently) ordered, but this `local' order need not induce a `global' order on the set of objects.

\begin{term}[Ordered simplicial complexes] A simplicial set $X$ is called an `ordered simplicial complex' if its unordering $\Unord X$ is an ordinary simplicial complex.
\end{term}

\nid Ordered simplicial complexes can be more explicitly described: a simplicial set $X$ is an ordered simplicial complex if each non-degenerate simplex in $K$ can be uniquely identified by its set of vertices (as a subset of $K[0]$).

\begin{rmk}[Non-degenerate simplices inject into complexes] A simplex $x \in K[m]$ in an ordered simplicial complex is non-degenerate if and only if the presheaf map $x : [m] \to K$ is an injection. A similar observation holds in the case of simplicial complexes.
\end{rmk}

\begin{term}[The category of ordered simplicial complexes] The `category of ordered simplicial complexes' $\OrdSCplx$ is the full subcategory of $\SSet$ consisting of ordered simplicial complexes. Its morphisms will be called `ordered simplicial maps'.
\end{term}

\begin{term}[Ordering unordered simplicial complexes and their maps] \label{term:orderings} An `ordering' of a simplicial complex $K$ (resp.\ of a simplicial map $F : K \to L$) is a choice of a preimage $\Order K$ of $K$ (resp.\ a preimage $\Order F : \Order K \to \Order L$ of $F$) under the unordering functor. For fixed orderings $\Order K$ and $\Order L$ of simplicial complexes $K$ and $L$, a simplicial map $F : K \to L$ is said to be `order-preserving' if it has a (necessarily unique) ordering $\Order F : \Order K \to \Order L$.
\end{term}

\nid The standard ordering of the unordered standard simplex $\unsimp m$ is the (ordered) standard simplex $[m]$.

\subsubsecunnum{The definition of framings}

We now introduce framings on ordinary simplicial complexes. Recall the definition of embedded frames $\cF$ of $m$-simplices $S$ from \autoref{defn:combinatorial-frame-embedded} which endows $S$ with an isomorphism $S \iso [m]$ together with an $n$-embedded frame $\cF$ of $[m]$ (given by an injection $\cF : \spine[m] \into \bnum n$). Recall also, for a $j$-face $f : T \into S$, such an $n$-embedded frame of $S$ restricts to an $n$-embedded frame of $T$ defined by the restricted isomorphism $\rest {(S \iso [m])} f : T \iso [j]$ and the restricted $n$-embedded frame $\rest \cF f$ (see \autoref{defn:frame-restriction-gen}).

\begin{defn}[Framings of simplicial complexes]  An \textbf{$n$-framing} $(\alpha, \cF)$ of a simplicial complex $K$ endows each $m$-simplex $x : S \into K$ with an $n$-embedded frame $(\alpha_x : S \iso [m], \cF_x)$ such that, for any $j$-face $f : T \into S$, the restriction of the chosen frame of $x$ to the face $f$ coincides with the chosen frame of $x \circ f$; that is, $\alpha_{x \circ f} = \rest {\alpha_x} f$ and $\cF_{x \circ f} = \rest {\cF_x} f$.
\end{defn}

\nid The following observation will simplify the data of framings.

\begin{obs}[Isomorphism data of framings is an ordering] \label{obs:framings-as-orderings} An $n$-framing $(\alpha,\cF)$ of a simplicial complex $K$ gives rise to an ordering of $K$; indeed, for each $m$-simplex $x : S \into K$ in $K$, an order on the vertices of $S$ is determined by the isomorphism $\alpha_x : S \iso [m]$ and together (since choices of $\alpha_x$ are compatible with faces) these orderings of vertices of simplices determine an ordering of $K$ itself. This has an inverse: any ordering of $K$ restricts to an ordering of each simplex $x : S \into K$, and thus yields (compatible) isomorphisms $S \iso [m]$.
\end{obs}

\nid We therefore obtain the following equivalent way of phrasing the notion of $n$-framings.

\begin{rmk}[Framings of simplicial complexes by orderings] \label{defn:framings-of-simp-cplx} To define an $n$-framing $\cF$ of a simplicial complex $K$ we may equivalently specify an ordering of $K$ together with an $n$-embedded frame $\cF_x$ of $[m]$ for each (order-preserving) $m$-simplex $x : [m] \into K$, such that, for any face $f : [k] \into [m]$, we have $\cF_{x \circ f} = \rest {\cF_x} f$.
\end{rmk}

\nid We will henceforth adopt this more concise reformulation of $n$-framings in terms of orderings (as recorded by \autoref{defn:framings-of-simp-cplx}). We refer to the pair $(K, \cF)$, of a simplicial complex and an $n$-framing on it, as an `$n$-framed simplicial complex'.  Note that the ordering of $K$ itself will be implicit in our notation, and we adopt the following convention.

\begin{conv}[Keeping orderings implicit] Given a framed simplicial complex $(K,\cF)$, when considering $K$ as an `ordered simplicial complex' (for instance, when considering maps $[m] \to K$) we will always assume this ordering to be the ordering of $K$ provided by the $n$-framing $\cF$.
\end{conv}

\nid Note that a simplicial complex $K$ cannot contain simplices of dimension greater than $n$ in order for it to admit an $n$-framing.

\begin{eg}[Framings and non-framings on simplicial complexes] \label{eg:compatible-choices-of-frames} In \autoref{fig:compatible-and-incompatible-choices-of-frames}, on top, we depict three simplicial complexes consisting of one, two and three 2-simplices respectively. Underneath each we depict an ordering, underneath which we depict several choices of frames for individual simplices (to depict these we use coordinate frame notation as introduced in \autoref{eg:combinatorial-frames-emb}). Not all choices lead to framings: for the indicated simplices in blue and in red the compatibility of frame choices with face restrictions is not met.
\begin{figure}[h!]
    \centering
    \def\svgwidth{1\columnwidth}
    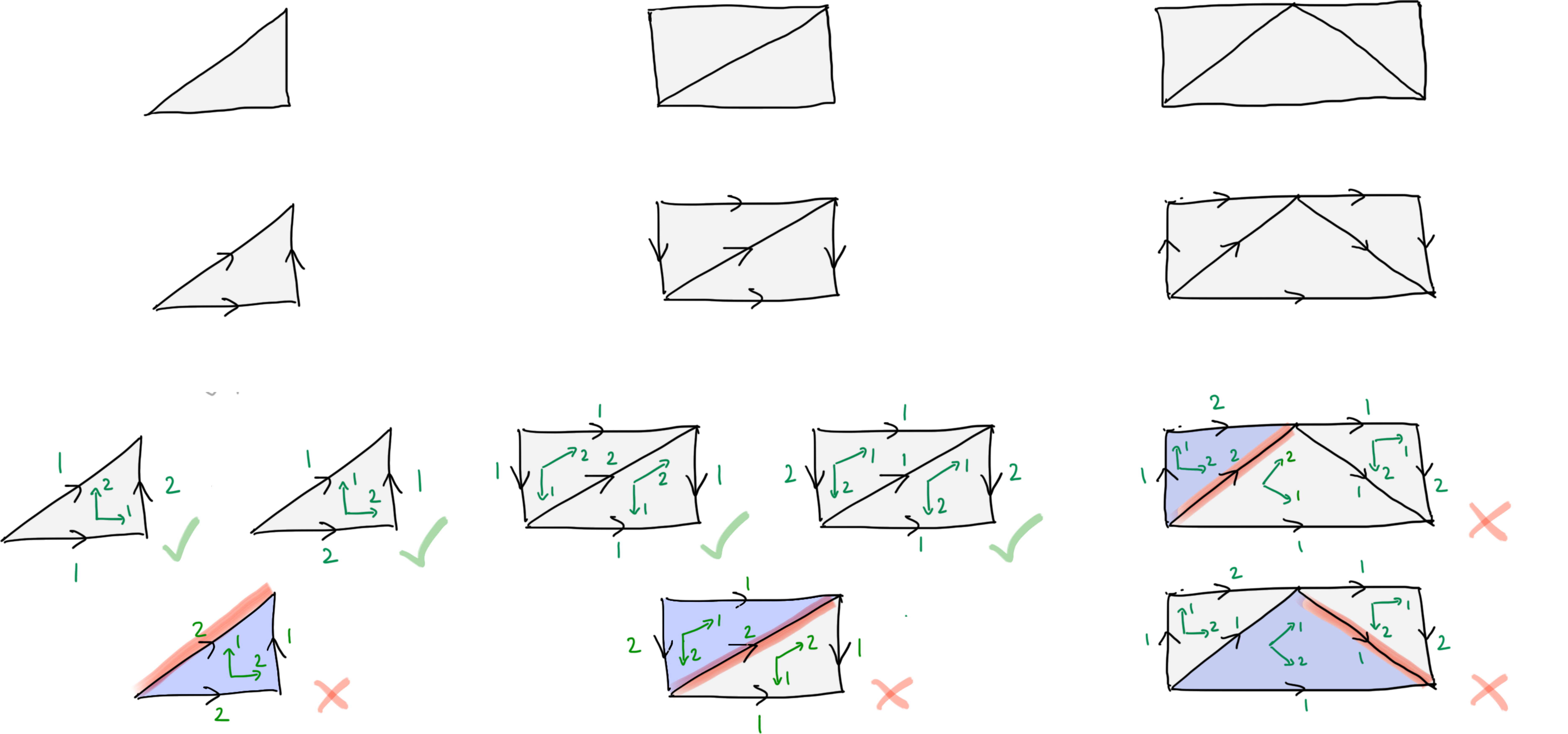

    \caption{Simplicial complexes with orderings, framings and non-framings.}
    \label{fig:compatible-and-incompatible-choices-of-frames}
\end{figure}
\end{eg}

\begin{rmk}[1-Skeletal notation] \label{rmk:1-skeletal-notation} Instead of depicting frames $\cF_x$ separately for each simplex $x$ in $K$ (as in the previous example) it suffices to give frame labels of 1-simplices: indeed, each simplex frame is determined the frames of its spine, and the compatibility condition in the definition of framed simplicial complexes ensures that if two simplices share a spine vector that vector must carry the same label in both of their frames.\footnote{The converse is not true: we cannot `only' label 1-simplices in a simplicial complex with frame labels in $\bnum n$ to define an $n$-framing---we must also check that such a labeling defines valid frames for each $k$-simplex in the complex, $k > 1$.}
\end{rmk}

\begin{eg}[More framings on simplicial complexes] \label{eg:1-skeletal-notation} We depict three framed simplicial complexes in \autoref{fig:further-examples-of-framed-simplicial-complexes}. Observe the difference between the second and third example: the second example realizes to the annulus $S^1 \times D^1$ while the their example realizes to the Mobius band. Note that the framing of the Mobius band, in a sense, is `singular' as it `flips' when traversing the band; we will later describe this as a failure of \emph{local flatness} of the framing (see \autoref{rmk:mobius-band-non-flat}).
\begin{figure}[ht]
    \centering
    \def\svgwidth{1\columnwidth}
    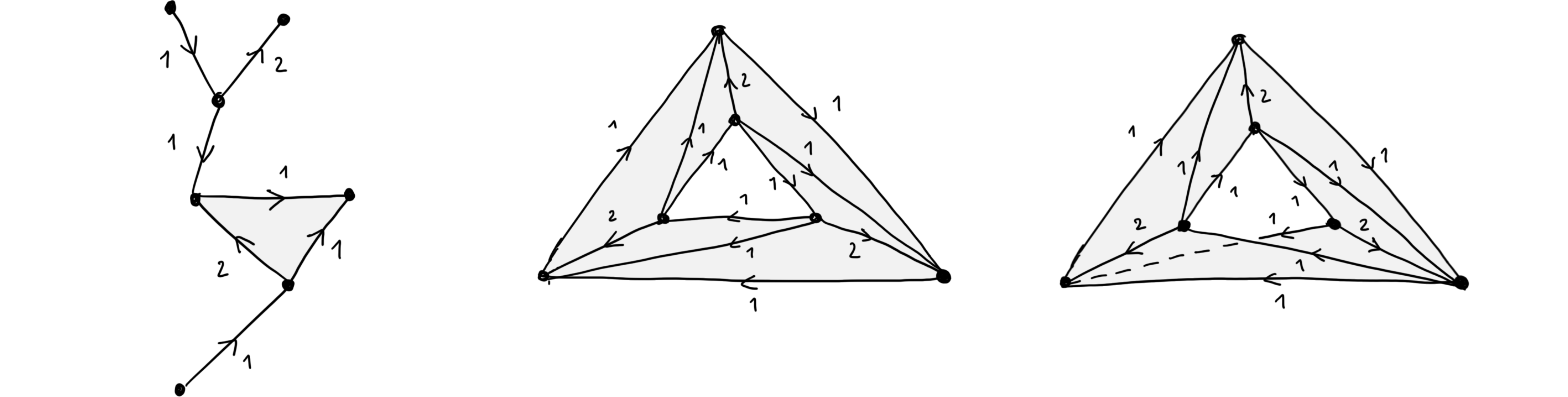

    \caption{Framed simplicial complexes depicted using 1-skeletal notation.}
    \label{fig:further-examples-of-framed-simplicial-complexes}
\end{figure}
\end{eg}

\begin{defn}[Restrictions] \label{rmk:framed-subcomplex} Given an $n$-framed simplicial complex $(K,\cF)$ and a simplicial subcomplex $L \into K$, the \textbf{restriction $\rest \cF L$ of $\cF$ to $L$} is the $n$-framing of $L$ obtained by restricting the ordering of $K$ to $L$, and setting $(\rest \cF L)_x = \cF_x$ for simplices $x : [k] \into L \into K$.
\end{defn}

\pauseae

Framings of simplicial complex can similarly be defined in the `partial' case.

\begin{defn}[Partial framings of simplicial complexes] A \textbf{partial $n$-framing} $(\alpha, \cF)$ of a simplicial complex $K$ endows each $m$-simplex $x : S \into K$ with an $n$-embedded partial frame $(\alpha_x : S \epi [k], \cF_x)$ such that, for any $j$-face $f : T \into S$, the restriction of the chosen $n$-embedded partial frame of $x$ to the face $f$ coincides with the chosen $n$-embedded partial frame of $x \circ f$; that is, $\alpha_{x \circ f} = \rest {\alpha_x} f$ and $\cF_{x \circ f} = \rest {\cF_x} f$.
\end{defn}

\nid Note that a `partial $n$-framed simplicial complex' $(K,\cF)$ may now have $m$-simplices of any dimension $m$. Note also that \autoref{obs:framings-as-orderings} no longer holds: a partial framing need not determine an ordering on a simplicial complex. Going forward, we will be mainly interested in the non-partial case of framings; nonetheless, all subsequent definitions have analogs in the partial case as well.

\begin{eg}[A familiar partial framing] An example of a partial $1$-framing was depicted and discussed earlier in \autoref{fig:the-saddle-singularity}.
\end{eg}

\pauseae

We next define maps of framed simplicial complexes. This straight-forwardly generalizes the notion of framed maps of framed simplices (see \autoref{defn:framed-map-of-framed-simplices}).

\begin{defn}[Framed maps of framed simplicial complexes] \label{defn:framed-maps-of-framed-scplx} Consider $n$-framed simplicial complexes $(K,\cF)$ and $(L,\cG)$. A \textbf{framed simplicial map} $F : (K,\cF) \to (L,\cG)$ (or simply, a `framed map') is a simplicial map $F : K \to L$ that restricts on all simplices $x : [k] \into K$ and $y = \im(F \circ x) : [l] \into L$ to a framed map $F : ([k],\cF_x) \to ([l],\cG_y)$ of $n$-embedded framed simplices.
\end{defn}

\nid In particular, note that any framed map $F : (K,\cF) \to (L,\cG)$ is order-preserving. When referring to the ordering of $F$, we will speak of the `ordered simplicial map' $F : K \to L$.

\begin{notn}[The category of framed simplicial complexes] The category of $n$-framed simplicial complexes and framed maps will be denoted by $\FrSCplx n$.
\end{notn}

\begin{defn}[Unframing framed simplicial complexes] The \textbf{unframing} functor $\Unframe : \FrSCplx n \to \SCplx$ takes a framed simplicial complex $(K,\cF)$ to the simplicial complex $K$, and a framed map $F : (K,\cF) \to (L,\cG)$ to the simplicial map $F : K \to L$.
\end{defn}

\begin{rmk}[Subframed maps of framed simplicial complexes] \label{rmk:subframed-maps-of-framed-simp-cplx} The notion of `subframed maps of framed simplices' described in \autoref{def:subframed} also generalizes to framed simplicial complexes. A `subframed maps of framed simplicial complexes' $F : (K,\cF) \to (L,\cG)$ is a simplicial map $F : K \to L$ that restricts on all simplices $x : [k] \into K$ and $y = \im(F \circ x) : [l] \into L$ to a subframed map $F : (\unsimp k \iso [k],\cF_x) \to (\unsimp l \iso [l],\cG_y)$. Note, unlike framed maps, subframed maps need not be order-preserving.
\end{rmk}

\subsubsecunnum{The definition of proframings} \label{ssec:combinatorial-proframings}

We next introduce proframings on simplicial complexes. The definition takes a `global route to framings': a proframing endows a simplicial complex with a sequence of simplicial surjections defined on the entire complex (we also refer to these surjections as `simplicial projections'). In particular, the definition of $n$-proframings differs from that of $n$-framings in that $n$-proframings do not merely require compatible choices of proframes on each simplex (however, such `local' proframes can be derived from the definition). The `global nature' of proframings will make them a useful tool for describing global properties of framings; most importantly, we will see that `globally flat' framings can be completely understood in terms of certain proframings, which will underlie their later `constructible classification' (in \autoref{ch:classification-of-framed-cells}).

\begin{defn}[Proframings of simplicial complexes] \label{defn:proframings-of-simp-cplx} An \textbf{$n$-proframing} of a simplicial complex $K$ is a ordering of $K$ together with a sequence $\cP$ of ordered simplicial surjections
    \begin{equation}
        K = K_n \xto {p_n} K_{n-1} \xto {p_{n-1}} ... \xto {p_2} [K_1] \xto {p_1} K_0 = [0]
    \end{equation}
such that on each simplex $x : [m] \into K$ the restricted sequence $\rest \cP x$ is an $n$-embedded proframe of $[m]$.
\end{defn}

\nid We will refer to the pair $(K, \cP)$ of a simplicial complex with an $n$-proframing $\cP = (p_n,p_{n-1},...,p_1)$ on it, as an `$n$-proframed simplicial complex'.

\begin{eg}[Proframings of simplicial complexes] In \autoref{fig:proframings-of-two-different-simplicial-complexes} we depict three sequences of simplicial surjections. Each projection $p_i$ is suggested as a geometric projection, but we also highlight affine kernels of $p_i$ on each simplex (in blue). While the first three sequences define proframings, the last sequence fails to do so as it fails to be an $2$-embedded proframe when restricted to its central 2-simplex.
\begin{figure}[h!]
    \centering
    \def\svgwidth{1\columnwidth}
    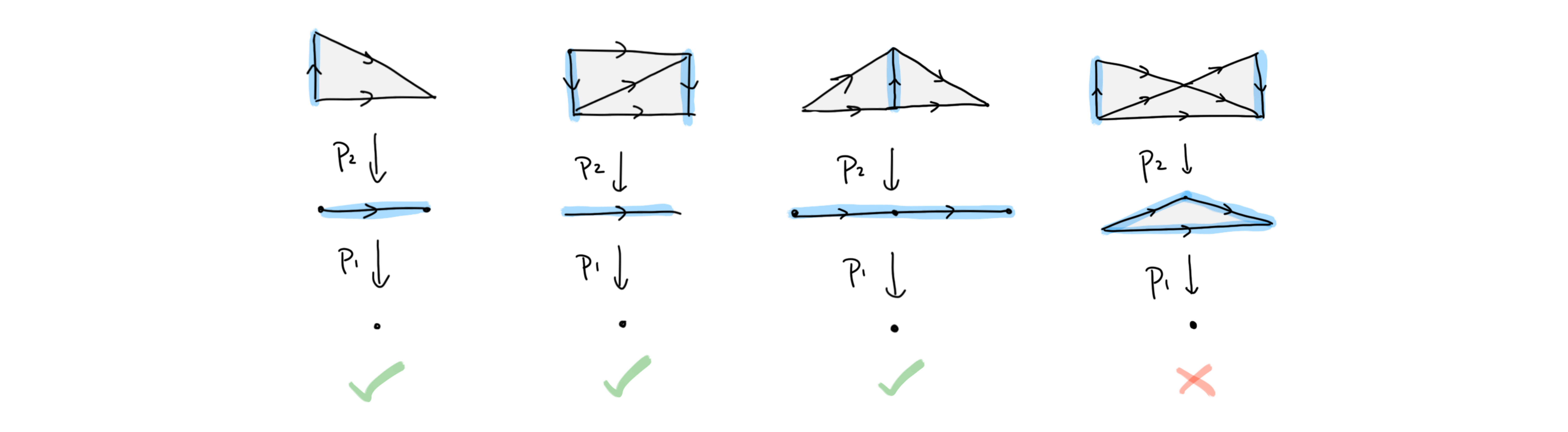

    \caption{Two 2-proframings on a simplicial complex with two 2-simplices.}
    \label{fig:proframings-of-two-different-simplicial-complexes}
\end{figure}
\end{eg}

\pauseae

As in the case of framings, we can further generalize the definition of proframings to allow for partial proframes.

\begin{defn}[Partial proframing] \label{defn:partial-proj-framings} An \textbf{partial $n$-proframing} of a simplicial complex $K$ is a sequence of surjective simplicial maps (all but the first of which are \emph{ordered} simplicial maps)
    \begin{equation}
        K \epi K_n \xto {p_n} K_{n-1} \xto {p_{n-1}} ... \xto {p_2} K_1 \xto {p_1} K_0 = [0]
    \end{equation}
such that, on each $m$-simplex $x : S \into K$ of $K$, the sequence restricts to an $n$-embedded partial proframe of $S$.
\end{defn}

\begin{eg}[A familiar partial proframing] The sequence of ordered simplicial maps depicted earlier in \autoref{fig:the-saddle-singularity} is, in fact, a partial 1-proframing.
\end{eg}

\nid While many subsequent definitions also apply (with evident adjustments) to partial proframings, we will henceforth be interested only in the non-partial case.

\begin{rmk}[Goodbye partiality] This will, in fact, be last time that we will mention partial (pro)framings in this book. Ultimately, neither notion will play a role for us here, and instead, our focus will lie with (combinatorial) geometric objects in which frames span all `tangential directions'. Nonetheless, partial (pro)framings provide a natural and useful generalization of (pro)framings---whose applications, however, go beyond the scope of this present work.
\end{rmk}

\pauseae

We next define maps of proframed simplicial complexes. This generalizes our earlier definition of proframed maps of embedded proframed simplices (see \autoref{defn:proframed-maps}).

\begin{defn}[Maps of proframings] \label{defn:maps-of-proj-framing} Given $n$-proframed simplicial complexes $(K,\cP)$ and $(L,\cQ)$, a \textbf{proframed simplicial map} $F : (K,\cP) \to (L,\cQ)$ is a map of sequences $(F_n,F_{n-1},...,F_1,F_0) : \cP \to \cQ$ (consisting of cellular maps $F_i$) which, on a simplex $x : [k] \into K$ with image $y = \im(F_n \circ x): [l] \into L$, restricts to a proframed map $F : \rest \cP x \to \rest \cQ y$ of $n$-embedded proframed simplices.
\end{defn}

\begin{notn}[Category of proframings] The category of $n$-proframed simplicial complexes and their proframed maps will be denoted by $\PFrSCplx n$.
\end{notn}

\begin{notn}[Truncations of proframings] \label{rmk:framed-cplx-truncation} Given an $n$-proframing $\cP = (K_n \xto {p_n} K_{n-1} \xto {p_{n-1}} ... \xto {p_1} K_0)$ of $K$, its \textbf{$i$-truncation} $\cP_{\leq i}$ (for $i \leq n$) is the $i$-proframing $(K_i \xto {p_i} K_{i-1} \xto {p_{i-1}} ... \xto {p_1} K_0)$ of the simplicial complex $\Unord K_i$. Similarly truncating maps, one obtains functors $(-)_{\leq i} : \PFrSCplx n \to \PFrSCplx i$.
\end{notn}

\begin{defn}[Restrictions of proframings]  Given an $n$-proframing $\cP$ of $K$ and a subcomplex $L \into K$, the \textbf{restriction $\rest \cP L$ of $\cP$} to $L$ is the $n$-proframing of $L$ obtained by restricting the sequence $\cP$ to $L$.
\end{defn}

\subsection{Gradient framings and integral proframings} \label{ssec:gradients-integrability-pro-framings}

\subsubsecunnum{Gradient framings}
From each proframed simplicial complex we can extract a framed simplicial complexes by the following `gradient' operation, which generalizes gradients of framed simplices as described in \autoref{defn:gradient-of-proframed-simplices}.

\begin{defn}[Gradients of proframed simplicial complexes] \label{defn:grad-of-profr-simp-cplx} Given an $n$-proframing $\cP$ of a simplicial complex $K$, the \textbf{gradient framing $\Gradfr \cP$} is the $n$-framing of $K$ with the same ordering as $\cP$, and with an $n$-embedded frame $\Gradfr \cP_x$ on each simplex $x : [m] \into K$ given by the gradient frame $\Gradfr \rest \cP x$.
\end{defn}

\nid The fact that choices of frames in the gradient framing $\Gradfr \cP$ are compatible with face restrictions follows from the compatibility of gradient frames with face restrictions (see \autoref{rmk:framed-faces-are-compatible-with-gradients-and-integration}).

\begin{defn}[Gradients of proframed maps] Given a proframed map $F = (F_n,F_{n-1}, ..., F_1,F_0) : (K,\cP) \to (L,\cQ)$ of $n$-proframed simplicial complexes, the \textbf{gradient framed map} $\Gradfr F : (K,\Gradfr \cP) \to (K,\Gradfr \cQ)$ is the framed map determined by the simplicial map $F_n : K \to L$.
\end{defn}

\begin{term}[The gradient framing functor] \label{notn:gradient-framing-functor} The construction of gradients on proframings and their maps yield the `gradient framing' functor
\begin{equation}
    \Gradfr : \PFrSCplx n \to \FrSCplx n. \qedhere
\end{equation}
\end{term}

\begin{eg}[Proframings and their gradients] \label{more-proframings} We depict three 2-proframings in \autoref{fig:less-simple-proframings} together with their gradients framings; note that the gradient framings recover our examples from \autoref{eg:1-skeletal-notation}. \begin{figure}[ht]
    \centering
    \def\svgwidth{1\columnwidth}
    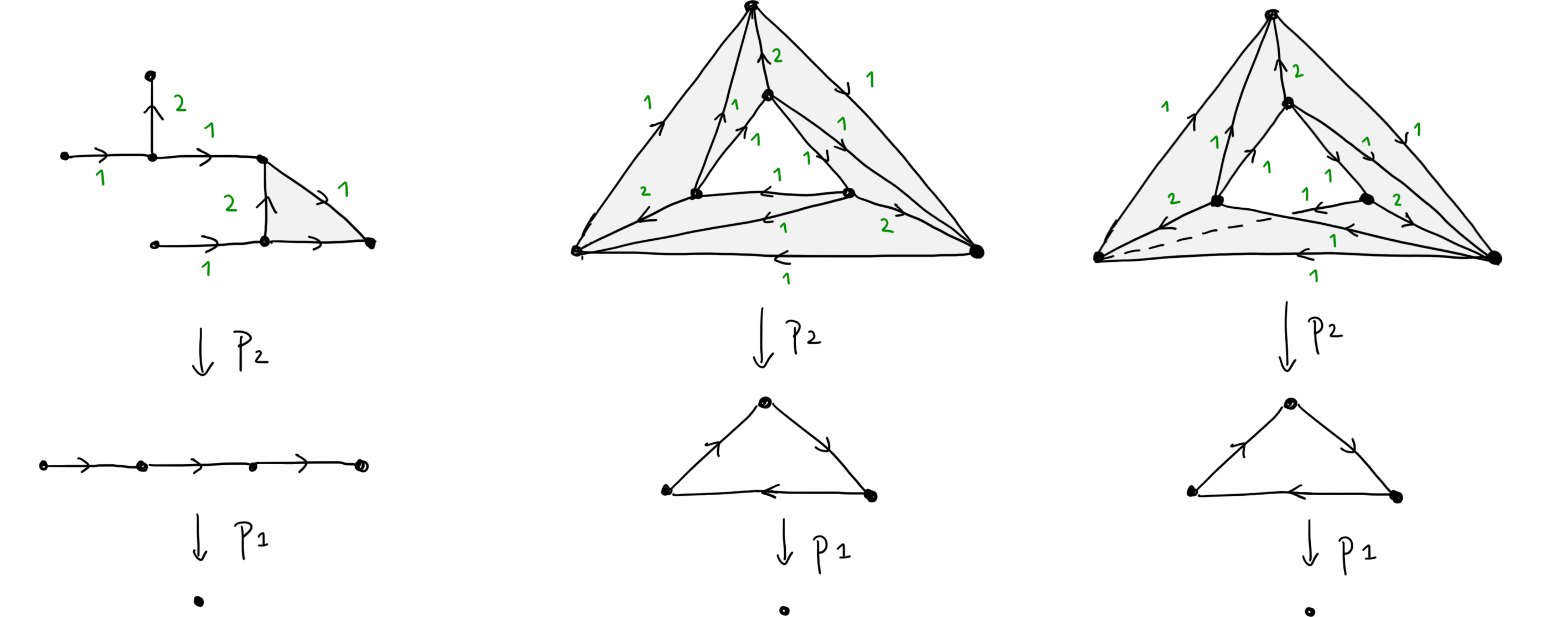

    \caption{Proframed simplicial complexes with their gradient framings.}
    \label{fig:less-simple-proframings}
\end{figure}
\end{eg}

\subsubsecunnum{Integral proframings}

The `converse' to taking gradients is taking integrals. However, not all framings are `integrable'---and if they are, they need not have unique integral proframings. Only in special cases will there be a correspondence between framings and proframings via gradients and integrals (noteworthily, in the case of `flat framings' which we will discuss shortly).

\begin{defn}[Integrating proframings] \label{defn:integrating-proframing} Given an $n$-framed simplicial complex $(K,\cF)$, an \textbf{integral proframing $(K,\cP)$} of $(K_n,\cF)$ is an $n$-proframing whose gradient $n$-framing $\Gradfr \cP$ recovers $\cF$.
\end{defn}

\begin{eg}[A framing with two proframings] \label{eg:proj-frame-inequivalent} In \autoref{fig:two-proframings-that-induce-the-same-choice-of-frames} we depict a 2-framed simplicial complex on the left; to its right we depict two different integral 2-proframings $\cP = (p_2,p_1)$ and $\cQ = (q_2,q_1)$.
\begin{figure}[h!]
    \centering
    \def\svgwidth{1\columnwidth}
    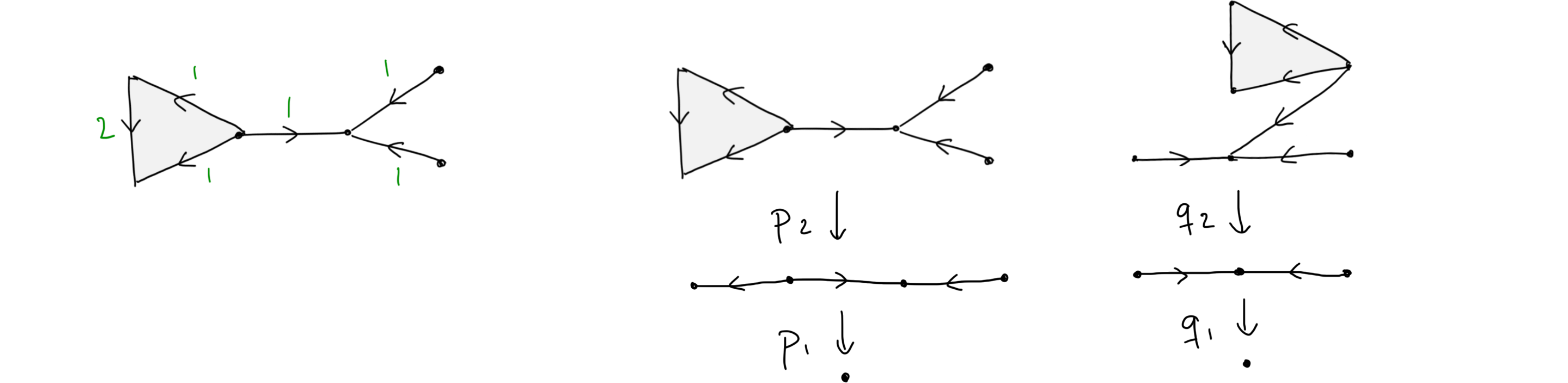

    \caption{A 2-framing with two different integral proframings.}
    \label{fig:two-proframings-that-induce-the-same-choice-of-frames}
\end{figure}
\end{eg}

\begin{eg}[A non-integrable framing] In \autoref{fig:a-framing-without-integrating-proframing} we depict two 2-framings of the boundary $\partial {\unsimp 2}$ of the unordered 2-simplex $\unsimp 2$: neither framing admits an integral proframing.
    \begin{figure}[ht]
    \centering
    \def\svgwidth{1\columnwidth}
    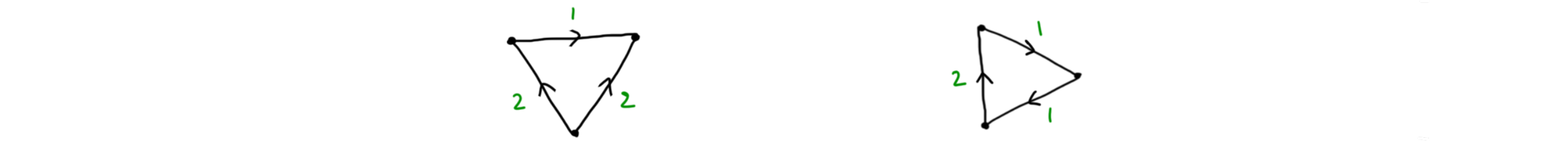

    \caption{Two framings without integral proframings.}
    \label{fig:a-framing-without-integrating-proframing}
\end{figure}
\end{eg}

\nid We remark that the failure of integrability in the previous example may be explained by the following more general observation (which we state without proof).

\begin{rmk}[Integrating simplex boundary framings] An $n$-framing $\cF$ of the unordered simplex boundary $\partial {\unsimp m}$ is integrable if and only if it is the restriction of some $n$-framing $\bar\cF$ of the unordered $m$-simplex $\unsimp m$; that is, $\cF = \rest {\bar\cF} {\partial {\unsimp m}}$.
\end{rmk}

\subsection{Flat framings and flat proframings} \label{ssec:combinatorial-flatness}

In this section we discuss `flat framings': these provide an important class of integrable framings for which, as it will turn out, the correspondence between framings and proframings will be fully restored. Intuitively, flatness can be thought of as requiring framed spaces to `trivialize' as framed subspaces of standard framed euclidean space.

\subsubsecunnum{Flat framings}

We start by generalizing the definition of framed realizations of $n$-embedded framed simplices to the case of framed simplicial complexes (see \autoref{defn:framed-real-emb}).

\begin{defn}[Framed realization] \label{defn:framed-embeddings} Given an $n$-framed simplicial complex $(K,\cF)$, a linear embedding $r : \abs{K} \into \lR^n$ (that is, an embedding that is linear on each simplex) is called a \textbf{framed realization} of $K$ if for each simplex $x : [m] \into K$ the restriction $r \circ \abs{x}$ is a framed realization of $([m],\rest \cF x)$.
\end{defn}

\begin{rmk}[Flat framings are determined by their framed realizations] \label{rmk:flat-framing-via-realization} Since any framed realization $r : \abs{S} \into \lR^n$ of an $n$-embedded framed simplex $(S \iso [m],\cF)$ determines both the isomorphism $S \iso [m]$ and the frame $\cF$ on $[m]$ uniquely, it follows that any framed realization of $r : \abs{K} \into \lR^n$ of an $n$-framed simplicial complex $(K,\cF)$ determines the framing $\cF$ of $K$ uniquely. We may therefore depict framed simplicial complexes by their framed realizations: this is illustrated in \autoref{fig:omitting-arrows-in-flat-framings}, where we depict a framed realization $r : \abs{K} \into \lR^2$ on the left, together with the induced $2$-framing $\cF$ on the right.
\begin{figure}[ht]
    \centering
    \def\svgwidth{1\columnwidth}
    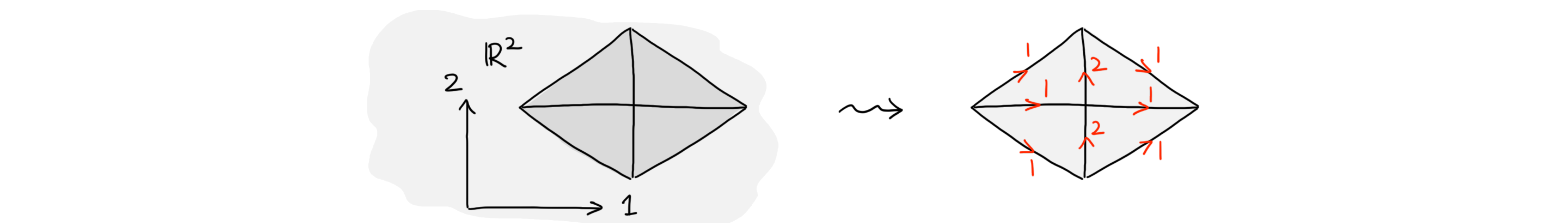

    \caption{Flat framings are determined by their framed realizations.}
    \label{fig:omitting-arrows-in-flat-framings}
\end{figure}
\end{rmk}

\nid We will be interested in framed realizations that are sufficiently regular in that they have `framed bounded' images. The appropriate framed notion of bounds is spelled out in the following definition.

\begin{term}[Framed half-spaces] Recall the projection $\pi_i : \lR^{i} \to \lR^{i-1}$ which forgets the last coordinate of $\lR^{i}$. Given continuous section $\gamma : \lR^{i-1} \to \lR^i$ of $\pi_i$, the `upper half-space' $\lR^n_{\geq \gamma}$ of $\lR^n$ is the subspace of $\lR^n$ of points $x$ such that $\pi_{>i} (x) \geq \gamma \circ \pi_{>i-1} (x)$ in the $\lR$-fiber of $\pi_i$ over $\pi_{>i-1} (x)$ (here, $\pi_{>i}$ abbreviates the composite $\pi_{i+1} \circ ... \circ \pi_{n-1} \circ \pi_n$). Similarly one defines the `lower half-space' $\lR^n_{\leq \gamma}$.
\end{term}

\begin{defn}[Framed bounded subspaces] \label{defn:framed-bounds} A \textbf{framed bounded subspace $S$} of $\lR^n$ is a compact subspace of the form
    \begin{equation}
        S = \bigcap_{i\leq n} (\lR^n_{\geq \gamma^-_i} \cap \lR^n_{\leq \gamma^+_i})
    \end{equation}
    where $\gamma^\pm_i$ are sections of $\pi_i$ such that $\gamma^-_i \leq \gamma^+_i$. We call $\gamma^-_i$ resp.\ $\gamma^-_i$ the `lower' resp.\ `upper $i$th bound' of $S$.
\end{defn}

\begin{eg}[Framed half-spaces and bounded subspaces]
In \autoref{fig:framed-half-spaces-and-bounded-subspace} we depict a set of framed bounds $\gamma^\pm_1$ and $\gamma^\pm_2$ together with their corresponding half-spaces. Intersecting these half-spaces yields the shown framed bounded subspace of $\lR^2$.
\begin{figure}[ht]
    \centering
    \def\svgwidth{1\columnwidth}
    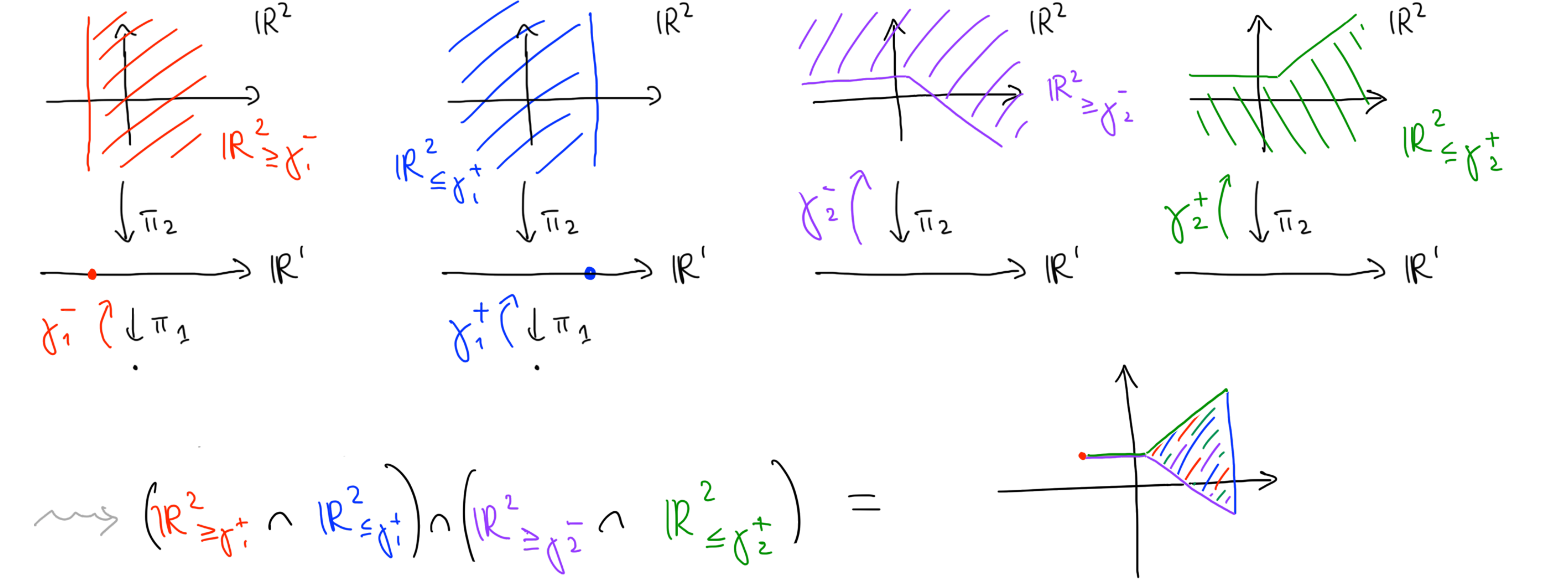

    \caption{Framed half-spaces bounding a framed bounded subspace of $\lR^2$.}
    \label{fig:framed-half-spaces-and-bounded-subspace}
\end{figure}
\end{eg}


\nid We may combine the notions of framed realizations and framed bounds into the following terminology.

\begin{term}[Bounded framed realization]  \label{term:framed-bounded-real} A framed realization of a framed simplicial complex is called a `framed bounded realization' if its image is framed bounded.
\end{term}

\begin{defn}[Flat framings] \label{defn:flat-framings} A framed simplicial complex $(K,\cF)$ is \textbf{flat} if it admits a framed bounded realization $r : \abs{K} \into \lR^n$.
\end{defn}

\begin{notn}[Category of flat framings] Denote the full subcategory of $\FrSCplx n$ consisting of flat framings by $\FlFrSCplx n$.
\end{notn}

\begin{eg}[Flat framings] \label{eg:flat-nonproj-framings} In the upper row of \autoref{fig:flat-framings-simple-examples} we depict examples of flat framings. The first framing admits a framed bounded realization with bounds as given in \autoref{fig:framed-half-spaces-and-bounded-subspace}; one checks that the other framings similarly admit framed bounded realizations. In the lower row of \autoref{fig:flat-framings-simple-examples} we depict non-examples of flat framings; the first framing fails to be flat since no framed realization into $\lR^1$ exists. The second example embeds in $\lR^2$ but its image is not bounded continuously in the sense of \autoref{defn:framed-bounds}. The third framing does not admit a framed realization in $\lR^2$ (note that directions of the left and right vertical 1-simplices are `reversed'). The fourth example again fails to admit a framed bounded realization; while the last example fails to admit a framed realization in $\lR^3$.
\begin{figure}[ht]
    \centering
    \def\svgwidth{1\columnwidth}
    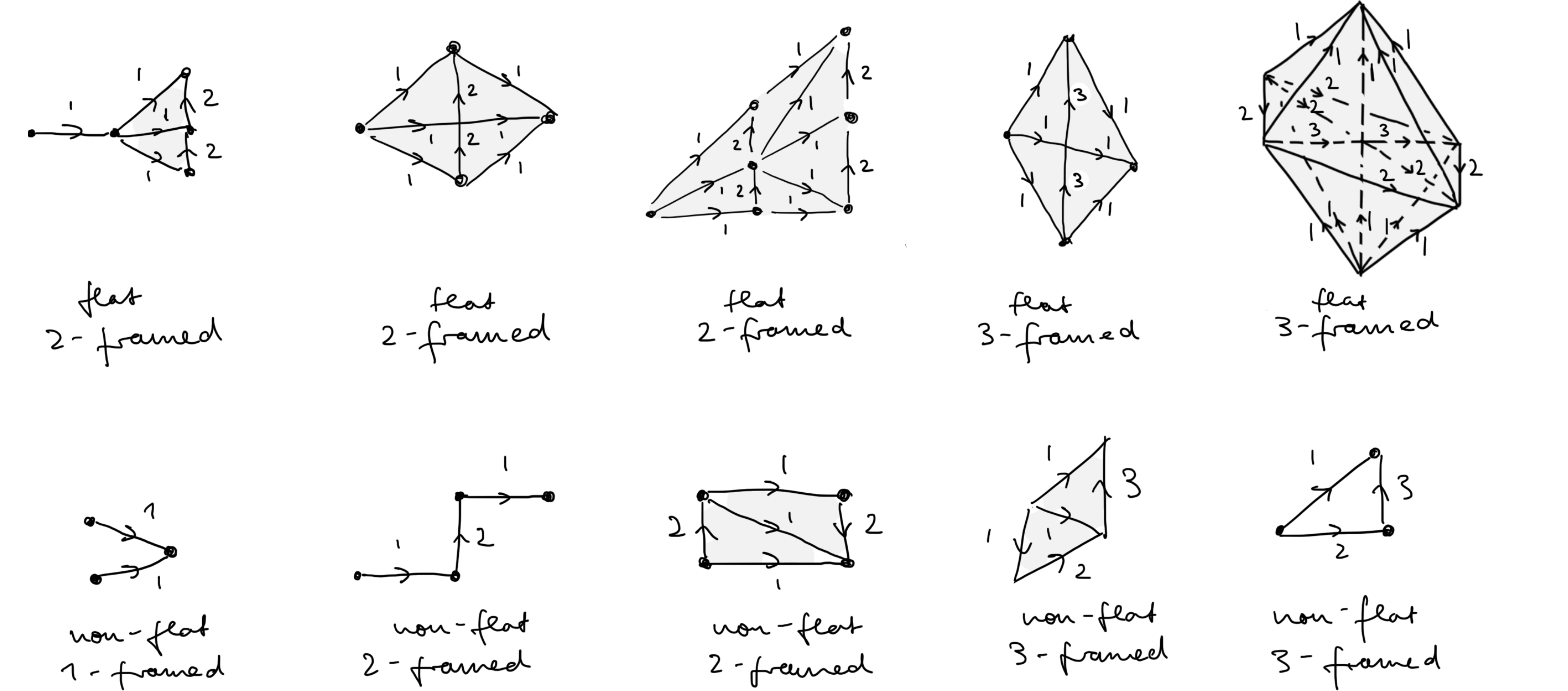

    \caption{Flat and non-flat framed simplicial complexes.}
    \label{fig:flat-framings-simple-examples}
\end{figure}
\end{eg}

\pauseae

We may further impose flatness \emph{locally}. This can be used to recover classical intuition about framings of manifolds as we will see.

\begin{term}[Stars] \label{rmk:stars} Recall the ordinary notion of `stars' $\sstar x$ of vertices $x$ in simplicial complexes $U$: this is defined as the minimal subcomplex $\sstar x$ of $U$ containing all simplices that have $x$ as a vertex.
\end{term}

\begin{defn}[Local flatness] We say an $n$-framing $(K,\cF)$ of simplicial complex $K$ is locally flat if for each vertex $x \in K$ the restricted $n$-framing $(\sstar x, \rest {\cF} {\sstar x})$ is flat.
\end{defn}

\begin{rmk}[Flatness implies local flatness] It is true that each flat framing is in particular locally flat.
\end{rmk}

\begin{eg}[Locally flat framed triangulated 1-manifolds] \label{eg:locally-flat-triangulations} In \autoref{fig:locally-flat-framed-triangulated-1-manifolds} on the left we depict two locally flat 1-framed simplicial complexes $(K,\cF)$ each triangulating a connected 1-manifold. On the right we depict two framed simplicial complex that are not locally flat framed.
\begin{figure}[ht]
    \centering
    \def\svgwidth{1\columnwidth}
    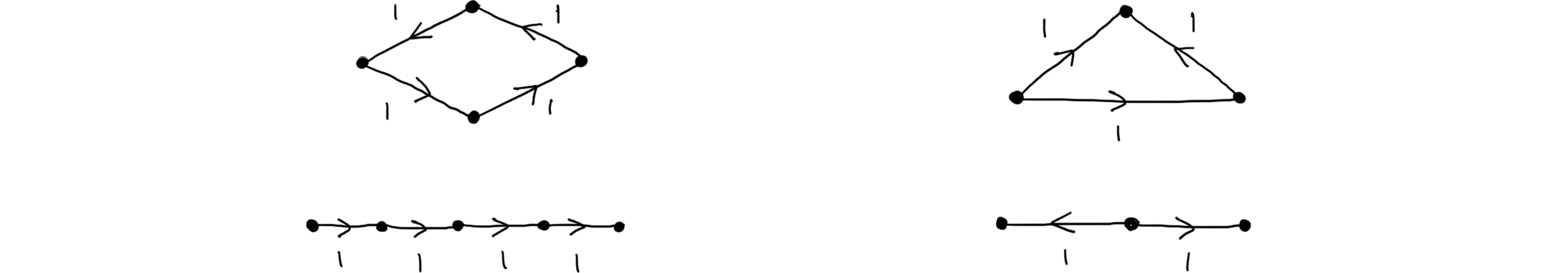

    \caption[Locally flat framed triangulated 1-manifolds.]{Examples and non-examples of locally flat framed triangulated 1-manifolds.}
    \label{fig:locally-flat-framed-triangulated-1-manifolds}
\end{figure}
\end{eg}

\nid We can now revisit our earlier example of the Mobius band, adding the following insight.

\begin{rmk}[Frameability of Mobius band revisited] \label{rmk:mobius-band-non-flat} Recall the last two framed simplicial complexes from \autoref{fig:further-examples-of-framed-simplicial-complexes}: the first example is locally flat, while the second (the `Mobius band') isn't. In fact, the Mobius band is not `locally flat 2-frameable', in the sense that no triangulation admits a locally flat framing.
\end{rmk}

\subsubsecunnum{Flat proframings} In this section we discuss an alternative characterization of flatness in the context of proframings. The two definitions of flatness will turn out to coincide in a precise sense, and this comparison will underlie the later classification of framed combinatorial structures in \autoref{ch:classification-of-framed-cells}.

\begin{term}[Section and spacers simplices] \label{term:sections-spacers-in-proframings} Let $(K,\cP)$ be an $n$-proframed simplicial complex, and pick a simplex $x : [m] \into K_i$ in some $K_i$, $1 \leq i \leq n$. Assume the restriction $\rest {\cP_{\leq i}} x$ of the truncation $\cP_{\leq i}$ is an $i$-proframe of $[m]$ of the form $(p_i,...,p_1)$. We call $x$ a `section simplex' of $\cP$ if $p_i = \id$, and a `spacer simplex' if $p_i \neq \id$.
\end{term}

\begin{term}[Upper and lower sections of spacers] \label{defn:framed-upper-lower-sections} Consider an $n$-proframed simplicial complex $(K,\cP)$ and a spacer simplex $x : [m] \into K_i$. The top projection $p_i$ of the restricted truncated proframe $\rest {\cP_{\leq i}} x = (p_i,...,p_1)$ equals a $j$th degeneracy map $s_j : [m] \to [m-1]$ ($1 \leq j \leq  m$). Precomposing with simplicial face maps $d_k : [m-1] \into [m]$, we obtain the $(m-1)$-simplices $\partial_- x := x \circ d_{j+1} : [m-1] \into K_i$, called the `lower section of $x$ in $\cP$', and $\partial_+ x := x \circ d_j : [m-1] \into K_i$, called the `upper section of $x$ in $\cP$'.
\end{term}

\nid Observe that, given a spacer simplex $x$ in a proframed simplicial complex $(K,\cP)$, its upper and lower sections $\partial_\pm x$ are indeed section simplices in $(K,\cP)$.

\begin{eg}[Sections and spacers] In \autoref{fig:sections-and-spacer-in-framed-simplicial-complexes} we illustrate 2-proframed simplicial complex $K = K_2 \to K_1 \to K_0$. We highlight two spacer simplices in $K_2$ (in blue), and their respective upper and lower sections (in red).
\begin{figure}[h!]
    \centering
    \def\svgwidth{1\columnwidth}
    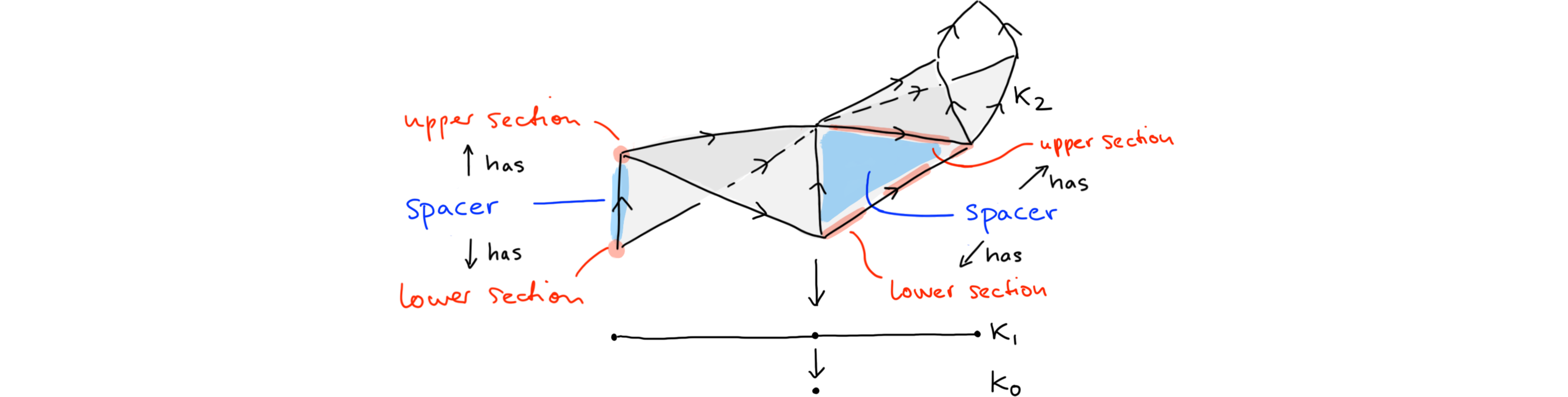

    \caption{Sections and spacers in a proframed simplicial complex.}
    \label{fig:sections-and-spacer-in-framed-simplicial-complexes}
\end{figure}
\end{eg}

\begin{term}[Fiber complexes and categories] \label{defn:fundamental-fiber-category} Consider an $n$-proframed simplicial complex $(K,\cP)$, with $\cP = (p_n,...,p_1)$, and a $k$-simplex $z : [k] \into K_{i-1}$.
    \begin{enumerate}
        \item The `fiber' $K_z$ over $z$ is the subset of non-degenerate simplices in $K_{i}$ that are mapped to $z$ by $p_{i}$.
        \item The `fiber complex' $\FG K_z$ of the fiber over $z$ is the ordered simplicial complex whose vertices are sections $x \in K_z$ over $z$ and whose 1-simplices $y : x_- \to x_+$ are spacers $y \in K_z$ with $\partial_\pm y = x_\pm$.
        \item The `fiber category' $\FC K_z$ is the free category on the complex $\FG K_z$ (with vertices as objects, and edges as generating morphisms). \qedhere
    \end{enumerate}
\end{term}

\begin{constr}[Transition functors of fiber categories] \label{rmk:compatibility-fiber-orders} For an $n$-proframed simplicial complex $(K,\cP)$, consider simplices $z : [k] \into K_{i-1}$ and $w : [l] \into K_{i-1}$ such that $w$ is a face of $z$ (that is, $w$ factors through $z$ as a face $[l] \into [k]$). Note that each simplex $x \in K_z$ restricts over $w$ to a simplex $\rest x {w \subset z} \in K_w$: this restriction takes sections to sections, while spacers restrict either to spacers or to sections. The restriction thus induces the `transition functor' $\rest{-}{w \subset z} : \FC K_z \to \FC K_w$.
\end{constr}

\begin{eg}[Fiber categories and transition functors] In \autoref{fig:fiber-categories-and-transition-functors}, for the  proframing given in the previous \autoref{fig:sections-and-spacer-in-framed-simplicial-complexes}, we colored each simplex in $K_1$, and depicted its fiber category in the corresponding color (with a point for each object, and an arrow for each morphisms). We also indicated the transition functors between them (by mapping arrows `$\mapsto$').
\begin{figure}[ht]
    \centering
    \def\svgwidth{1\columnwidth}
    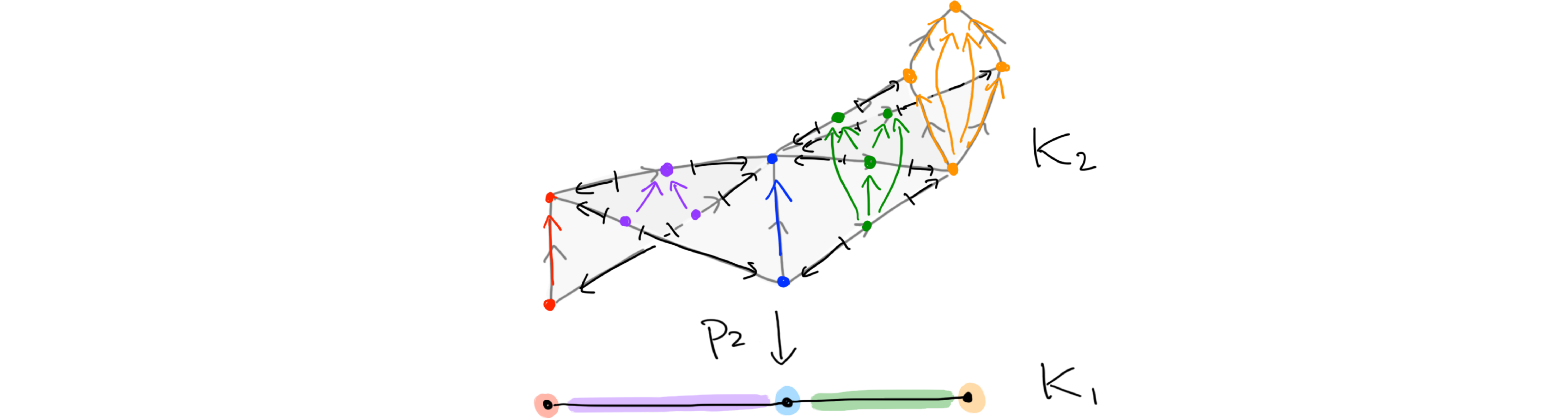

    \caption{Fiber categories and transition functors.}
    \label{fig:fiber-categories-and-transition-functors}
\end{figure}
\end{eg}

\nid Fiber transition functors relate fiber categories across fibers. In order to define flat proframings, we will be interested in the case where fibers are `linear', and transitions functors glue fibers onto fibers, as follows.

\begin{defn}[Flat proframings] \label{defn:flat-proframings}
    An $n$-proframed simplicial complex $(K,\cP)$ is said to be \textbf{flat} (or to have `flat proframing') if the following holds.
    \begin{enumerate}
        \item \textit{Fibers are linear}: For any simplex $z : [k] \into K_{i}$, the fiber category $\FC K_z$ is a total order.
        \item \textit{Fiber transition is endpoint-preserving}: For simplices $w \subset z$ in $K_{i}$ the transition functor $\rest{-}{w \subset z}$ is endpoint-preserving (meaning it preserves least and greatest elements as a map of total orders). \qedhere
    \end{enumerate}
\end{defn}

\begin{term}[Linear complexes] \label{rmk:lin-complex} The `finite linear simplicial complex with $j$ edges' $\lincplx j$ is the ordered simplicial complex with $j+1$ vertices $\Set{0,1,...,j}$ and $j$ directed 1-simplices $i \to i+1$.
\end{term}

\begin{obs}[Fiber complexes in flat proframings are linear] \label{rmk:linear-fiber-complexes-in-flat-proframings} The first condition in \autoref{defn:flat-proframings} is equivalent to the condition that all fiber complexes $\FG K_z$ are linear.
\end{obs}

\begin{notn}[Categories of flat proframings] Denote the full subcategory of $\PFrSCplx n$ consisting of flat proframings by $\FlPFrSCplx n$.
\end{notn}

\begin{eg}[Flat proframings] In \autoref{fig:euclidean-vs-non-euclidean-framed-simplicial-complexes} we depict six 2-proframings (in each case highlight in blue the kernels of projection). The first example is not flat since it has non-linear fiber categories; the third example is not flat since its transition functors are not endpoint-preserving; the fifth example is not flat since, again, it has non-linear fiber categories (this time already at level 1 of the sequence).
\begin{figure}[h!]
    \centering
    \def\svgwidth{1\columnwidth}
    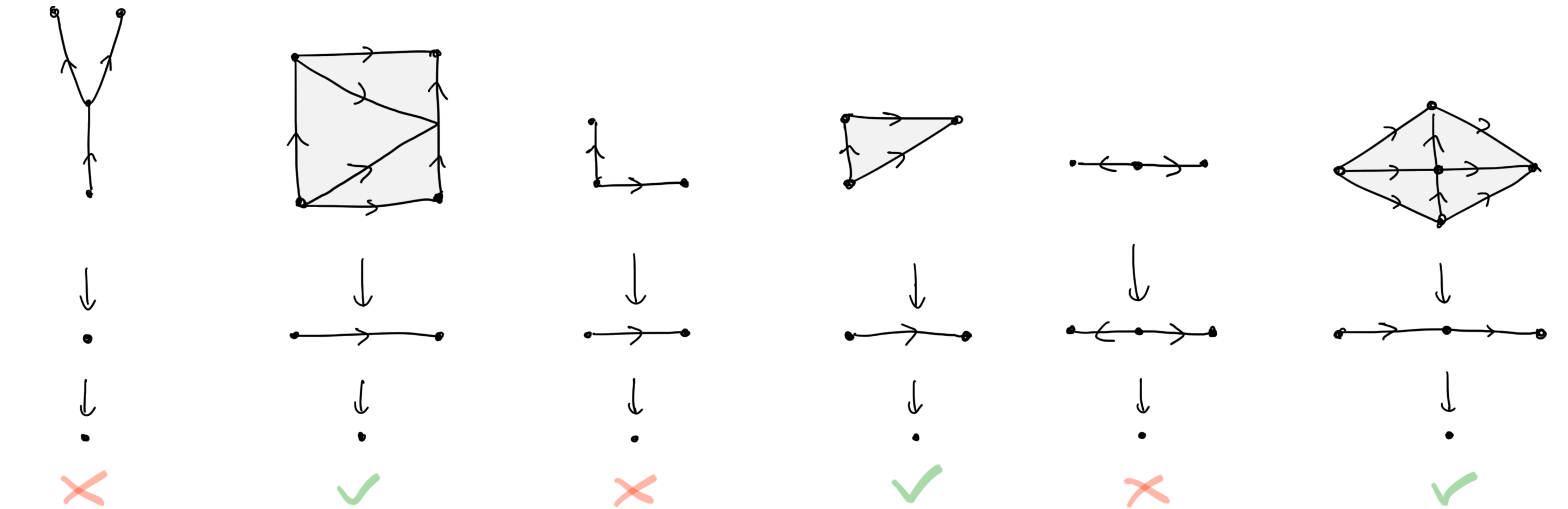

    \caption{Flat and non-flat proframings of simplicial complexes.}
    \label{fig:euclidean-vs-non-euclidean-framed-simplicial-complexes}
\end{figure}
\end{eg}

\pauseae

Finally, just as flat $n$-framed simplicial complexes have framed realizations, flat proframed simplicial complexes `proframed realizations'. These now embed in the standard euclidean $n$-proframe $\Pi = (\pi_n,\pi_{n-1},...,\pi_1)$ (see \autoref{defn:standard-proframe}).

\begin{rmk}[Flat proframings admit proframed realizations] \label{rmk:flat-proframing-via-realization} Any flat $n$-proframed simplicial complex $(K,\cP)$ admits a `proframed realization', by which we mean an embedding of sequences $r = (r_n,r_{n-1},...,r_1,r_0) : \abs{\cP} \into \Pi$, with components $r_i : \abs{K_i} \into \lR^i$, which restricts on each simplex $x : [m] \into K$ to a proframed realization in the sense of \autoref{defn:proframed-real-emb}. As in the case of framed realization, note that any proframed realization of $(K,\cP)$ determines the proframing $\cP$ of $K$ uniquely. A proframed realization of a $3$-proframed simplicial complex is illustrated in \autoref{fig:a-3-framed-simplicial-complex-with-flat-framing}.
\begin{figure}[h!]
    \centering
    \def\svgwidth{1\columnwidth}
    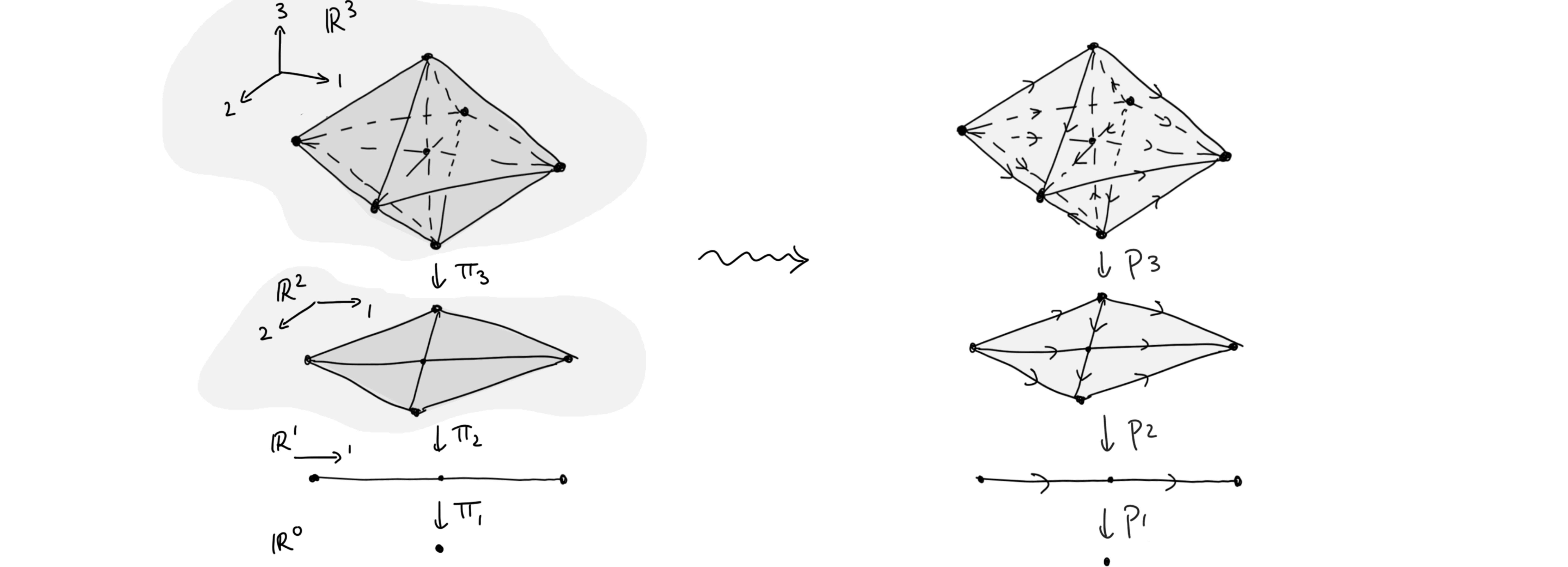

    \caption{Flat proframings are determined by their proframed realizations.}
    \label{fig:a-3-framed-simplicial-complex-with-flat-framing}
\end{figure}
\end{rmk}

\subsubsecunnum{Equivalence of flat framings and flat proframings}

Flatness fully restores the correspondence of framings and proframings. We record this in several statements below (all of which will be proven later in \autoref{ch:classification-of-framed-cells}).

\begin{prop}[Gradients of flat proframings are flat] \label{prop:grad-of-flat-proframing-is-flat-framing} Given a flat $n$-proframed simplicial complex $(K,\cP)$ its gradient $n$-framing $\Gradfr \cP$ is a flat $n$-framing.
\end{prop}

\begin{prop}[Flat framings have flat integral proframings] \label{prop:flat-framings-have-integrating-flat-proframing} Every flat $n$-framing $(K,\cF)$ has a essentially unique integral flat $n$-proframing $\int \cF$.
\end{prop}

\begin{thm}[Correspondence of flat proframings and flat framings] \label{obs:frames-vs-proj-frames-flat} The gradient framing functor restricts to an equivalence of categories
    \begin{equation}
        \Gradfr : \FlPFrSCplx n \eqv \FlFrSCplx n \quad .
    \end{equation}
The inverse to $\Gradfr$ will be called `\emph{integration}', and denoted by $\Intfr$.
\end{thm}

\nid The theorem parallels (and, in fact, generalizes) the correspondence of framed simplices and proframed simplices that we've seen in \autoref{cor:iso-of-frame-vs-proj-frame-simplex}. The proofs of \autoref{prop:grad-of-flat-proframing-is-flat-framing}, \autoref{prop:flat-framings-have-integrating-flat-proframing} and \autoref{obs:frames-vs-proj-frames-flat} will be deferred until \autoref{ssec:frames-vs-proj-frames-flat}.

\section{Framed regular cell complexes} \label{sec:framed-reg-complexes}

In the final section of this chapter, we will introduce $n$-framings on regular cell complexes. Our goal will be the definition of a category of `$n$-framed regular cell complexes' $\FrCCplx n$, fitting into the following diagram of categories (in which vertical arrows are fully faithful embeddings of categories, while horizontal arrows forget framing structures)
\begin{equation}
    \begin{tikzcd} [baseline=(W.base)]
        \FrSCplx n \arrow[r, "\Unframe"] \arrow[d, hook] & \SCplx \arrow[d, hook] \\
        \FrCCplx n \arrow[r, "\Unframe"]                 & |[alias=W]| \CCplx
    \end{tikzcd} \quad .
\end{equation}

\begin{figure}[ht]
    \centering
    \def\svgwidth{1\columnwidth}
    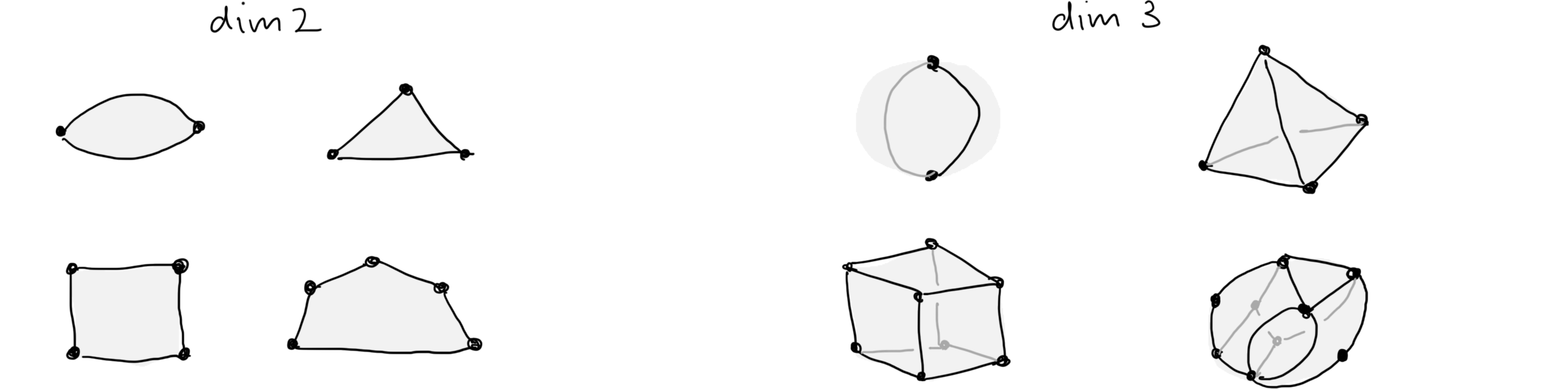

    \caption{An illustration of regular 2-cells and 3-cells.}
    \label{fig:an-illustrutation-of-regular-2-cells-and-3-cells}
\end{figure}

Regular cell complexes, as we recall in \autoref{ssec:reg-cell-complex}, are complexes of cells of general `polytopic shape': in \autoref{fig:an-illustrutation-of-regular-2-cells-and-3-cells} we illustrate shapes of regular cells in dimension 2 and 3.
Despite their generality, regular cell complexes, unlike their non-regular counterpart, are \emph{combinatorializable}. The fundamental property of regular cell complexes that enables this combinatorialization is the `homotopical triviality' of cell attachments (see \autoref{rmk:reg-cell-cplx-homotopically-trivial}). This entails that one can describe a class of so-called `cellular' posets, which are exactly the entrance path posets of regular cell complexes; in fact, geometric realizations of cellular posets recover the homeomorphism type of their corresponding regular cell complexes. The resulting translation between regular cell complexes and cellular posets provides the claimed combinatorialization of regular cell complexes. Crucially however, the combinatorialization is (computably) \emph{intractable} in the following sense: given a poset, there can be no general algorithm to determine if that poset is cellular. In particular, it is impossible to algorithmically write down a list classifying `all the shapes' of regular cells (such as those in \autoref{fig:an-illustrutation-of-regular-2-cells-and-3-cells}) up to some general bound in, say, the number boundary cells. By \emph{framing} regular cells this intractibility will find a natural resolution.

The definition of framed regular cells will directly rely on our previous work on framings of simplicial complexes. Namely, an $n$-framed regular cell will be an $n$-framing of the simplicial complex that underlies the cellular poset of the cell (with the additional condition that the framing is flat on the cell itself as well as on each subcell in its boundary). Our approach therefore combines the following two ingredients: firstly, use to the correspondence of regular cells and cellular posets to endow regular cells with canonical `piecewise linear' structure; secondly, require simplicial framings on cells to be flat. The second condition assures that framings of cells are `trivial' (in the sense that the framed cell is framed realizable as a subspace of $\lR^n$) which, heuristically, reflects that framed regular cells, just like framed simplices, will play the role of `small', and therefore `trivializable' framed pieces from which larger framed spaces will be built.

    In contrast to the case of nonframed regular cells, framed regular cells can now be computably recognized and classified. The classification will be constructed in \autoref{ch:classification-of-framed-cells}. At the same time, the generality of framed regular cells, as opposed to mere framed simplices, will be at the very heart of fundamental results in framed combinatorial topology---for instance, working with framed cells will enable the construction of `canonical cellulations' of flat framed stratifications, as we will explain in \autoref{ch:hauptvermutung} (see also \autoref{eg:generality-of-blocks} in this section). This will highlight that the passage from framed simplices to framed regular cells is not `by choice', but of central importance in framed combinatorial topology.

    We outline this section. After recalling regular cell complexes and related notions in \autoref{ssec:reg-cell-complex}, we will first introduce notions of both framings and proframings of regular cells in \autoref{ssec:framed-reg-cells}. In \autoref{ssec:framed-reg-cell-complex} we will obtain a notion of `framed space' by defining framings (resp.\ proframings) on regular cell complexes by compatibly endowing each cell with a framing (resp.\ considering sequences of surjective cellular maps). We will further discuss an important special case of framed regular cell complexes that are themselves flat, which will later resurface in the construction of canonical cellutation of flat framed stratifications. A diagram summarizing the categories defined in this chapter can be found at the end of the section in \autoref{fig:diagrams-chapter-1}.

\subsection{Regular cell complexes} \label{ssec:reg-cell-complex}

We start by recalling definitions and facts about regular cell complex. Our subsequent goal will be the introduction of a `combinatorial' category $\CCplx$ of regular cell complexes, generalizing that of simplicial complexes $\SCplx$.

\subsubsecunnum{Regular cell complexes as cellular posets}

\begin{defn}[Regular cell complexes] A \textbf{regular cell complex} is a cell complex in which the closure of each cell is a closed ball.
\end{defn}

\begin{term}[Cellulations] In analogy to the notion of `triangulation', we speak of `cellulation' when decomposing a space into a regular cell complex.
\end{term}

\nid When referring to cells in cell complexes we will say `closed cells' to mean closures of cells and `open cells' to mean cell interiors (that is, cells without their boundaries). A simplicial complex is a simple type of regular cell complex in which each $n$-dimensional cell has exactly $(n+1)$ faces of dimension $(n-1)$.

\begin{defn}[Entrance path posets of regular cell complexes] \label{defn:face-poset} The \textbf{entrance path poset} $\Entr X$ of a regular cell complex $X$ is the poset whose objects are the open cells $x$ in $X$, with an arrow $x \to y$ whenever the closure $\overline x$ contains $y$.
\end{defn}

\nid Entrance path posets of regular cell complexes are graded by dimension, that is, they admit a functor $\dim : \Entr X \to \lN\op$ with discrete preimages, mapping each cell to its dimension. Cells of a regular cell complex $X$ which are minimal elements in $\Entr X$ will be called `facets'---these are exactly cells which are not contained in any other cell's boundary.

It will further be useful to think of regular cell complexes as stratified (i.e.\ `singular') spaces.

\begin{defn}[Regular cell complexes as stratified spaces] \label{obs:reg-cell-cplx-as-strat} The \textbf{cell stratification} $\cellstrat X$ of a regular cell complex $X$ is the stratification of the space underlying $X$ whose strata are the open cells of $X$.
\end{defn}

\nid Note that the definition of entrance path posets of regular cell complexes is consistent with that of stratified spaces (see \autoref{defn:entr}), that is, $\Entr X = \Entr {(\cellstrat X)}$. In particular, the characteristic function $\cellstrat X : X \to \Entr X$ (see \autoref{term:characteristic_function}) takes points in the open cell $x$ of $X$ to the object $x \in \Entr X$. As an aside, we remark that regular cell complex have the following properties as stratifications.

\begin{rmk}[The homotopical triviality of cell attachments] \label{rmk:reg-cell-cplx-homotopically-trivial} Regular cell complexes are, as stratified spaces, \emph{conically} stratified (see \autoref{prop:regular-cell-cplx-conical}). They thus have `entrance path $\infty$-categories' (see \autoref{defn:tentr}) which is a higher categorical analog of entrance path posets. Importantly, reflecting the `homotopical triviality' of cell attachments, the entrance path $\infty$-categories of regular cell complexes are in fact $0$-truncated, and thus equivalent to their entrance path posets (see \autoref{lem:reg_use_Entrz}).
\end{rmk}

Understanding regular cell complexes as stratifications means they inherit the following notion of maps. Recall a stratified map is a map that maps strata into strata; equivalently, this means stratified maps factor through characteristic maps by a map of entrance path posets (see \autoref{defn:strat_maps}).

\begin{defn}[Maps of regular cell complexes] A \textbf{map of regular cell complexes} $F : X \to Y$ is a stratified map $F : \cellstrat X \to \cellstrat Y$.
\end{defn}

\begin{notn}[Regular cell complexes as stratified spaces] We henceforth notationally identify $X$ and $\cellstrat X$, thinking of regular cell complexes in terms of their cell stratifications.
\end{notn}

\begin{rmk}[Functoriality of entrance path posets] \label{rmk:entr-face-functorial} Note that entrance path poset construction is functorial: for any map of regular cell complexes $F : X \to Y$ we obtain a poset map $\Entr(F) : \Entr(X) \to \Entr(Y)$, mapping a cell of $X$ to the cell in $Y$ that contains the image $F(x)$.
\end{rmk}

\pauseae

While `entrance path posets' extract posets from regular cell complexes (and more generally, from stratifications) we can, conversely, turn posets into stratifications by the following notion of `classifying stratifications'.

\begin{term}[Upper and strict upper closures] Given a poset $P$ and an element $x \in P$, then the `strict upper closure' $P^{> x}$ of $x$ in $P$ is the full subposet with objects  $y \in P$ with $y > x$. Similarly, the `upper closure' $P^{\geq x}$ is the full subposet of objects $y \in P$ with $y \geq x$.
\end{term}

\begin{term}[Classifying spaces of posets] Recall, the `classifying space' $\abs{P}$ of a poset $P$ is the geometric realization of the nerve $NP$ of $P$. (Abusing terminology, we also refer to $\abs{P}$ itself as the `geometric realization' of $P$, and say $P$ `realizes' to $\abs{P}$.)
\end{term}

\begin{term}[Classifying stratification of posets] \label{term:classifying-stratifications} Given a poset $P$, the `classifying stratification' $\CStr {P}$ is the stratification of $\abs{P}$ whose strata are given subspaces $\abs{P^{\geq x}} \setminus \abs{P^{> x}}$ for $x \in P$; its characteristic function, denoted by $\CStr {P} : \abs{P} \to P$, takes points in $\abs{P^{\geq x}} \setminus \abs{P^{> x}}$ to $x$. \end{term}

\begin{term}[Classifying stratified maps of posets maps] \label{term:classifying-stratifications-maps} Given a poset map $F : P \to Q$, the `classifying stratified map' $\CStr F : \CStr P \to \CStr Q$ maps vertices $p \in \abs{P}$ to vertices $F(p) \in \abs{Q}$, and linearly extends this mapping to all other simplices in $\abs{P}$.
\end{term}

\nid An alternative phrasing of classifying stratifications, based on an explicit definition of classifying spaces, can be found in \autoref{constr:classstrat} (for maps, see also \autoref{constr:classstrat-functor}).

\begin{eg}[Classifying stratifications] We illustrate three classifying stratifications in \autoref{fig:the-classifying-stratification-of-the-569xsimplex}: from left to right, we depict the classifying stratifications $\CStr P$ of the 1-simplex $P = [1]$, the $2$-simplex $P = [2]$ and the product of two 1-simplices $P = [1] \times [1]$.
\begin{figure}[ht]
    \centering
    \def\svgwidth{1\columnwidth}
    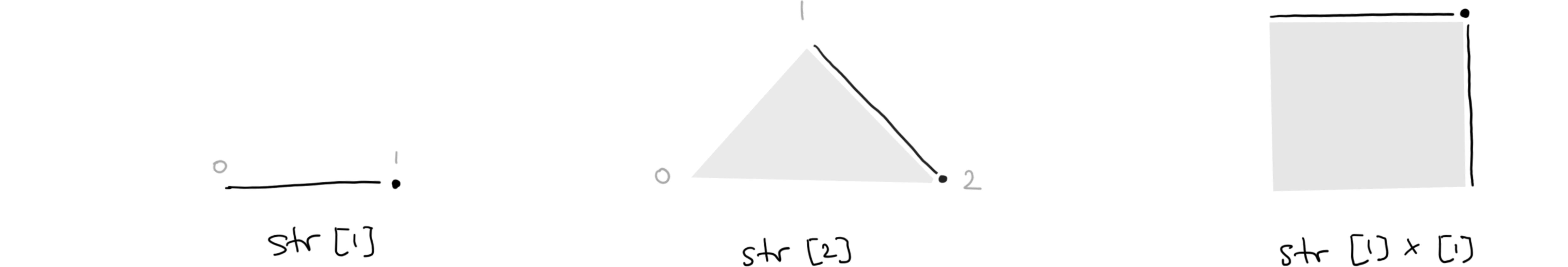

    \caption{Classifying stratifications of posets.}
    \label{fig:the-classifying-stratification-of-the-569xsimplex}
\end{figure}
\end{eg}

The class of posets that can be obtained as entrance path posets of regular cell complexes will be described by so-called `cellular posets'

\begin{defn}[Cellular posets] \label{defn:sph-poset} A poset $(X,\leq)$ is called \textbf{cellular} if the realization $\abs{X^{>x}}$ of the strict upper closure of any $x \in P$ is homeomorphic to a sphere.
\end{defn}

\begin{eg}[Cellular and non-cellular posets] In \autoref{fig:cellular-posets-and-non-cellular-posets} we depict cellular posets, as well as posets which fail to be cellular. Note, in particular, even if upper closures $P^{\geq x}$ realizes to topological balls it need not be the case that strict upper closures $P^{>x}$ realize to spheres.
\begin{figure}[h!]
    \centering
    \def\svgwidth{1\columnwidth}
    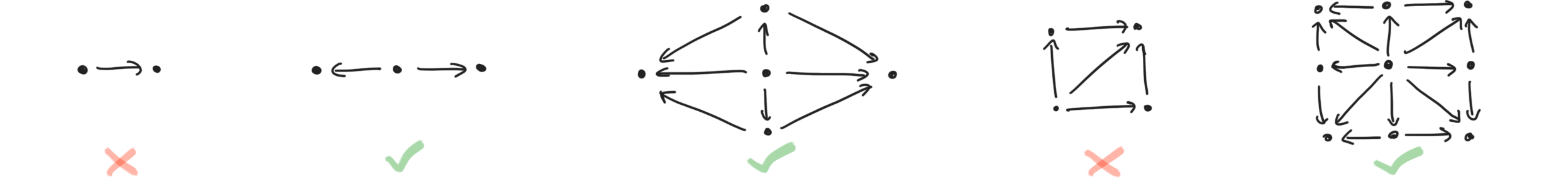

    \caption{Cellular and non-cellular posets.}
    \label{fig:cellular-posets-and-non-cellular-posets}
\end{figure}
\end{eg}

\nid Centrally, classifying stratifications of cellular posets are exactly regular cell complexes as recorded in the following result.

\begin{prop}[Regular cell complexes are classifying stratifications of cellular posets] \label{prop:cellular-posets-vs-reg-cell-cplx} Regular cell complexes are exactly classifying stratifications of cellular posets, in the following sense.
    \begin{enumerate}
        \item The classifying stratification of a cellular poset is a regular cell complex.
        \item The entrance path poset of a regular cell complex is a cellular poset.
        \item Every regular cell complex $X$ is stratified homeomorphic the classifying stratification of its entrance path poset, that is, $X \iso \CStr \Entr X$.\footnote{The isomorphism is `canonical up to homotopy'.}
        \item Every cellular poset $X$ is canonically homeomorphic to the entrance path poset of its classifying stratification, that is, $X \iso \Entr \CStr X$.
    \end{enumerate}
\end{prop}

\begin{proof} Statement (2) and (4) follow from the definitions. Statement (1) and (3) are discussed in \cite[\S 3]{bjorner1984posets}.
\end{proof}

\nid The correspondence of regular cell complexes (up to stratified homeomorphism) and cellular posets provides the \emph{combinatorialization} of regular cell complexes by cellular posets, as promised earlier.

\subsubsecunnum{The category of combinatorial regular cell complexes}

    Let us discuss how to extend the combinatorialization to include maps between regular cell complexes resp.\ between cellular poset. This well lead us to introduce a class of so-called `cellular' maps.

\begin{term}[Closure preservation for stratifications] A stratified map is said to be `closure preserving' if it maps closures of strata onto closures of strata.
\end{term}

\begin{defn}[Cellular maps of regular cell complexes] \label{defn:cellular-maps-of-rcc} A \textbf{cellular map of regular cell complexes} $F : X \to Y$ is a map of regular cell complexes that is closure preserving; that is, for each open cell $x$ in $X$ we have $F(\overline x) = \overline y$, where $y = \Entr F (x)$ is the open cell into which $x$ is mapped.
\end{defn}

\begin{notn}[The category of regular cell complexes and cellular maps] Denote by $\CStrat$ the category of regular cell complexes and their cellular maps.
\end{notn}


Let us next discuss cellular maps for cellular posets.

\begin{term}[Closure preservation for posets] \label{term:upper-closure-preservation} A map of posets $F : P \to Q$ is `upper-closure preserving' if for each $x \in P$ there exists $y \in Q$ such that the image $F P^{\geq x}$ equals $Q^{\geq y}$.
\end{term}

\begin{defn}[Cellular maps of cellular posets] \label{defn:cellular-maps-of-cellular-posets} A \textbf{cellular map of cellular posets} $P$ and $Q$ is a map of posets that is upper-closure preserving.
\end{defn}

\begin{notn}[The category of cellular posets and cellular maps] Denote by $\CPos$ the category of cellular posets and their cellular maps.
\end{notn}


\begin{eg}[Cellular and non-cellular maps] \label{eg:cell-complex-projections} In \autoref{fig:cellular-and-non-cellular-maps-of-cellular-posets} we depict cellular and non-cellular maps of regular cell complexes as well of their entrance path posets, which are cellular posets. In each case we indicate the mapping by coloring images and preimages in the same color.
\begin{figure}[ht]
    \centering
    \def\svgwidth{1\columnwidth}
    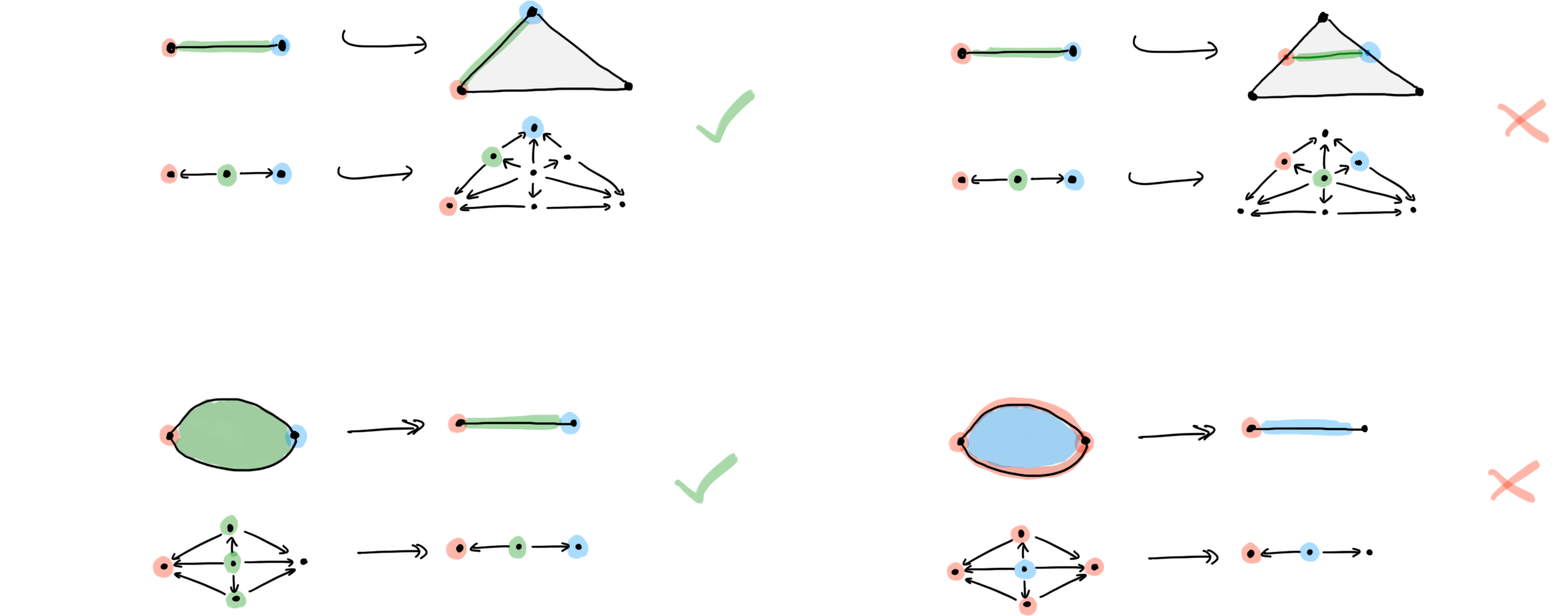

    \caption{Cellular and non-cellular maps.}
    \label{fig:cellular-and-non-cellular-maps-of-cellular-posets}
\end{figure}
\end{eg}

\nid The two definitions of cellular maps can be related by the following functors.

\begin{obs}[Entrance path and classifying functors] The entrance path poset construction (see \autoref{rmk:entr-face-functorial}), as well as the classifying stratification construction (see \autoref{term:classifying-stratifications} and \autoref{term:classifying-stratifications-maps}), yield functors
\begin{equation}
    \begin{tikzcd}[column sep=55pt]
        \CStrat \arrow[r, "\Entr{}", bend left=9]  & \CPos \arrow[l, "\CStr{}", bend left=9]
    \end{tikzcd} \quad .
\end{equation}
\end{obs}

\nid The central observation about cellular maps is the following.

\begin{rmk}[The equivalence of regular cell complexes and cellular posets] \label{rmk:cellular-map-yield-eqv} For sufficiently nice spaces, the entrance path poset construction in fact yields a functor of \emph{topological} categories $\Entr : \TStrat \to \TPos$ (see \autoref{rmk:entrz_cont}). Denote the topological subcategory of $\TStrat$ given by regular cell complexes and their cellular maps by $\TCStrat$. Similarly, one can consider the subcategory of $\TPos$ given by cellular posets and cellular maps: due to the cellularity condition on maps, this subcategory in fact has discrete hom spaces, recovering the `ordinary' category $\CPos$. The entrance path poset functor now restricts to a functor $\Entr : \TCStrat \to \CPos$. The functor provides a `weak equivalence' of topological categories (the `weak inverse' of $\Entr$ is given by the classifying stratification functor $\CStr{}$). The category $\TCStrat$ of regular cell complexes and their cellular maps is, in this sense, `combinatorializable'.
\end{rmk}

\nid While we not prove the preceding remark (nor will we use it in any form), it motivates that $\CPos$ is in fact a good model for the `combinatorial' category of regular cell complexes. We record this as a definition.

\begin{defn}[The combinatorial category of regular cell complexes] The \textbf{`combinatorial' category of regular cell complexes} $\CCplx$ is the category of cellular posets and their cellular maps.
\end{defn}

\nid Note that the adjective `combinatorial' is simply meant to highlight the combinatorial nature of the category's objects. We often refer to objects in $\CCplx$ as `combinatorial regular cell complexes', which emphasizes that objects can be interpreted both in combinatorial terms (as cellular posets) and a stratified-topological terms (as regular cell complexes obtained via the classifying stratification construction).

We will be particularly interested in `combinatorial regular cells' as recorded by the following terminology. Recall, the `depth' of an object $x$ in a poset measures the length $k$ of maximal chains $x_0 \to  x_1 \to ... \to x_k$ in the poset starting at that object, $x = x_0$.

\begin{term}[Regular cells] \label{term:regular-cells} A cellular poset $X \in \CCplx$ is called a `combinatorial regular $k$-cell' (or simply, a `regular $k$-cell' if no confusion arises) if $X$ has an initial object of depth $k$. We often denote initial elements of cells by $\ino$.
\end{term}

\nid Equivalently, in stratified topological terms, a regular cell complex $X$ is a regular $k$-cell if it is the closure of a single $k$-dimensional cell in the complex.

\pause

We next show that the `combinatorial' category of simplicial complexes $\SCplx$ fully faithfully embeds in the `combinatorial' category of regular cell complexes $\CCplx$. 

\begin{obs}[Simplicial complexes embed in regular cell complexes] \label{obs:entr-of-simp-cplx} Given a simplicial complex $K$, its `entrance path poset' $\Entr K$ is the poset whose objects of $\Entr K$ are simplices $x$ in $K$ with an arrow $x \to y$ whenever the simplex $y$ is a face of the simplex $x$. For a simplicial map of unordered simplicial complexes $F : K \to L$, we obtain a map of entrance path posets $\Entr F : \Entr K \to \Entr L$ mapping a simplex $x$ in $K$ to the simplex $Fx$ in $L$. This yields the functor
    \begin{equation}
        \Entr : \SCplx \to \CCplx
    \end{equation}
which, in fact, is a fully faithful embedding of categories. Note if we were to allow non-cellular poset maps in $\CCplx$ the latter claim would fail to hold.
\end{obs}

\pause

To end this section, let us briefly address the discrepancy between `topology' and `PL topology' which is, in fact, also visible at the level of `combinatorial regular cell complexes'. For this, we introduce the following PL analog of the definition of cellular posets.

\begin{defn}[PL cellular posets] \label{defn:PL-cellular-poset}
A poset $(X,\leq)$ is called \textbf{PL cellular} if the realization $\abs{X^{>x}}$ of the strict upper closure of any $x \in P$ is PL homeomorphic to a PL sphere.
\end{defn}

\begin{rmk}[Cellular is not always PL cellular] \label{rmk:cellular-not-PL-cellular} Note that, while `PL cellular' always implies `cellular' the converse is in general not true. Indeed, there exist triangulations of the sphere which are not PL spheres (see \cite{edwards1980topology}, \cite[Thm. 9.1]{bryant2002piecewise}). Adjoining a new minimal element to the entrance path poset of such a triangulation yields a poset that is cellular but not PL cellular. In contrast, we will later on find that in the framed setting the adjectives `cellular' and `PL cellular' can be used interchangeably (see \autoref{punch:cellular-vs-PL-cellular}).
\end{rmk}

\begin{term}[The combinatorial category of regular cell complexes] The `(combinatorial) category of regular PL cell complexes' $\CCplxPL$ is the category of cellular posets and their cellular maps.
\end{term}

\begin{rmk}[The discrepancy of PL and TOP and its resolution] \label{rmk:discrepancy-of-top-and-PL-cells} Note that $\CCplxPL \subsetneq \CCplx$. While our notation therefore needs to distinguish the PL from the topological case here, as we will see the difference disappears in the framed case. In other words, we will find that `framed' analogs $\FrTCCplx n$ and $\FrCCplx n$ of the above categories are in fact the same category, and we shall therefore not distinguish them notationally (see \autoref{punch:cellular-vs-PL-cellular}).
\end{rmk}


\subsection{Framings and proframings on regular cells} \label{ssec:framed-reg-cells}

In this section we define framings and proframings on regular cells. This will combine our discussion of framings on simplices and simplicial complex, with the notion of `combinatorial regular cell complexes', and their category $\CCplx$, as discussed in detail in the previous section. The following (abuse of) notation will reduce the amount of symbols needed in this and subsequent sections.

\begin{notn}[Structures related to cellular posets] \label{rmk:interpretations-of-combinatorial-regular-cell-complexes} Consider a cellular poset $X \in \CCplx$ (i.e.\ a `combinatorial regular cell complex').
    \begin{enumerate}[wide, labelwidth=!, labelindent=0pt]
        \item The regular cell complex $\CStr X$ classifying $X$ will, abusing notation, usually be referred to simply as the `regular cell complex $X$'.
        \item The simplicial complex $\Unord {NX}$ obtained by unordering the nerve of $X$ will, abusing notation, be referred to as the `(underlying) simplicial complex $X$'.
    \end{enumerate}
    The abuse of notation similarly applies to maps of cellular posets $F : X \to Y$ (and thus to subposets $X \into Y$), which context-dependently may be used to denote maps of corresponding regular cell complexes and of underlying simplicial complexes.
\end{notn}

\begin{figure}[ht]
    \centering
    \def\svgwidth{1\columnwidth}
    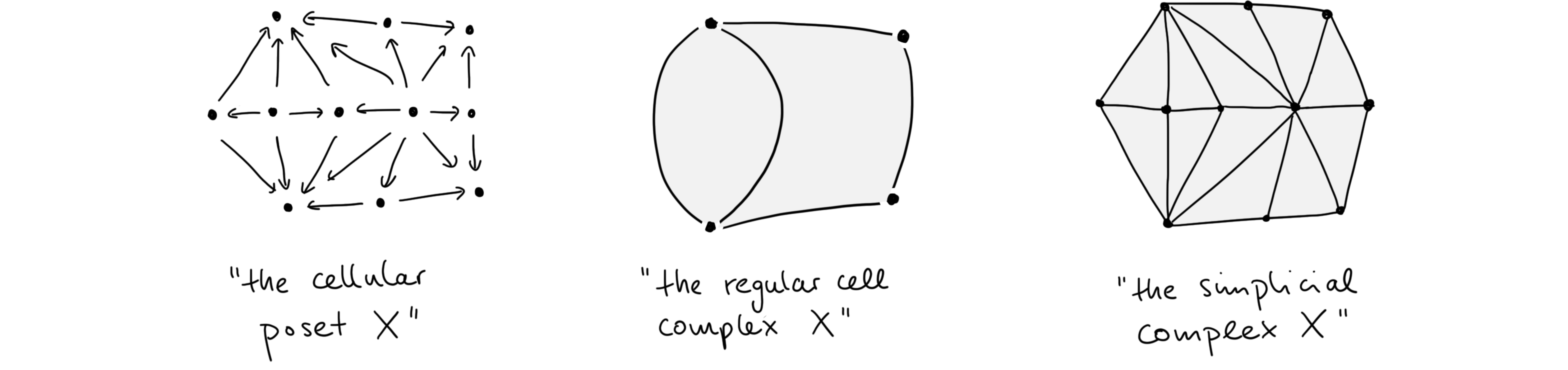

    \caption[Structures related to cellular posets.]{A cellular poset, together with its corresponding regular cell complex and underlying simplicial complex.}
    \label{fig:overloading-meaning-of-objects-in-ccplx}
\end{figure}

\nid In \autoref{fig:overloading-meaning-of-objects-in-ccplx} we illustrate a cellular poset $X$, together with its corresponding regular cell complex $X$ and its underlying simplicial complex $X$.

\subsubsecunnum{The definition of framed cells}

We now define framings on regular cells. Recall, a regular cell is a cellular poset with an initial element (see \autoref{term:regular-cells}).

\begin{defn}[Framed regular cells] \label{defn:framed-reg-cells} An \textbf{$n$-framing $\cF$ of a regular cell $X$} is an $n$-framing $\cF$ of the simplicial complex $X$ such that, for each $x \in X$, the framing restricts on the upper closure $X^{\geq x}$ of $x$ to a flat framing $\rest \cF {X^{\geq x}}$ of the simplicial subcomplex $X^{\geq x} \into X$.
\end{defn}

\nid We will refer to the pair $(X,\cF)$, of a regular cell $X$ together with an $n$-framing $\cF$ on it, as an `$n$-framed regular cell'.

\begin{notn}[The ordered simplicial complex $X$] \label{rmk:overloading-notation} Given an $n$-framed regular cell $(X,\cF)$, note that the $n$-framing $\cF$ endows the underlying simplicial complex $X$ with an ordering: this yields, abusing notation, the `ordered simplicial complex $X$' (which generally differs from the `nerve of $X$').
\end{notn}

\begin{notn}[Restricting framings to subcells] \label{notn:frames-on-subcells} Given an $n$-framed regular cell $(X,\cF)$ and $x \in X$, we abbreviate the restricted $n$-framing $\rest \cF {X^{\geq x}}$ by $\rest \cF x$.
\end{notn}

\begin{eg}[2-Framed regular cells] In \autoref{fig:framed-regular-cells-basic-examples} we illustrate several examples of 2-framed regular $k$-cell $(X,\cF)$, for $k = 1$ and $2$: in each case we depict three pieces of data, namely, the regular cell complex $X$ in the top row, the underlying simplicial complex $X$ in the second row, and the ordered simplicial complex $X$ together with the framing $\cF$ (indicated by frame labels) in the third row. In the fourth row, we exploit the flatness of $\cF$ to provide an alternative illustration of the third row by framed realizations in $\lR^2$ (see \autoref{rmk:flat-framing-via-realization}).
\begin{figure}[ht]
    \centering
    \def\svgwidth{1\columnwidth}
    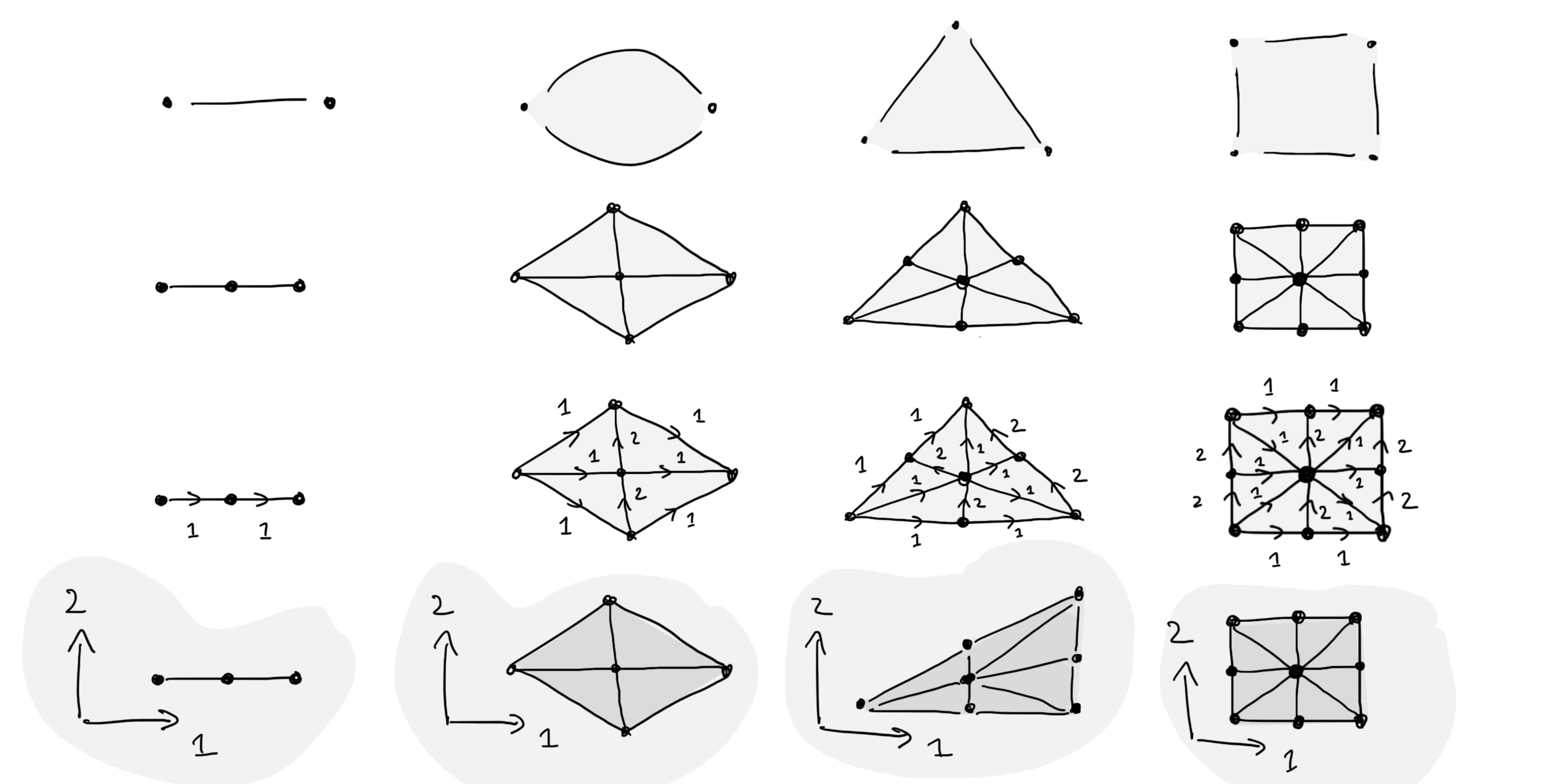

    \caption{Data of 2-framed regular cells.}
    \label{fig:framed-regular-cells-basic-examples}
\end{figure}
\end{eg}

\nid Usually, depictions of regular cell complexes $X$, of their underlying simplicial complexes, and of $n$-framings $\cF$ on them, will be condensed into a single picture as the next remark explains.

\begin{rmk}[Condensing cell structure, simplicial structure and framings] \label{rmk:condensing-framed-cell-depiction} We usually depict $n$-framed regular cells $(X,\cF)$ by embedding the regular cell $X$ in $\lR^n$: this, by passing to the underlying simplicial complex $X$ and then deriving an $n$-framing from the resulting framed realization in $\lR^n$, determines all data of $(X,\cF)$. We illustrate this with further examples in \autoref{fig:framed-regular-cells-further-examples}: note that the first row  re-illustrates exactly the four examples given in \autoref{fig:framed-regular-cells-basic-examples}.
\begin{figure}[ht]
    \centering
    \def\svgwidth{1\columnwidth}
    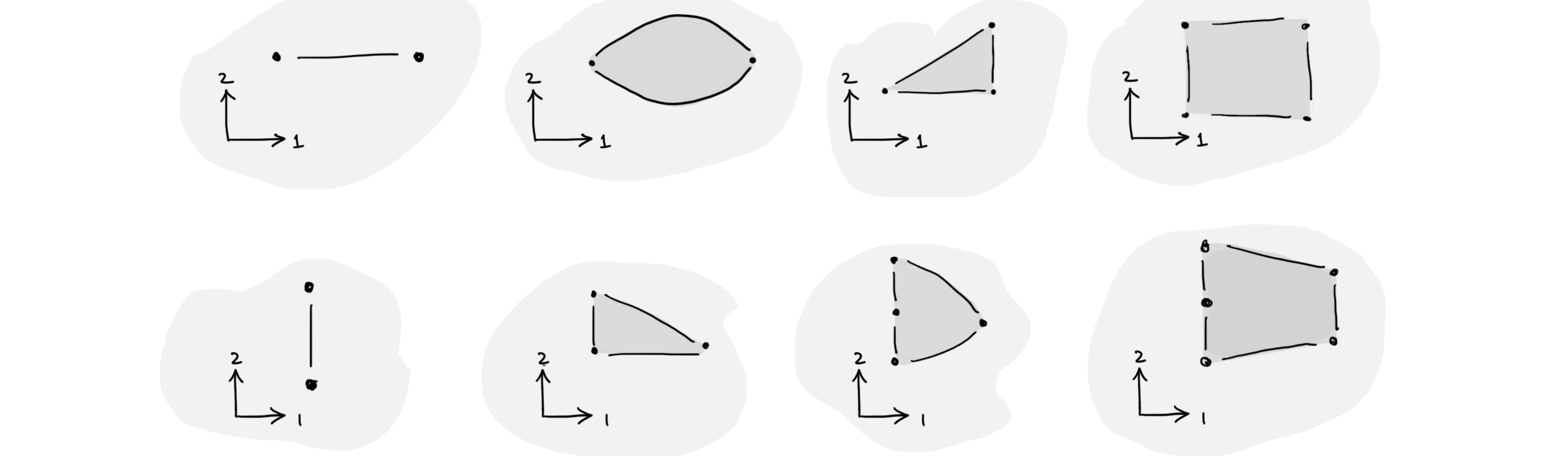

    \caption{2-Framed regular cells depicted via framed realizations.}
    \label{fig:framed-regular-cells-further-examples}
\end{figure}
\end{rmk}

\begin{eg}[Non-examples of framed regular cells] In \autoref{fig:framed-regular-cell-non-examples} we illustrate three non-examples of $n$-framings on regular $k$-cells $X$: in each case, we separately depict the regular cell $X$ (as a cell complex) on top, together with a framing $\cF$ on $X$ (as a simplicial complex) underneath it. In the leftmost example, the chosen framing fails to be flat. In the second example, the framing is flat, but fails to be flat when restricted to either of the two 1-cells. Conversely, the third framing is flat when restricted to any subcell in the boundary but not on the 2-cell itself.
\begin{figure}[ht]
    \centering
    \def\svgwidth{1\columnwidth}
    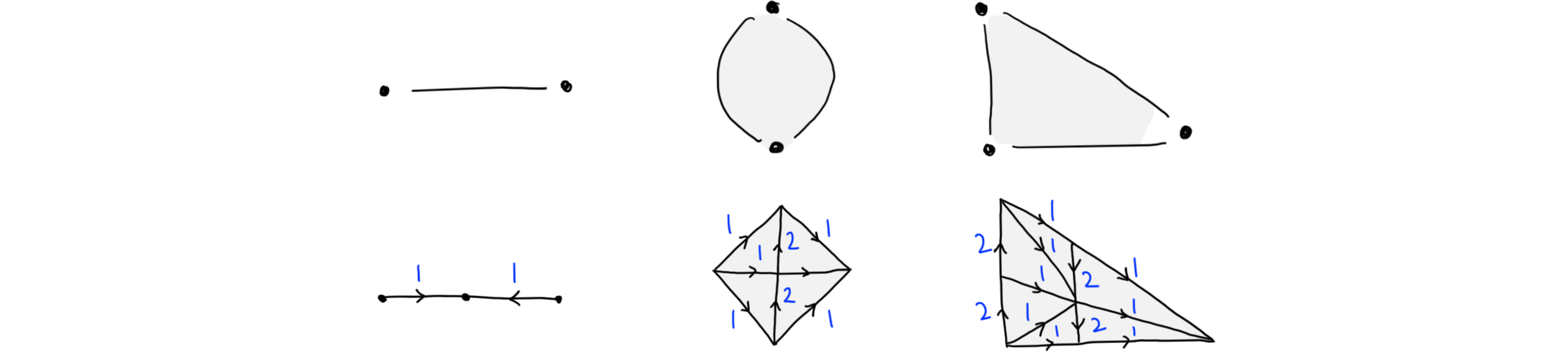

    \caption{Non-examples of framed regular cells.}
    \label{fig:framed-regular-cell-non-examples}
\end{figure}
\end{eg}

\begin{eg}[The simplest 3-framed regular 3-cell] In \autoref{fig:different-framings-of-the-hemspherical-2-sphere-complex} on the left, we depict a framed regular $3$-cell $(X,\cF)$. For clarity, we also separately depict the $3$-framing $\cF$ on the underlying simplicial complex $X$ on the right.
\begin{figure}[h!]
    \centering
    \def\svgwidth{1\columnwidth}
    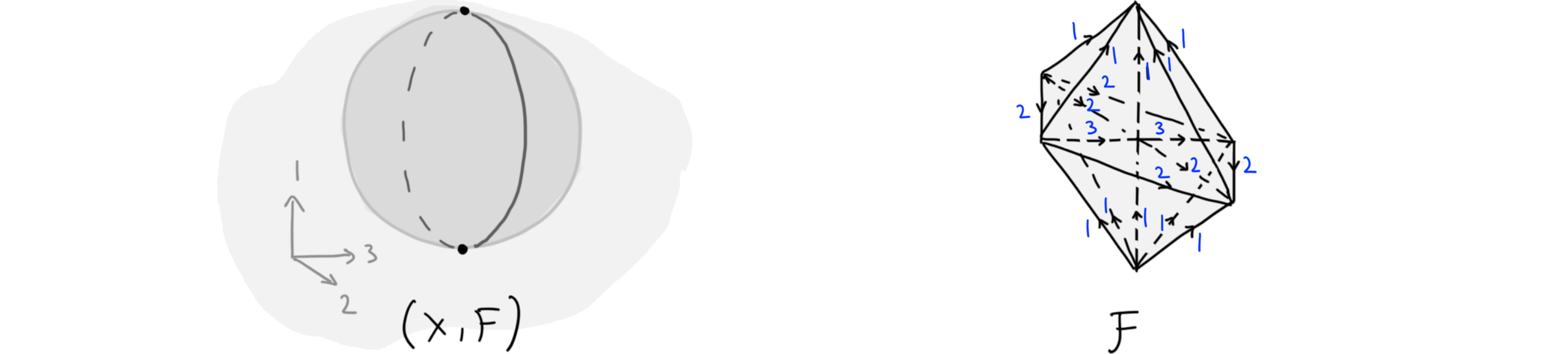

    \caption{The simplest 3-framed regular 3-cell.}
    \label{fig:different-framings-of-the-hemspherical-2-sphere-complex}
\end{figure}
\end{eg}

\nid The previous example is the simplest example of a 3-framed regular 3-cell (indeed, its underlying regular cell complex is the simplest regular cell complex of the 3-ball). In general, 3-framed regular 3-cells can be of various shapes as the next example illustrates.

\begin{eg}[More 3-framed regular 3-cells] \label{eg:more-3-blocks} We depict few more 3-framed regular 3-cells in \autoref{fig:3-framed-regular-3-cells-examples}. Note that cells in the boundary of these 3-cells yields various 3-framed k-cells for $k < 3$.
\begin{figure}[ht]
    \centering
    \def\svgwidth{1\columnwidth}
    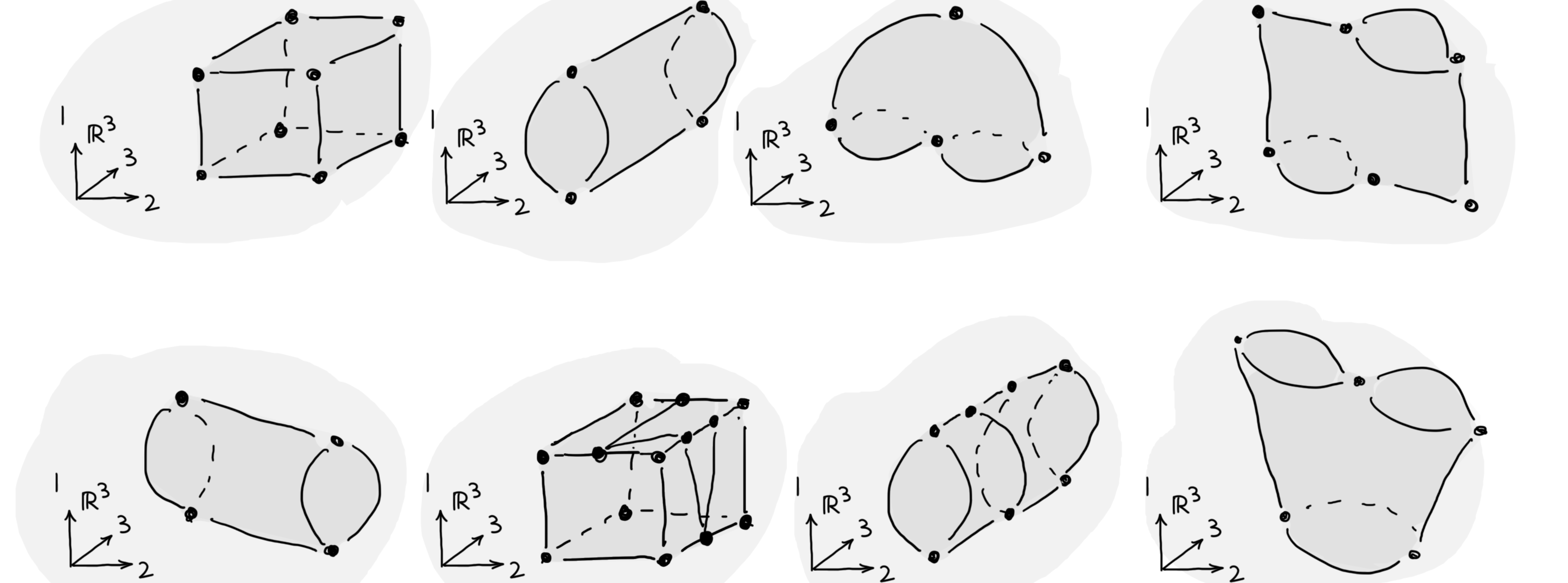

    \caption{3-framed regular 3-cells.}
    \label{fig:3-framed-regular-3-cells-examples}
\end{figure}
\end{eg}

\pauseae

Having introduced the notion of framed regular cells, let us briefly revisit two central punchlines. Firstly, recall our earlier \autoref{rmk:discrepancy-of-top-and-PL-cells} about the `discrepancy of topological cells and PL cells'. In the framed setting the discrepancy disappears.

\begin{rmk}[Unity of framed cells and PL cells] \label{punch:cellular-vs-PL-cellular} Topological and PL cells coincide in the framed setting, that is, given a framed regular cell $(X,\cF)$ then the framing forces $X$ to not only be a regular cell, but in fact a regular \emph{PL} cell: this will be shown in \autoref{lem:block-shellability}. There is consequently no need to consider notions of `framed cells' and `framed regular PL cell' separately, as both these notions coincide.
\end{rmk}

\nid Secondly, framed regular cells can (unlike nonframed regular cells) be classified.

\begin{rmk}[Tractability of framed regular cells] \label{punch:block-classification} We can recognize framed regular cell among framed posets (here, a `framed poset' is simply a poset together with a framing of its underlying simplicial complex). We will show this by providing a `constructive' classification of framed regular cells in \autoref{ch:classification-of-framed-cells}.
\end{rmk}

\nid In this way, framed regular cells provide a novel class of general shapes that is computationally tractable.

\begin{rmk}[Framed regular cells with degenerate faces] A yet more general, and still tractable class of shapes can be obtained by allowing faces of framed regular cells to `degenerate'. While these shapes naturally appear when considering general presheafs on framed regular cells, we will not consider them further here.
\end{rmk}

\subsubsecunnum{Framed cellular maps}

We next define framed maps of framed regular cells. Recall, in the case of framed simplices, we defined framed maps as maps that preserve frames on each `vector'; the setup crucially relied on vectors in $m$-simplices being generated by their $m$ spine vectors (which provided a `basis' for the affine space of the simplex). Unfortunately, there is no good analog for `spine vectors' in the context of regular cells in general; however, for \emph{framed} regular cells, we may recover a notion of `final frame vectors' for any given framed regular cell as explained below. This turns out to be just enough to define framed maps of framed regular cell.

\begin{term}[Final frame vectors of framed regular cells]  Given an $n$-framed regular $k$-cell $(X,\cF)$ with initial element $\ino \in X$, a `final frame vector' $v$ is a 1-simplex in the ordered simplicial complex $X$ containing $\ino$ as a vertex and whose frame label $\cF_v \in \bnum n$ is maximal among all frame labels of 1-simplices in $X$.
\end{term}

\begin{notn}[Final frame vectors] Given an $n$-framed regular $k$-cell $(X,\cF)$ and an object $x \in X$, we denote by $\svec (x) \subset X$ the subcomplex of the ordered simplicial complex $X$ spanned by the final frame vectors of the cell $(X^{\geq x}, \rest \cF x)$.
\end{notn}

\begin{eg}[Final frame vectors of framed regular cells] In \autoref{fig:final-vectors-of-framed-cells}, for several framed regular cells $(X,\cF)$ (which reproduce selected cells from both \autoref{fig:framed-regular-cells-further-examples} and \autoref{fig:3-framed-regular-3-cells-examples}), we depict the final frame vectors by red edges; in the first row, we also depict underneath each cell its corresponding ordered simplicial complex (but leave these complexes implicit in the second row). Note that, in each case, there are in fact exactly two final frame vectors: an edge that ends in the vertex $\bot$, and an edge that starts at the vertex $\bot$. We further indicate, by green and blue edges, the final frame vectors of a selection of subcells in the given cells boundary.
\begin{figure}[ht]
    \centering
    \def\svgwidth{1\columnwidth}
    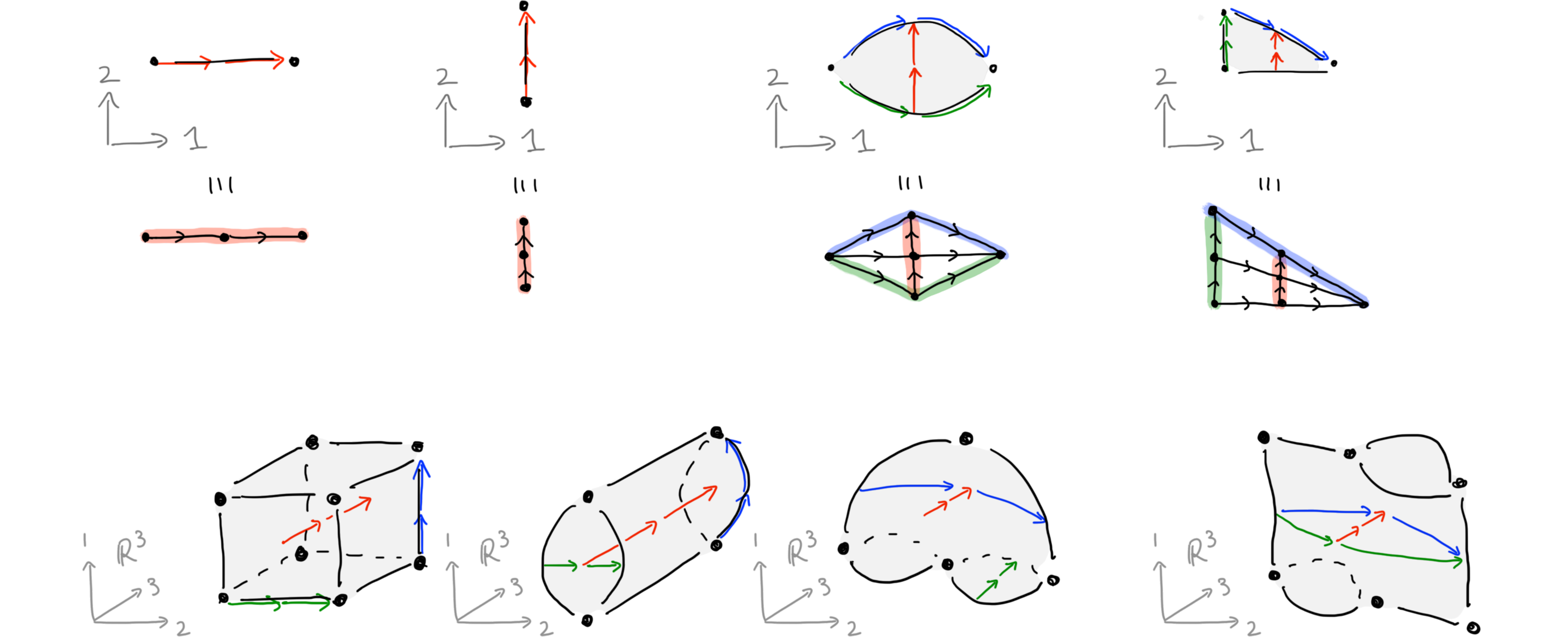

    \caption{Final frame vectors of framed regular cells.}
    \label{fig:final-vectors-of-framed-cells}
\end{figure}
\end{eg}

\begin{rmk}[Final frame vectors `form a vector'] \label{obs:final-frame-vec-is-vec} The observation in the preceding example generalizes: any framed cell $(X,\cF)$ has exactly two final frame vectors (if $X$ is non-trivial, that is, $\dim(X) > 0$). The subcomplex $\svec (\ino)$ of the ordered simplicial complex $X$ is isomorphic to the linear simplicial complex (\begin{tikzpicture}
\begin{scope}[decoration={markings,mark=at position 0.5 with {\arrow{>}}}]
    \draw[postaction={decorate}] (-.7,0)--(0,0);
    \draw[postaction={decorate}] (0,0)--(.7,0);
    \foreach \x/\name in {-.7,0,.7}{
        \path[draw,fill=black] (\x,0) circle[radius=1.8pt];
    }
\end{scope}
\end{tikzpicture}), while the subposet $\svec (\ino)$ of the poset $X$ (determined by the subcomplex $\svec (\ino) \subset X$) is isomorphic to the entrance path poset $(\bullet \ot \bullet \to \bullet)$ of a 1-simplex. In this sense the final frame vectors $\svec (\ino)$ of a framed regular cell may be regarded as forming a `single directed 1-simplex'. We will revisit and prove this claim in \autoref{cor:final-frame-vec-is-vec}.
\end{rmk}

\nid The definition of framed maps of framed regular cells now mirrors the definition of framed maps of framed simplices (see \autoref{defn:framed-map-of-framed-simplices}) as follows.

\begin{defn}[Framed maps of framed regular cells] Given $n$-framed regular cells $(X,\cF)$ and $(Y,\cG)$, a \textbf{framed cellular map} $F : (X,\cF) \to (Y,\cG)$ is a cellular map of cellular posets $F : X \to Y$, such that for all $x \in X$, either $F$ preserves the framing of final frame vectors $\svec x$, that is, $F$ restricts to a framed simplicial isomorphism $(\svec x, \rest \cF {\svec x}) \iso (F(\svec x), \rest \cG {F (\svec x)})$, or the final frame vectors are degenerated by $F$, i.e.\ $F(\svec x) \subset Y$ is a point.
\end{defn}

\begin{eg}[Framed cellular maps] In the upper row of \autoref{fig:framed-cellular-map-examples} we illustrate examples of framed maps of 2-framed regular cells; in each case we highlight image and preimage cells in the same color. Underneath each map we also depict the corresponding maps of simplicial complex as well as the orderings of those complexes. In the lower row of \autoref{fig:framed-cellular-map-examples} we similarly depict non-examples of framed maps. While all three depicted maps are cellular they fail preserve frame vectors in the required sense; the first fails to preserve final frame vectors of one of the two blue edges; the second fails to preserve final frame vectors of the blue 2-cells; the last fails to preserve the final frame vectors of the left red edge.\begin{figure}[ht]
    \centering
    \def\svgwidth{1\columnwidth}
    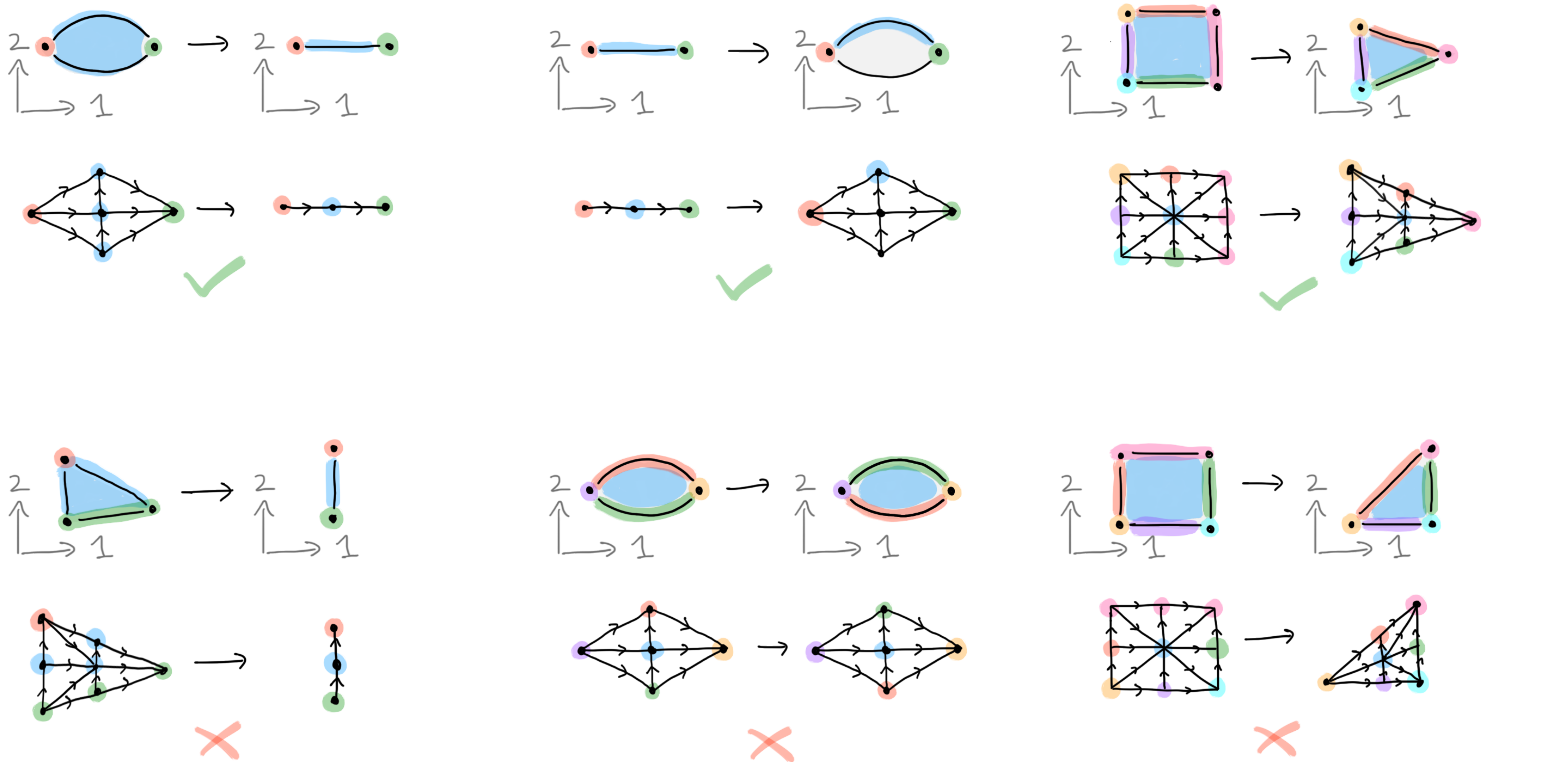

    \caption{Examples and non-examples of framed cellular maps.}
    \label{fig:framed-cellular-map-examples}
\end{figure}
\end{eg}

\begin{rmk}[Framed maps of cell are subframed on simplices] Note that in framed maps $F : (X,\cF) \to (Y,\cG)$ of framed regular cells, need not descend to framed maps of framed simplicial complexes $F : (X,\cF) \to (Y,\cG)$, but they do descend to \emph{subframed} maps of framed simplicial complexes $F : (X,\cF) \to (Y,\cG)$ (see \autoref{rmk:subframed-maps-of-framed-simp-cplx}); that is, $F$ need not preserve the simplicial ordering and may `specialize' frame labels of vectors. This is illustrated in \autoref{fig:a-framed-cellular-map-need-not-induce-a-framed-simplicial-map}.
\begin{figure}[ht]
    \centering
    \def\svgwidth{1\columnwidth}
    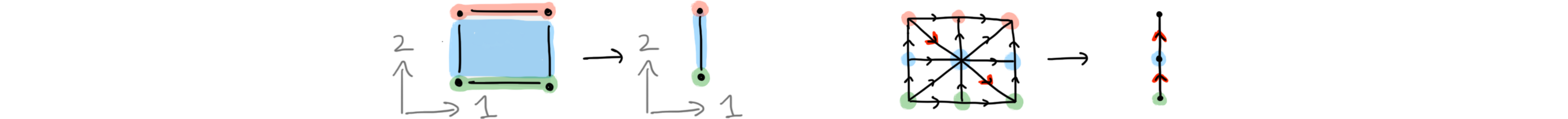

    \caption[Framed cellular maps are generally not framed simplicial.]{A framed cellular map need not induce a framed simplicial map, but does induce a subframed simplicial map.}
    \label{fig:a-framed-cellular-map-need-not-induce-a-framed-simplicial-map}
\end{figure}
\end{rmk}

\begin{notn}[The category of framed regular cells] \label{notn:frccell}  The category of $n$-framed regular cells and their framed cellular maps will be denoted by $\FrCCell n$.
\end{notn}

\subsubsecunnum{The definition of proframed cells} \label{sssec:proframings-on-cells}

Just has simplicial complexes may be endowed both with framings or proframings, we now introduce proframings of regular cells. Recall that proframed simplicial complexes were defined as certain sequences of surjective ordered simplicial maps (see \autoref{defn:proframings-of-simp-cplx}). The role of surjective simplicial maps will now be played by surjective cellular maps. The following terminology will further be useful.

\begin{term}[Ordering a sequence of simplicial maps] \label{rmk:ordering-pre-cell-proframe} Recall from \autoref{term:orderings}, that an `ordering' of a simplicial map $F : K \to L$ is an ordered simplicial map which, after unordering, recovers $F$. Similarly, given a sequence of simplicial maps, an `ordering' of it is a sequence of ordered simplicial maps which, after unordering, recovers the given sequence.
\end{term}

\begin{defn}[Proframed regular cells] \label{defn:proj-framed-reg-cells} An \textbf{$n$-proframing $\cP$ of a regular cell $X$} is a sequence of surjective cellular maps $X = X_n \xto {p_n} X_{n-1} \xto {p_{n-1}} ... \xto {p_1} X_0$, whose sequence of underlying simplicial maps is endowed with an ordering, such that the resulting sequence of ordered simplicial maps (denoted again by $\cP$) is an $n$-proframing of the simplicial complex $X$, and with the condition that, for each $x \in X$, the restriction $\rest \cP {X^{\geq x}}$ of $\cP$ to the subcomplex $X^{\geq x}$ is flat.
\end{defn}

\nid We will refer to the pair $(X,\cP)$, of a regular cell $X$ together with an $n$-proframing $\cP$ on it, as an `$n$-proframed regular cell'. Note that `$\cP$' refers, similarly to previous notation, to multiple structures simultaneously.

\begin{rmk}[Notation for data of proframes] \label{rmk:overloading-notation-proframes}  Abusing notation, we may refer to $\cP$ as a `sequence of cellular maps', a `sequence of simplicial maps', or a `sequence \emph{ordered} simplicial maps' all of which refer to different parts of the data of the proframed cell $(X,\cP)$ (see also \autoref{rmk:overloading-notation}).
\end{rmk}

\begin{notn}[Restricting framing to subcells] \label{notn:proframes-on-subcells} We abbreviate the cell restrictions $\rest \cP {X^{\geq x}}$ by $\rest \cP x$.
\end{notn}

\begin{eg}[2-proframed regular cells] \label{eg:proframed-cells} We illustrate four 2-proframed regular cells $(X,\cP)$ in \autoref{fig:2-proframed-regular-cell-examples}. In each case we depict the sequence of surjective cellular maps $X_2 \to X_1 \to X_0$ embedded in the standard proframe $\lR^2 \to \lR^1 \to \lR^0$. In analogy to our examples of framed regular cells, each cell embedding $X_i \into \lR^i$ represents an embedding $\abs{X_i} \into \lR^i$ of the realized underlying simplicial complex $X_i$, and requiring the latter embeddings to be the components of a proframed realization (see \autoref{rmk:flat-proframing-via-realization}) fully determines the proframing $\cP$.
\begin{figure}[ht]
    \centering
    \def\svgwidth{1\columnwidth}
    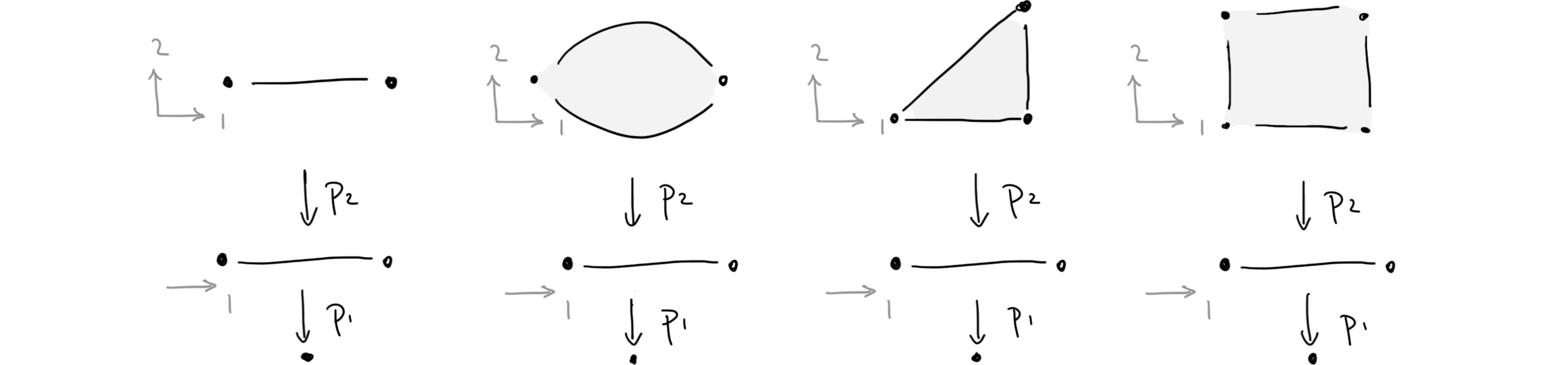

    \caption{2-Proframed regular cells.}
    \label{fig:2-proframed-regular-cell-examples}
\end{figure}
\end{eg}

\begin{eg}[3-proframed regular 3-cells] Recall the `simplest' 3-framed 3-cell from \autoref{fig:different-framings-of-the-hemspherical-2-sphere-complex}. The analogous 3-proframe is depicted in \autoref{fig:a-tower-of-surjective-cellular-poset-maps-realizes-to-a-tower-of-cell-complex-projections} (the `analogy' will be made precise in \autoref{obs:frames-vs-proj-frames-for-reg-cells}). We also illustrate four (less simple) 3-proframed regular 3-cells in \autoref{fig:more-3-blocks}.
\begin{figure}[h!]
    \centering
    \def\svgwidth{1\columnwidth}
    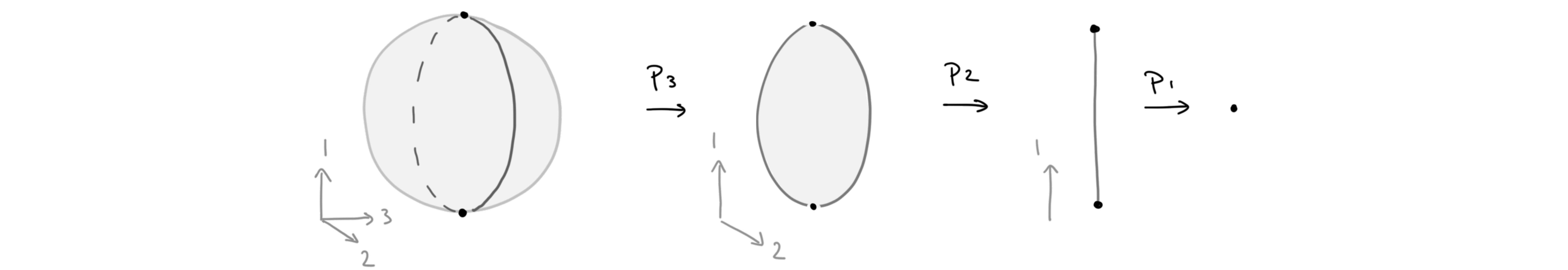

    \caption{The simplest 3-proframed regular 3-cell.}
    \label{fig:a-tower-of-surjective-cellular-poset-maps-realizes-to-a-tower-of-cell-complex-projections}
\end{figure}
\begin{figure}[h!]
    \centering
    \def\svgwidth{1\columnwidth}
    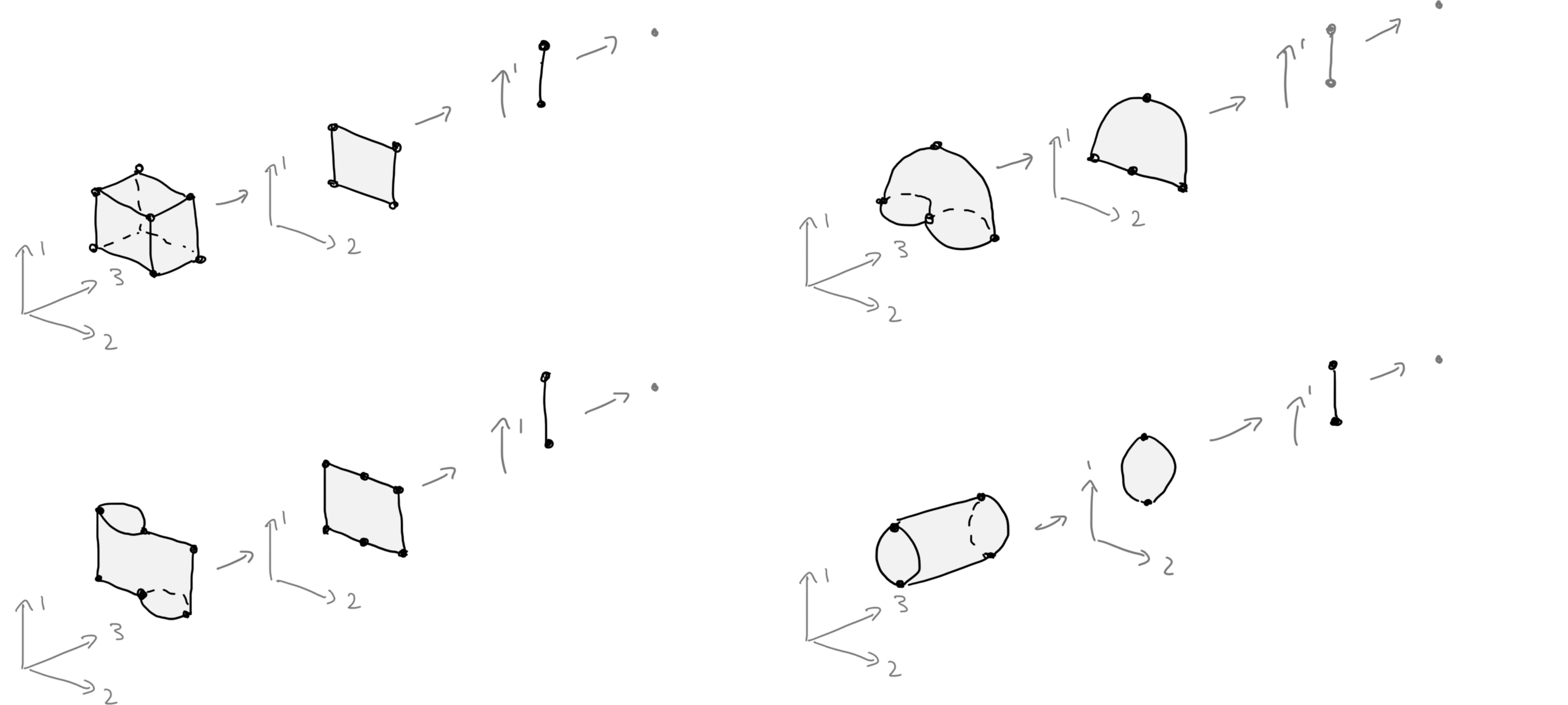

    \caption{3-proframed regular 3-cells.}
    \label{fig:more-3-blocks}
\end{figure}
\end{eg}

\subsubsecunnum{Proframed cellular maps}

\nid We define maps of proframed regular cells. Like in the case of framed maps of framed regular cells, we are confronted with the difficulty of regular cells not have `spine vectors'; as before, we will remedy this using the observation that framed regular cells have `final frame vectors'. The definition of proframed maps then takes the following form, which is parallel to the earlier definition of proframed maps of framed simplices (see \autoref{defn:proframed-maps}).

\begin{defn}[Proframed maps] \label{defn:proframed-cellular-map} Given $n$-proframed regular cells $(X, \cP = (p_n,...,p_1))$ and $(Y,\cQ = (q_n,...,q_1))$, a \textbf{proframed cellular map} $F : (X,\cP) \to (Y,\cQ)$ (or simply, a `proframed map') is a map of cellular poset sequences
\begin{equation}
    \begin{tikzcd}[baseline=(W.base)]
        {X = X_n} \arrow[d, "F_n"'] \arrow[r, "p_n"] & {X_{n-1}} \arrow[d, "F_{n-1}"'] \arrow[r, "p_{n-1}"] & {} \arrow[d, "\cdots", phantom] \arrow[r, "p_2"] & {X_1} \arrow[d, "F_1"] \arrow[r, "p_1"] & {X_0 = [0]} \arrow[d, "F_0"] \\
        {Y = Y_n} \arrow[r, "{q_n}"'] & {Y_{n-1}} \arrow[r, "{q_{n-1}}"'] & \cdots \arrow[r, "{q_2}"'] & {Y_1} \arrow[r, "{q_1}"'] & |[alias=W]| {Y_0 = [0]}
    \end{tikzcd}
\end{equation}
where each $F_i$ is a cellular map, and such that for every $x \in X$, either the proframing of the final frame vectors $\svec x$ is preserved, i.e.\ $F$ restricts to a proframed simplicial isomorphism $F : (\svec x, \rest \cP {\svec x}) \iso (F_n (\svec x), \rest \cQ {F_n (\svec x)})$, or the final frame vectors are degenerated, i.e.\ $F_n (\svec x) \subset Y_n$ is a point.
\end{defn}

\nid We will forego giving examples of proframed cellular maps here (earlier examples of framed maps of framed regular cells may be adapted appropriately).

\begin{notn}[The category of proframed cells] The category of $n$-proframed regular cells and their proframed cellular maps will be denoted by $\PFrCCell n$.
\end{notn}

\subsection{Framings and proframings on regular cell complexes} \label{ssec:framed-reg-cell-complex}

We now define framings and proframings of regular cell complexes, generalizing the notions of framed and proframed regular cells introduced in the previous section. In fact, this generalization only requires minor modifications of previous definitions, replacing `cells' by `cell complexes' and requiring `(pro)framed maps' to be `(pro)framed maps on each cell'.

\subsubsecunnum{The definition of framings on cell complexes}

Framings of regular cell complexes are framings of their simplicial complexes in which `all cells are flat'.

\begin{defn}[Framings on cell complexes] \label{defn:framed-reg-cell-cplx} Given a regular cell complex $X$, an \textbf{$n$-framing} $\cF$ of $P$ is an $n$-framing of the simplicial complex $X$ such that, for all $x \in X$, the framing restricts on the upper closure $X^{\geq x}$ to a flat framing $\rest \cF {X^{\geq x}}$ of the simplicial complex $X^{\geq x}$.
\end{defn}

\nid We will refer to the pair $(X,\cF)$, of a regular cell complex $X$ together with an $n$-framing $\cF$ on it, as an `$n$-framed regular cell complex'. As before, we will denote framings restricted to subcells $X^{\geq x}$ by $\rest \cF x$ (see \autoref{notn:frames-on-subcells}).

\begin{term}[Flat and locally flat framed regular cell complexes] An $n$-framed regular cell complex $(X,\cF)$ is called `flat' (resp.\ `locally flat') if the framing $\cF$ of the simplicial complex $X$ is flat (resp.\ `locally flat').
\end{term}

\begin{eg}[Non-flat framed regular cell complexes] \label{eg:framed-cell-complexes} Framed regular cell complexes can be thought of as `gluings' of framed regular cells: we illustrate this in two instances in \autoref{fig:block-complex-unit-counit-gluing}: in both cases we glue two 2-framed regular cells along their boundary is shown; the resulting frame structure of the cell complex is indicated by moving ambient frames `inside' their respective cells. Note that in the first case, the resulting cell complex realizes to a contractible space, while in the second case it realizes to the $2$-sphere. Neither of the framed regular cell complexes is complexes flat.
\begin{figure}[h!]
    \centering
    \def\svgwidth{1\columnwidth}
    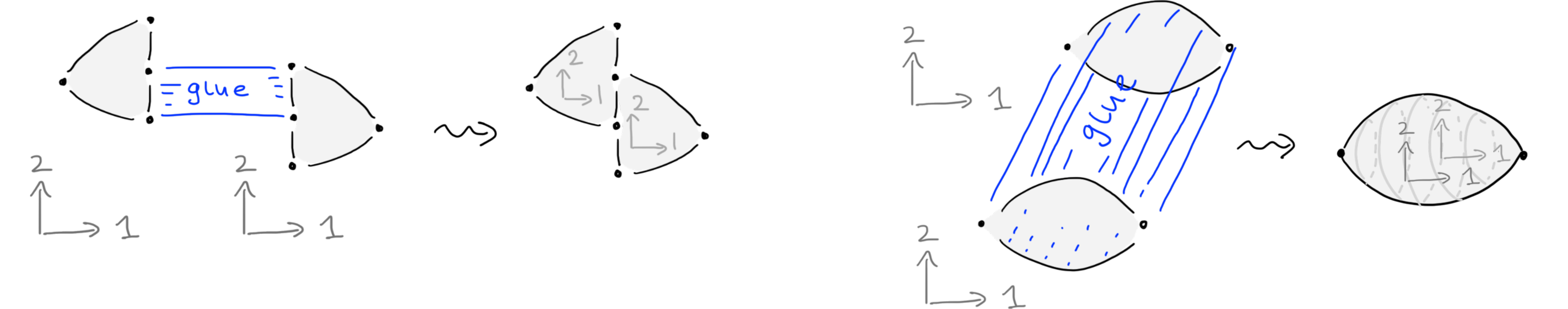

    \caption{Framed regular cell complexes obtained by gluing framed regular cells.}
    \label{fig:block-complex-unit-counit-gluing}
\end{figure}
\end{eg}

\begin{eg}[Locally flat cellulations of 1-manifolds] \label{eg:locally-flat-cellulations} In \autoref{fig:locally-flat-framed-cellulated-1-manifolds} on the left we depict several locally flat 1-framed regular cell complexes $(X,\cF)$ each cellulating a connected 1-manifold. On the right we depict two 1-framed regular cell complexes that are not locally flat framed.
\begin{figure}[ht]
    \centering
    \def\svgwidth{1\columnwidth}
    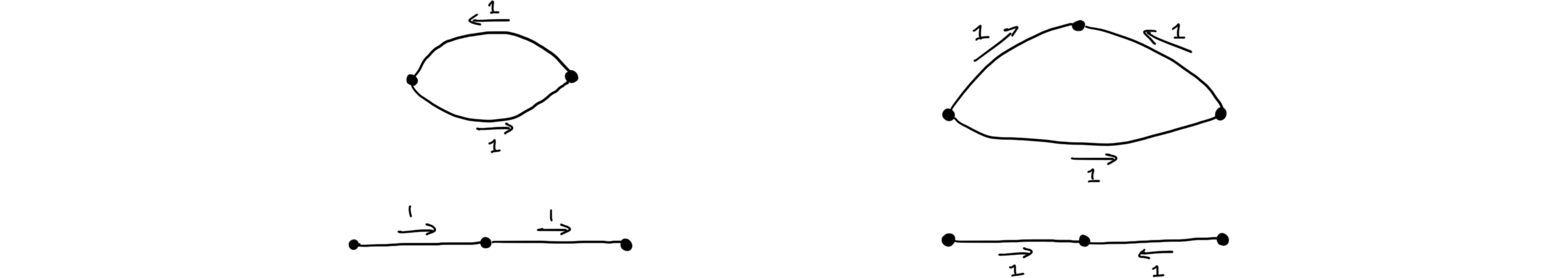

    \caption{Locally flat and non-locally flat 1-framed cellulated 1-manifolds.}
    \label{fig:locally-flat-framed-cellulated-1-manifolds}
\end{figure}
\end{eg}

\begin{eg}[Flat cellulation of the Hopf circle] \label{eg:generality-of-blocks} We motivate a slightly more complicated example. Consider the `Hopf circle' embedding of the circle $S^1$ into $\lR^3$, whose image projects along $\pi_3 : \lR^3 \to \lR^2$ to the figure eight, as shown in \autoref{fig:twisted-embedding-of-circle} below.
\begin{figure}[h!]
    \centering
    \def\svgwidth{1\columnwidth}
    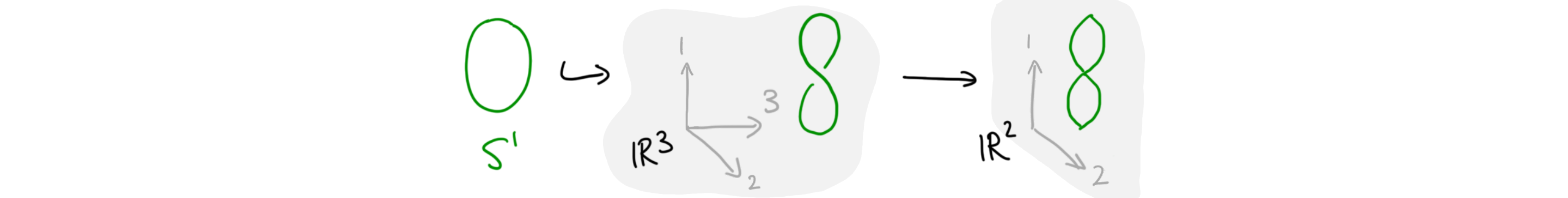

    \caption{The Hopf circle embedded in $\lR^3$.}
    \label{fig:twisted-embedding-of-circle}
\end{figure}
Take any compact 3-cube $I^3 \subset \lR^3$ containing the image of the Hopf circle embedding $S^1 \into \lR^3$. There exist a `canonical framed cellulation' of $I^3$ by a flat 3-framed regular cell complex containing the circle embedding in its 1-skeleton (here, the term `cellulation' is the cellular analog to the notion of `triangulation' by simplices). This complex is shown in \autoref{fig:triangulating-the-twisted-embedding-of-the-circle}: it consists of eight 3-framed 3-cells glued together as indicated (note all framings are determined by the same flat ambient 3-frame, which is indicated only once); we highlighted 0- and 1-cells cellulating the Hopf circle in blue.
\begin{figure}[h!]
    \centering
    \def\svgwidth{1\columnwidth}
    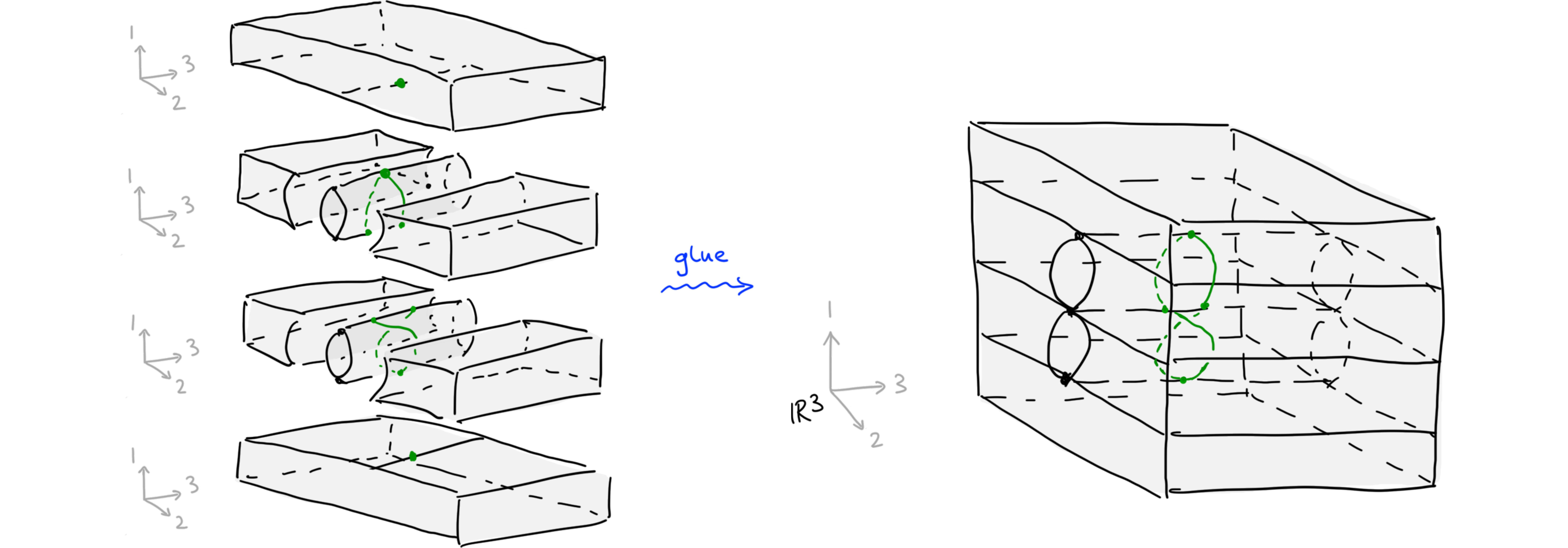

    \caption[Cellulation of the Hopf circle.]{Cellulation the Hopf circle embedding into eight 3-framed 3-cells.}
    \label{fig:triangulating-the-twisted-embedding-of-the-circle}
\end{figure}
We will formally show that this complex is canonical in the sense that any other cellulation by framed regular cells must subdivide it---the precise statement for the case of the Hopf circle embedding will be revisited in \autoref{rmk:resolving-the-twisted-circle-example}.
\end{eg}

\begin{eg}[Dualizing framed cell complexes] \label{eg:self-duality-of-blocks} In \autoref{fig:the-dual-of-the-twisted-embedding-of-the-circle} we depict another flat 3-framed regular cell complex made up of four framed 3-cell: this is the `geometric dual' to the framed cell complex given in the previous example.
\begin{figure}[h!]
    \centering
    \def\svgwidth{1\columnwidth}
    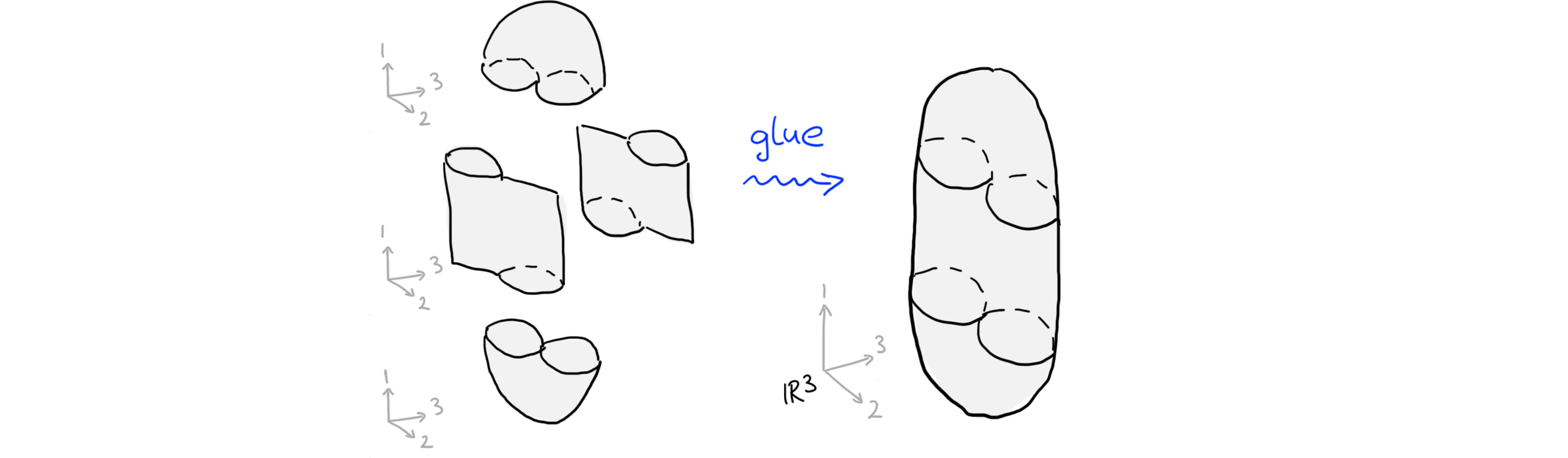

    \caption{Dualizing the cellulation of the Hopf circle.}
    \label{fig:the-dual-of-the-twisted-embedding-of-the-circle}
\end{figure}
The two complexes can be related by a `dualization operation' which switches the roles of dimension and codimension for each cell---we will reconsider this example more formally in \autoref{rmk:dualization-of-blocks}.
\end{eg}

\pauseae

We next define framed maps of framed regular cell complexes; the definition mirrors the definition of framed simplicial maps (see \autoref{defn:framed-maps-of-framed-scplx}).

\begin{defn}[Framed maps of framed regular cell complexes] \label{defn:framed-maps-of-reg-cplx} Consider $n$-framed regular cell complexes $(X,\cF)$ and $(Y,\cG)$. A \textbf{framed cellular map} $F : (X,\cF) \to (Y,\cG)$ is a cellular map $F : X \to Y$ that restricts on all cells $x \in X$ to a framed cellular map $F : (X^{\geq x}, \rest \cF x) \to (Y^{\geq F x},\rest \cG {Fx})$.
\end{defn}

\nid We often refer to `framed cellular maps' simply as `framed maps'.

\begin{notn}[The category of framed regular cell complexes] \label{notn:frccplx} The category of $n$-framed regular cell complexes and their framed maps will be denoted by $\FrCCplx n$. The full subcategory of $n$-framed regular cell complexes that are flat will be denoted by $\FrCDiag n$.
\end{notn}

\pauseae

Finally, let us compare framed regular cell complexes to framed simplicial complexes. Recall, we set out to construct a fully faithful embedding $\FrSCplx n\into \FrCCplx n$ which descends to the ordinary (nonframed) embedding $\SCplx \into \CCplx$ via functors that forget framings. The desired embedding will be given by `framed' version of the entrance path poset construction, which we outline as follows.

\begin{constr}[Framed entrance path posets] \label{constr:framed-entrance-path-poset} Given an $n$-embedded framed simplicial complex $(K,\cF)$, we construct a framed regular cell complex $(\Entr K,\FrEntr \cF)$, referred to as the `framed entrance path poset' of $(K,\cF)$, where $\Entr K$ is the entrance path poset of $K$ and $\FrEntr \cF$ is an $n$-framing inductively constructed as follows.

    First assume $K = S$ is an $m$-simplex. The case $m = 0$ is trivial. For $m > 0$, given $n$-embedded framed $m$-simplex $(S \iso [m],\cF)$ construct the integral proframing $\Intfr \cF$. Assume $\Intfr \cF$ is of the form $(\id,...,\id,p_j,p_{j-1},...,p_1)$ with $p_j \neq \id$. Set $\cP$ to be the proframing $(\id,...,\id,\id,p_{j-1},...,p_1)$ of the $(m-1)$-simplex $[m-1]$. Arguing inductively in $m$, construct the framing $\FrEntr \Gradfr \cP$ of $\Entr [m-1]$. The framing $\FrEntr \cF$ is now defined to label vectors in the kernel of $\Entr p_j : \Entr [m] \to \Entr [m-1]$ by $j \in \bnum n$, and all other vectors by the label that their image under $\Entr p_j$ is given in the framing $\FrEntr \Gradfr \cP$.

    Now assume $K$ is any simplicial complex. Then $\FrEntr \cF$ is the framing that, on each simplex $x : [m] \into K$, restricts on $\Entr x : \Entr [m] \into \Entr K$ to the framing $\rest {(\FrEntr \cF)} {\Entr x} = \FrEntr {(\rest \cF x)}$.
\end{constr}

\begin{constr}[The framed entrance path poset functor] \label{obs:cell-cplx-of-framed-complex} The \textbf{framed entrance path poset functor}
\begin{equation}
        \FrEntr : \FrSCplx n \to \FrCCplx n
\end{equation}
takes $n$-embedded framed simplicial complex $(K,\cF)$ to their framed entrance path posets $(\Entr K,\FrEntr \cF)$, and framed simplicial maps $F : (K,\cF) \to (L,\cG)$ to the framed cellular map determined by the cellular map $\Entr F : \Entr K \to \Entr L$.
\end{constr}

\nid We omit a detailed verification of \autoref{constr:framed-entrance-path-poset} and \autoref{obs:cell-cplx-of-framed-complex}. The punchline is recorded in the following observation.

\begin{obs}[Framed simplicial complexes are framed regular cell complexes] \label{obs:face-functor-fully-faithful} The functor $\FrEntr : \FrSCplx n \to \FrCCplx n$ is a fully faithful embedding of categories.
\end{obs}

\nid In this sense, framed simplicial complexes are specific framed regular cell complexes, i.e.\ the latter notion is a generalization of the former notion.

\subsubsecunnum{The definition of proframings on cell complexes}

We next discuss proframings of regular cell complexes. The definition is an almost verbatim generalization of the case of proframed regular cells. Recall the notion of `orderings' of sequences of simplicial maps (see \autoref{rmk:ordering-pre-cell-proframe}).

\begin{defn}[Proframed regular cells] \label{defn:proj-framed-reg-cell-complex} An \textbf{$n$-proframing $\cP$ of a regular cell complex $X$} is a sequence of surjective cellular maps $X = X_n \xto {p_n} X_{n-1} \xto {p_{n-1}} ... \xto {p_1} X_0$, whose sequence of underlying simplicial maps is further endowed with an ordering, such that the resulting sequence of ordered simplicial maps (denoted again by $\cP$) is an $n$-proframing of the simplicial complex $X$, and with the condition that, for each $x \in X$, the restriction $\rest \cP {X^{\geq x}}$ of $\cP$ to the subcomplex $X^{\geq x}$ is flat.
\end{defn}

\nid We will refer to the pair $(X,\cP)$, of a regular cell $X$ together with an $n$-proframing $\cP$ on it, as an `$n$-proframed regular cell'. As before, we usually abbreviate the cell restrictions $\rest \cP {X^{\geq x}}$ by $\rest \cP x$ (see \autoref{notn:proframes-on-subcells}).

\begin{term}[Flat proframed regular cell complexes] An $n$-proframed regular cell complex $(X,\cP)$ is called `flat' if the proframing $\cP$ is flat.
\end{term}

\begin{eg}[Non-flat 2-proframed regular cell complexes] In \autoref{fig:proframed-regular-cell-complex-examples} we depict two 2-proframed regular cell complexes. (Our notation is an immediate generalization of the case of proframed regular cells.) Neither proframing is flat.
\begin{figure}[ht]
    \centering
    \def\svgwidth{1\columnwidth}
    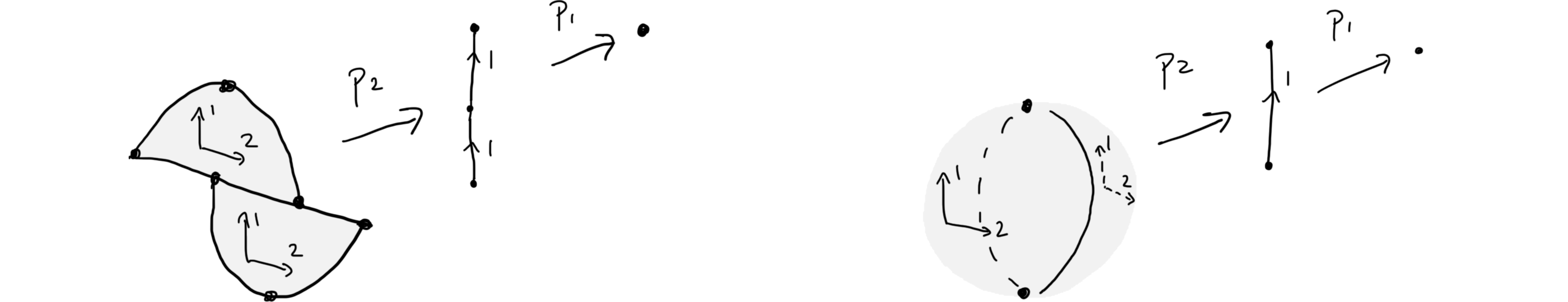

    \caption[Proframed cell complexes and their gradients.]{Two 2-proframed regular cell complexes together with their gradient framings.}
    \label{fig:proframed-regular-cell-complex-examples}
\end{figure}
\end{eg}

\begin{eg}[Flat 3-proframed regular cell complexes] In \autoref{fig:proframed-twisted-circle-triangulation} we depict two 3-proframed regular cell complexes. Both of proframings are flat (and they are depicted as proframed realized in the standard euclidean proframe).
\begin{figure}[ht]
    \centering
    \def\svgwidth{1\columnwidth}
    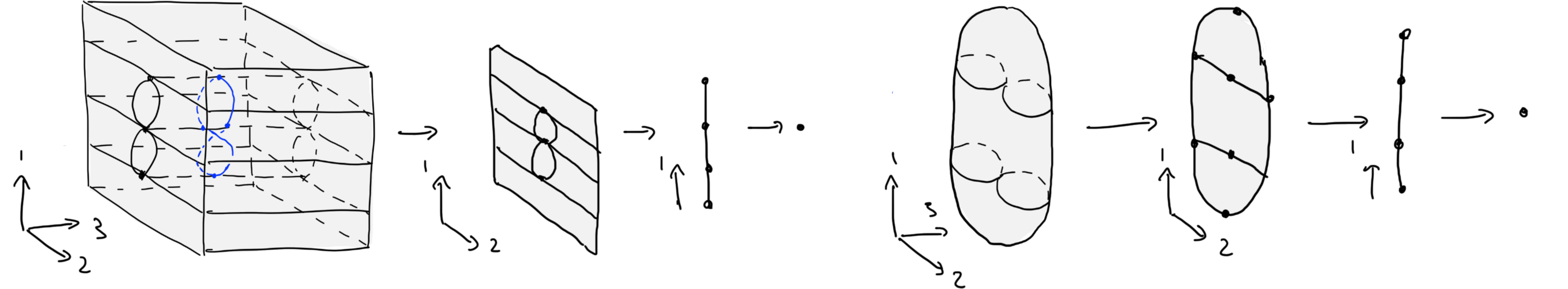

    \caption[Proframing of the Hopf circle cellulation.]{Proframing of the Hopf circle cellulation and a proframing of its dual.}
    \label{fig:proframed-twisted-circle-triangulation}
\end{figure}
\end{eg}

\pauseae

We next define framed maps of framed regular cell complexes; the definition is directly parallel to the definition of proframed simplicial maps (see \autoref{defn:maps-of-proj-framing}).

\begin{defn}[Proframed maps of proframed regular cell complexes] \label{defn:proframed-maps-of-reg-cplx} Consider $n$-proframed regular cell complexes $(X,\cP)$ and $(Y,\cQ)$. A \textbf{proframed cellular map} $F : (X,\cP) \to (Y,\cQ)$ is a map of sequences $F = (F_n,F_{n-1},...,F_1,F_0) : \cP \to \cQ$ (consisting of cellular maps $F_i$) which, for all $x \in X$, restricts to a proframed cellular map $F : (X^{\geq x}, \rest \cP x) \to (Y^{\geq F x},\rest \cQ {Fx})$.
\end{defn}

\nid We often refer to `proframed cellular maps' simply as `proframed maps'.

\begin{notn}[The category of proframed regular cell complexes] The category of $n$-proframed regular cell complexes and their proframed cellular maps will be denoted by $\PFrCCplx n$. The full subcategory consisting of flat $n$-proframings will be denoted by $\PFrCDiag n$.
\end{notn}

\subsubsecunnum{Gradient cellular framings and integral cellular proframings}

Parallel to the definition of gradient framings of proframed simplicial complexes (see \autoref{defn:grad-of-profr-simp-cplx}), proframed regular cell complexes too have gradient framings as follows.

\begin{defn}[Gradient framings of proframed regular cell complexes] Given an $n$-proframed regular cell complex $(X,\cP)$, its \textbf{gradient framing} is the $n$-framed regular cell complex $(X,\Gradfr \cP)$ whose framing is given by the gradient framing $\Gradfr \cP$ of the simplicial complex $P$. Further, given a proframed cellular map $F = (F_n,F_{n-1},...,F_1,F_0) : (X,\cP) \to (Y,\cQ)$ its \textbf{gradient map} is the framed cellular map $(X,\Gradfr \cP) \to (Y,\Gradfr \cQ)$ determined by the cellular map $F_n : X \to Y$.
\end{defn}

\begin{term}[Gradient framing functor] Gradients of proframings of regular cell complexes and their proframed maps assemble into a `gradient framing' functor as follows:
\begin{equation}
        \Gradfr : \PFrCCplx n \to \FrCCplx n . \qedhere
\end{equation}
\end{term}

\begin{eg}[Gradients of proframed regular cells] \label{eg:proframed-cell-gradients}
Passing to the gradient $2$-framings of the $2$-proframed cell depicted in \autoref{fig:2-proframed-regular-cell-examples} recovers the upper four examples of framings given in \autoref{fig:framed-regular-cells-further-examples}.
\end{eg}

\begin{eg}[Gradients of proframed regular cell complexes]
In our earlier \autoref{fig:proframed-regular-cell-complex-examples} we depicted two 2-proframed regular cell complexes; their respective gradient framings recover the 2-framed regular cell complexes from our earlier \autoref{fig:block-complex-unit-counit-gluing}. Similarly, in \autoref{fig:proframed-twisted-circle-triangulation} we depicted two 2-proframed regular cell complexes whose respective gradient framings recover the 2-framed regular cell complexes from our earlier \autoref{fig:triangulating-the-twisted-embedding-of-the-circle} and \autoref{fig:the-dual-of-the-twisted-embedding-of-the-circle}.
\end{eg}

\begin{defn}[Integral proframings] Given an $n$-framed regular cell complex $(X,\cF)$, an $n$-proframed regular cell complex $(X,\cP)$ is an \textbf{integral $n$-proframing} for $(X,\cF)$ if $(X,\cF)$ is the gradient framing of $(X,\cP)$.
\end{defn}

\nid The previous definitions leave open the question of whether each $n$-framed regular cell complexes has an integral $n$-proframing: unlike in the case of simplices, this is now even non-trivial in the case of framed cells, since the definition of proframed regular cells requires a `sequence of surjective cellular maps' as part of its data. Nonetheless, we have the following result.

\begin{thm}[Equivalence of framed and proframed regular cells] \label{obs:frames-vs-proj-frames-for-reg-cells} The gradient framing functor is an equivalence of categories of proframed regular cells and framed regular cells:
    \begin{equation}
        \Gradfr : \PFrCCell n \iso \FrCCell n .
    \end{equation}
The inverse to $\Gradfr$ will be called \emph{integration}, and denoted by $\Intfr$.
\end{thm}

\nid The proof will be given in \autoref{ssec:frames-vs-proj-frames-reg-cell}. The theorem may, in fact, be phrased in yet more general form: namely, parallel to the case of simplicial complexes, the notion of framings and proframings on regular cell complexes become equivalent when flatness is imposed.

\begin{thm}[Equivalence of flat framed and flat proframed regular cell complexes] \label{obs:frames-vs-proj-frames-for-reg-diags} The gradient framing functor is an isomorphism of categories of flat proframed regular cell complexes and flat framed regular cell complexes:
    \begin{equation}
        \Gradfr : \PFrCDiag n \iso \FrCDiag n
    \end{equation}
The inverse to $\Gradfr$ will be called \emph{integration}, and denoted by $\Intfr$.
\end{thm}

\nid The proof can be found in \autoref{ssec:frames-vs-proj-frames-reg-cell}.

\pause

\begin{figure}[!ht]
\centering
\adjustbox{scale=.85,center}{
\begin{tikzcd}
 &[+15pt] &[-20pt] \PFrDelta n \arrow[lld, "\FrEntr" description, pos=0.2, hook'] \arrow[r, "\Gradfr", shift left=2]{}[swap]{\text{\tiny(see \ref{cor:iso-of-frame-vs-proj-frame-simplex})}} \arrow[d, hook] &[+15pt] \FrDelta n \arrow[l, "\Intfr", shift left=2] \arrow[lld, "\FrEntr" description, pos=0.8, hook', crossing over] \arrow[d, hook] & \\[+25pt]
\PFrCCell n \arrow[d, hook] \arrow[r, "\Gradfr", shift left=2]{}[swap]{\text{\tiny(see \ref{obs:frames-vs-proj-frames-for-reg-cells})}} & \FrCCell n \arrow[l, "\Intfr", shift left=2] & \FlPFrSCplx n \arrow[r, "\Gradfr", shift left=2]{}[swap]{\text{\tiny(see \ref{obs:frames-vs-proj-frames-flat})}} \arrow[d, hook] \arrow[lld, "\FrEntr" description, pos=0.2, hook'] & \FlFrSCplx n \arrow[l, "\Intfr", shift left=2] \arrow[d, hook] \arrow[lld, "\FrEntr" description, pos=0.8, hook', crossing over] & \\[+40pt]
\PFrCDiag n \arrow[r, "\Gradfr", shift left=2]{}[swap]{\text{\tiny(see \ref{obs:frames-vs-proj-frames-for-reg-diags})}} \arrow[d, hook] & \FrCDiag n \arrow[l, "\Intfr", shift left=2] & \PFrSCplx n \arrow[r] \arrow[lld, "\FrEntr" description, pos=0.2, hook'] & \FrSCplx n \arrow[lld, "\FrEntr" description, pos=0.8, hook'] \arrow[r, "\Unframe"] & \SCplx \arrow[lld, hook'] \\[+40pt]
\PFrCCplx n \arrow[r] & \FrCCplx n \arrow[r, "\Unframe"'] & \CCplx & & \arrow[from=2-2, to=3-2, hook, crossing over] \arrow[from=3-2, to=4-2, hook, crossing over]
\end{tikzcd}
}\vspace{15pt}
\caption[Categories of framings and proframings.]{Categories of framings and proframings of simplicial complexes and cell complexes.}
\label{fig:diagrams-chapter-1}
\end{figure}
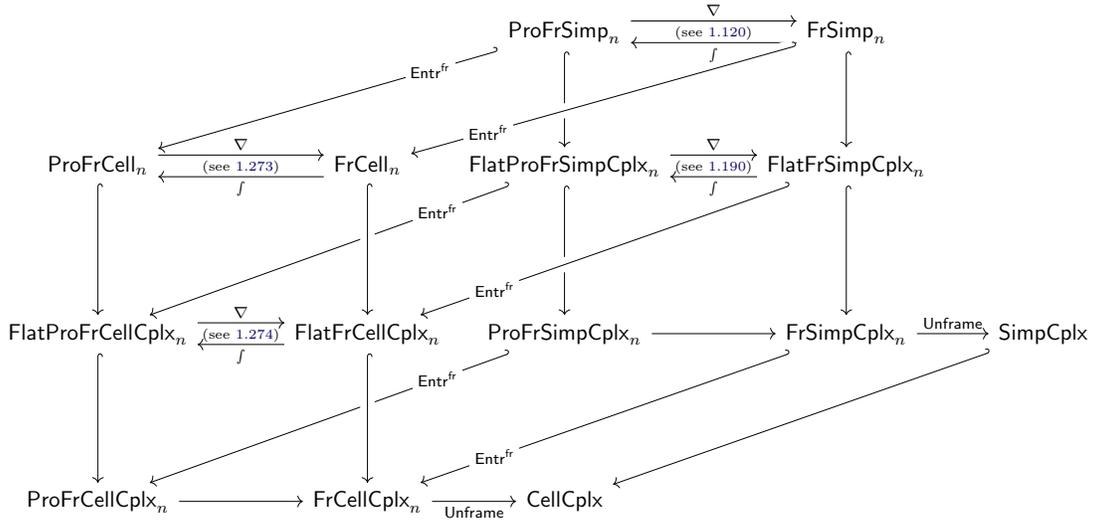

Both theorems, together with our earlier results and constructions, then organize into the diagram in \autoref{fig:diagrams-chapter-1}, which summarizes the categories defined in this chapter. (Note that except for the isomorphism $\PFrDelta n \iso \FrDelta n$, which was explained in \autoref{cor:iso-of-frame-vs-proj-frame-simplex}, all isomorphisms in \autoref{fig:diagrams-chapter-1} are yet to be constructed; we will construct the remaining isomorphisms in \autoref{sec:equiv-of-frames-and-proj-frames}).

\chapter{Constructible framed combinatorics: trusses} \label{ch:trusses}

In this chapter we develop the basic theory of `trusses': trusses are combinatorial structures that can be obtained by iteratively building towers of `constructible combinatorial bundles'. Trusses both underlie the tractable classification of framed regular cell complexes (as discussed in \autoref{ch:classification-of-framed-cells}) while also providing a link between framed combinatorics and framed topology (as discussed in \autoref{ch:meshes} and \autoref{ch:hauptvermutung}). We will start in \autoref{sec:1-trusses} by introducing 1-dimensional trusses, or `1-trusses' as well as their `bordism' and `bundles'. Subsequently, in \autoref{sec:truss-induction}, we describe an important ordering property of bundles of 1-trusses, which we will refer to as `truss induction'. This in turn will allow us to describe the theory of trusses in general dimension $n$ as discussed in \autoref{sec:n-trusses}.

\section{1-Trusses, 1-truss bordisms, and 1-truss bundles} \label{sec:1-trusses}

1-Trusses are `1-dimensional' combinatorial objects. As we will illustrate, their 1-dimensionality has a concrete geometric interpretation: 1-trusses can be obtained by taking the `\emph{framed} entrance path posets' of `\emph{framed} stratified intervals'. Here, the notion of `framed' entrance path poset generalizes the usual notion of entrance path posets of stratifications to the case of (appropriately) `framed' stratifications. While both these notions will be formalized only later in \autoref{ch:meshes} in the context of `constructible framed topology', thinking of 1-trusses as modelling stratified intervals will be nonetheless useful as a guiding intuition, and we will illustrate this intuition with several examples throughout this section.

We outline the section. In \autoref{ssec:1-trusses} we introduce the category of 1-trusses and their maps. We will discuss several (often geometrically motivated) notions and constructions for 1-trusses. In particular, we define an involutive `dualization' functor on 1-trusses. In geometric terms, this will correspond to dualizing the dimension of each stratum in a stratification. In \autoref{ssec:1-truss-bord} we then introduce a notion of 1-truss bordisms. In geometric terms, the notion models the behavior of `stratified intervals in $\lR$-parametrized families'. Understanding 1-truss bordisms, will in turn allow us to define 1-trusses bundles in \autoref{ssec:1-truss-bun}, which then provide a combinatorial model of general `constructible stratified line bundles'. 1-Truss bundles play a central role: by iterating such combinatorial constructible stratified line bundles we will later build combinatorial models of `$n$-dimensional' stratified space, yielding the notion of `$n$-trusses'.

\subsection{1-Trusses} \label{ssec:1-trusses}

\subsubsecunnum{1-Trusses as framed fences}

\nid Recall the classical combinatorial notion of fences (1-trusses will be a strengthening of this notion). The definition of fences finds a concise expression via the notion of a classifying space of a category: recall this is the space realizing the simplicial set given by the nerve of a category.

\begin{defn}[Fence] \label{defn:fences}
    A \textbf{fence} $X$ is a category of countable size, whose classifying space is a connected $k$-manifold where $k \in \{0,1\}$.
    \begin{enumerate}
        \item If $F$ realizes to a point we say $F$ is \textbf{trivial}.
        \item If $F$ realizes to an interval we say $F$ is \textbf{linear}.
        \item If $F$ realizes to a circle we say $F$ is \textbf{circular}. \qedhere
    \end{enumerate}
\end{defn}

\nid We call a fence $X$ a `finite fence' if it is finite as a category. For convenience, we often extend both the class of linear and circular fences to include the trivial fence as well (and speak of the `trivial linear' resp.\ `trivial circular' fence in this case). Note that a fence cannot contain composable non-identity morphisms.

\begin{eg}[Fences] \label{eg:fences} Finite fences of different types are shown in \autoref{fig:fences-of-different-types}; in each case we indicate a fence $X$ by its classifying space $\abs{X}$ (where directions of morphisms in $X$ are recorded by directing edges of $\abs{X}$).
\begin{figure}[ht]
    \centering
    \def\svgwidth{1\columnwidth}
    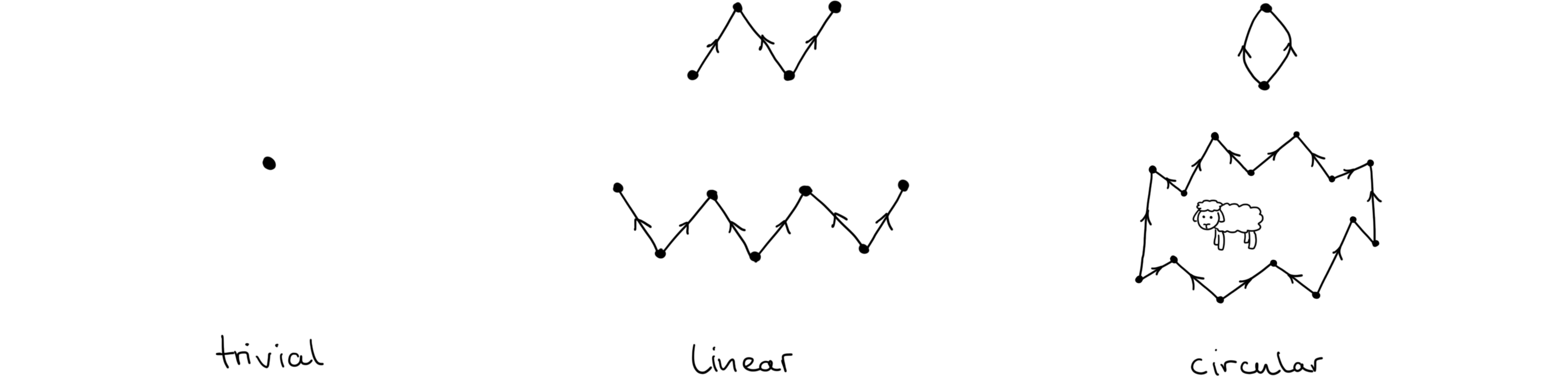

    \caption{Fences of different types.}
    \label{fig:fences-of-different-types}
\end{figure}
\end{eg}

\begin{obs}[Regular cell complexes] Given a finite fence $X$, then the classifying space $\abs{X}$ of $X$ naturally has the structure of a regular cell complex (whose 0-cells realize objects of $X$, and whose 1-cells realize morphisms of $X$).
\end{obs}

\begin{obs}[Height maps] Every non-trivial fence $X$ comes equipped with `height map' $X \to [1]$ given by the unique functor whose preimages are discrete categories (the map `folds the fence onto a single fence post'). Note that for trivial fences, the choice of height map $X \to [1]$ is ambiguous!
\end{obs}

\nid Conversely, regular cell complexes with height maps yield fences as follows.

\begin{obs}[Non-trivial fences are regular cell complexes with height map] Any regular cell complex $X$ cellulating a connected 1-manifold, together with a cellwise bijective map $X \to \abs{[1]}$, determines a non-trivial fence whose classifying space is $X$ and whose height map realizes to $X \to \abs{[1]}$.
\end{obs}

\nid The notion of 1-trusses strengthens that of fences in two ways: firstly, 1-trusses contain the data of a `framing' of the 1-manifolds that non-trivial fences cellulate (technically, this uses the combinatorial notion of locally flat 1-framings of regular cell complexes, see \autoref{defn:framed-reg-cell-cplx} and \autoref{eg:locally-flat-cellulations}); secondly, 1-trusses resolve the ambiguity of height maps for trivial fences by introducing them as additional structure (inverting variance, this now will be recorded by a map to $[1]\op$, which we call the `dimension map').

\begin{defn}[General 1-trusses] \label{defn:general-1-trusses} A \textbf{1-truss} $(T,\dim,\cF)$ is a finite fence $T$, together with a poset map $\dim : T \to [1]\op$ as well as a locally flat 1-framing $\cF$ of the corresponding regular cell complex.  A 1-truss is \textbf{trivial} resp.\  \textbf{linear} resp.\  \textbf{circular} whenever its underlying fence is.
\end{defn}

\begin{eg}[General 1-trusses] \label{eg:general-1-trusses} 1-Trusses of different types are shown in \autoref{fig:1-trusses-of-different-types}: in each case we depict the underlying fence as before, we depict the dimension map by coloring preimages of $0$ in red, and preimages of $1$ in blue, and we depict the locally flat 1-framing by indicating the 1-framing of each cell by a purple frame vector.
\begin{figure}[ht]
    \centering
    \def\svgwidth{1\columnwidth}
    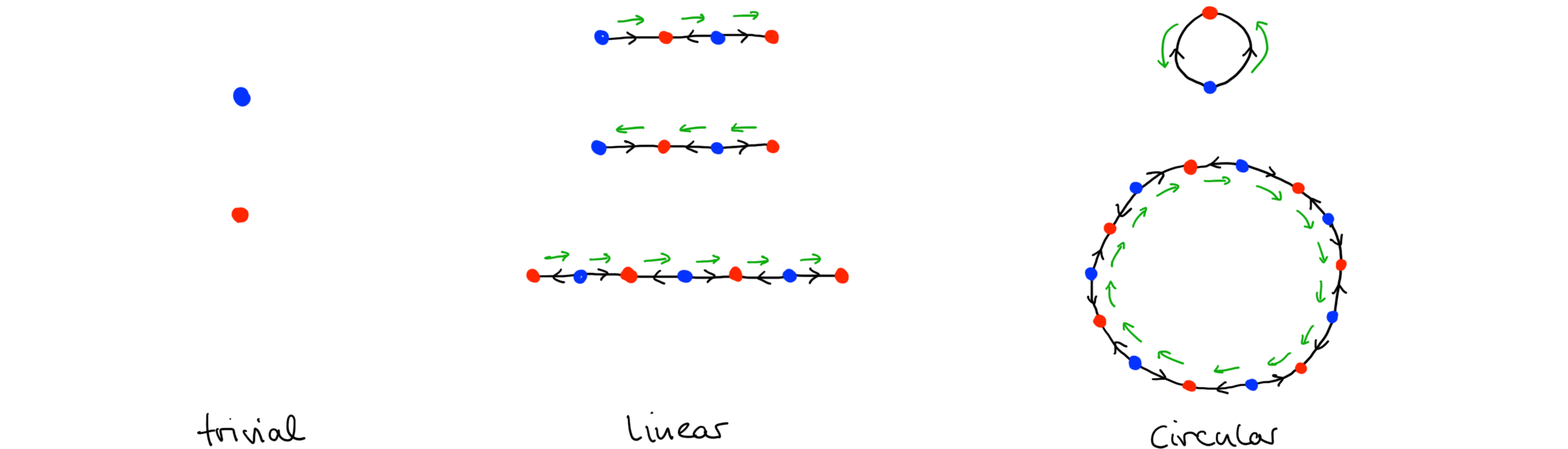

    \caption{1-Trusses of different types.}
    \label{fig:1-trusses-of-different-types}
\end{figure}
\end{eg}

\nid Going forward, we now make the following simplification: we will henceforth focus solely on the case of linear trusses. While much of the theory developed here, including notions of `higher-dimensional' $n$-trusses, does generalize to the case of general 1-trusses as well, our main interest and applications will ultimately lie with the linear case. We therefore adopt the following convention.

\begin{conv}[Linear 1-truss by default] We will use the term `1-truss' to refer to `linear 1-trusses' unless otherwise noted.
\end{conv}

\nid In the linear case, we may further reformulate the above definition in purely order-theoretic terms. First, note that if a fence $T$ is linear (which we take to include the trivial case), then it must be a poset.\footnote{\label{notn:preorder-as-categories} Every preorder (resp.\ poset) $(X,\leq)$ is equivalently a category with object set $X$ and a single morphism $x \to y$ whenever $x \leq y$; and every map of preorders is a functor of corresponding categories.} We usually denote the resulting poset order of linear fences $T$ either by $(T,\to)$ or $(T,\eleq)$ (the latter notation being convenient for expressing `strictly less/greater' relations using the symbols `$\eles$/$\egre$'). Secondly, a locally flat 1-framing $\cF$ of a linear fence $T$ may be equivalently recorded by a total order, which we usually denote by $\fleq$. This leads us to the following definition of 1-trusses (i.e.\ linear 1-trusses).

\begin{defn}[1-Trusses] \label{defn:1-trusses}
    A (linear) \textbf{1-truss} $(T,\eleq,\dim,\fleq)$ is a finite non-empty set $T$ together with the following structure.
\begin{enumerate}
\item a partial order $\eleq$ called the `face order' of $T$,
\item a poset map $\dim : (T,\eleq) \to [1]\op$ called the `dimension map' of $T$,
\item a total order $(T,\fleq)$ called the `frame order' of $T$, in which $a$ succeeds $b$ or $b$ succeeds $a$ iff either $a \eles b$ or $b \egre a$. \qedhere
\end{enumerate}
\end{defn}

\begin{notn}[1-Trusses] When working with 1-trusses, we will usually keep face orders, dimension maps as well as frame orders implicit; that is, we will abbreviate 1-trusses $(T,\eleq,\dim,\fleq)$ simply by $T$.
\end{notn}

\nid The fact that we refer to $\eleq$ as the `face order' of $T$ reflects the relation of 1-trusses and (framed) stratified intervals that we alluded to earlier: under this relation, elements $a$ of 1-trusses translate to open 1-cells of dimension $\dim(a)$, and arrows $a \eles b$ translate to `$b$ being a face of $a$'. The next example illustrates this `geometric translation' of 1-trusses.\footnote{Note that the resulting translation of 1-trusses to stratified intervals is fundamentally different to our earlier translation of fences to regular cell complexes.}

\begin{eg}[Translating 1-trusses to stratified intervals] \label{eg:1-trusses} In the upper row of \autoref{fig:1-trusses-and-their-realization-as-stratified-intervals} we depict two 1-trusses $T$. In each case, we depict the face order by black arrows, we color preimages of $1$ of the dimension map in blue and preimages of $0$ in red, and indicate the frame order by a purple coordinate axis.  Underneath each 1-truss, we depict a stratified interval: note that each object $x \in T$ corresponds to a stratum of dimension $\dim(x)$, and the face order $(T,\eleq)$ can be recovered as the entrance path poset of the stratified interval (see \autoref{defn:entr} for a definition of `entrance path posets').
\begin{figure}[ht]
    \centering
    \def\svgwidth{1\columnwidth}
    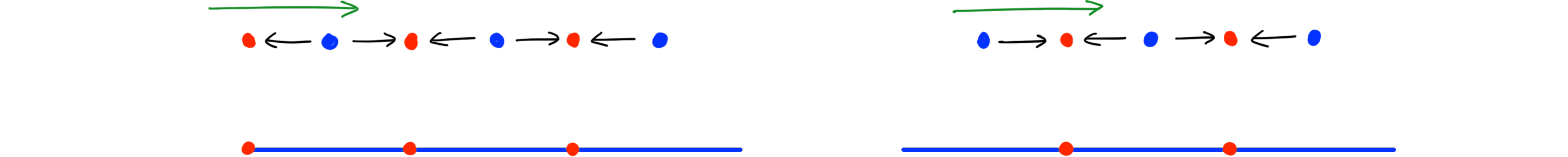

    \caption{1-Trusses and their translation to stratified intervals.}
    \label{fig:1-trusses-and-their-realization-as-stratified-intervals}
\end{figure}
\end{eg}

\begin{term}[Singular and regular objects] \label{term:sing_reg_objects} An element $p \in T$ of a 1-truss $T$ is called `singular' if $\dim(p) = 0$ and `regular' if $\dim(p) = 1$. We denote the subset of singular objects $T$ by $\sing(T)$ resp.\ the subset of regular objects by $\reg(T)$.
\end{term}

\begin{rmk}[Induced orders on singular/regular objects] \label{rmk:reg-sing-subposets} Note that, endowed with the face order both $(\sing(T),\eleq)$ and $(\reg(T),\eleq)$ are discrete orders, and endowed instead with the frame order both $(\sing(T),\fleq)$ and $(\reg(T),\fleq)$ are total orders.
\end{rmk}

\subsubsecunnum{Maps of 1-trusses} We next define maps between 1-trusses. In order to simultaneously keep track of face and frame orders in the definition of 1-trusses, the following abstraction will be useful.

\begin{defn}[Diposets and their maps] \label{defn:diposets} A \textbf{diposet} $(X,\eleq,\fleq)$ is a set $X$ with two orders $\eleq$ and $\fleq$. A \textbf{diposet map} $F : (X,\eleq,\fleq) \to (Y,\eleq,\fleq)$ is a map of sets $F : X \to Y$ that separately respects both orders, i.e.\ induces poset maps $F : (X,\eleq) \to (Y,\eleq)$ and $F : (X,\fleq) \to (Y,\fleq)$.
\end{defn}

\begin{defn}[1-Truss maps] \label{defn:1-truss-map} A \textbf{map of 1-trusses} $T \to S$ is a diposet map $(T,\eleq,\fleq) \to (S,\eleq,\fleq)$.
\end{defn}

\nid Note that the preceding definition of 1-truss maps does not impose any conditions on how maps interact with the dimension data of trusses; there are several sensible ways to impose such conditions.

\begin{defn}[Regular, singular, and balanced maps] \label{defn:reg-sing-cell-1-truss-map} Let $F :T \to S$ be a map of 1-trusses.
    \begin{enumerate}
        \item We call $F$ \textbf{singular} if it maps singular objects of $T$ to singular objects of $S$. That is, for all $x \in T$, $\dim(x) \geq \dim(Fx)$.
        \item We call $F$ \textbf{regular} if it maps regular objects of $T$ to regular objects of $S$. That is, for all $x \in T$, $\dim(x) \leq \dim(Fx)$.
        \item We call $F$ \textbf{balanced} if it is both regular and singular. That is, for all $x \in T$, $\dim(x) = \dim(Fx)$. \qedhere
    \end{enumerate}
\end{defn}

\begin{eg}[Maps of 1-trusses] In \autoref{fig:maps-of-1-trusses} we depict examples of singular, regular and balanced maps of 1-trusses.
\begin{figure}[h!]
    \centering
    \def\svgwidth{1\columnwidth}
    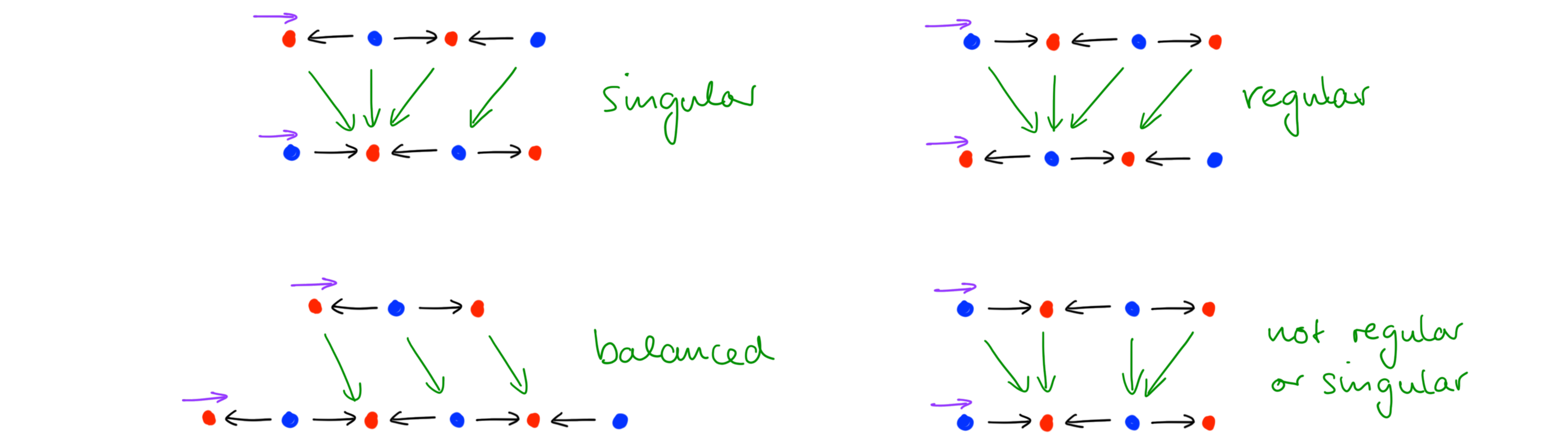

    \caption{Four different types of maps of 1-trusses.}
    \label{fig:maps-of-1-trusses}
\end{figure}
\end{eg}

\begin{notn}[Category of 1-trusses] \label{notn:1truss-category} The category of 1-trusses and their maps is denoted by $\truss 1$. The wide subcategory containing only regular respectively singular respectively balanced maps will be denoted by $\rtruss 1$ respectively $\struss 1$ respectively $\rstruss 1$.
\end{notn}

\nid To avoid identifying the trivial truss $T$ with $\dim(a) = 0$ (where $a \in T$ is the single element of $T$) with the trivial truss $T$ with $\dim(a) = 1$ we will default to the following notion of isomorphism.

\begin{term}[Balanced isomorphism by default] \label{rmk:1-truss-iso} Unless otherwise noted, the term `isomorphism of 1-trusses' will by refer to isomorphisms in the category $\rstruss 1$, that is, to balanced bijective truss maps.
\end{term}

\begin{rmk}[Balanced isomorphisms are unique] \label{rmk:1-trusses-skeletally} Observe that if two trusses are balanced isomorphic, then the balanced isomorphism must be unique. There is thus no harm in thinking about 1-trusses `skeletally', i.e.\ in terms of their balanced isomorphism classes.\footnote{In fact, note that (not necessarily balanced) isomorphisms between any two trusses are unique as well if they exists; the reason for working with balanced isomorphism classes roots in their relation to `1-truss bordism isomorphism classes', as discussed later in \autoref{rmk:1-trusses-skeletally-bordism}.}
\end{rmk}

\nid We now distinguish balanced isomorphism classes by their `endpoint type'.

\begin{defn}[Endpoints of 1-trusses] \label{term:endpoints} Given a truss $(T,\eleq,\dim,\fleq)$ we denote by $\ept_- T$ the minimal element of $(T,\fleq)$, called the \textbf{lower endpoint}, and by $\ept_+ T$ the maximal element, called the \textbf{upper endpoint} of $T$.
\end{defn}

\begin{term}[Endpoint types of 1-trusses] \label{rmk:trusses-by-boundary-type} Let $T$ be a 1-truss. Then $T$ falls in exactly one of the following six cases of isomorphism classes in $\rstruss 1$.
    \begin{enumerate}
        \item If $T$ has a single regular element, then $T$ is the \textbf{trivial open} 1-truss denoted by $\CTT_0$.

        \item If $T$ has a single singular element, then $T$ is the \textbf{trivial closed} 1-truss denoted by $\OTT_0$.
    \end{enumerate}

\nid For the next cases, $T$ has more than one element.

    \begin{enumerate}[resume]
        \item If both endpoints $\ept_\pm T$ of $T$ are regular, then $T$ is \textbf{open}. If $T$ has $2k + 1$ elements it is denoted by $\OTT_k$.

        \item If both endpoints $\ept_\pm T$ of $T$ are singular, then $T$ is \textbf{closed}. If $T$ has $2k + 1$ elements it is denoted  by $\CTT_k$.

        \item If $\ept_- T$ is regular and $\ept_+ T$ singular, then $T$ is \textbf{half-open half-closed}. If $T$ has $2k$ elements it is denoted by $\OCTT_{\!k}$.

        \item If $\ept_- T$ is singular and $\ept_+ T$ regular, then $T$ is \textbf{half-closed half-open}. If $T$ has $2k$ elements it is denoted by $\COTT_{\!k}$. \qedhere
    \end{enumerate}
\end{term}

\begin{eg}[Types of 1-trusses] In \autoref{fig:cellular-isomorphism-classes-of-1-trusses} we illustrate each type of 1-truss with an example.
\begin{figure}[h!]
    \centering
    \def\svgwidth{1\columnwidth}
    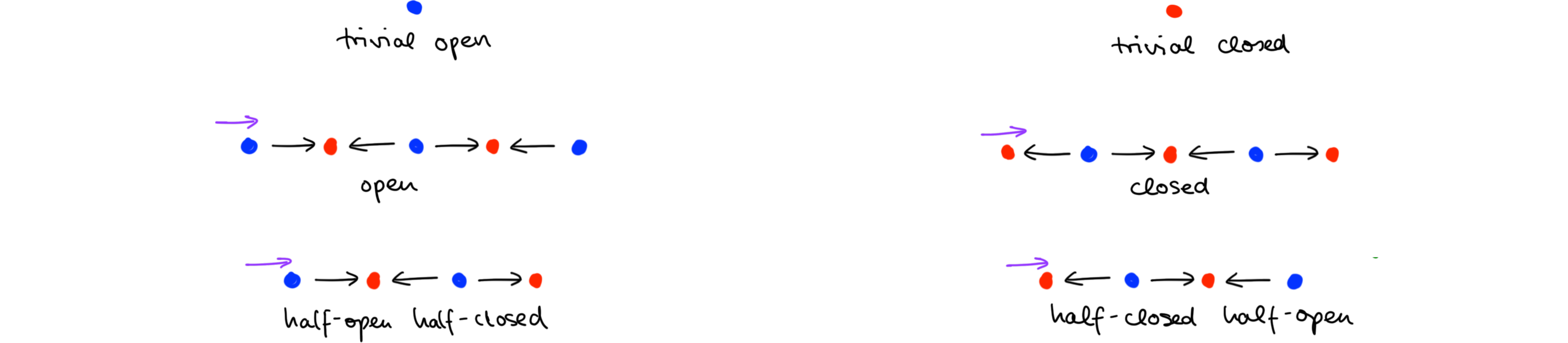

    \caption{Types of 1-trusses.}
    \label{fig:cellular-isomorphism-classes-of-1-trusses}
\end{figure}
\end{eg}

\nid We will be mainly be interested in the case of closed (and, dually, open) trusses, because of their relation to framed cellular geometry. We therefore introduce the following notation for their corresponding categories.

\begin{notn}[Open and closed 1-trusses] \label{notn:1truss-categories-open-closed} The subcategory of $\truss 1$ containing only open trusses (including the trivial open truss) will be denoted by $\otruss 1$, and the subcategory of $\truss 1$ containing only closed trusses (including the trivial closed truss) by $\ctruss 1$. Similarly to \autoref{notn:1truss-category}, we add a superscript `$\mathsf{s}$' resp.\ `$\mathsf{r}$' resp.\ `$\mathsf{rs}$'  to indicate the corresponding wide subcategories $\sotruss 1$ resp.\ $\rotruss 1$ resp.\ $\rsotruss 1$ (as well as $\sctruss 1$ resp.\ $\rctruss 1$ resp.\ $\rsctruss 1$) containing only singular resp.\ regular resp.\ balanced maps.
\end{notn}

\begin{obs}[Maps that are balanced] \label{rmk:automatic-cellularity-1} Note that any regular map of open trusses is balanced and, likewise, any singular map of closed trusses is balanced. In other words, the categories $\rotruss 1$ and $\rsotruss 1$ (resp.\ $\sctruss 1$ and $\rsctruss 1$) are identical.
\end{obs}

\subsubsecunnum{Dualization of 1-trusses} Trusses admit a natural dualization operation, which maps closed trusses to open trusses and vice-versa.

\begin{constr}[Dualization of 1-trusses] \label{constr:1-truss-duals} The \textbf{dualization functor}
    \begin{equation}
        \dagger : \truss 1 \iso \truss 1
    \end{equation}
    is an involutive functor defined as follows. Given a 1-truss $T \equiv (T,\eleq,\dim,\fleq)$ its dual is the 1-truss $T^\dagger \equiv (T,\eleq\op,\dim\op,\fleq)$: that is, face order of $T^\dagger$ is opposite to that of $T$; its dimension map is the opposite of the dimension map of $T$ (post-composed with the identification $[1] \iso [1]\op$); its frame order is equal to that of $T$. Similarly, the dual map $F : T^\dagger \to S^\dagger$ of a 1-truss map $F : T \to S$ is the map that equals $F$ as a map of underlying sets.
\end{constr}

\begin{eg}[Dualization] In \autoref{fig:cellular-isomorphism-classes-of-1-trusses}, 1-trusses in the left column dualize to 1-trusses shown in the right column; while in \autoref{fig:maps-of-1-trusses} the depicted regular map dualizes to the depicted singular map.
\end{eg}

\begin{obs}[Dualization of dimension] Given a 1-truss $T$ and an element $p \in T$, then $p$ is regular (resp.\ singular) in $T$ if and only if $p$ is singular (resp.\ regular) in $T^\dagger$. In particular, if $F : T \to S$ is a map of trusses, then $F^\dagger$ is regular (resp.\ singular) if and only if $F^\dagger$ is singular (resp.\ regular).
\end{obs}

\begin{obs}[Dualization of boundary types] The preceding observation entails that the dualization functor restricts to isomorphisms $\otruss 1 \iso \ctruss 1$ as well as $\rotruss 1 \iso \sctruss 1$ and $\sotruss 1 \iso \rctruss 1$.
\end{obs}

\subsection{1-Truss bordisms} \label{ssec:1-truss-bord}

\begin{figure}[h!]
    \centering
    \def\svgwidth{1\columnwidth}
    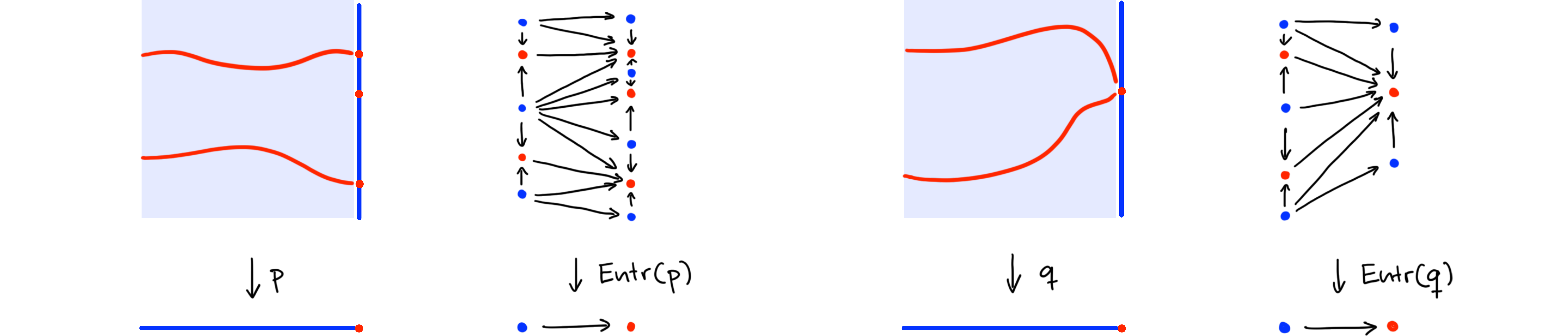

    \caption[Continuous families of open stratified intervals.]{Continuous families of open stratified intervals as stratified bundles over the stratified 1-simplex.}
    \label{fig:families-of-stratified-intervals-as-bundles}
\end{figure}

We now introduce 1-truss bordisms. While 1-trusses geometrically translate into (framed) stratified intervals as previously illustrated, 1-truss bordisms model changes that can occur in continuous families of such intervals: these families can be formalized in terms of (framed constructible) stratified bundles\footnote{Recall, stratified bundles generalize fiber bundles in that they allow (stratified) fibers to change when passing between strata in the base, see \autoref{defn:stratified-bundle}. The precise `framed constructible' variation alluded to here is introduced in \autoref{defn:1-mesh-bundle} under the name of `1-mesh bundles'.} over the stratified 1-simplex $\Set{1} \subset [0,1]$, which is illustrated in two instances in \autoref{fig:families-of-stratified-intervals-as-bundles}. In the first example, the generic fiber is a stratified open interval with two point strata; when reaching the special fiber, a third point stratum spontaneously appears. In the second example, the two point strata in the generic fiber converge into a single point stratum in the special fiber. In both cases, we also indicate the corresponding `poset bundle' obtained by passing to entrance path posets.
The notion of 1-truss bordisms will provide a unified combinatorial description of the `deformations between fibers' that may appear in such entrance path poset bundles (while appropriately keeping track of framings and dimension data of 1-trusses). The description of these fiber deformations plays a most central role: it will shortly enable us to define 1-truss bundles, which in turn, `by iteration', will later on lead us to the definition of $n$-truss bundles as towers of 1-truss bundles. Note that our chosen terminology hints at a relation of `truss bordisms' to `classical \emph{stratified} bordisms'.\footnote{Most visually, this becomes apparent with later examples of (higher-dimensional) trusses underlying stratifications such as those in \autoref{fig:tangles-singularities}.}

\subsubsecunnum{Recollection on Boolean profunctors} To streamline the definition of 1-truss bordisms it is convenient to encode their general structure in terms of Boolean profunctors between preorders. Recall, the category $\Bool$ of booleans has two objects `true' (also written `1') and `false' (also written `0'), with a single non-identity arrow from `false' to `true'. It will be useful to think of preorders as $\Bool$-enriched categories (which translate back to categories via the fully faithful inclusion $\Bool \into \SetCat$).

\begin{defn}[Boolean profunctors and (co)functoriality]
    Given two preorders $X$ and $Y$, a \textbf{boolean profunctor} $R : X \proto{} Y$ is a functor $R : X \op \times Y \to \Bool$. $R$ is called \textbf{functorial} if it is of the form $\Hom_Y(f-,-)$ (for a functor $f : X \to Y$) and \textbf{cofunctorial} if it is of the form $\Hom_X(-,f-)$ (for a functor $f : Y \to X$).\footnote{Frequently, (co)functorial profunctors are also called `(co)representable profunctors'; however, we find this terminology badly chosen since the notion does not relate to an actual representability condition. Moreover, the terminology does not `decategorify' well to the case of relations, where one often speaks of `(co)functional' relations as mentioned in \autoref{rmk:cofunctionality}.}
\end{defn}

\begin{rmk}[Composition of boolean profunctors] \label{rmk:boolean-profunctor-composition} Given Boolean profunctors $R: X \proto{} Y$ and $S: Y \proto{} Z$, then their composition $S \circ R: X \proto{} Z$ is determined by setting $Z(x \in X, z \in Z)$ to be true if and only if there is an element $y \in Y$ so that both $R(x,y)$ and $S(y,z)$ are true.
\end{rmk}

\begin{notn}[Category of boolean profunctors] Preorders and their boolean profunctors form a category denoted by $\BoolProf$.
\end{notn}

\begin{defn}[Underlying relations of boolean profunctors] \label{rmk:underlying-relations} The \textbf{underlying relation functor} $\URel : \BoolProf \to \Rel$ takes preorders to their object sets, and Boolean profunctors $R : X \proto{} Y$ to the relation $R\inv (1) \subset X \times Y$.
\end{defn}

\begin{rmk}[Boolean profunctors are `functorial relations']
    Given preorder $X$ and $Y$, a Boolean profunctor $R$ between them may be completely encoded by its underlying relation $\URel R \subset X \times Y$; in other words, the functor $\URel$ is faithful. However, the functor is not full: the relations $R \subset X \times Y$ in the image of $\URel$ are exactly the `functorial relations', namely those for which $\{a \to a', R(a',b)\} \imp R(a,b)$ and $\{b' \to b, R(a,b')\} \imp R(a,b)$.
\end{rmk}

\begin{rmk}[Boolean profunctors of discrete preorders and (co)functionality] \label{rmk:cofunctionality} A boolean profunctor $R : X \proto {} Y$ between discrete preorders is a relation between the corresponding sets in the ordinary sense. In this case, `(co)functorial profunctors' are usually called `entire (co)functional relations'---keeping `entireness' implicit, we will simply speak of (co)functional relations.
\end{rmk}

\begin{term}[Non-emptiness] We say a Boolean profunctor $R : X \proto{} Y$ is \textbf{non-empty} if there is a choice of $a \in X$, $b \in Y$ for which $R(a,b)$ holds.
\end{term}

\begin{term}[Restrictions of profunctors] \label{rmk:restricting-profunctors} Given a boolean profunctors $R : X \proto{} Y$ and subpreorders $U \into X$, $V \into Y$, then $R$ restricts to boolean profunctor between $U$ and $V$ which, abusing notation, we denote again by $R : U \proto{} V$.
\end{term}

\nid A final helpful notion for the definition of 1-truss bordisms is that of `bimonotonicity' of relations: bimotone relations of preoders do not transpose any two order-related elements.

\begin{defn}[Transpositions and bimonotonicity] Consider a relation $R \subset X \times Y$ between preorders $X$ and $Y$ (that is, $R$ is an ordinary relation between the object sets of $X$ and $Y$). A pair consisting of a non-identity arrow $x \to x'$ in $X$ and a non-identity arrow $y \to y'$ in $Y$ is called a \textbf{transposition of $R$} if both $R(x,y')$ and $R(x',y)$ are true. If the relation $R$ has no transpositions it is said to be \textbf{bimonotone}.
\end{defn}

\begin{eg}[Transpositions and bimonotonicity] In \autoref{fig:transpositions-in-relations} we depicted relations $R, S \subset X \times Y$ between posets $X \iso Y \iso [1]$. The relation $R$ has a transposition while $S$ has no transpositions and thus is bimonotone.
\begin{figure}[h!]
    \centering
    \def\svgwidth{1\columnwidth}
    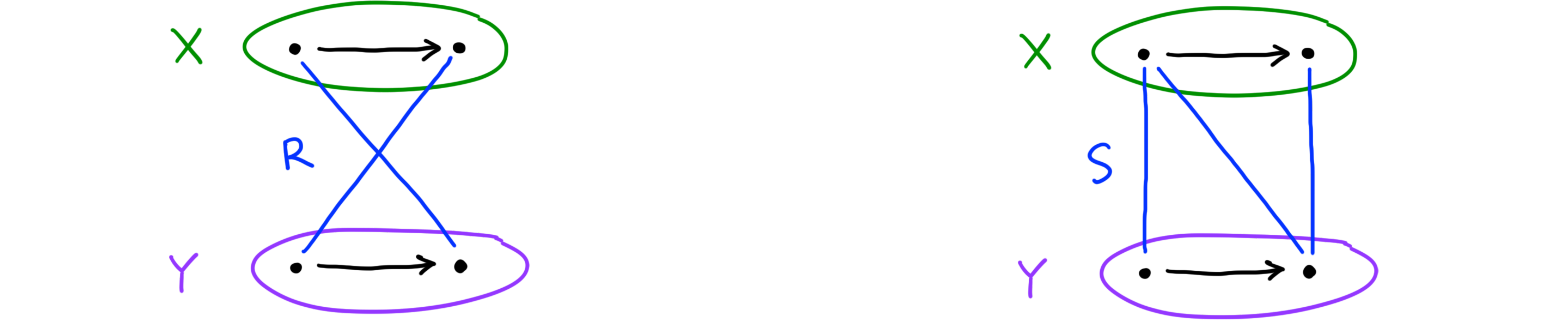

    \caption[Transpositions and bimonotonicity of relations.]{Transpositions and bimonotonicity of relations.}
    \label{fig:transpositions-in-relations}
\end{figure}
\end{eg}

\subsubsecunnum{1-Truss bordisms as bimonotone bifunctional profunctors} We finally turn to a definition of 1-truss bordisms. Recall that singular (resp.\ regular) objects in trusses $T$ yield discrete subposets $(\sing~T,\eleq)$ (resp.\ $(\reg~T, \eleq)$) of the face order $(T,\eleq)$ (see \autoref{rmk:reg-sing-subposets}).

\begin{defn}[1-truss bordisms of 1-trusses] \label{defn:1-truss-bord} A \textbf{1-truss bordism} $R : T \proto{} S$ of 1-trusses $T$ and $S$ is a non-empty boolean profunctor $R : (T,\eleq) \proto{} (S,\eleq)$ of the face orders of $T$ and $S$ subject to the following conditions.
    \begin{enumerate}
        \item \textit{Bifunctionality}: The restriction $R : (\sing~T,\eleq) \proto{} (\sing~S,\eleq)$ is functional, and the restriction $R : (\reg~T,\eleq) \proto{} (\reg~S,\eleq)$ is cofunctional. \vspace{4pt}

        \item \textit{Bimonotonicity}: The relation $\URel R$ of frame orders $(T,\fleq)$ and $(S,\fleq)$ is bimonotone. \qedhere
    \end{enumerate}
\end{defn}

\begin{rmk}[Unwinding bifunctionality] \label{rmk:unwinding-birepresentatbility} The bifunctionality condition can be made more explicit: it requires that for each singular object $a \in \sing~T$ there is a unique singular object $\singbif_R(a) \in \sing~T$ such that $R(a,\singbif_R(a))$, and for each regular object $d \in \reg~T$ there is a unique regular object $\regbif^R(d) \in \reg~T$ such that $R(\regbif^R(d),d)$. The resulting function $\singbif_R: (\sing~T,\eleq) \to (\sing~S,\eleq)$ is called the \textbf{singular function} of $R$. Similarly, the function $\regbif^R: (\reg~S,\eleq) \to (\reg~T,\eleq)$ is called the \textbf{regular function} of $R$.
\end{rmk}

\begin{figure}[h!]
    \centering
    \def\svgwidth{1\columnwidth}
    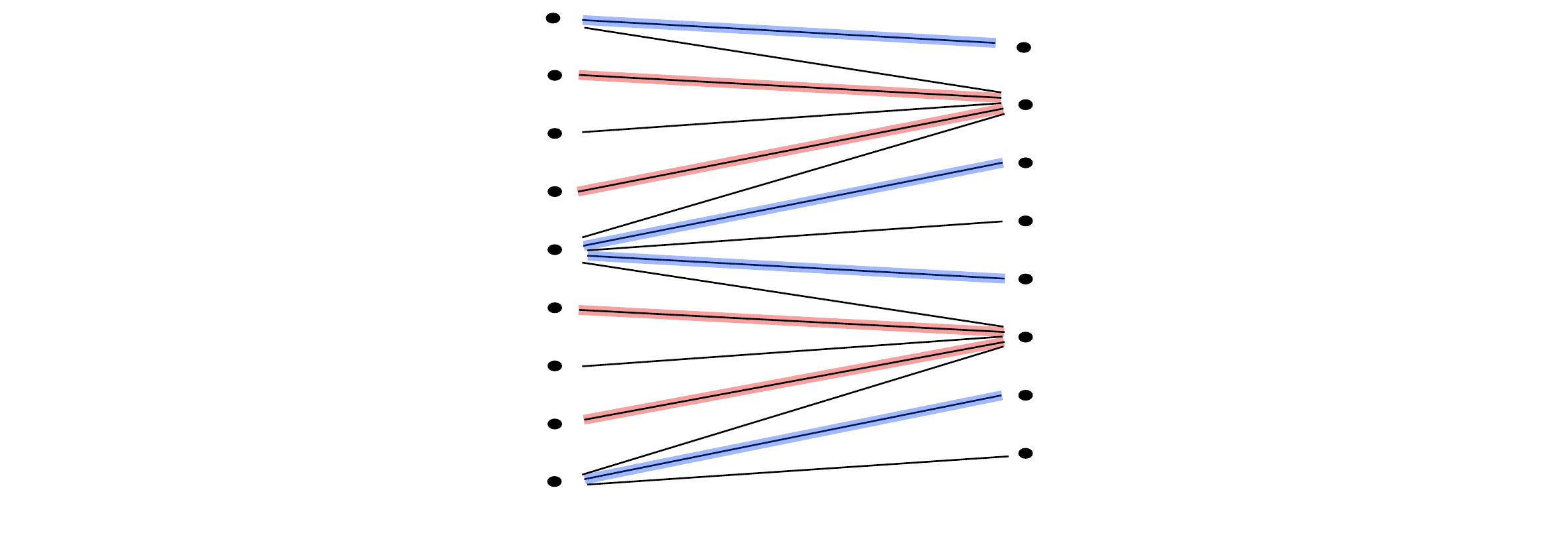

    \caption[A 1-truss bordism.]{A 1-truss bordism $R : T \proto{} S$.}
    \label{fig:1-truss-bordism}
\end{figure}

\begin{eg}[A 1-truss bordism]
    In \autoref{fig:1-truss-bordism} we depict a 1-truss bordism $R: T \proto{} S$. The domain 1-truss $T$ is drawn on the left, the codomain 1-truss $S$ is drawn on the right, and elements of the underlying relation of $R$ are indicated by black lines (the `edges' of $R$) between objects of $T$ and $S$. The fact that our choice of $R$ is bimonotone, i.e.\ has no transpositions, is witnessed by there being no crossings among its edges. The fact that $R$ is bifunctional is witnessed by the left-to-right function highlighted in red (providing the singular function $\singbif_R: \sing(T) \to \sing(S)$) and the right-to-left function highlighted in blue (providing the regular function $\regbif^R: \reg(S) \to \reg(T)$).
\end{eg}

\begin{rmk}[Relation to 1-truss maps] The relationship of 1-truss bordisms to 1-truss maps (as introduced in \autoref{defn:1-truss-map}) is analogous to the relationship of profunctors and functors. In particular, while every 1-truss map descends to a functor of underlying face orders, every 1-truss bordism will descend to a profunctor of underlying face orders.
\end{rmk}

\subsubsecunnum{Composition and dualization of 1-truss bordisms} The fact that 1-truss bordisms organize into a category follows from the following observation.

\begin{obs}[Bimonotonicity and bifunctionality compose] The properties of `bimonotonicity' and `bifunctionality' in \autoref{defn:1-truss-bord} are preserved under composition of Boolean profunctors. Indeed, given 1-truss bordisms, $R : T \proto{} T'$ and $R' : T' \proto{} T''$ then the Boolean profunctor $R' \circ R : T \proto{} T''$ is again bifunctional (its singular function is $\singbif_{R'} \circ \singbif_R$, and its regular function is $\regbif^{R'} \circ \regbif^R$), as well as bimonotone (any transposition of $R' \circ R$ with respect to the frame orders of $T$ and $T''$ would imply a transposition in at least on of $R$ or $R'$). Therefore, the Boolean profunctor $R' \circ R : T \proto{} T''$ in fact defines another 1-truss bordism, which we can take to be the composite of $R$ and $R'$.
\end{obs}

\nid Note further that `identity 1-truss bordisms' are given by ($\Bool$ enriched) hom functors $\Hom_{(T,\eleq)}(-,-)$, for any 1-truss $(T,\eleq,\dim,\fleq)$. Thus, 1-truss bordisms are morphisms of the following category.

\begin{notn}[The category of 1-truss bordisms] \label{notn:category-ebun}
    The \textbf{category of 1-truss bordisms} (whose objects are 1-trusses and whose morphisms are 1-truss bordisms) will be denoted by $\ttr 1$. The wide subcategories containing only open respectively closed 1-trusses will be denoted by $\ottr 1$ respectively by $\cttr 1$.
\end{notn}

\begin{obs}[The terminal and initial 1-trusses] \label{rmk:ebun-initial-terminal} The terminal object of $\ttr 1$ is the trivial closed 1-truss $\CTT_0$. The unique bordism $R : T \proto{} \CTT_0$ can be defined by setting $R(a,0)$ to hold for all $a \in T$.

The initial object of $\ttr 1$ is the trivial open 1-truss $\OTT_0$. The unique bordism $R : \OTT_0 \proto{} T$ can be defined by setting $R(0,a)$ to hold for all $a \in T$.
\end{obs}

\begin{rmk}[Invertible 1-truss bordisms are unique] \label{rmk:1-trusses-skeletally-bordism} Given two 1-trusses, any invertible 1-truss bordism between them is unique if one exists. There is thus no harm in working `skeletally', i.e.\ working with 1-trusses up to invertible bordisms. Further the classes of 1-trusses up to invertible 1-truss bordisms are exactly the classes of 1-trusses up to balanced 1-truss isomorphisms (see \autoref{rmk:1-trusses-skeletally}).
\end{rmk}

\nid Recall the dualization functor $\dagger$ on 1-trusses (and their maps) from \autoref{constr:1-truss-duals}. A similar dualization construction applies to the category of 1-truss bordism as follows.

\begin{constr}[Dualization of 1-truss bordisms] \label{constr:truss-trans-dualization} Let $R: T \proto{} S$ be a 1-truss bordism. We define its \textbf{dual} $R^\dagger : S^\dagger \proto{} T^\dagger$ to be the 1-truss bordism determined by transposing the underlying relation of $R$; that is, $R^\dagger$ is given by the Boolean profunctor
    \begin{equation}
        (S,\eleq\op)\op \times (T,\eleq\op) \iso (T,\eleq)\op \times (S,\eleq) \xto {R} \Bool
    \end{equation}
    One verifies that $R^\dagger$ satisfies both bimonotonicity and bifunctionality. Dualization gives rise to an involutive isomorphism of categories
    \[ \dagger : \ttr 1 \iso (\ttr 1)\op \]
As before this restricts to an equivalence $\dagger : \ottr 1 \iso (\cttr 1)\op$.
\end{constr}

\subsubsecunnum{Properties of 1-truss bordisms} We discuss several immediate properties of 1-truss bordisms.

\begin{lem}[Dimension monotonicity] \label{lem:truss-trans-dimensions} Given a 1-truss bordism $R: T \proto{} S$ and elements $a \in T$, $b \in S$ such that $R(a,b)$ holds, then $\dim(b) \leq \dim(a)$.
\end{lem}

\begin{proof} We show that if $R(a,b)$ holds and $\dim(a) = 0$, then we must have $\dim(b) = 0$. In other words, $R$ never relates singular elements of $T$ to regular elements in $S$.  Arguing by contradiction, assume $R(a,b)$ holds for singular $a$ and regular $b$. Since frame orders are total, and since $b$ is regular, we either have $b \fles \singbif_R(a)$ or $\singbif_R(a) \fles b$. Assume the former case holds (the latter case is symmetric). By bimonotonicity we must have $\regbif^R(b) \fles a$. Then there exists a face order arrow $a-1 \eles a$ in $(T,\eleq)$ (`$a-1$' denotes the predecessor of $a$ in the total order $(T,\fleq)$). By profunctoriality of $R$ this implies that $R(a-1,\singbif_R(a))$ holds. But, since we assumed $R(a,b)$ and since $b \fles \singbif_R(a)$, this contradicts bimonotonicity.
\end{proof}

Recall our definition and notation for endpoints of 1-trusses (see \autoref{term:endpoints}).

\begin{lem}[Endpoint preservation] \label{lem:truss-trans-endpoints} Any 1-truss bordism $R: T \proto{} S$ must relate lower (resp.\ upper) endpoints of $T$ to lower (resp.\ upper) endpoints of $S$.
\end{lem}

\begin{proof} We argue in the case of lower endpoints $\ept_- T$, $\ept_- S$ (the case of upper endpoints is symmetric). First, assume $\ept_- S$ is regular. Set $a := \regbif^R(\ept_- S)$. We must have $a = \ept_- T$ which shows the claim: indeed, if $a \neq \ept_-T$ then $R(a-1,\singbif_R(a-1))$ would hold, which would contradict bimonotonicity (since $\ept_- S \fles \singbif_R(a-1)$).

    Next, assume $\ept_- S$ is singular. If $\ept_- T$ is singular then a dual argument to the first case shows $\singbif_R(\ept_- T) = \ept_- S$ as claimed. If $\ept_- T$ is regular then we argue as follows. If both $T$ and $S$ contain a single object, then $R(\ept_-T,\ept_-S)$ holds since $R$ is assumed non-empty. Otherwise, we can assume that either $T$ or $S$ contain at least two objects. We argue in the former case (the argument for the latter case is similar). Since $T$ has at least two elements, its lower endpoint has a successor $a = \ept_-T + 1$. Set $b := \singbif_R(a)$. If $b = \ept_-S$, then $R(\ept_-T,\ept_-S)$ holds as claimed, since $\ept_-T \eles a$ and since $R$ is profunctorial.  On the other hand, if $b \fgre \ept_-S$, then we must have $\regbif^R(\ept_-S + 1) = \ept_- T$ by bimonotonicity. Since $\ept_-S \egre \ept_-S + 1$, profunctoriality of $R$ implies again that $R(\ept_-T,\ept_-S)$ holds as claimed.
\end{proof}


\nid As we've seen every 1-truss bordism $R : T \proto{} S$ determines a `singular function'  $\singbif_R : (\sing(T),\fleq) \to (\sing(S),\fleq)$ and a `regular function' $\regbif^R : (\reg(S),\fleq) \to (\reg(T),\fleq)$. By the previous \autoref{lem:truss-trans-endpoints}, it follows that $\singbif_R$ `preserves singular endpoints' in the following sense: whenever a (lower resp.\ upper) endpoint $\ept_\pm T$ is a singular in $T$, then the image $\singbif_R(\ept_\pm) = \ept_\pm S$ is a (lower resp.\ upper) endpoint in $S$, which is necessarily singular itself. Dually, $\regbif^R$ `preserves regular endpoints' in that, whenever a (lower resp.\ upper) endpoint $\ept_\pm S$ is regular in $S$, then the image $\regbif^R(\ept_\pm S) = \ept_\pm T$ is a (lower resp.\ upper) endpoint in $T$, which is necessarily a regular itself. In fact, frame order preserving maps of singular resp.\  regular elements in trusses, which preserve singular resp.\ regular endpoints in this sense, determine 1-truss bordisms as the next result records.

\begin{lem}[Singular and regular determined bordisms] \label{lem:truss-trans-representations} Let $T$ and $S$ be 1-trusses.
\begin{list}{}{
\setlength\leftmargin{18pt}
}
\item \textsc{Singular determined}: Given a function $f: (\sing(T),\fleq) \to (\sing(S),\fleq)$ that preserves singular endpoints, there is a unique 1-truss bordism $R: T \proto{} S$ with singular function $\singbif_R = f$.
\item \textsc{Regular determined:} Given a function $f: (\reg(S),\fleq) \to (\reg(T),\fleq)$ that preserves regular endpoints, there is a unique 1-truss bordism $R: T \proto{} S$ with regular function $\regbif^R = f$.
\end{list}
\end{lem}
\begin{proof}
For the first case, define $R(a,b)$ to hold if and only if either (1) the object $a$ is singular and $b = \singbif_R(a)$, or (2) the object $a$ is regular and both $b \leq \singbif_R(a+1)$ (whenever $a+1 \in T$) and $b \geq \singbif_R(a-1)$ (whenever $a-1 \in T$).  For the second case, define $R(a,b)$ to hold if and only if either (1) the object $b$ is regular and $a = \regbif^R(b)$, or (2) the object $b$ is singular and both $a \leq \regbif^R(b+1)$ (when $b+1 \in S$) and $a \geq \regbif^R(b-1)$ (when $b-1 \in S$).
\end{proof}

We say a relation (or Boolean profunctor) $R : T \proto{} S$ `fully matches objects' if for each $a \in T$ there exists $a' \in S$ with an $R(a,a')$ and for each $b \in S$ there exists $b' \in T$ with $R(b',b)$. The description of 1-truss bordisms in the previous lemma (and its proof) then has the following two corollaries.

\begin{cor}[Truss bordisms fully match objects] \label{lem:truss-trans-right-left-lifts} Any 1-truss bordism $R: T \proto{} S$ fully matches objects. \qed
\end{cor}

\begin{cor}[Dependence of full matching and functionality] \label{cor:full-matching} Consider 1-trusses $T$ and $S$, and a boolean profunctor $R : (T,\eleq) \proto{} (S,\eleq)$ such that $\URel R \subset (T,\fleq) \times (S,\fleq)$ is bimonotone and fully matches objects. If \emph{either} $R(u,v)$ is functional on singular objects \emph{or} cofunctional on regular objects, then $R$ is a 1-truss bordism.
\end{cor}
\begin{proof} The assumptions on $R$ imply that $R$ satisfies exactly the conditions in one of the two cases of \autoref{lem:truss-trans-representations}.
\end{proof}

\subsection{1-Truss bundles} \label{ssec:1-truss-bun}

We now define 1-truss bundles: these are bundles over posets whose fibers are 1-trusses, and whose fiber transitions are described by 1-truss bordisms.

\subsubsecunnum{1-Truss bundles as collections of 1-truss bordisms}

\begin{defn}[1-Truss bundle] \label{defn:1-truss-bundle}
    Consider posets $(B,\to)$ and $(T,\leq)$, equipped with a functor $\dim : (T,\eleq) \to [1]\op$ and a second order $(T,\fleq)$. A \textbf{1-truss bundle} is a map of diposets $p: (T,\eleq,\fleq) \to (B,\to,=)$ (where `$=$' is the discrete order) satisfying the following.
    \begin{enumerate}
        \item For every object $x$ in $B$, the datum $(T,\fleq,\dim,\fles)$ restricts on the fiber over $x$ to a 1-truss $p\inv(x) \equiv (p\inv(x),\fleq,\dim,\fles)$.
        \item For every arrow $x \to y$ in $B$, the fiber $p\inv(x \to y)$ is a 1-truss bordism, in the sense that the relation
            \begin{center}
                $a \sim b$ iff ($a \eleq b$ in $(T,\eleq)$ s.t. $p(a \eleq b) = x \to y$)
            \end{center}
is the underlying relation of a 1-truss bordism $p\inv(x) \proto{} p\inv(y)$.
    \end{enumerate}
We call $(T,\eleq)$ the \textbf{total poset} and $(B,\leq)$ the \textbf{base poset} of the 1-truss bundle $p$. We call $\eleq$ the \textbf{face order}, and $\fleq$ the \textbf{frame order} of the bundle.
\end{defn}

\begin{notn}[Face, dimension, and frame structure] When working with 1-truss bundles, we will usually keep face orders, dimension maps, frame orders as well as base poset orders implicit; that is, we will denote 1-truss bundles simply by maps $p : T \to B$. When needed, we will use either arrows `$\to$' or the relation symbol `$\eleq$' for face orders; `$\dim$' for dimension maps $\dim : (T,\eleq\op) \to [1]$; the relation symbol `$\fleq$' for frame orders; and arrows `$\to$' for base poset orders.
\end{notn}

\begin{rmk}[Frame orders are trivial across fibers] Observe that, given a 1-truss bundle $p : T \to B$, then objects $a, b \in T$ are related in the frame order of $p$ if and only if $a$ and $b$ live in the same fiber of $p$: indeed, no objects between different fibers can be frame order related since $p : (T,\fleq) \to (B,=)$ is a poset map, and frame orders are total on each fiber by the condition that fibers are 1-trusses.
\end{rmk}

\begin{eg}[1-Truss bundle] We illustrate a 1-truss bundle $p : T \to B$ in \autoref{fig:1-truss-bundle} (we indicate total frame orders on each fiber, obtained by restricting the frame order $\fleq$ to that fiber, by a green coordinate axis). Note that we chose all such axis to point in the same (`upwards') direction. Flipping this direction of frame orders on all fibers (i.e.\ letting axis point downwards) would lead to another valid example.
\begin{figure}[h!]
    \centering
    \def\svgwidth{1\columnwidth}
    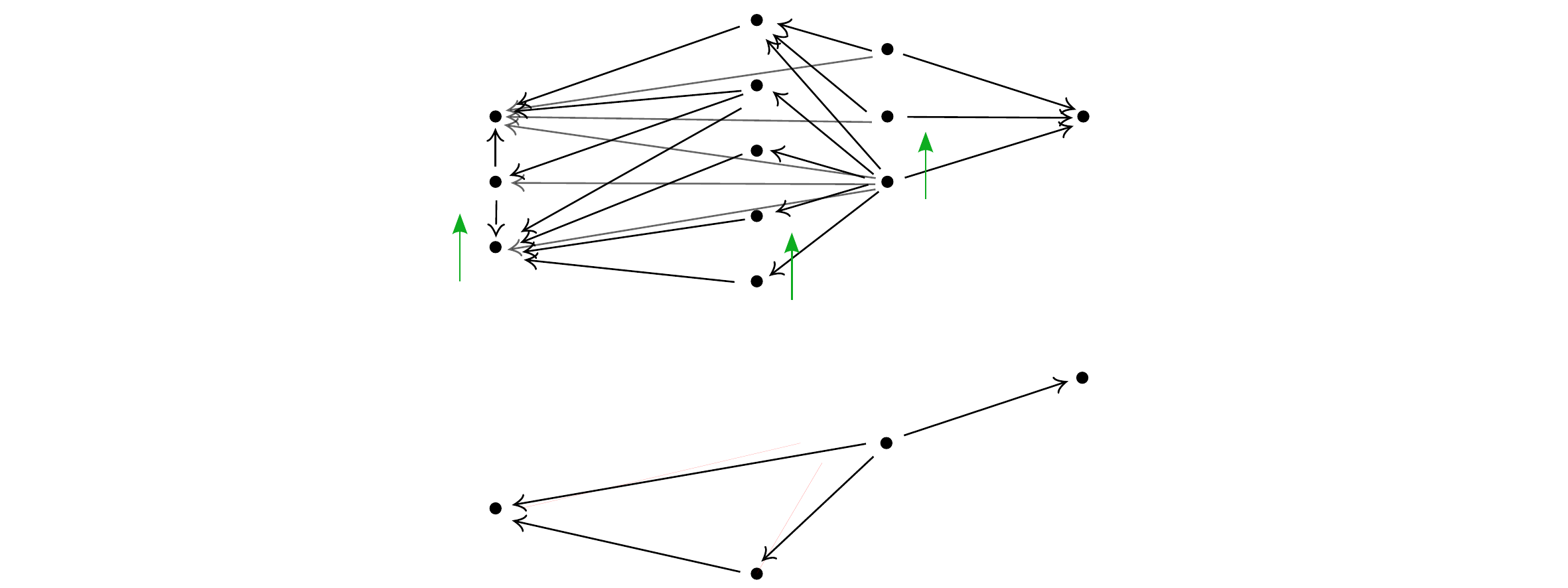

    \caption[A 1-truss bundle.]{A 1-truss bundle $p : T \to B$.}
    \label{fig:1-truss-bundle}
\end{figure}
\end{eg}

\begin{rmk}[$\lZ_2$ action on frame orders] \label{rmk:z2-ambiguity-1-truss-bundles} As observed in the preceding example, $\lZ_2$ acts on the frame order of 1-truss bundles by `flipping' the direction of total the frame order on every fiber. Most examples in this chapter do not change in nature under this action. Therefore, we usually depict frame orders of bundles only `up to' this $\lZ_2$ action; namely, we will parallely align all fibers of the bundle (as in \autoref{fig:1-truss-bundle}) and assume all fiber frame orders run in the `same direction', but we will usually not fix this direction.
\end{rmk}

\begin{term}[Open and closed bundles]  A 1-truss bundle in which all fibers are open (resp.\ closed) 1-trusses is called an `open' (resp.\ `closed') 1-truss bundle.
\end{term}

\begin{term}[Singular and regular objects in bundles] \label{term:sing_reg_objects_in_bundles} Given a 1-truss bundle $p : T \to B$ we call an element $a \in T$ `singular' if $\dim(a) = 0$ and `regular' if $\dim(a) = 1$. We denote by $\sing(T)$ respectively by $\reg(T)$ the full subposets of $(T,\eleq)$ containing all singular respectively all regular objects.
\end{term}

\begin{obs}[Underlying discrete (op)fibrations] \label{rmk:1-truss-bundle-fibrations} Recall the standard notion of discrete opfibrations: a functor $F : \iC \to \iD$ is called a discrete opfibration if for every object $c$ in $\iC$ and every morphism of the form $g : F(c) \to d$ in $\iD$ there is a unique morphisms $h : c \to c'$ such that $F(h) = g$. Dually, one defines a discrete fibration $F : \iC \to \iD$ to be a discrete opfibration $F\op : \iC\op \to \iD\op$ after dualization. From the bifunctionality of 1-truss bordism it follows that 1-truss bundles $p : T \to B$ restrict on singular objects to a discrete opfibration $p : \sing(T) \to B$ and, dually, on regular objects to a discrete fibration $p : \reg(T) \to B$.
\end{obs}

\nid 1-Truss bundles further have the following lifting property.

\begin{obs}[Lifts in 1-truss bundles] \label{rmk:1-truss-bundle-lifts} Consider a 1-truss bundle $p : T \to B$. For an object $a$ in $T$ and an arrow $x \to y$ in $B$ with $p(a) = x$ there exists at least one object $a'$ in $T$ with an arrow $a \to a'$ and such that $p(a \to a') = (x \to y)$. Conversely, for $b$ in $T$ and $x \to y$ in $B$ with $p(b) = y$, there exists at least one $b'$ in $T$ with $b' \to b$ and $p(b' \to b) = (x \to y)$. Both claims follow from \autoref{lem:truss-trans-right-left-lifts}.
\end{obs}

\subsubsecunnum{Maps of 1-truss bundles} Using the notion of maps of 1-trusses `fiberwise', we obtained the following notion of maps of 1-truss bundles.

\begin{defn}[Maps of 1-truss bundles] \label{defn:1-truss-bundle-map} For 1-truss bundles $p : T \to B$ and $q : S \to C$, a \textbf{1-truss bundle map} $F : p \to q$ is a diposet map $F : T \to S$ such that $q \circ F$ factors through $p$ by a (necessarily unique) map $G$; that is, the following square commutes:
    \begin{equation}
\begin{tikzcd}
T \arrow[r, "F"] \arrow[d, "p"'] & S \arrow[d, "q"] \\
B \arrow[r, "G"]                 & C
\end{tikzcd}
    \end{equation}
In the case where $G = \id_B$ we say the 1-truss bundle map $F$ is \textbf{base preserving}.
\end{defn}

\nid Note that a 1-truss bundle map $F : p \to q$ restricts on each fiber of $p$ to a 1-truss map: that is, for each $x \in B$, the bundle map $F$ restricts to a 1-truss map $F : p\inv(x) \to q\inv G(x)$.

\begin{term}[Singular, regular and balanced 1-truss bundle maps] \label{defn:sing-and-reg-1-truss-bun-map} Consider a 1-truss bundle map $F : p \to q$. If $F$ maps $\sing(T)$ to $\sing(S)$ we call $F$ a `singular' bundle map; similarly, if it maps $\reg(T)$ to $\reg(S)$ we call it a `regular' bundle map. If $F$ is both singular and regular, than we call it `balanced'.
\end{term}

\nid Note that a 1-truss bundle map $F : p \to q$ is singular if and only if $\dim F(x) \leq \dim x$ for all $x \in T$, and regular if and only if $\dim F(x) \geq \dim x$ for all $x \in T$.

\begin{notn}[Category of 1-truss bundles and their maps] The category of 1-truss bundles and their maps is denoted by $\trussbun 1$. The full subcategories of open resp.\ closed trusses are denoted by $\otrussbun 1$ resp.\ $\ctrussbun 1$. As before, restricting either of those categories to singular, resp.\ regular, resp.\ balanced maps will be denoted by superscripts `$\mathsf{s}$', resp.\ `$\mathsf{r}$', resp.\ `$\mathsf{rs}$'.
\end{notn}

\begin{notn}[Categories of 1-truss bundles over a fixed base] The subcategory of the category of 1-truss bundles $\trussbun 1$ of bundles over a fixed poset $B$ and base preserving maps will be denoted by $\trussbun 1 (B)$. Decorations for bundle fiber types and map types apply as before.
\end{notn}

\begin{defn}[Restrictions of 1-truss bundles] Given a truss bundle $p : T \to B$ and a subposet $C \into B$, the \textbf{restriction} of $p$ to $C$ is the 1-truss bundle $\rest p C : \rest T C \to C$ obtained by restricting $p$ to $\rest T C = p\inv(C) \into T$ (and all data on $T$ is restricted accordingly). This induces the restriction functor $\rest {-} C : \trussbun 1 (B) \to \trussbun 1 (C)$
\end{defn}

\subsubsecunnum{Classification and totalization for 1-truss bundles} 1-Truss bundles admit bundle classification and totalization constructions; their classifying category is the category of 1-truss bordisms $\ttr 1$. For instance, in the simple case of the base poset being the 1-simplex $[1]$, observe that 1-truss bundles $p : T \to [1]$ are in correspondence with 1-truss bordisms and thus with functors $\iF : [1] \to \ttr 1$. We now describe a more general correspondence of 1-truss bundles $p : T \to B$ over a general base poset $B$ with functors $\iF : B \to \ttr 1$ from $B$ into the category of 1-truss bordisms.

\begin{constr}[Classifying functors] \label{constr:classification-1}
    Let $p : T \to B$ be a 1-truss bundle over a poset $B$. We construct the \textbf{classifying functor} $\fcl {} p : B \to \ttr 1$ as follows. For each object $x$, the image $\fcl {} p (x)$ is the 1-truss $p\inv(x)$ and for each arrow $x \to y$, the image $\fcl {} p (x \to y)$ is the 1-truss bordism $p\inv(x \to y)$ (as given in \autoref{defn:1-truss-bundle}). The functoriality of $\fcl {} p$ follows from $(T,\eleq)$ being a poset and by the definition of Boolean profunctor composition.
\end{constr}

Conversely, one constructs `total bundles' as follows.

\begin{constr}[Totalizing bundles] \label{constr:totalization-1} Given a poset $B$,  consider a functor $\iF: B \to \ttr 1$. We construct the \textbf{total 1-truss bundle} $\fp {} \iF: \Tot{}\iF \to B$ as follows. The `total poset' $(\Tot{}\iF,\eleq)$ has as elements the pairs $(x \in B, a \in \iF(x))$. Morphisms $(x,a) \eleq (y,b)$ are given whenever the 1-truss bordism $\iF(x \to y)$ from the 1-truss $\iF(x)$ to the 1-truss $\iF(y)$ relates the elements $a$ and $b$. The frame order $(\Tot{}\iF,\fleq)$ is defined to relate $(x,a) \fleq (x,b)$ whenever $a \fleq b$ in $\iF(x)$. The dimension functor $\dim : (T,\eleq\op) \to [1]$ is equally defined fiberwise; the fact that this extends to a functor on all of $(T,\eleq\op)$ follows since truss bordisms are `dimension monotonic' by \autoref{lem:truss-trans-dimensions}.
\end{constr}

    In order to promote the preceding constructions to functors, we introduce the following notion of `bundle concordance'.

    \begin{defn}[Concordance of 1-truss bundles] For a poset $B$, a \textbf{1-truss bundle concordance} $u : p \Rightarrow q$ between 1-truss bundles $p : T \to B$ and $q : S \to B$ is a 1-truss bundle $u : U \to B \times [1]$ such that $\rest r {B \times \{0\}} = p$ and $\rest r {B \times \{1\}} = q$.
\end{defn}

\nid If the base $B$ is trivial, a 1-truss bundle concordance is simply a 1-truss bordism.

\begin{rmk}[Invertible concordance are unique] \label{rmk:1-trusses-skeletally-concordance} Generalizing \autoref{rmk:1-trusses-skeletally-bordism}, observe that, given two 1-truss bundles, any invertible concordance between them is unique if one exists. There is thus no harm in working `skeletally', i.e.\ working with 1-truss bundles up to invertible concordances. Further, the classes of 1-truss bundles up to invertible concordance are exactly the classes of 1-truss bundles up to base preserving balanced 1-truss bundle isomorphisms.
\end{rmk}

\nid Note that, given a 1-truss bundle concordance $u : p \Rightarrow q$, its classifying functor $\fcl {} u$ is a functor $B \times [1] \to \ttr 1$ restricting on $B \times \{0\}$ to $\fcl {} p$ and on $B \times \{1\}$ to $\fcl {} q$. Equivalently, $\fcl {} u$ is thus a natural transformation of classifying functors $\fcl {} p \Rightarrow \fcl {} q$. We will refer to $\fcl {} u$ as the `classifying natural transformation' of $u$.

\begin{defn}[Categories of 1-truss bundle concordances] For a poset $B$, define the \textbf{category of concordances of 1-truss bundles over $B$}, denoted by $\trussconc 1(B)$, to have 1-truss bundles over $B$ as objects, and 1-truss bundle concordances as morphisms. Composition of two truss concordances $u : p \Rightarrow q$ and $v : q \Rightarrow r$ is determined by the condition that $\fcl {} {u \circ v} = \fcl {} v \circ \fcl {} u$.\footnote{Explicitly, working skeletally (that is, identifying 1-truss bundles up to invertible concordances, see \autoref{rmk:1-trusses-skeletally-concordance}), this can be defined by setting $v \circ u = \fp {} {\fcl {} v \circ \fcl {} u}$.}
\end{defn}

\nid With these definitions in place, we obtain the following equivalence of categories.

\begin{constr}[Classification and totalization functors] \label{constr:class-tot-1} Given a poset, we define the \textbf{classification functor}
\[
    \fcl{}{-} : \trussconc 1(B) \to \Fun(B,\ttr 1)
\]
to map 1-truss bundles $p : T \to B$ to their classifying functor $\fcl {} p$, and concordances $u : p \Rightarrow q$ to their classifying natural transformations $\fcl{}u : \fcl{}p \Rightarrow \fcl{}q$. Conversely, we define the \textbf{totalization functor}
\[
    \fp{}{-} : \Fun(B,\ttr 1) \to \trussconc 1(B)
\]
to take functors $\iF: B \to \ttr 1$ to their total 1-truss bundles $\fp {} \iF$, and natural transformations $\alpha : B \times [1] \to \ttr 1$ to their total 1-truss bundles $\fp {} \alpha$.
\end{constr}

\begin{obs}[Classification and totalization are inverse] \label{obs:classification-vs-totalization-1} Classification and totalization functors provide an equivalence of categories.
\end{obs}

\begin{eg}[The composition of 1-truss bordisms as a truss bundle over the 2-simplex] \label{eg:compmorph}
Two composable 1-truss bordisms $R_1 : T_0 \to T_1$ and $R_2 : T_1 \to T_2$ together with their composite $R_2 \circ R_1 : T_0 \to T_2$ evidently define a functor $[2] \to \ttr 1$ from the poset $[2]$ to the category of 1-truss bordisms; by the above construction this situation can be equivalently considered as a 1-truss bundle over the 2-simplex. For example, in \autoref{fig:1-truss-bordisms-composition-as-a-bundle} we illustrate truss bordisms $R_1$ and $R_2$, along with the corresponding bundle $p : T \to [2]$ on the right.
\begin{figure}[h!]
    \centering
    \def\svgwidth{1\columnwidth}
    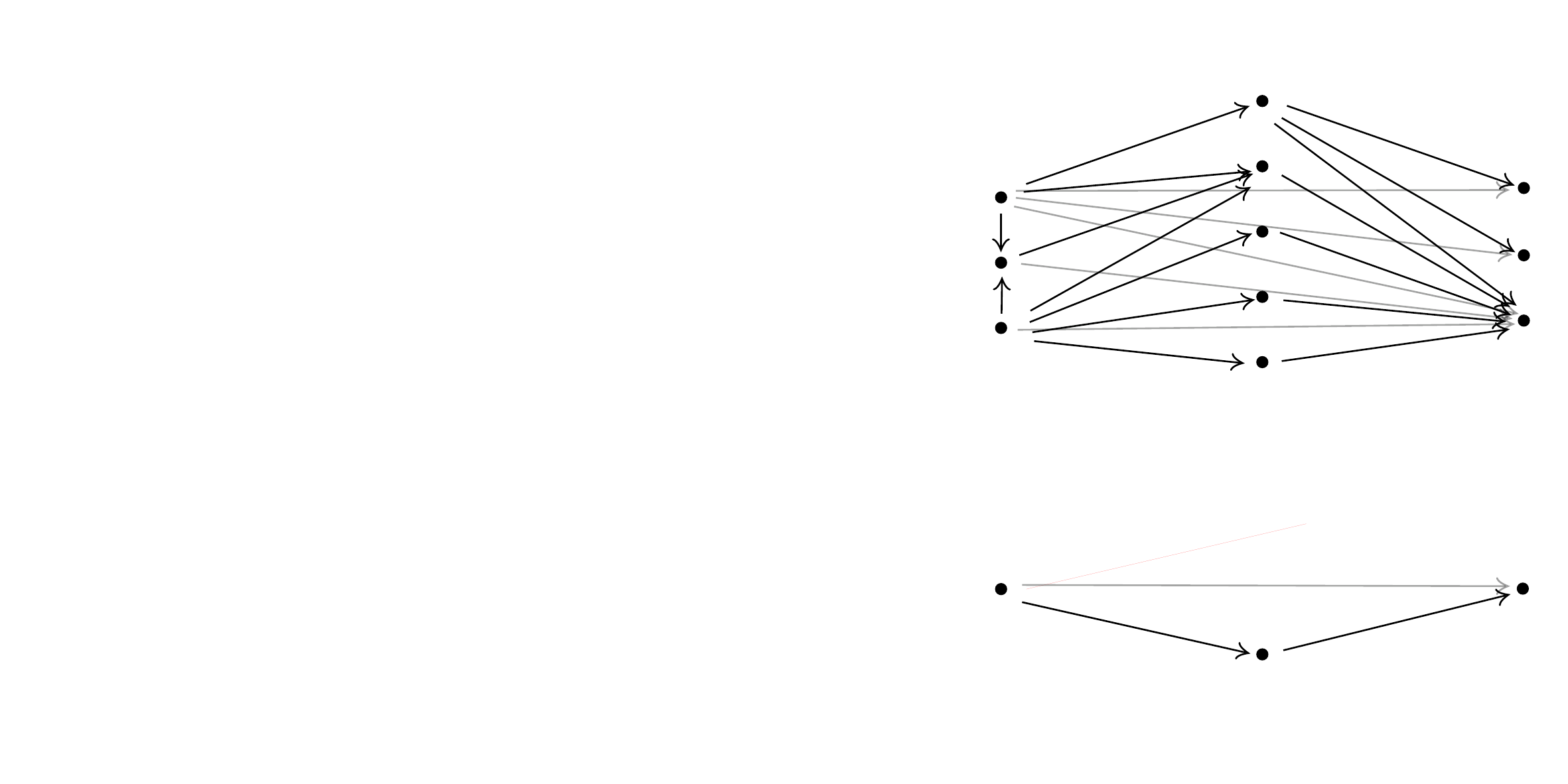

    \caption{1-Truss bordisms composition as a bundle over the 2-simplex.}
    \label{fig:1-truss-bordisms-composition-as-a-bundle}
\end{figure}
\end{eg}

\nid The classification and totalization constructions described here for 1-truss bundles, are of course analogous to many of the classical combinatorial instances of these constructions (see for instance \cite{benabou2000distributors}, \cite{street2001powerful}).

\subsubsecunnum{Pullback, dualization, and suspension of 1-truss bundles}

With notions of classification and totalization at hand, we now describe several further constructions on 1-truss bundles. Firstly, 1-truss bundles may be pulled back along poset maps into their base as follows.

\begin{defn}[Pullback 1-truss bundles] \label{constr:1-truss-bundle-pullback} Let $p : T \to B$ be a 1-truss bundle with classifying map $\fcl {} p : B \to \ttr 1$. Consider a map of posets $F : C \to B$. We define the \textbf{pullback $F^*p : F^*T \to C$ of $p$ along $F$} to be the total 1-truss bundle $\fp {} {\fcl {} p F}$ of the composite functor $\fcl {} p \circ F : C \to \ttr 1$. The canonical 1-truss bundle map $F^*p \to p$, which is fiberwise given by isomorphisms $(F^*p)\inv y \iso p\inv F(y)$, will be denoted by $\Tot{}F$.
\end{defn}

\begin{notn}[Pullback of 1-truss bundles] \label{notn:pullback-of-1-truss-bundles} A pullback of 1-truss bundles will often be indicated by a `pullback square'
\begin{equation}
\begin{tikzcd}
S \arrow[dr, phantom, "\lrcorner" , very near start, color=black] \arrow[r, "\Tot{}F"] \arrow[d, "q"'] & T \arrow[d, "p"] \\
C \arrow[r, "F"]                       & B
\end{tikzcd}
\end{equation}
where the 1-truss bundle $q$ is the pullback of the 1-truss bundle $p$ along the poset map $F$. Note, while `1-truss bundles' and `poset maps' do not live in the same category, we may forget part of the 1-truss bundle structure and think of $p$ and $q$ as poset maps instead: this then makes the above a pullback in the category of posets and poset maps.
\end{notn}

\begin{rmk}[Pullbacks generalize restrictions] \label{rmk:1-truss-bundle-pullback-restriction}) In the special case where $F : C \into B$ is a subposet of $B$, the pullback bundle $F^*p$ equals the restriction  bundle $\rest p C$.
\end{rmk}

\begin{eg}[Pullback 1-truss bundles] We depict the canonical pullback map $(\Tot{}F, F) : F^*p \to p$ of a 1-truss bundle pullback in \autoref{fig:pullback-1-truss-bundle} (in fact, this is a restriction in the sense of \autoref{rmk:1-truss-bundle-pullback-restriction})).
\begin{figure}[h!]
    \centering
    \def\svgwidth{1\columnwidth}
    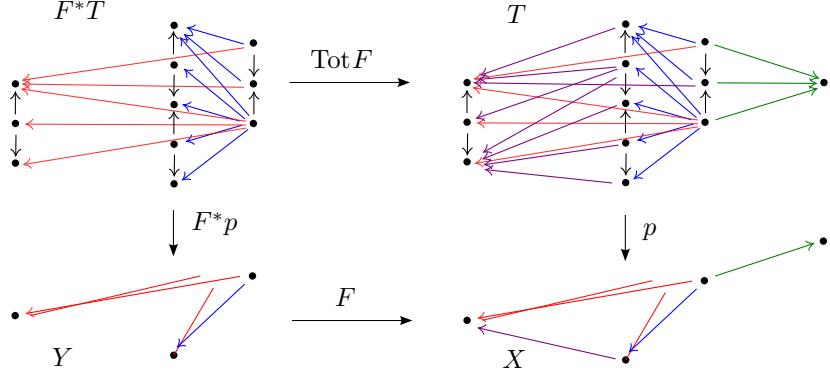

    \caption[The pullback of a 1-truss bundle along a base poset map.]{The pullback of a 1-truss bundle $p$ along a base poset map $F$.}
    \label{fig:pullback-1-truss-bundle}
\end{figure}
\end{eg}

Next, we generalize our dualization construction for 1-trusses to the case of 1-truss bundles; this simply `dualizes all fibers' of a given bundles.

\begin{defn}[Dual 1-truss bundles] \label{constr:1-truss-bundle-dualization} Given a 1-truss bundle $p : T \to B$ over $B$, its \textbf{dual bundle} $(p^\dagger : T^\dagger \to B\op)$ is the 1-truss bundle over $B\op$ whose total poset $(T,\eleq\op)$ is opposite to that of $p$, whose dimension functor is the composite $(T,\eleq\op) \xto {\dim\op} [1] \iso [1]\op$ (where $\dim$ is the dimension functor of $p$), and whose frame order $(T,\fleq)$ is equal to that of $p$.
\end{defn}

\begin{rmk}[Dualization via classifying functors] Using the classification and totalization equivalence, given a 1-truss bundle $p$, its dual $p^\dagger$ may equivalently be defined to have classifying functor
    \[
        (\fcl {} {p^\dagger} : B\op \to \ttr 1) = \big(B \xto {\fcl {} T} \ttr 1 \xto {\dagger} (\ttr 1)\op \big)\op .
    \]
    (The functor $\dagger : \ttr 1 \to (\ttr 1)\op$ was defined in \autoref{constr:truss-trans-dualization}.)
\end{rmk}

\begin{defn}[Dual 1-truss bundle maps] Given a 1-truss bundle map $F : (p : T \to B) \to (q : S \to C)$ one defines the \textbf{dual bundle map} $F^\dagger : p^\dagger \to q^\dagger$ by setting $F^\dagger = F\op : (T,\eleq\op) \to (S,\eleq\op)$.
\end{defn}

\nid Note that dualization maps open 1-truss bundles to closed 1-truss bundles and vice-versa, as well as singular bundle maps to regular ones and vice versa.

\begin{defn}[Dual 1-truss bundle concordances] Given a 1-truss bundle concordance $u : p \Rightarrow q$ one defines the \textbf{dual bundle concordance} $u^\dagger : p^\dagger \Rightarrow q^\dagger$ by setting $u^\dagger$ to be the dual bundle of $u$.
\end{defn}

\nid The preceding constructions now yield the following functors.

\begin{obs}[Dualization functors on 1-truss bundles] Generalizing the dualization functor on the category of 1-trusses, we obtain a dualization functor
    \[
        \dagger : \trussbun 1 \iso \trussbun 1.
    \]
    This recovers our earlier \autoref{constr:1-truss-duals} if we restrict to $\trussbun 1(\ast) \into \trussbun 1$.

    Similarly, generalizing the dualization functor on the category of 1-truss bordisms, we obtain, for any poset $B$, a dualization functor
    \[
        \dagger : \trussconc 1(B) \iso \trussconc 1(B)\op .
    \]
This recovers \autoref{constr:truss-trans-dualization} if we set $B = \ast$.
\end{obs}

As a final elementary construction, we now introduce `suspensions' of 1-truss bundles: this adds new initial and final elements to a given truss bundle. We start with the case of posets.

\begin{constr}[Suspension of posets] Let $X$ be a poset. Its \textbf{suspension} $\Sigma X$ is the poset obtained from $X$ by adjoining two elements $\top$ and $\bot$ and arrows $x \to \top$ and $\bot \to x$ for each $x \in X$.
\end{constr}

\begin{constr}[Suspension of 1-truss bundles] \label{constr:1-truss-bundle-suspension} Let $p : T \to B$ be a 1-truss bundle. The suspension bundle $\Sigma p : \Sigma T \to \Sigma B$ is the 1-truss bundle
    \begin{enumerate}
        \item whose base poset is the suspension $\Sigma B$ of $B$,
        \item whose total poset $(\Sigma T, \eleq)$ is the suspension $\Sigma (T,\eleq)$ of $(T,\eleq)$,
        \item whose dimension functor restricts on $T \into \Sigma T$ to the dimension functor of $p$ while mapping the new initial object $\bot$ in $\Sigma T$ to $0$ and terminal object $\top$ in $\Sigma T$ to $1$,
        \item whose frame order $(\Sigma T, \fleq)$ relates elements if and only if they are already related in $(T,\fleq)$. \qedhere
    \end{enumerate}
\end{constr}

\begin{rmk}[Suspension bundles via classifying functors] Given a 1-truss bundle $p$, its suspension bundle $\Sigma p$ is classified by the unique map $\fcl{}{\Sigma p} : \Sigma B \to \ttr 1$ which restricts on $B \into \Sigma B$ to $\fcl{}p$, maps $\bot$ to the initial object $\OTT_0$ in $\ttr 1$, and map $\top$ to the terminal object $\CTT_0$ in $\ttr 1$ (see \autoref{rmk:ebun-initial-terminal}).
\end{rmk}

\begin{eg}[Suspension bundles] In \autoref{fig:suspension-bundle} we depict a 1-truss bundle $p : T \to B$ on the left together with its suspension bundle $\Sigma p : \Sigma T \to \Sigma B$ on the right (for clarity we colored arrows in the total poset and in the base poset in corresponding colors).
\begin{figure}[h!]
    \centering
    \def\svgwidth{1\columnwidth}
    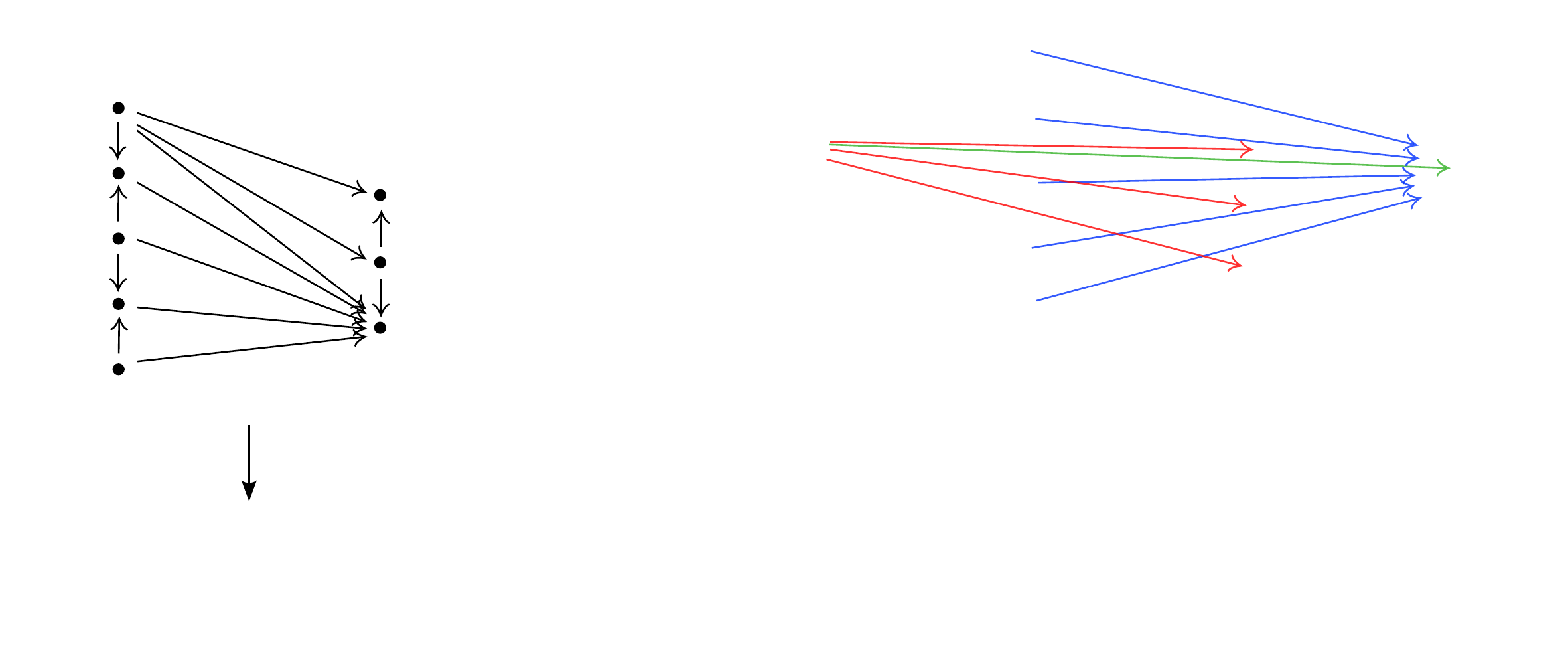

    \caption[Suspension bundle of a 1-truss bundle.]{Suspension bundle $\Sigma p$ of a 1-truss bundle $p$.}
    \label{fig:suspension-bundle}
\end{figure}
\end{eg}

\section{Truss induction along simplices in bundles} \label{sec:truss-induction}

We study simplices in total posets of 1-truss bundles and find natural linear orders on subsets of $k$-simplices that lie over the same base simplex; we will refer to this linear order as the `scaffold order' of the truss bundle. This observation will provide a useful tool in later combinatorial proofs, by allowing us to induct along a progression of simplices in the total poset of a bundle as determined by the scaffold order; we refer to this general principle as `truss induction'.

\subsection{Sections and spacers}

In analogy to notions of section and spacer simplices for proframings (see \autoref{term:sections-spacers-in-proframings}) we first introduce the following distinction of section and spacer simplexes in 1-truss bundles.

\subsubsecunnum{The definition of sections and spacers} Recall that a $k$-simplex $K : [k] \to P$ in a poset $P$ is called non-degenerate if it is injective on objects; otherwise we say it is degenerate.

\begin{defn}[Sections] \label{defn:sections} For a 1-truss bundle $p : T \to B$, a \textbf{$k$-section $K$ in $p$} is a non-degenerate $k$-simplex $K : [k] \into (T,\eleq)$ such that the composite map $p K : [k] \to B$ is a non-degenerate simplex in $B$.
\end{defn}

\begin{defn}[Spacers] \label{defn:spacers} For a 1-truss bundle $p : T \to B$, a \textbf{$k$-spacer $K$ in $p$} is a non-degenerate $k$-simplex $K : [k] \into (T,\eleq)$ such that $p K : [k] \to B$ is a degenerate simplex in $B$.
\end{defn}

\begin{term}[Simplices in bundles] Given a 1-truss bundle $p : T \to B$, we often say `a simplex in $p$' to mean a non-degenerate simplex $[k] \to (T,\eleq)$.
\end{term}

\nid Note that every simplex in 1-truss bundle $p$ is either a section or a spacer.

\begin{term}[Base projections] Let $K : [k] \into T$ be a $k$-simplex in a 1-truss bundle $p : T \to B$. We denote by $\im (pK) : [m] \into B$ the unique non-degenerate simplex in $B$ whose image equals that of $pK : [k] \into B$, and call it the `base projection' of $K$.
\end{term}

\nid In the special where the base poset $b$ is itself the $m$-simplex $B = [m]$, it further makes sense to ask for simplices in truss bundles to project to the `entire' base.

\begin{term}[Base-surjectivity] \label{defn:base-surjectivity} A $k$-simplex $K : [k] \into T$ in a 1-truss bundle $p : T \to B$ is called `base-surjective' if $\im(pK) = [m]$.
\end{term}

\begin{eg}[Base-surjective sections and spacers in 1-truss bundles] In \autoref{fig:sections-and-spacers-in-truss-bundles} we depict simplices in a 1-truss bundle $p : T \to B$: on the left we depict a selection of sections in $p$, on the right we depict a selection of spacers in $p$. All simplices are base-surjective.
\begin{figure}[h!]
    \centering
    \def\svgwidth{1\columnwidth}
    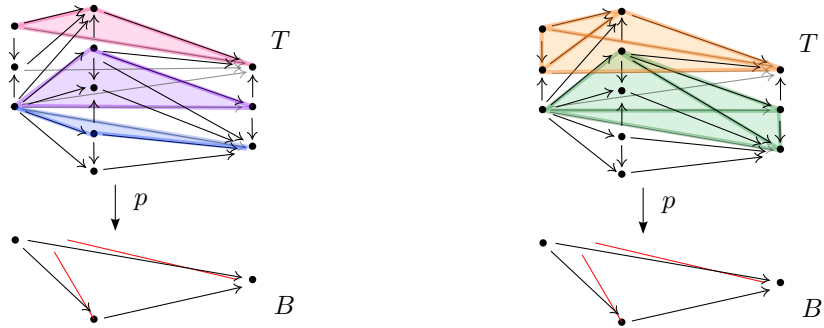

    \caption[Sections and spacers in a 1-truss bundle.]{Sections and spacers in a 1-truss bundle $p : T \to B$.}
    \label{fig:sections-and-spacers-in-truss-bundles}
\end{figure}
\end{eg}

\nid In the case of base-surjective simplices we may simplify their visualization by solely depicting their spine as explained by the following remark.

\begin{rmk}[Spine-only notation]
    \label{eg:base-surjective-simplices} On the left in \autoref{fig:base-surjective-simplices} we depict base-surjective sections (in blue and red) and spacers (in green and orange) of a 1-truss bundle over the 2-simplex $[2]$. On the right we depict the same simplices in the same bundle using a more convenient representation: namely, to depict $p : T \to [m]$, it suffices depict 1-truss bordisms lying over the spine of $[m]$ (note that this fully determines $p$), and similarly, to depict simplices in $p$, it suffices to depict their spines $T$.
\begin{figure}[h!]
    \centering
    \def\svgwidth{1\columnwidth}
    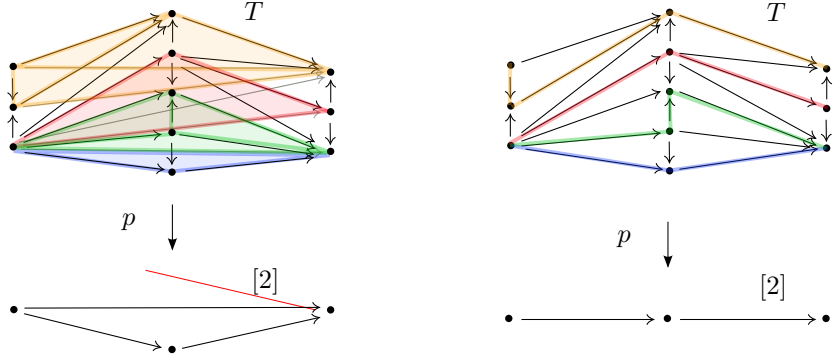

    \caption[Base-surjective sections and spacers.]{Base-surjective simplices may be represented by their spine which lies in the spine of the bundle.}
    \label{fig:base-surjective-simplices}
\end{figure}
\end{rmk}

\nid The case of 1-truss bundles over simplices and their base-surjective simplices is of particular importance: indeed, any simplex in any 1-truss bundle canonically factors as a base surjective simplex and a bundle inclusion as follows.

\begin{rmk}[Factoring simplices through base-surjective simplices] Let $K : [k] \into T$ be a $k$-simplex in a 1-truss bundle $p : T \to B$. Abbreviate its base projection $\im(pK) : [m] \into B$ by $F : [m] \into B$. Then $K : [k] \into T$ uniquely factors through the bundle pullback map $\Tot{}F : F^*p \into p$ as a $k$-simplex $F^*K : [k] \into F^*T$. Note $F^*K$ is base-surjective. In fact, the mapping $K \mapsto F^*K$ yields a bijection between simplices in the bundle $p$ whose base projection is $F$ and base-surjective simplices in the pullback bundle $F^*p$.
\end{rmk}

\nid In the rest of this section, we will almost exclusively work with bundles over simplices, and assume base-surjectivity by default.

\begin{conv}[Base-surjectivity by default] All sections and spacers will be assumed to be base-surjective unless otherwise noted.
\end{conv}


\subsubsecunnum{The spines of sections and spacers} Sections and spacers admit concrete combinatorial descriptions of their respective spines, as we will now explain. The fundamental observation is as follows.

\begin{obs}[Spines of simplices in truss bundles] \label{obs:spines-of-simplices} consider a 1-truss bundle $p : T \to [m]$ over the $m$-simplex and a $k$-simplex $K : [k] \to T$ in $p$. Recall, if $(i,a) \to (j,b)$ is a morphism in $T$, then we must have $\dim(a) \geq \dim(b)$ (see \autoref{lem:truss-trans-dimensions}). In particular, for each spine vector $(i-1 \to i)$ in $[k]$ we must have $\dim K(i-1) \geq \dim K(i)$. As a consequence, simplices in truss bundles fall into one of three categories.
    \begin{enumerate}
        \item The chain of morphisms $K(i-1 \to i)$ contains only regular objects $K(i)$ (i.e.\ $\dim K(i) = 1$). Since all non-identity arrows in fibers of $p$ run between regular and singular objects, $K$ must be a section simplex; we call $K$ a `purely regular' section in this case.
        \item The chain of morphisms $K(i-1 \to i)$ contains only singular objects $K(i)$ (i.e.\ $\dim K(i) = 0$). Similar to the previous case, $K$ must be a section simplex; we call it a `purely singular' section in this case.
        \item There is a unique morphism $K(j-1 \to j)$ whose domain is regular and whose codomain is singular. The index $j$ is called the `transition index' of $K$. All objects $K(i < j)$ are regular, and all objects $K(i \geq j)$ are singular. \qedhere
\end{enumerate}
\end{obs}

\nid The next two remarks spell out the preceding observation separately in the cases of sections resp.\ spacers. Both remarks use the following `tuple' notation to keep track of base projections of total poset elements in truss bundles.

\begin{notn}[Pair notation for 1-truss bundles] \label{notn:tuple-notation} Given a 1-truss bundle $p : T \to B$, we will sometimes redundantly denote objects $a \in T$ by pairs $(x,a)$: here, $x \in B$ is the object $p(a)$ in the base poset $B$, and thus $a \in p\inv(x)$ is lives in the 1-truss fiber over $x$.
\end{notn}

\begin{rmk}[Spines of sections] \label{obs:spines-of-sections} If $K$ is a section simplex in a 1-truss bundle $p : T \to [m]$, then the spine vectors $K(i \to i+1)$ of $K$ must be a chain of morphisms in $T$ of the form
\begin{equation}
    (0,a_0) \to (1,a_1) \to ... \to (j-1,a_{j-1}) \to (j,{b_j}) \to (j+1,{b_{j+1}}) \to... \to  (m,b_m)
\end{equation}
where each $a_i$ is regular, each $b_j$ is singular and $j$ is a unique index $0 \leq j \leq m + 1$. Corresponding to the first two cases of \autoref{obs:spines-of-simplices} note the following: if $j = 0$ then $K$ is purely singular, and conversely, if $j = m+1$ then $K$ is purely regular. In either of these boundary cases, we again refer to $j$ as the `transition index' of $K$.
\end{rmk}

\begin{rmk}[Spines of spacers] \label{obs:spines-of-spacers} If $L$ is a spacer in $p : T \to [m]$, then the spine vectors $L(i \to i+1)$ of $L$ must be a chain of morphisms in $T$ of the form
\begin{equation}
    (0,a_0) \to  (1,a_1) \to  ... \to (j,a_{j}) \to  (j,{b_j}) \to (j+1,{b_{j+1}}) \to... \to  (m,b_m)
\end{equation}
where each $a_i$ is regular, each $b_j$ is singular and $j$ is a unique index $0 \leq j \leq m$. Note that in particular we must have $\dim(K) = m+1$, that is, all filler simplices in bundles over the $m$-simplex are of dimension $(m+1)$. Corresponding to the third case in \autoref{obs:spines-of-simplices}, the index $j$ is called the transition index of $K$. Note that the fiber $p\inv(j)$ over $j$ contains \emph{two} vertices $\La(i)$ of $\La$ (a regular and a singular one), and that every spacer in particular has at least one regular and one singular vertex.
\end{rmk}

The boundary cases observed with sections (in which all elements are either regular or singular) can be resolved by passing to suspensions. Recall the construction of suspension bundles from \autoref{constr:1-truss-bundle-suspension}. The following notation will be useful.

\begin{notn}[Suspending simplices] \label{notn:suspending-simplices} For numeric convenience, we will usually identify the suspension $\Sigma [m]$ of the $m$-simplex $[m]$ with the poset $(-1 \to 0 \to 1 \to ... \to (m-1) \to m \to m+1)$.
\end{notn}

\begin{constr}[Suspending sections] \label{obs:suspending-sections}
    Let $p : T \to [m]$ be a 1-truss bundle over the $m$-simplex $[m]$, and consider a section $K : [m] \to T$ in $p$. Construct the suspended section $\Sigma K : \Sigma [m] \to \Sigma T$ by setting $\Sigma K(\bot) = \bot$ and $\Sigma K (\top) = \top$ (where $\bot$ respectively $\top$ denote newly adjoined initial respectively terminal objects); equivalently, using our numeric notational conventions (see \autoref{notn:tuple-notation} and \autoref{notn:suspending-simplices}), this may be written as $\Sigma K (-1) = (-1,0)$ and $\Sigma K (m+1) = (m+1,0)$. Observe that the mapping $K \mapsto \Sigma K$ establishes a 1-to-1 correspondence between sections in $p$ and in $\Sigma p$: indeed, the inverse mapping is obtained by restricting sections $\Sigma [m] \to \Sigma T$ along $[m] \into \Sigma [m]$ and $T \into \Sigma T$.
\end{constr}

\nid Note since $\Sigma K$ adjoint a new regular object to the start of $K$ and a new singular object to the end of $K$ it follows that $\Sigma K$ always contains at least one regular and one singular object (even if $K$ doesn't). In other words, suspended sections are never purely regular or purely singular.

    We now characterize sections and spacers in terms of certain morphism. We start with the case of sections in a 1-truss bundle $p$, which we will show correspond to so-called `jump morphism' in the suspension bundle $\Sigma p$.

\begin{defn}[Jump morphisms] Let $p : T \to [m]$ be a 1-truss bundle over the $m$-simplex. A \textbf{jump morphism} $f$ in $p$ is a morphism $(T,\eleq)$ whose domain $\dom(f)$ is regular, whose codomain $\cod(f)$ is singular and whose base projection $pf$ is a spine vector in $[i]$ (that is, $pf$ is of the form $i-1 \to i$ for $0 < i \leq m$).
\end{defn}

\begin{constr}[Correspondence of sections in $p$ and jump morphisms in $\Sigma p$] \label{obs:sections-and-jump-morphisms}  \label{constr:jumps-of-sections} Let $p : T \to [m]$ be a 1-truss bundle over the $m$-simplex $[m]$. Every section $K : [m] \to T$ in $p$ has a `associated jump morphism' $f$ in $\Sigma p$, given by $f = \Sigma K (j-1 \to j)$ where $j$ is the transition index of $K$. Conversely, let $f$ be a jump morphism in $\Sigma p$ lying over a spine vector $j-1 \to j$ in $\Sigma [m]$ (use \autoref{notn:suspending-simplices} for objects in $\Sigma[m]$). Then $f$ has an `associated section' $K : [m] \to T$ constructed by setting
    \begin{enumerate}
        \item for $i < j$, $\Ka(i) = \regbif^{\fcl{}p (i \to j-1)}(\dom f)$,
        \item for $i \geq j$, $\Ka(i) = \singbif_{\fcl{}p (j \to i)}(\cod f)$.
    \end{enumerate}
The constructions are mutually inverse, thus providing  a correspondence between sections in $K$ and jump morphisms in $\Sigma p$.
\end{constr}

\begin{eg}[Correspondence of sections and jump morphisms] In \autoref{fig:sections-vs-jump-morphisms} we depict a bundle $p : T \to [2]$, together with its suspension $\Sigma p$ (highlighted in purple). We indicate regular elements of $T$ in blue and singular elements in red. We colored the spines of several sections in $P$ (using the `spine-only' notation discussed in \autoref{eg:base-surjective-simplices}); in each case we then mark their associated jump morphism in $\Sigma p$ by a big dot of the same color.
\begin{figure}[h!]
    \centering
    \def\svgwidth{1\columnwidth}
    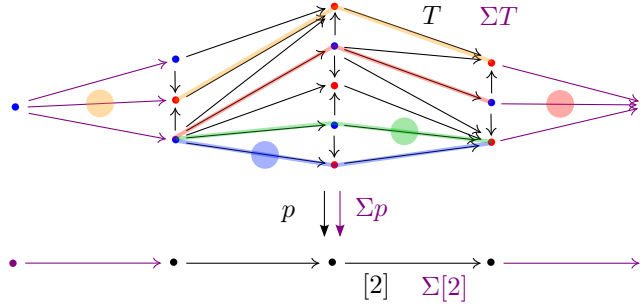

    \caption[Sections and their associated jump morphisms.]{Sections and their associated jump morphisms.}
    \label{fig:sections-vs-jump-morphisms}
\end{figure}
\end{eg}

Next, we establish a similar correspondence between spacers and so-called `fiber morphisms' defined as follows.

\begin{defn}[Fiber morphisms] Let $p : T \to [m]$ be a 1-truss bundle over the $m$-simplex. A \textbf{fiber morphism} $f$ in $p$ is a morphism $f$ in $(T,\eleq)$  whose domain $\dom(f)$ is regular, whose codomain $\cod(f)$ is singular, and whose base projection $pf$ is an identity morphism in $[m]$.
\end{defn}

\begin{constr}[Correspondence of spacers and fiber morphisms] \label{obs:spacers-vs-fiber-mor} Let $p : T \to [m]$ be a 1-truss bundle over the $m$-simplex $[m]$. Every spacer $L : [m + 1] \to T$ in $p$ has a `associated fiber morphism' $f$ in $p$, given by $f = L (j-1 \to j)$ where $j$ is the transition index of $L$. Conversely, let $f$ be a fiber morphism in $\Sigma p$ lying the object $j$ in $[m]$. Then $f$ has an `associated spacer' $L : [m + 1] \to T$ constructed by setting
\begin{enumerate}
    \item for $i \leq j$, $\La(i) = \regbif^{\fcl{}p (i \to j)}(\dom f)$,
    \item for $i > j$, $\La(i) = \singbif_{\fcl{}p (j \to i)}(\cod f)$.
\end{enumerate}
The constructions are mutually inverse, thus providing a correspondence between spacers in $K$ and fiber morphisms in $p$.
\end{constr}

\begin{eg}[Correspondence of spacers and fiber morphisms] In \autoref{fig:spacers-vs-fiber-morphisms} we colored the spines of several spacers in a bundle $p : T \to [2]$; we then mark their associated jump morphisms in $p$ by a big dot of the same color.
\begin{figure}[h!]
    \centering
    \def\svgwidth{1\columnwidth}
    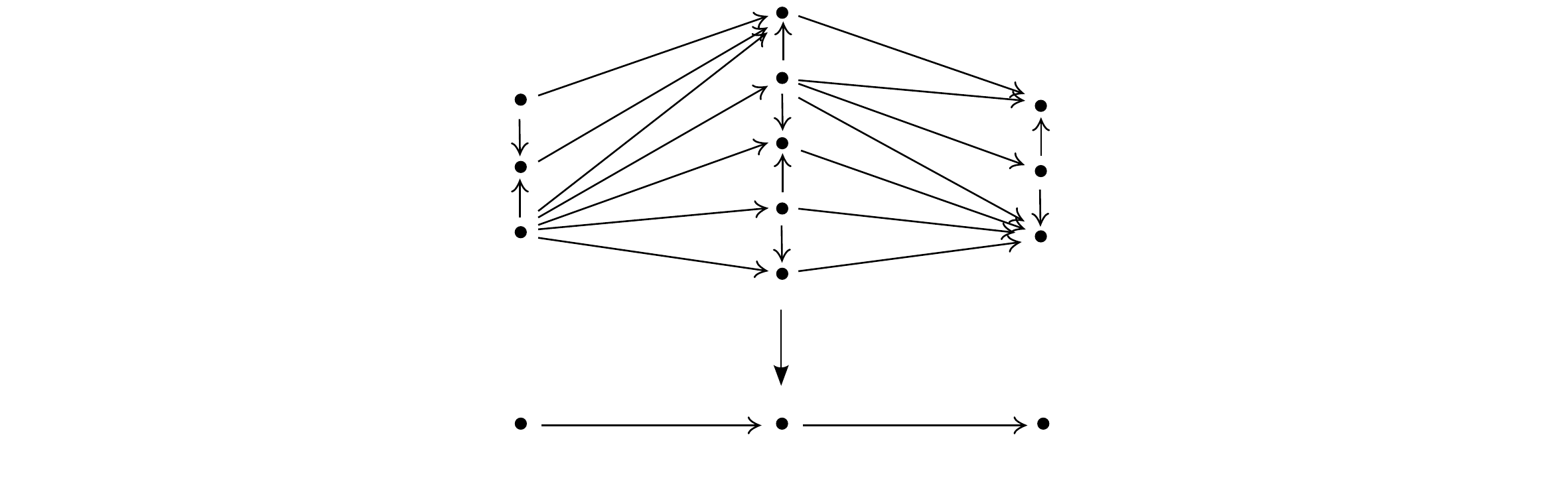

    \caption[Spacers and their associated fiber morphisms.]{Spacers and their associated fiber morphisms (indicated by dots with corresponding color).}
    \label{fig:spacers-vs-fiber-morphisms}
\end{figure}
\end{eg}

\subsection{The scaffold order}

We now construct canonical linear orders on the sets of sections and spacers in a 1-truss bundle over a simplex.

\begin{notn}[Set of sections and spacers] \label{notn:sections-and-spacers} Given a 1-truss bundle $p : T \to [m]$, we denote its sets of sections and spacers as follows.
\begin{align}
    \Gamma_p &= \{ \text{sections $K : [m] \to T$ of $p : T \to [m]$ } \}\\
    \Psi_p &= \{ \text{spacers $L : [m+1] \to T$ of $p : T \to [m]$} \} \qedhere
\end{align}
\end{notn}

\subsubsecunnum{The case of sections}

We first construct a total order on the set of sections $\Gamma_p$ of a 1-truss bundle $p : T \to [m]$ over the $m$-simplex: we call this order the `scaffold order of sections'. Using the correspondence of sections in $p$ and jump morphisms in the suspension $\Sigma p$ (see \autoref{obs:sections-and-jump-morphisms}) one may easily visualize the scaffold order as follows, and in fact, our formal construction of the scaffold order will directly rely on the idea underlying this visualization.

\begin{eg}[Scaffold orders] \label{eg:scaffold-norm-2} In \autoref{fig:scaffold-order-on-sections} we depict \emph{all} sections in a 1-truss bundle $p : T \to [2]$ over the 2-simplex (by coloring their spine). We also depict their corresponding jump morphisms in $\Sigma p$ (by correspondingly colored dots). The scaffold order on sections $\Gamma_p$, and thus on jump morphisms, is the total order indicated by a directed path (in red) passing through all jump morphisms in order. This directed path is uniquely determined as follows: it is (up to reversal) the unique directed path that passes through all jump morphisms of the suspended bundle once while intersecting exactly one fiber morphism in between any two consecutive jump morphisms.
\begin{figure}[h!]
    \centering
    \def\svgwidth{1\columnwidth}
    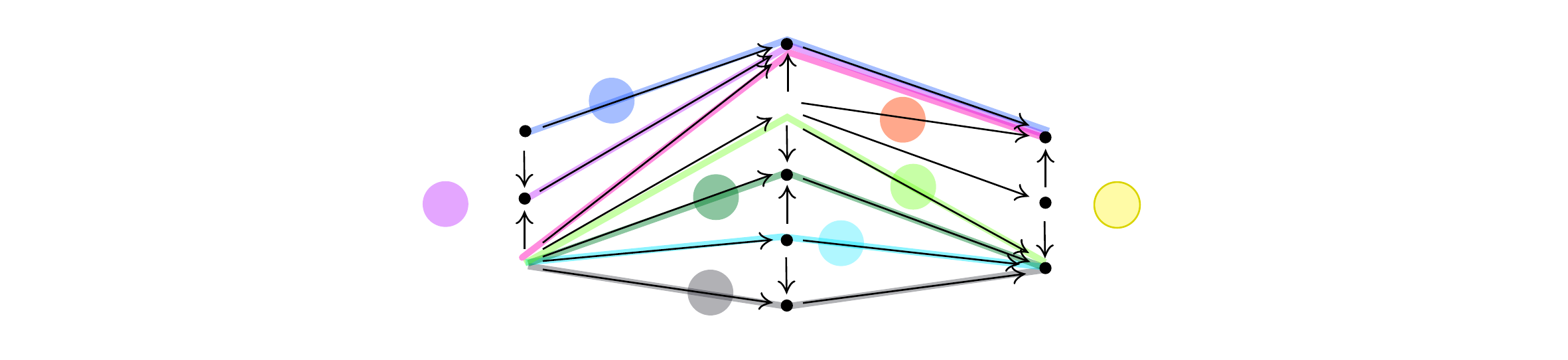

    \caption[The scaffold order on sections over 2-simplex.]{The scaffold order on sections in a 1-truss bundle over a 2-simplex indicated by a directed path through corresponding jump morphisms.}
    \label{fig:scaffold-order-on-sections}
\end{figure}
A slightly more complex example, using the same directed path notation to illustrate the scaffold order of jump morphisms (and thus of sections), is shown in \autoref{fig:scaffold-order-on-sections-2}: this depicts the scaffold order of sections in a bundle $p : T \to [3]$ over the 3-simplex.
\begin{figure}[h!]
    \centering
    \def\svgwidth{.8\columnwidth}
    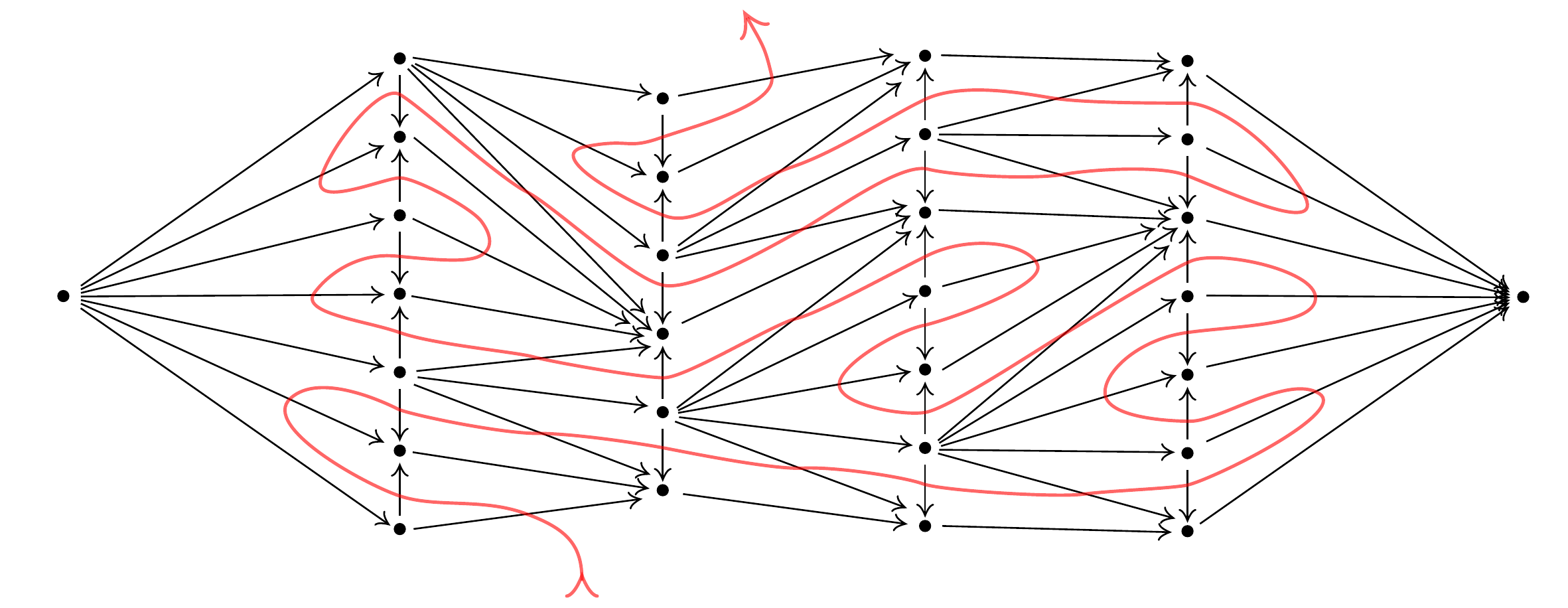

    \vspace{8pt}
    \caption[The scaffold order on sections over 3-simplex.]{The scaffold order on sections in a 1-truss bundle over a 3-simplex indicated by a directed path.}
    \label{fig:scaffold-order-on-sections-2}
\end{figure}
\end{eg}

We now formally construct the scaffold order on sections $\Gamma_p$. We first introduce a `scaffold norm'; this is a function which associates to each section $K \in \Gamma_p$ a natural number. We will show that, in fact, the scaffold norm is an injective function $\Gamma_p \into \lN$; endowing $\lN$ with the standard order of natural numbers, this will then induce the scaffold order on $\Gamma_p$ as required.

\begin{defn}[Scaffold norm of sections] \label{def:scaff-norm-sections}
    Consider a 1-truss bundle $p : T \to [m]$ and its set of sections $\Gamma_p$. The \textbf{scaffold norm} $\avg{-}$ is the function
\begin{align*}
\avg - : \Gamma_p  &\to \lN \\
\Ka & \mapsto \sum_{i \in [m]} \mathsf{num}(K(i))
\end{align*}
where $\mathsf{num} : (T,\fleq) \to (\lN,\leq)$ takes elements $a \in T$ to numbers $i-1 \in \lN$ iff $a$ is the $i$th element in the frame ordered fiber over $p(a)$.\footnote{In other words, $\mathsf{num} : (T,\fleq) \to (\lN,\leq)$ restricts on fibers $p\inv(i)$, $i \in [m]$, to the unique isomorphisms $\mathsf{num} : (p\inv(i),\fleq) \iso [n] \into (\lN,\leq)$ of frame ordered fibers with standard simplices.}
\end{defn}

\begin{obs}[Suspension preserves scaffold norm] \label{obs:suspension-scaff} Recall the suspension operation on sections defined in \autoref{obs:suspending-sections}. Observe that taking suspensions $\Sigma : \Gamma_p \to \Gamma_{\Sigma p}$ preserves scaffold norm, that is, $\avg{K} = \avg{\Sigma K}$ for $K \in \Gamma_p$.
\end{obs}

In order to describe the image of the scaffold norm, we first construct its `bottom' and `top sections'.

\begin{constr}[Minimal and maximal values of the scaffold norm] \label{defn:bottom-top-sections}
    Let $p : T \to [m]$ be a 1-truss bundle over the $m$-simplex $[m]$. We construct sections $\iK^-_p, \iK^+_p : [m] \to T$ of $p$, called the \textbf{bottom} respectively \textbf{top} sections of $p$, which are minimal respectively maximal sections with respect to the scaffold norm. Namely, we define $\iK^\pm_p$ to map $i \in [m]$, to the lower respectively upper endpoints of the fiber over $i$, that is, $\iK^\pm_p(i) = \ept_\pm p\inv(i)$. The fact that this defines valid sections follows since 1-truss bordisms relate endpoints (see \autoref{lem:truss-trans-endpoints}). Further, observe that the scaffold norm indeed attains its minimal and maximal value on (and only on) $\iK^\pm_p(i)$: namely, its minimal value is $\avg{\iK^-_p} = 0$ and its maximal values is $\avg{\iK^+_p} = \#T-\#B$ (where $\#T$ and $\#B$ are the numbers of elements in $T$ respectively $B$).
\end{constr}

\nid We will denote extremal values of the scaffold norm, achieved on the bottom resp.\ top sections of a bundle $p$, by $\scaff_\pm (p) := \avg {\iK^\pm_p}$. These values bound an interval of natural numbers, which described exactly the image of the scaffold norm.

\begin{lem}[Scaffold order of sections] \label{lem:section-induction}
    Given a 1-truss bundle $p$, the scaffold norm $\avg - : \Gamma_p  \to \lN$ on sections is a bijection with its image; this image consist of all integers between $\scaff_- (p)$ and $\scaff_+ (p)$. As a result, $\Gamma_p$ obtains a total order, called the `scaffold order of sections'.
\end{lem}

\begin{proof}[Proof of \autoref{lem:section-induction}] In \autoref{defn:bottom-top-sections} we showed that the scaffold norm has unique bottom and top sections $\iK^\pm_p$. We now construct for each $K \neq \iK^+_p$ a successor section $\successor(K)$ with scaffold norm $\avg{\successor(K)} = \avg{K} + 1$, and, conversely, for each section $K \neq \iK^-_p$ a predecessor section $\successor(K)$ with scaffold norm $\avg{\predecessor(K)} = \avg{K} - 1$. Showing that successor and predecessor section constructions are mutually inverse will prove the claim. We can further make the following assumption (which simplifies the treatment of the `boundary cases' in \autoref{obs:spines-of-sections}). Recall, by \autoref{obs:suspending-sections} we have an isomorphism $\Gamma_p \iso \Gamma_{\Sigma p}$, and by \autoref{obs:suspension-scaff} this isomorphism commutes with the scaffold norm. After potentially replacing $p$ with its suspension $\Sigma p$ in the lemma, we may therefore assume that all sections $K$ in $p : T \to [m]$ have transition index $j$ in the range $0 < j < m+1$ (in other words, $p$ does not contain purely regular or purely singular sections).

    Let $K$ be a section of $p$ with jump morphism $K(j-1) \to K(j)$ in $T$. Assume that $K \neq \iK^+_p$, i.e.\ $K$ is not the top section. We claim that then there is \emph{either} an arrow $K(j-1) + 1 \to K(j)$ \emph{or} an arrow $K(j-1) \to K(j) + 1$ in $T$ (but there can never be both, due to bimonotonicity of 1-truss bordisms). Note, by assumption on $K$, there is at least on index $l \in [m]$ such that $K(l)$ has a successor $a \equiv K(l) + 1$ in the fiber over $l$. We verify the claim in the case where $l < j$ (the case $l \geq j$ follows by the same argument after dualizing $p$). Note if $l < j$ then $K(l)$ is regular. Thus the successor $a$ is singular, and is mapped by the singular function of the 1-truss bordism $p\inv(l \to j-1)$ to a singular element $b$ in the fiber over $(j-1)$. Since $K(l) \fles a$, bimonotonicity implies that $K(j-1) \fles b$. In particular, the regular element $K(j-1)$ has a singular successor $K(j-1) + 1$ in the fiber over $j-1$ in this case. Now, again by bimonotonicity, the 1-truss bordism $p\inv(j-1 \to j)$ must relate the singular element $K(j-1) + 1$ either to $K(j)$, or otherwise to some singular $b$ with $K(j) \fles b$. Using the description of singular determined 1-truss bordisms (see \autoref{lem:truss-trans-representations}), this implies there is an arrow $K(j-1) \to K(j) + 1$, which verifies the claim.

    We now construct the successor section $\successor(K)$ of $K$ (still assuming $K \neq \iK^+_p$) by distinguishing the two cases established above, namely, whether $K(j-1) + 1 \to K(j)$ or $K(j-1) \to K(j) + 1$. Both case are illustrated in \autoref{fig:successor-construction-cases}. In each case we highlight singular elements in red, regular element in blue, and the jump morphism $K(j-1) \to K(j)$ of $K$ is marked by green dot. The purple dot marks the jump morphism of the successor section of $K$. In general, the successor section can now be constructed as follows.
\begin{figure}[h!]
    \centering
    \def\svgwidth{1\columnwidth}
    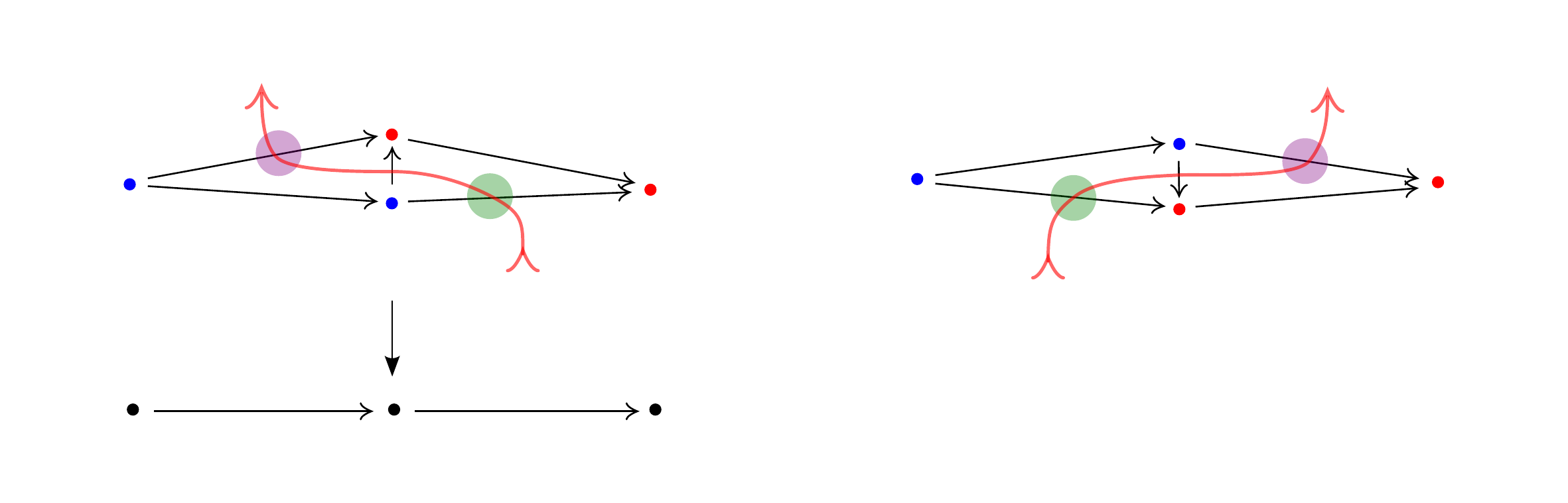

    \caption[The construction of successors.]{The construction of successor sections.}
    \label{fig:successor-construction-cases}
\end{figure}
    \begin{enumerate}
        \item[] \textbf{Case 1} If there is an arrow $K(j-1) + 1 \to K(j)$, we construct the successor section $\successor(K)$ by setting $\successor(K)(j-1) = K(j-1)+1$ and $\successor(K)(i) = K(i)$ if $i \neq j-1$. By assumption on $p$ we must have $j > 1$ (otherwise, $\successor(K)$ would be a purely singular section). Profunctoriality of the 1-truss bordism $p\inv(j-2 \to j-1)$ then implies there is an arrow $K(j-2) \to K(j-1) + 1$ in $T$, showing that $\successor(K)$ is a valid section.

    \item[] \textbf{Case 2} If there is an arrow $K(j-1) \to K(j) + 1$ in $T$, we construct the successor $\successor(K)$ by setting $\successor(K)(j) = K(j)+1$ and $\successor(K)(i) = K(i)$ if $i \neq j$. By assumption on $p$ we must have $j < m$ (otherwise, $\successor(K)$ would be a purely regular section). Profunctoriality of the 1-truss bordism $p\inv(j \to j+1)$ then implies there is an arrow $K(j)+1 \to K(j+1)$ in $T$, showing that $\successor(K)$ is a valid section.
    \end{enumerate}
    \nid This completes the construction of successors. The construction of predecessor sections $\predecessor(\Ka)$ of $\Ka$ (for $\Ka \neq \iK^-_p$ not equal to the bottom section) reduces the construction of successors after passing to the opposite frame order (i.e.\ flipping the frame order on each fiber of $p$, which amounts to reading the total posets in \autoref{fig:successor-construction-cases} upside down). One verifies that $\predecessor(\successor(\Ka)) = \Ka$ and similarly $\successor(\predecessor(\Ka)) = \Ka$, which completes the proof.
\end{proof}


\subsubsecunnum{The case of spacers}

Similar to the case of sections, we will now show that the set of spacers $\Psi_p$ in a 1-truss bundle $p : T \to [m]$ over the $m$-simplex carries a natural total order, called the `scaffold order of spacers'.

\begin{eg}[Scaffold order on spacers] Recall from \autoref{eg:scaffold-norm-2} that the scaffold norm of sections corresponds to a directed path that passes through all jump morphism of the suspended bundle once while intersecting exactly one fiber morphism in between any two consecutive jump morphisms---the order in which this path travels through fiber morphisms is exactly the scaffold order on spacers! We illustrate this in \autoref{fig:scaffold-norm-on-spacers} for a truss bundle over the 2-simplex, highlighting fiber morphism by colored dots.
\begin{figure}[h!]
    \centering
    \def\svgwidth{1\columnwidth}
    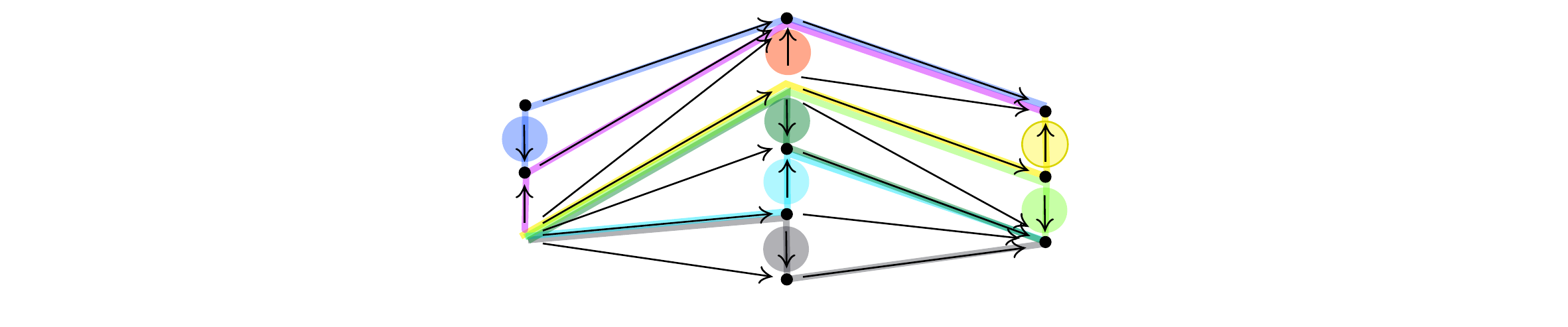

    \caption[The scaffold order on spacers.]{The scaffold order on spacers in a 1-truss bundle over the 2-simplex indicated by a directed path through fiber morphisms.}
    \label{fig:scaffold-norm-on-spacers}
\end{figure}
\end{eg}

Fixing a spacer and its corresponding fiber morphism in \autoref{fig:scaffold-norm-on-spacers}, then the `consecutive jump morphism' which the path intersects immediately before the given fiber morphisms correspond to so-called `lower and upper boundary sections' of the spacer. We formally define these as follows.

\begin{constr}[Upper and lower boundaries of spacer] \label{defn:upper-and-lower-sec} For a spacer $L : [m+1] \to T$ of a 1-truss bundle $p : T \to [m]$ we construct sections $\partial_- L, \partial_+ L : [m] \to T$ called the \textbf{lower boundary} respectively \textbf{upper boundary} of $L$. Let $j$ denote the transition index of $L$, that is, $L(j) \to L(j+1)$ is the fiber morphisms of $L$. If $L(j) \fles L(j+1)$, then we construct $\partial_- L$ as the $(j+1)$th face $Ld_{j+1}$ of $L$, and $\partial_+ L$ as the $j$th face $Ld_j$ of $L$. Otherwise, if $L(j+1) \fles L(j)$, then we construct $\partial_- L$ as the $j$th face $Ld_j$ of $L$, and $\partial_+ L$ as the $(j+1)$th face $Ld_{j+1}$.
\end{constr}

\begin{eg}[Upper and lower boundaries] We depict two spacers $L$ and $L'$ in a 1-truss bundle together with the lower and upper boundary sections $\partial_\pm L$ respectively $\partial_\pm L'$ in \autoref{fig:boundaries-of-spacers}.
\begin{figure}[h!]
    \centering
    \def\svgwidth{1\columnwidth}
    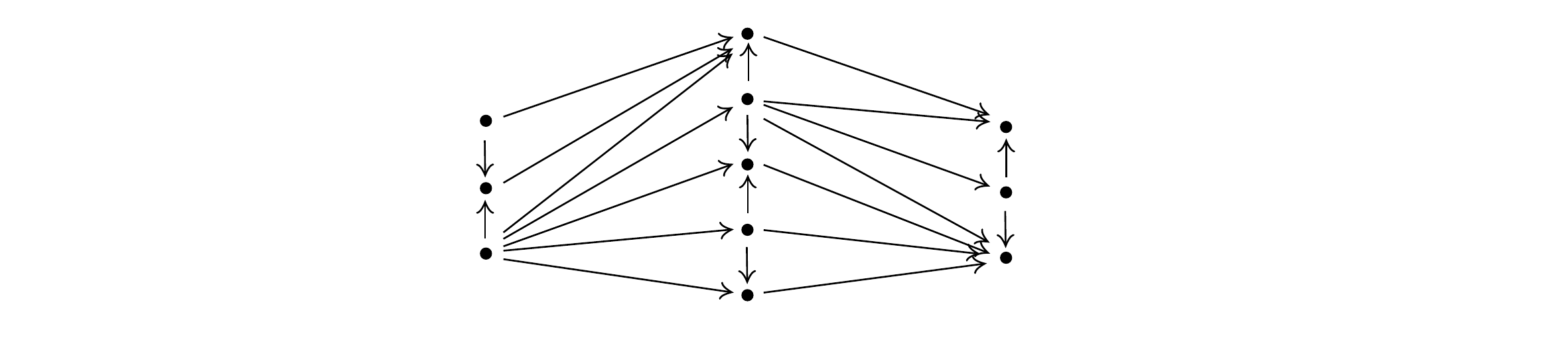

    \caption{Upper and lower boundaries of spacers.}
    \label{fig:boundaries-of-spacers}
\end{figure}
\end{eg}

\begin{rmk}[Uppers boundaries succeed lower boundaries] Note that the preceding construction ensures that we have $\avg{\partial_+ L} = \avg{\partial_- L} + 1$, that is, $\partial_+ L$ is the successor of $\partial_- L$ in the totally ordered set $(\Gamma_p,\fleq)$ of sections in $p$.
\end{rmk}

\nid In analogy to the definition of scaffold norms for sections (see \autoref{def:scaff-norm-sections}), we introduce the following scaffold norm of spacers.

\begin{defn}[Scaffold norm of spacers] \label{def:scaff-norm-spacers}
    Consider a 1-truss bundle $p : T \to [m]$ and its set of spacers $\Psi_p$ (see \autoref{notn:sections-and-spacers}). The \textbf{scaffold norm} $\avg{-}$ is the function
\begin{align*}
\avg - : \Psi_p  &\to \frac{1}{2}\lN\\
\La  &\mapsto  \frac{\avg{\partial_- \La} + \avg{\partial_+ \La}}{2}
\end{align*}
where $\avg{\partial_\pm L}$ denotes the scaffold norm of the sections $\partial_\pm L$.
\end{defn}


\nid In analogy to the construction of the scaffold order of sections in \autoref{lem:section-induction}, we now construct the scaffold order of spacers; as before, the order is induced by exhibiting $\Psi_p$ as a subset of a totally ordered set (in this case, the half integers $\frac{1}{2}\lN$) via the scaffold norm function.

\begin{lem}[Scaffold order for spacers] \label{lem:spacer-induction}
    Let $p : T \to [m]$ be a 1-truss bundle over the $m$-simplex $[m]$.  The scaffold norm $\avg - : \Psi_p  \to \frac{1}{2}\lN$ on spacers is a bijection with its image; this image consist of all half integers between $\scaff_- (p) + \frac{1}{2}$ and $\scaff_+ (p) - \frac{1}{2}$. As a result, $\Psi_p$ obtains a total order, called the `scaffold order of spacers'.
\end{lem}

\begin{proof}[Proof of \autoref{lem:spacer-induction}] Note that each spacer $\La$ in $p$ is uniquely determined by its boundary sections $\partial_\pm \La$. By the previous \autoref{lem:section-induction} it follows that $\avg - : \Psi_p  \to \frac{1}{2}\lN$ is injective. To see that it is also surjective, recall from the proof of \autoref{lem:section-induction} the construction of successors $\successor(K)$ of sections $K$ in $p$ that are not maximal in $(\Gamma_p,\fleq)$. Inspecting the two cases for the construction of $\successor(K)$, note that in case 1, the fiber morphism $\successor(K)(j-1) \to K(j-1)$ determines a spacer $L$ (see \autoref{obs:spacers-vs-fiber-mor}) with $\partial_- L = K$ and $\partial_+ L = \successor(K)$. Similarly, in case 2, the fiber morphism $K(j) \to \successor(K)(j)$ determines a spacer $L$ with $\partial_- L = K$ and $\partial_+ L = \successor(K)$. This shows any two consecutive sections are the boundary sections of a spacer, which thus proves surjectivity of the scaffold norm of spacers as claimed.
\end{proof}

Finally, as a first application of the results here, we will show that that truss bundles satisfy conditions analogous to those imposed on flat proframings (see \autoref{defn:flat-proframings}). This observation will prove useful for the later classification of framed regular cell complexes. The following mirrors the definition of fiber categories in proframed simplicial complexes (see \autoref{defn:fundamental-fiber-category}).

\begin{term}[Fiber categories in 1-truss bundles] Let $p : T \to B$ be a 1-truss bundle and consider a non-degenerate simplex $z : [m] \to B$. The `fiber category' $\FC p_z$ in $p$ (over $z$) is the free category whose objects are sections $K$ in $\Gamma_{z^*p}$ (that is, sections of the pullback bundle $z^*p$) together with a generating morphism $L : \partial_- L \to \partial_+ L$ for each spacer $L$ in $\Psi_{z^*p}$.
\end{term}

\begin{defn}[Transition functors of fiber categories]  Let $p : T \to B$ be a 1-truss bundle and consider non-degenerate simplices $z : [m] \to B$ and $y : [l] \to B$ such that $z$ is a face of $y$, that is, $z$ factors through $y$ by an injective map $[m] \into [l]$. We can restrict sections and spacers in the bundle $y^*p$ along the bundle inclusion $z^*p \into y^*p$ (note sections restrict to sections, while spacers restrict to sections or spacers). This induces a functor $\rest {-} {z \subset y} : \FC p_{y} \to \FC p_{z}$ called the \textbf{fiber transition} functor from $y$ to $z$ in $p$.
\end{defn}

\begin{obs}[`Flatness' of 1-truss bundles] \label{prop:1-truss-bun-euclidean} For all 1-truss bundles $p : T \to B$ the following holds.
    \begin{enumerate}
        \item All fiber categories in $p$ are total orders.
        \item All fiber transitions functors in $p$ are endpoint preserving (and thus surjective).
    \end{enumerate}
Indeed, the first statement follows from \autoref{lem:section-induction}
    and \autoref{lem:spacer-induction} together with the observation that $\avg{\partial_\pm L} = \avg{L} \pm \frac{1}{2}$. The second statement follows by construction of top and bottom sections (see \autoref{defn:bottom-top-sections}) as the reader can check.
\end{obs}

\section{$n$-Trusses, $n$-truss bordisms, $n$-truss bundles, and truss blocks} \label{sec:n-trusses}

We finally turn to the notion of $n$-trusses, which are `$n$-dimensional' generalizations of 1-trusses. Just as 1-trusses geometrically translate to (framed) stratified intervals, $n$-trusses geometrically translate to certain (framed) stratified subspaces of $\lR^n$. Before diving into the theory of $n$-trusses, we will define $n$-trusses and illustrate this geometric translation with examples. In fact, defining $n$-trusses is rather straight-forward: $n$-trusses are towers of 1-truss bundles over the trivial poset $[0]$ as follows.

\begin{defn}[$n$-Trusses] \label{defn:n-trusses} An \textbf{$n$-truss $T$} is a tower of 1-truss bundles
    \begin{equation}
        T_n \xto {p_n} T_{n-1} \xto {p_{n-1}} ... \xto {p_2} T_1 \xto {p_1} T_0 = [0]
    \end{equation}
where the total poset of each $p_i$ is the base poset of $p_{i+1}$. We call $T_n$ (with face order) the \textbf{total poset} of $T$.
\end{defn}

\begin{term}[Closed and open $n$-trusses] We call an $n$-truss $T$ \textbf{closed} respectively \textbf{open} if all its 1-truss bundles $p_i$ are closed respectively open.
\end{term}

Let us illustrate $n$-trusses by examples. In each case we also depict the `corresponding stratifications' in $\lR^n$: for now, this correspondence may be understood to identify entrance path posets of stratifications with face orders of the corresponding trusses. In \autoref{ch:meshes} we will make the correspondence fully precise, by also accounting for framing structures on both stratifications and trusses.

\begin{eg}[A closed $2$-truss] In \autoref{fig:a-first-2-truss-example} we illustrate a closed $2$-truss $T$. The bundle $p_1 : T_1 \to T_0$ has base post $T_0 = [0]$ and therefore its total poset is simply a 1-truss $T_1$. In contrast, the total poset $T_2$ of the bundle $p_2 : T_2 \to T_1$ looks more complicated---nonetheless it is `2-dimensional' in character as can be seen. Accordingly, the stratification that corresponds to $T$ can be embedded in $\lR^2$ as shown on the right in \autoref{fig:a-first-2-truss-example}---in fact, this stratification can also be understood as flat $2$-framed regular cell complex (see \autoref{defn:framed-reg-cell-cplx}). Note that the entrance path poset (see \autoref{defn:face-poset}) of this framed regular cell complex can be canonically identified with $T_2$.
\begin{figure}[ht]
    \centering
    \def\svgwidth{1\columnwidth}
    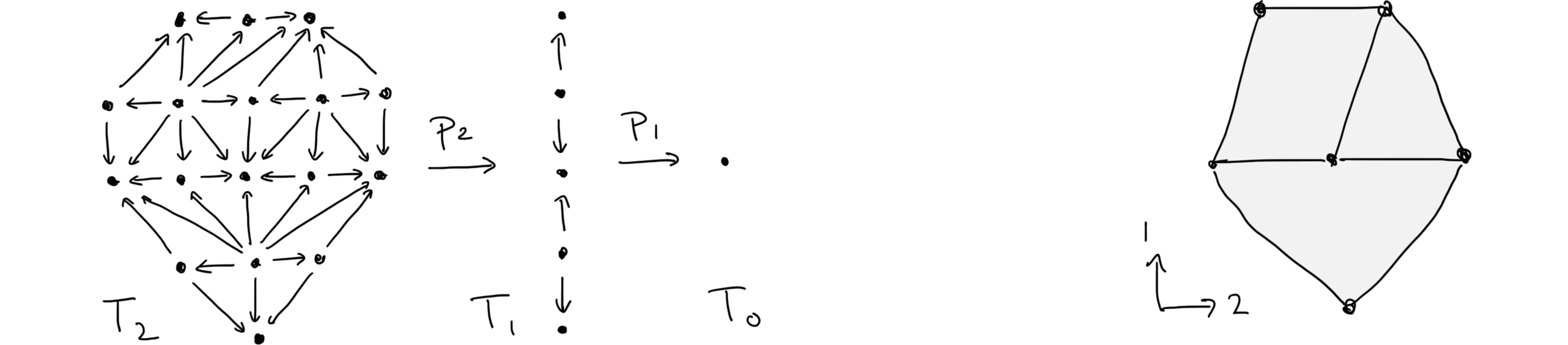

    \caption[An example of a closed 2-truss.]{A closed 2-truss together with its corresponding stratification.}
    \label{fig:a-first-2-truss-example}
\end{figure}
\end{eg}

\nid The previous example is an instance of a more general correspondence of \emph{closed} $n$-trusses and flat framed regular cell complexes: a proof of this correspondence will be the subject of \autoref{ch:classification-of-framed-cells}. The next example illustrates the case of an open 3-truss.

\begin{eg}[An open $3$-truss] Consider the open $3$-truss $T$ shown in \autoref{fig:a-first-3-truss-example}: note that, in order to simplify our notation, we depicted only `generating' arrows in $T_2$ and $T_3$ (with all other arrows being composites of the depicted ones). To the right of $T$ we illustrate its corresponding stratification. Note that this now yields a stratification of the \emph{open} 3-cube embedded in $\lR^3$. Note further, that objects in the total poset $T_3$ are again in correspondence with strata in the stratification realizing $T$, that is, the face order of $T_3$ can be identified with the entrance path poset of the stratification (see \autoref{defn:entr}).
\begin{figure}[ht]
    \centering
    \def\svgwidth{1\columnwidth}
    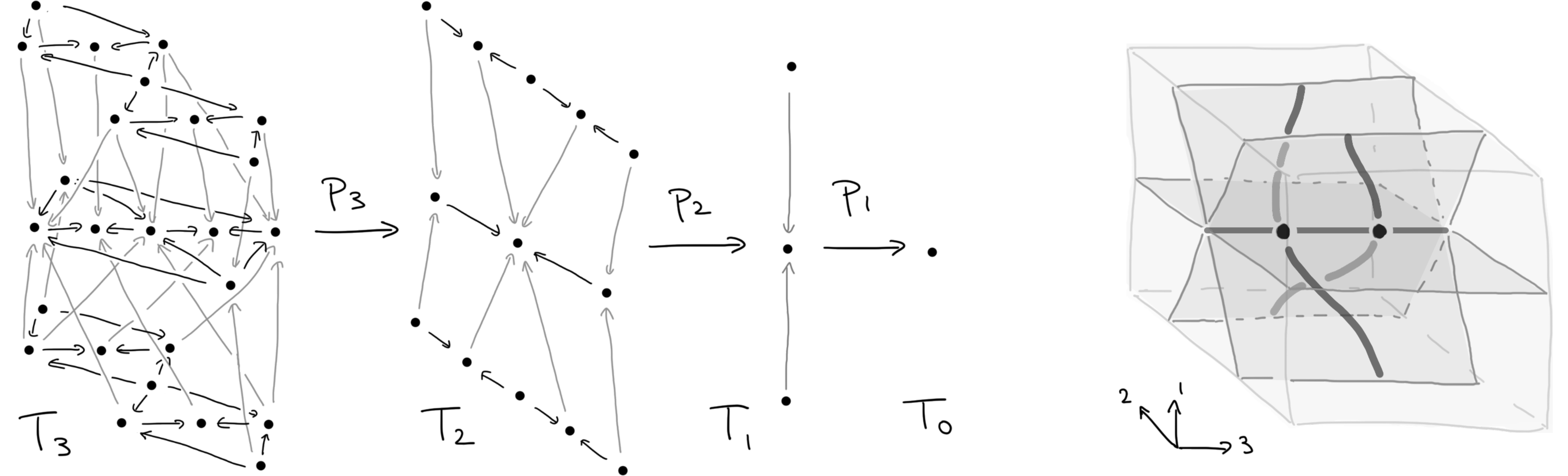

    \caption[An example of an open 3-truss.]{An open 3-truss together with its corresponding stratification.}
    \label{fig:a-first-3-truss-example}
\end{figure}
\end{eg}

\begin{preview}[Classifying the braid] \label{preview:stratified-trusses} The truss in the preceding example later plays a particularly interesting role: its corresponding stratification \emph{canonically} refines the braid\footnote{Here, the `braid' is the stratification of the 3-cube containing two `crossing' wires.}. This will formalized and generalized in \autoref{ch:hauptvermutung} to a correspondence of so-called `stratified trusses' and `flat framed stratifications'. We revisit the completed combinatorial description of the braid in \autoref{eg:revisiting-the-braid}.
\end{preview}

While our \autoref{defn:n-trusses} of $n$-trusses is short and simple, their combinatorial theory turns out to be rich. Moreover, properly dealing with notions such as `$n$-truss bordism' and `$n$-truss bundles' requires understanding a yet different aspect of truss theory, namely, of so-called `labelings' of bundles which we shall first introduce in the case of 1-truss bundles in \autoref{ssec:labeled-1-truss-bundle}. Labeled truss bundles then allow to categorically organize $n$-trusses in several interesting ways, as we will discuss in \autoref{ssec:n-trusses} (they also provide the correct structure for capturing our later definition of `stratified trusses'). In the final \autoref{ssec:block-and-brace-trusses}, we introduce the notion of `block' and `brace' trusses: these are the `building blocks' of $n$-trusses analogous to cell complexes being `built of cells'. More precisely, we will see that closed (and dually, open) trusses are gluings of blocks (resp.\ of braces). Importantly, considering (regular) presheaves on the category of such blocks will enable our classification of framed regular cell complexes in \autoref{ch:classification-of-framed-cells}

\subsection{Labeled 1-trusses, bordisms, and bundles} \label{ssec:labeled-1-truss-bundle}

Labeled 1-trusses are 1-trusses equipped with a `labeling' in a category $\iC$ as follows.

\begin{defn}[Labeled 1-trusses] \label{defn:labeled-1-trusses} Let $\iC$ be a `labeling' category. A \textbf{$\iC$-labeled 1-truss} $T$ is a pair consisting of an `underlying' truss $\und T$ together with a `labeling' functor $\lbl T : (\und T, \eleq) \to \iC$.
\end{defn}

\begin{rmk}[Labeled 1-trusses as spans] \label{rmk:labeled-1-trusses-as-spans}
    Given a $\iC$-labeled 1-truss $T \equiv (\und T, \lbl T)$ we may equivalently express its data as a span
    \begin{equation}
        \begin{tikzcd}
        \iC & \Totb{}{\und T} \arrow[r, "\fp {} {\und T}"] \arrow[l, "\lbl T"'] & {[0]}
        \end{tikzcd}
    \end{equation}
where the right leg $\fp {} {\und T}$ is a 1-truss bundle with underlying 1-truss $\und T : [0] \to \ttr 1$ (see \autoref{constr:totalization-1}), and the left leg is the labeling functor $\lbl T$ up to the canonical identification of 1-trusses $(\Totb{}{\und T},\eleq,\dim,\fleq) \iso (\und T,\eleq,\dim,\fleq)$.
\end{rmk}

\nid Many notions from the theory of 1-trusses as developed in \autoref{sec:1-trusses} naturally generalize to the labeled case. In particular, we now discuss a notion of `labeled 1-truss bordism', after which we define `labeled 1-truss bundles'.

\subsubsecunnum{The definition of labeled 1-truss bordisms}

We define labeled 1-truss bordism. In analogy to the definition of labeled 1-trusses, we want to equip 1-truss bordisms with functorial `labeling data' in a given category $\iC$. To make this precise, we will use that 1-truss bordisms $R$ correspond to 1-truss bundles $\fp{}{R} : \Totb{}{R} \to [1]$ over the 1-simplex $[1]$ via the classification-totalization correspondence for 1-truss bundles.

\begin{defn}[Labeled 1-truss bordisms] \label{defn:labeled-1-truss-bordism} For a category $\iC$, a \textbf{$\iC$-labeled 1-truss bordism} $R$ is a pair consisting of an `underlying' 1-truss bordism $\und R$ together with a `labeling' functor $\lbl R : (\Totb{}{\und R},\eleq) \to \iC$ from the total poset $\Totb{}{\und R}$ of $\und R$ into $\iC$.
\end{defn}

\nid In other words, a $\iC$-labeled 1-truss bordism $R \equiv (\und R,\lbl R)$ is a span of the form
    \begin{equation}
        \begin{tikzcd}
            \iC & \Totb{}{\und R} \arrow[r, "\fp {} {\und R}"] \arrow[l, "\lbl R"'] & {[1]}
        \end{tikzcd}
    \end{equation}
where $\fp {} {\und R}$ is the total bundle of $\und R : [1] \to \ttr 1$.

\begin{term}[Domain and codomain of labeled 1-truss bordisms] Given a $\iC$-labeled 1-truss bordism $R \equiv (\und R, \lbl R)$, the `domain' $\dom(R)$ (resp.\ `codomain' $\cod(R)$) of $R$ is the $\iC$-labeled 1-truss $(\rest {\Totb{}{\und R}} i, \lbl R : \rest {\Totb{}{\und R}} i \to \iC)$ obtained by restricting the totalization $\Tot{}{\und R}$ and the labeling $\lbl T$ to the fiber over the element $i = 0$ (resp.\ over $i = 1$) in the base $[1]$. We write this as $R : \dom(R) \xslashedrightarrow{} \cod(R)$.
\end{term}

\nid The interesting question concerning the definition of labeled 1-truss bordisms is whether labeled 1-truss bordisms have well-defined compositions.

\begin{defn}[Composition of labeled 1-truss bordisms] \label{defn:composition-witnesses-1} Given two $\iC$-labeled truss bordisms $R_{01}$ and $R_{12}$ which are composable (that is, $\cod(R_{01}) = \dom(R_{12})$), a \textbf{composition candidate} $R_{02}$ of $R_{12}$ and $R_{01}$ is a $\iC$-labeled truss bordism for which there exists `composition witness' $W \equiv (\und W : [2] \to \ttr 1, \lbl W : \Totb{}{\und W} \to \iC)$ which restricts on each arrow $(i \to j) : [1] \into [2]$ to $R_{ij}$; that is, $R_{ij}$ equals the $\iC$-labeled 1-truss bordism $(\rest {\Totb{}{\und W}}{i \to j}, \lbl W : \rest {\Totb{}{\und W}}{i \to j} \to \iC)$.
\end{defn}

\nid As we now show composition candidates and composition witnesses always uniquely exist. The proof will use `truss induction' as discussed in \autoref{sec:truss-induction}.

\begin{lem}[Existence and uniqueness of composition witnesses] \label{lem:comp-witnesses} Composition candidates and witnesses of composable $\iC$-labeled 1-truss bordisms exist uniquely.
\end{lem}

\begin{proof} Consider $\iC$-labeled 1-truss bordisms $R_{01}$ and $R_{12}$, composable at a $\iC$-labeled 1-truss $S = \cod(R_{01}) = \dom(R_{12})$. We construct their unique composition witness $W = (\und W, \lbl W)$ of (their composition candidate is determined by restricting $W$ over $(0 \to 2)$). First observe, the functor $\und W: [2] \to \ttr 1$ must be defined by setting $\und W(0 \to 1) = \und {R_{01}}$ and $\und W(1 \to 2) = \und {R_{12}}$. The task is to see that the labeling functor $\lbl W : \Tot{}{\und W} \to \iC$ exists uniquely.

    Since $W$ must restrict to $R_{01}$ over $(0 \to 1)$, and to $R_{12}$ over $(1 \to 2)$, it remains to define $\lbl W$ on arrows $f$ in $\Tot{}{\und W}$ lying over $0 \to 2$. Since 1-truss bordisms compose as their underlying relations, any arrow $f$ lying over $0 \to 2$ in $[2]$ can be written as a composite $f = gh$ for some arrow $g$ over $(0 \to 1)$ and some arrow $h$ over $(1 \to 2)$. Any candidate labeling $\lbl W : \Tot{}{\und W}  \to \iC$ will have to satisfy that
\begin{equation}
    \lbl W f = \lbl {R_{12}} h \circ \lbl {R_{01}} g
\end{equation}
(thus the labeling $\lbl W$ must be unique if it exists). Provided that this composition is independent of our choice of $g$ and $h$, the value $\lbl W f$ is well-defined, thus completing
the definition of $\lbl W$.

Abbreviate the composite morphism $\lbl {R_{12}} g \circ \lbl {R_{01}} h$ in $\iC$ by $\mathsf{lbl}(g,h)$. Given another pair of arrows $g'$ and $h'$ decomposing $f$, we need to verify that $\mathsf{lbl}(g,h) = \mathsf{lbl}(g',h')$. To prove this, we use truss induction (recall, in particular, the scaffold order $\fleq$ on sections, see \autoref{lem:section-induction}). The arrows $g$ and $h$ (resp.\ $g'$ and $h'$) are the spine vectors of a section simplex $K$ (resp.\ $K'$) in the bundle $\fp{}{\und W}$. Assume that $K \fleq K'$ (otherwise switch $K$ and $K'$ in the following). Let $K = K_0, K_1, ..., K_k = K'$ be a chain of successor sections starting at $K$ and ending in $K'$. Arguing inductively, we can assume $k = 1$, that is, $K' = K_1 = \successor(K)$. By \autoref{lem:spacer-induction} there is a spacer $L$ with lower boundary $\partial_- L = K$ and upper boundary $\partial_+ L = K'$. Note that the spine $L(2 \to 3) \circ L(1 \to 2) \circ L(0 \to 1)$ of $L$ composes both to the spine $K(1 \to 2) \circ K(0 \to 1) = h \circ g$ of $K$ and to the spine $K'(1 \to 2) \circ K'(0 \to 1) = h' \circ g'$ of $K'$. Functoriality of the labelings $\lbl {R_{12}}$ of and $\lbl {R_{12}}$ now implies
\begin{equation}
    \mathsf{lbl}(g,h) = \lbl {R_{12}} \La(2 \to 3) \circ \lbl {S} \La(1 \to 2) \circ \lbl {R_{01}} \La(0 \to 1) = \mathsf{lbl}(g',h')
\end{equation}
(where $\lbl {S}$ is the restriction of both $\lbl {R_{12}}$ and $\lbl {R_{01}}$ to the fiber over $1 \in [2]$) as required.
\end{proof}

\nid The preceding result shows that labeled 1-truss bordisms have a well-defined notion of composition, allowing us to introduce the following category.

\begin{defn}[Labeled 1-truss bordisms] \label{defn:labeled-1-truss-bord-category} Given a category $\iC$, the \textbf{category of $\iC$-labeled 1-truss bordisms $\lttr 1 \iC$} is the category whose objects are $\iC$-labeled 1-trusses (see \autoref{defn:labeled-1-trusses});
 whose morphisms are $\iC$-labeled 1-truss bordisms (see \autoref{defn:labeled-1-truss-bordism}); and whose composition is as given in \autoref{defn:composition-witnesses-1}.
 \end{defn}

\nid Noteworthily, the definition of the category of labeled 1-truss bordisms generalizes by allowing labelings to take values in `\infty-categories' as well, which we record in the following remark (without using it later on).

\begin{rmk}[1-Truss bordism labeled in an $\infty$-category] \label{altdefn:infty-truss-bord} Let $\cC$ be a quasicategory (i.e.\ a simplicial set in which inner horn fillers exist, see \cite{boardman2006homotopy}). We define the quasicategory $\qlttr 1 \cC$ of `$\cC$-labeled 1-truss bordisms' to be the simplicial set whose $k$-simplices are given by tuples consisting of a `$k$-simplex of bordisms' $S : [k] \to \ttr 1$ together with a `$\cC$-labeling' of its total space now given by an $\infty$-functor $\lbl S : \Totb{}{S} \to \cC$. Note that, for a simplicial map $f : [l] \to [k]$, we obtain an $l$-simplex in $\qlttr 1 \cC$ as the tuple consisting $S \circ f : [l] \to \ttr 1$ and the labeling $\lbl S \circ \Tot{}{f} : \Totb{}{S \circ f} \to \cC$. To see that the resulting simplicial set $\qlttr 1 \cC$ is indeed a quasicategory, i.e.\ as inner horn fillers, one uses truss induction for bundles over general $k$-simplices (conceptually following the proof of our earlier \autoref{lem:comp-witnesses}, which used truss induction only over the $2$-simplex).
\end{rmk}

\begin{rmk}[Unlabeled 1-truss bordisms are trivially labeled] \label{rmk:trivial-labels-again} Unlabeled 1-truss bordisms can be thought of as having `trivial' labelings, i.e.\ labelings in the terminal category $\ast$: indeed, the functor $\lttr 1 \ast \to \ttr 1$ mapping $\ast$-labeled 1-trusses $T$ (resp.\ bordisms $R$) to their underlying 1-truss $\und T$ (resp.\ their underlying 1-truss bordism $\und R$) is an isomorphism of categories.
\end{rmk}

The construction of $\lttr 1 \iC$ is, in fact, functorial in the category $\iC$, yielding the `labeled 1-truss bordism' functor defined below. Its behavior on functors of categories can be described as follows.

\begin{defn}[Relabeling along a functor] Consider a functor $F : \iC \to \iD$ between categories $\iC$ and $\iD$. The \textbf{$F$-relabeling functor} between categories of labeled 1-truss bordisms
    \begin{equation}
        \lttr 1 \iF : \lttr 1 \iC \to \lttr 1 \iD
    \end{equation}
is defined to map a $\iC$-labeled 1-truss $T$ to the $\iD$-labeled 1-truss with underlying truss $\und T$ and labeling $F \circ \lbl T$, and a $\iC$-labeled 1-truss bordism $R$ to the $\iD$-labeled 1-truss bordism with underlying 1-truss bordism $\und R$ and labeling $F \circ \lbl R$.
\end{defn}

\begin{term}[The label forgetting functor] \label{constr:forgetful-functor} Let $\iC$ be a category. Relabeling by the terminal functor $\iC \to \ast$ yields the functor
    \begin{equation}
        \lttr 1 {\iC \to \ast} : \lttr 1 {\iC} \to \lttr 1 {\ast} \iso \ttr 1
    \end{equation}
which we call the `forgetful functor' and usually abbreviate by an underline $(\und {-})$.
\end{term}

\begin{defn}[The labeled 1-truss bordism functor] \label{constr:lebun-functor} The \textbf{labeled 1-truss bordism functor} is the endofunctor
    \begin{equation}
        \lttr 1 - : \Cat \to \Cat
    \end{equation}
that takes a category $\iC$ to the category $\lttr 1 \iC$ of $\iC$-labeled 1-truss bordisms and a functor $\iF : \iC \to \iD$ to the $F$-relabeling $\lttr 1 \iF : \lttr 1 \iC \to \lttr 1 \iD$.
\end{defn}

\pauseae

We end our discussion of labeled 1-truss bordisms by highlighting an alternative definition of the notion in abstract categorical terms, namely, using the language of double categories. The reader without an inclination for abstract categorical structures may readily skip ahead to the next section.

\begin{term}[Vertical comma categories] For a pseudofunctor $H : \iD \to \Prof$ from a category $\iD$ into the bicategory of profunctors (which is the `horizontal category' of the double category of profunctors, see \cite[\S 3.1]{grandis1999limits}), and a category $\iC$, define the `vertical comma bicategory' $H_{\vslash \iC}$ as follows: objects are tuples $(d \in \iD,\iF : H(d) \to \iC)$; 1-morphisms $(d,\iF) \to (d',\iF')$ are tuples $(f : d \to d', \alpha : H(f) \Rightarrow \Hom_{\iC}(\iF-,\iF'-))$; 2-morphisms $(f,\alpha) \Rightarrow (f',\alpha')$ are natural isomorphisms $\beta : H(f) \Rightarrow H(f')$ such that $\alpha'\beta = \alpha$.\footnote{More generally, we could allow non-invertible 2-morphisms in $H_{\vslash \iC}$. However, for our purposes here, we want to restrict our attention to natural isomorphisms. }
\end{term}

\nid To apply the notion of vertical comma categories to the case of 1-truss bordism, we next define a functor $\iota : \ttr 1 \to \Prof$. Recall the fully faithful inclusion $\Bool \into \SetCat$, mapping 0 to the empty set $\emptyset$, and 1 to the singleton set $\Set{\ast}$.

\begin{term}[Inclusion of truss bordisms into ordinary profunctors] We define the `bordisms-as-profunctors' pseudofunctor
    \begin{equation}
        \iota : \ttr 1 \to \Prof
    \end{equation}
from category 1-truss bordisms to the bicategory of profunctors $\Prof$, by mapping 1-trusses $T$ to the category $(T,\eleq)$ and 1-truss bordisms $R$ to the post-composition of their underlying Boolean profunctor with the inclusion $\Bool \into \SetCat$.
\end{term}

\nid The fact that $\iota$ is a well-defined pseudofunctor is addressed in the following remark; we highlight that this claim relies on the rather special properties of 1-truss bordisms and a similar definition fails for general Boolean profunctors.

\begin{rmk}[The specialness of 1-truss bordisms] Ordinary profunctors compose by coends (see \cite[\S 5]{loregian2021co}). Importantly, given general Boolean profunctors $R : X \proto{} Y$ and $S : Y \proto{} Z$ between preorders it need \textit{not} be the case that $\iota (S \circ R) \iso (\iota S) \circ (\iota R)$; that is, the composition of Boolean profunctors need not coincide with their composition as ordinary profunctors. However, in the case of 1-truss bordisms $R$ and $S$, the isomorphism $\iota (S \circ R) \iso (\iota S) \circ (\iota R)$ always uniquely exists: this follows by explicitly evaluating the colimit defining the profunctor composite $(\iota S) \circ (\iota R)$ by following the arguments in the proof of \autoref{lem:comp-witnesses}.
\end{rmk}

\nid Combining the above notions, we now reach the following alternative description of labeled 1-truss bordisms.

\begin{obs}[Labeled 1-truss bordisms] \label{altdefn:lebun} The category of $\iC$-labeled 1-truss bordisms is equivalent to the vertical comma bicategory $\iota_{\vslash \iC}$ as defined above. In particular, note that the only $2$-morphisms in the bicategory $\iota_{\vslash \iC}$ are identities, and thus $\iota_{\vslash \iC}$ is in fact an ordinary 1-category.
\end{obs}

\nid We omit a verification of this observation, which includes showing that $\iota_{\vslash \iC}$ has `trivial' $2$-categorical structure (i.e.\ is an ordinary 1-category).

\subsubsecunnum{The definition of labeled 1-truss bundles} We next generalize 1-truss bundles to the labeled case; similar to the previous definitions of labeled 1-trusses and their bordism, this now endows the total poset of a 1-truss bundle with a labeling functor as follows.

\begin{defn}[Labeled 1-truss bundles] \label{defn:labeled-1-truss-bundle} Let $\iC$ be a `labeling' category and $B$ a `base' poset. A \textbf{$\iC$-labeled 1-truss bundle} $p$ over $B$ is a pair consisting of an `underlying' 1-truss bundle $p : T \to B$ over $B$ together with a `labeling' functor $\lbl p : (T,\eleq) \to \iC$.
\end{defn}

\nid In other words, a $\iC$-labeled 1-truss bundle $p \equiv (\und p,\lbl p)$ is a span of the form
    \begin{equation}
        \begin{tikzcd}
            \iC & T \arrow[r, "p"] \arrow[l, "\lbl p"'] & {B} \quad.
        \end{tikzcd}
    \end{equation}
Labeled 1-truss bundles support the following notion of maps, which mirrors the definition of maps of their unlabeled counterpart (see \autoref{defn:1-truss-map}).

\begin{defn}[Maps of labeled 1-truss bundles] \label{defn:labeled-1-truss-bundle-maps} For categories $\iC$ and $\iD$, let $p$ be a $\iC$-labeled 1-truss bundle, and $q$ an $\iD$-labeled 1-truss bundle. A \textbf{map of $\iC$-labeled 1-truss bundles} $F : p \to q$ is a pair consisting of an `underlying' 1-truss bundle map $\und F : \und p \to \und q$, as well as a `label category' functor $\lbl F : \iC \to \iD$ such that $\lbl F \circ \lbl p = \lbl q \circ \und F$.
\end{defn}

\nid Equivalently, one may think of labeled 1-truss bundle maps $F : p \to q$ as maps of spans $p \to q$, i.e. commuting diagrams
\begin{equation}
\begin{tikzcd}
\iC \arrow[d, "\lbl F"'] & T \arrow[r, "{\und p}"] \arrow[l, "\lbl p"'] \arrow[d, "\und F"] & B \arrow[d] \\
\iD                 & S \arrow[r, "{\und q}"'] \arrow[l, "\lbl q"]                & C
\end{tikzcd}
\end{equation}
in which $\und F$ is a 1-truss bundle map $\und p \to \und q$.

\begin{term}[Singular, regular, and balanced labeled bundle maps] A labeled 1-truss bundle map $F$ is said to be `singular', `regular' or `balanced' if its underlying 1-truss bundle map $\und F$ is
\end{term}

\nid Maps of labeled 1-truss bundles compose in the obvious way, yielding the following category.

\begin{notn}[The category of labeled 1-truss bundle] The category of labeled 1-truss bundles and their maps will be denoted by $\trussbunlbl 1$. The wide subcategories of regular resp.\ singular resp.\ balanced maps will be denoted by superscripts `$\mathsf{r}$' resp.\  `$\mathsf{s}$' resp.\  `$\mathsf{rs}$' as usual.
\end{notn}

\begin{rmk}[All 1-truss bundles are trivially labeled] \label{rmk:trivial-labels} Every 1-truss bundle $p : T \to B$ admits a trivial labeling $T \to \ast$ in the terminal category $\ast$. A map between trivially labeled 1-truss bundles is simply a map between their underlying 1-truss bundles. This provides a fully faithful inclusion
\begin{equation}
    \trussbun 1 \into \trussbunlbl 1
\end{equation}
and in this sense the notion of 1-truss bundles is properly generalized by the notion of labeled 1-truss bundles.
\end{rmk}

\begin{term}[Label and base preserving maps] A labeled 1-truss bundle map $(F,G,H)$ is called `label preserving' if $H = \id_\iC$ is the identity on the labeling category $\iC$, and `base preserving' if $F = \id_B$ is the identity on the base poset $B$.
\end{term}

\begin{defn}[Restrictions of labeled 1-truss bundles] Given a $\iC$-labeled 1-truss bundle $p = (\und p : T \to B, \lbl p : T \to \iC)$ and a subposet $C \into B$, then the \textbf{restriction} $\rest p C$ is the $\iC$-labeled 1-truss bundle $(\rest {\und p} C : \rest T C \to C, \lbl p : \rest T C \to \iC)$.
\end{defn}

\nid Note for any non-identity morphisms $f : [1] \into B$, we may think of the restriction $\rest p f$ as a $\iC$-labeled 1-truss bordism. In this sense, analogous to the unlabeled case, labeled 1-truss bordisms describe `fiber transitions' in labeled 1-truss bundles.

\subsubsecunnum{Classification and totalization for labeled 1-truss bundles} We construct a classification-totalization equivalence between $\iC$-labeled 1-truss bundles and their concordances and functors into the category of $\iC$-labeled 1-truss bordisms (analogous to the classification-totalization equivalence for unlabeled 1-truss bundles, see \autoref{constr:class-tot-1}). We start with the construction of `classifying functors'.

\begin{constr}[Classifying functors of labeled 1-truss bundles] \label{constr:classification-2} Consider a $\iC$-labeled 1-truss bundle $p$ over a poset $B$. Its \textbf{classifying functor} $\fcl {} p$ is the functor $B \to \lttr 1 {\iC}$ that maps an object $x : [0] \into B$ to the $\iC$-labeled 1-truss $\rest p x$ and a non-identity arrow $f : [1] \into B$ to the $\iC$-labeled 1-truss bordism $\rest p f$. The fact that this constructs a functor can be checked using the uniqueness of composition witnesses (see \autoref{lem:comp-witnesses}).
\end{constr}

Conversely, one constructs `total bundles' as follows.

\begin{constr}[Totalizing bundles of labeled 1-truss bundles] \label{constr:totalization-2}
    For a poset $B$ and a category $\iC$, consider a functor $\iF: B \to \lttr 1 \iC$. We define its \textbf{total bundle} $\fp{}{\iF} \equiv (\fp{}{\iF}, \flblo \iF)$ to be the $\iC$-labeled 1-truss bundle
    \begin{itemize}
        \item[-] whose underlying bundle $\fp{}{\iF}$ is the total 1-truss bundle $\fp{}{\und \iF}$ (of the composition of $\iF$ with the label forgetting functor $(\und {-}) : \lttr 1 \iC \to \ttr 1$),
        \item[-] and whose labeling functor $\flblo \iF : \Totb{}{\und \iF} \to \iC$ is determined on fibers over objects $x \in B$ as the labeling functor $\lbl {Fx} : \Totb{}{\und Fx} \to \iC$ (of the labeled 1-truss $Fx$) and on fibers over non-identity arrows $f : [1] \to B$ as the labeling functor $\lbl {Ff} : \Totb{}{\und Ff} \to \iC$ (of the labeled 1-truss bordism $Ff$).
    \end{itemize}
The fact that $\flblo \iF$ is indeed a functor can be seen to follow again using the construction of composition witness of labeled 1-truss bordisms (see \autoref{lem:comp-witnesses}).
\end{constr}

In order to promote the preceding constructions to functors we introduce a notion labeled 1-truss bundle concordances as follows.

\begin{defn}[Concordances of labeled 1-truss bordisms] For a poset $B$ and a `labeling' category $\iC$, a \textbf{$\iC$-labeled 1-truss bundle concordance} $u : p \Rightarrow q$ between $\iC$-labeled 1-truss bundles $p$ and $q$ over $B$ is a $\iC$-labeled 1-truss bundle over $B \times [1]$ such that $\rest u {B \times \{0\}} = p$ and $\rest u {B \times \{1\}} = q$.
\end{defn}

\begin{rmk}[Invertible concordances need not be unique] \label{rmk:1-trusses-skeletally-concordance-labeled} In contrast to \autoref{rmk:1-trusses-skeletally-concordance}, invertible labeled bundle concordances between fixed labeled 1-truss bundles need not be unique if they exist (since the labeling category $\iC$ may have non-trivial automorphisms).
\end{rmk}

Note that, given a $\iC$-labeled 1-truss bundle concordance $u : p \Rightarrow q$, its classifying functor $\fcl{}{u}$ is equivalently a natural transformation $\fcl{}{p} \Rightarrow \fcl{}{q} : B \to \lttr 1 \iC$; we often refer to this natural transformation as the `classifying natural transformation' of $u$.

\begin{defn}[Categories of 1-truss bundle concordances] For a poset $B$ and a category $\iC$, define the \textbf{category of concordances of $\iC$-labeled 1-truss bundles over $B$}, denoted by $\trussconclbl 1(B,\iC)$, to have $\iC$-labeled 1-truss bundles over $B$ as objects, and $\iC$-labeled 1-truss bundle concordances as morphisms. Composition of two such concordances $u : p \Rightarrow q$ and $v : q \Rightarrow r$ is determined by the condition that $\fcl{}{v \circ u} = \fcl{}{v} \circ \fcl{}{u}$.\footnote{Explicitly, working skeletally with underlying 1-truss bundles (that is, identifying unlabeled 1-truss bundles up to invertible concordances, see \autoref{rmk:1-trusses-skeletally-concordance}), we may define $v \circ u = \fp {} {\fcl {} v \circ \fcl {} u}$.}
\end{defn}

\nid Classification and totalization then organize into equivalences of categories as follows.

\begin{constr}[Classification and totalization functors] \label{constr:class-tot-lbl-1-truss}
Given a poset, we define the \textbf{classification functor}
\[
    \fcl{}{-} : \trussconclbl 1(B,\iC) \to \Fun(B,\lttr 1 \iC)
\]
to map labeled 1-truss bundles $p : T \to B$ to their classifying functor $\fcl {} p$, and labeled concordances $u : p \Rightarrow q$ to their classifying natural transformations $\fcl{}u : \fcl{}p \Rightarrow \fcl{}q$. Conversely, we define the \textbf{totalization functor}
\[
    \fp{}{-} : \Fun(B,\lttr 1 \iC) \to \trussconc 1(B,\iC)
\]
to take functors $\iF: B \to \lttr 1 \iC$ to their total 1-truss bundles $\fp {} \iF$, and similarly natural transformations $\alpha : B \times [1] \to \ttr 1$ to their total 1-truss bundles $\fp {} {\alpha}$.
\end{constr}

\begin{obs}[Classification and totalization are inverse] Classification of labeled 1-truss bundles and totalization of functors into labeled 1-truss bordisms provide an equivalence of categories.
\end{obs}

\subsubsecunnum{Pullback, dualization, and suspension of labeled 1-truss bundles}

We introduce pullbacks, duals, and suspensions of labeled 1-truss bundles. The first, in analogy to the unlabeled case (see \autoref{constr:1-truss-bundle-pullback}), can be defined as follows.

\begin{defn}[Pullbacks of labeled 1-truss bundles] \label{constr:labeled-1-truss-bundle-pullback} Given a $\iC$-labeled 1-truss bundle $p$ over a poset $B$ and a poset map $F : C \to B$, the \textbf{pullback $F^*p$ of $p$ along $F$} is the $\iC$-labeled 1-truss bundle with classifying functor $\fcl{}{F^*p} = \fcl{}{p} \circ F$.
\end{defn}

\nid In other words, $F^*p \equiv (\und {F^*p},\lbl {F^*p})$ can be thought of as the span obtained as the upper row in the following diagram
\begin{equation}
\begin{tikzcd}[row sep=5pt]
 &[+25pt] F^*T \arrow[ddr, phantom, "\lrcorner" , very near start, color=black] \arrow[ld, "\lbl p \circ \Tot{}F"' pos=.1] \arrow[dd, "\Tot{}F"'] \arrow[r, "F^*\und p"] & C \arrow[dd, "F"] \\
\iC & & \\
 & T \arrow[lu, "\lbl p" pos=.4] \arrow[r, "\und p"] & B
\end{tikzcd}
\end{equation}
where $\Tot{}{F}$ is the canonical pullback map of 1-truss bundles $\Tot{}F : F^*\und p \to \und p$.

\begin{rmk}[Restrictions] \label{notn:labeled-1truss-bundle-restriction} Pullbacks of labeled 1-truss bundles $p$ over $B$ along subposets $F : C \into B$ recover restrictions, that is, $F^*p = \rest p C$.
\end{rmk}

Next\nid  we consider dualization of labeled 1-truss bundles.

\begin{defn}[Duals of labeled 1-truss bundles] \label{constr:labeled-1-truss-bundle-dualization} Given a $\iC$-labeled 1-truss bundle $p \equiv (\und p,\lbl p)$, its \textbf{dual bundle} $p^\dagger$ is the $\iC\op$-labeled 1-truss bundle whose underlying bundle $\und {p^\dagger}$ is the dual $(\und p)^\dagger$ of the underlying bundle $\und p$ (see \autoref{constr:1-truss-bundle-dualization}) and whose labeling $\lbl {p^\dagger}$ is the opposite labeling $(\lbl p)\op$.
\end{defn}

\begin{defn}[Duals of labeled dual 1-truss bundle maps] Given a labeled 1-truss bundle map $F : p \to q$ one defines the \textbf{dual bundle map} $F^\dagger : p^\dagger \to q^\dagger$ to have underlying 1-truss bundle map is $\und{F^\dagger} = \und F\op$ (i.e. $\und F^\dagger$) and labeling functors $\lbl {F^\dagger} = (\lbl F)\op$.
\end{defn}

\begin{defn}[Duals of labeled 1-truss bundle concordances] Given a labeled 1-truss bundle concordance $u : p \Rightarrow q$ one constructs the \textbf{dual bundle concordance} $u^\dagger : q^\dagger \Rightarrow p^\dagger$ by setting $u^\dagger$ to be the dual bundle of $u$.
\end{defn}

\nid The preceding definitions now yield the following functors.

\begin{obs}[Dualization functors on 1-truss bundles] The preceding definitions construct a dualization functor of labeled $1$-truss bundles as follows:
    \[
        \dagger : \trussbunlbl 1 \iso \trussbunlbl 1.
    \]
Similarly, for a category $\iC$ and a poset $B$, we obtain a dualization functor of $\iC$-labeled $1$-truss bundle concordances as follows:
    \[
        \dagger : \trussconclbl 1(B,\iC) \iso \trussconclbl 1(B,\iC\op)\op .
    \]
If we set $B = \ast$, then this further specializes to an isomorphism of $\iC$-labeled $1$-truss bordisms:
    \[
        \dagger : \lttr 1 \iC \iso (\lttr 1 {\iC\op})\op . \qedhere
    \]
\end{obs}

\begin{rmk}[Dualization via classifying functors] Using the classification and totalization equivalence for labeled 1-truss bundles, given a $\iC$-labeled 1-truss bundle $p$, its dual $p^\dagger$ may equivalently be defined to have classifying functor $\fcl{}{p^\dagger}$ equal to the composite of $(\fcl{}{p})\op$ with the isomorphism $(\lttr 1 \iC)\op \iso \lttr 1 {\iC\op}$.
\end{rmk}

\begin{rmk}[Suspensions in the labeled case] \label{rmk:labeled-suspension-1-trs} Recall the definition of suspensions of 1-truss bundles from \autoref{constr:1-truss-bundle-suspension}. This again has a labeled analog, subject to the additional assumption that $\iC$ has both initial and terminal objects, which can be used as labels for the newly adjoined initial and terminal objects in the total poset of the suspension---we leave details to the reader.
\end{rmk}

\subsection{$n$-Truss bordisms and $n$-truss bundles} \label{ssec:n-trusses}

We now discuss the combinatorial theory of $n$-trusses, as well as their bordisms and their bundles. Recall from the introduction to this section, which defined $n$-trusses as tower of 1-truss bundles ending in the $0$-simplex $[0]$. We saw the definition of $n$-trusses already in \autoref{defn:n-trusses}; for a more uniform description of truss theory, it will be useful to generalize this to a notion of `labeled $n$-trusses' as follows.

\begin{defn}[Labeled $n$-trusses] \label{defn:n-trusses-labeled} Given a category $\iC$, an \textbf{$\iC$-labeled $n$-truss $T$} is a pair consisting of an `underlying' $n$-truss $\und T$ together with a `labeling' functor $\lbl T : T_n \to \iC$ (where $T_n$ is the total poset of $\und T$).
\end{defn}

\nid Unpacking the underlying $n$-truss $\und T$ as a tower of 1-truss bundles $p_i : T_i \to T_{i-1}$ we may equivalently consider an $n$-truss to be a diagram of the form
    \begin{equation}
        \iC \xot {\lbl T} T_n \xto {p_n} T_{n-1} \xto {p_{n-1}} ... \xto {p_2} T_1 \xto {p_1} T_0 = [0]
    \end{equation}
where $p_i$ are 1-truss bundles whose base poset is the total poset of $p_{i-1}$, and where $\lbl T : T_n \to \iC$ is a functor from the total poset $T_n$ of $p_n$ to the category $\iC$.

The goal of this section will be to study notions of (labeled) $n$-truss bordisms, $n$-truss bundles as well as their maps and concordances, which generalize previously introduced notions to `higher dimensions'. This will involve an interesting inductive iteration of ideas from previous sections.

\begin{rmk}[On towers and sequences] We use the terms `towers of maps' and `sequences of maps' largely synonymously, but prefer `tower' when maps have `bundle character'.
\end{rmk}

\subsubsecunnum{The definition $n$-truss bordisms} We introduce $n$-truss bordisms, which, in analogy to the case of 1-trusses of 1-truss bordisms, will describe `fiber transitions' in bundles whose fibers are $n$-trusses.

\begin{defn}[$n$-Truss bordisms] \label{defn:n-truss-bordism} An \textbf{$n$-truss bordism $R$} is a tower of 1-truss bundles
    \begin{equation}
        R_n \xto {p_n} R_{n-1} \xto {p_{n-1}} ... \xto {p_2} R_1 \xto {p_1} R_0 = [1] \quad .
    \end{equation}
    We refer to $R_n$ as the \textbf{total poset} of $R$.
\end{defn}

\nid Similarly, we obtain the following labeled analog of the notion.

\begin{defn}[Labeled $n$-truss bordism] \label{defn:labeled-n-truss-bordism} For a category $\iC$, an \textbf{$\iC$-labeled $n$-truss bordism} $R$ is a pair consisting of an `underlying' $n$-truss bordism $\und R$ and a `labeling' functor $\lbl R : R_n \to \iC$ (where $R_n$ is the total poset of $R$).
\end{defn}

\nid Unpacking the tower $\und R$, we may equivalently think of $R$ as a diagram of the form
\begin{equation}
        \iC \xot {\lbl R} R_n \xto {p_n} R_{n-1} \xto {p_{n-1}} ... \xto {p_2} T_1 \xto {p_1} R_0 = [1]
    \end{equation}
    where $p_i$ are 1-truss bundles, and $\lbl R : (R_n,\eleq) \to \iC$ is a functor. Clearly, \autoref{defn:labeled-n-truss-bordism} specializes to \autoref{defn:n-truss-bordism} if we label in the terminal category $\iC = \ast$, and we shall therefore work mainly in the more general case of labeled $n$-truss bordisms. Our main goal will be understanding how labeled $n$-truss bordisms organize (as morphisms) into a category.

\begin{term}[Domain and codomain of a labeled $n$-truss bordism] Given a $\iC$-labeled $n$-truss bordism $R \equiv (\und R, \lbl R)$, the `domain' $\dom(R)$ is the $\iC$-labeled $n$-truss $T$ whose underlying $n$-truss $\und T$ is obtained by (iteratively) restricting the tower of bundles $\und R$ to $0 \in [1]$, and whose labeling $\lbl T$ is the restriction of $\lbl R$ to total poset $T_n$ of $\und T$; in other words, $T$ is the determined by the upper row in the following diagram of 1-truss bundle pullbacks (all of which are restrictions, see \autoref{rmk:1-truss-bundle-pullback-restriction})
\begin{equation} \small
        \begin{tikzcd}
 & T_n \arrow[ddr, phantom, "\lrcorner" , very near start, color=black] \arrow[dl, "\lbl {T}"' pos=.2] \arrow[r] \arrow[dd, hook'] & T_{n-1} \arrow[ddr, phantom, "\lrcorner" , very near start, color=black] \arrow[r] \arrow[dd, hook'] & \cdots \arrow[r] \arrow[dd, "\cdots", phantom] & T_1 \arrow[ddr, phantom, "\lrcorner" , very near start, color=black] \arrow[r] \arrow[dd, hook'] & {[0]} \arrow[dd, hook', "0"] \\[-17pt]
            \iC & & & & & \\[-17pt]
& R_n \arrow[ul, "\lbl R" pos=.2] \arrow[r, "p_n"'] & R_{n-1} \arrow[r, "p_{n-1}"'] & \cdots \arrow[r, "p_2"'] & R_1 \arrow[r, "p_1"'] & {[1]}
\end{tikzcd}
\end{equation}
Similarly, one defines the `codomain' $\cod(R)$ by restricting $R$ to the object $1 \in [1]$.
\end{term}

\nid The following notion of composition immediately generalizes \autoref{defn:composition-witnesses-1}.

\begin{defn}[Composition of labeled $n$-truss bordisms] \label{defn:composition-witnesses-n} Given two $\iC$-labeled truss bordisms $R_{01}$ and $R_{12}$ which are composable (that is, $\cod(R_{01}) = \dom(R_{12})$), a \textbf{composition candidate} $R_{02}$ of $R_{12}$ and $R_{01}$ is a $\iC$-labeled truss bordism for which there exists a diagram $W$ of the form
    \begin{equation}
    \iC \xot {\lbl W} W_n \xto {p_n} W_{n-1} \xto {p_{n-1}} ... \xto {p_2} W_1 \xto {p_1} W_0 = [2]
    \end{equation}
where $p_i$ are 1-truss bundles and $\lbl W : R_n \to \iC$ is a functor, such that $\lbl W$ restricts on each non-identity arrow $(i \to j)$ in $[2]$ to the $\iC$-labeled $n$-truss bordism $R_{ij}$. In this case, we call the above diagram a \textbf{composition witness} of $R_{01}$ and $R_{12}$.
\end{defn}

\begin{lem}[Existence and uniqueness of $n$-truss bordism composition witnesses] \label{lem:comp-witnesses-n} Composition candidates and witnesses of labeled $n$-truss bordisms exist uniquely.
\end{lem}

\nid A direct proof of this statement can be given, but since the statement will be a corollary of subsequent constructions, we will defer a proof until \autoref{rmk:composition-witnesses-revisited}. The preceding result allows us to define the following category.

\begin{defn}[The category of $\iC$-labeled $n$-truss bordisms] \label{defn:labeled-n-truss-bordism-category} Given a category $\iC$, the \textbf{category $\nlttr {n} \iC$ of $\iC$-labeled $n$-truss bordisms} is the category whose objects are $\iC$-labeled $n$-trusses and whose morphisms are $\iC$-labeled $n$-truss bordisms endowed with a notion of composition induced by the existence and uniqueness of composition candidates.
\end{defn}

\nid Importantly, the category of $\iC$-labeled 1-truss bordisms admits an equivalent formulation, which highlights its `fundamentally iterative' nature. We introduce this formulation as follows.

\begin{defn}[The labeled $n$-truss bordism functor] The $n$-fold composite $\lttr 1 - \circ \lttr 1 - \circ ... \circ \lttr 1 - : \Cat \to \Cat$ of the  labeled 1-truss bordism functor $\lttr 1 - : \Cat \to \Cat$ is called the \textbf{labeled $n$-truss bordism functor}, and denoted by $\lttr {n} -$. If $n = 0$, we take this definition to mean $\lttr {0} - = \id$.
\end{defn}

\nid Applying the labeled $n$-truss bordism functor to a specific category $\iC \in \Cat$, we obtain the following.

\begin{defn}[The `iterative' category of $\iC$-labeled $n$-truss bordisms] \label{defn:labeled-n-truss-bordism-endo} Given a category $\iC$, the \textbf{`iterative' category of $\iC$-labeled $n$-truss bordisms} is the category $\lttr {n} \iC$.
\end{defn}

\begin{rmk}[Unwinding the iteration] Recall that the functor $\lttr 1 - : \Cat \to \Cat$, takes a category and maps it to the category of 1-truss bordisms labeled in that category. On the one hand, the `iterative' category of $\iC$-labeled $n$-truss bordisms can be obtained as the category
\[
    \lttr {n} \iC = \lttr {1} {\lttr {n-1} \iC}
\]
From this perspective, we can understand $\lttr {n} \iC$ as the category of 1-truss bordisms \emph{labeled} in the category of $\iC$-labeled $(n-1)$-truss bordisms. On the other hand, we may also write
\[
    \lttr {n} \iC = \lttr {n-1} {\lttr {1} \iC}
\]
This means we can also understand $\lttr {n} \iC$ as the category of $(n-1)$-truss bordisms \emph{labeled} in the category of $\iC$-labeled 1-truss bordisms. Moreover, setting the labeling category $\iC = \ast$ to be trivial, these observations specialize the case of unlabeled $n$-truss bordisms: they may either be understood as 1-truss bordisms labeled in $(n-1)$-truss bordisms, or as $(n-1)$-truss bordisms labeled in 1-truss bordisms.
\end{rmk}

\nid The fact that \autoref{defn:labeled-n-truss-bordism-category} and \autoref{defn:labeled-n-truss-bordism-endo} define the same category requires proof.

\begin{prop}[Equivalence of viewpoints on labeled $n$-truss bordism] \label{prop:labeled-n-truss-bordism} The categories $\nlttr {n} \iC$ and $\lttr {n} \iC$ are canonically equivalent.
\end{prop}
\nid Again, we defer a proof since the statement will be a corollary of subsequent constructions, and will be revisited in \autoref{rmk:equivalence-of-viewpoints-on-labeled-n-truss-bordisms}.

\subsubsecunnum{The definition of $n$-truss bundles}

We discuss $n$-truss bundles which generalize 1-truss bundles. We start, as before, with the unlabeled case.

\begin{defn}[$n$-Truss bundles] \label{defn:n-truss-bundle}
    An \textbf{$n$-truss bundle} $p$ over a poset $B$ is a tower of 1-truss bundles
\begin{equation}
    T_n \xto {p_n} T_{n-1} \xto {p_{n-1}} ~\cdots~ \xto{p_2} T_1 \xto{p_1} T_0 = B \quad .
\end{equation}
We call $T_n$ the \textbf{total poset} of $p$.
\end{defn}

\begin{defn}[Labeled $n$-truss bundles] \label{defn:n-truss-bundle-labeled}
    Given a category $\iC$, a \textbf{$\iC$-labeled $n$-truss bundle} $p$ over a poset $B$ is a pair consisting of an `underlying' $n$-truss bundles $\und p$ together with a `labeling' functor $\lbl p : T_n \to \iC$ (where $T_n$ is the total poset of $p_n$).
\end{defn}

A $\iC$-labeled $n$-truss bundle $p = (\und p, \lbl p)$ can equivalently be understood as a diagram of the form
    \begin{equation}
        \iC \xot {\lbl p} T_n \xto {p_n} T_{n-1} \xto {p_{n-1}} ~\cdots~ \xto{p_2} T_1 \xto{p_1} T_0 = B
    \end{equation}
    where $p_i$ are the 1-truss bundles in $\und p$ and $\lbl p : (T_n,\eleq) \to \iC$ is a functor.

\begin{rmk}[$n$-Trusses and $n$-truss bordisms as bundles] \label{rmk:trusses-as-bundles} Note that both $\iC$-labeled $n$-truss as well as $n$-truss bordisms are special instance of the preceding definition, which can be recover for base posets $B = [0]$ resp.\  $B = [1]$.
\end{rmk}

\begin{term}[Open and closed bundles] If, for all $i$, $i$th bundles $p_i$ of an $n$-truss bundle $p$ are open (respectively closed), then we will call $p$ itself `open' (respectively `closed').
\end{term}

\begin{defn}[Labeled $n$-truss bundle maps] \label{defn:n-truss-bun-map}
    Given labeled $n$-truss bundles $p$ and $q$, a \textbf{labeled $n$-truss bundle map} $F : p \to q$ consists of a commuting diagram
\begin{equation}
    \begin{tikzcd}
        \iC \arrow[d, "\lbl F"'] & T_n \arrow[l, "\lbl p"'] \arrow[r, "p_n"] \arrow[d, "F_n"'] & T_{n-1} \arrow[r, "p_{n-1}"] \arrow[d, "F_{n-1}"'] & \cdots \arrow[r, "p_2"] \arrow[d, "\cdots", phantom] & T_1 \arrow[r, "p_1"] \arrow[d, "F_1"] & T_0 \arrow[d, "F_0"] \\
        \iD & S_n \arrow[l, "\lbl q"] \arrow[r, "q_n"']                  & S_{n-1} \arrow[r, "q_{n-1}"']                      & \cdots \arrow[r, "q_2"']                             & S_1 \arrow[r, "q_1"']                 & S_0
\end{tikzcd}
\end{equation}
where, for $i > 0$, each $F_i : p_i \to q_i$ is a 1-truss bundle map (see \autoref{defn:1-truss-bundle-map}) and $\lbl F : \iC \to \iD$ is a functor, called the `label category' functor of $F$.
\end{defn}

\begin{term}[Label and base preserving] \label{term:label-and-base-preserving} An $n$-truss bundle map $F$ is `base preserving' if $F_0 = \id$, and `label preserving' if its label category functor $\lbl F = \id$ is the identity.
\end{term}

\begin{term}[Singular, regular, and balanced $n$-truss bundles maps] If for each $i$, the $i$th bundle map of an $n$-truss bundle map singular resp.\ regular resp.\ balanced (in the sense of \autoref{defn:sing-and-reg-1-truss-bun-map}) then we call the bundle map itself `singular', resp.\ `regular', resp.\ balanced.
\end{term}

\nid The preceding definition of maps of labeled $n$-truss bundle specializes both to the case `unlabeled' $n$-truss bundles as well as the case of `$n$-trusses'.

\begin{rmk}[The unlabeled case] For unlabeled $n$-truss bundles $p,q$ (which can be considered as labeled $n$-truss bundles labeled in the terminal category $\ast$) the labeling functor $\lbl F : \ast \to \ast$ of any labeled $n$-truss bundle map $F : p \to q$ must be trivial; in this case we call $F$ simply an `$n$-truss bundle map'. Moreover, every labeled $n$-truss bundle map $F : p \to q$ induces an $n$-truss bundle map $\und F : \und p \to \und q$ of underlying $n$-truss bundles, by setting $\und F_i = F_i$. We call $\und F$ the `underlying' $n$-truss bundle map of $F$.
\end{rmk}

\begin{rmk}[The case of $n$-trusses] For labeled $n$-trusses $T$ and $S$ (which can be regarded as labeled $n$-truss bundles over the trivial base poset $[0]$) we usually refer to labeled $n$-truss bundle maps $F : T \to S$ as `labeled $n$-truss maps'.
\end{rmk}

\nid Categories of $n$-trusses as well as $n$-truss bundles and their maps can now be introduced as follows.

\begin{notn}[$n$-Truss and $n$-truss bundle categories] We denote by $\trussbunlbl n$ the category of labeled $n$-truss bundles with the following commutative conventions for category names, decorations and superscripts.
    \begin{itemize}
        \item[-] We drop the suffix `$\mathsf{Bun}$' to indicate that we restrict to the full subcategory of labeled $n$-trusses (writing e.g\ $\trusslbl n \into \trussbunlbl n$).
        \item[-] We drop the prefix `$\mathsf{Lbl}$' to indicate that we restrict to a subcategory of unlabeled $n$-truss bundles (writing e.g.\ $\truss n \into \trusslbl n$).
        \item[-] To restrict to the full subcategories of `closed' resp.\ `open' $n$-truss bundles we add decorations `$\bar{\mathsf{T}}$' and $\mathring{\mathsf{T}}$ to the letter ${\mathsf{T}}$ (writing e.g.\ $\ctruss n \into \truss n$).
        \item[-] To restrict to the wide subcategories containing only singular or regular or balanced maps we add superscripts `$\mathsf{s}$' or `$\mathsf{r}$' or `$\mathsf{rs}$' (writing e.g.\ $\sctruss n \into \ctruss n$). \qedhere
    \end{itemize}
\end{notn}

\begin{defn}[Upper truncation] \label{defn:upper_truncations}
For $0 \leq k \leq n$, we define the \textbf{upper $k$-truncation functor}
    \begin{equation}
        (-)_{\geq k} : \trussbunlbl n \to \trussbunlbl {n-k}
    \end{equation}
    which maps a labeled $n$-truss bundle $p = (\und p, \lbl p)$ where $\und p = (p_n,p_{n-1},...,p_1)$ to the labeled $k$-truss bundles $p_{\geq k} = (\und p_{\geq k}, \lbl p)$ where $\und p_{\geq k} = (p_n,p_{n-1},...,p_{k+1})$.
\end{defn}

\begin{defn}[Lower truncation] \label{defn:lower_truncations}
For $k \leq n$, we define the \textbf{lower $k$-truncation functor}
    \begin{equation}
        (-)_{\leq k} : \trussbunlbl n \to \trussbun k
    \end{equation}
which maps a labeled $n$-truss bundle $p = (\und p, \lbl p)$ where $\und p = (p_n,p_{n-1},...,p_1)$ to the unlabeled $k$-truss bundle $p_{\leq k}$ consisting of 1-truss bundles $(p_k,p_{k-1},...,p_1)$.
\end{defn}

\begin{defn}[Labeled $n$-truss bundle restrictions] Given a $\iC$-labeled $n$-truss bundle $p$ over a poset $B$ with underlying bundle $\und p = (p_n,p_{n-1},...,p_1)$, as well as a subposet $C \into B$, we define the \textbf{restriction} $\rest p C$ of $p$ to $C$ to be the $\iC$-labeled $n$-truss bundle over $C$ obtained as the upper row in the following diagram
\begin{equation} \small
        \begin{tikzcd}
 & \rest {T_n} C \arrow[ddr, phantom, "\lrcorner" , very near start, color=black] \arrow[dl, "\lbl {\rest p C}"' pos=.2] \arrow[r, "\rest {p_n} C"] \arrow[dd, hook'] & \rest {T_{n-1}} C \arrow[ddr, phantom, "\lrcorner" , very near start, color=black] \arrow[r, "\rest {p_{n-1}} C"] \arrow[dd, hook'] & \cdots \arrow[r, "\rest {p_2} C"] \arrow[dd, "\cdots", phantom] & \rest {T_1}{C} \arrow[ddr, phantom, "\lrcorner" , very near start, color=black] \arrow[r, "\rest {p_1} C"] \arrow[dd, hook'] & {[0]} \arrow[dd, hook'] \\[-17pt]
            \iC & & & & & \\[-17pt]
& T_n \arrow[ul, "\lbl p" pos=.2] \arrow[r, "p_n"'] & T_{n-1} \arrow[r, "p_{n-1}"'] & \cdots \arrow[r, "p_2"'] & T_1 \arrow[r, "p_1"'] & {B}
\end{tikzcd}
\end{equation}
(note all pullbacks of 1-truss bundles in the diagram are in fact restrictions).
\end{defn}

\pauseae

We briefly remark on how to construct a class of `generating arrows' in the total poset of $n$-truss bundles, yielding a minimal set of arrows whose transitive closure recovers all non-identity arrows of the total poset. Note that, this class is sometimes referred to as the `covering relation' of a poset and can be defined for any poset; however, $n$-truss bundles admit an explicit construction of this class as follows.

\begin{constr}[Generating arrows] \label{rmk:generating-arrows} Given an $n$-truss bundle $p = (T_n \xto {p_n} T_{n-1} \xto {p_{n-1}} ... \xto {p_1} T_0 = B)$ over a base poset $B$, we inductively construct subsets of non-identity arrows $\mathrm{cov}(T_i) \subset \mor(T_i,\eleq)$, called the `generating arrows' of $T_i$.
    \begin{itemize}
        \item[-] If $i = 0$, define $\mathrm{cov}(B)$ to be the minimal subset of $\mor(B)$ that generates all non-identity arrows in $B$, i.e.\ its covering relation.
        \item[-] If $i > 0$, then $f \in \mathrm{cov}(T_i)$ if and only if one of the following holds.
            \begin{itemize}
                \item[$\cdot$] either $p_i f = \id$ (that is, $f$ lies in a fiber of $p_i$),
                \item[$\cdot$] or, $p_i f\in \mathrm{cov}(T_{i-1})$ and $f \in \reg(T_i) \cup \sing(T_i)$ (that is, $f$ lies over a generating arrow of $T_{i-1}$ and passes either between two regular or two singular objects, see \autoref{defn:sing-and-reg-1-truss-bun-map}).
            \end{itemize}
    \end{itemize}
The central observation about the generating arrow set $\mathrm{cov}(T_i)$ is that each non-identity arrow in $(T_i,\eleq)$ can be written as a composite of generating arrows (in other words, it is the covering relation of $T_i$). The observation will be particularly useful when illustrating $n$-trusses (and their bundles) in that we will often only depict their generating arrows.
\end{constr}

\begin{eg}[Generating arrows in $n$-truss bundles] In \autoref{fig:generating-arrows-in-truss-bundles} we illustrate a closed $2$-truss $T$ on the left, and depict generating arrows of $T_2$ on the right.
\begin{figure}[ht]
    \centering
    \def\svgwidth{1\columnwidth}
    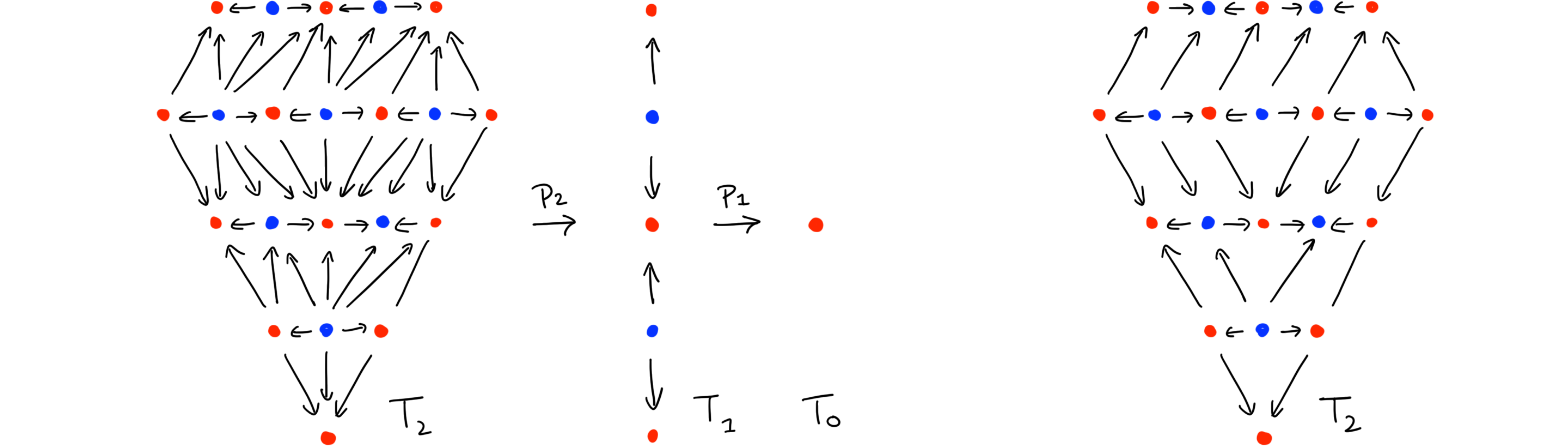

    \caption{A $2$-truss bundle and its generating arrows.}
    \label{fig:generating-arrows-in-truss-bundles}
\end{figure}
\end{eg}

\subsubsecunnum{Classification and totalization for labeled n-truss bundles} We discuss classification and totalization for labeled $n$-truss bundles, analogous to the case of labeled 1-truss bundles (see \autoref{constr:class-tot-lbl-1-truss}). Crucially, this now uses the iterative nature of the category $\lttr n \iC$ of $\iC$-labeled $n$-truss bordisms.

\begin{constr}[Classification of labeled $n$-truss bundles] \label{constr:grothendieck-labeled-n-dim} Let $p = (\und p, \lbl p)$ be an $\iC$-labeled $n$-truss bundle over a poset $B$ with $\und p = (p_n,p_{n-1},...,p_1)$. We construct the \textbf{classifying functor} $\fcl {} p : B \to \lttr n \iC$ of $p$ inductively as follows.
\begin{itemize}
    \item[-] Define a functor $\fcl n p := \lbl p : T_n \to \iC$.
    \item[-] Inductively in descending $i$, define a functor $\fcl {i-1} p : T_{i-1} \to \lttr {n-i} \iC$ to be the classifying map of the $\lttr {n-i+1} \iC$-labeled 1-truss bundles $(p_i, \fcl i p)$.
\end{itemize}
Then set $\fcl {} p = \fcl 0 p$.
\end{constr}

\begin{eg}[The labeled tower of a truss bundle] Here is an illustration of the form of a labeled tower of a $\iC$-labeled 3-truss-bundle $p$ with underlying 3-truss bundle $T_3 \xto {p_3} T_2 \xto{p_2} T_1 \xto {p_1} T_0$ and with label $\lbl p$, together with its classifying maps $\fcl i p$.
\begin{equation}
\begin{tikzcd}[column sep=35pt]
&[-27pt] \iC & T_3 \arrow[d, "p_3"] \arrow[l, "\fcl 3 T ~=~ \lbl p"']        \\
& \lttr 1 \iC & T_2 \arrow[d, "p_2"] \arrow[l, "\fcl 2 p"'] \\
    \lttr 1 {\lttr 1 \iC} \arrow[r,equal] & \lttr 2 \iC & T_1 \arrow[d, "p_1"] \arrow[l, "\fcl 1 p"']\\
\lttr 1 {\lttr 2 \iC} \arrow[r,equal] & \lttr 3 \iC & T_0 \arrow[l, "\fcl 0 p"']
\end{tikzcd}
\end{equation}
Each of the following subsets of the above diagrams determines the entire truss bundle $p$: the functor $\fcl 0 p$, or the bundle $p_1$ together with the functor $\fcl 1 p$, or the bundles $p_1$ and $p_2$ together with the functor $\fcl 2 p$, or the bundles $p_1$, $p_2$ and $p_3$ together with the functor $\fcl 3 p := \lbl p$.
\end{eg}

\nid Conversely, we may construct `total bundles' as follows. This simply inverts each inductive step in the preceding construction separately.

\begin{constr}[Totalization of labeled $n$-truss bundles] Given a functor $\iF : B \to \lttr n \iC$ we construct its \textbf{total $\iC$-labeled $n$-truss bundle} $\fp{} \iF = (\fpu{}{\iF}, \flblo \iF)$ inductively as follows.
\begin{itemize}
    \item[-] Define $\flbl 0 \iF = \iF$.
    \item[-] Inductively in ascending $i$, define $(\fp i \iF,~ \flbl i \iF)$ to be the $\lttr{n-i} \iC$-labeled total 1-truss bundle constructed from $\flbl {i-1} \iF : \Tot {i-1} \iF \to \lttr 1 {\lttr {n-i} \iC}$.
\end{itemize}
Then set $\fpu{}{\iF} = (\fp n \iF, \fp {n-1} \iF,..., \fp 1 \iF)$ and $\flblo \iF = \flbl n \iF$.
\end{constr}

To promote the preceding constructions to functors we once more introduce the notions of bundle concordances as follows.

\begin{defn}[Concordances of labeled $n$-truss bundles] For a poset $B$ and a `labeling' category $\iC$, a \textbf{$\iC$-labeled $n$-truss bundle concordance} $u : p \Rightarrow q$ between $\iC$-labeled $n$-truss bundles $p$ and $q$ over $B$ is a $\iC$-labeled $n$-truss bundle over $B \times [1]$ such that $\rest u {B \times \{0\}} = p$ and $\rest u {B \times \{1\}} = q$.
\end{defn}

\begin{rmk}[Uniqueness of $n$-truss bundle concordances] \label{rmk:n-trusses-skeletally-concordance} Fixing unlabeled $n$-truss bundles, any concordance between them must be unique if it exists (this follows by iteratively applying \autoref{rmk:1-trusses-skeletally-concordance}). In analogy to \autoref{rmk:1-trusses-skeletally-concordance-labeled}, this need not hold in the labeled case; but there is no harm in identifying labeled $n$-truss bundles up to concordance of their underlying (unlabeled) $n$-truss bundles.
\end{rmk}

Note that, given a $\iC$-labeled $n$-truss bundle concordance $u : p \Rightarrow q$, its classifying functor $\fcl{}{u}$ is equivalently a natural transformation $\fcl{}{p} \Rightarrow \fcl{}{q} : B \to \lttr n \iC$; we often refer to this natural transformation as the `classifying natural transformation' of $u$.

\begin{defn}[Categories of $n$-truss bundle concordances] For a poset $B$ and a category $\iC$, define the \textbf{category of concordances of $\iC$-labeled $n$-truss bundles over $B$}, denoted by $\trussconclbl 1(B,\iC)$, to have $\iC$-labeled $n$-truss bundles over $B$ as objects, and $\iC$-labeled $n$-truss bundle concordances as morphisms. Composition of two such concordances $u : p \Rightarrow q$ and $v : q \Rightarrow r$ is determined by the condition that $\fcl{}{v \circ u} = \fcl{}{v} \circ \fcl{}{u}$.\footnote{Explicitly, working skeletally with underlying $n$-truss bundles (that is, identifying unlabeled $n$-truss bundles up to invertible concordances, see \autoref{rmk:n-trusses-skeletally-concordance}), we may define $v \circ u = \fp {} {\fcl {} v \circ \fcl {} u}$.}
\end{defn}

\nid Classification and totalization then organize into equivalences of categories as follows.

\begin{constr}[Classification and totalization functors] \label{constr:class-tot-lbl-n-truss}
Given a poset, we define the \textbf{classification functor}
\[
    \fcl{}{-} : \trussconclbl n (B,\iC) \to \Fun(B,\lttr n \iC)
\]
to map labeled $n$-truss bundles $p$ to their classifying functor $\fcl {} p$, and labeled concordances $u : p \Rightarrow q$ to their classifying natural transformations $\fcl{}u : \fcl{}p \Rightarrow \fcl{}q$. Conversely, we define the \textbf{totalization functor}
\[
    \fp{}{-} : \Fun(B,\lttr 1 \iC) \to \trussconc 1(B,\iC)
\]
to take functors $\iF: B \to \lttr 1 \iC$ to their total $n$-truss bundles $\fp {} \iF$, and similarly natural transformations $\alpha : B \times [1] \to \ttr 1$ to their total $n$-truss bundles $\fp {} {\alpha}$.
\end{constr}

\begin{obs}[Classification and totalization are inverse] Classification of labeled $n$-truss bundles and totalization of functors into labeled $n$-truss bordisms provide an equivalence of categories.
\end{obs}

The correspondence between $n$-truss bundles and their classifying functors into the category $\lttr n \iC$ enables a simple proof of our earlier claim that the `non-iterative' and `iterative' definitions of labeled $n$-truss bordisms coincide. We record this by the following two observations.

\begin{obs}[Proof of \autoref{lem:comp-witnesses-n}] \label{rmk:composition-witnesses-revisited} In \autoref{lem:comp-witnesses-n} we claimed that the composition candidates and witnesses of composable $\iC$-labeled $n$-truss bordisms $R_{01}$ and $R_{12}$ always uniquely exist. Indeed, the composition witness (and thereby the candidate) can be easily constructed using totalization and classification as follows. First, composability implies $\fcl {} {\dom(R_{12})} = \fcl {} {\cod(R_{01})}$, and thus $\fcl {} {R_{01}} : [1] \to \lttr n \iC$ and $\fcl {} {R_{12}} : [1] \to \lttr n \iC$ are composable morphisms in $\lttr n \iC$. Define a functor $W : [2] \to \lttr n \iC$ by setting $W(0 \to 1) = \fcl {} {R_{01}}$ and $W(1 \to 2) = \fcl {} {R_{12}}$. The required composition witness is (up to unique invertible concordance of the underlying $n$-truss bundle) the totalization $\fp {} {W}$ of $W$.
\end{obs}

\begin{obs}[Proof of \autoref{prop:labeled-n-truss-bordism}] \label{rmk:equivalence-of-viewpoints-on-labeled-n-truss-bordisms} The equivalences between the categories $\nlttr {n} \iC$, as defined in \autoref{defn:labeled-n-truss-bordism-category}, and $\lttr n \iC$, as defined in \autoref{defn:labeled-n-truss-bordism-endo}, can be constructed using totalization and classification as follows. A $\iC$-labeled $n$-truss $T$ (resp.\ $n$-truss bordism $R$) in $\nlttr {n} \iC$ can be classified by functors $\fcl {} T : [0] \to \lttr {n} \iC$ (resp.\ $\fcl {} T : [1] \to \lttr {n} \iC$) thus yielding objects (resp.\ morphisms) in $\lttr {n} \iC$. Conversely, objects (resp.\ morphisms) in $\lttr {n} \iC$ correspond to functors $T : [0] \to \lttr n \iC$ (resp.\ functors $R : [1] \to \lttr {n} \iC$) which we can totalize to $\iC$-labeled $n$-trusses $\fp {} T$ (resp.\ $n$-truss bordisms $\fp {} R$) thus yielding objects (resp.\ morphisms) in $\nlttr {n} \iC$.
\end{obs}

\subsubsecunnum{Pullback, dualization, and suspension of labeled n-truss bundles}

Our usual constructions of pullbacks, duals, and suspensions immediately carry over from 1-trusses bundles to $n$-trusses bundles.

\begin{constr}[Pullbacks of labeled $n$-truss bundles] \label{constr:labeled-n-truss-bundle-pullback} Consider a $\iC$-labeled $n$-truss bundle $p = (\und p, \lbl p)$ over a poset $B$ with underlying bundle $\und p = (p_n,p_{n-1},...,p_1)$, as well as a poset map $F : C \to B$. We define the \textbf{pullback $F^*p \equiv (F^*\und p, \lbl {F^*p})$ of $p$ along $F$} to be the $\iC$-labeled $n$-truss bundle with underlying bundle $F^*\und p = (F^*\und p_n, F^*\und p_{n-1}, ..., F^*\und p_1)$ constructed and labeling $\lbl {F^*p}$ as follows.
\begin{itemize}
    \item[-] Define $\Tot 0 F := F$.
    \item[-] For ascending $i$, define $\Tot i F : F^*p_i \to p_i$ to be the 1-truss bundle pullback map of $p_i$ along the poset map $\Tot {i-1} F : F^*T_{i-1} \to T_{i-1}$ (where $T_{i-1}$ resp.\ $F^*T_{i-1}$ is the total poset of $p_{i-1}$ resp.\ of $F^*p_{i-1}$).
\end{itemize}
Finally, define the labeling $\lbl {F^*p}$ as the composite $\lbl {p} \circ \Tot n F$.
\end{constr}

\nid In other words, the pullback $F^*p$ of a labeled $n$-truss bundle $p$ along a base poset map $F$ can be obtained by the upper row in the following diagram
    \begin{equation} \small
        \begin{tikzcd}
 & F^*T_n \arrow[ddr, phantom, "\lrcorner" , very near start, color=black] \arrow[dl, "\lbl {F^*p}"' pos=.2] \arrow[r, "F^*p_n"] \arrow[dd, "\Tot n F"] &[+20pt] F^*T_{n-1} \arrow[ddr, phantom, "\lrcorner" , very near start, color=black] \arrow[r, "F^*p_{n-1}"] \arrow[dd, "\Tot {n-1} F"'] &[+20pt] \cdots \arrow[r, "F^*p_2"] \arrow[dd, "\cdots", phantom] &[+20pt] F^*T_1 \arrow[ddr, phantom, "\lrcorner" , very near start, color=black] \arrow[r, "F^*p_1"] \arrow[dd, "\Tot 1 F"'] &[+20pt] C \arrow[dd, "\Tot 0 F = F"'] \\[-10pt]
            \iC & & & & & \\[-10pt]
& T_n \arrow[ul, "\lbl p" pos=.4] \arrow[r, "p_n"'] & T_{n-1} \arrow[r, "p_{n-1}"'] & \cdots \arrow[r, "p_2"'] & T_1 \arrow[r, "p_1"'] & B
\end{tikzcd}
    \end{equation}
Observe that the poset maps $\Tot i F$ (together with the labeling category functor $\id : \iC \to \iC$) assemble into a  $\iC$-labeled $n$-truss bundle map $F^*p \to p$ which we call the `pullback bundle map'. As before, in the special case where $F : C \into B$ is a subposet of $B$, this recovers our earlier definition of restrictions, i.e. $F^*p = \rest F C$.

Next we define duals of $\iC$-labeled $n$-truss bundles.

\begin{defn}[Duals of labeled $n$-truss bundles] \label{constr:labeled-n-truss-bundle-dualization} Given a $\iC$-labeled $n$-truss bundle $p \equiv (\und p,\lbl p)$ with underlying $n$-truss bundle $\und p = (p_n,p_{n-1},...,p_1)$, its \textbf{dual bundle} $p^\dagger$ is the $\iC\op$-labeled $n$-truss bundle whose underlying $n$-truss bundle is $\und {p^\dagger} = (p_n^\dagger,p_{n-1}^\dagger,...,p_1^\dagger)$ (where $p_i^\dagger$ is the dual of $p_i$, see \autoref{constr:1-truss-bundle-dualization}), and whose labeling $\lbl {p^\dagger}$ is the opposite labeling $(\lbl p)\op$.
\end{defn}

\begin{defn}[Duals of labeled $n$-truss bundle maps] Given a labeled $n$-truss bundle map $F : p \to q$ one defines the \textbf{dual bundle map} $F^\dagger : p^\dagger \to q^\dagger$ to have underlying $n$-truss bundle map is $\und{F^\dagger} = \und F\op$ (i.e. $\und F^\dagger$) and labeling functors $\lbl {F^\dagger} = (\lbl F)\op$.
\end{defn}

\begin{defn}[Duals of labeled $n$-truss bundle concordances] Given a labeled $n$-truss bundle concordance $u : p \Rightarrow q$ one constructs the \textbf{dual bundle concordance} $u^\dagger : q^\dagger \Rightarrow p^\dagger$ simply as the dual bundle of $u$.
\end{defn}

\nid The preceding definitions now yield the following functors, which are (covariant or contravariant) involutive isomorphisms of categories.

\begin{obs}[Dualization functors on labeled $n$-truss bundles] \label{thm:labeled-n-truss-bundle-dualization} The preceding definitions construct a dualization functor of labeled $n$-truss bundles:
    \[
        \dagger : \trussbunlbl n \iso \trussbunlbl n.
    \]
    Similarly, for a fixed category $\iC$ and a poset $B$, we obtain a dualization functor of labeled $n$-truss bundle concordances:
    \[
        \dagger : \trussconclbl n(B,\iC) \iso \trussconclbl n(B,\iC\op)\op .
    \]
If we set $B = \ast$, then this further specializes to an isomorphism of labeled $n$-truss bordisms:
    \[
        \dagger : \lttr n \iC \iso (\lttr n {\iC\op})\op . \qedhere
    \]
\end{obs}

\begin{rmk}[Duality of closed and open trusses, and singular and regular maps] Dualization maps closed $n$-trusses to open $n$-truss bundles, and singular $n$-truss bundle maps to regular $n$-truss bundle maps; in particular, dualization of $n$-truss bundles specializes to an isomorphism
    \begin{equation}
        \dagger : \sctruss n \toot \rotruss n : \dagger
    \end{equation}
between closed $n$-trusses with singular maps and open $n$-truss with regular maps.
\end{rmk}

Finally, we construct suspension of (unlabeled) $n$-truss bundles.

\begin{constr}[Suspension of unlabeled $n$-truss bundles] \label{constr:n-truss-bundle-suspension} For an (unlabeled) $n$-truss bundle $p = (p_n,p_{n-1},...,p_1)$ its \textbf{suspension $\Sigma p$} is the $n$-truss bundle with 1-truss bundle maps $(\Sigma p_n, \Sigma p_{n-1}, ..., \Sigma p_1)$ (where $\Sigma p_i$ is the suspension bundle of the 1-truss bundle $p_i$, see \autoref{constr:1-truss-bundle-suspension}).\end{constr}

\nid In analogy to \autoref{rmk:labeled-suspension-1-trs}, the construction may further be generalized to the labeled case as long as the labeling category $\iC$ has initial and terminal objects---we leave details to the reader.

\subsubsecunnum{Face, embedding, coarsening, \& degeneracy truss maps and their rigidity} Trusses are `rigid' when considered with special classes of $n$-truss maps, that is, no two maps in that class admit a non-trivial natural transformation between them. As we will later see, the types of maps in question combinatorially mirror the properties of certain stratified topological maps, and we introduce them here with accordingly mirrored terminology.

\begin{term}[Subtrusses, faces, and embeddings of 1-trusses] \label{term:1-truss-inj} A map of 1-trusses $F : T \to S$ that is injective on objects is called an `injection'. An injection $F$ may further have the following properties.
    \begin{enumerate}
        \item \emph{Balanced}: If $F$ is balanced then $F$ is a `subtruss' inclusion map.
        \item \emph{Closed singular}: If $T$ and $S$ are closed and $F$ is singular then we call $F$ a `closed face' (or simply, a `face').
        \item \emph{Open regular}: If $T$ and $S$ are open and $F$ is regular then we call $F$ an `open embedding' (or simply, a `embedding'). \qedhere
    \end{enumerate}
\end{term}

\begin{obs}[Characterizing faces and embeddings] \label{rmk:automatic-cellularity-n} Note that closed faces are exactly subtrusses of closed 1-trusses, and open embeddings are exactly subtrusses of open 1-trusses.
\end{obs}

\nid Conversely, we introduce the following terminology for surjective 1-truss maps.

\begin{term}[Degeneracies and coarsenings of 1-trusses] \label{term:1-truss-surj} A map of 1-trusses $F : T \to S$ that is surjective on objects is called a `surjection'. A surjection $F$ may further have the following properties.
    \begin{enumerate}
        \item \emph{Singular}: If $T$ and $S$ have the same endpoint type\footnote{Recall that this means $\dim(\ept_\pm T) = \dim(\ept_\pm S)$, see \autoref{rmk:trusses-by-boundary-type}. Note also that any surjection must be endpoint preserving, and thus the condition implies that $F$ preserves dimensions of endpoints.} and $F$ is singular then $F$ is called a `degeneracy'.
        \item \emph{Regular}: If $T$ and $S$ have the same endpoint type and $F$ is regular, then $F$ is called a `coarsening'.
        \item \emph{Closed singular}: If $T$ and $S$ are closed and $F$ is singular then we call $F$ a `closed degeneracy'.
        \item \emph{Open regular}: If $T$ and $S$ are open and $F$ is regular then we call $F$ an `open coarsenings'.
    \end{enumerate}
Note that if a surjective 1-truss map $F$ is balanced then it must be an isomorphism.
\end{term}

\begin{obs}[Characterizing closed degeneracies and open coarsenings] Note that closed degeneracies are exactly degeneracies of closed 1-trusses, and open coarsenings are exactly coarsenings of open 1-trusses.
\end{obs}

\begin{eg}[Faces, embeddings, degeneracies, and coarsenings of 1-trusses] In \autoref{fig:faces-embeddings-degeneracies-and-coarsenings} we depict examples of `subtrusses', ` faces', `embeddings', as well as `(closed) degeneracies' and `(open) coarsenings' of 1-trusses.
\begin{figure}[ht]
    \centering
    \def\svgwidth{1\columnwidth}
    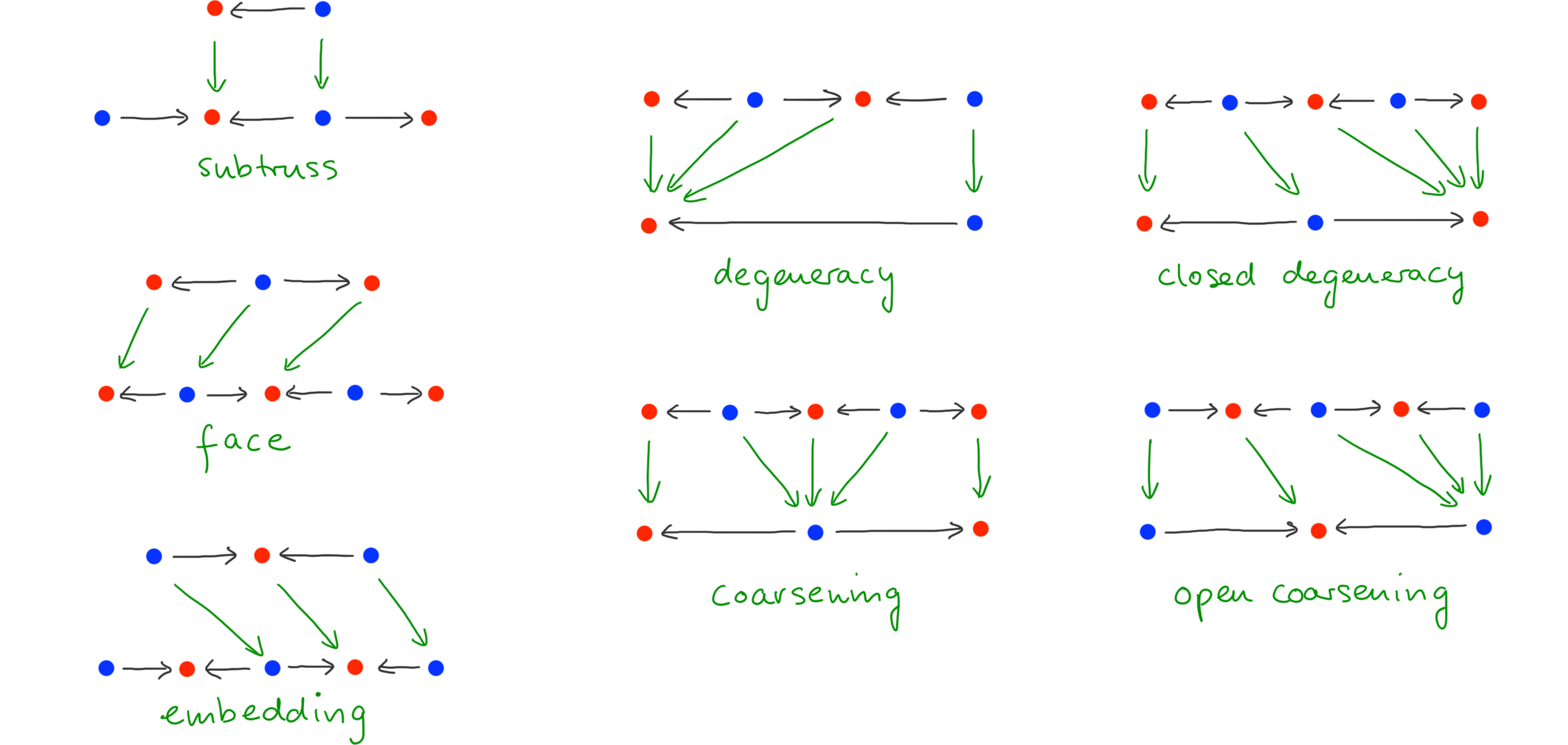

    \caption{Faces, embeddings, degeneracies, and coarsenings.}
    \label{fig:faces-embeddings-degeneracies-and-coarsenings}
\end{figure}
\end{eg}

We next generalize the preceding terminology to the case of (labeled) $n$-trusses as well as $n$-truss bundles.

\begin{term}[Faces, embeddings, degeneracies, and coarsenings of $n$-trusses]  \label{term:embeddings-subtrusses}  \label{term:coarsening-and-quotients} Given $n$-trusses $T = (p_n,p_{n-1},...,p_1)$ and $S=(q_n,q_{n-1},...,q_1)$, a truss map $F : T \to S$ is said to be a `injection' if each 1-truss bundle map $F_i : p_i \to q_i$ is fiberwise a 1-truss injection. Similarly, one says that $F$ is a `face', or an `embedding', or a `subtruss', or a `surjection', or a `(closed) degeneracy', or a `(open) coarsening' if each $F_i$ is fiberwise so.

    The terminology also applies to $n$-truss bundle maps $F : T \to S$ (we require $F$ to be base preserving), and further, to maps $F : T \to S$ of labeled $n$-trusses or labeled $n$-truss bundles (we require $F$ to be base preserving and label preserving).
\end{term}

\begin{notn}[Categories of truss coarsenings and truss degeneracies] \label{notn:truss-crs-and-deg} Denote by $\crstruss n$ the category of $n$-trusses and their coarsenings, and by $\degtruss n$ the category of $n$-trusses and their degeneracies.
\end{notn}

\begin{rmk}[Coarsenings vs refinements] \label{rmk:truss-crs-vs-ref} Given a coarsening of $n$-trusses  $F : T \to S$ we also call $F$ a `refinement' of $S$ by $T$. That is, we use the term `coarsening' and `refinement' synonymously but describing dual processes: a coarsening `coarsens' the domain, while a refinement, in opposite direction, `refines' the codomain.
\end{rmk}

The role of dual categories $\sctruss n$ of `closed trusses with singular maps' and $\rotruss n$ of `open trusses with regular maps' is special; these categories will be shown to combinatorially model framed regular cells (resp.\ its dual). In many ways, their properties reflect that of other `categories of shapes', such as the category of simplices $\Delta$. Our terminology here is meant to highlight this parallel: for instance, just as any morphisms in $\Delta$ factors into a degeneracy and face map (yielding its `(epi,mono)-factorization'), any singular $n$-truss map of closed trusses factors into degeneracy and a face map. Importantly, this existence of (epi,mono)-factorizations is special to the categories $\sctruss n$ and $\rotruss n$, and fails for general truss bundle maps as illustrated in the next example.

\begin{eg}[Failure of (epi,mono)-factorization in general] \autoref{fig:failure-of-(epi,mono)-factorization-of-a-2-truss-map} shows a map $F : T \to S$ of $2$-trusses: components $F_i : T_i \to S_i$ are highlighted in bold, and the mapping of the top components $F_2$ is determined by coloring its images corresponding preimages in the same color (note, as usual, singular elements are shown as red dots, while regular elements are shown as blue dots). The map $F$ cannot admit an (epi,mono) factorization since the image $\im(F_2) \subset S_2$ does not describe a subtruss of $S$.
\begin{figure}[ht]
    \centering
    \def\svgwidth{1\columnwidth}
    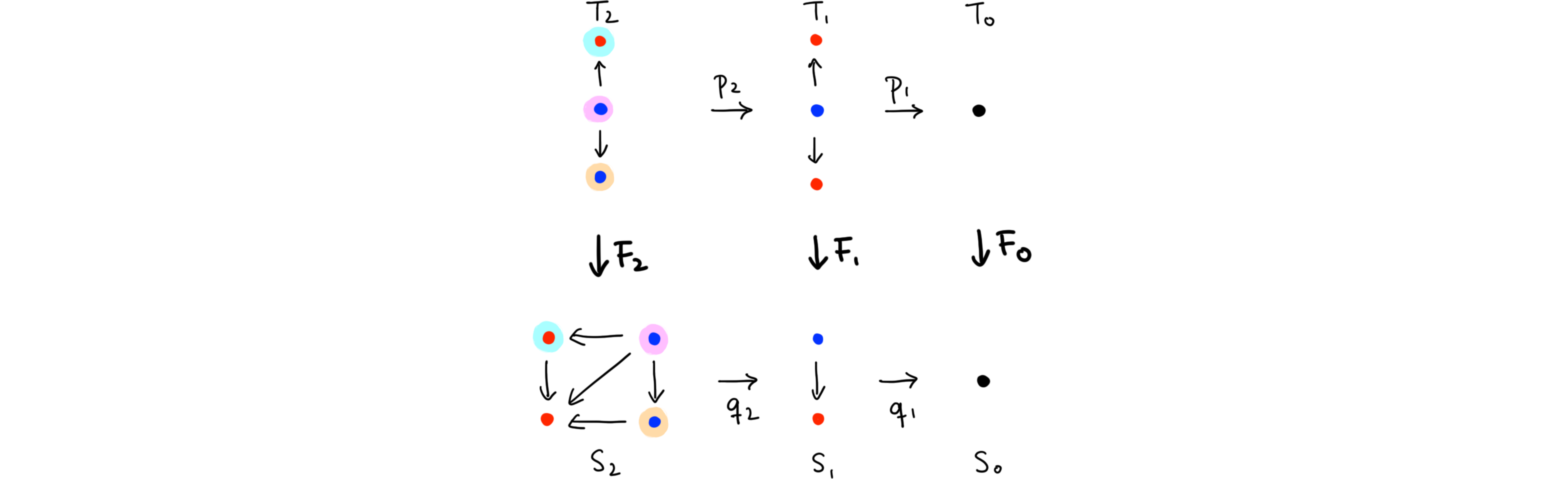

    \caption{Failure of (epi,mono)-factorization of a 2-truss map.}
    \label{fig:failure-of-(epi,mono)-factorization-of-a-2-truss-map}
\end{figure}
\end{eg}

\nid In contrast, singular maps of closed $n$-trusses, and, dually, regular maps of open $n$-trusses, do have (epi,mono)-factorizations. We record the existence of (epi,mono)-factorization for the categories $\sctruss n$ and $\rotruss n$ as follows.

\begin{lem}[(Epi,mono)-factorization] \label{lem:epi-mono-factorization}
    Any singular truss map $F$ of closed $n$-trusses factors uniquely into a degeneracy $F^{\epilabel}$ followed by a face $F^{\monolabel}$. Thus, after dualizing face posets, any regular truss map $F$ of open $n$-trusses factors uniquely into a coarsening $F^{\epilabel}$ followed by an embedding $F^{\monolabel}$.
\end{lem}
\begin{proof} In both cases the factorization $F = F^{\monolabel} F^{\epilabel}$ is determined by factoring $i$th poset maps $F_i = F^{\monolabel}_i F^{\epilabel}_i$ using the (epi,mono)-factorization in $\Pos$. We omit the verification that this determines $n$-truss bundle maps $F^{\epilabel}$ and $F^{\monolabel}$.
\end{proof}

\begin{eg}[(Epi,mono)-factorization of closed singular truss maps]  In \autoref{fig:(epi,mono)-factorization-of-closed-singular-2-truss-map} we depict a singular map $F : T \to S$ of closed $2$-trusses, together with its (epi,mono)-factorization $F = F^{\monolabel} \circ F^{\epilabel}$. As before, we indicate the mappings of the top components by coloring their images corresponding preimages in the same color.
\begin{figure}[ht]
    \centering
    \def\svgwidth{1\columnwidth}
    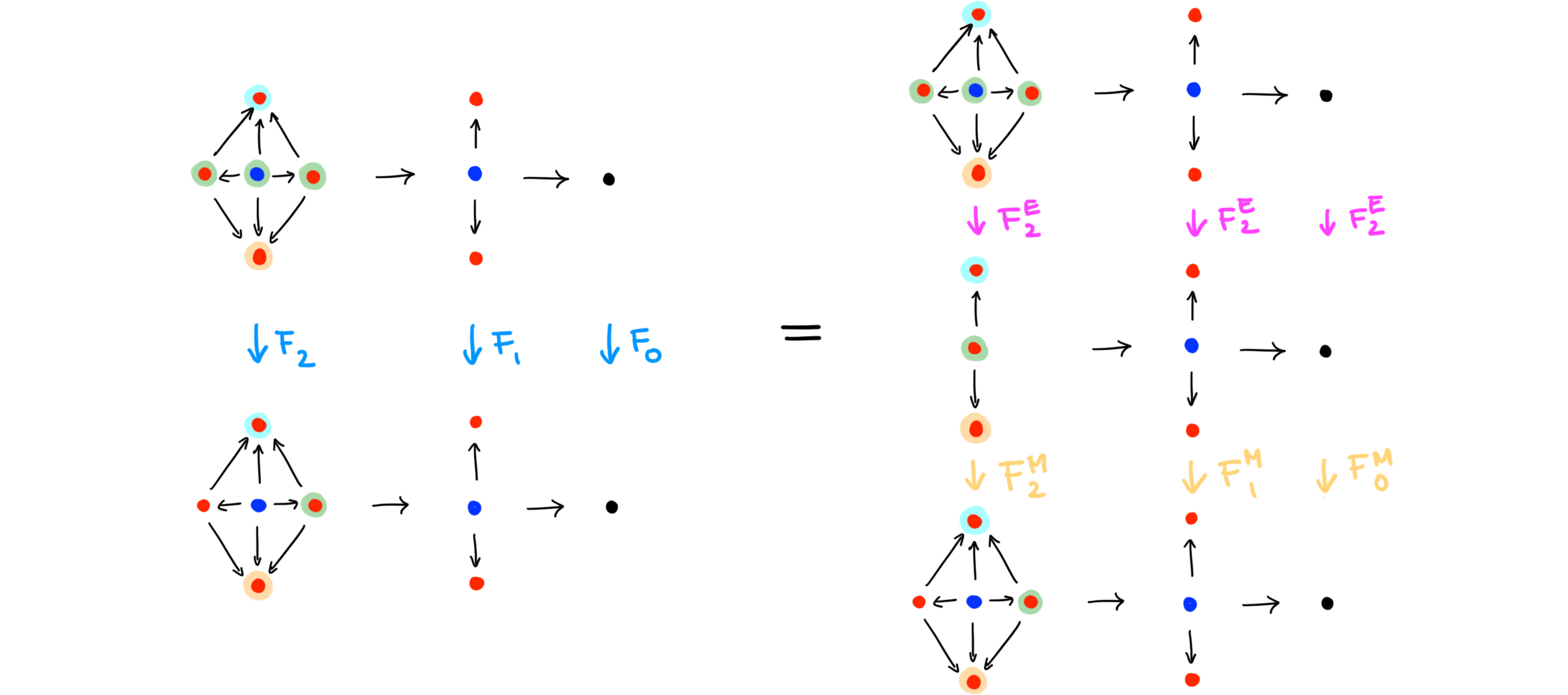

    \caption{(Epi,mono)-factorization of a closed singular 2-truss map.}
    \label{fig:(epi,mono)-factorization-of-closed-singular-2-truss-map}
\end{figure}
\end{eg}

Centrally, we now observe the following `rigidity of homs' in the categories $\sctruss n$ and $\rotruss n$, meaning that `hom posets' (consisting of maps and their natural transformations) are discrete. The result further applies to general `coarsening' and `degeneracies' as we record. We start with the case of 1-truss maps.

\begin{lem}[Rigidity of natural transformations for 1-trusses] \label{lem:nat-trafo-open-1truss} Consider 1-trusses $T$ and $S$, and 1-truss maps $E,F : T \to S$. Assume one of the following holds.
    \begin{enumerate}
        \item $T$ and $S$ are open and $E$ and $F$ are either embeddings, coarsenings, or (more generally) regular truss maps,
        \item $T$ and $S$ are closed and $E$ and $F$ are either faces, degeneracies, or (more generally) singular truss maps,
        \item $E$ and $F$ are coarsenings,
        \item $E$ and $F$ are degeneracies.
    \end{enumerate}
Then any natural transformation $\nu: E \Rightarrow F : (T,\eleq) \to (S,\eleq)$ must be the identity.
\end{lem}
\begin{proof} We argue in the first case of open 1-trusses and regular maps (the second case follows by duality, and the other cases follow similar arguments). Both embeddings $E$ and $F$ map the regular values $a$ in $T$ to regular values $E(a)$ resp.\ $F(a)$ in $S$; because there are no non-identity arrows between regular values in $S$, we must have $E(a) = F(a)$.  Because $T$ is assumed to be open, every singular value $b$ in $T$ has two adjacent regular values $b \pm 1$; since the maps $E$ and $F$ are functorial and monotone, both $E(b)$ and $F(b)$ must be the unique element $E(b) = F(b)$ in $S$ such that $E(b-1) \eleq E(b) \egeq E(b+1)$.  The functors $E$ and $F$ are thus identical and the natural transformation $\nu$ is necessarily trivial.
\end{proof}

\begin{lem}[Rigidity of natural transformations for $n$-trusses] \label{lem:nat-trafo-open-ntruss} Consider $n$-trusses $T$ and $S$, and $n$-truss maps $E, F : T \to S$. Assume one of the following holds
    \begin{enumerate}
        \item $T$ and $S$ are open and $E$ and $F$ are either embeddings, coarsenings, or (more generally) regular truss maps,
        \item $T$ and $S$ are closed and $E$ and $F$ are either faces, degeneracies, or (more generally) singular truss maps,
        \item $E$ and $F$ are coarsenings,
        \item $E$ and $F$ are degeneracies.
    \end{enumerate}
Then any natural transformation $\nu: E_n \Rightarrow F_n$ of poset maps $E_n, F_n : (T_n,\eleq) \to (S_n,\eleq)$ must be the identity.
\end{lem}

\begin{proof} We argue in the case of open $n$-trusses $T, S$ and regular maps $E,F$ (the second case follows by duality, and the other cases follow similar arguments). Write $T = (p_n,p_{n-1},...,p_1)$ and $S = (q_n,q_{n-1},...,q_1)$. Arguing inductively, we assume the statement holds for $(n-1)$-trusses (the bases case of $1$-trusses was shown in the previous lemma). Postcomposing $\nu$ with the bundle map $q_n$ yields a natural transformation $q_n \circ \nu : q_n \circ E_n \Rightarrow q_n \circ F_n$, which equivalently is a natural transformation $E_{n-1} \circ p_n \Rightarrow F_{n-1} \circ p_n$.  We must have $q_n \circ \nu = \nu_{n-1} \circ p_n$ for some natural transformation $\nu_{n-1} : E_{n-1} \Rightarrow F_{n-1}$.\footnote{In general, given poset maps $E,F : B \to C$ and $G : A \to B$, any natural transformation $\nu : E \circ G \Rightarrow F \circ G$ will be of the form $\mu \circ G$ for some $\Rightarrow$.}  Inductively, we deduce that $\nu_{n-1} = \id$. Applying the rigidity of 1-trusses (see \autoref{lem:nat-trafo-open-1truss}) to the transformation $\nu$ restricted to the fibers of $p_n$ and $q_n$ shows that $\nu$ is itself trivial as claimed.
\end{proof}

\nid The lemma further generalizes to the case of $E$ and $F$ being base and label preserving maps of labeled $n$-truss bundles.

\subsection{Truss blocks and truss block sets} \label{ssec:block-and-brace-trusses}

A truss block (or simply a `block') is a closed truss with an initial element; the name reflects the idea that all closed trusses can be `built from blocks', as we will make precise in this section. Just as closed trusses dualize to open trusses, blocks dualize to open trusses with terminal elements, which we will refer to as `braces'.

\subsubsecunnum{The definition of truss blocks} Recall that the depth of an object in a poset is the maximal length of chains starting in that object (for instance, a maximal object is of depth 0).

\begin{defn}[Truss blocks] \label{defn:trussblock} An \textbf{$n$-truss $k$-block} $T$ is a closed $n$-truss $T$ whose total poset $(T_n,\eleq)$ has an initial object, and that object is of depth $k$.
\end{defn}

\nid When referring to $n$-truss $k$-blocks, we often keep either or both of the dimensions $n$ and $k$ implicit, or, yet more simply, just refer to them as `blocks'.

\begin{rmk}[Block truncate] Note that, given an $n$-truss block $T = (p_n,p_{n-1},...,p_1)$, then each truncation $T_{\leq i} = (p_i,p_{i-1},...,p_1)$ is again a $i$-truss block; indeed, the minimal element in the total poset of $T_{\leq i}$ is simply the image of the initial object in the total poset of $T$ under the map $p_{>i} = p_{i+1} \circ p_{i+2} \circ ... \circ p_n$.
\end{rmk}

\begin{notn}[Initial elements of blocks] We often denote the initial element in the total posets of a block by $\ino$.
\end{notn}

\begin{eg}[A 2-truss 2-block] In \autoref{fig:a-2-block-with-geo-real} we illustrate a 2-truss 2-block on the left, together with its `geometric realization' on its right. In contrast to our earlier example of a general closed 2-truss in \autoref{fig:a-first-2-truss-example}, note that the realization now consists of a single 2-cell.
\begin{figure}[ht]
    \centering
    \def\svgwidth{1\columnwidth}
    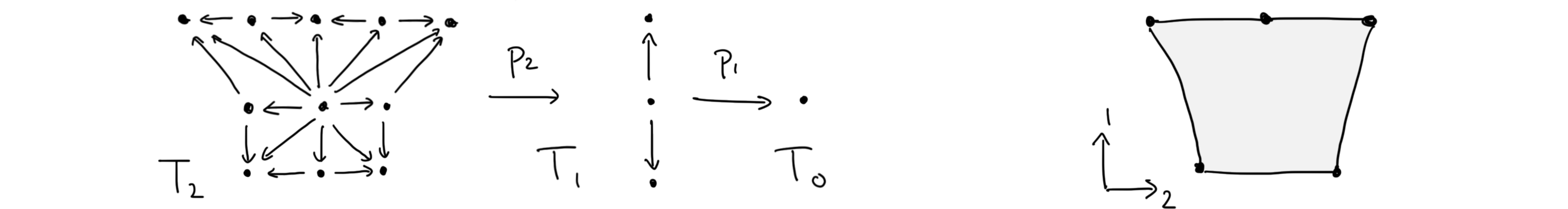

    \caption[A 2-truss 2-block.]{A 2-truss 2-block together with its geometric realization as a single framed 2-cell.}
    \label{fig:a-2-block-with-geo-real}
\end{figure}
\end{eg}

\begin{notn}[Category of blocks] The category $\blcat n$ of $n$-truss blocks and singular maps is defined to be the full subcategory of the category $\sctruss n$ of closed trusses and singular maps whose objects are $n$-truss blocks.
\end{notn}

Recall the notion of faces of closed trusses from \autoref{term:embeddings-subtrusses}. We construct faces in trusses that are blocks by taking closures of their elements, as follows.

\begin{constr}[Face blocks in closed trusses] \label{constr:blocks-in-trusses} Let $T = (p_n,p_{n-1},...,p_1)$ be a closed $n$-truss, and consider an element $x \in T_n$ in the total poset $T_n$ of $T$. We construct a subtruss inclusion $T^{\egeq x} \into T$ such that $T^{\egeq} = (p^{\eleq x}_n,p^{\eleq x}_{n-1},...,p^{\eleq x}_1)$ is a block, called the `face block of $x$'. For each $i\leq n$, denote by $x_i = p_{>i} x$ the image of $x$ under $p_{>i} = p_{i+1} \circ ... \circ p_n : T_n \to T_i$. Define $(T^{\eleq x}_i,\eleq)$ to be the subposet of $T_i$ given by the upper closure of $x$ in $(T_i,\eleq)$. Define $p^{\eleq x}_i : T^{\eleq x}_i \to T^{\eleq x}_{i-1}$ to be the restriction of $p_i : T_i \to T_{i-1}$ to these subposets. One verifies that we may uniquely endow $p_i$ with 1-truss bundle structure such that the subposet inclusions $T^{\eleq x}_i \into T_i$ induce 1-truss bundle maps $p^{\eleq x}_i \into p_i$. The subposet inclusions $T^{\eleq x}_i \into T_i$ also define the components of the claimed subtruss inclusion $T^{\egeq x} \into T$. Note further that $T^{\leq x}$ is a $k$-block where $k$ is the depth of $x$ in $T_n$.
\end{constr}

\begin{eg}[Blocks in closed trusses] In \autoref{fig:two-block-closures-of-elements-in-trusses} we depict a $2$-truss $T_2 \to T_1 \to T_0 = [0]$, and highlight two elements $x, y \in T_2$. For both elements we then construct their face block as shown. Identifying the total poset $T_2$ with the entrance path poset of the cell complex on the right, note that face blocks (obtained by taking `upper closures' in $T_2$) correspond closed cells in the cell complex (obtained by taking `topological closure).
\begin{figure}[ht]
    \centering
    \def\svgwidth{1\columnwidth}
    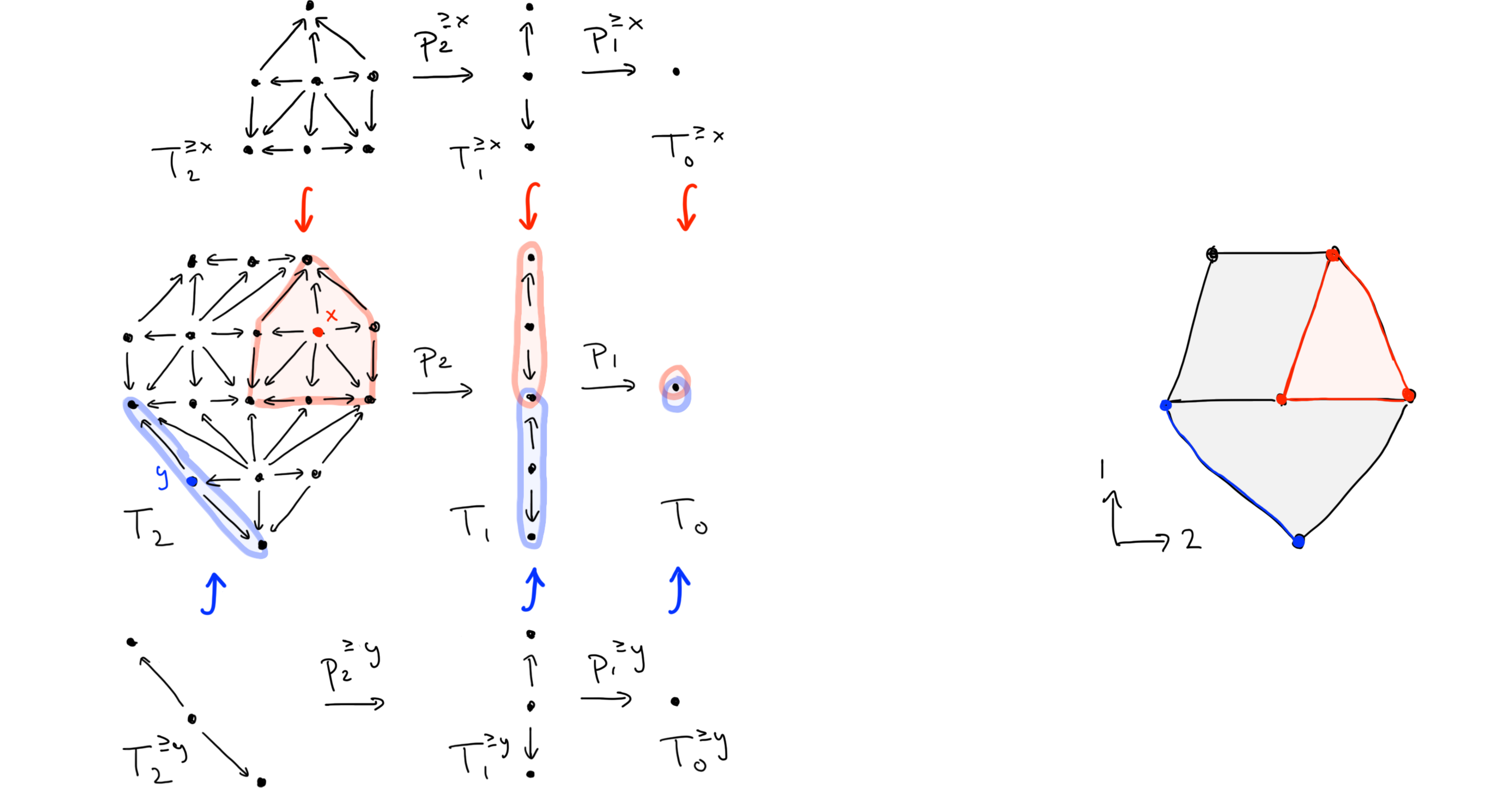

    \caption[The face blocks of elements in trusses.]{The face block of two elements $x$ (red) and $y$ (blue) in a $2$-truss $T$.}
    \label{fig:two-block-closures-of-elements-in-trusses}
\end{figure}
\end{eg}

\nid In fact, all subtrusses that are blocks are obtained by taking closures as in the previous construction.

\begin{rmk}[All subtruss blocks are face blocks] \label{rmk:block-faces-by-objects} Given an $n$-truss $T$, then subtrusses of $T$ which are $n$-truss blocks are in 1-to-1 correspondence with objects of $T_n$: indeed, every such subtruss block of $T$ determines an element of $T_n$ as the image of its initial element, and conversely objects in $T_n$ determined face blocks by taking closures as explained in \autoref{constr:blocks-in-trusses}.
\end{rmk}

\subsubsecunnum{The definition of  block sets}

In analogy with simplicial sets we now introduce the following notion of `$n$-truss block sets'.

\begin{defn}[Block sets] A \textbf{$n$-truss block set}, is a $\SetCat$-valued presheaf on the category $\blcat n$ of $n$-truss blocks.
\end{defn}

\nid We usually abbreviate `$n$-truss block set' simply to `block set', leaving the dimension $n$ implicit.

\begin{term}[Category of block sets] The `category of block sets and their maps' is the category of $\SetCat$-valued presheaves on the category $\blcat n$ of $n$-truss blocks.
\end{term}

\begin{term}[Faces, degeneracies, and non-degenerate blocks] Given a block set $X \in \blset n$ and $B \in \blcat n$, we call elements of $X(B)$ `blocks' of $X$ of `shape $B$'. Given a face map $F : C \to B$ in $\blcat n$ such that $(X(F))(b) = c$ for $b \in X(B)$, then we call $c \in X(C)$ a `face' (or `$F$-face') of $b$. Given a degeneracy map $F : C \to B$ in $\blcat n$ such that $(X(F))(b) = c$ for $b \in X(B)$, then we call $c \in X(C)$ a `degeneracy (or `$F$-degeneracy') of $b$. If $c \in X(C)$ is not an $F$-degeneracy for any non-identity $F$, then we call $c$ a `non-degenerate' block of $X$.
\end{term}

\begin{notn}[Representable block sets] Given $B \in \blcat n$ we often keep the Yoneda embedding notationally implicit, and write the representable presheaf $\blcat n(-,B)$ simply as $B$.
\end{notn}

\begin{constr}[Block nerve of trusses] Given a closed $n$-truss $T$ be can construct its `block nerve' to be the block set $\ctruss n(-,T)$ mapping a block $B$ to the hom set $\ctruss n(B,T)$. This construction is functorial and gives to the `block nerve functor' $\bnrv : \sctruss n \to \blset n$.
\end{constr}

\begin{rmk}[Building trusses from their blocks] We can now make precise our previous statement that each closed $n$-truss $T$ can be `built from their blocks'. Given a closed $n$-truss $T$, denote by $\blcat n \slash T$ the comma category of the inclusion $\blcat n \into \sctruss n$ over $T$ (that is, objects of $\blcat n \slash T$ are singular maps $B \to T$ from blocks into $T$, and morphisms are commuting block maps $B \into B'$). One verifies that $T$ is given by the colimit
    \begin{equation}
        T = \colim (\blcat n \slash T \to \sctruss n)
    \end{equation}
of the forgetful functor $\blcat n \slash T \to \sctruss n$, which maps $(B \to T)$ to $B$.
\end{rmk}
\nid The reader familiar with the general notion of nerves may also have observed that there are equivalent ways to state this remark: namely, it is equivalent to the observation that the functor $\blcat n \into \sctruss n$ is dense; it is also equivalent to the observation that the functor $\bnrv$ is fully faithful.

We next introduce a notion of `regular' block sets.

\begin{defn}[Regular block set] A block set $X \in \blset n$ is called \textbf{regular} if for each non-degenerate block $b \in X(B)$ the block set map $b : B \to X$ is a monomorphism.
\end{defn}

\begin{prop}[Block nerves are regular block sets] Given a closed $n$-truss $T$, its block nerve $\bnrv T$ is regular.
\end{prop}

\begin{proof} Non-degenerate blocks $b : B \to \bnrv T$ are of the form $\bnrv(F : B \into T)$ where $F : B \into T$ is a face block of $T$ (see \autoref{rmk:block-faces-by-objects}). One verifies that $\bnrv(F : B \into T)$ is a monomorphism of block sets as required.
\end{proof}

\begin{notn}[Category of regular block sets] The full subcategory of the category of block set $\blset n$ containing only regular block sets will be denoted by $\rblset n$.
\end{notn}


\subsubsecunnum{Truss braces and truss brace sets}

Finally, let us briefly record the duals of truss blocks, which we will call `truss braces'. Recall that open and closed $n$-trusses are related by duality functors $\dagger : \sctruss n \iso \rotruss n : \dagger$ which dualize face orders on each truss. The preceding discussion of truss blocks and blocks sets may be completely dualized in this way. We mirror the most central definitions of the discussion for convenience. The `height' of an object $x$ in a poset measures the length $k$ of maximal chains $x_{-k} \to  x_{-k+1} \to ... \to x_0$ in the poset ending at that object, $x = x_0$.

\begin{defn}[Truss braces] An \textbf{$n$-truss $k$-brace} $T$ is an open $n$-truss whose total poset $(T_n,\eleq)$ has a terminal object, and that object is of height $k$.
\end{defn}

\nid More concisely, we often refer to `$n$-truss $k$-braces' simply as `braces'. Dually to `face blocks' we then find a notion of `embedding braces' as follows.

\begin{constr}[Embedding braces in open trusses] Let $T$ be an open $n$-truss. For any $x$ in its total poset $T_n$, we construct the subtruss $T^{\eleq x} \into T$, called the \textbf{embedding brace of $x$}: it is the unique open subtruss of $T$ such that $T^{\eleq x}$ is a brace whose terminal element maps to $x$.
\end{constr}

\begin{notn}[Category of braces] The category $\brcat n$ of braces and regular maps is defined to be the full subcategory of the category $\rotruss n$ of open trusses and regular maps whose objects are braces.
\end{notn}

\begin{defn}[Brace sets] The \textbf{category of brace sets} $\brset n$ is the category of presheaves on the category $\brcat n$ of braces and regular maps. Objects of $\brset n$ will be called \textbf{brace sets} and morphisms \textbf{brace set maps}.
\end{defn}

\nid Note that the dualization isomorphism $\dagger : \sctruss n \iso \rotruss n$, induces an isomorphism of presheaf categories $\brcat n \iso \blcat n$ (by precomposition with $\dagger$).

\chapter{Constructibility of framed combinatorial structures}  \label{ch:classification-of-framed-cells}

The central theorems of this chapter will construct equivalences between the framed combinatorial structures introduced in \autoref{ch:framed-combinatorial-structures} and the iterated constructible combinatorial structures discussed in \autoref{ch:trusses}. The most elementary such equivalence states that $n$-framed regular cells are classified by $n$-truss blocks as follows.

\begin{thm}[Truss blocks classify framed regular cells] \label{thm:classification-of-cells} Framed regular cells are classified by truss blocks; that is, there is a canonical equivalence of categories
    \begin{equation}
       \begin{tikzcd}
            \FrCCell n  \arrow[r, "\kT", shift left] & \blcat n \arrow[l, "\kX", shift left]
        \end{tikzcd}
    \end{equation}
\end{thm}

\begin{figure}[ht]
    \centering
    \def\svgwidth{1\columnwidth}
    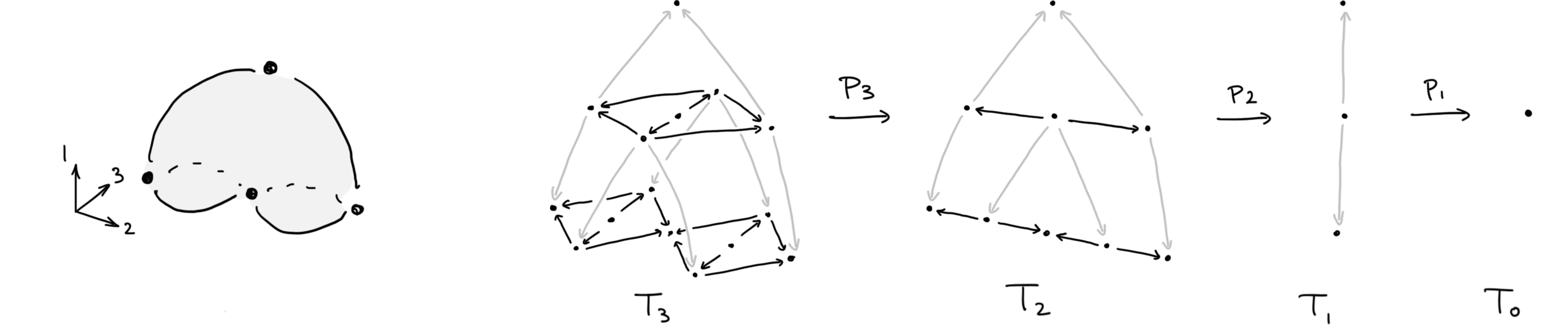

    \caption{Framed regular cells correspond to truss blocks.}
    \label{fig:framed-regular-cells-translate-to-blocks-and-vice-versa}
\end{figure}

\nid The equivalence is illustrated in \autoref{fig:framed-regular-cells-translate-to-blocks-and-vice-versa}: on the left, we depict a $3$-framed regular cell (note this cell also appeared in \autoref{fig:3-framed-regular-3-cells-examples}); and on the right, we depict a $3$-truss block $T$ (note that we only depict generating arrows in the posets $T_i$, see \autoref{rmk:generating-arrows}). The $3$-framed regular cell will be mapped to the $3$-truss by the functor $\cT$, and the inverse mapping will be described by the functor $\kX$, as given in the previous theorem.

The correspondence of framed regular cells and truss blocks can be generalized to a correspondence of \emph{flat} framed regular cell complexes and general trusses as follows.

\begin{thm}[Trusses classify flat framed regular cell complexes] \label{thm:classification-of-flat-framed-cell-cplx} Flat framed regular cell complexes are classified by closed trusses and their singular maps; that is, there is a canonical equivalence of categories
    \begin{equation}
       \begin{tikzcd}
            \FrCDiag n  \arrow[r, "\kT", shift left] & \sctruss n \arrow[l, "\kX", shift left]
        \end{tikzcd}
    \end{equation}
\end{thm}

\begin{figure}[ht]
    \centering
    \def\svgwidth{1\columnwidth}
    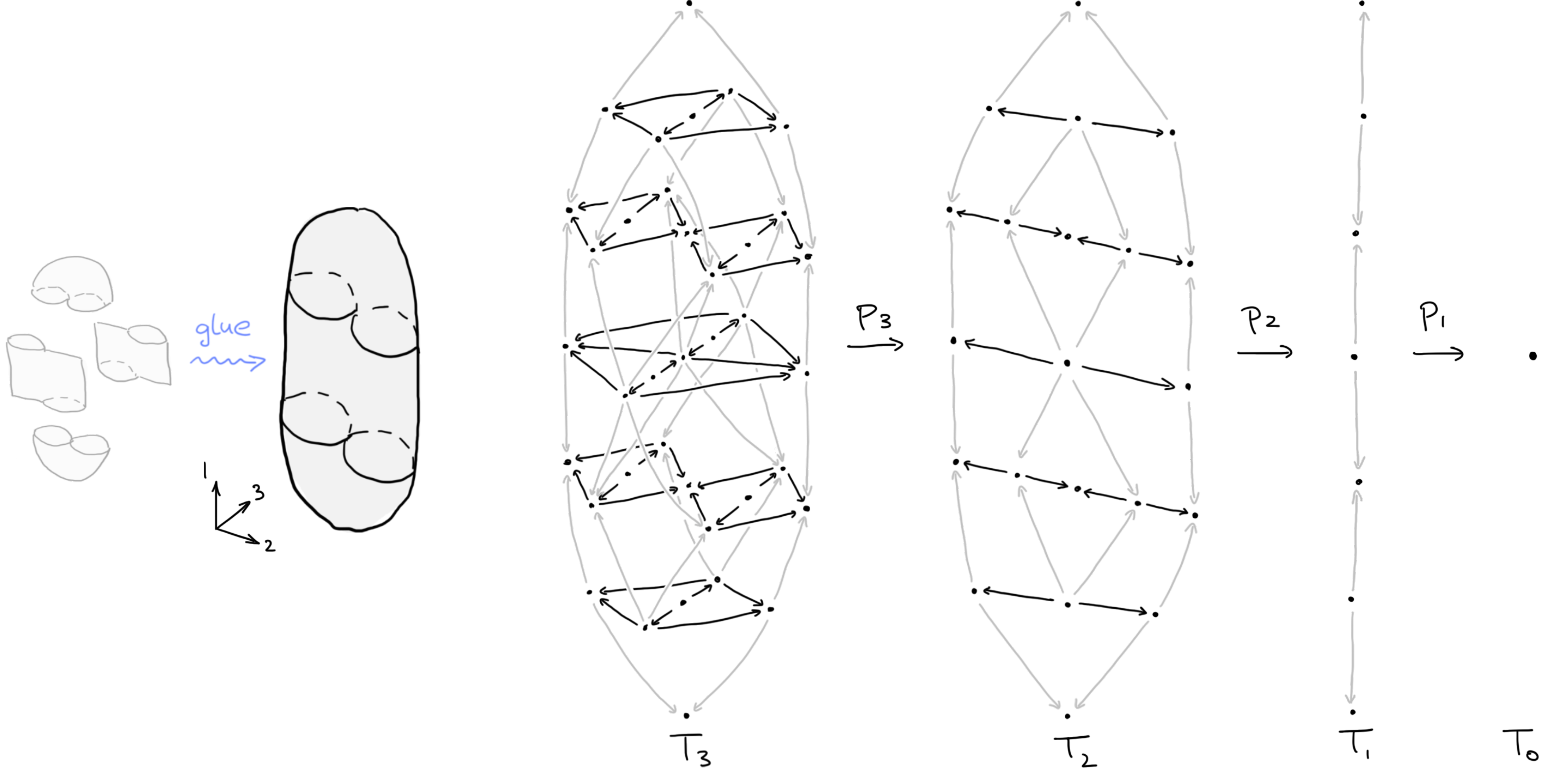

    \caption{Flat framed regular cell complexes correspond to closed trusses.}
    \label{fig:flat-framed-regular-cell-complexes-translate-to-closed-trusses-and-vice-versa}
\end{figure}

\nid In \autoref{fig:flat-framed-regular-cell-complexes-translate-to-closed-trusses-and-vice-versa} we illustrate an instance of this equivalence, depicting a flat $3$-framed regular cell complex on the left (note that this complex also appeared in \autoref{fig:the-dual-of-the-twisted-embedding-of-the-circle}). On the right, we depict its corresponding closed $3$-truss.l

Yet more generally, the preceding two equivalences are restrictions of the following equivalence.

\begin{figure}[ht]
    \centering
    \def\svgwidth{1\columnwidth}
    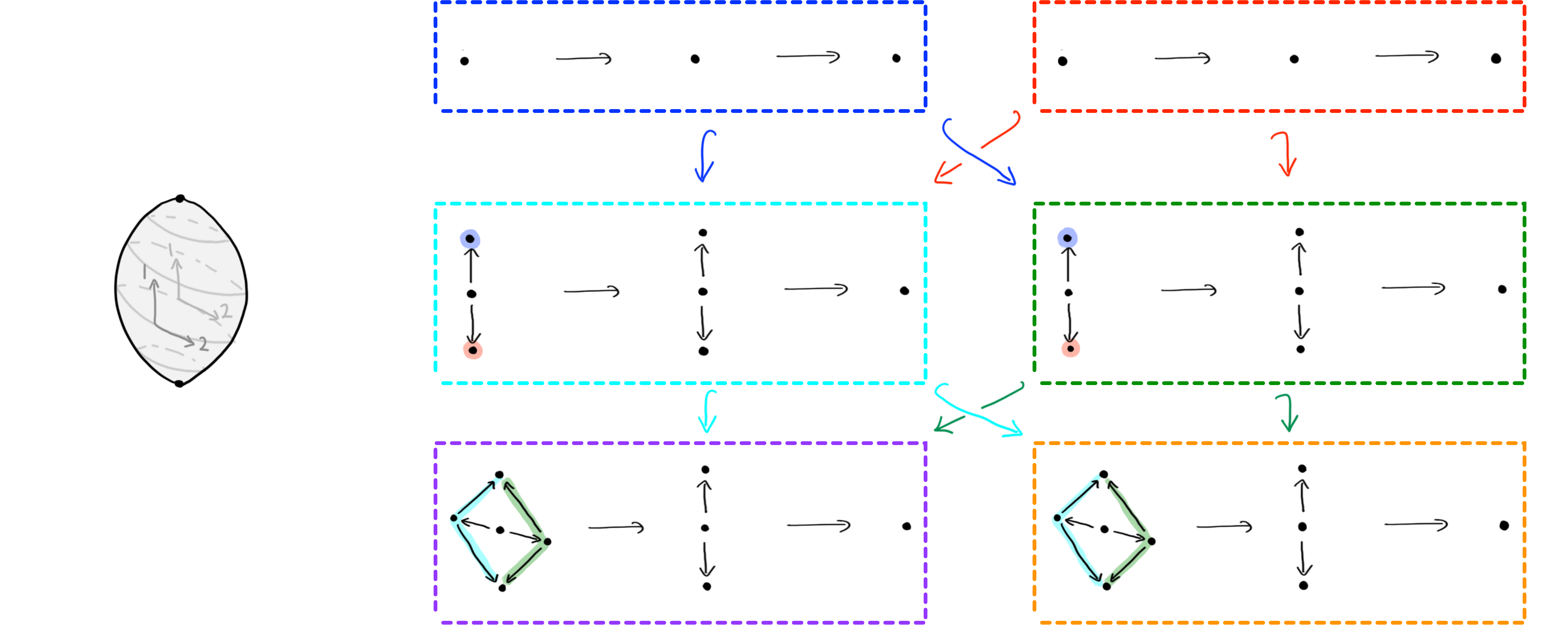

    \caption{Framed regular cell complexes correspond to regular truss block sets.}
    \label{fig:framed-regular-cell-complexes-translate-to-regular-block-sets-and-vice-versa}
\end{figure}

\begin{thm}[Truss block sets classify framed regular cell complexes] \label{thm:classification-of-framed-cell-cplx} Framed regular cell complexes are classified by regular truss block sets; that is, there is a canonical equivalence of categories
    \begin{equation}
       \begin{tikzcd}
            \FrCCplx n  \arrow[r, "\kT", shift left] & \rblset n \arrow[l, "\kX", shift left]
        \end{tikzcd}
    \end{equation}
\end{thm}
\nid In \autoref{fig:framed-regular-cell-complexes-translate-to-regular-block-sets-and-vice-versa} we illustrate an instance of this equivalence. On the left, we depict a $2$-framed regular cell complex, consisting of two 0-cells, connected by two 1-cells, between which we suspend two 2-cells (note that this complex also appeared in \autoref{fig:block-complex-unit-counit-gluing}). The regular truss block set on the right is indicated by its six non-degenerate truss blocks, each marked by its own box, which correspond to the cells on the left, and with face relations as indicated by colored inclusions.

\pauseae

In each of the above three equivalences, we call the functor $\kT$ the `truss translation', and the functor $\kX$ the `framed complex translation'. The latter functor is not hard to define. Namely, for a closed $n$-truss $T = (T_n \xto {p_n} T_{n-1} \xto {p_{n-1}} ... \xto {p_1} T_0)$ we obtain a framed complex translation $\kX T$ as follows. We take the underlying cellular poset to be $T_n$. For each (unordered) simplex $x : \Unord{[k]} \into T_n$, we can restrict the tower $T$ to $x$ to obtain a tower of (equally unordered) simplicial projections $\rest T x$. The frame orders of each 1-truss bundle in $T$ now determine a unique ordering of the tower $\rest T x$ (namely, such that kernel vectors of projections are morphisms in the frame order of fibers).  Ordering the tower $\rest T x$ defines a proframe on $\Unord [k]$, and passing to its gradient frames then defines a frame. Framing each (unordered) simplex in $T_n$ in this way defines a framing, and yields the framed regular cell complex $\kX T$.

Of course, to verify this construction one still needs to check, for instance, that the total posets $T_n$ is cellular and that the resulting framed complex is flat. Verifying these details (as well as describing the inverse functor $\kT$, which will be based on a series of technical observations about flat framed cell complexes) will take up the rest of this chapter. Readers should therefore feel to skip ahead to \autoref{ch:meshes}. For those interested in the technical details, we now outline the proofs in this chapter as follows.
\begin{enumerate}
    \item In \autoref{sec:equiv-of-frames-and-proj-frames} we establish the equivalence of framed and proframed regular cells, yielding an equivalence $\FrCCell n \iso  \PFrCCell n$. There is a similar equivalence between flat framed and flat proframed regular cell complexes. These equivalences allow us to work with proframed cells in place of framed cells. Since proframed cells (like trusses) are defined by towers of maps, this equivalence will facilitate the construction of the truss corresponding to a framed cell.
    \item In \autoref{sec:equiv-of-frccplx-and-trusses} we then construct the truss translation and framed complex translation functors, in particular verifying the correctness of the above sketch definition of the framed complex translation (the crucial step in this verification will be checking the cellularity of face posets of trusses, for which we employ a result relating `shellability' and `cellularity' of posets).
\end{enumerate}

\section{Equivalences of framings and proframings} \label{sec:equiv-of-frames-and-proj-frames}

Our goal in this section will be the construction of equivalences between framings and proframings in three instances: (1) flat framings and flat proframings of simplicial complexes, (2) framings and proframings of regular cells, and (3) flat framings and flat proframings of regular cell complexes.

\subsection{Equivalence of flat framed and flat proframed simplicial complexes} \label{ssec:frames-vs-proj-frames-flat}

We will prove that the gradient framing functor $$\Gradfr : \FlPFrSCplx n \to \FlFrSCplx n$$ (see \autoref{notn:gradient-framing-functor}) provides an equivalence of flat proframed simplicial complexes with flat framed simplicial complexes.

\subsubsecunnum{Gradients of flat proframings are flat}

We start with a proof of \autoref{prop:grad-of-flat-proframing-is-flat-framing}---recall, this stated that the gradient framing of a flat proframing is flat itself.

\begin{proof}[Proof of \autoref{prop:grad-of-flat-proframing-is-flat-framing}] Consider a flat $n$-proframing $\cP = (p_n,...,p_1)$ on a simplicial complex $K$. Since flatness of proframings truncates, we can argue inductively and assume that the gradient framing $\Gradfr \cP_{<n}$ of the truncation $\cP_{<n}$ of $\cP$ is flat: that is, it admits a framed bounded realization $e_{n-1} : \abs{K_{n-1}} \into \lR^{n-1}$. We now construct a geometric realization $\abs{p_n} : \abs{K_n} \to \abs{K_{n-1}}$ of the map $p_n$, and a framed bounded realization $e_n : \abs{K_n} \into \lR^n$ such that the following commutes
    \begin{equation}
        \begin{tikzcd}
        \abs{K_n} \arrow[r, hook, "e_n"] \arrow[d, "\abs{p_n}"'] & \lR^n \arrow[d, "\pi_n"] \\
    \abs{K_{n-1}} \arrow[r, hook, "e_{n-1}"] & \lR^{n-1}
        \end{tikzcd}
    \end{equation}
    We first define $e_n: \abs{K_n} \into \lR^n$ on vertices $x \in \abs{K_n}$ as follows. Assume $x$ lies over a base vertex $y = p_n x \in \abs{K_{n-1}}$. Define $e_n(x) = (e_{n-1}(y), e_y(x))$ where $e_x : \abs{K_y} \into \lR$ is some choice of framed realization of the linear fiber complex $K_y$ into $\lR$ (see \autoref{rmk:linear-fiber-complexes-in-flat-proframings}). Now extend the embedding $e_n$ linearly to all other simplices in $\abs{K_n}$. We then define $\abs{p_n}$ to be the composite $e_{n-1}\inv \circ \pi_n \circ e_n$. Using the definition of flat proframings, one check that $e_n$ defines a framed realization whose image is framed bounded (as required in the definition of flat framings), and that $\abs{p_n}$ realizes $p_n$ as claimed.
\end{proof}

\subsubsecunnum{Existence of integral proframings for simplicial complexes}

We next prove \autoref{prop:flat-framings-have-integrating-flat-proframing}---recall this stated that any flat framing has an essentially unique integral flat proframing. We will prove this statement in two steps, first observing existence and then essential uniqueness.

\begin{lem}[Existence of integral proframings of flat framings] \label{prop:existence-of-integrating-proj-framings} Every flat $n$-framed simplicial complex $(K,\cF)$ has an integral flat $n$-proframing $\cP$.
\end{lem}

\begin{proof} Our strategy is to construct $\cP$ inductively in $n$. For this, we construct a flat $(n-1)$-framing of $K_{n-1}$, denoted by $cF_{n-1}$, and uniquely determined by the following condition: $(K_{n-1},\cF_{n-1})$ admits a framed bounded realization $e_{n-1} : \abs{K_{n-1}} \into \lR^{n-1}$ (see \autoref{term:framed-bounded-real}) such that the composite $e_{n-1}\inv \circ \pi_n \circ e_n$ geometrically realizes a simplicial map $p_n : K_n \to K_{n-1}$ (where $K_n$ is the ordering of $K$ determined by its $n$-framing). The construction of $(K_{n-1},\cF_{n-1})$ relies on the framed boundedness of the embedding $e_n$ and we outline its steps as follows. By definition of flatness, $e_n \abs{K_n}$ is a framed bounded subspace of $\lR^n$ of the form $\bigcap_{i\leq n} (\lR^n_{\geq \gamma^-_i} \cap \lR^n_{\leq \gamma^+_i})$. Therefore, $\pi_n : \lR^n \to \lR^{n-1}$ projects $e_n \abs{K}$ to the bounded subspace $\bigcap_{i\leq n-1} (\lR^{n-1}_{\geq \gamma^-_i} \cap \lR^{n-1}_{\leq \gamma^+_i})$ of $\lR^{n-1}$ (which, in particular, is again framed bounded). The definition of framed realizations implies $\pi_n e_n$ maps simplices of $\abs{K}$ linearly to simplices in $\pi_n e_n \abs{K}$. Framed boundedness of $e_n$ implies that the collection of such image simplices in fact triangulates $\pi_n e_n \abs{K}$ (see \autoref{fig:the-projected-framing-of-a-flat-framing}). Moreover, one checks this triangulation is the embedding $e_{n-1} : \abs{K_{n-1}} \into \lR^{n-1}$ of a unique ordered simplicial complex $K_{n-1}$ subject to the condition that $e_{n-1}\inv \pi_n e_n : \abs{K} \to \abs{K_{n-1}}$ realizes an ordered simplicial map $p_n : K_n \to K_{n-1}$. Finally, endow $K_{n-1}$ with the unique framing $\cF_{n-1}$ making $e_{n-1}$ a framed realization. This completes the construction $(K_{n-1},\cF_{n-1})$.

    Now, arguing inductively, define $\cP_{<n}$ to be a flat $(n-1)$-proframing obtained as the integral proframing of $\cF_{n-1}$. Then define the $n$-proframing $\cP$ of $K$ by extending the tower $\cP_{<n}$ by the map $p_n$. Once more using the definition of framed realizations and boundedness, one verifies that the proframed simplicial complexes $(K,\cP)$ is indeed flat proframed.
\end{proof}

\begin{term}[Projected framings] \label{term:projected-framings} The $(n-1)$-framing $\cF_{n-1}$ of $(K_{n-1}$ constructed in the preceding proof will be referred to as the `projected framing' of the flat $n$-framing $\cF$ of $K$, and the map $p_n : K_n \to K_{n-1}$ as its `projection map'.
\end{term}

\begin{eg}[Projected framings] \label{eg:projected-framings}
In \autoref{fig:the-projected-framing-of-a-flat-framing}, we illustrate the projected framing $\cF_2$ for a flat $3$-framing $\cF$ of a simplicial complex $K$.
\begin{figure}[ht]
    \centering
    \def\svgwidth{1\columnwidth}
    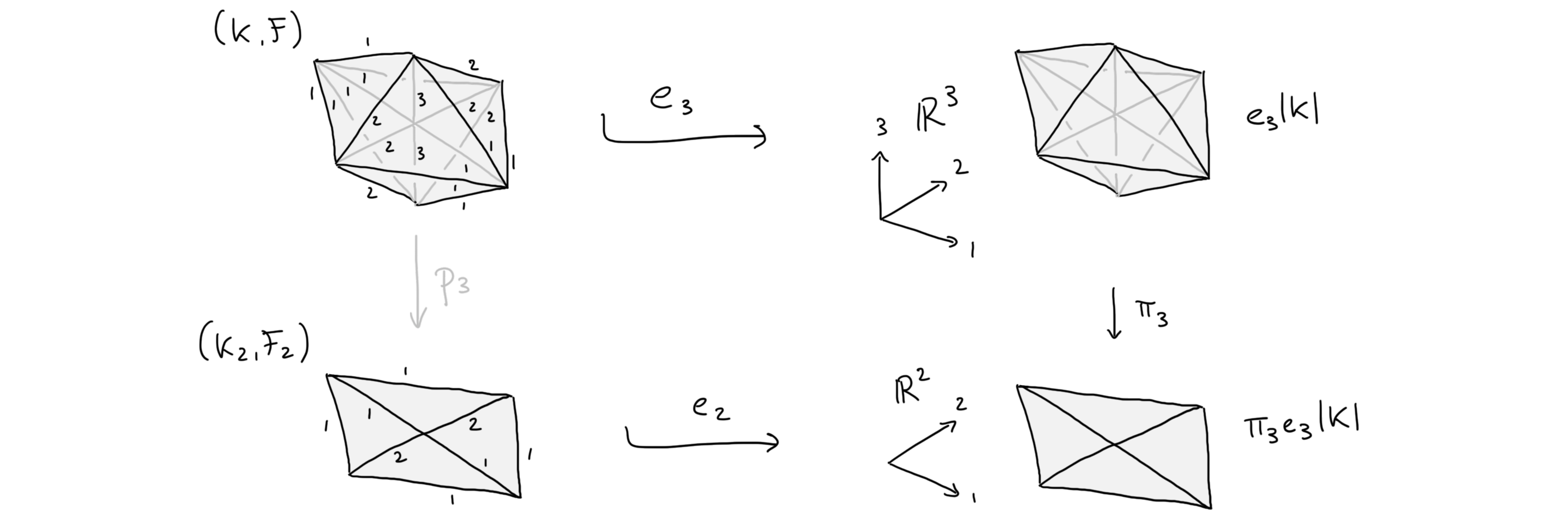

    \caption{The projected framing of a flat framing.}
    \label{fig:the-projected-framing-of-a-flat-framing}
\end{figure}
\end{eg}

\begin{lem}[Essential uniqueness of integral proframings of flat framings] \label{prop:uniqueness-of-integrating-proj-framings} Every flat $n$-framed simplicial complex $(K,\cF)$ has, up to unique isomorphism, a unique integral flat $n$-proframing.
\end{lem}

The following terminology will be useful.

\begin{term}[Simplex directions in proframings] \label{defn:kernels} Given an $n$-proframed simplicial complex $(K,\cF)$, we say 1-simplex $x : [1] \to K$ is `$i$-directed' if $p_{>i}(x)$ is a 1-simplex that is degenerated by $p_i$ (where $p_{>i}$ abbreviates $p_{i+1} \circ p_{i+2} \circ ... \circ p_n$ as before).
\end{term}

\begin{proof}[Proof of \autoref{prop:uniqueness-of-integrating-proj-framings}] Let $\cP = (p_n,...,p_1)$ and $\cQ = (q_n,...,q_1)$ be flat $n$-proframings of simplicial complexes $K$ resp.\ $L$ whose gradient $n$-framings $(K,\Gradfr \cP)$, $(L,\Gradfr \cQ)$ are isomorphic. We will prove that $\cP$ and $\cQ$ are isomorphic proframings related by a unique isomorphisms $\cP \iso \cQ$ whose $n$th component $K_n \iso L_n$ must equal the framed isomorphism $(K,\Gradfr \cP) \iso (L,\Gradfr \cQ)$. Our strategy is to argue by induction, this time, `from the bottom up'. We will first consider fibers of the map $p_{>1} : K_n \to K_1$. Note, since $\cP$ is flat, $K_1$ is a finite linear simplicial complex of the form $\lincplx j$ for some $j \in \lN$ (see \autoref{rmk:lin-complex}). Define `fiber' subcomplexes $K(i) := p_{>1}\inv(i)$ of $K$ for each vertex $i$ in $\lincplx j \iso K_1$. Note we can restrict the tower of projections $\cP$ to $K(i) \into K$, and drop the trivial last projection of the resulting restricted tower, to obtain an $(n-1)$-proframing $\cP(i)$ on $K(i)$. Since $\cP$ is flat it follows by definition of flatness that $\cP(i)$ is flat as well. Note that no 1-simplex in $K(i) \into K$ can be a 1-directed 1-simplex in $K$ (see \autoref{defn:kernels}). Conversely, any simplex from a vertex in $K(i)$ to a vertex $K(i')$, $i \neq i'$, must be 1-directed in $p$ and we must have $i' = i + 1$. There must also be at least one 1-simplex between $K(i)$ and $K(i+1)$ since $p_{>1}$ is assumed surjective on simplices (note, applying this observation inductively, one sees that $K_n(i)$ are in fact connected subcomplexes).

    We exhibited $K(i)$ as connected components of `$K$ without its 1-directed 1-simplices'; and further, 1-directed 1-simplices of $K$ fix an order $K(0)$, $K(1)$, $...$, $K(j)$ of these components, as they only run between $K(i)$ and $K(i+1)$. The assumed isomorphism $(K,\Gradfr \cP) \iso (L, \Gradfr \cQ)$ must therefore restrict to isomorphisms $(K(i),\Gradfr \cP(i)) \iso (L(i),\Gradfr \cQ(i))$ (as $(n-1)$-framed simplicial complexes). Arguing inductively, deduce unique isomorphism $\cP(i) \iso \cQ(i)$ for all $i$ such that $n$th components of these isomorphisms equal the isomorphisms $(K(i),\Gradfr \cP(i)) \iso (L(i),\Gradfr \cQ(i))$. One checks that this extends to the claimed unique isomorphism of $n$-proframings $(K,\cP) \iso (L,\cQ)$ as claimed. 
\end{proof}

\begin{proof}[Proof of \autoref{prop:flat-framings-have-integrating-flat-proframing}]
The proposition follows from \autoref{prop:existence-of-integrating-proj-framings} and \autoref{prop:uniqueness-of-integrating-proj-framings}.
\end{proof}

\nid As a corollary to the preceding observations, we can now record the following more explicit construction of the essentially unique integral proframings of flat framings.

\begin{constr}[The integration of flat framings via iterated projection] \label{obs:unwinding-integration-of-flat-fr} Combining the proofs of \autoref{prop:existence-of-integrating-proj-framings} and \autoref{prop:uniqueness-of-integrating-proj-framings}, it follows that the essentially unique integral proframing of a flat framed simplicial complex $(K,\cF)$ can be obtained by inductively constructing, for $i = n,n-1,...$, the projected framings $(K_i, \cF_i)$ and their projections $p_i : K_i \to K_{i-1}$ (see \autoref{term:projected-framings}) yielding a tower of simplicial maps $(p_n,...,p_1)$. Moreover, inspecting the proof of \autoref{prop:existence-of-integrating-proj-framings}, we see each $p_i$ acts by quotienting vertices of $K_i$ along 1-simplices $x$ with framing $\cF_i(x) = i$. We illustrate this in \autoref{fig:unwinding-the-essentially-unique-integration-of-flat-framings}.
\begin{figure}[ht]
    \centering
    \def\svgwidth{1\columnwidth}
    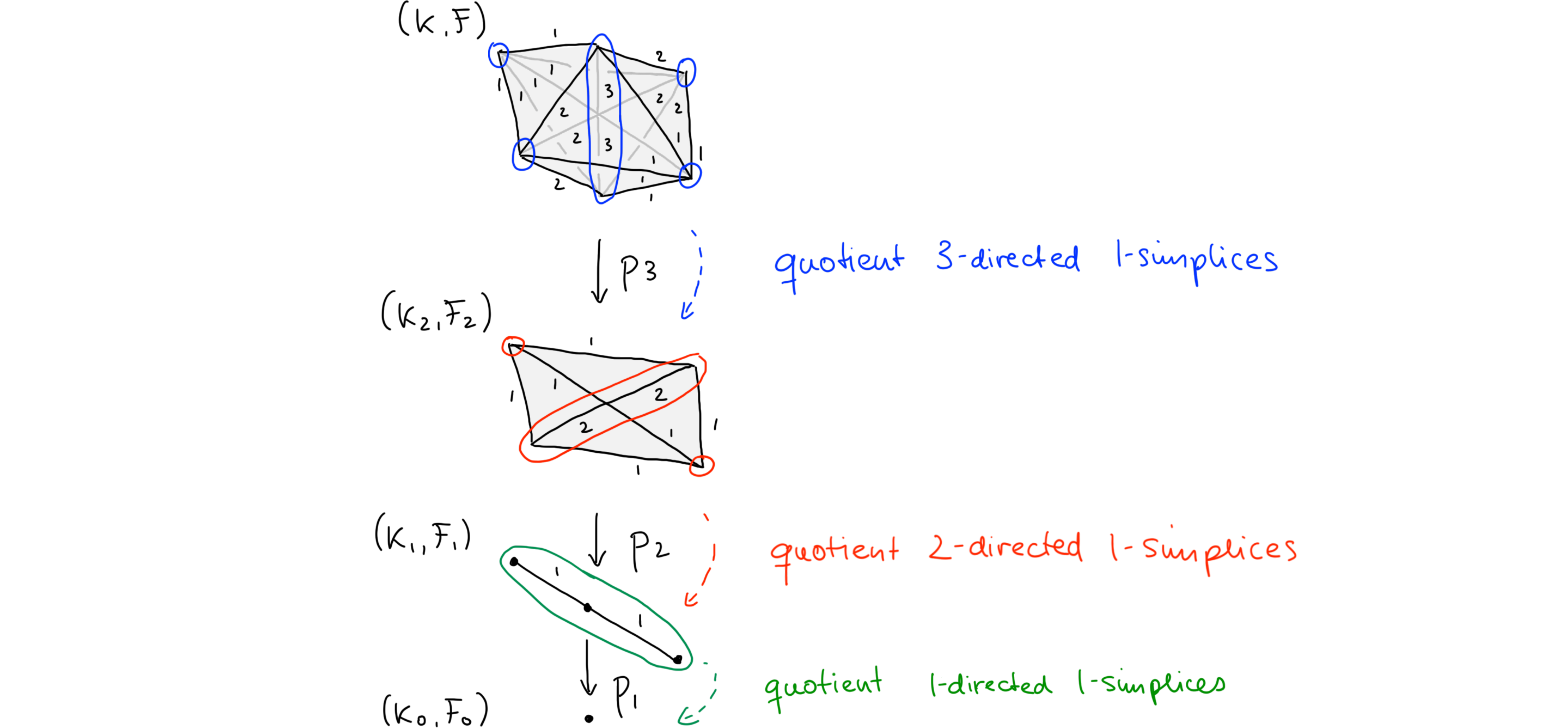

    \caption[The unique integration of flat framings.]{Constructing the essentially unique integration of flat framings as a tower of quotient maps.}
    \label{fig:unwinding-the-essentially-unique-integration-of-flat-framings}
\end{figure}
\end{constr}

\begin{rmk}[Equivalence of all framed realizations of flat framings] \label{rmk:framed-embeddings-of-flat-framings} The proof of \autoref{prop:existence-of-integrating-proj-framings} shows that any framed bounded realization $e = e_n : \abs{K} \into \lR^n$ of a flat $n$-framed simplicial complex $(K,\cF)$ descends through a tower of simplicial maps $p_i : K_i \to K_{i-1}$ to framed bounded realizations $e_{i} : \abs{K_i} \into \lR^i$ such that $\pi_i e_i = e_{i-1} \abs{p_i}$. The proof of \autoref{prop:uniqueness-of-integrating-proj-framings} then implies that all such framed bounded realizations are `equivalent' in that, given another framed bounded realization $e' = e'_n : \abs{K} \into \lR^n$, we have unique isomorphisms $\phi_i : e_i \abs{K_i} \iso e'_i \abs{K_i}$ (which are linear on each simplex) such that $\pi_i \phi_i = \phi_{i-1} \pi_i$. Namely, using notation from the proof of \autoref{prop:uniqueness-of-integrating-proj-framings}, the isomorphisms $\phi_i$ may be constructed by first inductively defining isomorphisms $\phi_i(j) : e_i \abs{K_i(j)} \to  e_i \abs{K_i(j)}$ on the `fibers of $p_{>1}$ over $j \in K_1$', and then extending these isomorphism of $p_{>1}$-fibers linearly to the remaining simplices of $(K_i,\cF_i)$ (those containing 1-directed edges). In fact, a similar argument applies to \emph{any} framed realization $e'$ of $(K,\cF)$ (without requiring boundedness of $e'$) which similarly shows equivalence of $e$ and $e'$ in the above sense. Note, in particular, that this entails that all framed realizations of flat framings are framed bounded.
\end{rmk}

\begin{constr}[The integration functor for flat framings] \label{constr:integration-constr-1} We construct the functor $\Intfr : \FrSCplx n \to \PFrSCplx n$, called the integration functor of flat framings. On objects, this functor maps flat $n$-framed simplicial complexes to their essentially unique integral flat $n$-proframed simplicial complexes, as constructed by the preceding two lemmas. On morphisms, construct the functor as follows. Let $F : (K,\cF) \to (L,\cG)$ be a framed map of flat $n$-framed simplicial complexes. Denote the integral proframings of $(K,\cF)$ resp.\ of $(L,\cG)$ by $(K,\Intfr \cF)$ resp.\ $(L,\Intfr \cG)$. Note that, since $F$ is a framed map, it must map $k$-directed 1-simplices in $K$ either to $0$-simplices or to $k$-directed 1-simplices in $L$. In particular, it preserves fibers of $p_n$ and $q_n$ which allows us to construct $F_{n-1} : K_{n-1} \to L_{n-1}$ such that $q_n F = F_{n-1} p_n$. Continuing inductively, this constructs a map of towers $\Intfr F : (K,\Intfr \cF) \to (L,\Intfr \cG)$. Comparing definitions, one checks that $\Intfr F$ is proframed since $F$ is framed.
\end{constr}

\nid Finally, we prove \autoref{obs:frames-vs-proj-frames-flat}, which stated that the gradient framing functor yields an equivalence of categories $\Gradfr : \PFrSCplx n \to \FrSCplx n $.

\begin{proof}[Proof of \autoref{obs:frames-vs-proj-frames-flat}] The inverse of the `gradient' functor $\Gradfr$ is the `integration' functor $\Intfr$ defined in \autoref{constr:integration-constr-1}.
\end{proof}


\subsection{Equivalence of flat framed and flat proframed cell complexes} \label{ssec:frames-vs-proj-frames-reg-cell}

We next show that the gradient framing functor, that takes proframed regular cells to framed regular cells, has an inverse. The inverse, as before, will be called `integration'. Recall, an $n$-framed regular $k$-cell $(X,\cF)$ consists of a combinatorial regular $k$-cell $X$ together with a flat $n$-framing $\cF$ of the simplicial complex $X$. Throughout this section, let $(X,\cF)$ denote an $n$-framed regular cell with $\ino$ being the initial element of $X$, and denote by $\cP = (X = X_n \xto {p_n} X_{n-1} \xto {p_{n-1}} ... \xto {p_1} X_0)$ the unique integral flat $n$-proframing of the \emph{simplicial complex} $X$ that integrates $\cF$ (as constructed in the previous section). Our goal will be to endow the simplicial complexes $X_i$ with cellular poset structures (i.e.\ to define cellular posets whose underlying simplicial complexes recover $X_i$), such that $\cP$ becomes an $n$-proframing of the \emph{regular cell} $X = X_n$. This will then provide the unique integral $n$-proframed regular cell for the $n$-framed regular cell $(X,\cF)$.

\subsubsecunnum{Upper and lower section cells} We now discuss a distinction of framed cells into `section' and `spacer' cells, such that spacer cells `fill the space' between their `upper' and `lower' section cells (the idea is analogous to the discussion of section and spacer simplices in the context of simplices in proframings, see \autoref{term:sections-spacers-in-proframings}).

\begin{prop}[Section and spacer cells] \label{obs:section-and-spacer-cells} For any $n$-framed regular cell $(X,\cF)$, one of the following holds true.
\begin{enumerate}
    \item Either, the fiber $X_{p_n(\ino)}$ over $p_n(\ino)$ is trivial (i.e.\ contains only $\ino \in X$) and the ordered simplicial map $p_n : X_n \to X_{n-1}$ is an isomorphism.
    \item Or, the fiber $X_{p_n(\ino)}$ is isomorphic to the 3-element linear complex $(\ino-1) \to \ino \to (\ino+1)$ with $\ino$ as its middle element. \qedhere
\end{enumerate}
\end{prop}

\begin{term}[Section and spacer cells] \label{term:section-and-spacer-cells} If a framed regular cell $(X,\cF)$ satisfies the first case of \autoref{obs:section-and-spacer-cells} then we call it a `section cell'. In the second case, we call it a `spacer cell', and call the elements $(\ino\pm1)$ the (upper resp.\ lower) `central fiber bounds' of the cell.
\end{term}

\begin{eg}[Section and spacer cells] We illustrate both cases of \autoref{obs:section-and-spacer-cells} in \autoref{fig:section-and-spacer-cells}, highlighting the poset structure $X$ of the framed cell $(X,\cF)$ in blue and central fiber elements in red.
\begin{figure}[ht]
    \centering
    \def\svgwidth{1\columnwidth}
    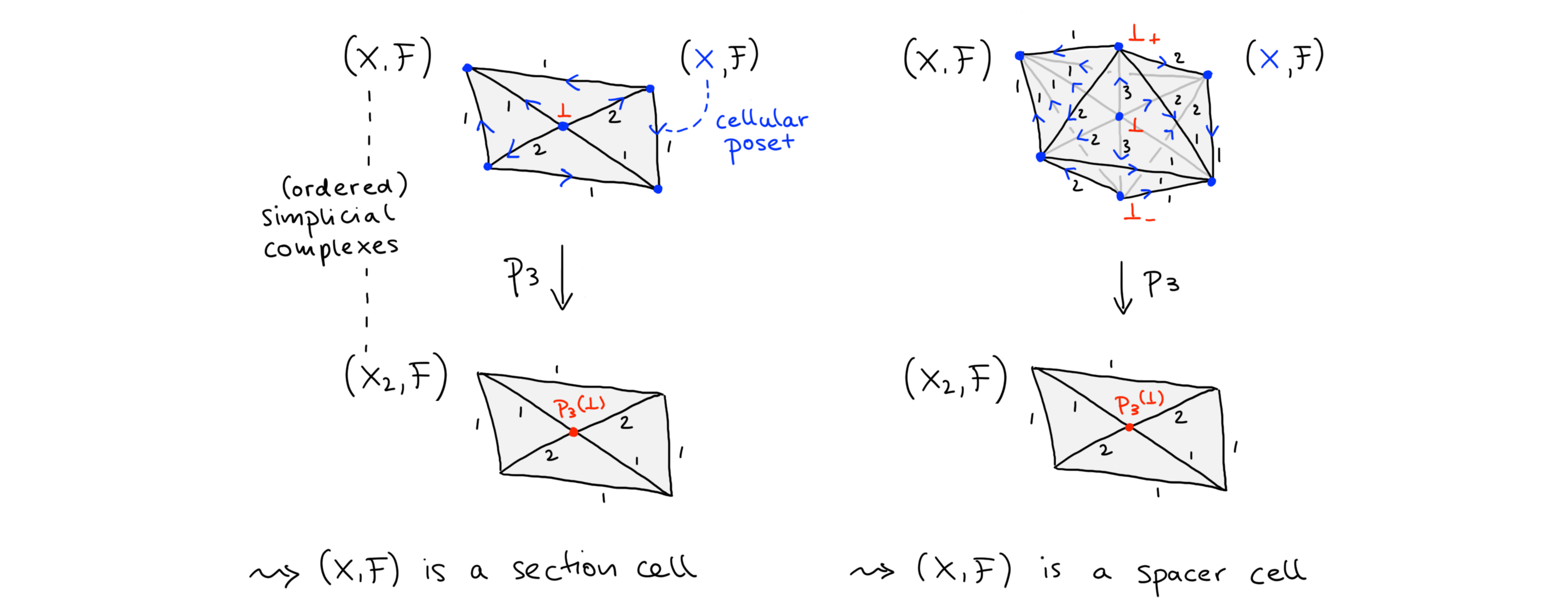

    \caption{Framed regular cells are either section or spacer cells.}
    \label{fig:section-and-spacer-cells}
\end{figure}
\end{eg}

\begin{proof}[Proof of \autoref{obs:section-and-spacer-cells}] We first note that the fiber $X_{p_n(\ino)}$ is a non-empty linear order with \emph{at most} three elements. Indeed, by definition of flat proframings, the ordered simplicial complex $X_{p_n(\ino)}$ is a linear complex of the form $\lincplx j$ for some $j \in \lN$, and, since $\ino$ is the initial object in the poset $X$, this implies that $X_{p_n(\ino)}$ can have at most two 1-simplices which must contain $\ino$ as a vertex (and thus at most three objects). Accordingly, we distinguish the following three cases.

    First, assume the fiber $X_{p_n(\ino)}$ has a single element $\ino$. We claim that $p_n : X_n \to X_{n-1}$ does not have spacer simplices (which, by surjectivity of $p_n$ and flatness of $\cP$, equivalently says that $p_n$ is an isomorphism). Arguing by contradiction, assume there is a spacer simplex $x$ of $X_n$, and pick such $x$ of maximal dimension. Since $\ino$ is initial in $X$, $x$ must contain $\ino$ in its vertices. Note $x$ cannot be the face of any other simplex (since this would have to be again a spacer simplex). Consequently, $\dim(x) = k$ equals dimension of the cell $X$. Denote by $z = p_n x$ the $(k-1)$-simplex in $X_{n-1}$ that $x$ projects to. Recall, the fiber category $\FC X_z$ over $z$ is the category consisting of section over $z$ with spacers in between section acting as generating morphisms (see \autoref{defn:fundamental-fiber-category}). Let $y$ be the initial section in the fiber category $\FC X_z$ (the terminal section works just as well). Note $\dim(y) = k-1$. Note $\ino$ is a vertex of $y$. We show that $y$ is the $(k-1)$-face of exactly one $k$-simplex in $X$. By our choice of $y$, it is the face of exactly one spacer (note each spacer has exactly two faces that are sections, its upper and lower section). Assume $y$ were the face of a section $y'$. Since $X$ is a $k$-cell, we must have $\dim(y') = k$ and $\FC X_{z'} = [0]$ where $z' = p_n y'$. But since $z$ has non-trivial fiber category, this contradicts the transition functor $\FC X_{z'} \to \FC X_{z}$ being endpoint-preserving. We deduce that $y$ is indeed the $(k-1)$-face of exactly one $k$-simplex in $X$. This, however, is impossible since $\ino$ is a vertex in $y$: indeed, it would imply that the point $\ino$ in $\abs{X}$ has no neighborhood homeomorphic to $\lR^k$, which contradicts $X$ being a regular $k$-cell. The contradiction proves that $p_n$ cannot have spacers, and thus $p_n$ is an isomorphism.

    Next, assume that the fiber $X_{p_n(\ino)}$ has two elements. We show this is impossible. Note $X_{p_n(\ino)}$ is of the form  $\ino \to (\ino+1)$ or $(\ino-1) \to \ino$. Argue in the former case (the latter case is similar). Since $p_n$ has at least one spacer simplex (namely in the fiber over $p_n(\ino)$), we can pick a spacer simplex $x$ of $X_n$ of maximal dimension, and, as before, we must have $\dim(x) = k$. By the same arguments as in the previous paragraph, construct an initial section $y$ over $p_n x$ and show that $\abs{\ino \in X}$ does not have euclidean neighborhood, thus contradicting that $X$ is a regular $k$-cell.

    This leaves only the case that $\FC X_{p_n(\ino)}$ has three elements, in which case it must be of the form $(\ino-1) \to \ino \to (\ino+1)$.
\end{proof}

\begin{term}[Upper and lower sections] \label{obs:upper-and-lower-sections-in-reg-cells} Define ordered simplicial maps $\gamma^- : X_{n-1} \to X_n$ resp.\ $\gamma^+ : X_{n-1} \to X_n$, called the `lower section' resp.\ `upper section' of $(X,\cF)$, mapping a $j$-simplex $z$ of $X_{n-1}$ to the $j$-simplex of $X_n$ that is initial resp.\ terminal in the fiber category $\FC X_z$ over $z$.
\end{term}

\nid Note that $\gamma^- = \gamma^+$ if and only if $(X,\cF)$ is a section cell. Indeed, if $(X,\cF)$ is spacer cell, then $(\ino\pm1) \in \im(\gamma^\pm)$. In fact, images of upper and lower sections are exactly the cells determined by $(\ino\pm1)$, as the following result shows.

\begin{lem}[Section images are cells] \label{lem:section-images-are-cells} If $(X,\cF)$ is a spacer cell, with lower and upper sections $\gamma^\pm$, then $X^{\geq (\ino\pm1)}$ are cells of dimension $(k-1)$, whose simplices are exactly those in the image of  $\gamma^\pm$.
\end{lem}
\begin{proof} We argue in the case of the lower section $\gamma^-$ (the case of the upper section is similar). First observe that any simplex in $X$ containing $(\ino-1)$ but not $\ino$ is a section simplex for $p_n$ (indeed, otherwise we could pick a spacer simplex containing $(\ino-1)$ but not $\ino$, and, by flatness of $\cP$ and since $(X,\cF)$ is a spacer $k$-cell, we could choose such spacer of dimension $k$, which is impossible).

    We show that $X^{\geq (\ino-1)}$ is a $(k-1)$-cell in $\partial X$. Pick any $x \in \partial X$ such that $X^{\geq x}$ is a $(k-1)$-cell and such that $(\ino-1) \in X^{\geq x}$. Since we assumed $(\ino-1) \in X^{\geq x}$, the framed cell $(X^{\geq x},\rest \cF x)$ must be a section cell (indeed, $X^{\geq x}$ will contain a $(k-1)$-simplex containing $(\ino-1)$ which, as we've just observed, must be a section simplex). In fact, each $(k-1)$-simplex in $X^{\geq x}$ must either contain $(\ino-1)$ or $(\ino+1)$: this follows, since taking the cone of the section $(k-1)$-simplices in $X^{\geq x}$ with cone point $\ino$ must yield spacer $k$-simplices. Observe that $X^{\geq x}$ cannot, however, contain both $(\ino-1)$ and $(\ino+1)$ (for instance, this would contradict the framed boundedness of all framed realizations of the flat framing $\rest \cF x$, see \autoref{rmk:framed-embeddings-of-flat-framings}). It follows that all $(k-1)$-simplices of the $(k-1)$-cell $X^{\geq x}$ must contain the vertex $(\ino-1)$. But this is only possible if $x = (\ino-1)$. 

    Finally, we check $\im(\gamma^-)$ contains the same simplices as $X^{\geq (\ino-1)}$. Certainly simplices of $X^{\geq (\ino-1)}$ are lowest sections and must lie in the image of $\im(\gamma^-)$. Conversely, since every simplex $X$ lies in the star of $\ino$ and thus every simplex in $X_{n-1}$ lies in the star of $p_n(\ino)$, which implies $\im(\gamma^-)$ must be contained in $X^{\geq (\ino-1)}$ as required.
\end{proof}

\nid Given a spacer cell $(X,\cF)$, we often refer to the subposets $X^{\geq (\ino\pm1)}$ of $X$, determined by the images of the sections $\gamma^\pm$, as the `upper' resp.\ `lower section cells' in $(X,\cF)$.

\subsubsecunnum{Projecting cells} We explain that spacer cells `project down' (resp.\ `up') onto their lower (resp.\ upper) section cells. This will provide the cellular projection maps necessary for the comparison of framed regular cells and proframed regular cells. The following notation will be useful. A framed regular cell $(X,\cF)$ restricts, by definition, to a framed cell $(X^{\geq x},\rest \cF x)$ on the subcell $X^{\geq x} \into X$, for any $x \in X$. We denote by $\cP^{x} = \Intfr {\rest \cF x}$ the integral proframing of the restricted $n$-framing $\rest \cF x$. Write $\cP^{x} = (p^{x}_n,p^{x}_{n-1}, ..., p^{x}_1)$ with $p^{x}_i : X^{x}_i \to X^{x}_{i-1}$. Importantly, while $X^{x}_n = X^{\geq x}_n$ is a subcomplex of $X_n$, the same need not a priori hold for $X^{x}_i$ and $X_i$ with $i < n$---the fact that it does hold, follows from the next result.

\begin{lem}[Projections restrict to cells] \label{lem:proj-restrict-on-cells} Given a framed regular cell $(X,\cF)$ with subcell $(X^{\geq x},\rest \cF x)$, to top simplicial projection $p^{x}_n$ of the integral proframing of $(X^{\geq x},\rest \cF x)$ is obtained by restricting the top simplicial projection $p_n$ of the integral proframing of $(X,\cF)$ along $X^{x}_n \into X_n$.
\end{lem}

\begin{proof} We can restrict any framed realization of $(X,\cF)$ to $X^{\geq x} \into X$ to obtain a framed realization of $(X^{\geq x},\rest \cF x)$ which, by \autoref{rmk:framed-embeddings-of-flat-framings}, is necessarily bounded. Tracing through the construction of integral proframings via framed bounded realizations in \autoref{prop:existence-of-integrating-proj-framings} then proves the claim.
\end{proof}


\begin{lem}[Projected cells] \label{obs:section-projections} For a framed cell $(X,\cF)$, with top simplicial projection $p_n : X = X_n \to X_{n-1}$ in its integral proframing, we can chose a unique poset structure on $X_{n-1}$ making $p_n$ a cellular map of posets.
\end{lem}

\begin{proof} If $(X,\cF)$ is a section cell, then $p_n$ is a simplicial isomorphism and thus we must have $p_n : X_n \iso X_{n-1}$ as cellular posets. If $(X,\cF)$ is a spacer cell, we define the poset structure on $X_{n-1}$ by setting $X_{n-1} \iso X^{\geq (\ino-1)}$, where we identify the corresponding simplicial complexes using the lower section map $\gamma^- : X_{n-1} \into X_n$. To see that $p_n : X = X_n \to X_{n-1} \iso X^{\geq (\ino-1)}$ is a cellular poset map note that $p_n$ maps each cell $X^{\geq x} \into X$ exactly to the cell $X^{\geq \gamma^- \circ p_n (x)} \into X^{\geq (\ino-1)} \iso X_{n-1}$ (which can be seen by combining \autoref{lem:proj-restrict-on-cells} and our earlier construction of section cells).
\end{proof}

\begin{rmk}[Isomorphism of upper and lower cells] \label{rmk:section-projections} If $(X,\cF)$ is a spacer cell, then the cellular poset map $p_n : X_n \to X_{n-1}$ must restrict on lower resp.\ upper section cells $X^{\geq (\ino\pm1)} \into X$ to cellular poset isomorphisms $p_n : X^{\geq (\ino\pm1)} \iso X_{n-1}$. In particular, $X^{\geq (\ino-1)}$ and $X^{\geq (\ino+1)}$ are canonically cellular isomorphic.
\end{rmk}


\begin{obs}[Framed cell structures project] \label{obs:fr-cell-str-project} Having constructed a cellular poset structure on the simplicial complex $X_{n-1}$ (which is a regular cell, since $X_{n-1} \iso X^{\geq (\ino+1)}$ is a regular cell) we may now further endow the poset $X_{n-1}$ with a framing. Recall our construction of the projected framing $\cF_{n-1}$ of $\cF$ (see \autoref{term:projected-framings}), which equips the simplicial complex $X_{n-1}$ with an $(n-1)$-framing. Note that this $(n-1)$-framed simplicial complex coincides (up to increasing the embedding dimension of frame labels by postcomposing with $\bnum {n-1} \into \bnum n$) with the $n$-framed simplicial complex $(X^{\geq (\ino\pm1)},\rest \cF {\ino\pm1})$. Therefore, considering $X_{n-1} \iso X^{\geq (\ino+1)}$ as cellular posets, since $\rest \cF {\ino\pm1}$ defines an $n$-framing of the latter cell, it follows $\cF_{n-1}$ defines an $(n-1)$-framing of the former cell.
\end{obs}

\begin{term}[Projected cells] \label{term:projected-cells} We call the $(n-1)$-framed regular cell $(X_{n-1},\cF_{n-1})$ the `projected cell' of the $n$-framed cell $(X,\cF)$ and the cellular map $p_n : X_n \to X_{n-1}$ its `cell projection'.
\end{term}

\nid Note that if $(X,\cF)$ is a section cell, then $X_{n-1} \iso X$ canonically. The next example illustrates the projected cell of a spacer cell.

\begin{eg}[Projected cells] \label{eg:projected-cells} In \autoref{fig:projected-cell-construction} we illustrate the projected cell of the 3-framed regular cell $(X,\cF)$ given in \autoref{fig:section-and-spacer-cells}, highlighting the `projected' poset structure on $X_{2}$ constructed by \autoref{obs:section-projections} in green.
\begin{figure}[ht]
    \centering
    \def\svgwidth{1\columnwidth}
    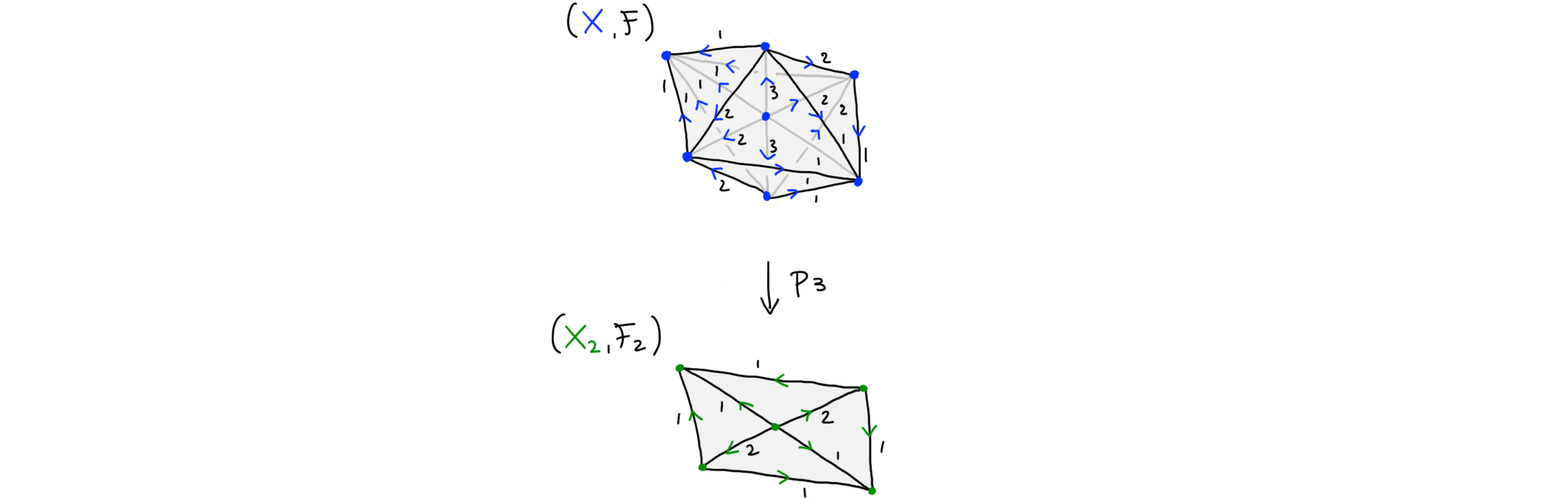

    \caption{The projected 2-framed cell of a 3-framed regular cell.}
    \label{fig:projected-cell-construction}
\end{figure}
\end{eg}

\subsubsecunnum{Existence of integral proframings for framed cells} Finally, we construct and verify uniqueness of integral proframings of framed cells.

\begin{lem}[Existence of integral $n$-proframings of framed cells] \label{prop:existence-of-integrating-proj-framings-reg-cell} Any $n$-framed regular cell $(X,\cF)$ has an integral $n$-proframing $(X,\cP)$.
\end{lem}

\begin{proof} Recall, the definition of $n$-proframed cells first requires the construction of a tower $\cP$ of cellular maps $X = X_n \xto {p_n} X_{n-1} \to ... \xto {p_1} X_0$. We inductively construct the cellular posets $X_i$ resp.\ the surjective cellular maps $p_i$ as `projected cells' and `cell projections' in the sense of \autoref{term:projected-cells}. Further, by the above constructions, the integral proframing $\Intfr \cF$ of the simplicial framing $\cF$ endows the tower $\cP$ with the structure of a simplicial proframing as required in the definition of proframed cell. Finally, the fact that $\cF$ is flat on each cell implies that $\cP$ is flat on each cell, which completes the definition of the integral proframing $\cP$ of the cell $X$.
\end{proof}

\begin{lem}[Essential uniqueness of integral $n$-proframings of framed cells] Integrating $n$-proframings of framed cells are essentially unique.
\end{lem}

\begin{proof} Let $(X,\cP)$ and $(Y,\cQ)$ be $n$-proframed regular cells with isomorphic gradient framed cells. Then the unique isomorphism $\cP \iso \cQ$ can be constructed as in \autoref{prop:uniqueness-of-integrating-proj-framings}. The fact that this simplicial isomorphism of towers is also a cellular isomorphisms of towers can be seen inductively using the relation of projected cells with lower (or upper) section cells (see \autoref{obs:fr-cell-str-project}).
\end{proof}


\begin{constr}[The integration functor of framed regular cells] \label{constr:integration-framed-reg-cell} We construct the functor $\Intfr : \FrCCell n \to \PFrCCell n$, called the integration functor of framed regular cells. On objects, this maps $n$-framed regular cells $(X,\cF)$ to their essentially unique integral $n$-proframed regular cells,  denoted by $(X,\Intfr \cF)$, which was shown to exist above. On morphisms, it maps a framed cellular map $F : (X,\cF) \to (Y,\cG)$ to the proframed cellular map $\Intfr F : (X,\Intfr \cF) \to (Y,\Intfr \cG)$ determined by setting the component $(\Intfr F)_n$ to equal $F$---the fact that this descents to a tower of cellular poset maps can be seen inductively the relation of projected cells with lower (or upper) section cells (see \autoref{obs:fr-cell-str-project}). Comparing definitions, one further checks that framedness of $F$ implies proframedness of $\Intfr F$ as required.
\end{constr}

Finally, our results assemble into the following proof of \autoref{obs:frames-vs-proj-frames-for-reg-cells}, which stated that the gradient functor $\Gradfr : \PFrCCell n \to \FrCCell n$ is an equivalence of categories.

\begin{proof}[Proof of \autoref{obs:frames-vs-proj-frames-for-reg-cells}] The inverse of the `gradient' functor $\Gradfr$ is given by the `integration' functor $\Intfr$ as defined in \autoref{constr:integration-framed-reg-cell}.
\end{proof}

As a consequence of the tools and observations described in this section, we may now revisit our earlier claim in \autoref{obs:final-frame-vec-is-vec} that the `final frame vector $\svec {\ino}$ is in fact a `vector', or more precisely, the entrance path poset of a 1-simplex.

\begin{cor}[Characterizing final frame vectors] \label{cor:final-frame-vec-is-vec} The final frame vectors $\svec {\ino}$ of any $n$-framed regular $k$-cell $(X,\cF)$, $k > 0$, span a subposet canonically isomorphic to the entrance path poset $\Entr [1] \iso (- \ot 0 \to +)$ of a 1-simplex.
\end{cor}
\begin{proof} The claim is trivial if $n = 1$. If $(X,\cF)$ is a spacer simplex then $\svec {\ino} \iso ((\ino-1) \ot \ino \to (\ino+1))$. If $(X,\cF)$ is a section, then $\svec {\ino} \iso \svec {\ino_{n-1}}$ where $\ino_{n-1}$ is the initial object of the projected cell $(X_{n-1},\cF_{n-1})$ and thus the statement follows inductively.
\end{proof}

\subsubsecunnum{Existence of integral proframings for flat framed cell complexes}

Building on the previous two sections, we now indicate the proof of \autoref{obs:frames-vs-proj-frames-for-reg-diags}: recall, this states that the gradient functor $\Gradfr : \PFrCDiag n \to \FrCDiag n$ is an equivalence of categories. Since the proof closely follows the corresponding arguments for simplicial complexes and for cells, we will only outline its main steps.

\begin{proof}[Proof of \autoref{obs:frames-vs-proj-frames-for-reg-diags}]  As before, we want to construct an inverse $\Intfr$ to $\Gradfr$. Let $(X,\cF)$ be a flat $n$-framed regular cell. Set $\cP$ to be the integral proframing of the framed simplicial complex $(X, \cF)$. Analogously to \autoref{obs:upper-and-lower-sections-in-reg-cells}, define upper and lower sections $\gamma^\pm : X_{n-1} \to X_n$. In analogy to \autoref{lem:section-images-are-cells}, one shows that images of lower (resp.\ upper) sections are in fact exactly simplices in an upper-closed subposets $X_-$ (resp.\ $X_+$) of $X$. In analogy to \autoref{lem:proj-restrict-on-cells} one finds that $p_n : X_n \to X_{n-1}$ induces a poset map $p_n : X \to X_{n-1}$, where we can identify $X_{n-1} \iso X_-$ via $\gamma^-$ (or similarly $X_{n-1} \iso X_+$ via $\gamma^+$). As before, the projected framing $\cF_{n-1}$ now yields another flat framed regular cell complex $(X_{n-1}, \cF_{n-1})$. Arguing inductively, we can construct an integral proframing for $(X_{n-1},\cF_{n-1})$, and augmenting it with the cellular map $p_n$ we obtain the integral proframing $(X,\cP)$ of $(X,\cF)$, as required. One checks essential uniqueness as before. The construction of $\Intfr$ on morphisms proceeds similarly to \autoref{constr:integration-framed-reg-cell}.
\end{proof}

\section[Equiv.\ of framed regular cell complexes \& regular block sets]{Equivalence of framed regular cell complexes \& regular block sets} \label{sec:equiv-of-frccplx-and-trusses}

We now prove the equivalence of framed regular cells and truss blocks, and, more generally, of flat framed regular cell complexes and trusses, and \emph{yet} more generally, of framed regular cell complexes and truss block sets. A first ingredient in these proofs is the equivalence of flat framed and flat proframed regular cell complexes, described in the previous section, which will allow us to think of flat framed regular cell complexes as `towers' of surjective cellular maps. A second important ingredient concerns the structure of trusses: namely, in \autoref{ssec:PL-cellularity-of-trusses} we will show that face posets of closed trusses are in fact cellular posets (moreover, they are \emph{PL} cellular posets). Equipped with both of these results, in \autoref{sec:truss-translation}, we first construct the correspondence of flat framed cell complexes and trusses, and then generalize this to a correspondence of framed cell complexes and truss block sets.

\subsection{Piecewise linear cellularity of closed trusses} \label{ssec:PL-cellularity-of-trusses}

In this section, we prove that face posets of trusses are PL cellular posets. Recall, this means that the strict upper closure of any element in a poset geometrically realizes to a PL sphere.

\begin{lem}[Cellularity of closed trusses] \label{lem:cellularity-of-closed-trusses} Given a closed $n$-truss $T$, each face poset $(T_i,\eleq)$ is a PL cellular poset.
\end{lem}

\nid While there are several direct ways of proving the result, our approach will be based on the `shellability' properties of trusses. Our proof will rely on the following result.

\begin{prop}[See {\cite[Prop. 4.5 ff.]{bjorner1984posets}}] \label{prop:bjoerner} The realization $\abs {X}$ of a poset $X$ is a regular cell complex PL homeomorphic to the PL $m$-sphere if $X$ is pure of dimension $m$, shellable, and thin. \qed
\end{prop}

\begin{term}[Purity, shellability, and thinness] \label{term:pure-shellable-thin} A simplicial complex is called `pure of dimension $m$' if its facets (that is, non-degenerate simplices which are not the face of any other non-degenerate simplex) are all of the same dimension $d$. Similarly, a poset $X$ is called pure of dimension $m$ if its nerve $NX$ is pure of dimension $m$.

    A poset $X$ is called `shellable' if the simplicial complex $NX$ is pure of dimension $m$ and its facets admit an ordering $K_0, K_1, K_2, ..., K_j$, such that, for all $0 < k \leq j$, the subcomplex $(\cup_{i < k} K_i) \cap K_k$ of $X$ (obtained by intersecting of the union $K_i$, $i < k$, with $K_k$) is a pure simplicial complex of dimension $(m-1)$.

Finally, a poset $X$ is called `thin' if for every non-refinable length-2 chain $x < y < z$ in $X$ there is exactly one $y' \neq y$ such that $x < y' < z$ (this is also sometimes called the `diamond property').
\end{term}

\begin{proof}[Proof of \autoref{lem:cellularity-of-closed-trusses}] Recall the following notation: given a poset $(X,\leq)$ and an element $x \in X$, the `upper closure' $X^{\geq x}$ (resp.\ `strict upper closure' $X^{> x}$) is the subposet of $X$ consisting of all elements greater than $x$ (resp.\ those that are strictly greater than $x$). Given an $n$-truss $T$ and an element $x \in T_i$ in a face poset $(T_i,\eleq)$ of $T$, for the proof of \autoref{lem:cellularity-of-closed-trusses} we need to show that the strict upper closure $T_i^{\egre x}$ realizes to a PL $m$-sphere. Inductively, we may assume this holds for $i < n$ and only consider strict upper closures $T_n^{\egre x}$ in the top level poset. Moreover, by passing to the face block of $T_n^{\egeq x}$ of $x$, we may assume $T$ is in fact an $n$-truss block, with initial element $\ino \equiv x$. What is left to show is that the truss block's `boundary' $\partial T_n \equiv T_n^{\egre \ino}$ realizes to a PL $m$-sphere. By \autoref{prop:bjoerner} this is equivalent to showing that $\partial T_n$ is pure, thin and shellable.

    Denote by $\ino_{n-1} = p_n(\ino)$ the projection of $\ino$ to $T_{n-1}$. Denote by $T_{<n}$ the $(n-1)$-truncation of $T$. Since $T$ is an $n$-truss $k$-block, note that $T_{<n}$ is an $(n-1)$-truss $l$-block (note $l \in \{k-1,k\}$). Arguing by induction, we assume the statement of the lemma to hold for $T_{<n}$; that is, we assume that the strict upper closure $\partial T_{n-1}$ of $\ino_{n-1}$ is thin, shellable and pure of dimension $l - 1$. There are now two basic cases.  Either, $\ino$ is singular in the fiber over $\ino_{n-1}$ or it is regular. In the first case, the 1-truss bundle $p_n : T_n \to T_{n-1}$ is an isomorphism of face posets (and $k = l$) so all claimed properties of $\partial T_n$ follow from that of $\partial T_{n-1}$. Thus, it remains to prove the lemma in the second case, that is, for $\ino$ being regular in the fiber over $\ino_{n-1}$. Note that initiality of $\ino$ implies the fiber $p_n\inv(\ino_{n-1})$ must be of the form $(\ino-1) \egre \ino \eles (\ino+1)$. We separately prove purity, shellability and thinness separately in this case.

    \emph{Purity}. We first show $\partial T_n$ is pure. Observe that facets of $T_n$ must be mapped to facets of $T_{n-1}$ by $p_n$ (this follows since $p_n$ is surjective on simplices, and since simplex fiber transition maps are surjective too, see \autoref{prop:1-truss-bun-euclidean}). Each facet in $T_{n-1}$ must contain the vertex $\ino_{n-1}$. Since the fiber over $\ino_{n-1}$ has spacers, so must the fiber of each facet in $T_{n-1}$. Consequently, facets in $T_n$ must themselves be spacers. Since all facets in $T_{n-1}$ have dimension $(k-1)$, it follows that facets of $T_n$ have dimension $k$. Note that this further implies that all facets of $\partial T_n$ are of dimension $k-1$.

    \emph{Shellability}. We next show $\partial T_n$ is shellable. Observe first that facets in $\partial T_n$ are either first or last sections lying over facets in $T_{n-1}$, or they are spacers lying over facets in $\partial T_{n-1}$. Inductively, construct a shelling $K_1, K_2, ..., K_{\#T_{n-1}}$ of facets of $\partial T_{n-1}$ (where $\#T_{n-1}$ is the number of facets in $\partial T_{n-1}$). This also induces a shelling $K_\bullet = (K^{\ino}_1, K^{\ino}_2, ... K^{\ino}_{\#T_{n-1}})$ of $T_{n-1}$ (where $K^{\ino}_i$ is obtained from $K_i$ by adjoining a new first vertex $\ino_{n-1}$). Now build a shelling $L_\bullet = (L_1,L_2, ... L_{\#T_n})$ of $\partial T_n$ (where $\#T_n$ is the number of facets in $\partial T_n$) in the following three steps.
    \begin{enumerate}
        \item \textit{Lower section shelling}: We define the first $\#T_{n-1}$ facets $$L_1, L_2, ... , L_{\#T_{n-1}}$$ in the sequence $L_\bullet$, by setting $L_i$ to be the bottom section (see \autoref{defn:bottom-top-sections}) lying over $K^{\ino}_i$. Note that this satisfies $L_i(0) = (\ino-1)$.

        \item \textit{Side shelling}:  We next define the subsequence $$L_{\#T_{n-1}+1},L_{\#T_{n-1}+2}, ... ,L_{\#T_m - \#T_{n-1}}$$ of $L_\bullet$ to be the sequence $L_{(1,1)}, L_{(1,2)}, ... L_{(1,j_0)}$, $L_{(2,1)},L_{(2,2)}, ..., L_{(2,j_1)}$ , $...$ , $L_{(\#T_{m-1},1)}, ..., L_{(\#T_{m-1},j_1)}$, where $L_{(i,j)}$ is the $j$th spacer (in scaffold order) lying over $K_i$.

        \item \textit{Upper section shelling}: Finally, we define the last $\#T_{n-1}$ facets $$L_{\#T_m - \#T_{n-1} + 1}, L_{\#T_m - \#T_{n-1} + 2}, ... , L_{\#T_{n-1}}$$ in the sequence $L_\bullet$ by setting $L_{\#T_m - \#T_{n-1} + i}$ to be the top section lying over $K^{\ino}_i$. Note that this satisfies $L_{\#T_m - \#T_{n-1} + i}(0) = (\ino+1)$.

    \end{enumerate}

    This constructs a shelling of $\partial T_n$. We illustrate an instance of the construction in \autoref{fig:truss-shelling-example}: for the 3-truss $T$ shown on the left, a shelling $(K_1,K_2,...,K_6)$ of $\partial T_2$ is depicted on the (lower) right. Above it, we depict the shelling $(L_1,L_2,...,L_{16})$ of $\partial T_3$ as constructed by the above procedure.
\begin{figure}[ht]
    \centering
    \def\svgwidth{1\columnwidth}
    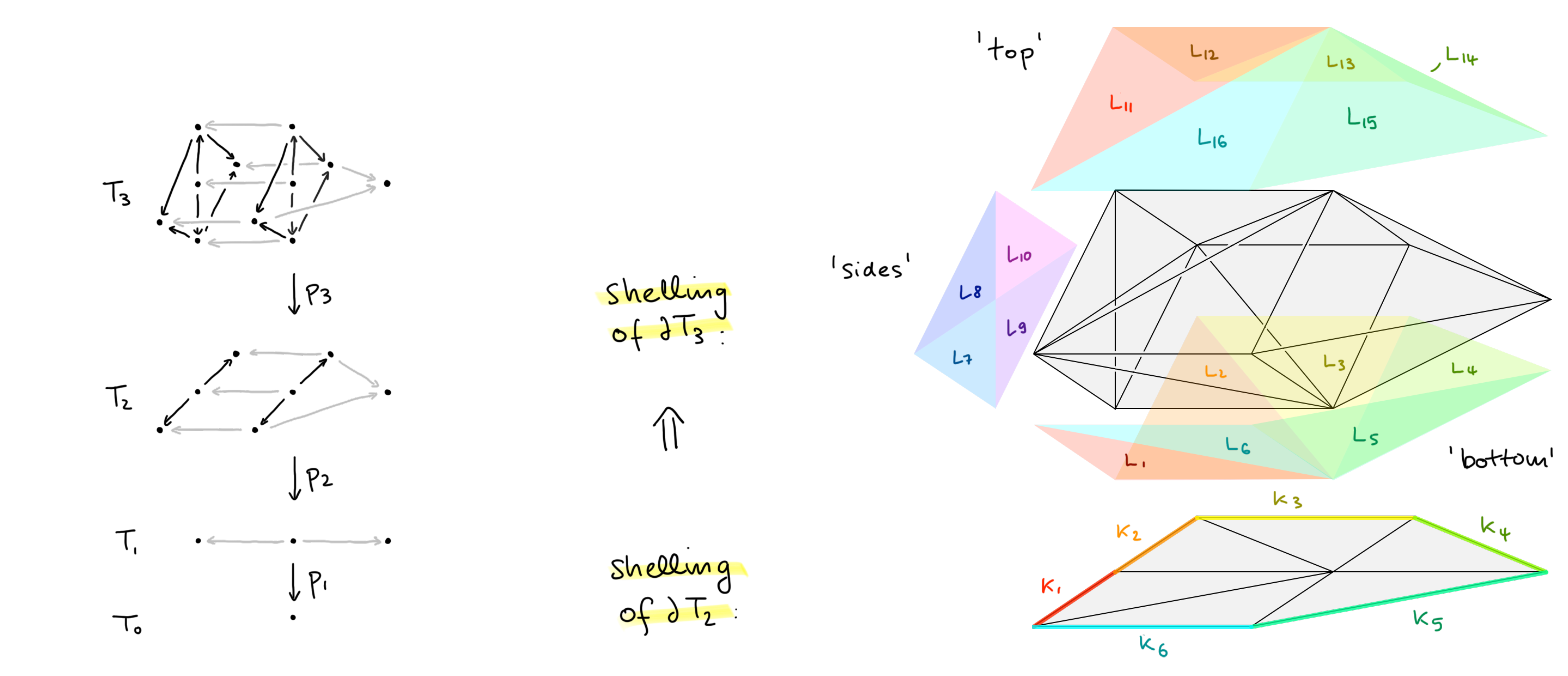

    \caption[Inductive construction of truss shellings.]{Inductive shelling of a face poset boundary $\partial T_n$, by sequentially shelling its bottom, sides, and top part.}
    \label{fig:truss-shelling-example}
\end{figure}

\emph{Thinness}. It remains to check that $\partial T_n$ is thin. We show that $T_n$ is thin which implies thinness of $\partial T_n$ (since $T_n$ adjoins an initial element to $\partial T_n$). Inductively assume $T_{n-1}$ is thin. Take a 2-simplex $K : [2] \to T$ such that the chain $K = (x \to y \to z)$ is non-refinable. We distinguish two cases based on the dimension of the base projection $\im(p_nK) : [j] \into T_{n-1}$, with $j \in \{1,2\}$.

First assume $j = 2$. Then the base projection $J := \im(p_nK)$ is a chain $(x_{n-1} \to y_{n-1} \to z_{n-1})$ in $T_{n-1}$. Note $J$ must be non-refinable (otherwise, $L$ would be refinable). Thinness of $T_{n-1}$ implies there is exactly one other non-refinable chain $J' : x_{n-1} \to y'_{n-1} \to z_{n-1}$. Since the 1-truss bordism lying over the $J$ compose to the same 1-truss bordism as the 1-truss bordism lying over $J'$, there must be at least one chain $K'$ from $x$ to $z$ lying over $J'$. Moreover, there cannot be a third chain $K''$ from $x$ to $z$, since that would have to lie over either $J$ or $J'$; assume, w.l.o.g., that it lies over $J$ and that $K \fles K''$ in the scaffold order of sections over $J$ (see \autoref{sec:truss-induction}); all spacers over $J$ between $K$ and $K''$ must now have fiber morphisms in the fiber over $y_{n-1}$. Thus the $3$-spacer containing the $2$-section $K$ as its lower section, has a spine that refines $K$; contradicting that $K$ was non-refinable. It follows that a third chain $K''$ cannot exist.

Next assume $j = 1$. In this case, the base projection $J := \im(p_nK)$ is a 1-simplex $(x_{n-1} \to z_{n-1})$ in $T_{n-1}$. Thus, $K$ must be a spacer over $J$. Arguing by truss induction on the 1-truss bundle $p_n$ over $J$, we find that exactly two non-refinable chains from $x$ to $z$ must exist. Namely, either the lower or the upper section of $K$ must have a jump morphism that lies over $J$ (this follows from the arguments in the proof of \autoref{lem:section-induction}, but it can be graphically understood when thinking of the section order as a directed path through jump morphisms, see \autoref{fig:scaffold-order-on-sections}). In the former case the two non-refinable chains are given by $K$'s spine and its predecessor's spine, in the latter case, by $K$'s spine and its successor's spine.

This completes the verification that $\partial T_n$ is pure of dimension $k-1$, shellable and thin. \end{proof}

\subsection{Correspondence of flat framed regular cell complexes and trusses} \label{sec:truss-translation}

\subsubsecunnum{From flat framed regular cell complexes to trusses}

We start with a construction of truss translation in the case of flat \emph{proframed} regular cell complexes. We then consider that construction together with our earlier definition of integration, yielding the truss translation of flat \emph{framed} regular cell complexes.

\begin{constr}[Truss translation of flat proframed regular cell complexes] \label{constr:truss-translation-proframed-case} Given a flat $n$-proframed regular cell complex $(X,\cP)$ with projections $\cP = (p_n,...,p_1)$, we construct a closed $n$-truss $\kT(X,\cP)$. Recall, that each projection $p_i : X_i \to X_{i-1}$ is both a cellular poset map as well as an ordered simplicial map (see \autoref{rmk:overloading-notation-proframes}). Denote the tower of (to be constructed) 1-truss bundles $\kT(X,\cP)$ by $(T_n \xto {p_n} T_{n-1} \xto {p_{n-1}} ... \xto {p_1} T_0)$. In order to construct the 1-truss bundles $p_i : T_i \to T_{i-1}$ we need to: (1) define poset maps $p_i : (T_i,\eleq) \to (T_{i-1},\eleq)$ of face orders, (2) define a dimension functor $\dim : (T_i,\eleq\op) \to [1]$, (3) define a frame order $(T_i,\fleq)$, (4) verify that fibers of $p_i$ over objects are closed 1-trusses, and (4') verify that fibers over morphisms are 1-truss bordism.
\begin{enumerate}[wide, labelwidth=!, labelindent=12pt]
        \item We set $(T_i,\eleq) = X_i$ and define $p_i : (T_i,\eleq) \to (T_{i-1},\eleq)$ to be $p_i : X_i \to X_{i-1}$.
        \item We define $\dim : (T_i,\eleq\op) \to [1]$ to map $x \in T_i$ to $0$ if $X^{\geq x}$ is a section cell and to $1$ if $X^{\geq x}$ is a spacer cell in $X$ (see \autoref{term:section-and-spacer-cells}). Since section cell can only border on other section cells, this defines a poset map $(T_i,\eleq\op) \to [1]$ as required.
        \item We define two elements $x, y$ in $T_i$ to be related in the frame order $(T_i,\fleq)$ such that $x \fles y$ if and only if they lie in the same fiber of $p_i$ and there is a sequence of 1-simplices $x \to \cdots \to y$ in the simplicial complex $X_i$.
        \item We check that the structures $\eleq$, $\dim$, and $\fleq$ restrict on fibers $p_i\inv(z)$ over objects $z$ in $T_{i-1}$ to yield closed 1-trusses $T_z = (p_i\inv(z),\eleq,\dim,\fleq)$. Since $\cP$ is flat, the fiber over $z$ is a linear subcomplex $x_0 \to x_1 \to ...\to x_k$ of the simplicial complex $X_i$; this shows, $(p_i\inv(z),\fleq)$ is a total order. Moreover, if $X^{\geq x_i}$ is a spacer cell, then $X^{\geq x_{i \pm 1}}$ must be its upper and lower section cells; thus $(p_i\inv(z),\eleq)$ is a linear fence. Our definition of the functor $\dim$ then implies that $T_z$ is a closed 1-truss as claimed (spacer cells have both upper and lower section cells, and thus cannot lie at fiber endpoints).
        \item[(4')] It remains to verify that fibers over morphisms are 1-truss bordisms. Let $f : z \to w$ be an arrow in $T_{i-1}$. Denote by $R: T_z \xslashedrightarrow{} T_w$ the relation $p_i\inv(f)$. Since $T_i = X_i$ is a poset, $R$ is a boolean profunctor. Since the proframing $\cP$ is flat, it follows that $R \subset (T_z,\fleq) \times (T_w,\fleq)$ is bimonotone. Since the transition functors are surjective, it follows that $R$ fully matches objects. Finally, if $x \in \sing(T_z)$ there is a unique $y \in \sing(T_w)$ such that $R(x,y)$: indeed, $p_n$ restricts on sections cells $X^{\geq x}$ to poset isomorphisms $p_x : X^{\geq x} \iso X^{\geq z}$, and thus $R(x,y)$ holds if and only if $X^{\geq y} = p_x\inv(X^{\geq w})$. The statement that $R$ is a 1-truss bordism now follows from \autoref{cor:full-matching}.
\end{enumerate}
This completes the construction of the closed $n$-truss $\kT(X,\cP)$.
\end{constr}

\begin{constr}[Truss translation of proframed maps] \label{constr:truss-translation-maps} Given a proframed cellular map $F : (X,\cP) \to (Y,\cQ)$ of flat proframed regular cell complexes we construct the singular truss map $\kT F : \kT(X,\cP) \to \kT(Y,\cQ)$. The $i$th component of $\kT F$ is defined by the $i$th component $F_i : X_i \to Y_i$ of $F$. We need to check that $F$ being a proframed cellular map translates to $\kT F$ being a singular map, i.e.\ $\kT F$ maps singular objects to singular objects. Indeed, singular objects $x \in T_i$ yield section cells $X^{\geq x}_i$ in $X_i$ and the definition of proframed maps implies that $F_i$ must map section cells (in $(X_i,\cP_{\leq i})$) to section cells (in $(Y_i,\cQ_{\leq i})$). It follows that $\kT F_i$ preserves singular objects as required.
\end{constr}

\begin{term}[Truss translation functors] \label{constr:truss-translation-functor} \label{defn:truss-translation-flat-framed} \label{defn:truss-translation-framed-cell} We denote the `truss translation functor' of flat proframed regular cell complexes, given by \autoref{constr:truss-translation-proframed-case} and \autoref{constr:truss-translation-maps}, by $\kT : \PFrCDiag n \to \sctruss n$.

    We define the `truss translation of flat framed regular cell complexes' to be the functor
    \begin{equation}
        \kT : \FrCDiag n \to \sctruss n
    \end{equation}
obtained as the composite of  the functors $\Intfr : \FrCDiag n \to \PFrCDiag n$ and $\kT : \PFrCDiag n \to \sctruss n$.

    We define the `truss translation of framed regular cells' to be the functor
    \begin{equation}
        \kT : \FrCCell n \to \blcat n
    \end{equation}
obtained by restricting the functor $\kT : \FrCDiag n \to \sctruss n$ to the full subcategory $\FrCCell n \into \FrCDiag n$.
\end{term}

\subsubsecunnum{From trusses to flat framed regular cell complexes} \label{sec:framed-complex-translation}

Inverse to truss translation, we next construct the framed complexes translations of $n$-truss blocks and closed $n$-trusses.

\begin{constr}[Proframed complex translation of closed trusses] \label{constr:proframed-complex-translation-truss} Given a closed $n$-truss $T = (T_n \xto{p_n} T_{n-1} \xto{p_{n-1}} ... \xto{p_1} T_0)$ we construct a flat $n$-proframed regular cell complexes $\kX T$. We write this proframed complex as $(X,\cP)$ with $X$ a cellular poset and $\cP = (p_n,...,p_1)$ a tower of (both cellular and ordered simplicial) projections. In order to define $(X,\cP)$, we need to: (1) define the sequence $\cP$ of cellular maps $p_i : X_i \to X_{i-1}$ (with $X = X_n$), (2) define an ordering of the underlying simplicial sequence $\cP$, making $\cP$ a simplicial proframing, and (3) verify that $\cP$ is both a flat proframing itself and flat on each cell of $X$.
    \begin{enumerate}[wide, labelwidth=!, labelindent=12pt]
        \item We define the poset map $p_i : X_i \to X_{i-1}$ to be $p_i : (T_i,\eleq) \to (T_{i-1},\eleq)$. The fact that each $X_i$ is cellular was proven in \autoref{lem:cellularity-of-closed-trusses}. The fact that the bundle maps $p_i$ are cellular maps follows since 1-truss bundle have lifts (see \autoref{rmk:1-truss-bundle-lifts}).
        \item We define the ordered simplicial maps $p_i : X_i \to X_{i-1}$. Inductively, we assume to have defined the ordered simplicial complex $X_{i-1}$. To obtain the ordered simplicial complex $X_i$, we need to consistently order vertices of all $k$-simplices $x : \unsimp k \into X_i$ in the unordered simplicial complex $X_i$. Take $k = 1$. Either the 1-simplex $x$ lies over an object $y$ of $X_{i-1}$, or it lies over a (ordered) 1-simplex $z : z(0) \to z(1)$ in $X_{i-1}$. In the first case, order $x = x(0) \to x(1)$ such that $x(0) \fles x(1)$ in the frame order $(T_i,\fleq)$. Otherwise, if $x$ lies over a 1-simplex $y : y(0) \to y(1)$, then order $x = x(0) \to x(1)$ such that its first vertex $x(0)$ lies over $y(0)$ while its second vertex $x(1)$ lies over $y(1)$. This determines a vertex order on all 1-simplices, which extends an ordering of $X_i$. The construction also entails that the poset map $p_i : X_i \to X_{i-1}$ is an ordered simplicial map $p_i : X_i \to X_{i-1}$.
        \item Finally, we verify that $\cP$ is both flat itself and flat on each cell of $X$. Recall, flatness requires fiber categories to be linear and transition functors to be endpoint-preserving. But these are exactly the conditions verified via truss induction in \autoref{prop:1-truss-bun-euclidean}. Thus, $\cP$ is flat. Similarly it follows that all cell restricted proframings $(X^{\geq x},\rest \cP x)$ are flat (namely, by applying \autoref{prop:1-truss-bun-euclidean} to the 1-truss bundles in the $n$-truss block $T^{\egeq x}$). \qedhere
    \end{enumerate}
\end{constr}

\begin{constr}[Proframed complex translation of singular maps] \label{constr:proframed-complex-translation-maps} Given a singular truss map $F : T \to S$ we construct a proframed cellular map $\kX F : \kX T \to \kX S$. Write the complex $\kX T$ as $(X,\cP)$ (consisting of a cellular poset $X$ and a proframing $\cP)$ and similarly write the complex $\kX S$ as $(Y,\cQ)$. The $i$th component $(\kX F)_i : X_i \to Y_i$ is defined by the $i$th component $F_i : T_i \to S_i$. We need to check that $\kX F$ is cellular and proframed. Both can be seen by an inductive argument. For cellularity, inductively assume $F_{i-1}$ is cellular. To see that $F_i$ is cellular, we check that, for each $x \in X_i$, the cell $X^{\geq x}_i$ in $X_i$ maps onto some cell $Y^{\geq y}_i$ in $Y_i$. If $x$ is singular (i.e.\ $X^{\geq x}_i$ is a section cell) then this claim follows from our inductive assumption that $F_{i-1}$ is cellular, and since $F_i$ preserves singular objects. If $x$ is regular (i.e.\ $X^{\geq x}_i$ is a spacer cell) the claim follows from the singular case since spacer cells are `sandwiched' between their lower and upper section cells.

    Next, to verify that $\kX F$ is proframed, we may inductively assume that the truncation $\kX F_{<n} : \kX T_{<n} \to \kX S_{<n}$ is proframed. Since $F$ is a singular truss map it follows that $\kX F_n$ maps section cells to sections cell, which (together with the inductive assumption) implies that $\kX F$ preserves final frame vectors as required in the definition of cellular proframed maps.
\end{constr}

\begin{term}[Proframed complex translation] \label{defn:framed-complex-translation-truss} \label{defn:framed-complex-translation-block}
We denote the `proframed complex translation functor', as given by \autoref{constr:proframed-complex-translation-truss} and \autoref{constr:proframed-complex-translation-maps}, by $\kX : \sctruss n \to \PFrCDiag n$.

We define the `framed complex translation of closed trusses' to be the functor
    \begin{equation}
        \kX : \sctruss n \to \FrCDiag n
    \end{equation}
obtained by taking the composite of the functor $\kX : \sctruss n \to \PFrCDiag n$ and the functor $\Gradfr : \PFrCDiag n \to \FrCDiag n$.

We define the `framed complex translation of truss blocks' to be the functor
    \begin{equation}
        \kX : \blcat n \to \FrCCell n
    \end{equation}
obtained by restricting the functors $\kX : \sctruss n \to \FrCDiag n$ to the full subcategory $\blcat n \into \sctruss n$.
\end{term}

\pauseae

As a last formality, we record that truss translation and framed complex translation assemble into the following equivalences of categories.

\begin{proof}[Proof of \autoref{thm:classification-of-cells} and \autoref{thm:classification-of-flat-framed-cell-cplx}] The functor $\kT$ defined in \autoref{defn:truss-translation-flat-framed} and the functor $\kX$ defined in \autoref{defn:framed-complex-translation-truss} are inverse; moreover, there is a unique choice for the natural isomorphisms $\id \iso \kT\circ\kX$ and $\id \iso \kX\circ\kT$. The equivalence further restricts to an equivalence of the subcategories $\FrCCell n \into \FrCDiag n$ and $\blcat n \into \sctruss n$.
\end{proof}

\begin{proof}[Proof of \autoref{thm:classification-of-framed-cell-cplx}] We show that the category of regular truss block sets (i.e.\ of presheafs on truss blocks which are regular) is equivalent to the category of framed regular cell complexes. Having shown the equivalence of categories of truss blocks $\blcat n$ and of framed regular cells $\FrCCell n$, it will suffice to show that the category of framed regular cell complexes is equivalent to a category of `regular' presheafs on framed regular cells, which can be seen as follows. Given a presheaf $W \in \PSh(\FrCCell n)$, a cell $y \in W(Y,\cG)$ is said to be non-degenerate if the map $y : (Y,\cG) \to W$ does not factor through a framed cellular surjection $(Y,\cG) \epi (Y',\cG')$ (except the identity). We say the presheaf $W$ is regular if each non-degenerate cell $y \in W(Y,\cG)$ includes back into $W$ by an injective presheaf map $y : (Y,\cG) \into W$. Every regular presheafs $W$ defines a framed regular cell complex $(X_W,\cF_W)$, whose cellular poset $(X_W,\leq)$ is the poset of all non-degenerate cells $y \in W(Y,\cG)$ such that $y \leq y'$ whenever $y$ factors through $y'$ by a framed cellular injection, and whose framing $\cW$ restricts on each cell $X_W^{\geq y} \iso Y$ to the framing $\cG$ of the framed regular cell $(Y,\cG)$. Conversely, every framed regular cell complexes $(X,\cF)$ defines a regular presheaf $W_{X,\cF}$ such that $W_{X,\cF}(Y,\cG)$ is the set of all framed cellular inclusion $(Y,\cG) \into (X,\cF)$. The constructions are inverse up to canonical isomorphism. The correspondence canonically extends to a correspondence of maps between regular presheafs and framed cellular maps which provides the claimed equivalence of categories.
\end{proof}

\begin{cor}[Framed regular cells are shellable PL cells] \label{lem:block-shellability} Given an $n$-framed regular cell $(X,\cF)$ then cellular poset $X$ is PL cellular, and both $X$ and $\partial X$ are shellable.
\end{cor}
\begin{proof} By \autoref{lem:cellularity-of-closed-trusses} we know that face posets of $\kT(X,\cF)$ are PL cellular. Since the poset $X$ is isomorphic to the top face poset of $\kT(X,\cF)$ (see \autoref{constr:truss-translation-proframed-case}) it follows that $X$ itself is PL cellular.
\end{proof}

\chapter{Constructible framed topology: meshes} \label{ch:meshes}


In this chapter we introduce the stratified topological notion of `meshes', which are towers of constructible stratified bundles whose fibers are points or framed stratified 1-manifolds. We begin in \autoref{sec:1-to-n-meshes} by describing these fibers, called `1-meshes', and their constructible bundles, before then generalizing these to notions of $n$-meshes and $n$-mesh bundles. In \autoref{sec:meshes-and-trusses} we discuss how meshes are the `topological counterpart' of trusses, by constructing the \emph{entrance path $n$-trusses} of $n$-meshes and the \emph{classifying $n$-meshes} of $n$-trusses, and showing that these constructions exhibit a `weak equivalence' between meshes and trusses.

\section{1-Meshes, 1-mesh bundles, and $n$-meshes}
\label{sec:1-to-n-meshes}

We introduce 1-meshes and $n$-meshes as certain structured stratifications (a recollection of basic ideas in the theory of stratified spaces can be found in \autoref{app:stratifications}). 1-Meshes will be defined as a stratified $0$- or 1-manifolds endowed with the structure of a flat framing, while $n$-meshes will be defined by iterating `constructible stratified bundles' of 1-meshes. Here, the notion of `stratified bundle' is a generalization of the ordinary notion of fiber bundles to the context of stratified topology (see \autoref{defn:stratified-bundle}). The term `constructible' bundle refers, generally, to bundles which can be (re)constructed from functors on the `fundamental category' of their base: while in the context of stratified topology this fundamental category is usually the entrance path $\infty$-category of the base stratification (see \autoref{defn:tentr}) in case of meshes constructibility holds already at the level of entrance path posets (which will follow from the comparison of meshes and trusses in \autoref{sec:meshes-and-trusses}).

The section will be organized as follows. In \autoref{ssec:1-meshes} we will introduce meshes and their maps; our definitions will closely mirror the definition of 1-trusses and 1-truss maps in \autoref{ch:trusses}. We then discuss bundles of 1-meshes and introduce a condition for their constructibility in \autoref{ssec:1-mesh-bundles}. Finally, in  \autoref{ssec:n-meshes}, we define $n$-meshes and briefly comment how this further gives rise to a notion of `$n$-mesh bundles'.

\subsection{1-Meshes} \label{ssec:1-meshes}

\subsubsecunnum{1-Meshes as 1-framed stratified manifolds}

We start with a notion of 1-framing of manifolds. Classically, a `tangent framing' of a smooth manifold is a trivialization of its tangent bundle. In the absence of smooth structures, we can similarly require trivializations of tangent microbundles. However, independent of whether we work in the smooth or topological case, the following holds true: every (sufficiently nice) codimension-$k$ submanifold of a framed manifold inherits a framing of its $k$-stabilized tangent bundle. Specifically, in the case $n = 1$, this motivates the following definitions.

\begin{term}[Manifolds] In the following, `manifold' will mean connected topological manifold with or without boundary. $\lR$ will mean the `standard euclidean 1-space' and $S^1$ the `oriented standard circle'.
\end{term}

\begin{defn}[1-Framings on manifolds] \label{defn:1-framings-on-mflds} A \textbf{1-framing of a manifold $M$} is an embedding $\gamma : M \into S^1$.
\end{defn}

\begin{rmk}[Dimension of 1-framed manifolds] Given a 1-framed manifold $(M,\gamma)$ if $\dim(M) = 0$ then $M$ obtains a `normal 1-frame' from the ambient $S^1$ and otherwise a `tangent 1-frame'.
\end{rmk}

\nid We will be particular interested in the case of `flat framings' where the ambient manifold is standard $\lR^n$.

\begin{defn}[Flat 1-framings on manifolds] \label{defn:flat-1-framings-on-mflds} A \textbf{flat 1-framing of a manifold $M$} is a bounded embedding $\gamma : M \into \lR$.
\end{defn}

\begin{rmk}[Flat 1-framings are framings] Any orientation preserving embedding $\lR \into S^1$ can be used to translate flat framings into framings by post-composition.
\end{rmk}

\begin{rmk}[The space of framings of a 1-manifold] Given a 1-manifold $M$, the subspace $\{ \gamma ~|~ (M,\gamma) \text{ is a 1-framing}\}$ of the mapping space $\Map(M,S^1)$ is homotopy equivalent to $\lZ_2$. That is, `up to homotopy' there are two framings on a 1-manifold. The same observation applies to flat 1-framings.
\end{rmk}

A mesh is a nicely stratified manifold with 1-framing: its strata are manifolds (without boundary).

\begin{defn}[General 1-meshes] \label{defn:1-mesh} A \textbf{1-mesh} $(M,f,\gamma)$ is a manifold $M$ with a finite stratification $f$, stratifying $M$ by manifolds (without boundary), together with a 1-framing $\gamma$ on $M$.
    \begin{enumerate}
        \item If $M$ has a single stratum, we call the 1-mesh \textbf{trivial}.
        \item If $M$ is a contractible 1-manifold, we call the 1-mesh \textbf{linear}.
        \item If $M$ is the circle, we call the 1-mesh \textbf{circular}. \qedhere
    \end{enumerate}
\end{defn}

\begin{eg}[General meshes] 1-Meshes of different types are shown in \autoref{fig:general-meshes}. In each case, we color 0-dimensional strata in red, and 1-dimensional strata in blue. For linear 1-meshes we depict the ambient euclidean space $\lR$; more commonly we will simply be indicate by `coordinate arrow' (indicated in green) which fixes the framing up to contractible choice. Similarly, we indicated the framing of circular meshes by a green arrow (describing the orientation of $S^1$).
\begin{figure}[ht]
    \centering
    \def\svgwidth{1\columnwidth}
    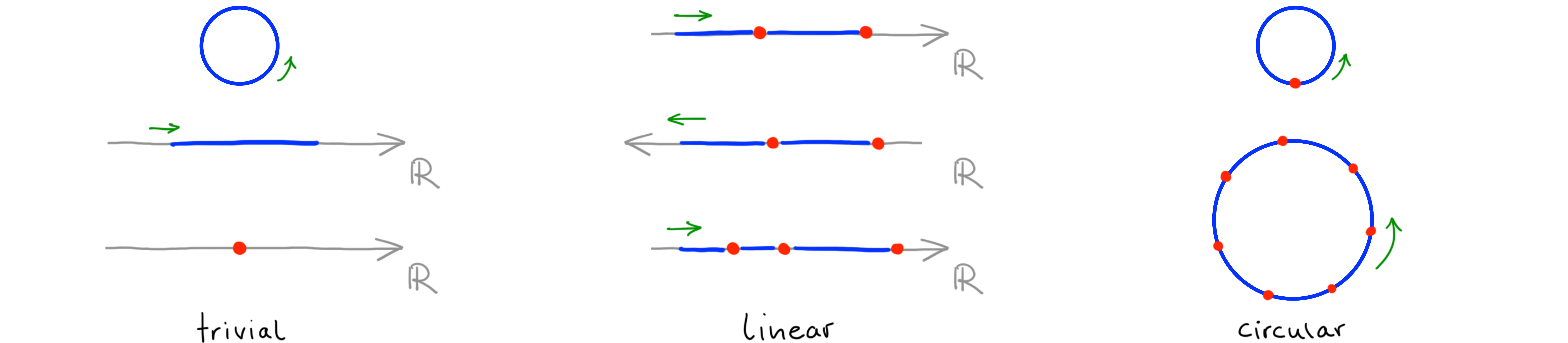

    \caption{1-Meshes of different types.}
    \label{fig:general-meshes}
\end{figure}
Note that, we can distinguish three types of trivial 1-meshes: the `trivial 0-dimensional mesh', the `trivial linear 1-dimensional mesh' and the `trivial circular 1-dimensional mesh'.\footnote{Note, in contrast, we had described only two trivial 1-trusses in \autoref{fig:1-trusses-of-different-types}. This indicate that the trivial 1-truss with one element of dimension 1 should have two distinct combinatorial incarnations: the `trivial linear' and the `trivial circular' 1-truss. We will not bother further with rectifying the combinatorial situation, since we are ultimately interested only in the linear case.}
\end{eg}

\nid While much of the theory of meshes can be developed in parallel for the `linear' and `circular' cases, our interest will (analogous to the case of 1-trusses) ultimately lie with the linear case. Going forward, we will therefore adopt the following convention.

\begin{conv}[Linear 1-meshes by default] Going forward, we will refer to `linear 1-meshes' simply as `1-meshes', and always endow 1-meshes with flat framings (i.e.\ embeddings into $\lR$).
\end{conv}

\begin{notn}[Flat 1-framing bounds] \label{notn:framing-closures} For a 1-mesh $(M,f,\gamma)$ with flat 1-framing $\gamma : M \into \lR$ denote by $\gamma^-(M)$ and $\gamma^+(M)$ the lower resp.\ upper bound of the end $\gamma(M) \subset \lR$.
\end{notn}

\begin{term}[Open and closed meshes] A 1-mesh $(M,f,\gamma)$ is called `closed' resp.\ `open' if the image $\gamma(M) \subset \lR$ is closed respectively open in $\lR$.
\end{term}

\subsubsecunnum{Maps of 1-meshes}

Let us first discuss framed maps of 1-framed manifolds (in the sense of \autoref{defn:1-framings-on-mflds}).

\begin{term}[Framed maps of standard framed spaces] A `framed map' of the standard circle $S^1$ is an orientation preserving map $S^1 \to S^1$ (by which we mean map of the form $e^{ix} \mapsto e^{i\phi(x)}$ where $\phi : \lR \to \lR$ is monotone).  A `framed map' of standard $\lR$ is an orientation preserving map $\lR \to \lR$ (by which we mean a monotone map $\lR \to \lR$). The notions similarly apply to connected subspaces of $S^1$ and $\lR^1$.
\end{term}

\begin{defn}[Framed maps of 1-framed manifolds] \label{defn:maps-of-1-meshes} For (flat) 1-framed manifolds $(M,\gamma)$ and $(N,\rho)$, a \textbf{framed map} $F : (M,\gamma) \to (N,\rho)$ is a continuous map $F : M \to N$ which induces a framed map $\tilde F : \gamma(M) \to \rho(N)$ between the images of their 1-framings, such that $\tilde F \circ \gamma = \rho \circ F$.
\end{defn}

\nid A 1-mesh map preserves both framing and stratification structure as follows.

\begin{defn}[Maps of 1-meshes] A \textbf{map of 1-meshes} $F : (M,f,\gamma) \to (N,g,\rho)$ is a continuous map $F : M \to N$ that is both a stratified map $F : (M,f) \to (N,g)$ as well as a framed map $F : (M,\gamma) \to (N,\rho)$.
\end{defn}

\nid Fully parallel to our earlier definitions of `regular', `singular', and `balanced' maps of trusses (see \autoref{defn:reg-sing-cell-1-truss-map}), we introduce the following terminology for mesh maps.

\begin{term}[Singular, regular, and balanced maps] \label{term:1-mesh-map-reg-sing} A map of 1-meshes $F : (M,f,\gamma) \to (N,g,\rho)$ is called `singular' if it maps point strata to point strata, `regular' if it maps interval strata into interval strata, and `balanced' if it is both singular and regular.
\end{term}

\nid In analogy to our earlier definitions of `subtrusses', `truss degeneracies' and `truss coarsenings' (see \autoref{term:1-truss-inj} and  \autoref{term:1-truss-surj}), we further introduce the following properties of mesh maps.

\begin{term}[Submeshes of 1-meshes] \label{term:1-submeshes} A map of 1-meshes $F : (M,f,\gamma) \to (N,g,\rho)$ is called a `submesh' if $F : (M,f) \to (N,g)$ is a substratification such that $\gamma$ is a restriction of $\rho$, that is, the induced map $\tilde F : \im(\gamma) \to \im(\rho)$ (see \autoref{defn:maps-of-1-meshes}) is a subspace inclusion (of subspaces in $\lR$).
\end{term}

\nid Note that any submesh inclusion is necessarily balanced.

\begin{term}[Mesh degeneracies and coarsenings of 1-meshes] \label{term:1-mesh-deg-crs} A map of 1-meshes $F : (M,f,\gamma) \to (N,g,\rho)$ is called a `mesh degeneracy' if it is a surjective singular mesh map that maps interval strata either homeomorphically onto their image stratum, or `degenerates' them into a point stratum. Dually, $F$ is called a `mesh coarsening' if it is a mesh map that descents to a coarsening of stratifications $F : (M,f) \to (N,g)$. (Note that, by definition of stratified coarsenings, see \autoref{defn:coarsenings-and-refinements}, mesh coarsenings are necessarily surjective regular mesh maps, and homeomorphisms on underlying spaces.)
\end{term}

\begin{eg}[Maps of 1-meshes] In \autoref{fig:1-meshes-and-their-maps} we depict singular, regular, balanced and `mixed' 1-mesh maps (whose mappings we indicate by green arrows). Note that the first example is in fact a mesh degeneracy, the second a mesh coarsening, and the third a submesh.
\begin{figure}[ht]
    \centering
    \def\svgwidth{1\columnwidth}
    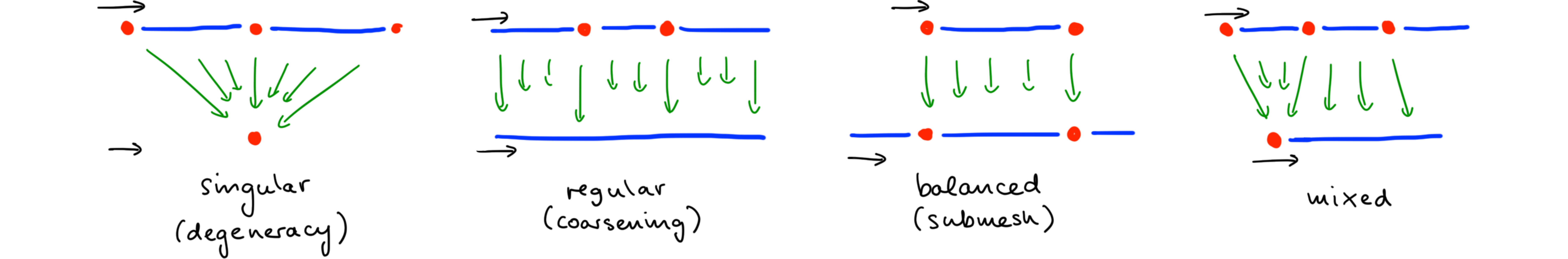

    \caption{1-Mesh maps.}
    \label{fig:1-meshes-and-their-maps}
\end{figure}
\end{eg}

\subsection{1-Mesh bundles}
\label{ssec:1-mesh-bundles}

\subsubsecunnum{The definition of constructible 1-mesh bundles}

We now introduce 1-mesh bundles and the condition for 1-mesh bundles to be constructible. Recall the notion of stratified bundles (see \autoref{defn:stratified-bundle}) and our earlier notation for the `bounds' of a flat 1-framing (see \autoref{notn:framing-closures}). We assume all stratifications to be finite (i.e.\ to have finitely many strata).

\begin{defn}[1-Mesh bundles] \label{defn:1-mesh-bundle} A \textbf{1-mesh bundle} $(p,\gamma)$ is a stratified bundle $p : (M,f) \to (B,g)$ together with a bundle embedding $\gamma : M \into B \times \lR$ into the trivial bundle $\pi : B \times \lR \to B$ which restricts on each fiber over $b \in B$ to a 1-mesh $(p\inv(b),f,\gamma)$ such that the maps $b \to \gamma^\pm p\inv(b)$ are continuous.
\end{defn}

\nid We will refer to the embedding $\gamma : M \into B \times \lR$ as the `fiber 1-framing' of the 1-mesh bundle $(p,\gamma)$. To lighten notation, we usually keep fiber 1-framings of 1-mesh bundles implicit, simply writing $p$ in place of $(p,\gamma)$.

\begin{term}[Fiber 1-framing bounds] Given a 1-mesh bundle $p :(M,f) \to (B,g)$ with fiber 1-framing $\gamma$, the maps $b \to \gamma^\pm p\inv(b)$ assemble into continuous sections of the projection $\pi : B \times \lR \to B$, which we denote by $\gamma^\pm : B \to B \times \lR$, and refer to as the `fiber 1-framing bounds'.
\end{term}

\begin{eg}[A 1-mesh bundle] \label{eg:1-mesh-bundle} In \autoref{fig:1-mesh-bundle-over-non-cellular-base} we depict a 1-mesh bundle whose base stratification is the classifying stratification $(B,g) = \CStr{P}$ of the poset $P = (a \ot b \ot c \to d)$ (this poset also appeared as the base poset in \autoref{fig:1-truss-bundle}). Note that fibers of the stratum $\cstratum(c) \subset \CStr {P}$ are open 1-meshes, fibers of the stratum $\cstratum(a)$ and $\cstratum(d)$ are closed 1-meshes, while the fiber over $\cstratum(b)$ is neither open nor closed. Note that the bundle map is depicted as a (restriction of the) projection $B \times \lR \to B$ (with the orientation of $\lR$ indicated by a green arrow) which provides the fiber 1-framing of the bundle.
\begin{figure}[ht]
    \centering
    \def\svgwidth{1\columnwidth}
    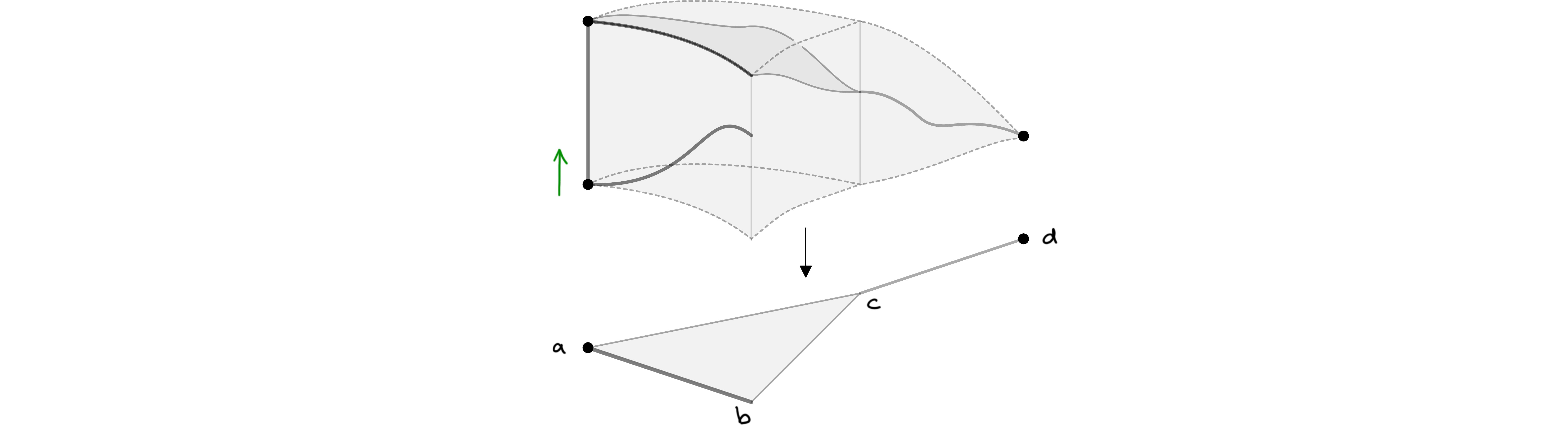

    \caption{A 1-mesh bundle.}
    \label{fig:1-mesh-bundle-over-non-cellular-base}
\end{figure}
\end{eg}

\begin{rmk}[Omitting orientations] When depicting bundle maps of 1-mesh bundles as projections $B \times \lR \to B$ we will (analogously to our depictions of 1-truss bundles, see \autoref{rmk:z2-ambiguity-1-truss-bundles}) sometimes omit indicating the standard orientation of $\lR$, which leaves a $\lZ_2$-ambiguity for such a choice.
\end{rmk}

\begin{term}[Open and closed 1-mesh bundles] \label{term:open-closed-1-mesh-bundles} A 1-mesh bundle $p$ is called open resp.\ closed if all its fibers are open resp.\ closed 1-meshes.
\end{term}

\begin{defn}[Maps of 1-mesh bundles] \label{defn:1-mesh-bundle-maps} Given 1-mesh bundles $p : (M,f) \to (B,g)$ and $p' : (M',f') \to (B',g')$, a \textbf{1-mesh bundle map} $(F,G) : p \to p'$ consists of stratified maps $F : (M,f) \to (M',f')$ and $G : (B,g) \to (B',g')$ commuting with $p$ and $p'$, such that $F$ restricts on each fiber to a 1-mesh map.
\end{defn}

\begin{rmk}[Bundle maps commute with fiber framings] \label{rmk:mesh-bundle-maps-commute-w-framings} Every 1-mesh bundle map $(F,G) : p \to p'$ induces a commutative diagram of continuous maps (where $\gamma$ resp.\ $\gamma'$ are fiber 1-framings of $p$ resp.\ of $p'$)
\begin{equation}
\begin{tikzcd}[baseline=(W.base)]
\im(\gamma) \arrow[d, "\tilde F"'] \arrow[rr, "\pi", bend left=12] & M \arrow[l, "\gamma"] \arrow[r, "p"'] \arrow[d, "F"'] & B \arrow[d, "G"] \\
\im(\gamma') \arrow[rr, "\pi"', bend right=12] & M' \arrow[l, "\gamma'"'] \arrow[r, "p'"] & |[alias=W]| B'
\end{tikzcd} . \qedhere
\end{equation}
\end{rmk}

\begin{term}[Singular, regular and balanced 1-mesh bundles maps] We call a 1-mesh bundle map $(F,G) : p \to p'$ singular resp.\ regular resp.\ balanced if it is fiberwise singular resp.\ regular resp.\ balanced in the sense of \autoref{term:1-mesh-map-reg-sing}.
\end{term}

\pauseae

We next introduce an additional condition on 1-mesh bundles that will guarantee their `constructibility'. The following notation and terminology will be helpful.

\begin{term}[Restrictions over base strata] Let $p : (M,f) \to (B,g)$ be a 1-mesh bundle. Over each stratum $s$ in the base $(B,g)$, $p$ restricts to a 1-mesh bundle denoted by $\rest p s : (p\inv(s), f) \to s$.
\end{term}

\begin{obs}[Trivialization over base strata] \label{obs:trivializing-1-mesh-bundles} Since automorphism spaces of 1-meshes are contractible, the restricted 1-mesh bundle $\rest p s$ is (non-uniquely) isomorphic to a trivial 1-mesh bundle which we will denote as follows
\begin{equation}
\begin{tikzcd}
    {(p\inv(s), f)} \arrow[rr, "\iso"] \arrow[rd, "\rest p s"'] &   & s \times \fib {s} \arrow[ld] \\
 & s &
\end{tikzcd} .
\end{equation}
Here, $\fib {s}$ is a 1-mesh and the trivialization restricts to a 1-mesh isomorphism on all fibers of $\rest p s$.
\end{obs}

\begin{term}[Regular and singular strata]  Let $p : (M,f) \to (B,g)$ be a 1-mesh bundle and $r$ a stratum in $f$ lying over a stratum $s$ in $g$. Up to trivializing the restricted bundle $\rest p s$, the stratum $r$ is of the form $r \equiv s \times u$ where $u$ is a stratum $\fib s$. We call $r$ `singular' if $u$ is a point stratum and `regular' if $u$ is an interval stratum. 
\end{term}

\nid Recall, a (`formal') entrance path $s \to r$ between strata $s$ and $r$ in a stratification $(X,f)$ is given if the intersection $\overline s \cap r$ of the closure of $s$ in $f$ and $r$ is non-empty (see \autoref{defn:fep}). Denote by $f^{\geq s}$ the `entrance path closure' of $s$ in $(X,f)$, given by the union of all strata $r$ for which there is a sequence of formal entrance paths $s = r_0 \to r_1 \to ... \to r_k = r$ (in other words, $f^{\geq s}$ is the preimage of the upper closure $\Entr(f)^{\geq s}$ under the characteristic map $f : X \to \Entr(f)$). The condition for constructibility of 1-mesh bundles is that `formal entrance paths lift uniquely to singular strata' in the following sense.

\begin{defn}[Constructible 1-mesh bundles] \label{defn:1-mesh-bundles} A 1-mesh bundle $p : (M,f) \to (B,g)$ is \textbf{constructible} if for any singular stratum $s$ of $(M,f)$ lying over a stratum $b = p(s)$ in the base $(B,g)$, the homeomorphism $p : s \to b$ extends to a homeomorphisms of entrance path closures $p : f^{\geq s} \to g^{\geq b}$.
\end{defn}

\begin{eg}[Constructible and non-constructible 1-mesh bundles] In \autoref{fig:constructible-and-non-constructible-1-mesh-bundles} we illustrate both a constructible 1-mesh bundle (on the left) as well as a non-constructible 1-mesh bundle (on the right). In each case, singular strata are colored in red and regular strata in blue. Note that the 1-mesh bundle on the right fails the constructibility condition on \emph{both} its singular 1-strata: firstly, for the upper singular 1-stratum the entrance path closure contains the regular 1-stratum as well as the singular 0-stratum, and, secondly, for the lower singular 1-stratum the entrance path closure is empty.\footnote{The two illustrated failures of the constructibility condition are somewhat different in nature, and weaker versions of the notion of `constructible bundle' (for which, in particular, fibers may be empty) may, in fact, allow the second type of failure to happen.}
\begin{figure}[ht]
    \centering
    \def\svgwidth{1\columnwidth}
    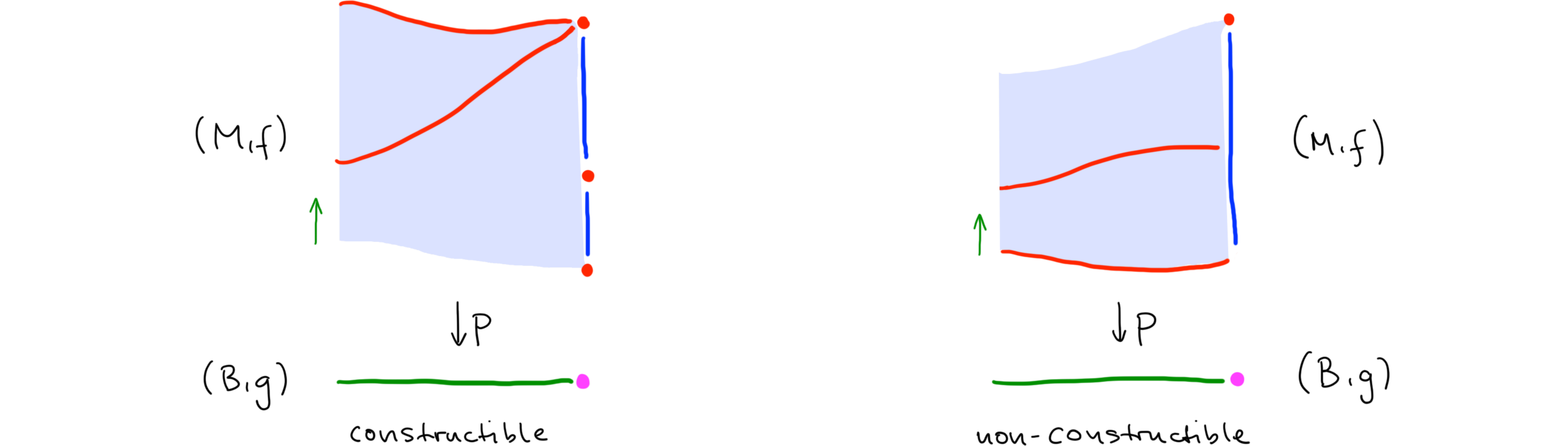

    \caption{Constructible and non-constructible 1-mesh bundles.}
    \label{fig:constructible-and-non-constructible-1-mesh-bundles}
\end{figure}
\end{eg}

\begin{rmk}[Constructibility in the case of conical bundles] For a conical stratification $(X,f)$ the entrance path closure $f^{\geq s}$ of a stratum $s$ equals the ordinary closure $\overline s$ of $s$ in $X$. Therefore, if both its base and total space stratification are conical, we may replace entrance path closures by topological closures of strata in the definition of constructible 1-mesh bundle. Importantly, note that the constructibility condition is still non-trivial even in the `better behaved' conical case: in \autoref{fig:a-conical-non-constructible-1-mesh-bundle} we depict a 1-mesh bundle which, despite its base and total stratification being conical, is non-constructible.
\begin{figure}[ht]
    \centering
    \def\svgwidth{1\columnwidth}
    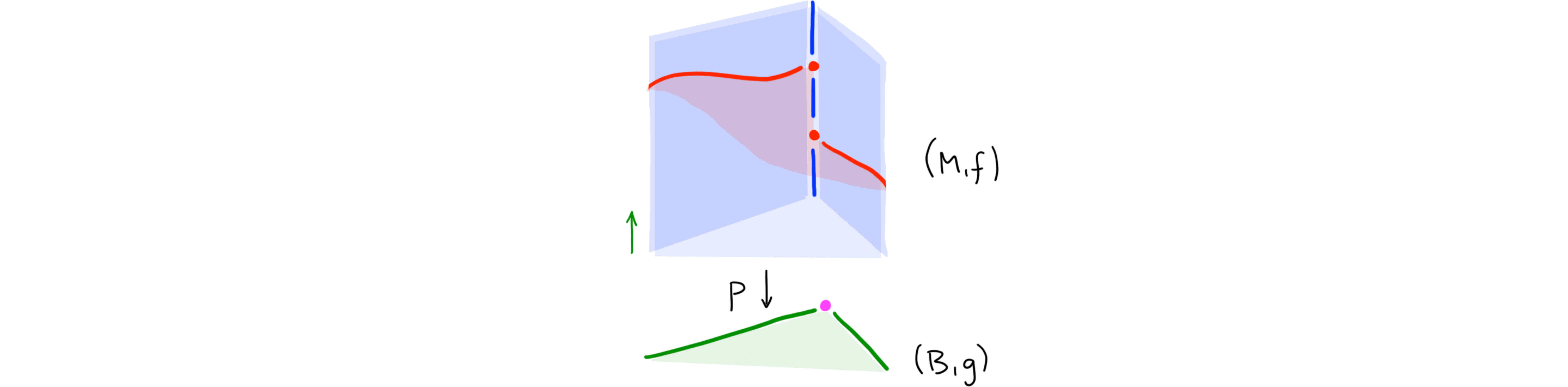

    \caption{A conical but non-constructible 1-mesh bundle.}
    \label{fig:a-conical-non-constructible-1-mesh-bundle}
\end{figure}
\end{rmk}

\begin{term}[Boundary-constructible stratifications] \label{term:boundary-constructibility} A stratification $(X,f)$ is called `boundary-constructible' if the entrance path closure $f^{\geq s}$ of each stratum $s$ equals the topological closure $\overline s$ of $s$ in $X$. (Boundary-constructibility may be phrased in several ways; for instance, $(X,f)$ is boundary constructible if and only if its characteristic maps is an open map, see \autoref{lem:boundary-constructibility-characteric-map}.)
\end{term}

\nid Importantly, the property of `boundary-constructibility' lifts in constructible 1-mesh bundles as follows.

\begin{obs}[Boundary-constructibility lifts along constructible bundles] \label{obs:boundary-constructibility-lifts} Let  $p : (M,f) \to (B,g)$ be a constructible 1-mesh bundle. If $(B,g)$ is boundary-constructible, it follows that $(M,f)$ itself must be boundary-constructible.
\end{obs}

Finally, let us briefly relate the condition for constructibility given above to the general idea of constructibility outlined earlier.

\begin{preview}[Constructibility] The fact that our definition of constructible 1-mesh bundles guarantees constructibility of 1-mesh bundles in the ordinary sense (namely, in that constructible bundles can be `classified up to homotopy equivalence by functor on the fundamental category of the base stratification') will be a corollary of the equivalence of meshes and trusses: since truss bundles are `constructible' (in that they are classified by functors on their base posets, as we have seen in \autoref{ch:trusses}) the equivalence will in particular entail constructibility of 1-mesh bundles.
\end{preview}

\subsubsecunnum{Pullbacks and compactifications of 1-mesh bundles} We discuss two important constructions on 1-mesh bundles: firstly, 1-mesh bundles can be fiberwise compactified and secondly, they can be pulled back along stratified maps in the base.

\begin{constr}[Fiberwise compactifications of 1-mesh bundle] \label{constr:fiberwise_compactification} Given a 1-mesh bundle $p : (M,f) \to (B,g)$ with fiber 1-framing $\gamma : M \into B \times \lR$, we construct a 1-mesh bundle $\overline p : (\overline M,\overline f) \to (B,g)$ called the `fiberwise compactification' of $p$. Define $\overline M$ to be the closure of $\gamma(M)$ in $B \times \lR$, and obtain $\overline p : \overline M \to B$ by restricting $\pi : B \times \lR \to B$ to $\overline M$. The stratification $(\overline M,\overline f)$ is defined to have strata to be images $\gamma(r)$ of strata $r$ in $f$, or images $\gamma^\pm(s)$ of strata $s$ in $g$.
\end{constr}

\begin{obs}[Constructibility and compactifications] \label{obs:constructibility-in-compactifications} Let $p : (M,f) \to (B,g)$ be a 1-mesh bundle, and $\overline p :  (\overline M,\overline f) \to (B,g)$ its compactification.
    \begin{enumerate}
        \item If $p$ is constructible, then so is $\overline p$.
        \item If $\overline p$ is constructible, and the inclusion $\sing(\Entr f) \into \sing(\Entr \overline f)$ preserves upper closures (where $\sing(\Entr f)$ is the full subposet of $\Entr f$ consisting of singular strata, and similarly for $\overline f$), then $p$ is itself constructible. \qedhere
    \end{enumerate}
\end{obs}

\begin{eg}[1-Mesh bundle compactification] In \autoref{fig:the-compactification-of-a-constructible-1-mesh-bundle} we depict the compactification of a 1-mesh bundle $p$. Note that, since $p$ is constructible so is its compactification $\overline p$.
\begin{figure}[ht]
    \centering
    \def\svgwidth{1\columnwidth}
    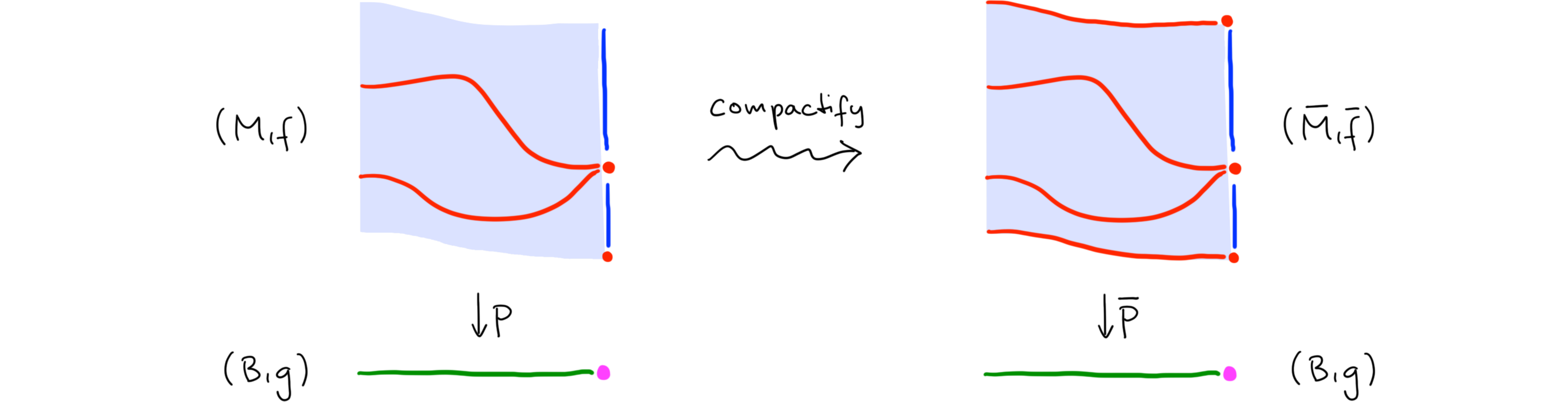

    \caption{The compactification of a constructible 1-mesh bundle.}
    \label{fig:the-compactification-of-a-constructible-1-mesh-bundle}
\end{figure}
\end{eg}

\begin{constr}[Pullbacks of 1-mesh bundles] \label{constr:pullback_of_constr_diagrams_bundles} Given a 1-mesh bundle $p : (M,f) \to (B,g)$ and a stratified map $G : (B',g') \to (B,g)$ we construct a 1-mesh bundle $G^*p : (G^*M, G^*f) \to (B',g')$ called the `pullback of $p$ along $G$'. As a map of stratified spaces, construct $G^*p$ by a pullbacks (see \autoref{term:strat-pullback}) as follows.
\begin{equation}
\begin{tikzcd}
    (G^*M, G^*f)
 \arrow[dr, phantom, "\lrcorner" , very near start, color=black]
 \arrow[r, "\Tot{}G"] \arrow[d, "G^*p"'] & (M,f) \arrow[d, "
 p"] \\
(B',g') \arrow[r, "G"]                      & (B,g)
\end{tikzcd}
\end{equation}
Note that $G^*p$ is a stratified bundle since $p$ is. Construct a fiber 1-framing $G^*\gamma : G^*M \into B' \times \lR$ for $G^*p$ by defining, for $x \in G^*M$, $G^*\gamma(x) \in B' \times \lR$ to have first component $ G^*p (x)$ and second component equal to the second component of $\gamma \circ \Tot{}G (x)$.
\end{constr}

\begin{obs}[Pullbacks preserve constructibility] Consider a 1-mesh bundle $p : (M,f) \to (B,g)$ and a stratified map $G : (B',g') \to (B,g)$. Then the pullback 1-mesh bundle $G^*p$ is constructible whenever $p$ is (this follows since stratified maps map closures of strata into closures of strata).
\end{obs}

\begin{obs}[Pullbacks preserve fiberwise compactifications] Consider a 1-mesh bundle $p : (M,f) \to (B,g)$ and a stratified map $G : (B',g') \to (B,g)$. Then the pullback $G^*{\overline p}$ of the fiberwise compactification of $p$ is the fiberwise compactification $\overline{G^*p}$ of the pullback $G^*p$.
\end{obs}

\subsubsecunnum{1-Mesh bundles over regular cells} We next discuss an important property of 1-mesh bundles which, roughly speaking, states that `strata in 1-mesh bundles over regular cells must themselves be regular cells'. We will work with the following slight generalization of the notion of regular cell complexes, which drops the requirement for the boundaries of each cell to be `complete' spheres, but instead requires them to be `completable'. Recall the notion of `constructible substratification' which is a substratification obtained as a preimage of a subposet in the entrance path poset of a stratification (see \autoref{defn:constr_substrat}).

\begin{defn}[Regular cell stratifications] \label{defn:regular-cell-strat} A stratification $(B,f)$ is called an \textbf{(open) regular cell stratification} if is an (open) constructible substratification of some regular cell complex.
\end{defn}

\nid Note that regular cell stratifications are necessarily conical (see \autoref{rmk:regular-cell-implies-conicality}).

\begin{prop}[Open regular cell 1-mesh bundles] \label{prop:constr-implies-cellular} Let $p : (M,f) \to (B,g)$ be a constructible closed (or open) 1-mesh bundle. If the base $(B,g)$ is open regular cell then $(M,f)$ is open regular cell.
\end{prop}

\nid If both base and total stratification of a 1-mesh bundle bundle $p$ are (open) regular cell, then we call $p$ itself `(open) regular cell'.

\begin{rmk}[General regular cell 1-mesh bundles] \label{rmk:regular-cell-for-general-1-mesh-bun} Working with `regular cell' stratifications in place of `open regular cells' stratifications, the preceding statement still holds. In fact, equally applies to general 1-mesh bundles. The additional assumption on stratifications to be \emph{open} regular cell, however, simplifies the proof and covers all cases of interest to us.
\end{rmk}

\nid The proof of \autoref{prop:constr-implies-cellular} will use the following general observation about cell of `quotient-cylindrical shape'.

\begin{lem}[Quotient-cylindrical cells] \label{obs:closed_int_bund_w_deg_over_disks} Let $D^m$ be the closed $m$-disk, $S^{n-1}$ its boundary, and let $p : X \to D^m$ be a subbundle of the projection $\pi : D^m \times \lR \to D^m$ whose fibers over $x \in D^m$ are subsets of $\lR$ of the form $[\gamma^-(x),\gamma^+(x)]$, where $\gamma^\pm : D^m \to D^m \times \lR$ are continuous sections. If $\gamma^-(x) < \gamma^+(x)$ everywhere, except possibly when $x \in S^{n-1}$, then $X$ is a closed $(m+1)$-disk.
\end{lem}

\begin{proof} Define the `trivial fiber locus' of $p$ to be the closed subset $\mathsf{triv}(p) = \{x ~|~ \gamma^-(x) = \gamma^+(x)\}$ of $D^m$. If $\mathsf{triv}(p)  = S^{m-1}$, then $X$ is the quotient of $D^m \times [-1,1]$ by $S^{m-1} \times [-1,1]$ which is an $m$-disk. If $\mathsf{triv}(p) \subsetneq S^{m-1}$ is a proper subset, the claim holds by the following standard argument. Define $V \subset S^{m-1} \times \lR$ to be the subspace with points $(x,y)$ where $x \in S^{m-1}$ and $y \in [-\mathrm{dist}(x,\mathsf{triv}(p)),+\mathrm{dist}(x,\mathsf{triv}(p))] \subset \lR$ (where $\mathrm{dist}$ is some metric on $S^{m-1}$). Let $W$ be the convex hull of $V$ in $\lR^{m+1}$. This projects to the convex hull of $S^{m-1}$ in $\lR^{m}$, which is $D^m$. By construction, the projection $q : W \to D^m$ is bundle isomorphic to $p : X \to D^m$, for instance, by a linearly identifying fibers $q\inv(z) \iso p\inv(z)$ for all $z \in D^m$. Since compact convex sets with non-empty interior are disks, we have $W \iso D^{m+1}$, and thus $X \iso D^{m+1}$ as claimed.
\end{proof}

\begin{proof}[Proof of \autoref{prop:constr-implies-cellular}] We first argue in the case of constructible \emph{closed} 1-mesh bundles $p : (M,f) \to (B,g)$ with open regular cell base. By definition, $(B,g)$ is a open constructible substratification of a regular cell complex $(Y,e)$ (where $e$ stratifies $Y$ by its cells). By removing cells from $Y$ that are not in the closure of $B$, we may further assume $\overline B = Y$. Let $\mathsf{triv}(p) \subset B$ be the closed subset of $B$ over which $p$ has point fibers, and denote by $\mathsf{triv}(p)^*$ the union of $\mathsf{triv}(p)$ and $Y \setminus B$. We build a bundle $q : (X,d) \to (Y,e)$ with $(X,d)$ being a regular cell complex containing $(M,f)$ as an open constructible substratification. The underlying map $q : X \to Y$ is the subbundle of the projection $\pi : Y \times \lR \to Y$ with fibers over $x$ given by $[-\mathrm{dist}(x,\mathsf{triv}(p)^*),+\mathrm{dist}(x,\mathsf{triv}(p)^*)]$.\footnote{Note every regular cell complex $Y$ is metrizable, for instance by using the `barycentric metric' of its barycentric subdivision. The same need not hold for non-regular cell complexes!} Denote by $\rest q B$ the restriction of $q$ to $B \subset Y$ in the base. Observe that we can identify bundles $\rest q B \iso p$ by linearly identifying their fibers $q\inv(z) \iso p\inv(z)$, $z \in B$ (note $p\inv(z)$ inherits linear structure from the fiber 1-framing $\gamma$ of $p$, which embeds $p\inv(z) \into \lR$). In particular, this identifies $q\inv(B) \iso M$. Now, stratify $X$ by defining strata $s$ of $(X,d)$ to be either strata $s$ in $(M,f)$ or of the form $r \times \{0\}$ where $r$ is a stratum of $(Y \setminus B,e)$. Verify that, by constructibility of $p$, this makes $q$ itself a constructible closed 1-mesh bundle. Using \autoref{obs:closed_int_bund_w_deg_over_disks} as well as the constructibility of $q$, one further verifies that $(X,d)$ is a regular cell complex. The constructible substratification $(M,f) \into (X,d)$ then exhibits $(M,f)$ as a open regular cell stratification. We illustrate the construction in \autoref{fig:closed-1-mesh-bundles-over-open-regular-cell-base-are-open-regular-cell}.
\begin{figure}[ht]
    \centering
    \def\svgwidth{1\columnwidth}
    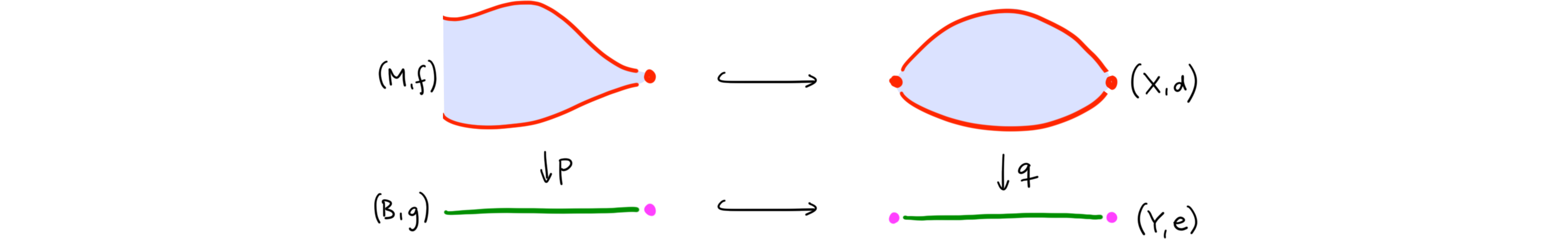

    \caption[Closed mesh bundle over open regular cell base.]{Closed 1-mesh bundles over open regular cell base have open regular cell total stratification.}
    \label{fig:closed-1-mesh-bundles-over-open-regular-cell-base-are-open-regular-cell}
\end{figure}

    Finally, to see the statement in the case of open 1-mesh bundles $p : (M,f) \to (B,g)$ over open regular cell stratifications $(B,g)$, one first compactifies $p$ to a closed 1-mesh bundle $\overline p : (\overline M, \overline f) \to (B,g)$ using \autoref{constr:fiberwise_compactification}. The previous argument shows that $(\overline M, \overline f)$ is open regular cell; since $(M,f)$ is an open constructible substratification of $(\overline M, \overline f)$ it follows that $(M,f)$ is open regular cell as claimed.
\end{proof}

\subsection{$n$-Meshes}
\label{ssec:n-meshes}

\subsubsecunnum{$n$-Meshes as constructible towers}

We now generalize the definition of 1-meshes to $n$-meshes.

\begin{defn}[$n$-Meshes] \label{defn:n-meshes} An \textbf{$n$-mesh} $M$ is a tower of constructible 1-mesh bundles $(M_n, f_n) \xto {p_n} (M_{n-1}, f_{n-1}) \xto {p_{n-2}} ... \xto {p_2} (M_1,f_1) \xto {p_1} (M_0,f_0) = \mathrm{pt}$.
\end{defn}

\nid Just as 1-meshes have flat 1-framings (embedding a 1-mesh in $\lR$), $n$-meshes have `flat $n$-framings' (embedding an $n$-mesh in $\lR^n$). This can be constructed from the individual fiber 1-framings of the 1-mesh bundles in an $n$-mesh as follows.  Recall the standard euclidean $n$-proframe $\Pi = (\pi_n,\pi_{n-1},...,\pi_1)$ consisting of projections $\pi_i : \lR^i \to \lR^{i-1}$ (see \autoref{defn:standard-proframe}).

\begin{constr}[Flat $n$-framings of $n$-meshes] \label{rmk:n-framings-of-n-meshes} Given an $n$-mesh $M$ consisting of 1-mesh bundles $p_i : (M_i,f_i) \to (M_{i-1},f_{i-1})$, recall that each $p_i$ comes equipped with a fiber 1-framing $M_i \into M_{i-1} \times \lR$. This allows us to inductively construct a map of towers of spaces
    \begin{equation}
        \begin{tikzcd}
            M_n \arrow[r, "p_n"] \arrow[d, "\gamma_n"'] & M_{n-1} \arrow[r, "p_{n-1}"] \arrow[d, "\gamma_{n-1}"'] & \cdots \arrow[r, "p_2"] \arrow[d, "\cdots", phantom] & M_1 \arrow[r, "p_1"] \arrow[d, "\gamma_1"] & M_0 = \mathrm{pt} \arrow[d] \\
\lR^n \arrow[r, "\pi_n"'] & \lR^{n-1} \arrow[r, "\pi_{n-1}"'] & \cdots \arrow[r, "\pi_2"'] & \lR^1 \arrow[r, "\pi_1"'] & \lR^0
        \end{tikzcd}
    \end{equation}
    where $\gamma_i$ is obtained by postcomposing the fiber 1-framing $M_i \into M_{i-1} \times \lR$ with the product map $\gamma_{i-1} \times \lR : M_{i-1} \times \lR \into \lR^{i-1} \times \lR$. We refer to the embedding of towers $\gamma : M \into \Pi$ as the \textbf{flat $n$-framing} of the $n$-mesh $M$. Note that $\gamma$ is fully determined by its top component $\gamma_n : M_n \into \lR^n$ and, abusing terminology, we sometimes refer the embedding $\gamma_n$ itself as the `flat $n$-framing' of $M$.
\end{constr}

\nid Given an $n$-mesh $M$ with flat $n$-framing $\gamma$, note that components of $\gamma$ may be considered either as subspace embeddings $\gamma_i : M_i \into \lR^i$ or as stratified maps $\gamma : (M_i,f_i) \to \lR^i$ (which are then coarsenings onto their images in $\lR^i$).

\begin{eg}[A mixed 2-mesh]  In \autoref{fig:a-mixed-2-mesh} we depict a 2-mesh via its flat $n$-framing $\gamma : M \into \Pi$ into the standard euclidean tower $\Pi$ (we indicate this by providing coordinate axis for the ambient euclidean space $\lR^i$). Note that the 2-mesh is `mixed' in that it is neither open nor closed.
\begin{figure}[ht]
    \centering
    \def\svgwidth{1\columnwidth}
    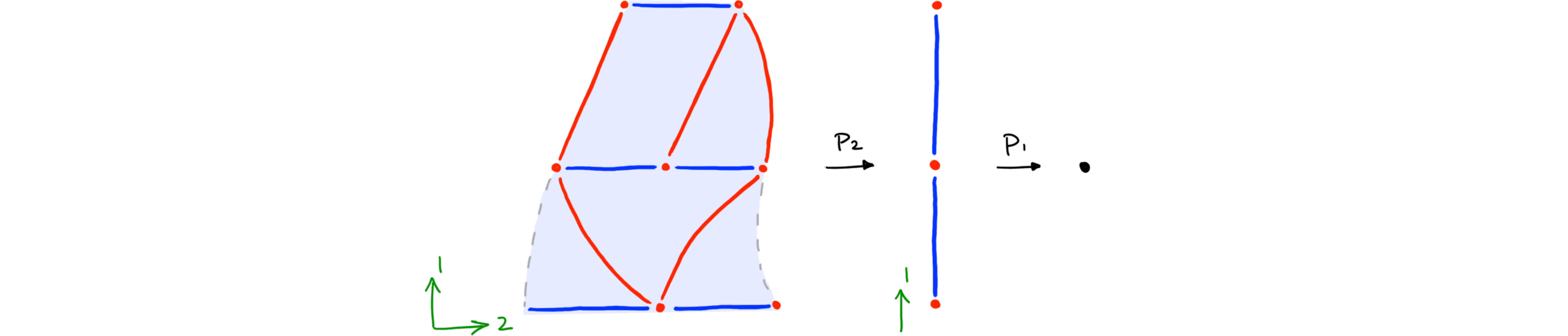

    \caption{A mixed 2-mesh.}
    \label{fig:a-mixed-2-mesh}
\end{figure}
\end{eg}

\begin{eg}[Closed and open $n$-Meshes] \label{eg:n-diag-and-codiag} In \autoref{fig:a-closed-and-open-meshes-in-dimension-2-and-3} we depict closed and open meshes in dimension 2 and 3. In each case we again depict the mesh via its flat $n$-framing $\gamma : M \into \Pi$.
\begin{figure}[ht]
    \centering
    \def\svgwidth{1\columnwidth}
    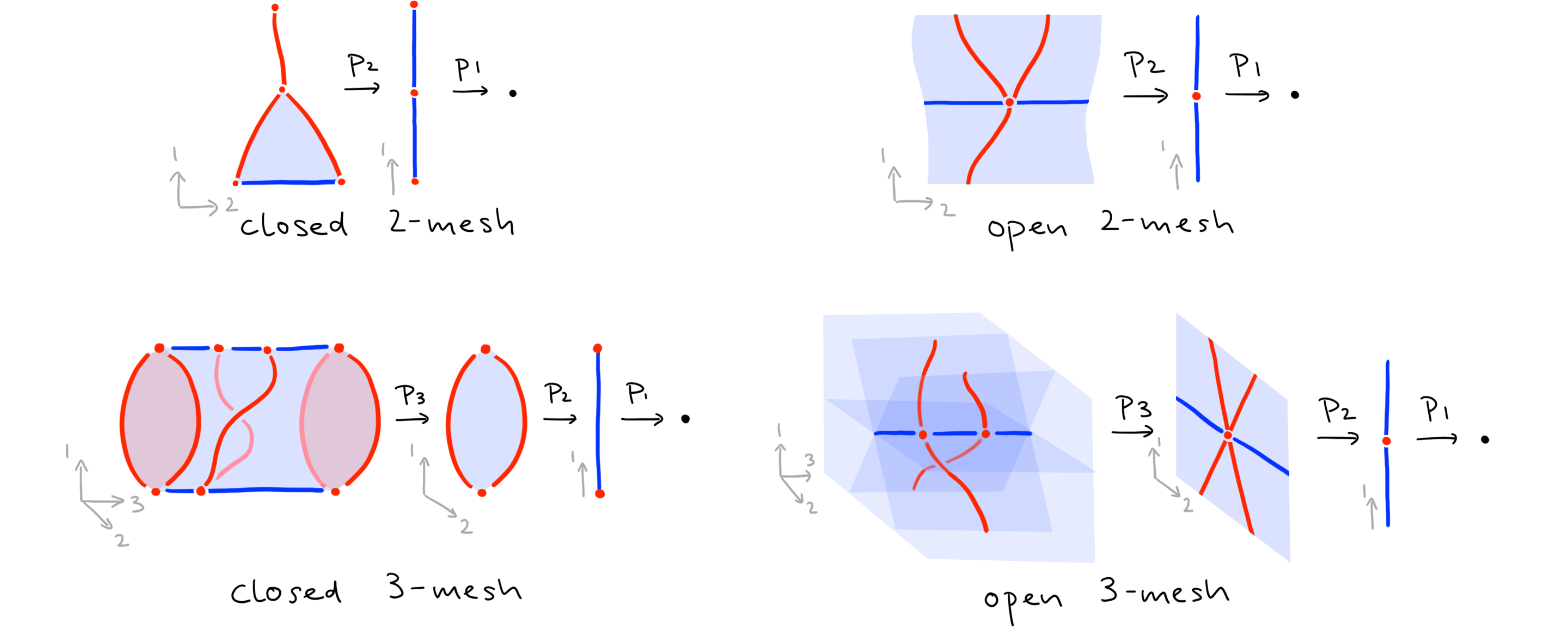

    \caption{A closed and open meshes in dimension 2 and 3.}
    \label{fig:a-closed-and-open-meshes-in-dimension-2-and-3}
\end{figure}
\end{eg}

\begin{term}[Closed or open $n$-meshes] A mesh is called `closed' (resp.\ `open') if each 1-mesh bundle in its tower is closed (resp.\ open).
\end{term}

\begin{obs}[$n$-Meshes are boundary-constructible] \label{rmk:meshes-are-boundary-constructible} Our earlier \autoref{obs:boundary-constructibility-lifts}, that `boundary-constructibility lifts in constructible $n$-mesh bundles', implies that for $n$-meshes $M$ each stratification $(M_i,f_i)$ is boundary-constructible.
\end{obs}

\begin{obs}[Closed and open $n$-meshes are regular cell] \label{rmk:closed-open-meshes-are-open-regular-cell} Our earlier \autoref{prop:constr-implies-cellular}, that `open regular cellularity lifts in closed resp.\ open constructible 1-mesh bundles', implies that for closed resp.\ open $n$-meshes $M$ each $(M_i,f_i)$ is an open regular cell stratification, and thus, in particular, conical.
\end{obs}

\begin{rmk}[General $n$-meshes are also regular cell] \label{rmk:general-meshes-are-open-regular-cell} The previous observation generalizes: general $n$-meshes are regular cell stratifications (while not necessarily open ones) and thus, in particular, conical. We omit a proof of this more general statement, since we will be mainly interested in the case of open and closed meshes. 
\end{rmk}

\subsubsecunnum{Maps of $n$-meshes}

\begin{defn}[Maps of $n$-meshes] \label{defn:maps_of_n_diagrams} Consider $n$-meshes $f$ and $g$. An \textbf{$n$-mesh map} $F : (M,f) \to (N,g)$ is a map between the towers $(M,f)$ and $(N,g)$ of the form
\begin{equation}
\begin{tikzcd}
    (M_n,f_n) \arrow[r, "p_n"] \arrow[d, "F_n"] & (M_{n-1},f_{n-1}) \arrow[r, "p_{n-1}"] \arrow[d, "F_{n-1}"] & \cdots \arrow[d, "\cdots", phantom] \arrow[r, "p_2"] & (M_1,f_1) \arrow[r, "p_1"] \arrow[d, "F_1"] & (M_0,f_0) \arrow[d, "F_0"] \\
    (N_n,g_n) \arrow[r, "q_n"] & (N_{n-1},g_{n-1}) \arrow[r, "q_{n-1}"] & \cdots \arrow[r, "q_2"] & (N_1,g_1) \arrow[r, "q_1"] & (N_0,g_0)
\end{tikzcd}
\end{equation}
such that each $(F_i,F_{i-1})$ is a 1-mesh bundle map $p_i \to q_i$.
\end{defn}

\begin{rmk}[$n$-Mesh maps commute with flat $n$-framings] \label{rmk:n-mesh-maps-commute-with-n-framing} In analogy to \autoref{rmk:mesh-bundle-maps-commute-w-framings}, given $n$-meshes $(M,f)$ and $(N,g)$ with flat $n$-framing $\gamma$ resp.\ $\rho$, an $n$-mesh map $F : (M,f) \to (N,g)$ map induces a map of the image towers $\tilde F : \im(\gamma) \to \im(\rho)$ (where both $\im(\gamma)$ and $\im(\rho)$ are subtowers of $\Pi$) such that $\tilde F \circ \gamma = \rho \circ F$.
\end{rmk}

\nid It will be useful to record terminology for the type of map of spaces that underlie the components of $n$-mesh maps as follows.

\begin{term}[Framed maps] \label{term:framed-maps} Given subspaces $Z \subset \lR^n$ and $W \subset \lR^n$, a map $F : Z \to W$ is called a `framed map' if for each $0 < i < n$, $F$ descends along the standard projection $\pi_{>i} = \pi_{i+1} \circ ... \pi_{n-1} \circ \pi_n : \lR^n \to \lR^i$ to a map $F_i : \pi_{>i} Z \to \pi_{> i} W$.
\end{term}

\nid In particular, given an $n$-mesh map $F : M \to N$, then each component $\tilde F_k$ of the corresponding map $\tilde F : \im(\gamma) \to \im(\rho)$ of flat $n$-framings (see \autoref{rmk:n-mesh-maps-commute-with-n-framing}) is a framed map in the preceding sense.

\begin{term}[Singular, regular and balanced maps] \label{term:reg-sing-balanced-n-mesh-map} We call an $n$-mesh map $F$ `singular' resp.\ `regular' resp.\ `balanced' if its is fiberwise resp.\ regular resp.\ balanced, that is, for all $1 \leq i \leq n$, the bundle maps $(F_i,F_{i-1})$ restrict on 1-mesh fibers to singular resp.\ regular resp.\ balanced 1-mesh maps in the sense of \autoref{term:1-mesh-map-reg-sing}.
\end{term}

\begin{term}[Submeshes] \label{term:submeshes} We call an $n$-mesh map $F$ a `submesh' inclusion of $n$-meshes if it is fiberwise a submesh inclusion of 1-meshes in the sense of \autoref{term:1-submeshes}.
\end{term}

\nid Note, in particular, any submesh inclusion induces a tower of substratifications.

\begin{term}[Mesh degeneracies and coarsenings] \label{term:mesh-crs-and-deg} An $n$-mesh map will be called a `mesh degeneracy' resp.\ a `mesh coarsenings' if it is fiberwise a 1-mesh degeneracy resp.\ a 1-mesh coarsening in the sense of \autoref{term:1-mesh-deg-crs}.
\end{term}

\nid Note, in particular, any mesh coarsening induces a tower of homeomorphisms on underlying spaces, and is thus a tower of coarsenings of stratifications in the usual sense (see \autoref{defn:coarsenings-and-refinements}).

\begin{rmk}[Mesh refinements] Similar to \autoref{rmk:truss-crs-vs-ref}, we synonymously refer to `mesh coarsenings' as `mesh refinements' (of their codomain).
\end{rmk}

\begin{notn}[Ordinary categories of $n$-meshes] The category of $n$-meshes and their maps will be denoted by $\mesh n$.
    \begin{enumerate}
        \item The full subcategories of closed resp.\ open $n$-meshes are denoted by $\cmesh n$ resp.\ $\omesh n$. Superscripts `$\mathsf{s}$' resp.\ `$\mathsf{r}$' resp.\ `$\mathsf{rs}$' will indicate that we restrict these categories to singular resp.\ regular resp.\ balanced maps as morphisms.
        \item The wide subcategories containing degeneracies resp.\ coarsenings will be denoted by $\degmesh n$ and $\crsmesh n$. \qedhere
    \end{enumerate}
\end{notn}

\nid The category of meshes is naturally $\Top$-enriched as follows.

\begin{notn}[Topological categories of $n$-meshes] The topologically enriched category of $n$-meshes $\tmesh n$ is obtained by topologizing hom sets $\mesh n (M,N)$ as subspaces of the product
\begin{equation}
\Map(M_n,N_n) \times \Map(M_{n-1},N_{n-1}) \times ... \times \Map(M_0,N_0) \quad .
\end{equation}
This further restricts to subcategories yielding the following topological categories (where hom subsets are given subspace topology).
\begin{enumerate}
    \item One obtains the topological categories closed meshes $\ctmesh n$ resp.\ open meshes $\otmesh n$ (superscripts `$\mathsf{s}$', `$\mathsf{r}$' and `$\mathsf{rs}$' apply as before).
    \item One similarly obtains the topological categories $\degtmesh n$ of meshes and degeneracies, and $\crstmesh n$ of meshes and coarsenings. \qedhere
\end{enumerate}
\end{notn}

\nid Like $n$-trusses, $n$-meshes `truncate' and truncation is continuous.

\begin{term}[Truncations]\label{constr:top-diag-truncations} Given an $n$-mesh $M = ((M_n, f_n) \xto {p_n} ... \xto {p_1} (M_0,f_0))$ its \textbf{$i$-truncation} $M_{\leq i}$ is the $i$-mesh $((M_i, f_i) \xto {p_n} ... \xto {p_1} (M_0,f_0))$ obtained by truncating the tower $(M,f)$ to its lower $i$ degrees.
\end{term}

\begin{rmk}[Truncating is continuous] \label{rmk:truncating-topological} $i$-Truncation induces a $\Top$-enriched functor $(-)_{\leq i} : \tmesh n \to \tmesh {i}$ of topological categories of meshes.
\end{rmk}

\begin{constr}[Inductive pullbacks of $n$-meshes] \label{constr:pullback-of-diagrams} Consider an $n$-mesh $M = ((M_n, f_n) \xto {p_n} ... \xto {p_1} (M_0,f_0))$, an $(n-1)$-mesh $N = ((N_{n-1}, g_{n-1}) \xto {q_n} ... \xto {q_1} (N_0,N_0))$ as well as an $(n-1)$-mesh map $G : N \to M_{\leq n-1}$. Using \autoref{constr:pullback_of_constr_diagrams_bundles} we can pull back $p_n$ along $G_{n-1}$ to obtain the constructible 1-mesh bundle bundle $q_n : (G_{n-1}^*M_n, G_{n-1}^*f_n) \to (N_{n-1},g_{n-1})$ as follows:
\begin{equation}
\begin{tikzcd}[baseline=(W.base)]
    (G_{n-1}^*M_n, G_{n-1}^*f_n)
 \arrow[dr, phantom, "\lrcorner" , very near start, color=black]
 \arrow[r, "\Tot{} G_{n-1}"] \arrow[d, "q_n"'] & (M_n,f_n) \arrow[d, "p_n"] \\
 (N_{n-1},g_{n-1}) \arrow[r, "G_{n-1}"]                      & |[alias=W]| (M_{n-1},f_{n-1})
\end{tikzcd} .
\end{equation}
Extending the tower $N$ by the constructible 1-mesh bundle $q_n$ defines the $n$-mesh $G^*M$, called the \textbf{inductive pullback of $M$ along $G$}. This comes together with a canonical $n$-mesh map $G^*M \to M$ obtained by extending $G$ with the map $\Tot{}G_{n-1}$.
\end{constr}

\begin{rmk}[Inductive pullbacks can be used inductively] Consider an $n$-mesh $M = ((M_n, f_n) \xto {p_n} ... \xto {p_1} (M_0,f_0))$, an $i$-mesh $N = ((N_{i}, g_{i}) \xto {q_n} ... \xto {q_1} (N_0,N_0))$ as well as an $i$-mesh map $G : N \to M_{\leq i}$. We can use the preceding construction inductively, first pulling back $(M_{\leq i+1},f_{\leq i+1})$ along $G$, and then, for $i < k < n$, pulling back $(M_{\leq k+1},f_{\leq k + 1})$ along the pullback map $(G^*M_{k},G^*f_k) \to (M_{k},f_{k})$.
\end{rmk}

\subsubsecunnum{$n$-Mesh bundles}

Finally, our discussion of $n$-meshes would not be complete if we were not to mention that they themselves admit a natural notion of bundles. Namely, the previous definition of $n$-meshes generalizes to a definition of `constructible bundles' of $n$-meshes over some (sufficiently nice) base stratification as follows.

\begin{defn}[Constructible $n$-mesh bundles] \label{defn:constr-n-mesh-bun} A \textbf{constructible $n$-mesh bundle} over a conical stratification $(B,g)$ is a tower of constructible 1-mesh bundles $p_i : (M_i, f_i) \to (M_{i-1},f_{i-1})$, $1 \leq i \leq n$, ending in $(M_0,f_0) = (B,g)$.
\end{defn}

\nid The definition of maps of $n$-meshes carries over verbatim to the case of $n$-mesh bundles and so does the definition of topologies on hom sets, which allows us to introduce the following.

\begin{notn}[Topological category of $n$-mesh bundles] The topological category of $n$-mesh bundles and their maps will be denoted by $\tmeshbun n$.
\end{notn}

\nid Note, requiring `conicality' of the base stratification in \autoref{defn:constr-n-mesh-bun} is a way of guaranteeing that the base stratification has a fundamental category (namely, its `entrance path $\infty$-category', see \autoref{defn:tentr}). In fact, from our definition of constructibility and the weak equivalence of $n$-meshes by $n$-trusses discussed in the next section, it follows that constructible $n$-mesh bundles over a conical base are always classified by functors from the entrance path \emph{poset} of their base (to the 1-category of $n$-truss bordisms). However, we will omit a detailed discussion of mesh bundles here, and focus on the case of meshes from now on.


\section{Weak equivalence of meshes and trusses} \label{sec:meshes-and-trusses}

In this section we prove the equivalence of $n$-meshes and $n$-trusses. The equivalence will be witnessed by the `entrance path truss' functor $\ETrs : \mesh n \to \truss n$ and, conversely, by the `classifying mesh' functor $\CMsh - : \truss n \to \mesh n$. As the names suggest, the former functor is a variation of the entrance path poset functor, while the latter functor is a variation of the classifying stratification functor. Namely, the entrance path truss functor will take an $n$-mesh, given by a tower $M$ of 1-mesh bundles, to a tower $T$ of 1-truss bundles defined by applying the entrance path poset functor to the tower $M$ and endowing fibers in $T$ with 1-truss structure obtained from the 1-mesh structure of fibers in $M$. Conversely and fully analogously, the classifying mesh functor will take an $n$-truss $T$ to a tower $M$ of 1-mesh bundles obtained by applying the classifying stratification functor to the tower $T$ and endowing fibers in $M$ with 1-mesh structure obtained from the 1-truss structure of fibers in $T$ (note, however, we will see that additional care has to be taken if $T$ is not closed). In the special case of closed resp.\ open $n$-meshes and $n$-trusses with singular resp.\ regular maps, we will obtain following central theorem.

\begin{thm}[Weak equivalence of meshes and trusses] \label{thm:diagram_classification} The entrance path truss and classifying mesh functor restrict to weak equivalences
\begin{equation}
    \begin{tikzcd}[column sep=55pt]
    \sctmesh {n} \arrow[r, "\ETrs", bend left=9] & \sctruss n  \arrow[l, "\CMsh", bend left=9] & \rotmesh n \arrow[r, "\ETrs", bend left=9]  & \rotruss n  \arrow[l, "\CMsh", bend left=9]
\end{tikzcd}
\end{equation}
between the topological categories of closed $n$-meshes $\sctmesh n$ and closed $n$-trusses $\sctruss n$ (both with singular maps) respectively between open $n$-meshes $\rotmesh n$ and open $n$-trusses $\rotruss n$ (both with regular maps).
\end{thm}

\nid The proof of this theorem will take up the rest of this section. After first defining the entrance path truss functor in \autoref{ssec:ff_entr}, we will demonstrate its `conservativity' in \autoref{sssec:conservative}, as well as its `weak faithfulness' in \autoref{sssec:faithful}. We then define the classifying mesh functor in \autoref{ssec:geo-real-trusses}, first in the (simpler) case of closed trusses and then for general trusses. Finally, the proof of \autoref{thm:diagram_classification} will be assembled in \autoref{sec:constructibility-of-meshes}. This last section will also discuss several immediate applications of \autoref{thm:diagram_classification}, which we now briefly outline.

Most immediately, composing the equivalence of closed $n$-trusses and flat $n$-framed regular cell complexes established in \autoref{ch:classification-of-framed-cells} with the weak equivalence of \autoref{thm:diagram_classification}, we obtain the following corollary.

\begin{cor}[Equivalence of closed meshes and flat framed cell complexes] \label{cor:mesh-complex-cell-mesh-equiv} The equivalences of \autoref{thm:diagram_classification} and \autoref{thm:classification-of-flat-framed-cell-cplx} compose to weak equivalences $\kX \circ \ETrs : \sctmesh n \eqv \FrCDiag n : \CMsh \circ \kT$.
\end{cor}

\nid We will denote the composite $\CMsh \circ \kT$ by $\CellMsh$, called the `cell mesh' functor, and the composite $\kX \circ \ETrs$ by $\MshCplx$, called the `mesh complex' functor.

Using the notion of mesh refinements, the translation of framed cells into meshes will allow us to define `framed subdivisions' of framed cells. We will then deduce the following result, which shows that the notion of framed subdivision is `combinatorializable'.
\begin{cor}[Classifying subdivision of flat framed cells] \label{cor:classifying-subdiv-1} Up to homotopy, framed subdivisions of any framed cell $(X,\cF)$ are in correspondence with truss refinements of the truss $\kT(X,\cF)$.
\end{cor}

Finally, the equivalence of trusses and meshes in \autoref{thm:diagram_classification} will also allow us to `transfer' the duality of closed and open $n$-trusses to a dualization equivalences between closed $n$-meshes and open $n$-meshes.
\begin{cor}[Dualization of meshes] \label{thm:dualization_functors} Dualization of trusses induces weak equivalences of topological categories $\dagger : \sctmesh n \eqv \rotmesh n : \dagger$.
\end{cor}
\nid This last result is central. While the `cellular' world of framed regular cell complexes provided us with (combinatorial) topological realizations of \emph{closed} trusses, a similar topological realization of \emph{open} trusses wasn't readily available. The stratified topological notion of meshes accommodates topological realizations of both closed and open trusses as well as of their duality (namely, the resulting weak equivalence of \autoref{thm:dualization_functors} is a `geometric duality' in that it dualizes dimensions of strata). This observation underlies later, more powerful dualization operations of so-called `flat framed stratifications', and provides a general bridge between the dual worlds of `cellular' geometry and `string' geometry (see also \autoref{sec:looking-ahead}).

\subsection{Entrance path trusses of meshes} \label{ssec:ff_entr}

We formally construct the entrance path truss functor $\ETrs : \mesh n \to \truss n$. The construction is inductive, and start with the case of constructible 1-mesh bundles. We will assume that all 1-mesh bundles to have boundary-constructible base and total space stratifications (see \autoref{term:boundary-constructibility}). Note that this holds true for all 1-mesh bundles in $n$-meshes by \autoref{rmk:meshes-are-boundary-constructible}.

\begin{obs}[Restriction 1-mesh bundles to closures of strata] \label{prop:bounding-sections-for-strata} Given a constructible 1-mesh bundle $p : (M,f) \to (B,g)$, and a stratum $t$ in $f$ lying over a stratum $s = p(t)$ in $g$, then $p$ restricts on the closure $\overline t$ to a 1-mesh bundle $\rest p {\overline t}: \overline t \to \overline s$, such that fiber 1-framing $\gamma : M \into B \times \lR$ of $p$ restricts to fiber 1-framing $\rest \gamma {\overline t}$ of $\rest p {\overline t}$ whose upper and lower bounds will be denoted by $\gamma^\pm_{\overline t} : \overline s \to \overline s \times \lR$. If $t$ is singular then upper and lower bounds coincide, namely, $\gamma^\pm_{\overline t} \circ \rest p {\overline t} = \rest \gamma {\overline t}$.
\end{obs}



\begin{constr}[Entrance path 1-truss bundles] \label{prop:constr_bun_entrance_paths} Given a constructible 1-mesh bundle $p : (M,f) \to (B,g)$ we endow the entrance path poset map $\Entr(p) : \Entr(f) \to \Entr(g)$ with structure of a 1-truss bundle, yielding the `entrance path 1-truss bundle' $\ETrs(p)$.

    We first endow fibers of $\Entr(p)$ with 1-truss structure. Trivializing $p\inv(s) \iso s \times \fib {s}$ over a base stratum $s$, where $\fib(s)$ is a 1-mesh (see \autoref{obs:trivializing-1-mesh-bundles}), note that $\Entr(p)\inv(s) = \Entr(p\inv(s)) \iso \Entr(\fib s)$ canonically. Endow $\Entr(p)\inv(s)$ with a frame order $\fleq$ ordering strata using the order determined the flat 1-framing of the 1-mesh $\Entr(\fib s)$, and with a dimension map $\dim : \Entr(p)\inv(s) \to [1]\op$ mapping strata $t$ in $p\inv(s)$ to the dimension of the corresponding stratum in $\fib s$.

    It remains to check, given an entrance path $s \to r$ in the base $(B,g)$, the Boolean profunctor $R := p\inv(s \to r) : \Entr(p)\inv(s) \proto{} \Entr(p)\inv(r)$ defines a 1-truss bordisms. Pick $t \in \Entr(p)\inv(s)$. If $t$ is singular then the restriction $p : \overline t \to \overline s$ is a homeomorphism by definition of constructible mesh bundles and thus there is a unique $u \in \Entr(p)\inv(r)$ such that $R(t,u)$ holds. Otherwise, if $t$ is regular then we can apply \autoref{prop:bounding-sections-for-strata} to see that $R(t,u)$ holds if and only if $\gamma_{\overline t}^-(r) \leq \gamma(u)$ and $\gamma(u) \leq \gamma_{\overline t}^+(r)$. Comparing the two cases to the construction of `singular determined' 1-truss bordisms in \autoref{lem:truss-trans-representations} verifies that $R$ is a 1-truss bordism as claimed.
\end{constr}

\begin{eg}[Entrance path 1-truss bundles] The entrance path 1-truss bundle of the 1-mesh bundle depicted in \autoref{fig:1-mesh-bundle-over-non-cellular-base} recovers the 1-truss bundle depicted in \autoref{fig:1-truss-bundle}.
\end{eg}

\begin{defn}[Entrance path truss of a mesh] \label{defn:entrzff} Given an $n$-mesh $M$ consisting of 1-mesh bundles $p_i : (M_i,f_i) \to (M_{i-1},f_{i-1})$, its \textbf{entrance path truss} $\ETrs M$ is the $n$-truss defined by the tower of 1-truss bundles
\begin{equation}
    \Entr(f_n) \xto{\ETrs(p_n)} \Entr(f_{n-1}) \xto{\ETrs(p_{n-1})} ~ \cdots ~ \xto{\ETrs(p_2)} \Entr(f_1) \xto{\ETrs(p_1)} \Entr(f_0)
\end{equation}
where $\ETrs(p_i)$ is the entrance path 1-truss bundle of $p_i$ (see \autoref{prop:constr_bun_entrance_paths}).
\end{defn}


\nid The construction of entrance path trusses further extends to mesh maps as follows.

\begin{defn}[Entrance path truss maps of mesh maps] Given an $n$-mesh map $F : M \to N$ be a map of $n$-meshes $M$ and $N$ with components $F_i : (M_i,f_i) \to (N_i,g_i)$. The \textbf{entrance path truss map} $\ETrs F : \ETrs M \to \ETrs N$ is $n$-truss map with components $(\ETrs F)_i = \Entr(F_i) : \Entr(f_i) \to \Entr(g_i)$.
\end{defn}

\nid Note that $\ETrs F$ is regular resp.\ singular iff $F$ is regular resp.\ singular.

\begin{notn}[The entrance path truss functor]
The previous definitions yields the `entrance path truss functor', denoted by $\ETrs : \mesh n \to \truss n$, from the (ordinary) category of $n$-meshes to the category of $n$-trusses.
\end{notn}

\nid The entrance path truss functor preserves pullbacks in the following sense.

\begin{rmk}[$\ETrs$ preserves inductive pullbacks] \label{rmk:entr_preserves_pullback} Recall from \autoref{constr:pullback-of-diagrams} the notion of inductive pullbacks $G^*M$ of $n$-meshes $M$ along $(n-1)$-mesh maps $G : N \to M_{<n}$. Unwinding definitions, and writing $\ETrs {G^*f : G^*M \to M}$ as $F : S \to T$, it follows that the 1-truss bundle $q_n : S_n \to S_{n-1}$ in $S$ recovers the pullback of the 1-truss bundle $p_n : T_n \to T_{n-1}$ in $T$ along the map $F_{n-1}$ (which equals $\Entr(G_{n-1})$) in the sense of \autoref{constr:1-truss-bundle-pullback}. In this sense $\ETrs$ preserves `inductive pullbacks' of $n$-meshes.
\end{rmk}

\pauseae

    Let us now discuss the following observation: the entrance path truss functor preserves `higher structure' but, without sufficient care, fails to be an equivalence, as follows. First, form the $\Pos$-enriched category of trusses of whose objects are trusses, and whose morphisms are truss maps that are compared by (componentwise, commuting) natural transformations. By endowing these hom posets with their specialization topology (see \autoref{rmk:specialization-order}) we obtain the topological category $\ttruss n$.

\begin{prop}[$\ETrs$ as a topological functor] \label{prop:general-top-etrs} The entrance path truss functor induces topological functor $\ETrs : \tmesh n \to \ttruss n$.
\end{prop}
\begin{proof}
    We need to check that the assignment $F \mapsto \ETrs F$ on hom spaces is continuous on each hom space. This follows from continuity of $\Entr : \TStrat \to \TPos$ with respect to the specialization topology on hom posets as established in \autoref{rmk:entrz_cont}.
\end{proof}

\nid Now, importantly, the topological functor $\ETrs : \tmesh n \to \ttruss n$ is not an equivalence. This follows since natural transformations equip $\ttruss n$ with non-invertible 2-morphisms, while $\tmesh n$ as a topological category only has invertible $k$-morphisms for $k > 1$. This discrepancy in higher structure goes back to an earlier `mistake': really, the category of (conical) stratifications should not be defined as an $\infty$-categories, but as an $(\infty,2)$-categories.\footnote{Indeed, since conical stratifications have fundamental $(\infty,1)$-categories, their category should have $(\infty,2)$-categorical structure. This is similar to spaces having fundamental $(\infty,0)$-categories, and thus their category having $(\infty,1)$-categorical structure.} A natural approach in establishing the equivalence $\tmesh n \eqv \ttruss n$ would therefore be to construct $\tmesh n$ as an $(\infty,2)$-category. However, there will be no need to go down that route for us: this is because, in the cases of interest to us, `2-structure' disappears. Namely, in the case of closed trusses with singular maps $\sctruss n$, resp.\ of open trusses with regular maps $\rotmesh n$, we may topologize hom sets in either category discretely and still obtain the following topological functor.

\begin{prop}[Continuity of $\ETrs$] \label{prop:entrff-is-functor} The entrance path truss functor induces topological functors $\ETrs : \sctmesh n \to \sctruss n$ and $\ETrs : \rotmesh n \to \otruss n$.
\end{prop}

\begin{proof} This follows from the fact that singular resp.\ regular maps of closed resp.\ open trusses do not admit non-trivial natural transformations, see \autoref{lem:nat-trafo-open-ntruss}.
\end{proof}

\nid This observation has the following further variation. Namely, note that $\ETrs$ maps mesh degeneracies to truss degeneracies, and mesh coarsenings to truss coarsenings, this yields ordinary functors $\ETrs : \degmesh n \to \degtruss n$ and $\ETrs : \crsmesh n \to \crstruss n$ (see \autoref{notn:truss-crs-and-deg}). Now, again by rigidity of hom posets (this time in the case of truss degeneracies and coarsenings, see \autoref{lem:nat-trafo-open-ntruss}) we may observe the following.

\begin{prop}[More continuity of $\ETrs$] \label{prop:entrff-is-functor-on-deg-and-crs} The entrance path truss functor induces topological functors $\ETrs : \degtmesh n \to \degtruss n$ and $\ETrs : \crstmesh n \to \crstruss n$. \qed
\end{prop}

\nid In contrast to the general entrance path truss functor in \autoref{prop:general-top-etrs}, the restricted functors in \autoref{prop:entrff-is-functor} (and similarly in \autoref{prop:entrff-is-functor-on-deg-and-crs}) provide the first half of `weak equivalences' of the respective topological categories. In the next sections we will discuss the two main ingredients of proving this equivalence, which are `conservativity' and `weak faithfulness'.

\subsection{Conservativity of the entrance path truss functor} \label{sssec:conservative}

For $n$-meshes $M,M'$, the functoriality of $\ETrs$ implies that if $M \iso M'$ then $\ETrs M = \ETrs {M'}$ (note we write equality in place of isomorphism since isomorphisms between closed resp.\ between open trusses are necessarily unique). In this section we will see that the converse holds as well; namely, the functor $\ETrs$ (restricted to closed resp.\ open $n$-meshes) is `conservative' as recorded in the following result.

\begin{prop}[Conservativity] \label{prop:conservativity} For closed resp.\ open $n$-meshes $M,M'$, if $\ETrs M = \ETrs {M'}$ then $M \iso M'$.
\end{prop}

\nid Note that the statement in particular implies the conservativity of the functors $\ETrs : \scmesh n \to \sctruss n$ and $\ETrs : \romesh n \to \rotruss n$. The proof of \autoref{prop:conservativity} will take up the rest of this section. Furthermore, the properties of the construction of the isomorphism $M \iso M'$ (namely, its `continuity in families') will be reused in our later proof of `weak faithfulness' in the next section.

\begin{rmk}[More conservativity] \label{prop:conservativity-crs-deg} Note that isomorphisms of $n$-meshes are both mesh degeneracies and mesh coarsening. The preceding proof of conservativity will therefore immediately imply conservativity of the functors $\ETrs : \degmesh n \to \degtruss n$ and $\ETrs : \crsmesh n \to \crstruss n$ when restricting these functors to open and closed $n$-meshes. With enough care, the proof of conservativity of these functors can be extended to include all other $n$-meshes as well (however, since our main interest lies with the open and closed case, we will omit details of this extension).
\end{rmk}

\subsubsecunnum{Reduction to 1-mesh bundles}

As a first step we show that the statement of \autoref{prop:conservativity} may, by an inductive argument, be reduced to a statement about 1-mesh bundles. Consider closed (respectively open) $n$-meshes $M$ and $M'$ given by towers of 1-mesh bundles $p_i : (M_i,f_i) \to (M_{i-1},f_{i-1})$ resp.\ $p'_i : (M'_i,f'_i) \to (M'_{i-1},f'_{i-1})$. Assume $\ETrs M = \ETrs M'$. This implies $(\ETrs f)_{<n} = (\ETrs  M')_{<n}$, and thus we can inductively construct the isomorphism $G : M_{<n} \iso M'_{<n}$ of truncated meshes. To complete the proof of \autoref{prop:conservativity} it now remains to construct a 1-mesh bundle isomorphisms as the following remark records.

\begin{rmk}[Inductive reduction to 1-mesh bundle isomorphism] \label{rmk:kappa_reduction} Define the $n$-mesh $G^*M'$ to be the inductive pullback of $M'$ along $G$  (see \autoref{constr:pullback-of-diagrams}); write the top-level 1-mesh bundle of $G^*M'$ as $\tilde p_n : (\tilde M_n,\tilde f_n) \to (M_{n-1},f_{n-1})$ and denote the canonical map $(\tilde M_n,\tilde f_n) \to (M'_n, f'_n)$ by $F$, as shown in the diagram below (where we use \autoref{notn:abbreviating-stratifications} to abbreviate stratified spaces $(X,f)$ by $f$).
\begin{equation}
\begin{tikzcd}
     \tilde f_n \arrow[r, "F"] \arrow[d, "\tilde p_n"']
\arrow[dr, phantom, "\lrcorner" , very near start, color=black]
 & f'_n \arrow[d, "p'_n"] \\
    \lvl {n-1} f \arrow[r, "G"]                           & f'_{n-1}
\end{tikzcd}
\end{equation}
Note, since $G$ is assumed to be an isomorphism, $F$ too is an isomorphism. To prove \autoref{prop:conservativity} it remains to construct a bundle isomorphism $\kappa^{p}_{\tilde p}$ as shown below.
\begin{equation}
\begin{tikzcd}[baseline=(W.base)]
f_n \arrow[rr, "\kappa^{p_n}_{\tilde p_n}"]{}[swap]{\sim} \arrow[rd, "p_n"'] &  & {\tilde f_n} \arrow[ld, "\tilde p_n"] \\
 & |[alias=W]| \lvl {n-1} f &
\end{tikzcd} \qedhere
\end{equation}
\end{rmk}

\nid Recall from \autoref{rmk:closed-open-meshes-are-open-regular-cell} that closed (respectively open) $n$-meshes are open regular cell stratifications at each level. The `inductive step' for \autoref{prop:conservativity} is provided by the following result.

\begin{prop}[Isomorphisms between constructible 1-mesh bundles] \label{prop:kappa_constr} For constructible closed (resp.\ open) 1-mesh bundles $p : (M,f) \to (B,g)$ and $\tilde p : (\tilde M,\tilde f) \to (B,g)$ whose base $(B,g)$ is open regular cell and whose entrance path 1-truss bundles $\ETrs(p) = \ETrs(\tilde p)$ coincide, there is a bundle isomorphism $\kappa^p_{\tilde p} : f \iso \tilde f$.
\end{prop}

\nid We will construct the isomorphism $\kappa^p_{\tilde p} : f \iso \tilde f$ only in the case of closed 1-mesh bundles. The case of open 1-mesh bundles may be reduced to the closed case by fiberwise compactification (see \autoref{constr:fiberwise_compactification}). The construction of  $\kappa^p_{\tilde p}$ will need care, and as a motivation for our approach we will first discuss how \emph{not} to construct $\kappa^p_{\tilde p}$.

When defining the bundle isomorphism $\kappa^p_{\tilde p} : f \iso \tilde f$ fiberwise, note that when traveling along an entrance paths between two strata $r \to s$ in the base $(B,g)$, new singular strata can appear in the `special' fiber over $s$, which were not present in the `generic' fiber over $r$. This implies that, for instance, we generally cannot define $\kappa^p_{\tilde p}$ by fiberwise mapping point strata to point strata and extending this mapping of points linearly to interval strata (note that fibers in 1-mesh bundle inherit linear structure from the standard linear structure of $\lR$ via the fiber 1-framing embedding). We illustrate the failure of this approach in the following example.

\begin{eg}[Failure of continuity in fiberwise linear interpolating] To illustrate the failure of continuity that can happen in attempting to fiberwise linearly map between consider the bundle $p : (M,f) \to (B,g)$, $\tilde p : (\tilde M,\tilde f) \to (B,g)$ as shown in \autoref{fig:conservativity-inductive-step-iso} whose entrance path 1-truss bundles coincide (which we depicted, via their fiber 1-framings, as embedded in $B \times \lR$). We indicated an entrance path in $(B,g)$, together with a generic fiber as well as a special fiber in both $p$ and $\tilde p$. If we were to build a bundle isomorphism fiberwise, by first identifying point strata of fibers (as indicated by the mappings on the right in \autoref{fig:conservativity-inductive-step-iso}) and then linearly interpolating these mappings on interval strata, we would end up with a \emph{discontinuous} bundle isomorphism between $p$ and $\tilde p$.
\begin{figure}[ht]
    \centering
    \def\svgwidth{1\columnwidth}
    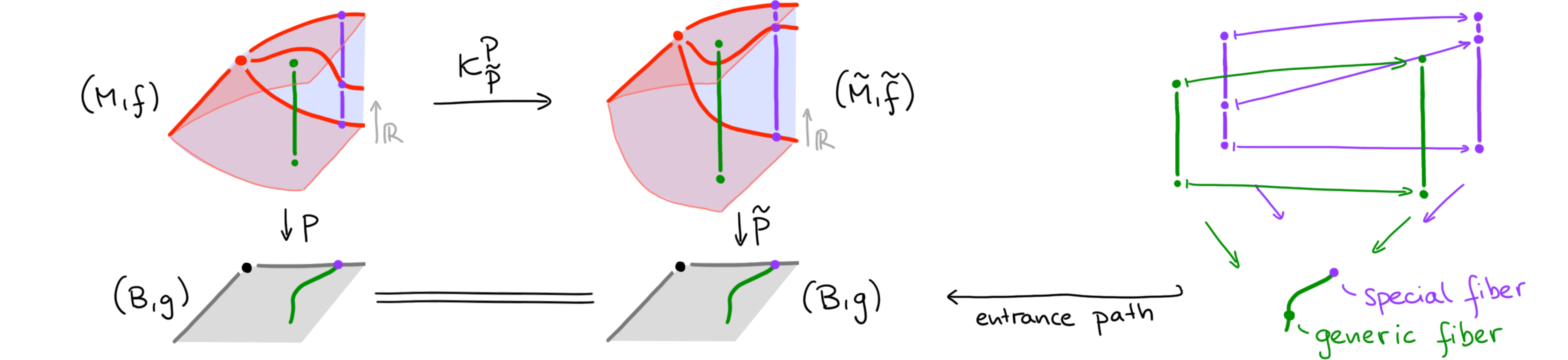

    \caption[Inductive step for conservativity of the entrance path truss functor.]{The isomorphism needed for the inductive step in the proof of conservativity of the entrance path truss functor generally cannot be constructed by simple fiberwise interpolation.}
    \label{fig:conservativity-inductive-step-iso}
\end{figure}
\end{eg}

\nid Our strategy to ensure continuity in our construction of $\kappa^p_{\tilde p}$ will be to use `affine combinations' of maps on generic and special fibers whenever we are getting `critically close' to special fibers (this will be referred to as fibers in the `critical region' of a stratum), and use simple linear interpolation otherwise.

\subsubsecunnum{Constructing critical regions}

We start with the construction of `critical regions' in the base stratification $(B,g)$, over which fiberwise maps will be defined as certain affine combinations of `general' and `special fiber maps' later on. Recall we assume the base stratification $(B,g)$ to be an open regular cell stratification.

\begin{term}[Interval contour] Let $p : (M,f) \to (B,g)$ be a constructible closed 1-mesh bundle, and consider a regular stratum $s$ of $f$ lying over a stratum $r = p(s)$ of $g$. The `interval contour under $s$', denoted by $\bdreg s$, is the subset of $\partial r = \overline r \setminus r$, over which the restricted bundle $\rest p {\overline s} : \overline s \to \overline r$ (see \autoref{prop:bounding-sections-for-strata}) has interval fibers.
\end{term}

\begin{figure}[ht]
    \centering
    \def\svgwidth{1\columnwidth}
    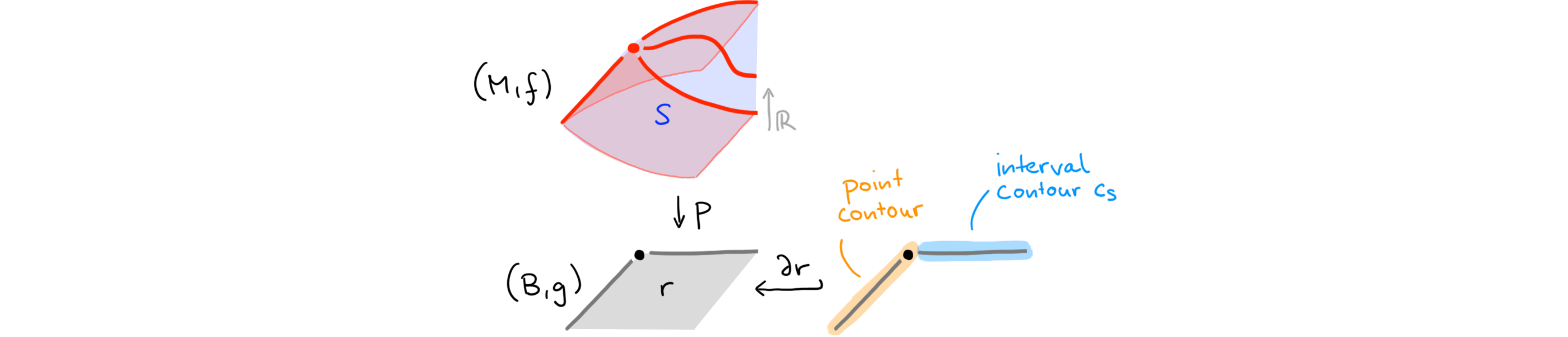

    \caption{Point and interval contour.}
    \label{fig:the-non-degeneracy-boundaries-under-2-different-strata}
\end{figure}

\nid The complement of $\bdreg s$ in $\partial r$ will also be called the `point contour under $s$'. For the bundle $p$ from \autoref{fig:conservativity-inductive-step-iso}, we highlight the point and interval contour of a chosen stratum $s$ in \autoref{fig:the-non-degeneracy-boundaries-under-2-different-strata} below

\begin{rmk}[Contours only depend on truss structure] \label{rmk:proj_bound_depend_on_truss_struct}
    Importantly, if $\tilde p : (\tilde M, \tilde f) \to (B,g)$ is another constructible closed 1-mesh bundle over $(B,g)$ and $\ETrs p = \ETrs \tilde p$  (which in particular allows us to identify strata $s$ of $f$ with strata $\tilde s$ of $\tilde f$) then $\bdreg s$ and $\bdreg {\tilde s}$ coincide as subspaces of $B$.
\end{rmk}

We will now thicken the interval contour by `pushing it into' the interior of the base stratum $r$ under $s$. This thickening will provide the `critical region' for our later definition of $\kappa^p_{\tilde p}$.

\begin{constr}[Critical regions] \label{constr:critical-regions} For a constructible closed 1-mesh bundle $p : (M,f) \to (B,g)$ with open regular cell base and a regular stratum $s$ lying over $r$ in $f$, we construct a `critical region' $\crreg s \into r$ as follows.
    First, take $Y$ to be a regular cell complex containing $(B,g)$ as an open constructible substratification. Denote by $\bdregclos s$ the closure of the interval contour in the boundary $S^{k-1} \iso \partial r$ of the open $k$-cell $r$ in $Y$, and set $\bdregbd s = \bdregclos s \setminus \bdreg s$.
    Now `thicken' the interval contour's interior by defining $\crregclos s$ to be the quotient of $\bdregclos s \times [0,1]$ by $\bdregbd s \times [0,1]$. Define $\crreg s \into \crregclos s$ to be the subspace $\bdreg s \times (0,1)$.
    Form the gluing $\iG_s$ by gluing $\overline r \iso D^n$ (where the closure is taken in $Y$) and $\crreg s$ along $\bdregclos s \into \overline r$ resp.\ $\bdregclos s \iso \bdregclos s \times \Set{0} \into \crregclos s$.
    A `choice of critical region for $s$' is a choice of an isomorphism $\iR_s : \iG_s \iso \overline r$ which maps $\bdregclos s \iso \bdregclos s \times \Set{1} \into \iG_s$ identically to $\bdregclos s \into \partial r$, and $\partial r \setminus \bdreg s \into \iG_s$ identically to $\partial r \setminus \bdreg s \into r$.
    Note such isomorphism can always be chosen since $\bdreg s \into \partial r = S^{k-1}$ is an open subset. Finally, restricting $\iR_s : \iG_s \iso \overline r$ to $\crreg s$ yields the `critical region' $\crreg s \into \iG_s$
\end{constr}

\begin{figure}[ht]
    \centering
    \def\svgwidth{1\columnwidth}
    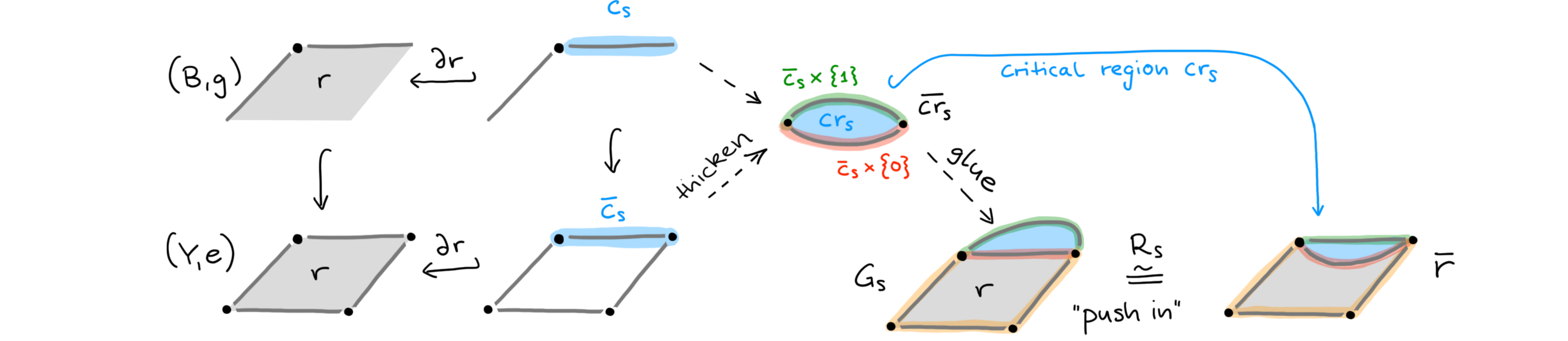

    \caption[Choosing a critical region for a regular stratum.]{Choosing a critical region in $r$ for a regular stratum $s$.}
    \label{fig:choosing-a-critical-region-for-a-regular-stratum}
\end{figure}

\nid For our previous example of $p$ and choice of $s$ in \autoref{fig:the-non-degeneracy-boundaries-under-2-different-strata} we illustrate a choice of critical region $\iR_s$ for $s$ in \autoref{fig:choosing-a-critical-region-for-a-regular-stratum}.

\nid Finally, note a valid choice of critical regions only depends on the entrance path truss structure of $p$, as recorded by the following remark.

\begin{rmk}[Critical regions only depend on 1-truss structure over the base stratification] \label{rmk:interpolation_reg_depend_on_truss_struct}
Importantly, if $\tilde p : (\tilde X, \tilde f) \to (Y,g)$ is another constructible closed 1-mesh bundle and $\ETrs(p) = \ETrs(\tilde p)$ then any choice of critical region for $p$ is also choice of critical regions for $\tilde p$.
\end{rmk}

\subsubsecunnum{Constructing the bundle isomorphism}

Having a notion of critical regions at hand, we can proceed to construct the bundle isomorphism $\kappa^p_{\tilde p}$, which will yield proofs of both \autoref{prop:kappa_constr} and \autoref{prop:conservativity}.

\begin{constr}[The bundle isomorphism $\kappa$] \label{constr:kappa} Consider constructible closed 1-mesh bundles $p : (M,f) \to (B,g)$ and $\tilde p : (\tilde M,\tilde f) \to (B,g)$ with open regular cell base $g$ and such that $\ETrs(p) = \ETrs(\tilde p)$. Fix choices of critical regions $\iR_s$ for each regular stratum in $p$ (see \autoref{constr:critical-regions}) which, by \autoref{rmk:interpolation_reg_depend_on_truss_struct} equally provides a choice of critical regions for $\tilde p$. We will define the bundle isomorphism $\kappa^p_{\tilde p}$ fiberwise by maps $\kappa^p_{\tilde p}(x,-) : p\inv(x) \to \tilde p\inv(x)$ over points $x \in r$ where $r$ is a stratum in $g$.

    We argue inductively in $\dim(r)$. If $\dim(r) = 0$ then $\kappa^p_{\tilde p}(x,-)$ is simply defined by mapping point strata of $p\inv(x)$ monotonicly to point strata of $\tilde p\inv(x)$ (note, since $\ETrs(p) = \ETrs(\tilde p)$ both in particular have the same number of point strata) and extending this mapping linearly to interval strata. Next assume $\dim(r) > 0$. Again, we define $\kappa^p_{\tilde p}(x,-)$ is simply defined by mapping point strata of $p\inv(x)$ monotonicly to point strata of $\tilde p\inv(x)$ (note, since $\ETrs(p) = \ETrs(\tilde p)$ both in particular have the same number of point strata). Note that interval strata $s_x$ in $p\inv(x)$ are restrictions of regular strata $s$ in $(M,f)$, and correspond to interval strata $\tilde s_x$ in $p\inv(x)$ obtained by restricting the corresponding regular strata $\tilde s$ in $(\tilde M,\tilde f)$. We define the restricted map $\kappa^p_{\tilde p}(x,-) : s_x \to \tilde s_x$. If $x \notin \crreg s$ lies outside the choice of critical region for $s$, define $\kappa^p_{\tilde p}(x,-) : s_x \to \tilde s_x$ by linearly extending the mapping $\kappa^p_{\tilde p}(x,-)$ on the endpoints of the interval stratum $s_x$. Otherwise, if $x \in \crreg s$, by \autoref{constr:critical-regions} we have $\crreg s = \bdreg s \times (0,1)$ and can thus write $x$ as a pair $(u \in \bdreg s,t \in (0,1))$. This yields points $x_0 = \iR_s(u,0) \in r$ and $x_1 = \iR_s(u,1) \in \partial r$. Define $\kappa^p_{\tilde p}(x,-) : s_x \to \tilde s_x$ by the affine combination
    \begin{equation}
        \kappa^p_{\tilde p}(x,-) = t \kappa^p_{\tilde p}(x_1,-) + (1-t)\kappa^p_{\tilde p}(x_0,-)
    \end{equation}
    where the first term on the right is defined since $x_1$ lies in a stratum of dimension lower than $r$, and the second term is defined since $x_0$ lies in $r$ but outside the critical region. One checks that this continuously extends the domain of definition of $\kappa^p_{\tilde p}$ to all fibers over $r$, and thus completes our inductive construction of the bundle isomorphism $\kappa^p_{\tilde p}$.
\end{constr}

The construction of $\kappa^p_{\tilde p}$ now makes proofs of earlier statement immediate.

\begin{proof}[Proof of \autoref{prop:kappa_constr}] The preceding construction defines $\kappa^{p}_{\tilde p}$ as required in the statement of \autoref{prop:kappa_constr}.
\end{proof}

\begin{proof}[Proof of \autoref{prop:conservativity}] By \autoref{rmk:kappa_reduction}, this follows as a corollary to \autoref{prop:kappa_constr}.
\end{proof}

\subsubsecunnum{Continuity properties of $\kappa$ construction} We mention two important properties of the construction of $\kappa^p_{\tilde p}$ in \autoref{constr:kappa}. Firstly, the construction `preserves identities' as recorded by the following remark.

\begin{obs}[Constructing $\kappa$ on identical bundles] \label{rmk:kappa_identity} Consider constructible closed (resp.\ open) 1-mesh bundles $p : (M,f) \to (B,g)$ and $\tilde p : (\tilde M,\tilde f) \to (B,g)$ with open regular cell base $g$ and such that $\ETrs(p) = \ETrs(\tilde p)$. Identify $M$ (resp.\ $\tilde M$) as a subspace of $B \times \lR$ using its fiber 1-framing $\gamma : M \into B \times \lR$ (resp.\ $\tilde \gamma : \tilde M \into B \times \lR$). If $(M,f)$ and $(\tilde M,\tilde f)$ are identical stratifications then the inductive construction of $\kappa^p_{\tilde p}$ (in \autoref{constr:kappa}) outputs the bundle identity map $\id : p = \tilde p$.
\end{obs}

\nid Secondly, the construction is `continuous in families'. To make this precise, we first introduce a notion of families of 1-mesh bundles.

\begin{defn}[Families of 1-mesh bundles] For a connected topological space $Z$ and a open regular cell stratification $(B,g)$, a \textbf{$Z$-family} of constructible closed (resp.\ open) 1-mesh bundles over $g$ is a constructible closed (resp.\ open) 1-mesh bundle $p : (M,f) \to Z \times (B,g)$.  For $z \in Z$ denote by $p_z : (M_z,f_z) \to (B,g)$ the restriction of $p$ to the subspace $B \iso \Set{z} \times B \into Z \times B$, called the \textbf{$z$-slice} of $p$.
\end{defn}

\nid Note that, using \autoref{prop:constr-implies-cellular}, all bundles $p_z$ in a $Z$-family of bundles are open regular cell bundles (however, $p$ itself may not be open regular cell).

\begin{term}[$\ETrs$-constant families] A $Z$-family of constructible closed (resp.\ open) 1-mesh bundles $p : (M,f) \to Z \times (B,g)$ is call `$\ETrs$-constant' if the canonical inclusion $\ETrs(p_z) \into \ETrs(p)$ is a 1-truss bundle isomorphism for each $z \in Z$.
\end{term}

\begin{rmk}[Constancy of entrance path truss bundle in families] \label{rmk:family_entr_invariance} Given a $Z$-family $p$ for a path-connected space $Z$, then $p$ is automatically $\ETrs$-constant.
\end{rmk}

\nid The `continuity in families' of our construction of $\kappa$ now takes the following form.

\begin{obs}[Constructing $\kappa$ for families] \label{rmk:kappa_continuity} For an open regular cell stratifications $(B,g)$, consider $\ETrs$-constant $Z$-families of constructible closed (resp.\ open) 1-mesh bundles $p : (M,f) \to Z \times (B,g)$ and $\tilde p : (\tilde M,\tilde f) \to Z \times (B,g)$, such that $\ETrs(p) = \ETrs(\tilde p)$. Pick any $z_0 \in Z$, and a chose critical regions for $p_{z_0}$. This choice equally provides critical regions to all bundles $p_z$ and $\tilde p_z$ for all $z \in Z$, and we may thus construct the 1-mesh bundle isomorphisms $\kappa^{p_z}_{\tilde p_z} : p_z \iso \tilde p_z$ using \autoref{constr:kappa}. The construction of these bundle isomorphisms is now `continuous in $z \in Z$' in that we obtain a continuous bundle isomorphisms $\kappa^{p}_{\tilde p} : (M,f) \to (\tilde M,\tilde f)$ whose restriction to $(M_z,f_z) \to (\tilde M_z,\tilde f_z)$ is defined to be $\kappa^{p_z}_{\tilde p_z}$.
\end{obs}

Finally, let us mention a particular way of constructing $\ETrs$-constant $Z$-families of 1-mesh bundles. Namely, using a restricted version of cartesian closedness in the context of stratified spaces, $\ETrs$-constant $Z$-families of constructible 1-mesh bundles may be constructed by pullback along `$\Entr$-constant $Z$-families of stratified maps' as the following remark explains.

\begin{rmk}[Families of bundles from pullbacks along families of maps]  \label{rmk:families_by_pullback}
    Consider open regular cell stratifications $(B,g)$ and $(B,\tilde g)$ and let $F : Z \to \TStrat(g,\tilde g)$ be a continuous map from a space $Z$ to the space of stratified maps between $g$ and $\tilde g$, such that $F$ is constant on entrance paths (that is, $\Entr \circ F : Z \to \TPos(\Entr(g), \Entr(\tilde g))$ is constant). By $\Top$-tensoredness of stratified spaces (see \autoref{rmk:strat_cc}), we can consider $F$ as a stratified map $F : Z \times g \to \tilde g$. Now, given a constructible closed (resp.\ open) 1-mesh bundle $\tilde p : (\tilde M, \tilde f) \to (\tilde B, \tilde g)$ we can construct a $\ETrs$-constant $Z$-family of constructible closed (resp.\ open) 1-mesh bundle as the pullback $F^*\tilde p$ of $\tilde p$ along $F$ (see \autoref{constr:pullback_of_constr_diagrams_bundles}).
\end{rmk}

\subsection{Weak faithfulness of the entrance path truss functor} \label{sssec:faithful} Our next goal will be to show that the entrance path truss functor $\ETrs$ is (in a weak sense) a `faithful' functor of topological categories in the following sense.

\begin{prop}[Weak faithfulness of entrance path truss functor] \label{prop:faithfulness} Given two closed (resp.\ open) $n$-meshes $M$ and $M'$, the hom space map $\ETrs : \sctmesh n (M,M') \to \sctruss n (\ETrs M,\ETrs M')$ (resp.\ $\ETrs : \rotmesh n (M,M') \to \rotruss n (\ETrs M,\ETrs M')$) has empty or weakly contractible preimages.
\end{prop}

\nid The proof of \autoref{prop:faithfulness} will take up the rest of this section. We will give the proof only in the case of closed $n$-meshes, noting that the case of open $n$-meshes is fully analogous (in fact, the proof for the closed case applies verbatim to the open case when replacing `closed' by `open' meshes resp.\ trusses, and `singular' by `regular' maps).

\begin{rmk}[Entrance path truss functor is `weakly fully faithful']
Once we have constructed the weak inverse of $\ETrs$ it will follow that fibers of the hom space maps of $\ETrs$ are, in fact, never empty.
\end{rmk}

\begin{rmk}[More weak faithfulness] \label{prop:faithfulness-crs-deg} Similar to \autoref{prop:conservativity-crs-deg}, the proof of \autoref{prop:faithfulness} may be adapted to show that the functors $\ETrs : \degtmesh n \to \degtruss n$ and $\ETrs : \crstmesh n \to \crstruss n$ are weakly faithful functors as well (the proof goes through without change in the closed resp.\ open case, and needs a little more care in the case of general meshes).
\end{rmk}

As a first step we show that the statement of \autoref{prop:faithfulness} may, by an inductive argument, be reduced to a statement about 1-mesh bundles. For $m \in \lN$, denote by $D^{m+1}$ the closed $(m+1)$-ball and by $S^{m}$ its boundary. Recall, a topological space $U$ is weakly contractible if for all maps $\zeta : S^m \to U$ there exists a map $\theta : D^{m+1} \to U$ such that $\rest \theta {S^m} : S^m \to U$ equals $\zeta$. We call $\theta$ a `filler' for $\zeta$.

Let $M$ and $M'$ be closed $n$-meshes consisting of 1-mesh bundles $p_i : (M_i,f_i) \to (M_{i-1},f_{i-1})$ resp.\ of $p'_i : (M'_i,f'_i) \to (M'_{i-1},f'_{i-1})$. Consider a map $\zeta : S^m \to  \sctmesh n (M,{M'})$ such that $\ETrs \circ \zeta$ is constant (in other words, $\zeta$ maps into a single fiber of $\ETrs$). Note that, by rigidity of singular truss maps of closed trusses (see \autoref{lem:nat-trafo-open-ntruss}) this constancy condition is satisfied automatically except when $m = 0$. Recall that truncation of meshes is a topological functor (see \autoref{rmk:truncating-topological}). Truncating $\zeta$ to degrees below $n$, we obtain the map $\zeta_{<n}: S^m \to \sctmesh {n-1} (M_{<n}, M'_{<n})$ which we denote by $\beta$. Arguing inductively, $\beta$ has a filler which we denote by $\eta : D^{m+1} \to \sctmesh {n-1} (M_{<n}, M'_{<n})$. Using the $\Top$-tensoredness of stratified spaces (see \autoref{rmk:strat_cc}), we may consider $\zeta$ as a stratified map $S^m \times f_n \to {f'_n}$, $\beta$ as a stratified map $S^m \times f_{n-1} \to f'_{n-1}$, and $\eta$ as a stratified map $D^{m+1} \times f_{n-1} \to f'_{n-1}$. To show that $\zeta$ has a filler it will be sufficient to prove the following.

\begin{prop}[Lifting fillers in closed 1-mesh bundles] \label{prop:filler-lifts-in-1diag-bun} Consider constructible closed 1-mesh bundles $p : (M,f) \to (B,g)$ and $\tilde p : (\tilde M,\tilde f) \to (\tilde B,\tilde g)$ with open regular cell base, and maps $\zeta : S^m \times f \to \tilde f$, $\beta : S^m \times g \to \tilde g$ such that, for each $e \in S^m$, $(\zeta(e,-),\beta(e,-)) : p \to \tilde p$ is a 1-mesh bundle map (if $m = 0$, further assume $\ETrs(\zeta(e,-),\beta(e,-))$ is independent of $e \in S^m$). Then any filler $\eta : D^{m+1} \times g \to \tilde g$ of $\beta$ `lifts' to a filler $\theta : D^{m+1} \times f \to \tilde f$ of $\zeta$ such that, for each $e \in D^{m+1}$, $(\theta(e,-),\eta(e,-)) : p \to \tilde p$ is a 1-mesh bundle map.
\end{prop}

\nid It will be convenient to consider $D^{m+1}$ as the quotient of $[0,1] \times S^m$ by the subset $\Set{1} \times S^m$. In particular, we will regard $\eta$ as a map $\eta : [0,1] \times S^m \times g \to \tilde g$ such that $\eta(1,-)$ is constant.

To construct the filler $\theta$ of $\zeta$ that lifts the filler $\eta$ of $\beta$ as claimed in \autoref{prop:filler-lifts-in-1diag-bun}, we will proceed in two steps. In the first step, by `pulling back along $\eta$', we construct `homotopy 1' $\theta_1 : [0,1] \times S^m \times f \to \tilde f$ which homotopes $\theta_1(0,-) = \zeta$ into a map $\theta_1(1,-) : S^m \times f \to \tilde f$ that descends to a map of base stratifications $S^m \times g \to \tilde g$ that is constant in $S^m$.  The second step, using `fiberwise contractions', constructs `homotopy 2' $\theta_2 : [0,1] \times S^m \times f \to \tilde f$, which homotopes $\theta_2(0,-) = \theta_1(1,-)$ into a map $\theta_2(1,-) : S^m \times f \to \tilde f$ that is constant on fibers as well and thus fully constant in $S^m$. Concatenating the homotopies $\theta_1$ and $\theta_2$ will provide us with the required filler $\theta$ of $\zeta$.

\begin{proof}[Proof of \autoref{prop:filler-lifts-in-1diag-bun}]
Define a constructible closed 1-mesh bundle $\widehat p$ by pulling back $p$ along $\beta$ and define a map $\widehat \zeta$ as the factorization of $\zeta$ through this pullback as shown below.
\begin{equation}
\begin{tikzcd}
S^m \times \empty f \arrow[d, "S^m \times p"'] \arrow[r, "{\widehat \zeta}"'] \arrow[rr, "\zeta", bend left=33] & \beta^*  \tilde f \arrow[dr, phantom, "\lrcorner" , very near start, color=black] \arrow[d,"\widehat p"] \arrow[r] & \tilde f \arrow[d, "\tilde p"] \\
S^m \times g \arrow[r, "\id"', equal] & S^m \times g \arrow[r, "\beta"'] & \tilde g
\end{tikzcd}
\end{equation}
Note that $\widehat p$ is an $\ETrs$-constant $S^m$-family of constructible closed 1-mesh bundles over $g$ (see \autoref{rmk:families_by_pullback}). We may trivially turn this into a $\ETrs$-constant $([0,1] \times S^m)$-family by taking the product $[0,1] \times -$: the resulting constructible closed 1-mesh bundle $[0,1] \times \widehat p$ is bundle isomorphic to the constructible closed 1-mesh bundle $\widecheck p$ defined by the pullback on the right below.
\begin{equation}
\begin{tikzcd} {{[0,1] \times \beta^*{\tilde f}}} \arrow[r, "\kappa"]{}[swap]{\sim}  \arrow[d, "{{[0,1] \times \widehat p}}"'] & \eta^*\tilde f \arrow[dr, phantom, "\lrcorner" , very near start, color=black] \arrow[r] \arrow[d, "\widecheck p"] & \empty \tilde f \arrow[d, "\tilde p"] \\
{[0,1]\times S^m \times g} \arrow[r, "\id"'] & {[0,1]\times S^m \times g} \arrow[r, "\eta"']  & \tilde g
\end{tikzcd}
\end{equation}
Indeed, $\kappa$ can be constructed using \autoref{rmk:kappa_continuity} since $g$ is assumed open regular cell. Our `homotopy 1' map $\theta_1 : [0,1] \times S^m \times f \to \tilde f$ is now simply defined as the composite
\begin{equation}
    [0,1] \times S^m \times \empty f \xto {[0,1] \times {\widehat \zeta}} [0,1] \times \beta^*{\tilde f} \xto {\kappa} \eta^*{\tilde f} \to \tilde f \quad.
\end{equation}
Since the $\kappa$ construction preserves identities (see \autoref{rmk:kappa_identity}) and since $\beta = \eta (0,-)$, we find that $\kappa(0,-)$ is the identity on $\beta^*{\tilde f}$. Thus, homotopy 1 satisfies $\theta_1(0,-) = \zeta$, and `lifts ' $\eta$ in the sense that
\begin{equation}
\begin{tikzcd}
{[0,1]\times S^m \times f} \arrow[r, "\theta_1"] \arrow[d, "{[0,1] \times S^m \times p}"'] & \empty \tilde f \arrow[d, "\tilde p"] \\
{[0,1]\times S^m \times g} \arrow[r, "\eta"']                                              & \tilde g
\end{tikzcd}
\end{equation}
This completes the construction of the first part of the homotopy.

It remains to construct the `homotopy 2' map $\theta_2 : [0,1] \times S^m \times \empty f \to \empty g$ such that $\theta_2(0,-) = \theta_1(1,-)$. We may construct the homotopy $\theta_2$ `fiberwise', by convexly combining fiberwise maps $\theta_1(1,e,x) : (p\inv(y),f) \to (\tilde p \inv \eta(1,e,y),\tilde f)$ on fiber over points $y \in B$ (note that $\eta(1,e,y)$ is in fact independent of $e \in S^m$). Namely, pick any $e_0 \in S^m$, and, for $t \in [0,1]$, $e \in S^m$ and $y \in B$, define the restriction of $\theta_2(t,e,-)$ to the fiber over $y$ to be the map
\begin{align}
    \theta_2(t,e,-) :  (p\inv(y),f) &\to (\tilde p \inv \eta(1,e,y),\tilde f)\\
x &\mapsto (1-t)\cdot \theta_1(1,e,x) + t\cdot\theta_1(1,e_0,x)
\end{align}
Note that $\theta_2(t,e,-)$ is indeed a 1-mesh map for all $t$ and $e$ since $\ETrs \theta_1(1,e,x) = \ETrs \theta_1(1,e_0,x)$ induce the same maps on 1-trusses (which in turn follows since, by assumption, $\ETrs(\zeta(e,-),\beta(e,-))$ is independent of $e \in S^m$). Note also at $t = 1$, the map $\theta_2(t,e,-)$ becomes independent of $e \in S^m$. We can chain the homotopies $\theta_1, \theta_2$ into a single homotopy
\begin{equation}
\theta := \theta_1 \ast \theta_2 : [0,1] \times S^m \times f \to \tilde f
\end{equation}
which defines the required filler $\theta$ of $\zeta$ lifting the filler $\eta$ of $\beta$ as required.
\end{proof}

\begin{proof}[Proof of \autoref{prop:faithfulness}] Since the statement of \autoref{prop:faithfulness} inductively reduces to the statement of \autoref{prop:filler-lifts-in-1diag-bun}, the former statement is now proven.
\end{proof}

\subsection{Classifying meshes} \label{ssec:geo-real-trusses}

Every $n$-truss $T$ has a `classifying mesh' $\CStr T$, and the construction of classifying meshes of trusses will provide an inverse to the construction of entrance path trusses of meshes (discussed in \autoref{ssec:ff_entr}). The situation is analogous to the classifying stratification functor $\CStr : \Pos \to \Strat$ being (left) inverse to entrance path poset functor $\Entr : \Strat \to \Pos$. However, only in the case of closed trusses will we obtain the classifying mesh construction as a direct analog of the classifying stratification construction. In the case of general trusses, more care needs to be taken to ensure that dimensions of strata do not incorrectly `degenerate'.

To illustrate this difficulty in constructing classifying meshes of general trusses consider the following basic case. The classifying stratification of the poset underlying the closed 1-truss with one element is a point: this thus correctly produces the underlying space of the closed 1-mesh with one stratum. However, the classifying stratification of the poset underlying the \emph{open} 1-truss with one element is again a point, and this does not equal the underlying space of the \emph{open} 1-mesh with one stratum (which is an open interval). As we will see, the problem of constructing classifying meshes of general trusses may be reduced to the case of closed trusses, by first `compactifying' trusses, and then applying the classifying mesh functor for closed trusses.

\subsubsecunnum{Classifying meshes of closed trusses}

We construct classifying meshes of closed $n$-trusses. We start by discussing the case of closed 1-truss bundles. Recall the construction of classifying stratifications of posets (\autoref{term:classifying-stratifications} and \autoref{term:classifying-stratifications-maps}).

\begin{constr}[Classifying 1-mesh bundles of closed 1-truss bundles] \label{prop:geo_closed_1_truss_bun} Given a closed 1-truss bundle $p : T \to X$, we endow the classifying stratified map $\CStr p : \CStr T \to \CStr X$ with the structure of a closed 1-mesh bundle, yielding the `classifying 1-mesh bundle' $\CMsh p$ of $p$. Recall, the classifying stratification $\CStr T$ stratifies the geometric realization $\abs T$. We need to define a fiber 1-framing $\gamma : \abs T \into \abs X \times \lR$ for $\CStr p$. Order-preservingly identify objects in the (frame ordered) fiber $(p\inv(x),\fleq)$ over $x \in X$ with objects $i$ in a total order $[m_x] = (0 \to ... \to m_x)$. Since objects in $T$ correspond to vertices in $\abs{T}$, we then define $\gamma$ to map vertices $i \in \abs{p}\inv(x) \subset \abs{T}$ to vertices $(x,i) \in \abs{X}\times \lN \into \abs{X} \times \lR$, and further extend this mapping on vertices linearly to simplices in $\abs{T}$. Using truss induction (see \autoref{sec:truss-induction}) one verifies that $\gamma$ is indeed an embedding such that the bounding maps $\gamma^\pm : \abs{X} \to \abs{X} \times \lR$ are continuous. This endows $\CStr p$ with the structure of a closed 1-mesh bundle.
\end{constr}

\begin{eg}[Classifying 1-mesh bundles of 1-truss bundles] In \autoref{fig:classifying-1-mesh-bundles-of-1-truss-bundles} we depict a closed 1-truss bundle $p : T \to X$ on the left (note that we only depict generating arrows, see \autoref{rmk:generating-arrows}), and on the right its classifying 1-mesh bundle $\CMsh p : \CStr T \to \CStr X$.
\begin{figure}[ht]
    \centering
    \def\svgwidth{1\columnwidth}
    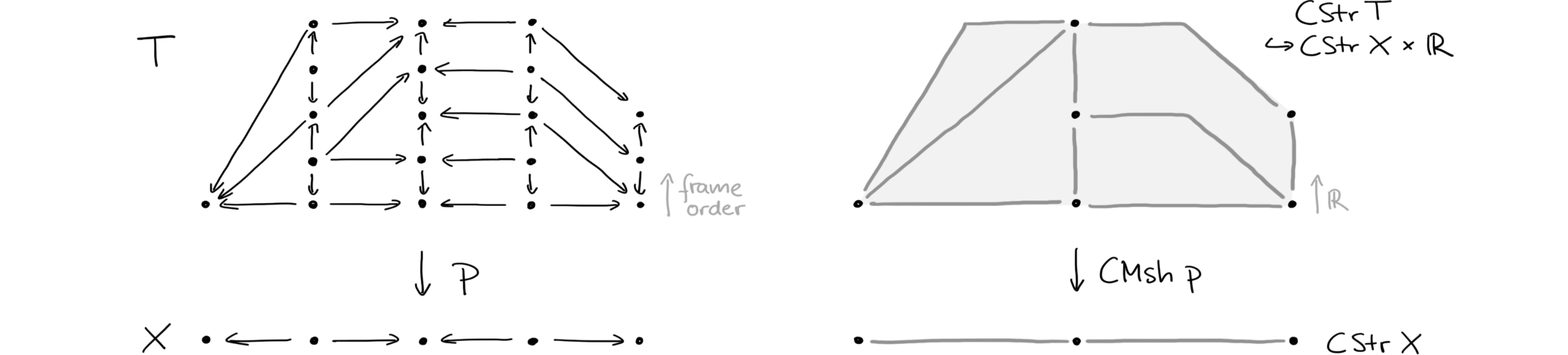

    \caption{Classifying 1-mesh bundle of a 1-truss bundle.}
    \label{fig:classifying-1-mesh-bundles-of-1-truss-bundles}
\end{figure}
\end{eg}

\begin{defn}[Classifying meshes of closed $n$-trusses] \label{constr:geo_closed} Given a closed $n$-truss $T$, consisting of 1-truss bundle $p_i : T_i \to T_{i-1}$, its \textbf{classifying mesh} $\CMsh T$ is the closed $n$-mesh defined by the tower of 1-mesh bundles
\begin{equation}
\CStr {T_n} \xto {\CStr {p_n}} \CStr {T_{n-1}} \xto {\CStr {p_{n-1}}} ... \CStr {T_1} \xto {\CStr {p_1}} \CStr {T_{0}}
\end{equation}
where $\CMsh(p_i)$ is the classifying 1-mesh bundle of $p_i$ (see \autoref{prop:geo_closed_1_truss_bun}).
\end{defn}

\nid The definition extends to maps of closed trusses as follows.

\begin{defn}[Classifying mesh maps of closed $n$-truss maps] Given a map of closed $n$-trusses $F : T \to S$ with components $F_i : T_i \to S_i$, the \textbf{classifying mesh map} $\CMsh F : \CMsh T \to \CMsh S$ is the $n$-mesh map with components given by classifying stratified maps $\CStr {F_i} : \CStr {T_i} \to \CStr {S_i}$.
\end{defn}

\begin{notn}[Classifying mesh functor of closed trusses] The previous definitions assemble into the `classifying mesh functor of closed trusses', denoted by $\CMsh : \ctruss n \to \ctmesh n$, from the category of closed $n$-trusses to the category of $n$-meshes.
\end{notn}

\nid Our earlier construction of entrance path trusses is right inverse to the classifying mesh functor as follows.

\begin{obs}[Entrance path trusses of classifying meshes] \label{rmk:geo_inverts_entr_closed} Unwinding definitions, one checks that the composite $\ETrs \circ \CMsh {T}$ is the identity on the category of closed $n$-trusses (up to a unique natural isomorphism).
\end{obs}

\subsubsecunnum{Classifying meshes of general trusses} We now turn to the construction of classifying meshes of general $n$-trusses. The construction needs care because, as discussed, the `naive' geometric realization of posets (by passing to their nerve, and using the usual geometric realization of simplicial sets) does not dualize dimensions correctly.

    The correct construction of classifying meshes will reduce the case general trusses to that of closed trusses by first `compactifying' general trusses to obtain closed trusses. The basic idea for truss compactifications is a combinatorial analog to fiberwise compactifications of 1-mesh bundles that were defined in \autoref{constr:fiberwise_compactification}. However, when employing fiberwise compactifications `inductively' to a \emph{tower} of bundles, we are left with a choice on how to extend bundles to compactifications of their base. This leads to several possible notions of compactifications. We mention the following universal choice, which we call `cubical compactification' of a truss, that admits a useful explicit construction. Rather than working only with trusses, we will define cubical compactifications in the more general setting of truss bundles.

\begin{term}[Dense subposet and truss subbundles] A subposet $P \into Q$ is called `dense' if the upward closure of $P$ is all of $Q$. A base preserving truss subbundle $F : T \into S$ is similarly called dense if all $i$-level maps $F_i$ are dense (recall, `base preserving' means $F_0 = \id$).
\end{term}

\begin{defn}[Cubical compactification] \label{defn:truss-compactifications} Given an $n$-truss bundle $T$ over a poset $X$, its \textbf{cubical compactification} $\overline T$ (or simply `compactification') is the unique closed $n$-truss bundle over $X$ satisfying the following.
    \begin{itemize}
        \item[-] \emph{Structure}. $T$ is a dense subbundle $\cint : T \into \overline T$ (called \textbf{inclusion map}). Conversely, $\overline T$ `retracts' to $T$ by a truss bundle surjection $\cret : \overline T \to T$ (called \textbf{retraction map}) such that $\cret \circ \cint = \id_T$.
        \item[-] \emph{Universal property}. For any pair of $n$-truss bundle maps $F : T \toot S : G$, where $S$ is a closed $n$-truss, $F$ is a dense subtruss, $G$ a truss bundle surjection, and $G \circ F = \id_T$, there is a unique truss bundle surjection $H : S \epi \overline T$ such that $H \circ F = \cint$ and $\cret \circ H = G$. \qedhere
    \end{itemize}
\end{defn}

\begin{eg}[Cubical compactifications] In \autoref{fig:globular-and-cubical-compactifications} we indicate the inclusion $\cint : T \into \overline T$ of an open $2$-truss $T$ (in black) into its cubical compactification $\overline T$ (extending $T$ by the red structure).
\begin{figure}[ht]
    \centering
    \def\svgwidth{1\columnwidth}
    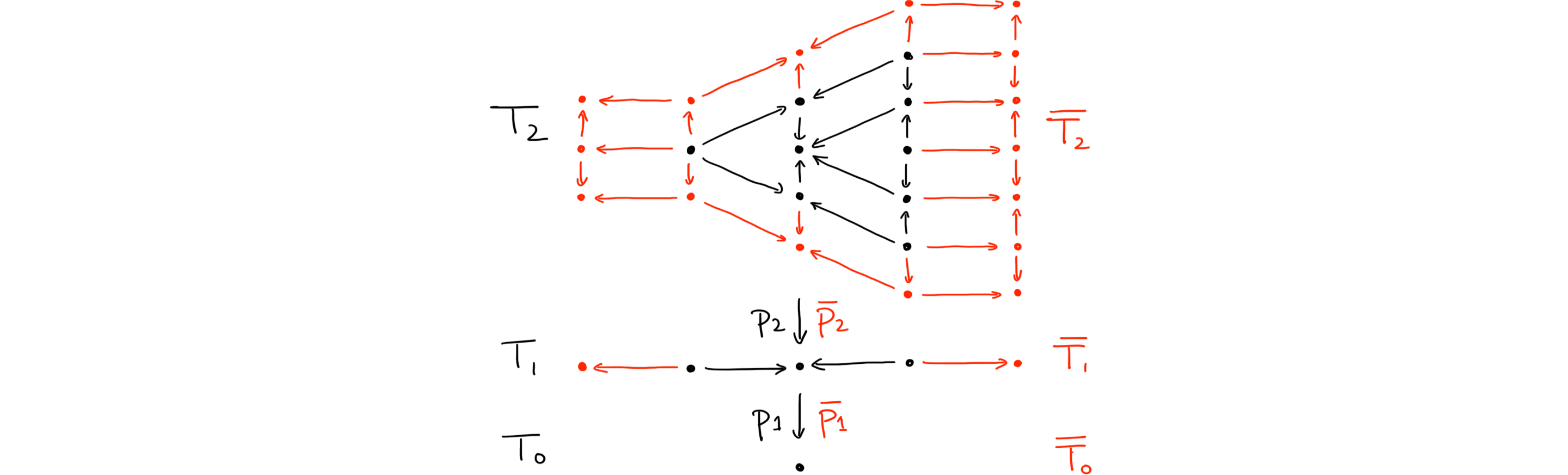

    \caption{Cubical compactifications of two 2-trusses.}
    \label{fig:globular-and-cubical-compactifications}
\end{figure}
\end{eg}

\nid The fact that cubical compactifications always exist is established by the following sequence of constructions, which explicitly construct cubical compactification with increasing generality first for 1-trusses, then for 1-truss bundles, and finally for $n$-truss bundles.

\begin{constr}[Cubical compactifications of 1-trusses]\label{constr:truss-compactifications-1-truss} Given a 1-truss $T$, its cubical compactification $\cint : T \toot \overline T : \cret$ is the unique closed 1-truss containing $T$, which satisfies the conditions in \autoref{defn:truss-compactifications} (in the case $n = 1$ and $X = [0]$): explicitly, $\overline T$ is obtained from $T$ by adjoining a new upper (resp.\ lower) singular endpoint if the upper (resp.\ lower) endpoint of $T$ is regular.
\end{constr}

\begin{constr}[Cubical compactifications of 1-truss bundles]\label{constr:truss-compactifications-1-truss-bun} Given a 1-truss bundle $p : T \to X$, its cubical compactification $\cint : p \toot \overline p : \cret$ is the unique closed 1-truss bundle containing $T$ and satisfying the conditions in \autoref{defn:truss-compactifications} (in the case $n = 1$): explicitly, $\overline p : \overline T \to X$ is obtained from $p$ by compactifying each fiber of $p$ using \autoref{constr:truss-compactifications-1-truss}, and then extending truss bordisms to compactified fibers in the unique endpoint preserving way.
\end{constr}

\begin{constr}[Cubical compactification for $n$-truss bundles] \label{constr:truss-compactifications} Let $T$ be an $n$-truss bundle over a base poset $X$ with bundle maps $p_i : T_i \to T_{i-1}$. Inductively construct its $(n-1)$-truncated compactification $\cint : T_{<n} \toot \overline T_{<n} : \cret$ (the case $n = 1$ is \autoref{constr:truss-compactifications-1-truss-bun}).

    First pull back $p_n : T_n \to T_{n-1}$ along $\cret_{n-1} : \overline T_{n-1} \to T_{n-1}$ to obtain a 1-truss bundle $\cret_{n-1}^*p_n$ over $\overline T_{n-1}$. Pulling back again, now by $\cint_{n-1}$, recovers $p_n$, that is, $\cint_{n-1}^*\cret_{n-1}^*p_n = p_n$. We thus have two 1-truss bundle maps: the 1-truss subbundle $\Totb{}{\cint_{n-1}} : p_n \into \cret_{n-1}^*p_n$ and the 1-truss bundle map $\Totb{}{\cret_{n-1}} : \cret_{n-1}^*p_n \to p_n$. The two 1-truss bundle maps form an inclusion-retraction pair $p_n \toot \cret_{n-1}^*p_n$. Composing this pair with the inclusion-retraction pair $\cret_{n-1}^*p_n \toot \overline {\cret_{n-1}^*p_n}$ constructed using \autoref{constr:truss-compactifications-1-truss-bun}, then defines the inclusion-retraction pair $\cint_n : p_n \toot \overline p_n : \cret_n$ (where $\overline p_n := \overline {\cret_{n-1}^*p_n}$).

    We define the cubical compactification $\overline T$ to be the truss bundle obtained by augmenting $\overline T_{<n}$ with the bundle map $\overline p_n$, and obtain an inclusion-retraction pair $\cint : T \toot \overline T : \cret$ as required.
\end{constr}

\nid We now give a definition of classifying meshes of general trusses. Recall the notion of `submeshes' of $n$-meshes from \autoref{term:submeshes} (which are mesh inclusions the $n$-framing of the codomain restricts to the $n$-framing of the domain). We will also use the notion of `constructible substratifications' (see \autoref{defn:constr_substrat}): a substratification $F : (Y,g) \into (X,f)$ is constructible if $Y = f\inv(Q)$ for some subposet $Q \into \Entr(f)$ of the entrance path poset of $f$. Constructible substratifications $F$ are thus fully determined by their entrance path poset maps $\Entr F$.

\begin{defn}[Classifying meshes of general $n$-trusses] \label{constr:geo_open} For an $n$-truss $T$, the \textbf{classifying $n$-mesh} $\CMsh T$ of $T$ is the submesh $\CMsh {\cint} : \CMsh T \into \CMsh \overline T$ given by constructible substratifications $\CMsh {\cint_i} : \CMsh T_i \into \CMsh \overline {T}_i$ with entrance path poset map $\Entr \CMsh {\cint_i} = \cint_i : T_i \into \overline {T}_i$.
\end{defn}

\nid (The verification that $\CMsh T$ indeed defines an $n$-mesh proceeds inductively, using the inductive construction of cubical compactifications.)

Note that if $T$ is closed, then the preceding definition of $\CMsh T$ specializes to our earlier \autoref{constr:geo_closed}.

\begin{eg}[Classifying meshes of general $n$-trusses] Recall the open $2$-truss $T$ from \autoref{fig:globular-and-cubical-compactifications}. In \autoref{fig:classifying-meshes-of-general-n-trusses} we depict the closed classifying mesh $\CMsh \overline T$ of the compactification $\overline T$ of $T$, together with the resulting open classifying mesh $\CMsh T$ of $T$. The latter mesh includes into $\CMsh \overline T$ as a submesh $\CMsh \cint : \CMsh T \into \CMsh \overline T$ as indicated.
\begin{figure}[ht]
    \centering
    \def\svgwidth{1\columnwidth}
    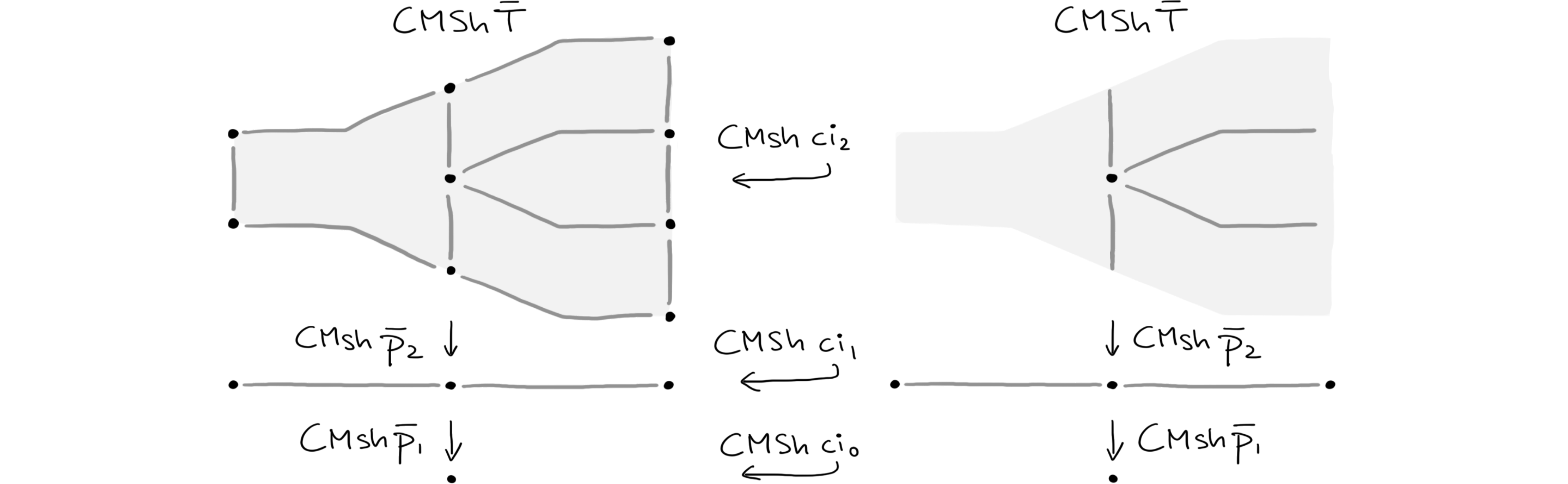

    \caption[The classifying mesh of an open $2$-truss.]{The classifying mesh of an open $2$-truss constructed via the classifying mesh of its compactification.}
    \label{fig:classifying-meshes-of-general-n-trusses}
\end{figure}
\end{eg}

\begin{obs}[Geometric and simplicial geometric realization of trusses] For an $n$-truss $T$, unwinding definitions, one checks that $\ETrs \CMsh T = T$.
\end{obs}

\nid Given an $n$-truss $T$, then classifying stratifications $\CStr{T_k}$ and classifying meshes $\CMsh{T_k}$ are in general distinct stratifications (the simplest example being obtained by taking $T$ to be the open 1-truss with one element). However, the former stratification canonically includes into the latter stratification as follows.

\begin{rmk}[Classifying meshes vs. classifying stratifications]  \label{obs:geo_vs_simp_geo_open} For an $n$-truss $T$, the classifying stratified map $\CStr{\cint_k} : \CStr{T_k} \into \CStr{\overline T_k} = \CMsh {\overline T_k}$ is a substratification that factors through the constructible substratification $\CMsh {T_k} \into \CMsh {\overline T_k}$ by a substratification $\CStr{T_k} \into \CMsh T_k$. This satisfies that
\begin{equation}
\CStr {T_k} \into \CMsh {T_k} \xto {\CMsh {\cint}} \CMsh {\overline T_k} = \CStr{\overline T_k} \xto {\CStr {\cret_k}} \CStr {T_k}
\end{equation}
composes to the identity.
\end{rmk}

\nid It remains to define classifying mesh maps for general truss maps.

\begin{defn}[Geometric realization of maps of closed trusses] \label{constr:geo_maps_open} Given a map of $n$-trusses  $F : T \to S$ with components $F_i : T_i \to S_i$, its \textbf{classifying mesh map} $\CMsh F : \CMsh T \to \CMsh S$ is the $n$-mesh map with components $\CMsh {F_i}$ defined by the composite
\begin{equation}
\CMsh {T_k} \xto {\CMsh {\cint}} \CMsh {\overline T_k} \xto {\CStr {\cret_k}} \CStr {T_k} \xto {\CStr {F_k}} \CStr{S_k} \into \CMsh {S_k}
\end{equation}
where the last map is the canonical inclusion from \autoref{obs:geo_vs_simp_geo_open} (the observation also implies that the association $F \mapsto \CMsh F$ is functorial).
\end{defn}

\nid Note that $\CMsh F$ is a singular resp.\ regular $n$-mesh map whenever $F$ is a singular resp.\ regular $n$-truss map.

\begin{notn}[Classifying mesh functor] The previous definitions assemble into the `classifying mesh functor', denoted by $\CMsh : \truss n \to \tmesh n$, from the category of trusses to the category of $n$-meshes.
\end{notn}

\begin{obs}[Geometric realization is inverse to entrance path truss construction] \label{rmk:geo_inverts_entr_comaps} Unwinding definitions, one checks that $\ETrs \circ \CMsh - \iso \id_{\truss n}$ by a unique natural isomorphism, and thus $\CMsh$ is `right inverse' to $\ETrs$.
\end{obs}

\subsubsecunnum{Classifying mesh maps of coarsenings} Finally, we discuss a more specific construction of classifying mesh maps in the case of coarsening. Note that the construction will not be needed for the proofs of the main theorems in this and later chapters (but it will be useful for several other results). The construction is motivated by the following observation.

\begin{obs}[Classifying mesh maps of degeneracies and coarsenings] \label{obs:class-crs-asymmetry} Given an $n$-truss degeneracy map $F : T \to S$, then the classifying mesh map $\CMsh F : \CMsh T \to \CMsh S$ is in fact a mesh degeneracy, and thus the classifying mesh functor restricts to a functor $\CMsh : \degtruss n \to \degmesh n$, which provides a `right inverse' to the functor $\ETrs \degmesh n \to \degtruss n$.

    However, if $F$ is an $n$-truss coarsening, then the analogous observation fails: that is, $\CMsh F$ need not be an $n$-mesh coarsening.
\end{obs}

\nid Addressing this asymmetry, in this section we will construct a `homotopical replacement' $\CrsMsh F$ of $\CMsh F$ (which is an $n$-mesh coarsening and homotopic to $\CMsh F$). We call $\CrsMsh F$ the `classifying mesh coarsening' of the $n$-truss coarsening $F$. We start with the construction of classifying mesh coarsenings in the case of closed trusses.

\begin{constr}[Classifying mesh coarsenings of closed truss coarsenings]
\label{constr:PL-realization-coarsening}
Given closed $n$-trusses $T = (p_n,...,p_1)$ and $S = (q_n,...,q_1)$ and a coarsening $F : T \to S$, we construct an $n$-mesh coarsening $\CrsMsh F : \CMsh T \to \CMsh S$ with the following two properties. Firstly, components $\CrsMsh F_i : \CMsh{T_i} \to \CMsh{S_i}$ are \emph{linear} on each simplex of $\CMsh{T_i} = \CStr {T_i}$. Secondly, we have $\ETrs \CrsMsh F = F$.

Arguing inductively in $n$, assume to have constructed $\CrsMsh F_{<n} : \CMsh {T_{<n}} \to \CMsh {S_{<n}}$ with the claimed properties. We define the top component $\CrsMsh F_n : \CMsh {T_n} \to \CMsh {S_n}$ on vertices $x \in T_n$, by picking any image point $ \CrsMsh F_n(x)$ in the fiber of the 1-truss bundle $\CMsh(q_n)$ over $\CrsMsh F_{n-1} \circ \CMsh(p_n) (x)$ subject to the following conditions: firstly, $\CrsMsh F_n(x) < \CrsMsh F_n(x')$ (in the linear order of the fiber) whenever $x \fles x'$ (in the frame order of $T_n$); secondly, $\CrsMsh F_n(x)$ lies in the stratum $F_n(x) \in S_n$ of $\CMsh {S_n}$. We then extend $\CrsMsh F_n$ linearly to all other simplices in $\CMsh {T_n} = \CStr {T_n}$. One checks that this makes $\CrsMsh F_n$ a coarsening, and thus $\CrsMsh F$ an $n$-mesh coarsening as claimed.
\end{constr}

\nid We next generalize our construction of classifying mesh coarsenings to the case of general trusses. The generalization is based on the observation that `cubical compactification is functorial on coarsenings of trusses' in the following sense.

\begin{obs}[Cubical compactification is functorial on coarsenings] \label{obs:cubical-compactification-functorial} Given a coarsening of $n$-trusses $F : T \to S$, there is a unique `cubically compactified coarsening' $\overline F : \overline T \to \overline S$ between the respective cubical compactifications of $T$ and $S$, such that the following two squares commute
\begin{equation}
\begin{tikzcd}[column sep=50pt, baseline=(W.base)]
\overline T \arrow[r, "\overline F"] \arrow[d, "\cret"' anchor=east, bend right] & \overline S \arrow[d, "\cret"' anchor=east, bend right, pos=.6] \\
T \arrow[r, "F"'] \arrow[u, "\cint"' anchor=west, bend right]           & |[alias=W]| S \arrow[u, "\cint"' anchor=west, bend right]
\end{tikzcd}.
\end{equation}
(The construction of $\overline F$ proceeds inductively, similar to the sequence of constructions defining cubical compactifications.)
\end{obs}

\begin{constr}[Classifying mesh coarsenings of truss coarsenings] \label{constr:PL-realization-coarsening-general} Given a coarsening of $n$-trusses $F : T \to S$, we construct an $n$-mesh coarsening $\CrsMsh F : \CMsh T \to \CMsh S$ with the following properties. Firstly, components $\CrsMsh F_i : \CMsh{T_i} \to \CMsh{S_i}$ are coarsenings which are linear each open simplex in $\CMsh{T_i} \into \CMsh{\overline T_i}$. Secondly, $\ETrs \CrsMsh F = F$.

    We may apply \autoref{constr:PL-realization-coarsening} and \autoref{obs:cubical-compactification-functorial} to construct a coarsening $\CrsMsh {\overline F} : \CMsh {\overline T} \to \CMsh {\overline S}$. This coarsening descends to a coarsening $\CrsMsh F : \CMsh  T \to \CMsh S$ along the submeshes $\CMsh {\cint} :  \CMsh T \into \CMsh {\overline T}$ and $\CMsh {\cint} : \CMsh S \into \CMsh {\overline S}$. One checks that this gives an $n$-mesh coarsening with the claimed properties.
\end{constr}

\nid Note that by \autoref{prop:faithfulness-crs-deg}, since $\ETrs \CrsMsh F = F = \ETrs \CMsh F$, we find that the classifying mesh coarsening $\CrsMsh F$ and the classifying mesh map $\CMsh F$ are in fact homotopic mesh maps, and in this sense $\CrsMsh F$ is a `homotopical replacement' of $\CMsh F$, as claimed earlier.

\subsection{The theorem and its applications} \label{sec:constructibility-of-meshes}

We finally prove \autoref{thm:diagram_classification}, which claimed that the classifying mesh functor and the entrance path truss functor provide a weak equivalence pair between (certain) topological categories of meshes and trusses.

\begin{proof}[Proof of \autoref{thm:diagram_classification}] Recall, we set out to show that the topological functors $\ETrs : \sctmesh {n} \to \sctruss {n}$ and $\ETrs : \rotmesh {n} \to \rotruss {n}$ are weak equivalences of topological categories. We argue in the closed case (the argument for the open case is fully analogous). We need to check the following (see \cite[Def. 1.1.3.6]{lurie2009higher}).
\begin{enumerate}
\item For each closed $n$-truss $T \in \ctruss {n}$ there exists a closed $n$-mesh $f \in \sctmesh {n}$ whose entrance path truss $\ETrs f$ is weakly equivalent to $T$.
\item The mapping of hom spaces $\ETrs : \sctmesh {n}(f,g) \to \sctruss {n}(\ETrs f, \ETrs g)$ is a weak equivalence of topological spaces.
\end{enumerate}
The first statement is a consequence of the construction of classifying meshes $\CMsh T$ (see \autoref{constr:geo_closed}) which satisfies that $\ETrs \CMsh T = T$ (see \autoref{rmk:geo_inverts_entr_comaps}). The second statement follows from the weak faithfulness of the entrance path truss functor (see \autoref{prop:faithfulness}) together with the observation that fibers of $\ETrs$ are never empty (since all elements of $\sctruss {n}(T,S)$ can be geometrically realized by \autoref{constr:geo_maps_open}).
\end{proof}

\noindent In particular, it follows that the topological category $\sctmesh n$ (respectively $\rotmesh n$) is 1-truncated (and thus it can be thought of as an ordinary 1-category, cf. \cite[Prop. 2.3.4.18]{lurie2009higher}). We next discuss applications of the theorem.

\subsubsecunnum{Meshes of framed cells and their subdivisions}

As a first application of the equivalence of meshes and trusses, we describe the equivalence of meshes and framed regular cells. Namely, flat framed regular cell complexes realize as so-called `cell meshes', and conversely, closed meshes translate to framed regular cell complexes called `mesh complexes'.

\begin{term}[Cell meshes] The composite of the truss translation functor $\kT : \FrCDiag n \to \sctruss n$ with classifying mesh functor $\CMsh : \sctruss n \to \sctmesh n$ will be denoted by $\CellMsh : \FrCDiag n \to \sctmesh n$ and called the `cell mesh functor'.
\end{term}

\begin{term}[Mesh cells] Meshes in the image of cell mesh functor $\CellMsh$ restricted to the subcategory of $n$-framed regular cells will be called `$n$-mesh cells'.
\end{term}

\nid Conversely to cell meshes, we define mesh complexes as follows.

\begin{term}[Mesh complexes] \label{term:mesh-complex} The composite of the entrance path truss functor $\ETrs : \scmesh n \to \sctruss n$ with the framed complex translation functor $\kX : \sctruss n \to \FrCDiag n$ will be denoted by $\MshCplx : \scmesh n \to \FrCDiag n$ and called the `mesh complex functor'.
\end{term}

As recorded in \autoref{cor:mesh-complex-cell-mesh-equiv} these functors (regarded as functors of topological categories) yield weak equivalences $\MshCplx : \sctmesh n \eqv \FrCDiag n : \CellMsh$. Restricting the equivalence further to subcategories of framed regular cells resp.\ mesh cell, we also record the following.

\begin{cor}[Equivalence of flat framed cells and mesh cells] \label{cor:eqv-mesh-cells-and-framed-reg-cells} The cell mesh and mesh complex functors establish a weak equivalence between $n$-framed regular cells and $n$-mesh cells.
\end{cor}

\pause

The translation of framed regular cells into mesh cells provides a `geometric realization' of framed regular cell as framed stratified spaces. We can leverage that stratified spaces come with a natural notion of `refinement map' (while cellular posets do not), to turn this into a definition of `framed subdivisions' of framed regular cells. We start with unframed case of (combinatorial) regular cell complexes.

\begin{term}[Subdivision of regular cells] Let $X$ be a combinatorial regular cell (i.e.\ a cellular poset with initial object), and $Y$ a combinatorial regular cell complex. A `subdivision' of $X$ by $Y$ is a stratified coarsening $F : \CStr Y \to \CStr X$ of between their respective classifying stratifications. Given $y \in Y$, the stratified map $F$ restricts on the cell $\CStr {(Y^{\geq y})} \into \CStr Y$ to the `cell inclusion' $\rest F y : \CStr {(Y^{\geq y})} \into \CStr X$.\footnote{Note that the cell inclusions $\rest F y$ generally are non-cellular maps (see \autoref{defn:cellular-maps-of-rcc}), and, in this sense, they are not immediately combinatorializable (see \autoref{rmk:cellular-map-yield-eqv}).}
\end{term}

\nid If $(X,\cG)$ is a \emph{framed} regular cell, then its classifying stratification $\CStr X$ is the underlying space of the cell mesh $\CellMsh (X,\cG)$, and this allows us to keep track of framing structures: a `framed subdivision' is a subdivision that restricts on each cell to a framed map, as follows.

\begin{defn}[Framed subdivisions of framed regular cells] \label{defn:framed-subdiv} For an $n$-framed regular cell $(X,\cF)$ and an $n$-framed regular cell complex $(Y,\cG)$, a \textbf{framed subdivision} $F : \CStr Y \to \CStr X$ is a subdivision such that, for each $y \in Y$, the cell inclusion $\rest F y : \CellMsh (Y^{\geq y},\rest \cG y) \into \CellMsh (X,\cF)$ is (the top component of a) framed map of meshes.
\end{defn}

\begin{eg}[Framed subdivisions of framed regular cells] In \autoref{fig:framed-subdivisions-of-framed-regular-cells} we illustrate the framed subdivision $F$ of a framed regular cell by a framed regular cell complex, as well as the condition that $F$ restrict on each framed regular cell to the top component of a framed map of meshes.
\begin{figure}[ht]
    \centering
    \def\svgwidth{1\columnwidth}
    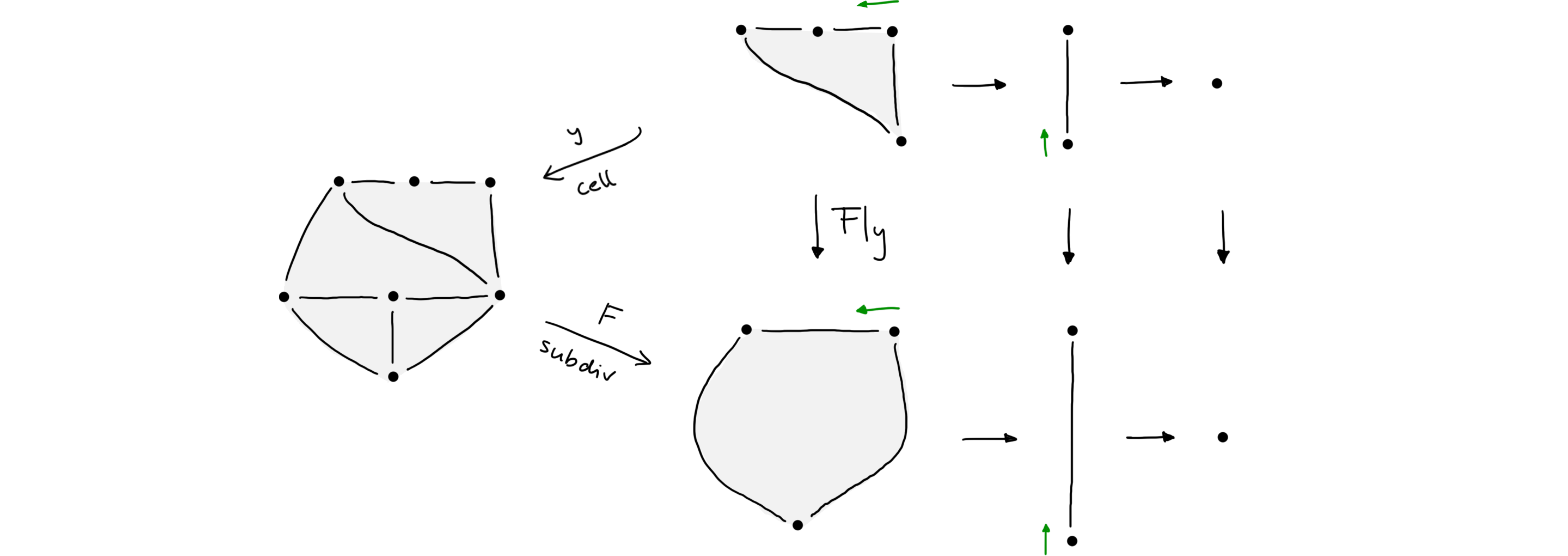

    \caption{Framed subdivisions of framed regular cells}
    \label{fig:framed-subdivisions-of-framed-regular-cells}
\end{figure}
\end{eg}

\begin{lem}[Framed subdivisions are flat] Given a framed subdivision $F : (Y,\cG) \epi (X,\cF)$ of an $n$-framed regular cell $(X,\cF)$ by an $n$-framed regular cell complex $(Y,\cF)$, then $(Y,\cG)$ is a flat.
\end{lem}

\nid The proof follows the description of framed regular cells as `section' and `spacer' cell from \autoref{ssec:frames-vs-proj-frames-reg-cell}, and we only provide an outline.

\begin{proof} We argue inductively in $n$. First, assume $(X,\cF)$ is a section cell. Then all cells in $(Y,\cG)$ must be section cells, and $F$ induces a framed subdivision $F : (Y,\cG_{n-1}) \to (X,\cF_{n-1})$ where $(X,\cF_{n-1})$ is the projected cell of $(X,\cF)$, an $(Y,\cG)_{n-1}$ is the complex of projected cells $(Y^{\geq y},(\rest \cG x)_{n-1})$ (see \autoref{term:projected-cells}). By induction, this makes $(Y,\cG_{n-1})$ a flat framed cell complex, which implies $(Y,\cG)$ too is a flat framed cell complex.

Next, assume $(X,\cF)$ is a spacer cell and denote its lower section cell by $\partial_- X$. Denote by $\partial_- Y \subset Y$ the preimage of $\partial_- X$ under $F$ on entrance path posets (that is, under the map $\Entr F : Y \to X$). Note that $F$ restricts to a framed subdivision $F : \CStr {\partial_- Y} \to \CStr {\partial_- X}$ (of the section cell $(\partial_- X,\rest \cF {\partial_- X})$ by the framed complex $(\partial_- Y,\rest \cG {\partial_- Y})$). Consider the cell projection $p_n : \CStr X \to \CStr X_{n-1}$ of $X$ to its projected cell $X_{n-1}$ (note $X_{n-1} \iso \partial_- X$, see \autoref{rmk:section-projections}). The top 1-mesh bundle in the $n$-mesh $\CellMsh(X,\cF)$ endows $p_n$ with the structure of an 1-mesh bundle. We may refine both the domain and codomain of $p_n$ by $F$ to obtain another 1-mesh bundle $q_n$ (with the same fiber 1-framing as $p_n$) as shown:
\begin{equation}
\begin{tikzcd}
\CStr Y \arrow[r, "F"] \arrow[d, "q_n"', dashed] & \CStr X \arrow[d, "p_n"] \\
\CStr \partial_- Y \arrow[r, "F"]                & \CStr \partial_- X
\end{tikzcd}
\end{equation}
(The construction of $q_n$ uses that spacer and section cells must be alternatingly `layered' over cells in the lower section $\CStr \partial_- Y \into \CStr Y$). Augmenting the inductively constructed $(n-1)$-mesh $\CellMsh(\partial_- Y, (\rest \cG {\partial_- Y})_{n-1})$ by the 1-mesh bundle $q_n$ yields an $n$-mesh, which we denote by $M$: by construction this satisfies $\MshCplx M \iso (Y,\cG)$ and thus $(Y,\cG)$ is a flat framed regular cell complex as claimed.
\end{proof}

\nid Having established that all subdividing framed regular cell complexes $(Y,\cG)$ must be flat (which enables us to construct the cell mesh $\CellMsh(Y,\cG)$), we may now equivalently characterize framed subdivisions as follows.

\begin{cor}[Framed subdivision are mesh coarsenings] \label{cor:fr-subdiv-are-mesh-crs} For an $n$-framed regular cell $(X,\cF)$ and cell complex $(Y,\cG)$, a stratified coarsening $F : \CStr Y \to \CStr X$ is a framed subdivision if and only if $F$ is a $n$-mesh coarsening $F : \CellMsh (Y,\cG) \to \CellMsh (X,\cF)$. \qed
\end{cor}

\begin{term}[Space of subdivisions] For an $n$-framed regular cell $(X,\cF)$ and cell complex $(Y,\cG)$, the space of `framed subdivision' $\mathsf{SubDiv}(Y,\cG;X,\cF)$ is the space of mesh coarsenings $\crstmesh n (\CellMsh (Y,\cG),\CellMsh (X,\cF))$.
\end{term}

\pauseae

\nid The notion of framed subdivisions of framed regular cells can be `combinatorialized'.

\begin{notn}[Set of truss coarsening] For an $n$-block $B$ and a closed $n$-truss $T$, denote the set of truss coarsening $T \to B$ by $\crstruss n (T,B)$.
\end{notn}

\begin{cor}[Framed subdivisions are up to homotopy truss coarsenings] \label{cor:classifying-subdiv-2} The map $\ETrs : \mathsf{SubDiv}(Y,\cG;X,\cF) \to \crstruss n (\kT(Y,\cG),\kT(X,\cF))$ is a weak homotopy equivalence.
\end{cor}

\nid In other words, up to contractible choice in $\mathsf{SubDiv}(-;X,\cF)$, framed subdivisions of $(X,\cF)$ correspond exactly to truss refinements of $\kT(X,\cF)$.

\begin{proof} We need to show that entrance path truss functor $\ETrs$, mapping framed subdivisions $\mathsf{SubDiv}(Y,\cG;X,\cF)$ to truss coarsenings $\truss n (\kT(Y,\cG),\kT(X,\cF))$, has weakly contractible preimages. This follows from \autoref{prop:faithfulness-crs-deg}.
\end{proof}

\begin{rmk}[Framed subdivisions can be made piecewise linear] Combining the corollary with our earlier construction of `piecewise linear' classifying mesh coarsenings (see \autoref{constr:PL-realization-coarsening-general}), it follows that any framed subdivision is homotopic to a piecewise linear framed subdivision.
\end{rmk}

\nid The corollary not only provides a combinatorialization (up to contractible choice) of framed subdivisions, but the resulting combinatorialization is also computationally tractable: namely, since truss refinements of $\kT(X,\cF)$ can be algorithmically listed, so can (up to contractible choice) framed subdivisions of the framed regular cell $(X,\cF)$. This stands in further contrast to the classical observation that it is \emph{impossible} to algorithmically decide when a simplicial complex is piecewise linearly homeomorphic to the simplex.

Finally, the fact that framed subdivisions can combinatorialized in terms of morphisms in the category of trusses (living in the same category as faces and degeneracies), is by itself rather remarkable property of framed regular cells. In contrast, for most classes of combinatorial shapes, notions subdivisions cannot be given in terms of combinatorial `morphism' (for instance, in the case of simplices, a combinatorial description of subdivision requires rather different techniques than the definition of simplicial maps, see \cite{pachner1991pl}).

\subsubsecunnum{Dualization of meshes}

As a second, and central, application of \autoref{thm:diagram_classification} we now construct the dualization functors between the categories of closed meshes with singular maps and open meshes with regular maps. This proves our earlier \autoref{thm:dualization_functors}.

\begin{proof}[Proof of \autoref{thm:dualization_functors}] The `mesh dualization'  functors
\begin{equation}
\dagger : \sctmesh n \eqv \rotmesh n : \dagger
\end{equation}
are defined by the respective composites in
\begin{equation}
\begin{tikzcd}
\sctmesh n \arrow[r, "\ETrs", shift left] & \sctruss n \arrow[r, "\dagger", shift left] \arrow[l, "\CMsh", shift left] & \rotruss n \arrow[l, "\dagger", shift left] \arrow[r, "\CMsh", shift left] & \rotmesh n \arrow[l, "\ETrs", shift left]
\end{tikzcd}
\end{equation}
where the central arrows (labeled by `$\dagger$') are the dualization functors of trusses, see \autoref{thm:labeled-n-truss-bundle-dualization}. Since each functor in the above composite is an equivalence, so are the mesh dualization functors.
\end{proof}

\begin{eg}[Meshes and their duals] In \autoref{fig:meshes-and-their-duals} we depict meshes $M$ together with their duals $M^\dagger$ (note in particular, this turns open meshes into closed meshes and vice versa). Note that the depicted closed 3-mesh contains \emph{two} 3-cells (corresponding to the two 0-cells in the depicted open 3-mesh); this is a subcomplex of the gluing depicted earlier in \autoref{fig:the-dual-of-the-twisted-embedding-of-the-circle}.
\begin{figure}[ht]
    \centering
    \def\svgwidth{1\columnwidth}
    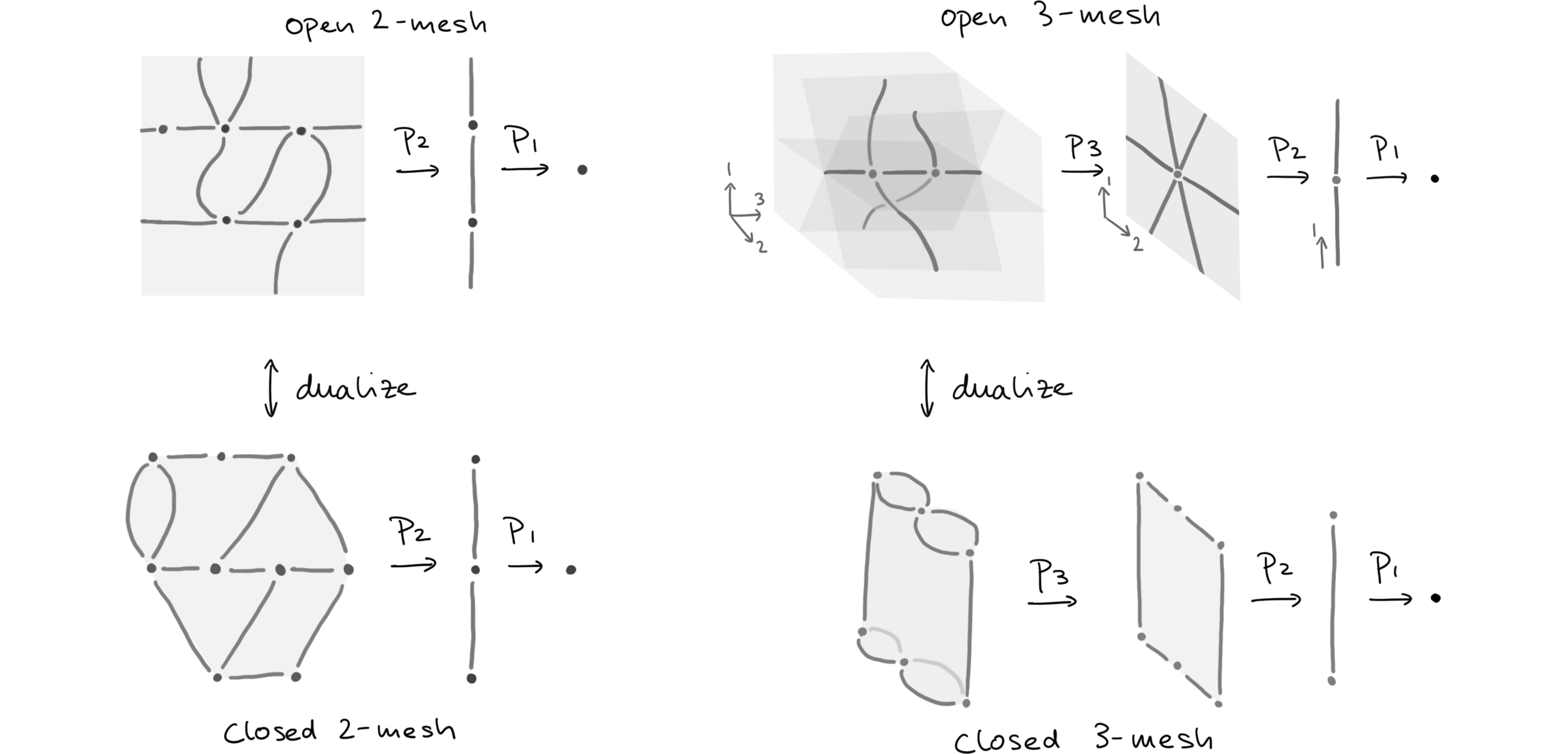

    \caption{Meshes and their duals.}
    \label{fig:meshes-and-their-duals}
\end{figure}
\end{eg}

\begin{rmk}[Dualization of flat framed regular cell complexes] \label{rmk:dualization-of-blocks} While flat framed regular cell complexes don't dualize to flat framed regular cell complexes (they correspond to closed meshes, which to open meshes), they do so `up to compactification' (see \autoref{defn:truss-compactifications}). Namely, given a flat $n$-framed regular cell complex $(X,\cF)$ define its `dual' to be the flat $n$-framed regular cell complex $\kX \overline {\kT(X,\cF)^\dagger}$ (that is, the framed complex translation of the compactification of the dual of $\kT(X,\cF)$). This is the `duality' used in an earlier example: namely, the framed complex in \autoref{fig:the-dual-of-the-twisted-embedding-of-the-circle} `dualizes' to the framed complex in \autoref{fig:triangulating-the-twisted-embedding-of-the-circle}.
\end{rmk}

\nid The fact that meshes admit a notion of dualization is fundamental, and provides the key to many ideas in the interplay of framed combinatorial topology and geometric higher category theory (such as those outlined in \autoref{intro:outlook}).

\chapter{Combinatorializability of flat framed stratifications} \label{ch:hauptvermutung}

%
%


In the final chapter of this book we will begin our study of framed combinatorial \emph{stratified} topology, or more precisely, of its model stratifications given by so-called \emph{flat} framed stratifications. Flat $n$-framed stratifications provide a general class of stratifications defined by the property of being `meshable', i.e.\ admitting a framed refinement by some $n$-mesh. Their `flat $n$-framing' is induced by a framed embedding into standard $n$-framed $\lR^n$. Flat framed stratifications provide `local models' for more general, `global' framed stratified spaces (much as framed regular cells provide local models for framed regular cell complexes).  While we will not discuss `global' framed stratifications here, the basic machinery of flat framed stratifications developed in this chapter will be tailored towards future applications in that and other related directions.

We outline the main results of this chapter. Recall, a stratified map is a `refinement of its image' if its underlying map of spaces is a homeomorphism. The definition of flat framed stratifications can be formally given as follows.

\begin{defn}[Flat framed stratifications] \label{defn:flat-framed-strat}
A \textbf{flat $n$-framed stratification} is a stratification $(Z,f)$ of a subspace $Z \subset \lR^n$, for which there exists an $n$-mesh whose $M$ flat $n$-framing (as a stratified map into $M_n \into \lR^n$) refines $(Z,f)$.
\end{defn}

\begin{figure}[ht]
    \centering
    \def\svgwidth{1\columnwidth}
    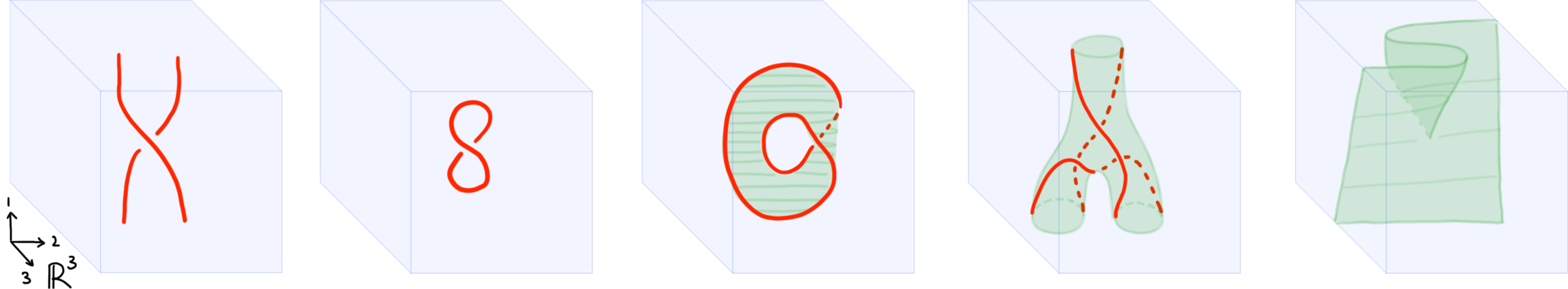

    \caption[Examples of flat 3-framed stratifications.]{Examples of flat 3-framed stratifications of the open 3-cube representing (from left to right) a `braid', a `Hopf circle', a `Mobius band', a `pair of pants with defects', and a `cusp singularity'.}
    \label{fig:flat-framed-stratifications}
\end{figure}

\nid In \autoref{fig:flat-framed-stratifications}, we illustrate several examples of flat $n$-framed stratifications $(I^3,f)$ of the open 3-cube $I^3 \subset \lR^3$. All examples can be seen to admit a refinement by the flat 3-framing of an open 3-mesh. (Note that, in general, the underlying space $Z \subset \lR^n$ of a flat framed stratification $(Z,f)$ need not be open, as it can be the image of the $n$-framing of any mesh.)


We will prove that flat framed stratifications, despite their definition in stratified topological terms, have an unexpected combinatorial counterpart. While we've shown in the previous chapter that meshes are `weakly' equivalent to constructible combinatorial objects (namely, to trusses), the above definition of flat framed stratifications makes reference merely to the \emph{existence} of some mesh refinement without singling out a particular such refinement. It therefore remains entirely unclear that flat framed stratifications can themselves be combinatorialized as claimed.

For the combinatorialization of meshes to transfer to that of flat framed stratifications, we will have to understand the class of all refining meshes in combinatorial terms. This will be achieved by the notion of `stratified trusses'. The notion is a specialization of our earlier definition of `labeled $n$-trusses' (see \autoref{defn:n-trusses-labeled}) requiring the labeling to be the characteristic map of a stratification (here, we regard posets as topological spaces via their specialization topology, see \autoref{conv:spec-top}). 

\begin{defn}[Stratified $n$-truss] \label{defn:norm-strat-trusses} A \textbf{stratified $n$-truss} $T$ is a labeled $n$-truss $T$ whose labeling $\lbl T$ is the characteristic map of a stratification on the total poset $T_n$ of $T$.
\end{defn}

\nid The correspondence between $n$-trusses and $n$-meshes can now be extended to a correspondence between \emph{stratified} $n$-trusses and \emph{refinements} of flat framed stratifications by $n$-meshes (up to `framed stratified homeomorphism'). Namely,  given a refinement of a flat framed stratification by an $n$-mesh, we can first consider the $n$-truss corresponding to the $n$-mesh, and then produce a labeling of that $n$-truss by precomposing the characteristic map of the stratification with the refinement's entrance path poset map. To obtain \emph{unique} combinatorial representatives of flat framed stratifications, we will be interested in the following subclasses of stratified trusses, consisting of those stratified trusses that do not admit a non-identity `label preserving truss coarsening' (see \autoref{term:label-and-base-preserving} and \autoref{term:coarsening-and-quotients}).

\begin{defn}[Normalized stratified $n$-trusses] \label{defn:norm-classes} We say a stratified $n$-truss $T$ is \textbf{normalized} if any label preserving truss coarsening of $T$ is an identity.
\end{defn}

Our main theorem in this chapter will show that flat framed stratifications are indeed `combinatorializable', by producing a correspondence of flat framed stratifications (up to framed stratified homeomorphism) with normalized stratified trusses (up to `stratified truss isomorphisms'), as follows.

\begin{thm}[Normalized stratified trusses classify flat framed stratifications] \label{thm:flat-fr-strat-are-norm-strat-trusses} Framed stratified homeomorphism classes of flat $n$-framed stratifications are in correspondence with isomorphisms classes of normalized stratified $n$-trusses.
\end{thm}

\nid The theorem's proof will be based on the following observation: while, a priori, flat framed stratifications may have many refining meshes, there is, in fact, always one canonical (namely, `coarsest') refining mesh. The correspondence of the theorem will take a given flat $n$-framed stratification to the stratified $n$-truss representing the coarsest $n$-mesh refinement of that stratification. The proof will therefore hinge on the construction of coarsest refining meshes which, in turn, will be based on a key property of meshes: the existence of so-called `mesh joins'. As suggested by the name, mesh joins are finest mutual coarsenings of any two given meshes (with identical support). The proof of the existence of mesh joins and the construction of coarsest refining meshes will take up all of \autoref{sec:join-stability}. This will provide the technical groundwork for \autoref{ssec:combinatorializability-of-fl-fr-str} where we will spell out the proof of \autoref{thm:flat-fr-strat-are-norm-strat-trusses}.

That theorem has several key corollaries, as we now outline. Most immediately, the theorem implies that (up to homotopically unique homeomorphism) any flat framed stratification may be replaced by a canonical piecewise linear one; here, a `piecewise linear stratification' means a stratification that admits a triangulation by linear simplices in $\lR^n$.

\begin{cor}[Flat framed stratification are piecewise linear] \label{cor:flat-framed-strat-are-PL} Any flat framed stratification is framed stratified homeomorphic to a canonical piecewise linear flat framed stratification.
\end{cor}

\nid The canonical piecewise linear flat framed stratification can be obtained by constructing the `classifying stratification' of the normalized stratified trusses corresponding to the given flat framed stratification. This corollary, in particular, entails that flat framed stratifications are better behaved than mere topological stratifications: indeed, there are bounded stratifications in $\lR^n$ that are not topologically (nonframed) homeomorphism to any piecewise linear stratifications (for instance, take any embedding of the $E_8$ manifold in euclidean space).

Not only are flat framed stratifications framed stratified homeomorphic to piecewise linear ones, but the notions of `framed stratified homeomorphism' and `framed stratified PL homeomorphism' coincide. 

\begin{cor}[Flat framed Hauptvermutung] \label{cor:flat-framed-hpt-vmtg} Any framed stratified homeomorphism between piecewise linear flat framed stratifications is homotopic to some framed stratified piecewise linear homeomorphism.
\end{cor}

\nid In fact, the homotopy is `unique up to contractible choice'. As we will discuss, the classical nonframed analog of the statement fails to hold: namely, there are bounded piecewise linear stratifications in $\lR^n$ that are topologically stratified homeomorphic but not piecewise linearly stratified homeomorphic.

Finally, we also have the following `converse' to \autoref{cor:flat-framed-strat-are-PL}.

\begin{prop}[Piecewise linear stratifications are flat framed] \label{prop:PL-strat-are-flat-frmd} Any bounded piecewise linear stratification in $\lR^n$ is (up to adding an `ambient' stratum) a flat framed stratification.
\end{prop}

\nid The previous three results will be content of \autoref{ssec:polyhedrality-of-fl-fr-str}. Together, they establish the `polyhedrality' of flat framed stratifications, and they can be summarized as follows.

\begin{rmk}[Framed, framed piecewise linear, and piecewise linear categories] \label{rmk:polyhedrality-summarized} \autoref{cor:flat-framed-strat-are-PL} shows that the functor from `piecewise linear flat framed stratifications' into `flat framed stratifications' is surjective on homeomorphism classes, and \autoref{cor:flat-framed-hpt-vmtg} shows that it is injective on homeomorphism classes. Further, \autoref{prop:PL-strat-are-flat-frmd} shows that the functor from `piecewise linear flat framed stratifications' to `piecewise linear stratifications' is surjective on homeomorphism classes. However, observe that that functor is far from being injective since framed stratified PL homeomorphism is a much finer equivalence relation than (nonframed) stratified PL homeomorphism, and the two categories are far from being equivalent.
\end{rmk}

A second important area of applications of the combinatorializability of flat framed stratifications concerns their `computability' properties.  We will highlight the following two results.

\begin{cor}[Canonical mesh refinement are computable] \label{cor:coarsest-mesh-comp} Given a flat framed stratification in $\lR^n$, its coarsest refining mesh can be algorithmically computed.
\end{cor}

\nid The proof of this result will translate the process of `coarsening mesh refinements' into a combinatorial notion of `reduction' of stratified trusses. Any chain of reductions of a given stratified truss eventually ends in the same normalized stratified truss (which is a consequence of the existence of coarsest mesh refinements) and this normalized stratified truss can be algorithmically computed. It follows that coarsest refining meshes also can be algorithmically computed. 

A second corollary will be the following.

\begin{cor}[Decidability of framed stratified homeomorphism] \label{cor:decidability-of-iso} Given two flat framed stratifications, presented as piecewise linear stratifications in $\lR^n$, we can algorithmically decide whether they are framed stratified homeomorphic.
\end{cor}

\nid The preceding two results will be discussed in \autoref{ssec:computability-of-fl-fr-str}.

The results in this chapter provide a first set of core properties of flat framed stratifications, in particular establishing their `combinatorial and computational tractability'. We will end the chapter in \autoref{sec:looking-ahead} with a short of overview of further steps in the program of framed combinatorial topology.

\section{Flat framed stratifications have canonical meshes} \label{sec:join-stability}

Recall from \autoref{defn:flat-framed-strat} that a flat $n$-framed stratification $(Z,f)$ is a stratification of a subspace $Z \subset \lR^n$ for which there is an $n$-mesh $M$, given by a tower of 1-mesh bundles $p_i : (M_i,f_i) \to (M_{i-1},f_{i-1})$ and $n$-framing components $\gamma_i : M_i \into \lR^i$, such that $\gamma_n : (M_n,f_n) \to (Z,f)$ is a stratified refinement of $(Z,f)$ by $(M_n,f_n)$.  We will refer to $M$ as a `refining $n$-mesh' of $(Z,f)$. Our goal in this section, will be to show that there is a canonical such refining $n$-mesh: namely, we will construct the `coarsest' refining $n$-mesh of a given flat framed stratification. The following notational convention will be helpful.

\begin{notn}[Keeping $n$-framings implicit] \label{notn:keeping-framings-implicit} We will usually keep the framing maps $\gamma_i : M_i \into \lR^i$ implicit, and instead think of the spaces $M_i$ as \emph{subspaces} $M_i \subset \lR^i$ of euclidean space. Thus $M$ is a refining mesh of a flat framed stratification $(Z,f)$ exactly if $(M_n,f_n) \to (Z,f)$ is a stratified refinement (whose underlying map is the identity on $M_n = Z \subset \lR^n$).
\end{notn}

\subsection{Definition of mesh joins}

A key step to prove the existence of coarsest mesh refinements will be the construction of so-called `mesh joins'. A `join' of two stratifications on the same space is the finest mutual coarsening of both stratifications; in general, this may be a prestratification as the following definition records (recall, a prestratification allows cycles in the formal entrance path relation, see \autoref{defn:prestrat-and-strat}).

\begin{defn}[Joins of stratification] Given stratifications $f$ and $g$ of a space $X$, the \textbf{join} $f \vee g$ is the unique prestratification of $X$ that coarsens both $f$ and $g$, and such that, for any other prestratification $h$ of $X$, if $h$ coarsens both $f$ and $g$ then it must be finer that $f \vee g$.
\end{defn}

\nid Explicitly, the join may be constructed as follows.

\begin{constr}[Joins of stratifications] \label{constr:join}
    Given a space $X$, and stratifications $(X,f), (X,g)$, let $\sim$ denote the (transitive closure of the) relation on the union of strata of both $f$ and $g$ given by
\begin{equation}
s  \sim t  \iff (s \cap t \neq \emptyset)
\end{equation}
Then the join $f \vee g$ is the prestratification of $X$ given by the decomposition of $X$ into non-empty connected subspaces $\bigcup_{s \in \is} s$, where $\is$ is an equivalence class under the relation $\sim$.
\end{constr}

\begin{rmk}[Joins as pushouts] \label{rmk:joins-as-pushouts} Joins of (pre)stratifications $(X,f)$ and $(Y,g)$ may be equivalently defined as pushouts in the category of prestratifications of the span of stratified coarsenings
    \begin{equation}
        (X,f) \ot (X,\mathsf{discr}(X)) \to (X,g)
    \end{equation}
where $\mathsf{discr}(X)$ is the discrete stratification on $X$ (see \autoref{rmk:discrete-and-indiscrete-strat}). Note that this pushout is preserved when passing to entrance path preorders, yielding a pushout in the category of preorders.
\end{rmk}

\begin{notn}[Equivalence classes and strata in joins] \label{notn:strata-in-joins} Given stratifications $(X,f)$ and $(X,g)$, abusing notation we often denote their strata by $\is \equiv \bigcup_{s \in \is} s$, where $\is$ is an equivalence class of strata---that is, we consider $\is$ both as an equivalence class of strata (for instance, we write $r \in \is$ to mean a member $r$ of that class) and as a stratum of $f \vee g$ (for instance, we write $x \in \is$ to mean a point in that stratum). Further, we denote by $\is_f$ the subclass of $\is$ consisting of strata that lie in $f$, and by $\is_g$ the subclass of strata in $\is$ that lie in $g$.
\end{notn}

\nid Note that the disjoint unions $\bigsqcup_{s \in \is_f} s$ and $\bigsqcup_{s \in \is_g} s$ in $M_n$ both equal the stratum $\is$.

\begin{eg}[Joins of stratifications] In \autoref{fig:joins-of-stratifications} we depict the joins of two stratifications. Note that in the second case, the join is in fact a prestratification, and not a stratification (using \autoref{rmk:joins-as-pushouts}, this reflects that pushouts of diagrams of posets in the category of preorders need not themselves land in the subcategory of posets).
\begin{figure}[ht]
    \centering
    \def\svgwidth{1\columnwidth}
    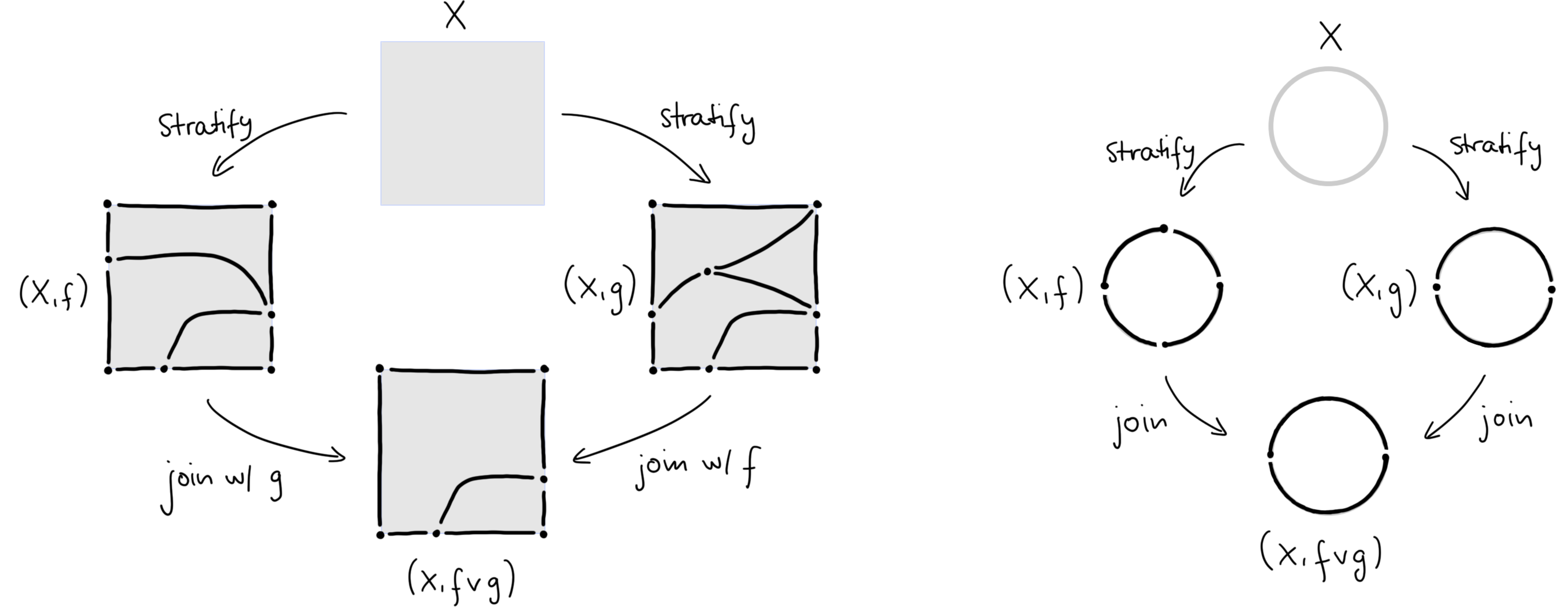

    \caption{Joins of stratifications.}
    \label{fig:joins-of-stratifications}
\end{figure}
\end{eg}

\nid For fixed domain and codomain, the join operations extend to stratified maps as the following observation records.

\begin{obs}[Joins of stratified maps] \label{obs:joins-of-maps} Given stratifications $(X,f)$, $(X,g)$, $(X',f')$, $(X',g')$, and stratified maps $F : (X,f) \to (X',f')$ and $G : (X,g) \to (X',g')$ that are identical as maps of underlying spaces $X \to X'$, then the latter map of spaces also underlies a stratified map $F \vee G : (X, f \vee g) \to (X', f' \vee g')$. We call $F \vee G$ the `join' of the stratified maps $F$ and $G$.
\end{obs}

Specifically, we will be interested in the joins of meshes. Using \autoref{notn:keeping-framings-implicit}, in the following we will identify meshes $M$ with the stratified image of their flat $n$-framings $\gamma : M \into \Pi$ in the standard euclidean tower $\Pi = (\lR^n \to \lR^{n-1} \to ... \to \lR^0)$. We call the underlying tower $M_n \to M_{n-1} \to ... \to M_0$ of subspaces $M_i \subset \lR^i$ the `support' of the mesh.

\begin{defn}[Mesh joins] Given two $n$-meshes $M$ and $N$ with identical support in $\Pi$ and consisting of projections $p_i : (M_i,f_i) \to (M_{i-1},f_{i-1})$ resp.\ $q_i : (N_i,f_i) \to (N_{i-1},f_{i-1})$ (where $M_i = N_i \subset \lR^i$ by assumption), then the \textbf{mesh join} $M \vee N$ is the tower of stratified subspaces in $\lR^i$ consisting of projections $p_i \vee q_i : (M_i, f_i \vee g_i) \to (M_{i-1}, f_{i-1} \vee g_{i-1})$.
\end{defn}

\begin{eg}[Mesh joins] In \autoref{fig:the-join-of-two-open-meshes} we depict the mesh join of two open 2-meshes. Note that mesh join is in fact again a 2-mesh.
\begin{figure}[ht]
    \centering
    \def\svgwidth{1\columnwidth}
    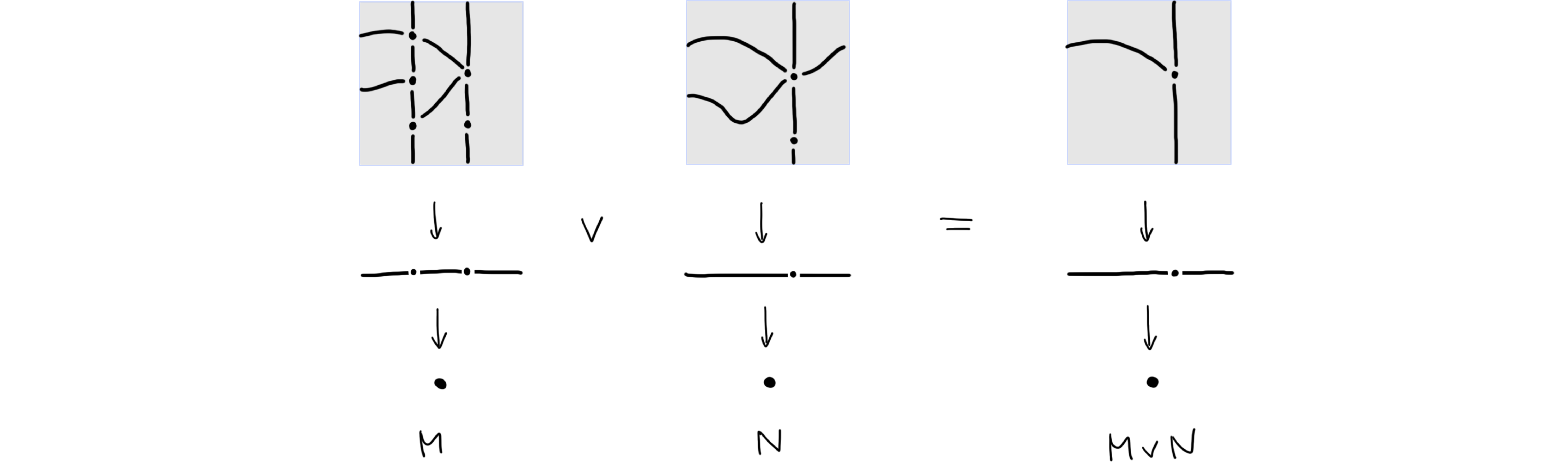

    \caption{The join of two open meshes.}
    \label{fig:the-join-of-two-open-meshes}
\end{figure}
\end{eg}

\nid The crucial fact, illustrated by the preceding example, is the following: mesh joins are meshes (automatically endowed with an $n$-framing by virtue of being subtowers of $\Pi$). The proof of this fact will take up the next section.

\subsection{Construction of mesh joins}

This section will be concerned with the proof of the following `key lemma'.

\begin{klem}[Join stability of open meshes] \label{thm:joins-of-open-meshes} Given $n$-meshes $M$ and $N$ with identical support, then their mesh join $M \vee N$ is an $n$-mesh (with the same support).
\end{klem}

Note that, in particular, if $M$ and $N$ are closed (resp.\ open) meshes then $M \vee N$ will be a closed (resp.\ open) mesh itself. Throughout this section we will denote the 1-mesh bundles in $M$ by $p_i : (M_i,f_i) \to (M_{i-1},f_{i-1})$ (or simply $p_i : f_i \to f_{i-1}$), and in $N$ by $p_i :(N_i,g_i) \to (N_{i-1},g_{i-1})$ (or simply $p_i : g_i \to g_{i-1}$), and their joins by $p_i :(M_i = N_i,f_i \vee g_i) \to (M_{i-1} = N_{i-1},f_{i-1} \vee g_{i-1})$ (or simply $p_i : f_i \vee g_i \to f_{i-1} \vee g_{i-1}$). Note, the notation reflects that all three stratified maps are identical as maps of underlying spaces $M_i = N_i$. We further usually write $M_i$ in place of $N_i$.

\begin{proof}[Proof of \autoref{thm:joins-of-open-meshes}] The proof is by induction in $n$. Inductively, the mesh join $M_{<n} \vee N_{<n}$ of the $(n-1)$-truncations $M_{<n}$ and $g_{<n}$ is an $n$-mesh; we write this $(n-1)$ mesh as $(M \vee N)_{<n}$, and its projections as $p_i : (f \vee g)_{i} \to (f \vee g)_{i-1}$. It remains to show that the stratified map $p_n : f_n \vee g_n \to (f \vee g)_{n-1}$ is in fact a constructible 1-mesh bundle. We do so in the following steps.
    \begin{enumerate}
        \item In \autoref{lem:join-bundle-surjective} we show that image of any stratum in $f_n \vee g_n$ under $p_n$ is exactly a stratum in $(f \vee g)_{n-1}$.
        \item In \autoref{lem:join-bundle-sections} we show that each stratum $\is$ in $f_n \vee g_n$ whose image under $p_n$ is the stratum $\ir$ in $(f \vee g)_{n-1}$, admits continuous sections $\gend^\pm_{\is} : \ir \to \ir \times \lR$ which fiberwise bound $\is$ from above and below.
        \item In \autoref{lem:join-bundle-constructible} we then assemble these observations into a proof that $p_n$ is a constructible 1-mesh bundle as required. \qedhere
    \end{enumerate}
\end{proof}

\begin{lem}[Joined strata project onto joined strata] \label{lem:join-bundle-surjective} For meshes $M$ and $N$ with identical support as before, the image $p_n \is$ of any stratum $\is$ of $f_n \vee g_n$ under $p_n$ is exactly a stratum in $(f \vee g)_{n-1}$.
\end{lem}

\begin{proof} First note $p_n (\is)$ is contained in some stratum $\ir$ of $(f \vee g)_{n-1}$ (see \autoref{obs:joins-of-maps}). Now, since $p_n : f_n \to f_{n-1}$ is 1-mesh bundle, images of strata in $f_n$ are exactly strata of $f_{n-1}$. Consequently, there is a subclass $\ir'_f$ of $\ir_f$ (see \autoref{notn:strata-in-joins}) such that the union of strata in $\ir'_f$ is exactly the image $p_n \is$. To proof the lemma, we show $\ir'_f = \ir_f$.  If $\ir'_f \neq \ir_f$, then we could find a stratum $r_g$ in the class $\ir_g$ which intersects both $\ir'_f$ and $\ir_f \setminus \ir'_f$. Pick a stratum $r_f$ in $\ir'_f$ that intersects $r_g$. By definition of $\ir'_f$, there is a stratum $s_f$ in the class $\is_f$ lying over $r_f$. Since $M$ and $N$ have identical support, and since $r_f$ intersects $r_g$, there is at least one stratum $s_g$ in $g_n$ lying over $r_g$ that intersects $s_f$. This implies $s_g$ is a member of the class $\is$, and thus of $\is_g$. Consequently, $p_n(s_g) = r_g$ must lie inside the subspace $\ir'_f$ contradicting our choice of $r_g$, and completing the proof of the lemma.
\end{proof}

We next want to show that the following `upper and lower fiber bounds' are in fact images of continuous sections. Note, since we take the tower $M$ (and similarly $N$) to be a subtower of the standard euclidean tower $\Pi$, the projection $p_n : M_n \to M_{n-1}$ is in particular a restriction of the projection $M_{n-1} \times \lR \to M_{n-1}$.

\begin{notn}[Upper and lower fiber bounds] Given a stratum $\is$ of $f_n \vee g_n$ lying over a stratum $\ir$ of $f_{n-1} \vee g_{n-1}$ (as shown in \autoref{lem:join-bundle-surjective}), for each point $x \in \ir$ we denote by $\is_x$ the restriction of $\is$ to the fiber $\Set{x} \times \lR$ of the projection $\ir \times \lR \to \ir$. We further denote by $\gend^-_{\is}(x)$ resp.\ $\gend^+_{\is}(x)$ the `lower resp.\ upper bound' of the subspace $\is_x$ in the fiber $\Set{x} \times \lR \iso \lR$.
\end{notn}

\begin{obs}[Upper and lower boundaries lie in singular strata] \label{rmk:pullbackjoin_bounds} We emphasize the basic but important observation that, for each $x \in \ir$, the point $\gend^\pm_{\is}(x) \in \ir \times \lR$ either equals $\gamma^\pm_n(x)$ (where $\gamma^\pm_n$ are the 1-fiber framing bounds of the top 1-mesh bundle $p_n$ of both $M$ and $N$) or it lies in singular strata of both $f_n$ and $g_n$. Indeed, it cannot lie in regular strata of either $f_n$ or $g_n$, since regular strata intersect fibers of $M_{n-1} \times \lR \to M_{n-1}$ in open intervals.
\end{obs}

\begin{lem}[Joined strata are bounded by continuous sections] \label{lem:join-bundle-sections} Given $n$-meshes $M$ and $N$ as before, and a stratum $\is$ of $f_n \vee g_n$ lying over a stratum $\ir$ of $(f \vee g)_{n-1}$, then map $\ir \to \ir \times \lR$ taking $x \mapsto \gend^\pm_{\is}(x)$ is continuous.
\end{lem}

\nid The resulting continuous maps will be denoted by $\gend^\pm_{\is} : \ir \to \ir \times \lR$.

\begin{proof} Recall meshes are regular cell stratifications (see \autoref{rmk:closed-open-meshes-are-open-regular-cell} and \autoref{rmk:general-meshes-are-open-regular-cell}). Since, by induction, we assumed that $M_{<n} \vee N_{<n}$ is an $(n-1)$-mesh, we may assume $\ir$ to be an open disk of dimension $k$.

    First observe that, for each stratum $v$ in the equivalence class $\ir$, the mapping $\gend^\pm_{\is}$ restricts to a continuous map on $v$. Indeed, assume $v \in \ir_f$ (or equally, $v \in \ir_g$). Then the class $\is$ contains a (non-empty) subclass $\is_v$ of consisting of strata lying over $v$, and $\is_v$ is exactly the intersection of $\is$ and $\ir \times \lR$. The image $\gend^\pm_{\is}(v)$ of the restricted mapping $\gend^\pm_{\is} : v \to v \times \lR$ is therefore either equal to $\gamma^\pm_n(v)$ or to some singular stratum in $f$ lying over $v$. In either case, $\gend^\pm_{\is}$ is continuous on $v$.

    To see that $\gend^\pm$ is continuous on all of $\ir$ we argue by contradiction. Assume there is a stratum $u_0 \in \ir$ such that $u_0$ contains a point $x$ at which $\gend^\pm$ is not continuous (that is, $\gend^\pm$ is non-continuous in any neighborhood of $x$ in the stratum $\ir$). By our first observation, we must have $\dim(u_0) < k$ where $k$ is the dimension of the open disk stratum $\ir$.

    We show that our assumption of discontinuity at $x$ implies that $\gend^\pm$ is discontinuous at all points of $u_0$.  Indeed, discontinuity at $x$ implies there is a sequence of points $x_i$ in $\ir$ converging to $x$ such that $\gend^\pm_{\is}(x_i)$ does not converge to $\gend^\pm_{\is}(x)$. Assume $u_0 \in \ir_f$ (the argument is symmetric if $u_0 \in \ir_g$). Since $f_{n-1}$ is a locally finite stratification (by being a finite regular cell stratification), we can pass to a subsequence $x_i$ converging to $x$ such that each $x_i$ lies in the same stratum $u'_0$ of $\ir_f$; the stratum $u'_0$ must contain $u_0$ in its closure (by boundary-constructibility of $f_{n-1}$). As initially observed, the subspaces $\gend^\pm_{\is}(u_0)$ (resp.\ $\gend^\pm_{\is}(u'_0)$) are either singular strata in $f_n$ or of the form $\gamma^\pm_n(u_0)$ (resp.\ $\gamma^\pm_n(u'_0)$). It follows by constructibility of $p_n : f_n \to f_{n-1}$ and continuity of $\gamma^\pm_n$, that either $\gend^\pm_{\is}(u'_0)$ contains \emph{all} or else \emph{none} of $\gend^\pm_{\is}(u_0)$ in its closure. The latter case must be true, since we assumed that $\gend^\pm_{\is}(x_i)$ does not converge to a point in $\gend^\pm_{\is}(u_0)$. This implies that, for any sequence $y_i \in u'_0$ converging to $y \in u_0$, the sequence $\gend^\pm_{\is}(y_i)$ cannot converge to $\gend^\pm_{\is}(y)$. We deduce that $\gend^\pm_{\is}$ is discontinuous at all points $y \in u_0$ as claimed.

    Now, let $\ir_g^{u_0}$ denote the subclass of strata in $\ir_g$ which intersect $u_0 \in \ir_f$. Note that the union of strata in $\ir_g^{u_0}$ covers $u_0$, but it cannot cover $u_0$ exactly (as this would imply $\ir = u_0$ which cannot be the case since we observed $\dim(u_0) < k$). Thus we can pick a stratum $u_1$ from $\ir_g^{u_0}$ which contains points outside of $u_0$. Since $u_1$ intersects $u_0$ non-trivially, note that $u_1$ has points at which $\gend^\pm_{\is}$ is discontinuous (in particular, again, we must have $\dim(u_1) < k$). By the same argument as before (with the roles of $f$ and $g$ reversed), this implies $\gend^\pm_{\is}$ is discontinuous at all points of $u_1$. Denote by $\ir_f^{u_0,u_1}$ the strata of $\ir_f$ which non-trivially intersect $u_0 \cup u_1$. Again, the union of strata $\ir_f^{u_0,u_1}$ cannot cover $u_0 \cup u_1$ exactly (since $\dim(u_i) < k$); thus we can pick a stratum $u_2 \in \ir_f^{u_0,u_1}$ which contains points outside of $u_0 \cup u_1$. Once more, $u_2$ inherits `discontinuity of $\gend^\pm_{\is}$' at all its points. Repeating the argument in this way, we obtain a sequence of strata $u_i$, such that each $u_i$ has points outside of $u_1 \cup ... \cup u_{i-1}$ and such that $\gend^\pm_{\is}$ is discontinuous on (all points of) $u_i$. The fact that the sequence is non-repeating contradicts the finiteness of the equivalence class $\ir$ and thus completes the argument.
\end{proof}

\begin{lem}[Constructibility of joined bundles] \label{lem:join-bundle-constructible} Given $n$-meshes $M$ and $N$ as before, the joined map $p_n : f_n \vee g_n \to (f \vee g)_{n-1}$ is a constructible 1-mesh bundle.
\end{lem}

\begin{proof} We first check that $p_n : f_n \vee g_n \to (f \vee g)_{n-1}$ is a 1-mesh bundle. Given a stratum $\is$ of $f_n \vee g_n$ lying over a stratum $\ir$ of $(f \vee g)_{n-1}$ denoted by $\gend^\pm_{\is}(ir)$ the respective images of $\ir$ under the mappings $\gend^\pm_{\is} : \ir \to \ir \times \lR$. Constructibility of both $p_n : f_n \to f_{n-1}$ and $p_n : g_n \to g_{n-1}$ implies that either $\gend^\pm_{\is}(\ir) = \gamma^\pm_n(\ir)$ or else $\gend^\pm_{\is}(\ir)$ is a union of singular strata in both $f$ and $g$. Note further we either have $\gend^-_{\is}(\ir) = \gend^+_{\is}(\ir)$ or else $\gend^-_{\is}(\ir)$ and $\gend^+_{\is}(\ir)$ are disjoint. In the former case $\is$ is a `singular stratum', that is, a section of $p_n : M_n \to M_{n-1}$ (restricted to $\ir$ in the base). In the latter case $\is$ is a `regular stratum' over $\ir$, fiberwise bounded by the continuous sections $\gend^\pm_{\is}$. This shows that $p_n : f_n \vee g_n \to (f \vee g)_{n-1}$ is stratified bundle. The bundle further inherits a fiber 1-framing since $M_n \subset \lR^n$, and its resulting bounding sections are continuous (namely, they equal the bounding sections $\gamma^\pm_n$ of $p_n : f_n \to f_{n-1}$ and $p_n : g_n \to g_{n-1}$). One verifies this endows fibers with the structure of 1-meshes as required.

    Finally, to see constructibility of $p_n : f_n \vee g_n \to (f \vee g)_{n-1}$, pick a singular stratum $\is$ over $\ir$. Note that, as required, the homeomorphism $p_n : \is \to \ir$ extends to the closure $p_n : \overline {\is} \to \overline{\ir}$, since $\is$ is a union of singular strata of $f$ and since both $p_n : f_n \to f_{n-1}$ and $p_n : g_n \to g_{n-1}$ are constructible by assumption.
\end{proof}

This completes our construction of mesh joins as meshes, and thus the proof of \autoref{thm:joins-of-open-meshes}.

\subsection{The coarsest refining mesh} \label{ssec:meshability-of-fl-fr-str}

We now define and then prove the existence of coarsest refining meshes of flat framed stratifications.

\begin{defn}[Coarsest refining meshes] Given a flat framed stratification $(Z,f)$, a refining $n$-mesh $M$, with refinement $(M_n,f_n) \to (Z,f)$, is the \textbf{coarsest refining mesh} of $(Z,f)$ if for any other refining $n$-mesh $N$, with refinement $(N_n,g_n) \to (Z,f)$, there is a refinement $(N_n,g_n) \to (M_n,f_n)$ factoring the latter refinement (of $f$ by $g_n$) through the former refinement (of $f$ by $f_n$).
\end{defn}

\nid The fact that coarsest refining meshes exists can be proven as follows.

\begin{thm}[Construction of coarsest refining meshes] \label{thm:minimal-meshes} Any flat $n$-framed stratification $(Z,f)$ has a (necessarily unique) coarsest refining $n$-mesh.
\end{thm}

\begin{proof} Given any two refining $n$-meshes $M$ and $N$ of $(Z,f)$, note that the mesh join $M \vee N$ (as constructed in  \autoref{thm:joins-of-open-meshes}) yields another $n$-mesh refining $(Z,f)$, which is coarser than both $M$ and $N$. Since meshes are finite stratifications (which excludes the possibility of infinite chains of coarsenings starting at any given mesh) it follows that there must exists a unique coarsest mesh through which all other meshes factor as claimed.
\end{proof}

\nid Importantly, the construction of coarsest refining meshes is also compatible with framed stratified homeomorphisms of flat framed stratifications, in the following sense.\footnote{Yet more generally, this observation in fact also holds for `weak equivalences' of flat framed stratifications, i.e.\ weakly invertible maps of flat framed stratifications, in an appropriate sense.}

\begin{term}[Framed stratified homeomorphism and mesh isomorphism] \label{term:fr-str-homeo} A `framed stratified homeomorphism' $(Z,f) \to (W,g)$ of flat framed stratifications is a stratified homeomorphism whose underlying map $Z \to W$ is framed, meaning that it descends along the standard euclidean projection $\pi_i : \lR^n \to \lR^i$ (see \autoref{term:framed-maps}). Note that, in contrast, we continue to speak of `$n$-mesh isomorphisms' to mean isomorphisms in the category $\tmesh n$ (since meshes are not just mere stratifications, but towers of 1-mesh bundles).
\end{term}

\begin{lem}[Homeos of flat framed stratifications induce isos of their coarsest refining meshes] \label{lem:iso-mesh} Let $(Z,f)$ and $(W,g)$ be flat $n$-framed stratifications with coarsest refining meshes $M$ resp.\ $N$. If there is a framed stratified homeomorphism $F : f \iso g$ then $F$ induces an $n$-mesh isomorphism $F : M \to N$ between their respective coarsening refining meshes.
\end{lem}

\begin{proof} Since $F$ is a framed stratified homeomorphism, note that there is a `push-forward' mesh $FM$ refining $g$: the mesh $FM$ is determined by an $n$-mesh isomorphism $M \iso FM$ whose top component is $F$. Conversely, there is a `pull back' mesh $F\inv N$ refining $f$ determined by $N \iso F\inv N$ having top component $F\inv$. Note that, since $N$ is the coarsest refining mesh of $g$, the refining mesh $FM$ must be finer than $N$. If $FM$ were strictly finer than $N$ (that is, $FM \to N$ is a non-identity coarsening of $FM$), then $M = F\inv F M \to F\inv N$ would be a non-identity coarsening of $M$. This is impossible since $M$ is the coarsest refining mesh of $g$. We deduce that $FM = N$, and thus we obtain an induced $n$-mesh isomorphism $M \iso N$ with top component $F : Z \to W$. Abusing notation, we denote this mesh isomorphism again by $F : M \iso N$.
\end{proof}

Having shown the existence of coarsest refining meshes (and their compatibility with framed stratified homeomorphism of flat framed stratification), we now give several illustrative examples. The slogan for these illustrations is that `coarsest refining meshes are refining meshes that record exactly where changes in framed stratified homeomorphism type of a flat framed stratification happen'. (Heuristically, they bear visual resemblance with `critical level set' decompositions from (higher) Morse theory.)

\begin{eg}[Coarsest refining meshes of the circle, the braid and Hopf circle] In \autoref{fig:coarsest-refining-meshes-of-circle-and-hopf-circle} we depict three flat framed stratifications (together with their respective coarsest refining meshes): these include flat 2-framed stratification of the open 2-cube obtained by an embedded circle and its complement, the flat 3-framed stratification of the open 3-cube given by the braid and its complement, and the flat 3-framed stratification given by the `Hopf circle' in the 3-cube and its complement. Observe that the last two examples recover the first two examples in \autoref{fig:flat-framed-stratifications} up to a rotation of the ambient $\lR^3$.
\begin{figure}[ht]
    \centering
    \def\svgwidth{1\columnwidth}
    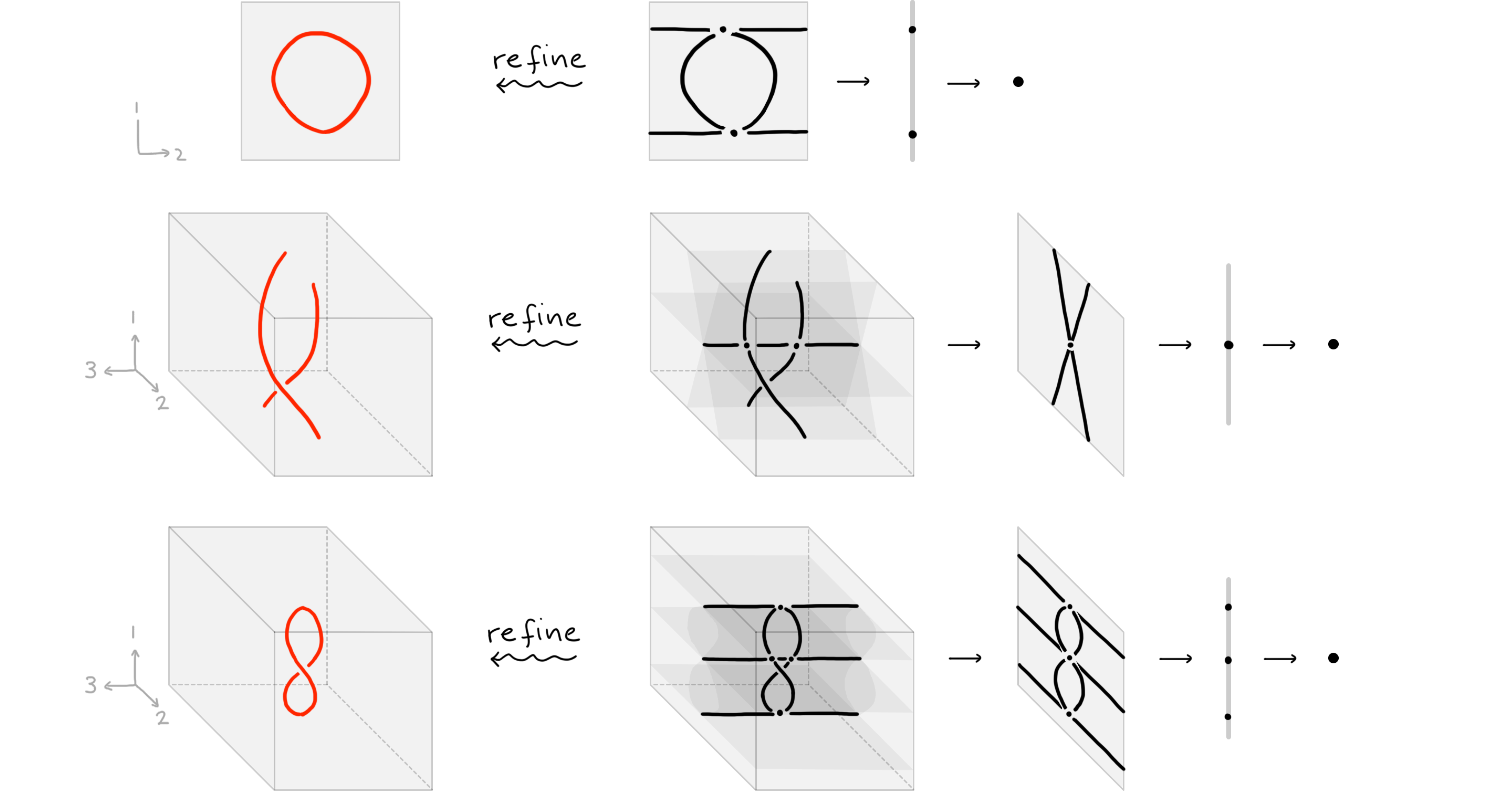

    \caption{Coarsest refining meshes of the circle, the braid, and Hopf circle.}
    \label{fig:coarsest-refining-meshes-of-circle-and-hopf-circle}
\end{figure}
\end{eg}

\begin{eg}[Non-coarsest refining meshes of the circle] In \autoref{fig:non-coarsest-refining-meshes-of-the-circle} we depict two refining $2$-meshes of a flat framed stratification of the open 2-cube obtained by an embedded circle and its complement. Neither of the two refining meshes is the coarsest refining mesh.
\begin{figure}[ht]
    \centering
    \def\svgwidth{1\columnwidth}
    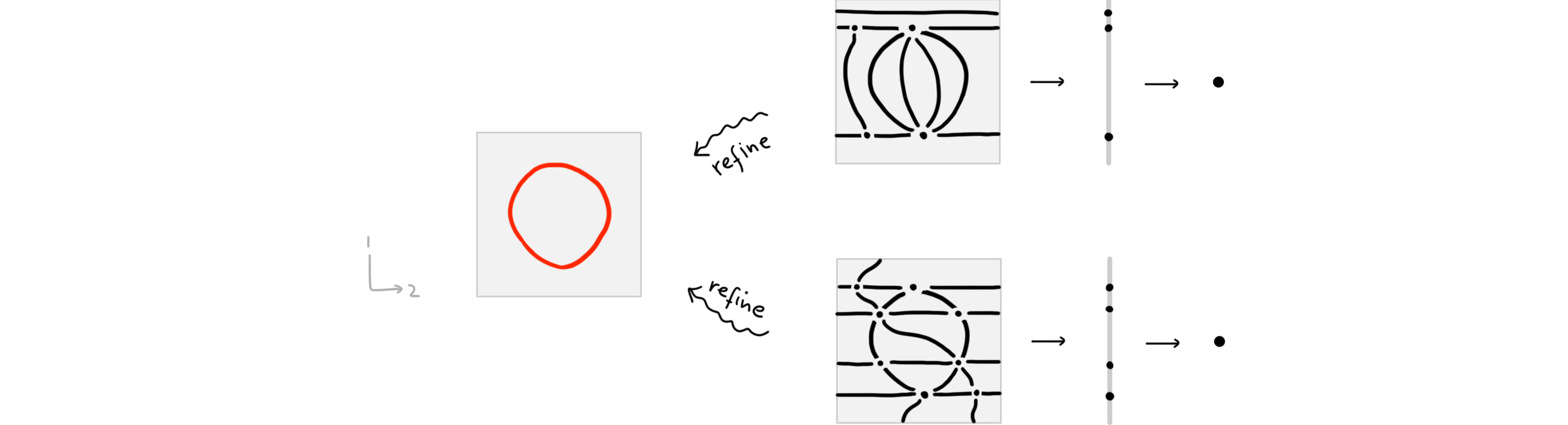

    \caption{Non-coarsest refining meshes of the circle.}
    \label{fig:non-coarsest-refining-meshes-of-the-circle}
\end{figure}
\end{eg}

\begin{eg}[Coarsest refining mesh of the cusp] In \autoref{fig:coarsest-refining-mesh-of-the-cusp} we depict the coarsest refining mesh of the cusp singularity (which appeared as the last example in our earlier \autoref{fig:flat-framed-stratifications}).
\begin{figure}[ht]
    \centering
    \def\svgwidth{1\columnwidth}
    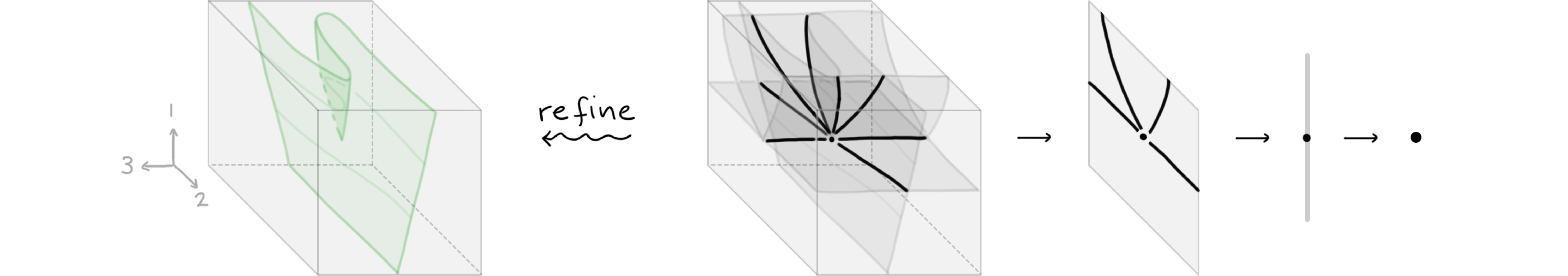

    \caption{Coarsest refining mesh of the cusp.}
    \label{fig:coarsest-refining-mesh-of-the-cusp}
\end{figure}
\end{eg}

\begin{rmk}[Twisted circle embedding] \label{rmk:resolving-the-twisted-circle-example} We revisit the Hopf circle embedding from \autoref{eg:generality-of-blocks}. As a flat framed stratification of the open $3$-cube, its coarsest refining mesh is an open 3-mesh $M$, as shown in \autoref{fig:coarsest-refining-meshes-of-circle-and-hopf-circle}. As a stratification of the closed cube, its coarsest refining mesh $\overline M$ is a closed 3-mesh (and equals, up to isomorphism, the classifying mesh $\CMsh \overline {\ETrs M}$ of the compactification $\overline {\ETrs M}$ of the entrance path truss $\ETrs M$ of $M$). The mesh complex of $\overline M$ (see \autoref{term:mesh-complex}) recovers exactly the flat framed regular cell complex given earlier in \autoref{eg:generality-of-blocks}.
\end{rmk}

\section{Tractability of flat framed stratifications} \label{sec:tractability-of-fl-fr-str}

In this section we discuss the `combinatorializability' of flat framed stratifications and its corollaries, as outlined in the beginning of the chapter. We prove the central theorem about the combinatorializability of flat framed stratifications in \autoref{ssec:combinatorializability-of-fl-fr-str}. A first set of corollaries relating to the `polyhedrality' of flat framed stratifications will be discussed in \autoref{ssec:polyhedrality-of-fl-fr-str}. Finally, the theorem also implies the `computability' of flat framed stratification as we will explain in \autoref{ssec:computability-of-fl-fr-str}.

\pagebreak
\subsection{Combinatorializability} \label{ssec:combinatorializability-of-fl-fr-str}

\subsubsecunnum{Stratified trusses} Recall from \autoref{defn:norm-strat-trusses}, that a stratified $n$-truss is a poset labeled $n$-truss $T$ (that is, a `underlying' $n$-truss $\und T = (T_n \to ... \to T_0)$ together with a `labeling' $\lbl T : T_n \to P$ in some poset $P$) with the condition that the labeling $\lbl T$ is the characteristic map of a stratification on $T_n$. This last condition can be more concretely rephrased as follows.

\begin{rmk}[Characteristic maps are connected-quotient maps]  A `quotient map' of posets is a surjective poset map for which subposets in the codomain are open, i.e.\ downward closed, if and only if their preimages are open in the domain. A `connected-quotient map' of posets is a quotient map of posets whose preimages are connected. A truss labeling $\lbl T : T_n \to P$ in a poset $P$ is the characteristic map of a stratification on $T_n$ if and only if it is a connected-quotient map. This is shown (in slightly more general form) in \autoref{lem:cquot_are_charact}.
\end{rmk}

\nid To highlight that a labeled truss $T = (\und T, \lbl T)$ is a \emph{stratified} truss we will usually denote the poset of labels by $\Entr(T)$ (and thus write $\lbl T : T_n \to \Entr(T)$ for the labeling of $T$). Before discussing the notion of stratified $n$-trusses further, let us give a first example.

\begin{eg}[Stratified trusses] In \autoref{fig:stratified-truss-examples} we depict a stratified open 2-truss.
\begin{figure}[ht]
    \centering
    \def\svgwidth{1\columnwidth}
    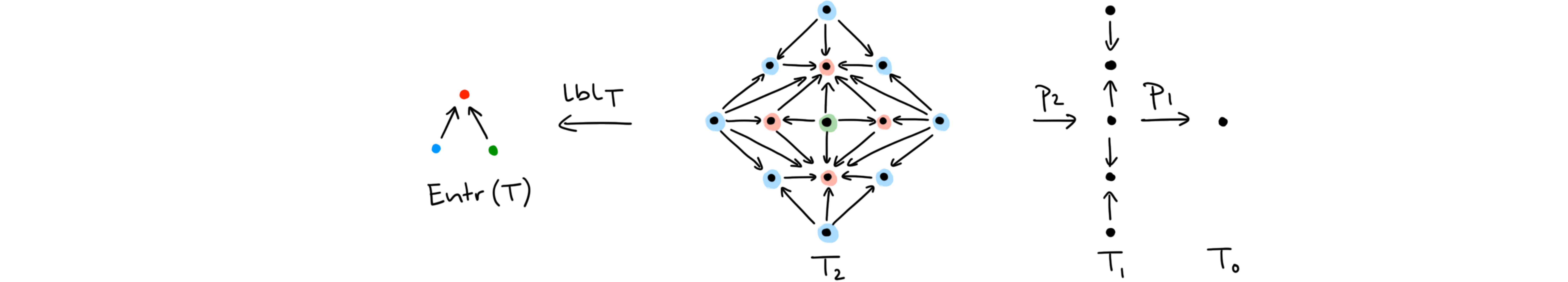

    \caption{A stratified open 2-truss.}
    \label{fig:stratified-truss-examples}
\end{figure}
\end{eg}

\begin{term}[Entrance path posets and strata of stratified trusses] Given a stratified $n$-truss $T \equiv (\und T, \lbl T : T_n \to \Entr(T))$ we call $\Entr(T)$ the `entrance path poset' of $T$. The `strata' of $T$ are the (by definition, connected) subposets of $T_n$ given by the preimages of $\lbl T$.
\end{term}

\nid Note that, up to isomorphism of entrance path posets, any decomposition of $T_n$ into connected subposets fully determines a stratified $n$-truss $T$. One may further obtain stratified trusses from `$P$-labelings' as follows.

\begin{rmk}[Poset-labeled truss yields stratified trusses] Note that any poset labeled truss $T = (\und T,\lbl T : T_n \to P)$ yields a stratified $n$-truss $\tilde T$, by defining strata of $\tilde T$ to be the connected components of the non-empty preimages of $\lbl T$. (This `connected component splitting' construction is formalized, in more general form, in \autoref{prop:conn_comp_split}.)  While we will not formally need this construction here, it turns out to be somewhat useful when illustrating examples: it allows us to reduce the number of colors needed in the depiction of stratified trusses, by replacing a given connected-quotient map with some (non-unique choice of) poset map that has smaller domain and whose connected component splitting recovers the original connected-quotient map. 

    For instance, in \autoref{fig:stratified-truss-examples} we may replace the labeling by a map to the poset $(0 \to 1)$, mapping both the blue and green stratum to $0$, and the red stratum to $1$.
\end{rmk}

\pauseae

The notion of maps for stratified $n$-trusses carries over from that of labeled $n$-trusses, but we may further make the following terminological distinctions.

\begin{term}[Stratified maps and coarsenings] \label{term:strat-truss-coarsenings} A `stratified map' of stratified $n$-trusses $T$ and $S$ is a labeled $n$-truss map $F \equiv (\und F, \lbl F) : T \to S$ (that is, a map of $n$-trusses $\und F : \und T \to \und S$, together with a map of labelings $\lbl F : \Entr(T) \to \Entr(S)$ such that $\lbl F \lbl T = \lbl S F_n$ commutes).
    \begin{enumerate}
        \item We call $F$ a `(truss preserving) label coarsening' if $\und F = \id$ is the identity truss map and $\lbl F$ is a connected-quotient map (see \autoref{defn:cquot-map-of-posets}).
        \item We call $F$ `(label preserving) truss coarsening' if $\und F$ is an $n$-truss coarsening (see \autoref{term:coarsening-and-quotients}) and $\lbl F = \id$.
        \item We call $F$ a `coarsening' if $\und F$ is an $n$-truss coarsening and $\und F$ is a connected-quotient map. \qedhere
    \end{enumerate}
\end{term}

\nid Note that the definition of `label coarsenings' reflects that, in the context of topological stratifications, connected-quotients of entrance path posets to poset describe exactly the coarsenings of a given stratification (see \autoref{rmk:coarsenings_as_cquotients}).

\begin{eg}[A truss coarsening] We depict a label preserving truss coarsening $F : T \to S$ of stratified $2$-trusses $S$ and $T$ in \autoref{fig:a-truss-coarsening-of-stratified-trusses}.
\begin{figure}[ht]
    \centering
    \def\svgwidth{1\columnwidth}
    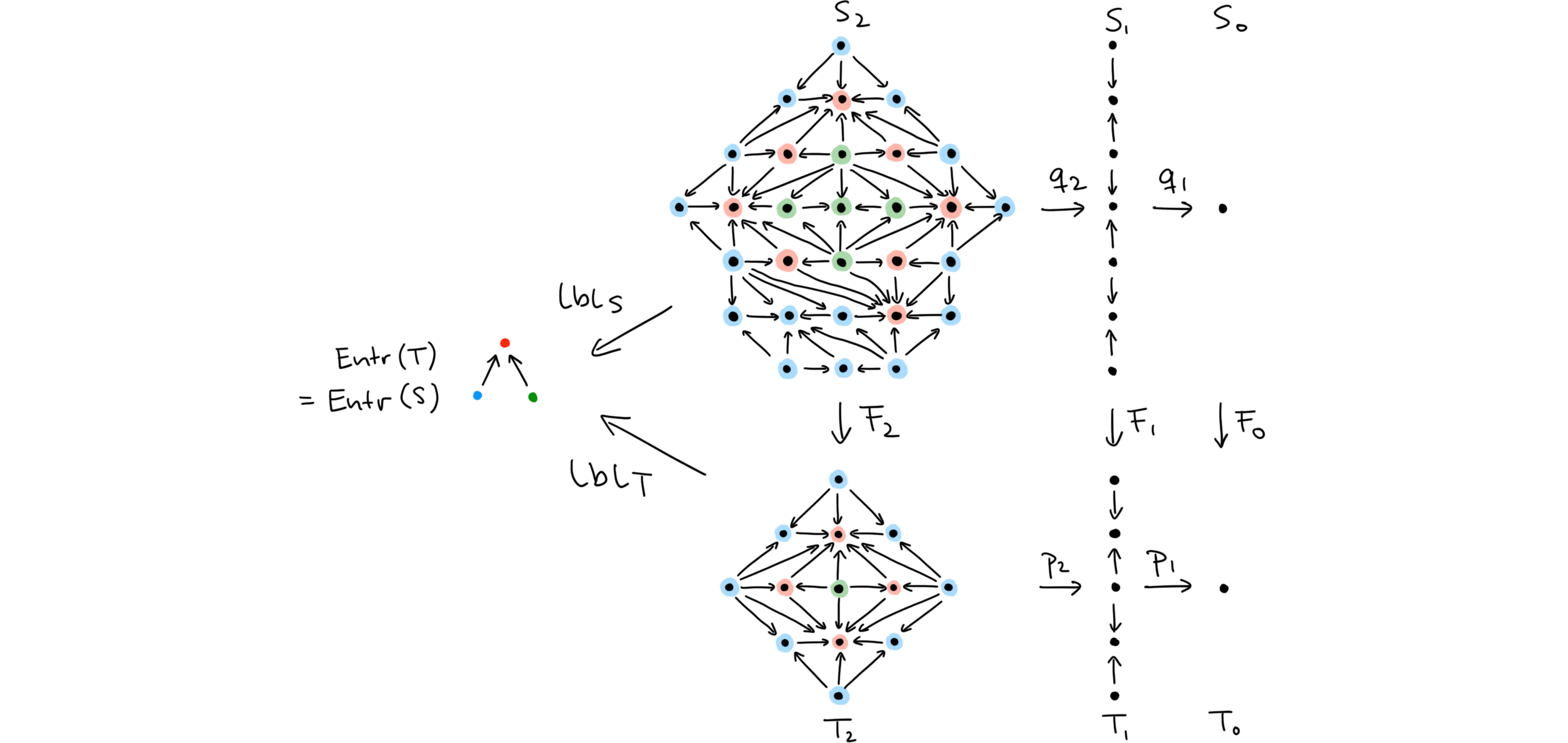

    \caption{A truss coarsening of stratified trusses.}
    \label{fig:a-truss-coarsening-of-stratified-trusses}
\end{figure}
\end{eg}

\nid Note that every coarsening can be written both as a unique composite of a truss coarsening followed by a label coarsening, and as a unique composite of a label coarsening followed by a truss coarsening.

\subsubsecunnum{Stratified meshes} Stratified trusses (as well as their maps) can be obtained from flat framed stratifications together with chosen refining meshes. We call the latter structure a `stratified mesh'.

\begin{defn}[Stratified mesh] \label{defn:stratified-meshes} A flat $n$-framed stratification $(Z,f)$ together with a refining mesh $M$ of $(Z,f)$ will be called a \textbf{stratified mesh} and denoted by a tuple $(M,f)$ (note that $Z$ may be recovered as the subspace $M_n \subset \lR^n$).
\end{defn}

\begin{defn}[Stratified mesh map] A \textbf{map of stratified meshes} $F : (M,f) \to (N,g)$ is a mesh map $F : M \to N$ whose top component is a map of stratifications $F_n : f \to g$.
\end{defn}

\begin{term}[Maps of stratified meshes and their coarsenings] Given a map of stratified meshes $F : (M,f) \to (N,g)$ we make the following terminological distinction.
    \begin{enumerate}
        \item If $F_n : f \iso g$, and if $N$ is a coarser refining mesh of $f$ than $M$ (i.e.\ $F : M \to N$ is an $n$-mesh coarsening) then we call $F : (M,f) \to (N,f)$ a `mesh coarsening' of stratified meshes.
        \item If $F : M \iso N$ and $g$ is coarser than $f$  (i.e.\ $F_n : f \to g$ is a coarsening), then we call $F : (M,f) \to (N,f)$ a `stratification coarsening' of stratified meshes.
        \item If both $F : M \to N$ and $F_n : f \to g$ are coarsenings, we call $F : (M,f) \to (N,f)$ simply a `coarsening' of stratified meshes.
    \end{enumerate}
\end{term}

\begin{eg}[A mesh coarsening] We depict a mesh coarsening $F : (M,f) \to (N,g)$ of stratified $2$-meshes $(M,f)$ (with $M = (q_2,q_1)$) and $(N,g)$ (with $N = (p_2,p_1)$) in \autoref{fig:a-mesh-coarsening-of-stratified-meshes}.
\begin{figure}[ht]
    \centering
    \def\svgwidth{1\columnwidth}
    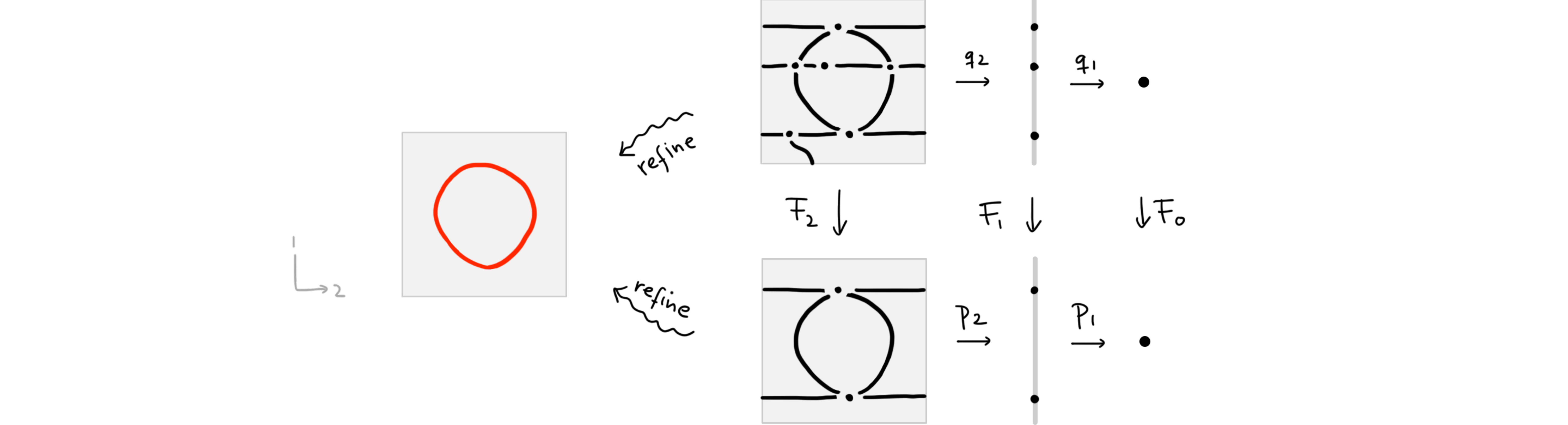

    \caption{A mesh coarsening of stratified meshes.}
    \label{fig:a-mesh-coarsening-of-stratified-meshes}
\end{figure}
\end{eg}

\subsubsecunnum{Stratified entrance path trusses} Every stratified mesh has a `stratified entrance path truss'.

\begin{defn}[Stratified entrance path trusses] Given a stratified $n$-mesh $(M,f)$ define the \textbf{stratified entrance path truss} $\ETrs (M,f)$ to be the stratified $n$-truss $(\ETrs M, \lbl {M,f})$ whose underlying $n$-truss is the entrance path truss $\ETrs M$ of $M$, and whose labeling $\lbl {M,f} : T_n \to \Entr(f)$ is the entrance path poset map $\Entr(f_n \to f) : \Entr f_n \to \Entr f$ of the stratified coarsening $(M_n,f_n) \to (Z,f)$.
\end{defn}

\begin{defn}[Stratified entrance path truss maps] Given a stratified mesh map $F : (M,f) \to (N,g)$, define the \textbf{stratified entrance path truss map} $\ETrs F : \ETrs (M,f) \to \ETrs (N,g)$ to be the map of stratified $n$-trusses whose underlying map of trusses $\ETrs (M \to N)$ is the entrance path truss map of the mesh map $F : M \to N$, and whose label map $\Entr(f \to g)$ is the entrance path poset map of the stratified map $F_n : f \to g$.
\end{defn}

\begin{obs}[Translations of notions of coarsenings] \label{obs:transl-notion-of-coarsening} Notions of coarsenings of stratified meshes translate to the corresponding notions of coarsenings of stratified trusses. Namely, given a mesh coarsening, a stratified coarsening, a coarsening of stratified meshes $(M,f) \to (N,f)$, then $\ETrs((M,f) \to (M,g)) : \ETrs(M,f) \to \ETrs (M,g)$ is a truss coarsening, resp.\ a label coarsening, resp.\ a coarsening of stratified trusses.
\end{obs}

\begin{eg}[Translating mesh coarsenings to truss coarsenings] The mesh coarsening in \autoref{fig:a-mesh-coarsening-of-stratified-meshes} translates to the truss coarsening in \autoref{fig:a-truss-coarsening-of-stratified-trusses}.
\end{eg}

\subsubsecunnum{Stratified classifying meshes} The stratified entrance path truss construction has a converse construction as follows. Recall that stratified coarsenings of a given stratification are exactly determined by connected-quotient maps of the stratification's entrance path poset (see \autoref{rmk:coarsenings_as_cquotients}).

\begin{defn}[Stratified classifying meshes of stratified trusses] \label{defn:stratified-class-mesh} Let $T$ be a stratified $n$-truss. Its \textbf{stratified classifying mesh} $\CMsh T$ is the stratified mesh $(\CMsh \und T, \sff_{T})$ whose underlying mesh $\CMsh \und T$ is the classifying mesh of the underlying $n$-truss $\und T$ of $T$, and whose stratification $\sff_{T}$ is the stratification of $Z = M_n \subset \lR^n$ obtained by coarsening $(M_n,f_n)$ by the connected-quotient map $\lbl T : T_n \to \Entr(T)$ on entrance path posets.
\end{defn}

\nid Note in particular, that $(Z = M_n,\sff_{T,f})$ is a flat framed stratification. The stratified classifying mesh construction extends to maps of stratified trusses as follows.

\begin{defn}[Stratified classifying mesh maps of stratified trusses maps] Given a map $F = (\und F,\lbl F) : T \to S$ of stratified trusses $T = (\und T,\lbl T)$ and $S = (\und S, \lbl S)$, we define the \textbf{stratified classifying mesh map} $\CMsh F : \CMsh T \to \CMsh S$ of stratified $n$-meshes to be the classifying $n$-mesh map $\CMsh \und F : \CMsh \und T \to \CMsh \und S$ which, by definition of stratified classifying meshes, can be seen to descend to a stratified map $\CMsh F_n : \sff_T \to \sff_S$ as required.
\end{defn}

\nid The construction are mutually (weakly) inverse in the following sense.

\begin{obs}[Mutual invertibility] \label{obs:stratified-cmsh-etrs-invertibility} Let $T$ be a stratified truss, and $(M,f)$ be a stratified mesh.
    \begin{enumerate}
        \item There is a unique isomorphism $T \iso \ETrs \CMsh (M,f)$. This follows since $\und T \iso \ETrs \CMsh M$ uniquely, and from the definitions of stratified entrance path trusses resp.\ stratified classifying meshes.
        \item There is a (up to contractible choice of homotopy unique) isomorphism of stratified meshes $(M,f) \iso \CMsh \ETrs (M,f)$. This follows since there is a (up to contractible choice of homotopy unique) mesh isomorphism $M \iso \CMsh \ETrs M$, and from the definitions of stratified classifying meshes resp.\ stratified entrance path trusses. \qedhere
    \end{enumerate}
\end{obs}

\begin{obs}[Translations of notions of coarsenings] \label{obs:transl-notion-of-coarsening-converse} Conversely to \autoref{obs:transl-notion-of-coarsening}, notions of coarsenings of stratified trusses translate to the corresponding notions coarsenings of stratified meshes. Namely, given a truss coarsening (resp.\ label coarsening resp.\ coarsening) of stratified trusses $T \to S$, then $\ETrs(T \to S) : \ETrs T \to \ETrs S$ is a mesh coarsening (resp.\ stratification coarsening resp.\ coarsening) of stratified meshes.
\end{obs}

\nid While not used for the main result of the chapter, the case of coarsenings will be of separate importance to us, and, as discussed before, it needs additional care (see \autoref{obs:class-crs-asymmetry}). Recall that truss coarsenings can be realized as mesh coarsenings, via the `classifying mesh coarsening' construction (see \autoref{constr:PL-realization-coarsening-general}).

\begin{defn}[Stratified classifying mesh coarsening of stratified trusses coarsening] Given a coarsening of stratified trusses $F : T \to S$, we define the \textbf{stratified classifying mesh coarsening} $\CMsh F : \CMsh T \to \CMsh S$ of stratified $n$-meshes to be the classifying $n$-mesh coarsening $\CrsMsh \und F : \CMsh \und T \to \CMsh \und S$ which, by definition of stratified classifying meshes, can be seen to descend to a stratified map $\CrsMsh F_n : \sff_T \to \sff_S$ as required.
\end{defn}

\subsubsecunnum{Normalization and combinatorializability}

Recall the definition of `coarsest refining meshes $M$' of flat framed stratifications $(Z,f)$: this is precisely a stratified meshes $(M,f)$ to which no non-identity mesh coarsening applies. Combinatorially, this motivates our earlier \autoref{defn:norm-classes} of `normalized' stratified trusses: namely, a normalized stratified truss $T$ is precisely a stratified truss to which no non-identity truss coarsening applies. Our earlier \autoref{obs:stratified-cmsh-etrs-invertibility} about the correspondence of stratified meshes and stratified trusses now specializes to a correspondence of coarsest stratified meshes (meaning stratified meshes $(M,f)$ such that $M$ is the coarsest refining mesh of $f$) and normalized stratified trusses up to isomorphism. This `almost' proves \autoref{thm:flat-fr-strat-are-norm-strat-trusses}, claiming the correspondence of flat framed stratifications and normalized stratified trusses up to isomorphism. The only missing ingredients lies in the comparison of framed stratified homeomorphisms of flat framed stratifications and isomorphisms of coarsest stratified meshes; but this comparison was exactly provided by our earlier \autoref{lem:iso-mesh}.

\begin{proof}[Proof of \autoref{thm:flat-fr-strat-are-norm-strat-trusses}]
    Given a normalized stratified truss $T$, define its corresponding flat framed stratification to be $\sff_T$ where $(\CMsh T,\sff_T)$ is the stratified classifying mesh of $T$. Note that, given another normalized stratified truss $S$ and $T \iso S$ by a label preserving balanced truss isomorphism, then $(\CMsh T,\sff_T) \iso (\CMsh S,\sff_S)$ as stratified meshes which descends to a framed stratified homeomorphism $\sff_T \iso \sff_S$ of flat framed stratifications. Conversely, given a flat framed stratification $(Z,f)$, define its corresponding normalized stratified truss $T$ to be the stratified entrance path truss $\ETrs (M,f)$ where $M$ is the coarsest refining mesh of $f$. Note that, given another flat framed stratification $(W,g)$, with $f \iso g$, we obtain an isomorphism of stratified meshes $(M,f) \iso (N,g)$ (where $N$ is the coarsest refining mesh of $g$) by \autoref{lem:iso-mesh}, and thus a label preserving balanced truss isomorphism $\ETrs (M,f) \iso \ETrs (N,g)$. Now, \autoref{obs:stratified-cmsh-etrs-invertibility} implies that both constructions are mutually inverse to each other, which proves the claim of the theorem.
\end{proof}

\begin{cor}[Normalization vs coarsest mesh refinements] Let $(M,f)$ be a stratified mesh and $T$ a stratified truss such that $T \iso \ETrs (M,f)$ (or equivalently, $(M,f) \iso \CMsh T$). Then $T$ is normalized if and only if $M$ is the unique coarsest mesh refining $f$. \qed
\end{cor}

We will use the rest of this section to briefly illustrate the correspondence established by \autoref{thm:flat-fr-strat-are-norm-strat-trusses}.

\begin{eg}[The normalized stratified truss of the circle] The normalized stratified truss corresponding to the flat framed stratification of the circle in the open 2-cube (see the first example in \autoref{fig:coarsest-refining-meshes-of-circle-and-hopf-circle}) is the stratified truss given in \autoref{fig:stratified-truss-examples}.
\end{eg}

\begin{eg}[The normalized stratified truss of the braid] \label{eg:revisiting-the-braid} The normalized stratified truss corresponding to the braid (see the second example in \autoref{fig:coarsest-refining-meshes-of-circle-and-hopf-circle}) is the stratified 3-truss depicted in \autoref{fig:the-normalized-stratified-truss-of-the-braid} (note that for readability we depict only generating arrows, see \autoref{rmk:generating-arrows}). Note that this fulfils our earlier promises in \autoref{preview:stratified-trusses}.
\begin{figure}[ht]
    \centering
    \def\svgwidth{1\columnwidth}
    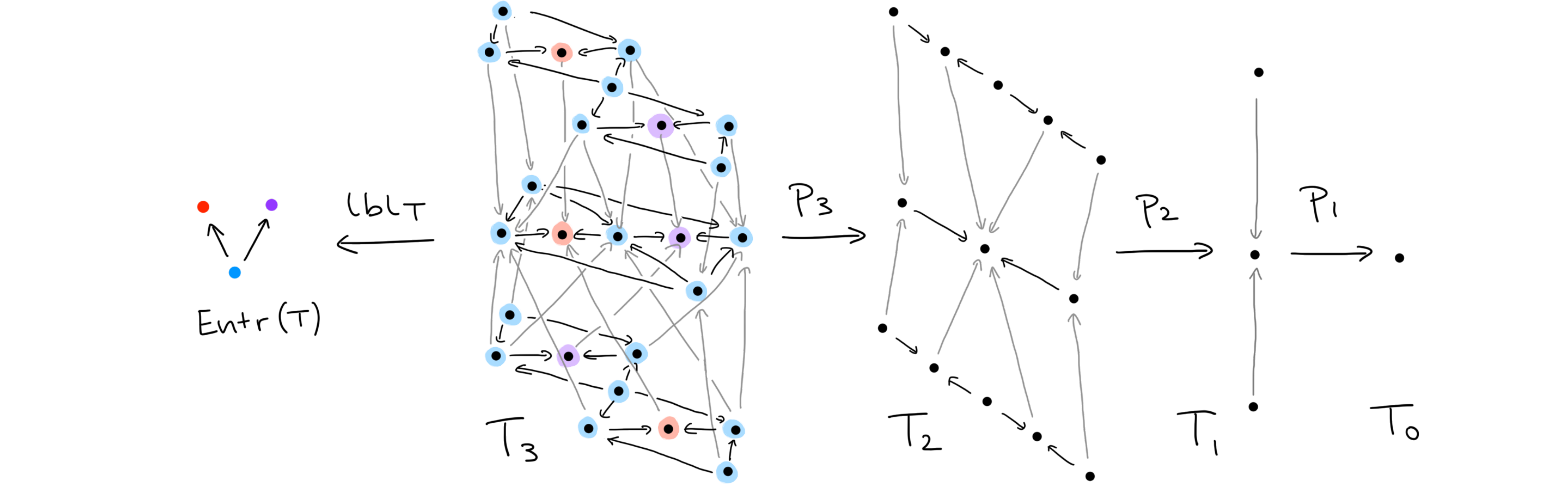

    \caption{The normalized stratified truss of the braid.}
    \label{fig:the-normalized-stratified-truss-of-the-braid}
\end{figure}
\end{eg}

\begin{eg}[The normalized stratified truss of the cusp] The normalized stratified truss corresponding to the cusp singularity (see \autoref{fig:coarsest-refining-mesh-of-the-cusp}) is the stratified 3-truss depicted in \autoref{fig:the-normalized-stratified-truss-of-the-cusp} (note again, that for readability we depict only generating arrows).
\begin{figure}[ht]
    \centering
    \def\svgwidth{1\columnwidth}
    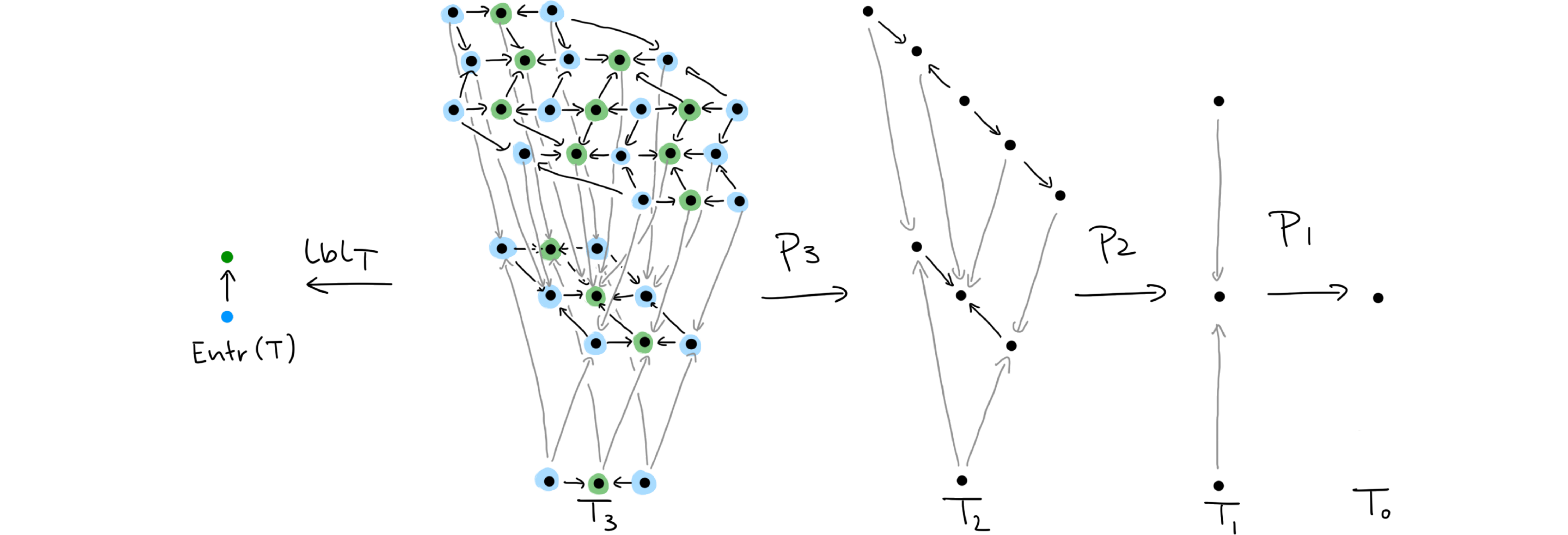

    \caption{The normalized stratified truss of the cusp.}
    \label{fig:the-normalized-stratified-truss-of-the-cusp}
\end{figure}
\end{eg}

\subsection{Polyhedrality} \label{ssec:polyhedrality-of-fl-fr-str}

In this section we discuss the relationship of flat framed stratifications and `piecewise linear stratifications'. We introduce the latter notion below. Note that all simplicial complexes $K$ in this section will be assumed to be compact and `linearly realized' in $\lR^n$, i.e. they come equipped with an embedding $\abs{K} \into \lR^n$ of their geometric realization with compact image, and such that the embedding is linear on each simplex (in particular, note that the underlying space of a simplicial complex, as a subspace of $\lR^n$, is a polyhedron in the usual sense, see \cite[Defn. 1.1]{rourke2012introduction}).

\begin{defn}[Piecewise linear stratifications] A \textbf{piecewise linear stratification} is a stratification $(Z,f)$, $Z \subset \lR^n$, for which there exists a `triangulating' simplicial complex $K$ such that each stratum of $f$ is a union of open simplices in $K$.
\end{defn}

\nid Note that we do not require $Z$ to have the same support as $K$ (as a subspace in $\lR^n$). In particular, $Z$ itself may be non-compact. We further record terminology for those piecewise linear stratifications whose underlying space $Z$ is that of a mesh.

\begin{term}[Mesh-supported piecewise linear stratifications] \label{term:mesh-supported} A piecewise linear stratification $(Z,f)$ is called `mesh-supported' if $Z \subset \lR^n$ is the space underlying some $n$-mesh $M \subset \lR^n$.
\end{term}

\nid Note that every piecewise linear stratification $(Z,f)$ can be trivially extended to a mesh-supported one, by the following construction. Pick $r \in \lR$ large enough such that the open $n$-cube $(-r,r)^n$ contains $Z$ and such that $(-r,r)^n \setminus Z$ has one connected component. Since the open cube $(-r,r)^n$ underlies an open $n$-mesh (for instance, the terminal open $n$-mesh), we obtain a stratification $(Z^+,f^+)$ where $Z^+ = (-r,r)^n$ and $f^+$ contains strata of $f$ together with $(-r,r)^n \setminus Z$ as a new `ambient' stratum. One checks that $(Z^+,f^+)$ is again a piecewise linear stratification (in fact, this also follows from our later \autoref{constr:simplicial-meshes-of-polyhedra}). Therefore, `up to adding an ambient stratum', we may think of all piecewise linear stratification as being mesh-supported. Note, for notational convenience (and `up to scaling') we will usually assume $r = 1$, and work in the open unit $n$-cube $\II^n := (-1,1)^n$.

As a first, and most immediate corollary to the combinatorializability of flat framed stratifications, let us observe that any flat framed stratification is in fact framed stratified homeomorphic to a piecewise linear stratification (which is necessarily `mesh-supported'). This proves our earlier claim in \autoref{cor:flat-framed-strat-are-PL}.

\begin{proof}[Proof of \autoref{cor:flat-framed-strat-are-PL}] Given a flat framed stratification $(Z,f)$, let $M$ be its coarsest refining mesh. Then $f$ is framed stratified homeomorphic to the piecewise linear stratification $\sff_{\ETrs (M,f)}$ obtained from the stratification of the stratified classifying mesh $\CMsh \ETrs (M,f)$ (note that the fact that $\sff_{\ETrs (M,f)}$ is piecewise linear follows from our construction of classifying meshes, which realizes truss posets linearly in $\lR^n$, see \autoref{constr:geo_closed} and \autoref{constr:geo_open}).
\end{proof}

We will now prove two further corollaries of the combinatorializability of flat framed stratifications in this section. Firstly, we address a `converse' to \autoref{cor:flat-framed-strat-are-PL}: namely, we will show that every piecewise linear stratification is, in fact (up to ambient extension) a flat framed stratification, which proves the claim of our earlier \autoref{prop:PL-strat-are-flat-frmd}. Secondly, we will discuss a variation of a classical question posed by the so-called `Hauptvermutung': this asks whether topological homeomorphism of piecewise linear structures can be, up to homotopy, replaced by piecewise linear homeomorphism. As we will recall, the conjecture fails to hold in the classical setting of piecewise linear stratifications in $\lR^n$. We then show the conjecture does hold when appropriately adapted to the flat framed setting, proving our earlier \autoref{cor:flat-framed-hpt-vmtg}.

\subsubsecunnum{Piecewise linear stratifications are flat framed}

We prove the claim that any piecewise linear stratification (up to passing to its `ambient extension') defines a flat framed stratification in $\lR^n$. The proof will rely on the following observations.

\begin{obs}[Intersection refinements] \label{rmk:intersection-refinement} Let $K$ be a finite simplicial complex, and $F : \abs{K} \to \lR^n$ a simplex-wise linear (not-necessarily injective) map. By standard results, there exists a simplex-wise linear embedding $G : \abs{L} \into \lR^{n-1}$ of some simplicial complex $L$, such that the image $F\abs{x}$ of each simplex $x$ in $K$ is a union of images $G\abs{y}$ of simplices $y$ in $L$ (see \cite[Thm. 2.15]{rourke2012introduction}). We call the simplicial complex $L$ (linearly realized in $\lR^n$ via $G$) an `intersection refinement' of $F$.
\end{obs}

\begin{obs}[Refining meshes of simplicial complexes] \label{constr:simplicial-meshes-of-polyhedra} Let $K$ be a simplicial complex linearly realized in $\lR^n$. Up to scaling, we may assume $K$ lives in the open unit $n$-cube $\II^n = (-1,1)^n$. Consider $K$ as a piecewise linear stratification in which each open simplex is a stratum, and pass to its ambient extension $K^+$ in $\II^n$ (whose strata are open simplices of $K$, as well as the complement $\II^n \setminus K$). We construct an open $n$-mesh $M$ that refines $K^+$. First, restrict the projection $p_n : \II^n \to \II^{n-1}$ to $K$, obtaining a map $p_n : K \to \II^{n-1}$ and, via \autoref{rmk:intersection-refinement}, construct an intersection refinement $L \into \II^{n-1}$ for it. Observe that the projection $p_n : K^+ \to L^+$ (with codomain the ambient extension of $L$ in $\II^{n-1}$) is `almost' a constructible 1-mesh bundle (but generally, fails to be a stratified map): indeed, if we define $K^+_L$ to be the refinement of $K^+$ whose strata are connected components of intersections $p_n\inv(s) \cap r$ where $s$ and $r$ are strata in $L^+$ resp.\ $K^+$ (note, if $r$ is an open simplex in $K$, then this intersection is necessarily connected) then, since $L$ is an intersection refinement, one checks that $p_n : K^+_L \to L^+$ is now a constructible 1-mesh bundle. 

    Now, arguing inductively, we construct an $(n-1)$-mesh $M_{<n}$ which refines $L^+$, consisting of 1-mesh bundles $p_i : (M_i,f_i) \to (M_{i-1},f_{i-1})$ (where $M_i = \II^i$). To construct the $n$-mesh $M$, we augment $M_{<n}$ by a constructible 1-mesh bundle $p_n : (M_n,f_n) \to (M_{n-1},f_{n-1})$ defined as the pullback bundle of the constructible 1-mesh bundle $p_n : K^+_L \to L^+$ along the refinement $(M_{n-1},f_{n-1}) \to L^+$. Since $(M_n,f_n) \to K^+_L \to K^+$ is a refinement, this completes the construction of the mesh $M$ refining $K$.
\end{obs}

\begin{eg}[Piecewise linear stratifications are flat framed] In \autoref{fig:a-piecewise-linear-stratification-and-its-mesh} we depict a simplicial complex $K$ linearly realized in $\lR^2$ (and bounded by the 2-cube $\II^2$) together with a 2-mesh $M$ refining it, as it may be constructed via \autoref{constr:simplicial-meshes-of-polyhedra}.
\begin{figure}[ht]
    \centering
    \def\svgwidth{1\columnwidth}
    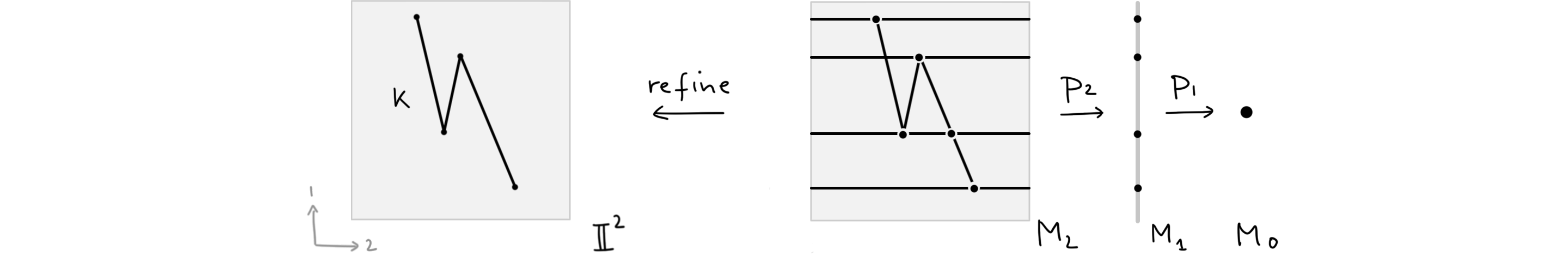

    \caption[A simplicial complex with refining mesh.]{A simplicial complex in $\lR^2$ with refining mesh.}
    \label{fig:a-piecewise-linear-stratification-and-its-mesh}
\end{figure}
\end{eg}

The preceding observations assemble into a proof of \autoref{prop:PL-strat-are-flat-frmd}, which claimed that all piecewise linear stratifications are flat framed, as follows.

\begin{proof}[Proof of \autoref{prop:PL-strat-are-flat-frmd}] Let $(Z,f)$ be a piecewise linear stratification with triangulating simplicial complex $K$. Up to scaling, we may assume $(Z,f)$ and $K$ lie inside the open unit $n$-cube $\II^n = (-1,1)^n$. Using \autoref{constr:simplicial-meshes-of-polyhedra} we can construct an $n$-mesh $M$ refining the ambient extension $K^+$, and thus also refining the ambient extension $(\II^n,f^+)$. This shows that (up to ambient extension) $(Z,f)$ is a flat framed stratification as claimed.
\end{proof}

\begin{obs}[Refining meshes of mesh-supported piecewise linear stratification] \label{rmk:refining-meshes-of-mesh-supported-strat} Let $(Z,f)$ be a piecewise linear stratification. Assume that $(Z,f)$ is mesh-supported (in the sense of \autoref{term:mesh-supported}). Then, as in the preceding proof, we can construct a refining mesh $M$ of the ambient extension $(\II^n,f^+)$ of $(Z,f)$ to the open $n$-cube $\II^n$. Now, we may restrict the $n$-mesh $M$ (and its tower of constructible 1-mesh bundles $p_i : (M_i,f_i) \to (M_{i-1},f_{i-1})$) to the constructible substratification $(M^Z_n, f^Z_n)$ of $(M_n,f_n)$ of strata lying in the subspace $Z$ of $\II^n$ (in particular, $M^Z_n = Z$ as subspaces in $\II^n$). The assumption that $(Z,f)$ is mesh supported implies that this defines a tower of constructible 1-mesh bundle, and thus an $n$-mesh $M^Z$. This mesh refines $(Z,f)$, which shows that $(Z,f)$ itself is in fact a flat framed stratification (without the need for ambient extension to $\II^n$).
\end{obs}

Finally, the following condition on refining meshes of piecewise linear stratifications will be useful to record. Recall from \autoref{constr:geo_open}, that the classifying mesh $\CMsh T$ of an $n$-truss $T$ is obtained as a constructible substratification of the classifying mesh $\CMsh \overline T$ of its compactification $\overline T$, with the latter mesh being isomorphic to the classifying stratification $\CStr {\overline T}$. An `open simplex' in $\CMsh T$ will refer to an open simplex in the simplicial complex $\abs{\overline T_n}$ (triangulating the stratification $\CMsh {\overline T_n}$) that lies in $\CMsh T_n \into \CMsh \overline T_n$.

\begin{term}[Linear meshes] An $n$-mesh $M$ is called a `linear mesh' if, setting $T = \CMsh M$, there is an isomorphism $\CMsh T \iso M$ that is linear on each open simplex in $\CMsh T$.
\end{term}

\begin{obs}[Linear refining meshes] \label{obs:linear-refining-meshes} Every finite simplicial complex $K$ (and similarly, every piecewise linear stratification) has, up to ambient extension, a linear refining $n$-mesh $M$. Indeed, this follows inductively, if in the inductive step of \autoref{constr:simplicial-meshes-of-polyhedra} we chose $M_{<n}$ to be a linear $(n-1)$-mesh. Tracing through the rest of the construction one checks that $M$ as constructed is necessarily linear as well. As an example, note that the refining 2-mesh in \autoref{fig:a-piecewise-linear-stratification-and-its-mesh} is linear.
\end{obs}

\subsubsecunnum{Flat framed Hauptvermutung} \label{sssec:hauptvermutung}

We next prove that notions of framed stratified homeomorphisms and piecewise linear framed stratified homeomorphism coincide, as recorded in the `flat framed Hauptvermutung'. We begin with a brief recollection of the classical Hauptvermutung (see e.g.\ \cite{ranicki1996hauptvermutung}).
 
\begin{disconj}[Hauptvermutung] Every homeomorphism $P \iso Q$ between polyhedra $P$, $Q$ is homotopic to a PL homeomorphism.
\end{disconj}

\nid The Hauptvermutung, and its related conjectures, have been famously disproven in various ways: not only are there polyhedra with different PL structures (where `PL structure' means `PL homeomorphism class'), but there are also topological manifolds with different PL structures and even topological manifolds that don't admit a PL structures at all (that is, the inclusion of PL manifolds into topological manifolds is not surjective on homeomorphism classes). The failure of the Hauptvermutung may be understood as a symptom of the `wildness' of topological homeomorphism which finds no counterpart in combinatorial topology; `taming' this wildness is possible (making, in particular, the Hauptvermutung a true statement) while usually  rather technically cumbersome (cf. \cite{shiota2014minimal} \cite{shiota2013minimal}). Our goal in this section will be to show that a `flat framed' variation of the Hauptvermutung holds naturally in the setting of flat framed stratifications. An adequate classical analog to this statement (for the case of stratifications embedded in $\lR^n$) may be phrased as follows. Note that we take a `stratified homotopy' to mean a homotopy of stratified maps that is constant on entrance path posets.

\begin{disconj}[Ambient stratified Hauptvermutung] Every stratified homeomorphism between (mesh-supported) piecewise linear stratifications is stratified homotopic to a piecewise linear stratified homeomorphism.
\end{disconj}

\nid Like its non-ambient counterpart the ambient Hauptvermutung fails to hold in general.\footnote{We are indebted to Mark Powell for outlining the given counter-example to us.}

\begin{proof}[Disproof of the ambient Hauptvermutung] We construct homeomorphic mesh-supported piecewise linear stratification $f$ and $g$ by piecewise linear embeddings $P_0 \into \II^n$ and $P_1 \into \II^n$ of polyhedra (where $\II^n$ is the open unit cube, stratified by $P_i$ and its complement in each case) as follows.

    \textit{Step 1}: Take simplicial complexes $K$ and $L$ that are both PL homotopy 5-tori but not PL homeomorphic.  Then both $K$ and $L$ are homeomorphic to the 5-torus $T^5$. (\textit{References for step 1}: Hsiang and Shaneson showed that there are non-standard PL homotopy 5-tori \cite{hsiang1970fake}. Later, Hsiang and Wall showed that all PL homotopy 5-tori are in fact homeomorphic to the 5-torus \cite{hsiang1969homotopy}.)

    \textit{Step 2}: PL embed both $K$ and $L$ in $\II^n$; the image of the embeddings define the polyhedra $P_0 \into \II^n$ and $Q_0 \into \II^n$. We can assume these embeddings are locally flat. (\textit{References for step 2}: Rourke and Sanderson construct the required embeddings in \cite[Theorem 5.5]{rourke2012introduction}; the following remark says they are locally flat. In the case of homotopy 5-tori the construction requires $n = 10$.)

    \textit{Step 3}: If we enlarge to $n = 12$, then the topological embeddings $P_0 \into \II^n$ and $Q_0 \into \II^n$ become ambient isotopic by a compactly supported isotopy since $P_0$ are $Q_0$ are both the 5-torus up to homeomorphism. (\textit{References for step 3}: This follows from \cite[Thm 4.4.2]{daverman2009embeddings}; see also Corollary 4.4.3 loc.\ cit.  This states homotopic topological embeddings of high codimension into a PL manifold, with 1-LCC images, are compactly supported ambient isotopic.  Here 1-LCC stands for `1-local-co-connected' and is implied by the embeddings being locally flat (see Prop. 1.3.1 loc.\ cit.).)

    \textit{Step 4}: The last time slice of the isotopy provides the topological homeomorphism of pairs.  A PL homeomorphism of pairs cannot exist by initial choice of PL structures.
\end{proof}

In contrast, the `flat framed' variation of the Hauptvermutung is true. Recall, as stated in \autoref{cor:flat-framed-hpt-vmtg}, this claims that any framed stratified homeomorphism between piecewise linear flat framed stratifications is stratified homotopic to some framed piecewise linear homeomorphism. In fact we will see that the homotopy is unique up to contractible choice (among stratified homotopies of framed stratified homeomorphisms).

\begin{proof}[Proof of \autoref{cor:flat-framed-hpt-vmtg}] Let $(Z,f)$ and $(W,g)$ be piecewise linear flat framed stratification that are framed stratified homeomorphic by a framed stratified homeomorphism $F : (Z,f) \iso (W,g)$. By definition of piecewise linearity, $(Z,f)$ and $(W,g)$ have triangulations $K$ resp.\ $L$. Using \autoref{rmk:refining-meshes-of-mesh-supported-strat} construct $n$-meshes $M$ and $N$ that refine the triangulations $K$ resp.\ $L$ (restricted to the subspaces $Z$ resp.\ $W$). In particular, $M$ and $N$ refine $(Z,f)$ resp.\ $(W,g)$. By \autoref{obs:linear-refining-meshes}, we may choose $M$ and $N$ to be linear $n$-meshes, that is, $\CMsh T \iso M$ and $\CMsh S \iso N$ piecewise linearly  where $T = \ETrs M$ and $S = \ETrs N$. Let $Q$ be the coarsest refining mesh of $(Z,f)$, and set $R = \ETrs Q$. Since coarsest refining meshes are compatible with framed stratified homeomorphism (see \autoref{lem:iso-mesh}) it follows that $F : Q \to FQ$ yields the coarsest refining mesh $FQ$ of $(W,g)$. In particular, $Q$ is coarser than $M$, and $FQ$ is coarser than $N$, and thus $R$ is coarser than both $T$ and $S$ by the coarsenings $\ETrs(M \to Q)$ resp.\ $\ETrs(N \to FQ)$. Recall the construction of `classifying mesh coarsenings' (see \autoref{constr:PL-realization-coarsening-general}) and note that these are linear on each open simplex in their domain. This allows us to construct piecewise linear coarsenings $\CrsMsh(T \to R)$ and $\CrsMsh(S \to R)$. We obtain the maps in the upper row of the following diagram
\begin{equation}
\begin{tikzcd}[baseline=(W.base)]
    (M_n,f_n) \arrow[d] \arrow[rrrrd, "\sim", phantom] & \CMsh T_n \arrow[r] \arrow[l] & \CMsh R_n & \CMsh S_n \arrow[r] \arrow[l] & (N_n,g_n) \arrow[d] \\
{(Z,f)} \arrow[rrrr, "F"'] & & & &|[alias=W]| {(W,g)}
\end{tikzcd}.
\end{equation}
This is a zig-zag of piecewise linear homeomorphisms (namely, the outer two arrows) and piecewise linear coarsenings (namely, the inner two arrows). Since all maps are piecewise linear homeomorphisms of underlying spaces, the upper row defines piecewise linear homeomorphism $G : Z \to W$. By construction, the map also induces a stratified homeomorphism $G : (Z,f) \to (W,g)$ which therefore provides a piecewise linear stratified homeomorphism between the two stratifications. 

Finally, note that $F$ and $G$ induce $n$-mesh maps $F, G : Q \to FQ$ on the respective minimal meshes of $(Z,f)$ and $(W,g)$ which are identical on entrance path trusses. Thus (by faithfulness of the entrance path functor construction, see \autoref{prop:faithfulness} and \autoref{prop:faithfulness-crs-deg}) we deduce that $F$ and $G$ are in fact homotopic as $n$-mesh maps (in fact, uniquely so up to contractible choice). It follows that they are homotopic as stratified maps $F, G : (Z,f) \to (W,g)$ as well (in fact, again uniquely so up to contractible choice). This proves that every framed stratified homeomorphism of piecewise linear stratification is homotopic to a piecewise linear framed stratified homeomorphism as claimed.
\end{proof}

\begin{rmk}[Equivalence of all triangulations] Let $(Z,f)$ be a piecewise linear stratification with two different triangulations $K$ and $L$. Running the preceding proof (setting $(W,g) = (Z,f)$ and $F = \id$) shows that all triangulations are piecewise linearly equivalent.
\end{rmk}

\begin{rmk}[Contractibility of framed structure groups] The fact that the Hauptvermutung fails in the nonframed setting but holds in the framed setting is underpinned by the following intuition: while the automorphism groups $\mathrm{Aut_{TOP}}(\lR^n)$ and $\mathrm{Aut_{PL}}(\lR^n)$ differ, the groups of framed automorphisms $\mathrm{Aut_{TOP}^{fr}}(\lR^n)$ and $\mathrm{Aut_{PL}^{fr}}(\lR^n)$ (see \autoref{term:framed-maps}) are both contractible and thus trivially equivalent. In fact, one would expect the same should hold using classical notions of framings, that is, the notions of `framed topological' and `framed PL' homeomorphisms (up to details of their definitions) should coincide.
\end{rmk}

\subsection{Computability} \label{ssec:computability-of-fl-fr-str}

Finally, we briefly address two basic computability questions about the theory of flat framed stratifications. First, we will see that `coarsest refining meshes' may in fact be algorithmically constructed. Secondly, as a corollary we will find that it can always be decided whether or not two framed stratifications are framed stratified homeomorphic or not. The idea is to reduce both problems to combinatorial problems by passing to entrance path posets.

\begin{prop}[Truss coarsenings are confluent] \label{prop:truss-reduction-is-confluent} Let $T$ be a stratified $n$-truss. Any chain of non-identity truss coarsenings of $T$ eventually ends in the same normalized $n$-truss $\NF{T}$, and the truss coarsening $T \to \NF{T}$ is unique.
\end{prop}

\begin{term}[Normal forms] We call $\NF{T}$ the `normal form' of $T$.
\end{term}

\begin{proof} Let $F_1 : T \to T_1$ and $F_2 : T \to T_2$ be non-identity truss coarsenings of $T$. The underlying maps of trusses are coarsenings $\und F_i : \und T \to \und T_i$ of $\und T$. Construct the classifying mesh coarsenings $\CrsMsh \und F_i : \CMsh \und T \to \CMsh \und T_i$ using \autoref{constr:PL-realization-coarsening-general} (up to pulling back $\CMsh \und T_i$ along $\CrsMsh \und F_i$ we will assume that $\CrsMsh \und T_i$ and $\CrsMsh \und T$ all have the same support in $\lR^n$). By \autoref{thm:joins-of-open-meshes}, we may take the join of $\CMsh {\und T_1}$ and $\CMsh {\und T_2}$ obtain the mesh $\CMsh {\und T_1} \vee \CMsh {\und T_2}$, and abbreviate it by $M$. Construct the stratified classifying mesh $(\CMsh \und T, \sff_T)$ of $T$. Note that, since $F_1$ and $F_2$ are truss coarsenings, both $\CMsh {\und T_1}$ and $\CMsh {\und T_2}$ refine $\sff_T$. Thus their join $M$ refines $\sff_T$ as well. Passing to stratified entrance path trusses, this yields a diagrams of truss coarsenings
    \begin{equation}
        \begin{tikzcd}[row sep=0pt, column sep=50pt]
         & T_i \arrow[rd] & \\
        T \arrow[ru, "F_1" description] \arrow[rd, "F_2" description] & & S \\
         & T_2 \arrow[ru] &
        \end{tikzcd}
    \end{equation}
    where $T_i = \ETrs (\CMsh {\und T_2},\sff_T)$ and $S := \ETrs(M,\sff_T)$. Since any chain of truss coarsenings of $T$ must ends eventually in some normalized truss, to above implies all chains must end in the same normalized $n$-truss $\NF{T}$. Note also that the truss coarsening $T \to \NF{T}$ is unique (otherwise, if there were two truss coarsening $F_1$ and $F_2$, repeat the above steps to contradict that $\NF{T}$ is normalized).
\end{proof}

\begin{rmk}[Mesh joins as truss pushouts] While the preceding proof of \autoref{prop:truss-reduction-is-confluent} crucially relies on our earlier construction of mesh joins, one may also prove the statement in purely combinatorial terms. Namely, the constructed diagram of truss coarsening in the proof of \autoref{prop:truss-reduction-is-confluent} is a pushout of the $n$-truss coarsenings $\und F_1 : \und T \to \und T_1$ and $\und F_2 : \und T \to \und T_2$ in the category of $n$-trusses and their maps, and these pushouts can be shown to exist by purely combinatorial arguments. The construction of pushouts of truss coarsenings may be regarded as a combinatorial counterpart to the construction of joins of meshes (see \autoref{rmk:joins-as-pushouts}).
\end{rmk}

\begin{obs}[Normal forms are computable] Given a stratified truss $T$, there is a finite set of label preserving surjective truss maps $F : T \to S$ from $T$ to another stratified truss $S$. Searching this set for truss coarsenings, and selecting the truss coarsening with smallest codomain, provides an algorithm to compute the normal form coarsening $T \to \NF T$.\footnote{Note that there are better algorithms than the brute-force search described here.}
\end{obs}

\nid Recall \autoref{cor:coarsest-mesh-comp}, which stated that coarsest refining meshes of flat $n$-framed stratifications $(Z,f)$ can be algorithmically computed.

\begin{proof}[Proof of \autoref{cor:coarsest-mesh-comp}] Let $(Z,f)$ be a flat framed stratification. Since it is flat framed it must have some refining mesh $M$ (we exclude the non-constructive existence of refining meshes). Compute its entrance path truss $\ETrs (M,f)$ and the resulting normal form $\NF{\ETrs (M,f)}$. This yields a coarsening of $M$ (determined by the coarsening $\ETrs (M,f) \to \NF{\ETrs (M,f)}$ on entrance path posets), that coarsens $M$ to the coarsest refining mesh of $(Z,f)$.
\end{proof}

\nid Immediately, we deduce that the question of framed stratified homeomorphism between flat framed stratifications is decidable, as stated \autoref{cor:decidability-of-iso}.

\begin{proof}[Proof of \autoref{cor:decidability-of-iso}] Given flat framed stratifications $(Z,f)$ and $(W,g)$ construct their coarsest refining meshes $M$ and $N$. $(Z,f)$ and $(W,g)$ are framed stratified homeomorphic if and only if the corresponding stratified trusses $\ETrs(M,f)$ and $\ETrs(N,g)$ are related by a label preserving balanced isomorphism (note it must be unique if it exists). Since the latter isomorphism problem can be algorithmically solved (e.g., by brute force) so can the former homeomorphism problem.
\end{proof}

\nid The above results establish the `computational tractability of flat framed stratifications'. Note that the existence of coarsest refining meshes (which, conceptually, is `opposite' to the classical quest for constructing mutually refining triangulations in the classical Hauptvermutung) plays a most fundamental role in this story.

\section{Towards transversality, tangles, singularities, and stability} \label{sec:looking-ahead}

In this final section, among a large pool of further topics that deserve discussion, we select and briefly discuss two future directions in the program of framed combinatorial topology, which we plan to address in subsequent work.
\begin{enumerate}

    \item In \autoref{ssec:transversality} we will describe a class of flat framed stratifications whose strata are `transversal'. We call such stratifications manifold diagrams. Manifold diagrams solve the long-standing problem of generalizing string diagrams and surface diagrams to higher dimensions.

    \item In \autoref{ssec:tangles} we describe a class of flat framed stratification that model tame tangles. We discuss how these tangles may be `perturbed', and how `perturbation-stable' tangles singularities provide a combinatorial model for singularity theory.

\end{enumerate}
\nid In \autoref{ssec:problems} we mention four specific problem for potential future research.

\begin{rmk}[Geometric higher category theory] 
We also emphasize that there is (at least) one central omission in our discussion of future work here. One may of course `globalize' our discussion of flat framed structures to more general (non-flat) framed combinatorial spaces and stratified spaces. This direction of investigation provides novel geometric models for various objects endowed with `directed higher structure'. The fundamental idea is that we can think of truss blocks and presheaves on blocks in a categorical manner: $k$-blocks are shapes of $k$-morphisms, and presheaves on blocks (or more precisely `sheaves' in an appropriate sense) can be used to model higher categories. This model of higher categories has some striking differences from other presheaf models \cite{leinster2004higher} \cite{cheng2004higher} \cite{bergner2020survey}. Traditional `filler' or `contractibility' conditions are absent in this model, and moreover the model allows us to work with `small examples' of higher categories. For instance, the block model provides an $n$-category that has a single generating 2-endmorphism on an object: higher coherence data is encoded in the geometry of pastings of cell diagrams. At the same time, due to the tight correspondence of combinatorics and stratified topology built into our model, this approach allows the formalization of several various folklore intuitions about the interaction of higher category theory and stratified topology.
\end{rmk}

\subsection{Transversality and manifold diagrams} \label{ssec:transversality} 
\begin{figure}[ht]
    \centering
    \def\svgwidth{1\columnwidth}
    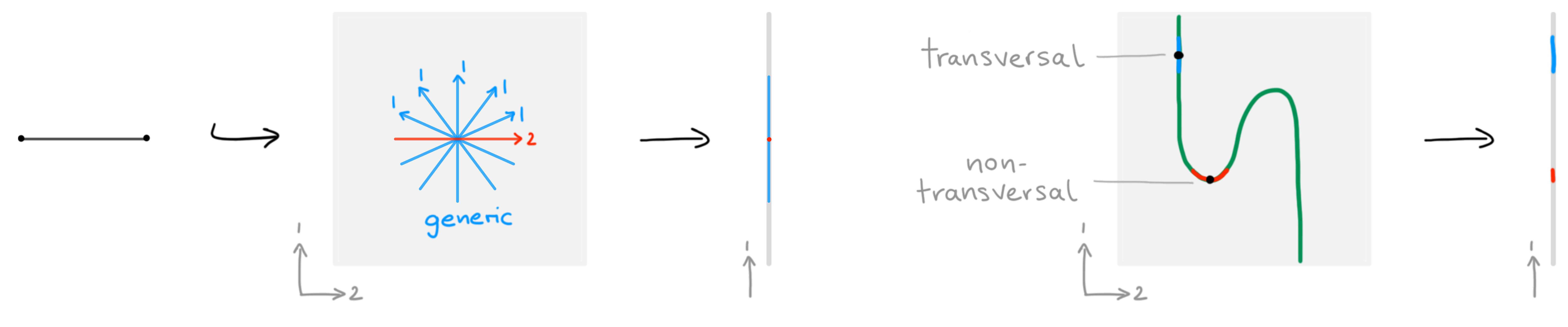

    \caption[Genericity of inherited framings of strata.]{Genericity of framings of 1-strata in $\lR^2$.}
    \label{fig:genericity-of-framings-of-strata-in-flat-framed-stratification}
\end{figure}

Strata in flat $n$-framed stratifications may inherit different types of framings by restricting the ambient flat framing of $\lR^n$ to a given stratum (the resulting framing may of course include `singularities'). For an unordered 1-simplex, consider linear embeddings of the 1-simplex in $\lR^2$ as shown in \autoref{fig:genericity-of-framings-of-strata-in-flat-framed-stratification} on the left. Requiring such an embedding to be a framed realization uniquely determines a 2-embedded frame on the 1-simplex (see \autoref{defn:framed-real-emb}). While there are two 2-embedded frames of the 1-simplex (with frame label `1' resp.\ `2' as indicated), only one of them is `generic' in the sense that it determines an open dense subset in the space of all linear embeddings.
This genericity may be expressed in more familiar terms as a \emph{transversality} condition: a linear embedding $e : \abs{S} \into \lR^2$ of a 1-simplex $S$ is `transversal' if it is transversal to the fibers of the projection $\pi_1 : \lR^2 \to \lR^1$. This transversality condition can alternatively be phrased as requiring $\pi_1 \circ e : \abs{S} \to \lR^1$ to be a homeomorphism onto its image.

The idea of transversality applies more generally as follows. We may say a 1-stratum in $\lR^2$ is `transversal' at a point if the projection $\pi_1 : \lR^2 \to \lR^1$ restricts, in a small neighborhood of the point in the stratum, to a homeomorphism onto its image. We illustrate this on the right in \autoref{fig:genericity-of-framings-of-strata-in-flat-framed-stratification} for two points in a 1-stratum; at one point the stratum is transversal while at the other it is not. More generally, we call a $k$-dimensional stratum $s$ in $\lR^n$ `transversal at a point' $p \in s$ if the projection $\pi_{>k} : \lR^n \to \lR^k$ restricts to a homeomorphism from a neighborhood of $p$ in the stratum $s$ onto its image.

\begin{figure}[ht]
    \centering
    \def\svgwidth{1\columnwidth}
    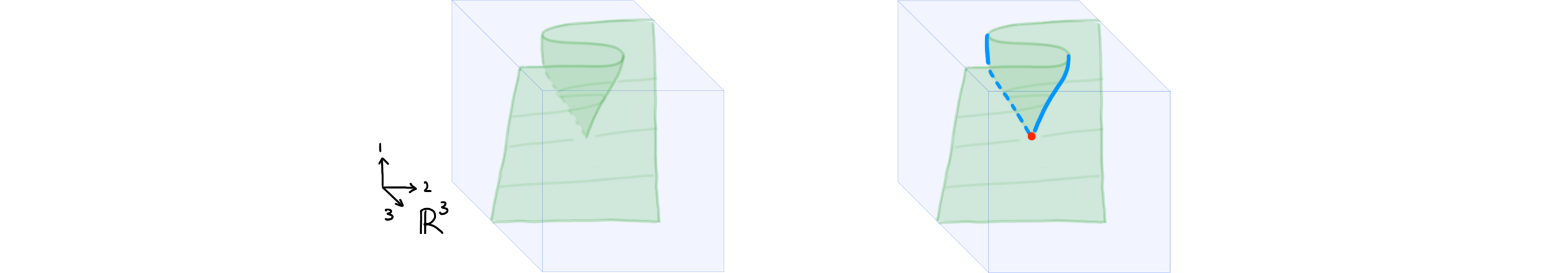

    \caption{A non-transversal stratum and its transversal refinement.}
    \label{fig:the-non-genericity-of-the-cusp-and-a-generic-refinement}
\end{figure}

In \autoref{fig:the-non-genericity-of-the-cusp-and-a-generic-refinement} we illustrate two flat $3$-framed stratifications: the left one (consisting of the green 2-dimensional stratum and the two connected components of its complement) is refined by the right one. While the right stratification is transversal at all points, the left stratification fails to be transversal. Indeed, the above transversality condition is not satisfied at any point that, after passing to the refinement, lies in either the 1- or 0-dimensional strata (indicated in blue and red).

A \emph{manifold $n$-diagram} is a flat framed stratification of the open unit $n$-cube $\II^n$ in which every stratum is transversal, i.e.\ transversal at all points.\footnote{\label{foot:conical-mfld-diag} Technically, we also want manifold $n$-diagrams to be \emph{conical} stratifications, with the additional condition that conical neighborhoods $C \times U$ of strata $s$ can be chosen such that the stratum $s$ is `transversal' to the cone $C$. We defer a more detailed discussion of this condition and the resulting notion of `transversal conicality'.} The flat framed stratification depicted on the right in \autoref{fig:the-non-genericity-of-the-cusp-and-a-generic-refinement} is an example of a manifold $3$-diagram.
This definition of manifold $n$-diagrams in terms is given in terms of topological transversality condition, but manifold diagrams admit a fully combinatorial representation. This representation is based on the correspondence of flat framed stratifications to normalized stratified trusses, together with the observation that transversality of strata in flat framed stratifications can be phrased as a combinatorial condition on strata in the corresponding stratified trusses. Alltogether, we obtain a correspondence between manifold $n$-diagrams and so-called normalized \emph{transversally} stratified trusses. Note further, the combinatorial condition of transversality is computably verifiable, so we can algortihmically recognize transversally stratified trusses among all stratified trusses.

This formalization of manifold $n$-diagrams in all dimensions resolves a long-standing open problem. Previously, stratified topological definitions of `manifold diagrams' had only been given in low dimensions, initially in dimension $2$ \cite{joyal1991geometry}, later in dimension $3$ \cite{hummon2012surface}, \cite{barrett2012gray} (see also \cite{trimble1999}), and in restricted form in dimension $4$ \cite{carter1996diagrammatics} \cite{carter1998knotted}. The importance of manifolds diagrams derives in part from their dual connection to cellular pastings and from their usefulness for `diagrammatic reasoning' via geometric deformations of the underlying stratifications \cite{selinger2010survey}. The approach outlined here provides formal underpinnings for such techniques in all dimensions. In particular, it allows us to formalize manifold $n$-diagram deformations as specific instances of manifold $(n+1)$-diagrams (namely, those that do not have point singularities).

\begin{figure}[ht]
    \centering
    \def\svgwidth{1\columnwidth}
    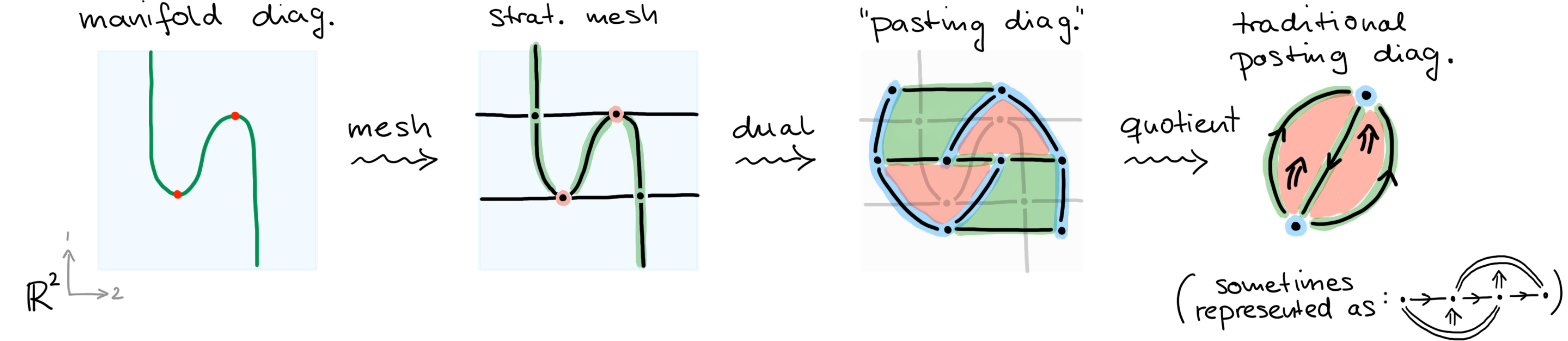

    \caption{The pasting diagram dual to a manifold diagram.}
    \label{fig:the-pasting-diagram-dual-to-the-snake}
\end{figure}

The relationship of `manifold' and `cellular' geometry relies on the dualization operations on meshes and trusses described in previous chapters. Given a manifold $n$-diagram $(Z,f)$, By our results on flat framed stratifications (see \autoref{thm:minimal-meshes}), we may first pass to its coarsest refining open mesh $M$, obtaining a stratified mesh $(M,f)$ (see \autoref{defn:stratified-meshes}). Dualizing the mesh $M$, we obtain a closed mesh $M^\dagger$, and the stratification $f$ induces a stratification $(M^\dagger,f^\dagger)$ on this mesh. (Formally, this stratification is the stratified classifying mesh of the dual $\ETrs(M,f)^\dagger$ of the labeled truss $\ETrs(M,f)$.) We call $(M^\dagger,f^\dagger)$ a `pasting diagrams of cells with degeneracies'. (Note, in particular, that $M^\dagger$ is a framed regular cell complex, see \autoref{term:mesh-complex}.) To obtain a `traditional pasting diagram' (without degeneracies) corresponding to the manifold diagram $(Z,f)$, we may further `quotient' $M^\dagger$ by collapsing to a single cell all the cells in each stratum of $f^\dagger$.
The process is illustrated in \autoref{fig:the-pasting-diagram-dual-to-the-snake}.

\subsection{Tangles, perturbative stability, and singularity theory} \label{ssec:tangles}

\begin{figure}[ht]
    \centering
    \def\svgwidth{1\columnwidth}
    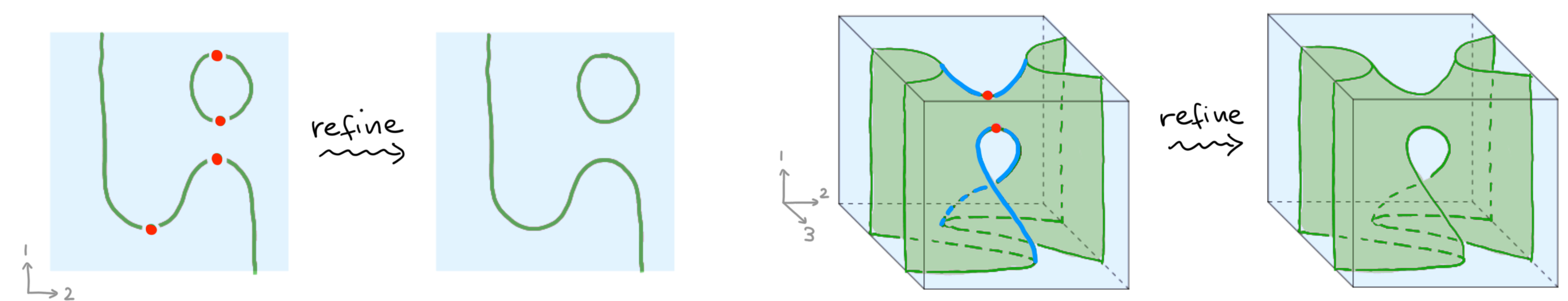

    \caption{Tangles and their refining manifold diagrams.}
    \label{fig:two-transversal-tangles}
\end{figure}

Having discussed transversality and the class of `manifold $n$-diagrams', we now introduce another important and related class of flat framed stratifications, namely `tame tangles'. Tame tangles provide a combinatorializable counterpart to the classical notion of tangles given by manifolds with corners embedded in the closed unit cube. Formally, we define a \emph{tame $k$-tangle in codimension $(n-k)$} to be a flat framed stratification of the open unit $n$-cube $\II^n$ consisting of an open $k$-manifold $M$ in $\II^n$ (with strata being the connected components of $M$ and $\II^n \setminus M$) that can be refined by an $n$-manifold diagram. As in the case of both flat $n$-framed stratifications and manifold $n$-diagrams, tangles may be captured in purely combinatorial terms by their normalized entrance path trusses, which we refer to as \emph{tangle trusses}.\footnote{It is, moreover, reasonable to ask that tangles $M \subset \II^n$ to not just be topological manifolds but combinatorial manifolds with piecewise linear structure induced by their combinatorialization as stratified trusses; this avoids `wildness' phenomena of \cite{cannonshrinking}.} We illustrate two (tame) tangles together with their coarsest refining manifold diagram in \autoref{fig:two-transversal-tangles}.
Note that, to obtain a `classical' tangle as a manifold with corners we may always compactify the open cube stratification to a closed cube stratification. (Formally, this compactification process uses our earlier combinatorial definition of `cubical compactifications', see \autoref{defn:truss-compactifications}; compactification establishes an equivalence between `open' and `closed' tame tangles, and we may work with either perspective.)

\begin{figure}[ht]
    \centering
    \def\svgwidth{1\columnwidth}
    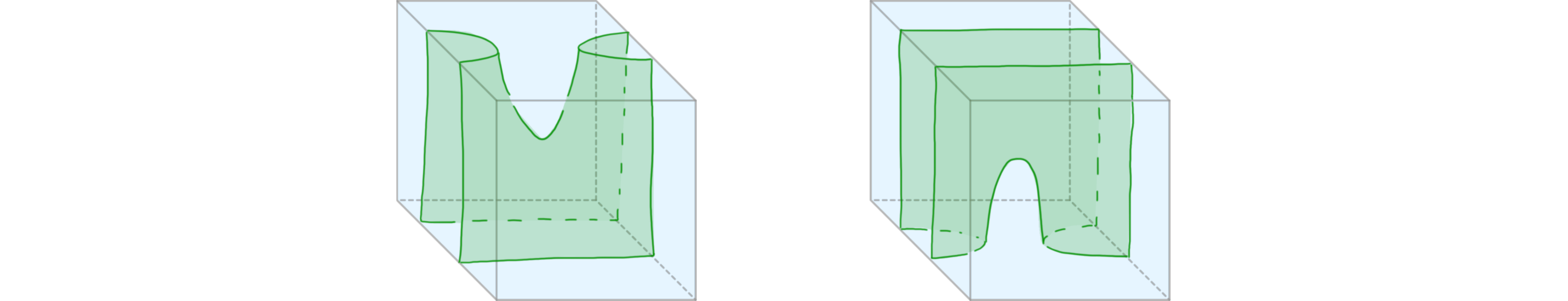

    \caption{Tangle singularities.}
    \label{fig:tangles-singularities}
\end{figure}

If the coarsest refining manifold diagram of a tangle contains a point stratum, then we call a small open neighborhood around this point (on which the tangle restricts to a `subtangle') a `tangle singularity'. For instance, the tangle depicted on the right of \autoref{fig:two-transversal-tangles} has two (red) point singularities in its coarsest refining manifold diagram; passing to the `subtangles' around those points yields the tangle singularities depicted in \autoref{fig:tangles-singularities}.

\begin{figure}[ht]
    \centering
    \def\svgwidth{1\columnwidth}
    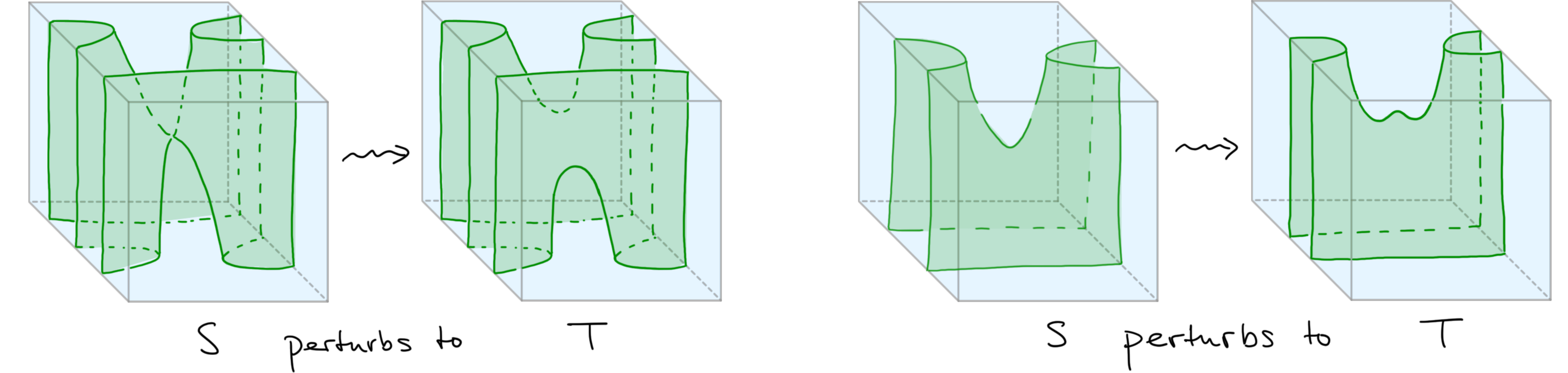

    \caption{Perturbations of tangle singularities.}
    \label{fig:perturbations-of-tangle-singularities}
\end{figure}

Not all tangles are created equal, and the space of all tangles carries a natural stratification. This stratification is, roughly speaking, a measure of the complexity of a tangle's singularities. More precisely, there is a notion of `perturbing a tangle $S$ to a tangle $T$': using the machinery of truss bundles, we may more define this perturbation by considering `bundles of tangles' over the 1-simplex $[1] = (0 \to 1)$ whose generic fiber (over $0$) is the tangle $T$, and whose special fiber is the tangle $S$. We may then say that a tangle singularity $S$ is `perturbation-stable' if, for any perturbation to another tangle $T$, the resulting singularities in $T$ look no less `complex' than the singularities in $S$ that they perturb. (Owing to the combinatorial nature of tangles, the `complexity' of a singularity can be measured in elementary combinatorial terms, for instance by counting the size of the singularity's normalized stratified entrance path truss). We illustrate two instances of perturbations of tangle singularities in \autoref{fig:perturbations-of-tangle-singularities}. In the first case, a `monkey saddle' (see  \cite{milnor1963morse}) is perturbed into two saddles, which are strictly `less complex'; thus the monkey saddle is not perturbation-stable. In contrast, in the second case, a saddle is perturbed into three tangle singularities (two saddles and a maximum) and these have the same complexity as the saddle we started with; indeed, the saddle is a perturbation-stable tangle singularity.

Perturbation-stable tangle singularities will be referred to as `elementary tangle singularities'. The notion of perturbation-stability moreover extends from tangle singularities to general tangles, and the resulting subspace of `perturbation-stable tangles' forms the deepest stratum in the stratification of the space of tangles, In particular, such tangles can only contain elementary singularities. Roughly speaking, the rest of the stratification may then be constructed by defining its entrance paths to be complexity-reducing perturbations. 

\begin{figure}[ht]
    \centering
    \def\svgwidth{1\columnwidth}
    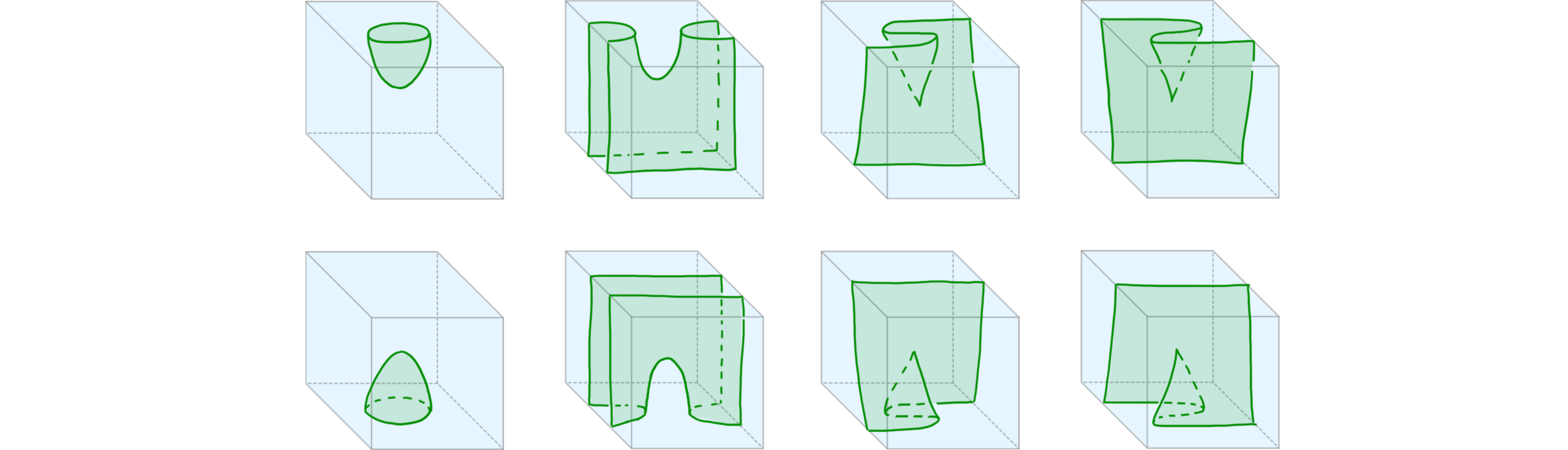

    \caption[Elementary 2-tangle singularities.]{Perturbation stable 2-tangles (in codimension-1).}
    \label{fig:perturbation-stable-2-tangles}
\end{figure}

Elementary tangle singularities turn out to be a particularly interesting class of tangles, and are intimately related to classical singularity theory. In \autoref{fig:perturbation-stable-2-tangles} we depict \emph{all} elementary 2-tangles in codimension-1.\footnote{Here, the fact that we emphasize tangles in codimension-1 reflects the fact that classical singularity theory often studies singularities in ($r$-parameter families of) $n$-variable functions $\lR^n \to \lR$, whose graphs (resp.\ families of graphs) yield hypersurfaces in $\lR^{n+1}$ (resp.\ in $\lR^{n+r+1}$).}
It happens that these 2-tangle singularities exactly recover a collection of classical 2-dimensional singularities, namely  namely the Morse-type singularities (saddles, minima, and maxima) and Cerf-type singularities (`cusps' of 1-parameter families of 1-variable functions). For codimension-1 $k$-tangle singularities with $k > 2$, the parallel with classical singularities initially appears to continue into higher dimensions. Some elementary tangle singularities maybe be associated to known classical singularities, including some of the `elementary singularities' that were described in the work of Thom \cite{thom2018tructural} and Arnold \cite{arnol1975critical}, such as the `swallowtail' singularities (succeeding the cusp in Arnold's list of $A_k$ singularities) whose corresponding elementary tangle singularity was depicted already in \autoref{fig:yang-baxterator-and-the-swallotail}.


However, with increasing number of parameters, the classical machinery of smooth singularity theory eventually breaks down (namely, for parameter dimensions above $5$, see \cite[\S7.6]{poston2014catastrophe}): the  `dimensions of jet spaces outgrow the dimensions of the  structure groups', which causes the smooth equivalence relation of singularities to become too fine, leading one to encounter uncountably many equivalence classes of singularities. In stark contrast, the combinatorial machinery of elementary tangle singularities works in all dimensions, and this leads to a profound discrepancy between the two approaches in higher dimensions. We mention one concrete example of the divergence between the two approaches: among the elementary singularities in Thom's classification one finds the so-called `elliptic umbilic' singularity (in Arnold's classification, it is called the $D^-_4$ singularity); in combinatorial terms, this umbilic singularity is no longer `elementary', in that it may be perturbed into a tangle consisting of less complex tangle singularities. The question of how to classify the combinatorial patterns underlying elementary tangle singularities appears fundamental, and will be revisited in the next section.

\subsection{Problems} \label{ssec:problems}

Wrapping up the discussion of the previous sections, we list four areas of open problems.

\subsubsecunnum{The classification of elementary singularities} In the last section we introduced `elementary', i.e.\ perturbation-stable, tangle singularities. We saw that for 2-tangles in codimension-1 there are exactly eight such singularities. In higher dimensions, elementary tangle singularities quickly become more complicated.

\begin{prob}[Classifying elementary singularities] \label{prob:tangle-sing} For tangles of a given dimension and codimension, can we list all elementary tangle singularities?
\end{prob}

\nid It is moreover reasonable to expect a \emph{finite} list of elementary tangle singularities in each dimension and codimension. As explained, while elementary tangle singularities reproduce singularities from classical singularity theory in low dimensions, in general, there is a discrepancy between the two approaches, and a more precise quantification of this discrepancy would be desirable. One can further consider tangles with structure (for instance, by tangentially framing tangles embedded in $\lR^n$), and similarly ask how the classification of tangle singularities changes in the presence of these additional structure. Via the generalized tangle and cobordism hypotheses (see \cite{baez1995higher} \cite{lurie2009classification}), such classification problems are closely related to understanding `(structured) dualizability laws' of objects in higher category theory. A better understanding of elementary tangle singularities may therefore provide insights both for the study of classical singularity theory as well as higher dualizable structures.

Let us also mention a related problem about the recognition of general (not necessarily perturbation-stable) tangles.
\begin{prob}[Recognizing tangle trusses] Can we recognize tangle trusses among all stratified trusses?
\end{prob}
\nid Even if we require tangles to be combinatorial manifolds (with piecewise linear structure induced by their combinatorialization as stratified trusses), it is not immediately clear how a recognition algorithm could work: such an algorithm would require checking the sphericality of certain links of strata in the stratified truss, which, as discussed in \autoref{punch:block-classification}, is undecidable in general. However, links in stratified trusses are better behaved than general simplicial complexes, and thus the unrecognizability of spheres need not necessarily imply the unrecognizability of tangle trusses; for instance, tangle trusses \emph{would} be algorithmically recognizable if it turned out that all links of tangle singularities were \emph{shellable} (see \autoref{term:pure-shellable-thin}). Other algorithms are also conceivable (which, for instance, could exploit a finite classification of elementary tangle singularities).

\subsubsecunnum{The classification of elementary homotopies} Given a manifold $n$-diagram we may `continuously deform' it: namely, a `homotopy' of manifold $n$-diagrams is a manifold $(n+1)$-diagram which itself does not contain any $0$-dimensional strata. Examples of such deformations were given earlier: the `braid' (depicted as the first example in \autoref{fig:flat-framed-stratifications}) is a homotopy that continuously deforms a manifold $2$-diagram consisting of two point singularities by rotating them around each other; one dimension up, we have also seen the `Reidemeister III move' (depicted in \autoref{fig:yang-baxterator-and-the-swallotail}) which continuously deforms an arrangement of three braids.

Homotopies of manifold diagrams encode coherences laws of higher algebraic structures. For instance, the `braid' witnesses the commutativity law of elements in 2nd (and higher) homotopy groups of spaces. In fact, the braid moreover is the simplest homotopy of manifold 2-diagrams in the following sense: other homotopies, such as the `triple braid', can be perturbed slightly into a composite of braids; in contrast, the braid itself cannot be perturbed any further. We call the braid an `elementary homotopy' and, more generally, we say that a homotopy is `elementary' if it cannot be perturbed into simpler homotopies. (The formalization of this notion mirrors our earlier discussion of elementary singularities.) Analogous to the question of classifying elementary tangle singularities, we may now ask the following.

\begin{prob}[Classification of elementary homotopies] Can we list all the elementary homotopies of manifold $n$-diagrams?
\end{prob}

\nid It is tempting to conjecture that the list of elementary homotopies is finite in any given dimension. However, in reality, little is known about elementary homotopies in higher dimensions (and codimension). In low dimension, such elementary deformations have been considered in the context of `knotting surfaces' \cite{carter1997combinatorial} \cite{carter1998knotted}. More generally, elementary homotopies relate to the question of how `Gray categories' (a notion of 3-categories with natural geometric semantics \cite{gordon1995coherence} \cite{hummon2012surface}) can be generalized to higher dimensions (a discussion of dimension 4 can be found in \cite{crans2000braiding} \cite{bar2017data}). 

\subsubsecunnum{Higher Morse theory} Morse theory is a powerful tool in the study of manifolds \cite{milnor1963morse}; more generally but less ubiquitously, stratified Morse theory is a similarly powerful tool in the study of (sufficiently nice) stratified spaces \cite{goresky1988stratified}. Both theories allow for the introduction of invariants associated to manifolds (resp.\ stratified spaces) based on the decomposition of space that is encoded by a (stratified) 1-Morse function. This leads, for instance, to the construction of Morse homology and intersection homology. Classical (stratified) Morse functions are, however, only the `first degree' in a more general picture which studies manifolds and stratified spaces via `generic' functions to $\lR^n$, that is, via `higher Morse functions'. The idea of 2-Morse functions has been largely realized in dimension 2 (see \cite{cerf1970stratification}) and used, for instance, in the classification of extended field theories in dimension $2$ in terms of generators and relations \cite{schommer2009classification}. However, defining $n$-Morse functions for $n>2$ turns out to be a technically difficult task---at least, when attempting to do so in traditional differential terms.

Using the notions of flat $n$-framed stratifications developed here, the idea of stratified $n$-Morse function finds a simple instantiation as follows.  A flat framed stratifications comes with an embedding into $\lR^n$, and, post-composing this embedding with the projections $\lR^n \to \lR^m$, yields functions to $\lR^m$, $m < n$. We may think of the resulting maps as \emph{instances} of `$m$-Morse functions'. Note that increasing $m$ gradually, each $m$-Morse function contains more information than its `predecessor' $(m-1)$-Morse function. While this discussion of higher Morse functions relies on `embedded framings' $M \into \lR^n$, it may be possible to leverage our notion of `partial framings' to obtain a more direct description of such functions (as an illustrative example see the earlier \autoref{fig:the-saddle-singularity}).

\begin{prob}[Formalization and application of $n$-Morse theory] Can we use the combinatorial framework developed here (including the notion of partial framings) to formalize a combinatorial theory of $n$-Morse functions? Can we derive stronger manifold invariants from $n$-Morse decompositions than from 1-Morse decompositions?
\end{prob}

\nid One may expect the answer to the second question to be `yes' for the following heuristic reason. Classical Morse functions `miss' exotic $m$-spheres from a singularity-theoretic point of view, in that exotic spheres may have 1-Morse functions with a single maximum and a single minimum (just as the standard sphere); instead the information of the smooth structure is hidden in the exotic diffeomorphism of fiber spheres. However, the same cannot be true for $n$-Morse functions, $n \geq m$, of exotic $m$-spheres: since now fibers are discrete, no `exotic' automorphisms of fibers are possible. Instead, the exotic structure needs to be recorded purely in terms of (elementary) tangle singularities and homotopies.

\subsubsecunnum{The relation of framed piecewise linear and smooth topology} Our final question concerns a more precise relation of framed piecewise linear and smooth topology. We phrase this in two conjectures as follows. (For consistency with earlier definitions, note that we identify $\lR^n$ with the open unit cube $\II^n$ by some framed homeomorphism.)

\begin{conj}[Framed stratified homeomorphisms implies diffeomorphism] \label{conj:smooth-struct} From any two smooth embeddings $e : M \into \lR^n$ and $e' : M \into \lR^n$ that correspond to the same normalized tangle truss in $\lR^n$, we can produce a diffeomorphism $M \iso M'$.
\end{conj}

\nid A proof of this conjecture could possibly use a `smooth' version of our earlier proof of conservativity of the entrance path truss functors (see \autoref{prop:conservativity}) to produce a smooth `ambient' isomorphism of the respective coarsest refining meshes of $M$ and $M'$, which would then restrict to a diffeomorphism $M \iso M'$.

\begin{conj}[Smooth embedded manifolds are generically tame tangles] \label{conj:smooth-struct-2} Given a compact smooth $k$-manifold $M$, any smooth embedding $e : M \into \lR^n$ has an arbitrarily small perturbation such that the image of the perturbed embedding $e' : M \into \lR^n$ defines a tame $k$-tangle in $\lR^n$.
\end{conj}

\nid In particular, since any smooth manifold can be embedded in $\lR^n$, the two conjectures together imply that any smooth structure on a compact manifold $M$ can be represented combinatorially by a tangle truss. Again, heuristically, this claim is not unexpected since, by arguments similar to those outlined in the previous section, combinatorializable`$n$-Morse embeddings' should be able to detect smooth structures.

If these conjectures hold, the resulting correspondence of `framed combinatorial manifolds $M$' (i.e. tame tangles with underlying manifold $M$) and smooth structures on $M$ would realize a goal similar to that of MacPherson's program of `combinatorial differential manifolds' \cite{macpherson1991combinatorial}, namely, the faithful combinatorial representation of smooth structures, and the ability to work `smoothly' without direct reference to the continuum.

\changelocaltocdepth{1}

\appendix

\chapter{Linear and affine frames} \label{app:frames}

We discuss classical notions of frames and trivializations of both linear and affine spaces, as well as various notions of `generalized' frames and trivializations.

In \autoref{appsec:linear-fr} we will start with the following two fundamental observations about linear frames. Firstly, we will see that orthonormal frames have reformulations both in terms of sequences of linear subspaces (yielding a notion of `indframes') and in terms of sequences of linear projections (yielding a notion of `proframes'); importantly, while `orthonormality' requires euclidean structure on the vector spaces, the definitions of indframes and proframes do not, and will therefore allow us to define a notion of `orthoequivalence' for general linear frames. Secondly, the notion of linear frames itself generalizes as follows: instead of considering trivializations by isomorphisms $V \toiso \lR^m$, we may also consider projections $V \epi \lR^k$ (yielding a notion of `partial trivialization') or subspaces $V \into \lR^n$ (yielding a notion of `embedded trivialization'), or, yet more generally, general linear maps $V \to \lR^n$ (yielding a notion of `embedded partial trivialization').

    In \autoref{appsec:linear-aff} we will then see that notions of linear frames on vector spaces, as well as their generalizations, can be adapted to notions of `affine frames' on affine spaces. This perspective will allow us to describe our earlier definition of combinatorial frames in purely linear algebraic terms.

\section{Linear frames} \label{appsec:linear-fr}

\subsecunnum{Linear trivializations, frames, indframes, and proframes}

We introduce notions of `linear indframes' and `linear proframes'. For euclidean vector spaces, both are equivalent to the ordinary notion of `orthonormal frames'. In the case of general vector spaces, they provide a notion of `orthoequivalence' of linear frames. The idea underlying both notions (translated appropriately into combinatorial terms) will be pervasive throughout our work, and in particular will play a crucial role in classifying framed combinatorial structures in later chapters. For consistency, let us first (re-)introduce linear trivializations and linear frames.

\begin{defn}[Linear trivializations] A \textbf{linear trivialization} of an $m$-dimensional vector space $V$ is a linear isomorphism $V \toiso \lR^m$.
\end{defn}

\nid Preimages of the standard basis vectors $e_i \in \lR^n$ under the linear trivialization map define an ordered list of `frame vectors' $v_i \in V$. Every linear trivialization therefore determines and is determined by a linear frame in the following sense.

\begin{defn}[Linear frames] A \textbf{linear frame} of an $m$-dimensional vector space $V$ is an ordered list of linearly independent vectors $\{v_1, v_2, \ldots, v_m\} \subset V$.
\end{defn}

We now want to compare the structure of linear trivializations (and equivalently of linear frames) on vector spaces to the following two structures.

\begin{defn}[Linear indframes] A \textbf{linear indframe} of a vector space $V = V_m$ is a sequence of inclusions of vector spaces $V_i$, with $\dim(V_i) = i$:
\[
    0 = V_0 \into V_1 \into V_2 \into \cdots \into V_{m-1} \into V_m = V. \qedhere
\]
\end{defn}

\begin{defn}[Linear proframes] \label{defn:linear-proframes} A \textbf{linear proframe} of a vector space $V = V^m$ is a sequence of projections of vector spaces $V^i$, with $\dim(V^{i}) = i$:
\[
    V = V^m \epi V^{m-1} \epi V^{m-2} \epi \cdots \epi V^{1} \epi V^0 = 0. \qedhere
\]
\end{defn}

\begin{obs}[Equivalence of indframes and proframes] \label{obs:linear-eqv-proframe-indframe} Note that linear indframes and proframes define the same structure on a vector space. For a linear indframe $\{V_i \into V_{i+1}\}_{0 \leq i < m}$ on $V$, the corresponding proframe is determined by the cokernels of the sequence of inclusions into the total vector space: $(V \epi V^{m-i}) = \coker(V_i \into V)$.  Conversely, for a linear proframe $\{V^i \epi V^{i-1}\}_{0 < i \leq m}$ on $V$, the corresponding indframe is determined by the kernels of the sequence of projections from the total vector space: $(V_i \into V) := \ker (V \epi V^{m-i})$.
\end{obs}

\nid There are two important standard instances of indframes and proframes.

\begin{term}[The standard euclidean indframe] \label{defn:standard-indframe} The `standard euclidean indframe' of $\lR^n$ is the sequence of subspace inclusions
\[
    0 \into \lR \into \lR^2  \into \cdots \into \lR^{n-1} \into \lR^n
\]
where $\lR^{i-1} \into \lR^i$ includes into the last $i-1$ coordinates of $\lR^i$.
\end{term}

\begin{term}[The standard euclidean proframe] The `standard euclidean proframe' of $\lR^n$ is the sequence of projections
\[
    \lR^n \epi \lR^{n-1} \epi \lR^{n-2} \epi \cdots \epi \lR^{1} \epi \lR^0
\]
where $\lR^i \to \lR^{i-1}$ forgets the last component of vectors in $\lR^i$.
\end{term}

\nid Observe that the complement of the image of each standard indframe inclusion $\lR^{i-1} \into \lR^{i}$ has two components $\lR^{i} \setminus \lR^{i-1} = \eps^-_i \sqcup \eps^+_i$, where the `negative' component $\eps^-_i$ resp.\ `positive' component $\eps^+_i$ consist of points with $i$th negative resp.\ positive coordinate. We call the assignment of positive and negative signs to those components an `orientation structure' on the standard indframe, and more generally introduce the following notions.

\begin{term}[Oriented indframe] An `oriented indframe' on a vector space $V$ is an indframe $\{V_i \into V_{i+1}\}$ together with an association of signs to the connected components of the complement of the image of each inclusion: $V_i \setminus V_{i-1} = \nu^-_i \sqcup \nu^+_i$. 
\end{term}

\nid An orientation structure on an indframe is equivalent to requiring each $V_i$ to be an oriented vector space.

Observe similarly that the fiber $\pi_i\inv(0)$ over $0 \in \lR^{i-1}$ of each projection $\pi_i : \lR^i \to \lR^{i-1}$ is $\lR$, and thus, $\pi_i\inv(0)\setminus 0$ has again a `negative' component $\eps^i_- = \lR_{<0} \subset \pi_i\inv(0)$ and a `positive' component $\eps^i_+ = \lR_{>0} \subset \pi_i\inv(0)$. We call this an `orientation structure' on the standard proframe, and more generally have the following notion.

\begin{term}[Oriented proframe] An `oriented proframe' on a vector space $V$ is a proframe $\{p_i : V^i \epi V^{i-1}\}$ together with an association of signs to the connected components of the complements $p_i\inv(0) \setminus 0 = \nu^i_- \sqcup \nu^i_+$.
\end{term}

\nid An orientation structure on a proframe is equivalent to requiring each $V^i$ to be an oriented vector space.  Our earlier correspondence of indframes $\{V^i \into V^{i-1}\}$ and proframes $\{p_i : V^i \epi V^{i-1}\}$ on an $m$-dimensional vector space $V$ extends to the oriented case: writing $p_{>i}$ for the composite projection $p_{i+1} \circ ... \circ p_{n-1} \circ p_m : V \epi V^i$, an orientation structure on the indframe determines an orientation structure on the corresponding proframe, and vice versa, by setting $\nu^i_\pm = p_{>i} (\nu^\pm_i)$, and conversely $\nu^\pm_i = p_{>i}\inv (\nu^i_\pm)$.

Any indframe on a vector space $V$ can be obtained by `pulling back' the standard indframe of $\lR^n$ along a trivialization $V \toiso \lR^m$. This similarly applies to the standard proframe in $\lR^n$. With a view towards later generalizations of our discussion here, it will be useful to formalize this using the following more general constructions.

\begin{term}[Pullback sequences] \label{term:pb-seq} Given a sequence of $j$ vector space inclusion $\{W_i \into W_{i+1}\}$ and a map $F : V \to W_j$ as shown below, then we obtain a `pullback' sequence of vector space inclusions $\{V_i \into V_{i+1}\}$ by iterated pullback as shown:
\begin{equation}
\begin{tikzcd}[baseline=(W.base)]
0 \arrow[r, dashed, hook] \arrow[dr] & V_0 \arrow[dr, phantom, "\lrcorner" , very near start, color=black]
\arrow[r, dashed, hook] \arrow[d, dashed] &
V_1 \arrow[dr, phantom, "\lrcorner" , very near start, color=black]
\arrow[r, dashed, hook] \arrow[d, dashed] &
V_{2} \arrow[dr, phantom, "\lrcorner" , very near start, color=black]
\arrow[r, dashed, hook] \arrow[d, dashed] & \cdots \arrow[r, dashed, hook] \arrow[d, "\cdots", phantom] &
V_{j-1} \arrow[dr, phantom, "\lrcorner" , very near start, color=black]
\arrow[r, dashed, hook] \arrow[d, dashed] & V_j = V \arrow[d,"F"] \\
& 0 \arrow[r, hook] & W_1 \arrow[r, hook] & W_2 \arrow[r, hook] & \cdots \arrow[r, hook] & W_{j-1} \arrow[r] & |[alias=W]| W_j
\end{tikzcd} . \qedhere
\end{equation}
\end{term}

\begin{term}[Restriction sequences] \label{term:fac-seq} Given a sequence of vector space projections $\{W^i \epi W^{i-1}\}$ and a map $F : V \to W_j$ as shown below, then we obtain a `restriction sequence' of vector space projections $\{V^{i+1} \epi V^i\}$ constructed by iterated image factorization as shown
\begin{equation}
\begin{tikzcd}[baseline=(W.base), column sep=20pt, row sep=20pt]
    V \arrow[r, two heads, dashed] \arrow[dr,"F"'] & V^j \arrow[r, two heads, dashed] \arrow[d, dashed, hook] & V^{j-1} \arrow[r, two heads, dashed] \arrow[d, dashed, hook] & V^{j-2} \arrow[r, two heads, dashed] \arrow[d, dashed, hook] & \cdots \arrow[r, two heads, dashed] \arrow[d, "\cdots", phantom] & V^1 \arrow[r, two heads, dashed] \arrow[d, dashed, hook] & V^0 = 0 \arrow[d, dashed, hook] \\
 & W^j \arrow[r, two heads] & W^{j-1} \arrow[r, two heads] & W^{j-2} \arrow[r, two heads] & \cdots \arrow[r, two heads] & W^1 \arrow[r, two heads] & |[alias=W]| 0
\end{tikzcd} . \qedhere
\end{equation}
\end{term}

\begin{obs}[Trivializations induce oriented indframes and proframes] A trivialization $F : V \toiso \lR^n$ induces both an oriented indframe and a oriented proframe on $V$ as follows. The indframe on $V$ is defined as the pullback sequence of the standard indframe of $\lR^n$ by $F$. A proframe on $V$ is defined as the restriction sequence of the standard proframe of $\lR^n$ by $F$. Note that the resulting proframe corresponds to the resulting indframe on $V$. Thus, to endow both with orientation it suffices to endow the indframe with an orientation $\nu$, which is canonically determined by requiring $F(\nu^\pm_i) = \eps^\pm_i$.
\end{obs}

Note well that different linear trivializations (and thus linear frames) on $V$ may induce the same oriented indframe resp.\ proframe. This leads to the following equivalence relation on trivializations (or frames).

\begin{term}[Orthoequivalence] We say two linear trivializations (or, equivalently, two linear frames) on $V$ are `orthoequivalent' if they induce the same indframe (and thus the same proframe).
\end{term}

\nid The name of this equivalence relation is derived from the following observation.

\begin{obs}[Orthonormal frames, oriented indframes, and oriented proframes are equivalent] \label{obs:ortho-corr} If $V$ is an $m$-dimensional euclidean vector space, then an `orthonormal frame' $(v_1,v_2, ..., v_m)$ is a linear frame consisting of orthonormal vectors. Orthonormal frames determine exactly those trivializations $V \toiso \lR^n$ which are isometries. Note that any oriented indframe (and similarly, any oriented proframe) on $V$ is induced by exactly one isometry $F : V \toiso \lR^n$, namely the one determined by the property that $F(\nu^\pm_i) = \eps^\pm_i$. Thus, orthonomal frames (resp.\ isometric trivializations) are in correspondence with oriented indframes and proframes.
\end{obs}

\begin{figure}[ht]
    \centering
    \def\svgwidth{1\columnwidth}
    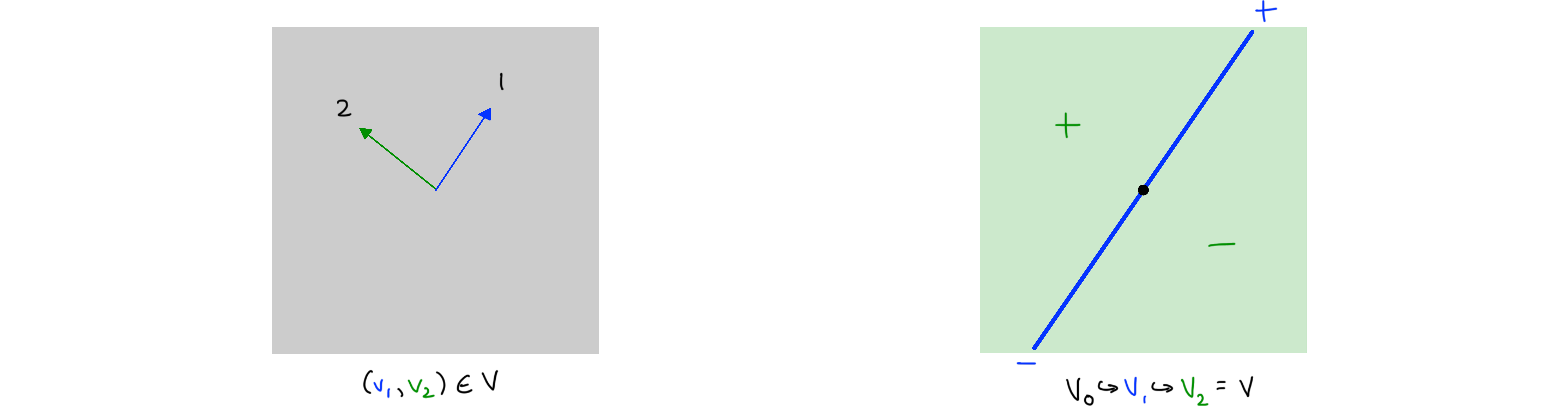

    \caption[An orthonormal frame and indframe]{An orthonormal frame and its corresponding indframe.}
    \label{fig:indframed-and-ortho-frame}
\end{figure}

\nid In other words, each orthoequivalence classes of trivializations of an euclidean vector spaces has a unique orthonormal representative. In \autoref{fig:indframed-and-ortho-frame} we depict an orthonormal frame $(v_1,v_2)$ in an euclidean vector space $V$, together with its corresponding oriented indframe (via \autoref{obs:ortho-corr}).

The preceding discussion suggests the following.

\begin{rmk}[Orthoequivalence as generalization of orthonormality] \label{punch:orthoequivalence} The notion of `trivializations up to orthoequivalence' can be regarded as a generalization of the notion of `orthonormal frames'. In the absence of euclidean structure, we may work with the former structure in place of the latter structure.
\end{rmk}

\subsecunnum{Partial and embedded trivializations and frames}

We now generalize the notions of linear trivializations and frames. Namely, instead of considering linear isomorphisms $V \toiso \lR^m$, we will consider projections $V \epi \lR^k$, subspaces $V \into \lR^n$, and general linear maps $V \to \lR^n$. This leads to a more subtle interplay of `trivializations' and `frames'. The ideas discussed here (translated into appropriate combinatorial form), especially the notion of `embedded' frames, will later be crucial to relate and `glue' frames of objects of different dimensions.

We start with the case of `partial' trivializations and frames.

\begin{defn}[Linear partial trivializations] A \textbf{linear $k$-partial trivialization} of a vector space $V$ is a linear projection $V \epi \lR^k$.
\end{defn}

\begin{defn}[Linear partial frames] A \textbf{linear $k$-partial frame} of a vector space $V$ is an ordered list $\{v_1, v_2, \ldots, v_k\}$ of $k$ linearly independent vectors in $V$.
\end{defn}

\begin{rmk}[Relation of partial frames and trivializations] \label{rmk:relation-partial} If $V$ is euclidean then linear partial trivializations and frames define the same structure on $V$. Indeed, any linear projection $V \epi \lR^k$ has a canonical splitting $\lR^k \into V$ (obtained as the inverse of the restriction of $V \epi \lR^k$ to its orthogonal kernel complement), whose images on standard basis vector $e_i$ determine a linear $k$-partial frame. Conversely, a linear $k$-partial frame $\{v_1, v_2, \ldots, v_k\}$ determines an inclusion $\lR^k \into V$, mapping $v_i$ to $e_i$, which splits to give a $k$-partial trivialization $V \epi \lR^k$. However, in the absence of euclidean structure, we say that a linear $k$-partial frame $\{v_1, v_2, \ldots, v_k\}$ is `compatible' with a linear $k$-partial trivialization $V \epi \lR^m$ if the projection is split by the inclusion $\lR^k \into V$ of the frame. In general, however, this compatibility relation does not yield a 1-to-1 correspondence. 
\end{rmk}

Next, instead of projections $V \epi \lR^k$, we may consider an `embedded trivialization' given rather by an inclusion $V \into \lR^n$.

\begin{defn}[Linear embedded trivialization] A \textbf{linear $n$-embedded trivialization} of an $m$-dimensional vector space $V$ is an inclusion $V \into \lR^n$.
\end{defn}

\begin{defn}[Linear embedded frame] A \textbf{linear $n$-embedded frame} of an $m$-dimensional vector space $V$ is an ordered list of $n$ vectors, $\{v_1, v_2, \ldots, v_n\} \subset V$, exactly $m$ of which are nonzero, and such that all nonzero vector are linearly independent.
\end{defn}

\begin{rmk}[Relation of embedded frames and trivializations] \label{rmk:relation-emb} As before we can translate embedded frames into embedded trivializations by mapping frame vectors to standard basis vectors: that is, every linear $n$-embedded frame $\{v_1, v_2, \ldots, v_n\} \subset V$ induces a linear $n$-embedded trivialization $V \into \lR^n$, defined to map nonzero vectors $v_i$ to $e_i$.  However, the translation is far from bijective (certainly not all linear embedded trivialization are induced by an embedded frame in this way).
\end{rmk}

Finally, there is a common generalization of partial trivializations and embedded trivializations, and an analogous notion of frames, as follows.

\begin{defn}[Linear embedded partial trivializations] A \textbf{linear embedded partial trivialization} of a vector space $V$ is a linear map $V \to \lR^n$.
\end{defn}

\begin{defn}[Linear embedded partial frames] A \textbf{linear $n$-embedded $k$-partial frame} of $V$ is an ordered list of $n$ vectors $\{v_1, v_2, \ldots, v_n\} \subset V$, exactly $k$ of which are nonzero, and such that all nonzero vectors are linearly independent.
\end{defn}

\begin{rmk}[Relation of embedded partial frames and trivializations] \label{rmk:relation-emb-partial}  

As in the previous two cases, the structures do not correspond `one-the-nose'.
\end{rmk}

Importantly, the failure of correspondence of trivialization and frame structures that we observed in the preceding three remarks can be remedied by working with `orthonormal' frames and trivializations up to `orthoequivalence'. As in the case of ordinary linear frames, orthoequivalence will be based on (now generalized) notions of indframes and proframes. Given a partial trivialization $V \epi \lR^k$, an $n$-embedded trivialization $V \into \lR^n$, or an $n$-embedded partial trivialization $V \to \lR^n$, we may form the pullback sequence of the standard euclidean indframe along the trivialization.  Sequences obtained in this way respectively describe notions of `partial', `embedded', and `partial embedded' indframes on $V$.

\begin{defn}[Linear partial, embedded, and embedded partial indframes] A \textbf{linear $k$-partial indframe} on $V$ is a sequence of the form (where $\dim(V_i) = i$):
\[
    0 \into V_{m-k} \into V_{m-k+1} \into \cdots \into V_{m-1} \into V_m = V.
\]
A \textbf{linear $n$-embedded indframe} on $V$ is a sequence of the form (where $\dim(V_{m_i}) = m_i$, and, for each $i$, either $m_{i+1} = {m_i}+1$ or $m_{i+1} = {m_i}$):
\[
    0 = V_0 = V_{m_0} \into V_{m_1} \into V_{m_2} \into \cdots \into V_{m_{n-1}} \into V_{m_n} = V_m = V.
\]
A \textbf{linear $n$-embedded $k$-partial indframe} is a sequence of the form (where $\dim(V_{k_i}) = k_i$, and for each $i$, either $k_{i+1} = {k_i}+1$ or $k_{i+1} = {k_i}$):
\[
    0 \into V_{m-k} = V_{k_0} \into V_{k_1} \into V_{k_2} \into \cdots \into V_{k_{n-1}} \into V_{k_n} = V_m = V.
\]
One further defines `orientations' $\nu$ as before by associating signs $V_{k_i} \setminus V_{k_{i-1}} = \nu^-_{k_i} \sqcup \nu^+_{k_i}$ to the connected components of complements of images of the above inclusions when these complements are non-empty (and similarly in the other two cases of indframes).
\end{defn}

Similarly, given a partial trivialization $V \epi \lR^k$, an $n$-embedded trivialization $V \into \lR^n$, or an $n$-embedded partial trivialization $V \to \lR^n$, we may form the restriction sequence of the standard euclidean proframe along the trivialization. Sequences obtained in this way respectively describe notions of `partial', `embedded' and `partial embedded', proframes on $V$.

\begin{defn}[Linear partial, embedded, and embedded partial proframes] A \textbf{linear $k$-partial proframe} on $V$ is a sequence of the form (where $\dim(V^i) = i$)
\[
    V = V^m \epi V^k \epi V^{k-1} \epi V^{k-2} \epi \cdots \epi V^{0} = 0.
\]
An \textbf{linear $n$-embedded proframe} on $V$ is a sequence of the form (where $\dim(V^{m_i}) = m_i$, and, for each $i$, either $m_{i-1} = {m_i}-1$ or $m_{i-1} = {m_i}$)
\[
    V = V^m = V^{m_n} \epi V^{m_{n-1}} \epi V^{m_{n-2}} \epi \cdots \epi V^{m_1} \epi V^{m_0} = 0.
\]
An \textbf{linear $n$-embedded $k$-partial proframe} is a sequence of the form (where $\dim(V^{k_i}) = k_i$, and, for each $i$, either $k_{i-1} = k_i-1$ or $k_{i-1} =k_i$)
\[
    V = V^m \epi V^k = V^{k_n} \epi V^{k_{n-1}} \epi V^{k_{n-2}} \epi \cdots \epi V^{k_1} \epi V^{k_0} = 0. \qedhere
\]
\end{defn}

\begin{obs}[Correspondence of indframes and proframes] In each of the above three cases (adding the adjectives `partial', `embedded', or `embedded partial'), the notions of indframes and proframes define identical structures on a vector space $V$: that is, one can be constructed from the other as before by taking cokernels and conversely kernels.
\end{obs}

\begin{obs}[Trivializations induce corresponding indframes and proframes] By passing to pullback sequences resp.\ restriction sequences (and adding the adjectives `partial', `embedded', or `embedded partial'), trivializations of a vector space $V$ induce oriented indframes resp.\ oriented proframes of $V$. Given a trivialization its induced oriented indframe corresponds to its induced oriented proframe.
\end{obs}

\begin{term}[Orthoequivalence] Two (partial, embedded, or embedded partial) trivializations are `orthoequivalent' if they induce the same oriented indframe, or equivalently, the same oriented proframe.
\end{term}

\nid Orthoequivalence allows us to relate trivializations and `orthonormal' frames; in the following, a list of vectors is called `orthonormal' if all its nonzero vectors are orthonormal.

\begin{obs}[Orthoequivalence classes of orthonormal embedded frames] \label{obs:ortho-emb-corr} Let $V$ be an $m$-dimensional euclidean vector space. The following auxiliary notion will be helpful: a `standard isometric' $n$-embedded trivialization $V \into \lR^n$ is an $n$-embedded trivialization which is an isometry containing $m$ standard vectors in its image. Note that, orthonormal $n$-embedded frames induce (in the sense of \autoref{rmk:relation-emb}) exactly standard isometric $n$-embedded trivialization, which yields a 1-to-1 correspondence between the two notions. We will show each orthoequivalence classes of an embedded trivialization contains exactly one orthonormal embedded frame.

Consider an oriented $n$-embedded indframe $\{V_{m_i} \into V_{m_{i+1}}\}$ of $V$ with orientation $\nu$. Note there is a unique isometry $F : V \toiso \lR^m$ such that $F(\nu^\pm_{m_i}) = \eps^\pm_{m_i}$. Define an isometric inclusion $\lR^m \into \lR^n$ by mapping basis vectors $e_{m_j}$ to $e_j$ if $m_j > m_{j-1}$. Composing the two maps yields a standard isometric $n$-embedded trivialization $V \into \lR^n$ and thus determines an orthonormal $n$-embedded frame. Thus orthonormal embedded frames (resp.\ their standard isometric embedded trivializations) are in correspondence with oriented embedded indframes and thus with orthoequivalence classes of embedded trivializations as claimed.
\end{obs}

\begin{figure}[ht]
    \centering
    \def\svgwidth{1\columnwidth}
    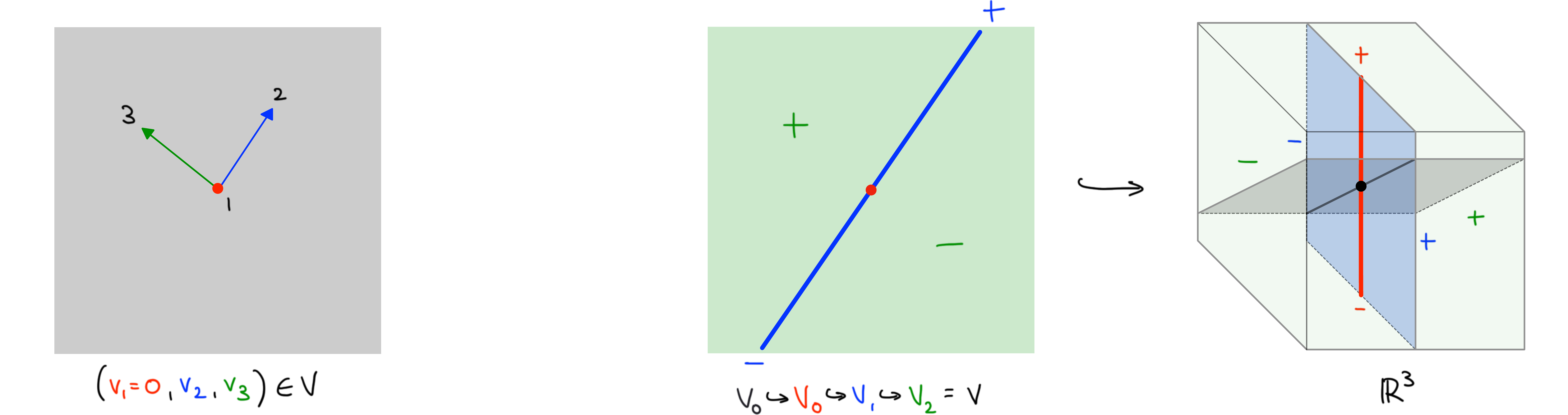

    \caption[An embedded orthonormal frame and indframe]{An $3$-embedded orthonormal frame and its corresponding $3$-embedded indframe.}
    \label{fig:embedded-indframed-and-ortho-frame}
\end{figure}

\nid In \autoref{fig:embedded-indframed-and-ortho-frame} we depict a $3$-embedded orthonormal frame $(v_1,v_2,v_3)$ on a 2-dimensional euclidean vector space $V$, together with its corresponding oriented $3$-embedded indframe (via \autoref{obs:ortho-emb-corr}): note that the embedded indframe determines exactly and orthoequivalence class of embedded trivializations (namely, those embedded trivializations that pull back the oriented standard indframe of $\lR^3$ to the depicted oriented embedded indframe of $V$).

We may further generalize the correspondence to the case of \emph{partial} embedded frames as follows.

\begin{obs}[Orthoequivalence classes of orthonormal embedded partial frames] \label{obs:ortho-emb-part-corr} Let $V$ be a euclidean vector space. The following auxiliary notion will be helpful: a `standard isometric' $n$-embedded $k$-partial trivializations is a subspace $W \into V$ together with a standard isometric embedded trivialization $W \into \lR^n$ in the sense of \autoref{obs:ortho-emb-corr}. Given an orthonormal $n$-embedded $k$-partial frame define $W \into V$ to be the subspace spanned by its nonzero vectors. This canonically splits by a projection $V \epi W$, which maps the orthonormal $n$-embedded $k$-partial frame of $V$ to an orthonormal $n$-embedded frame of $W$. Thus, by \autoref{obs:ortho-emb-corr}, orthonormal $n$-embedded $k$-partial frames are in correspondence with standard isometric $n$-embedded $k$-partial trivializations. Adapting the arguments of \autoref{obs:ortho-emb-corr}, we observe that each orthoequivalence classes of an embedded partial trivialization contains exactly one orthonormal embedded partial frame.
\end{obs}

\nid In other words, orthonormal embedded (partial) frames provide unique representatives orthoequivalence classes of embedded (partial) trivializations of an euclidean vector spaces. Importantly, in the absence of euclidean structure, we do not have access to a notion of `orthonormality' of frames any more, but we still do have access the notion of  `orthoequivalence' of trivializations; this suggests the the following.

\begin{rmk}[Orthoequivalence as generalization of orthonormality] \label{punch:orthoequivalence-emb} We may regard the notion of (embedded partial) `trivializations up to orthoequivalence' as a generalization of the (embedded partial) `orthonormal frames', and work with the former structure in place of the latter structure when euclidean structure is absent.
\end{rmk}

\section{Affine frames} \label{appsec:linear-aff}

\subsecunnum{Affine trivializations, frames, indframes, and proframes}

In the previous sections we discussed the correspondence of frames, indframes, and proframes in the setting of vector spaces. We now discuss how these ideas carry over to the case of affine spaces. This case is important as a parallel to the affine combinatorics of framed simplices.

An affine space $V$ is a space of `points' freely and transitively acted upon by a vector space $\vec V$ (the `associated vector space'); the vectors of $\vec V$ are called `translations'. A map of affine spaces $F : V \to W$ is a continuous map such that for a necessarily unique linear map $\vec F : \vec V \to \vec W$ (the `associated vector space map') we have $\vec F(v - v') = F(v) - F(v')$. Denote the category of affine spaces and affine maps by $\mathsf{Aff}$. This comes with an `associated vector space' functor $\vec - : \mathsf{Aff} \to \mathsf{Vect}$ to the category of vector spaces.

\begin{rmk}[Simplices realize to affine spaces] \label{rmk:simplices-as-affine-sp} Given an $m$-simplex $S$, its geometric realization $\abs{S}$ is the subspace of the free vector space $\lR\avg{S}$ (on the set of vertices of $S$) that consists of affine combinations of the standard basis. Therefore, the realization $\abs{S}$ is contained in an affine hyperplane $V(S)$ of $\lR\avg{S}$, and carries `partial' affine structure (`partial' in the sense that the action by translations in $\vec V(S)$ is partial). Given another affine space $W$, an `affine map' $\abs{S} \to W$ is a map that is the restriction of an affine map $V(S) \to W$.\footnote{In the main text, an `affine map' $\abs{S} \to W$ is, abusing terminology, referred to as a `linear map'.}
\end{rmk}

\begin{notn}[Standard geometric simplices] We denote the geometric realizations of the `standard' $m$-simplex $S$ with vertices $\{0,1,...,m\}$ by $\Delta^m$, and refer to it as the `standard geometric $m$-simplex'.
\end{notn}

\begin{notn}[Space of affine vectors] \label{term:affine-vectors} Given an affine space $V$, the space of affine embeddings $e : \Delta^1 \into V$ of the realized 1-simplex $\Delta^1$ into $V$ is the `space of affine vectors' in $V$, and is denoted by $\hat V$. Note that the space of affine vectors is itself an affine space (it has an action by $\vec V \oplus \vec V$). A `basepoint forgetting' map that takes an affine vector $e : \Delta^1 \into V$ to the translation $e(1) - e(0)$; we denote this map by $\textsf{unbase} : \hat V \to \vec V$. Note that any affine map $F : V \to W$ induces (by postcomposition) a map of affine vectors $\hat F : \hat V \to \hat W$; this induces an `affine vector' functor $\hat - : \mathsf{Aff} \to \mathsf{Aff}$.
\end{notn}

\nid Note that the notation similarly applies to vector spaces $V$ (which, by acting on themselves, are in particular affine spaces with $\vec V \equiv V$).

The notions of frames, indframes, and proframes carry over to the affine case, as follows.

\begin{defn}[Affine trivializations, frames, indframes, and proframes] \label{defn:affine-triv} An \textbf{affine trivialization} resp.\ an \textbf{affine frame} of an affine space $V$ is a linear trivialization resp.\ a linear frame of its associated vector space $\vec V$. Similarly an \textbf{affine indframe} resp.\ an \textbf{affine proframe} is defined as a linear indframe resp.\ proframe of $\vec V$.
\end{defn}

\nid The definitions of affine trivializations, frames, indframes, and proframes can further be translated into `affine terms', i.e.\ expressed in terms of structures on/maps of affine spaces and not of their associated vector spaces. Importantly, we will encounter that the definitions of indframes and proframes behave unequally under this translation. We start with the case of affine trivializations (and, equivalently, frames).

\begin{rmk}[Affine perspective on affine trivializations and frames] \label{obs:affine-frames} Given an affine space $V$ and a trivialization $\vec V \toiso \lR^m$, precomposing with $\textsf{unbase} : \hat V \to \vec V$ yields a map $\hat V \to \lR^m$ that trivializes `vector spaces of affine vectors based at $x$' for any $x \in V$. In particular, preimages of standard vectors $e_i \in \lR^n$ under this map provide a set of frame vectors $v^x_i : \Delta^1 \to V$ at each point $v^x_i(0) = x$ of $V$.
\end{rmk}

\nid Now, to express affine indframes and proframes in affine terms, the following observation needs to be taken into account.

\begin{obs}[Asymmetry of affine projections and subspaces] \label{obs:affine-asymmetry} Given an affine space $V$, and a projection of associated vector spaces $\vec V \epi \vec W$, then this canonically induces an `affine projection' $V \epi W$ where $W$ is constructed as the quotient of $V$ by the action of $\ker(\vec V \epi \vec W)$ (the associated vector space of $W$ is $\vec W$).
    In contrast, given an inclusion $\vec W \into \vec V$ then there is no canonical corresponding injective affine map $W \into V$. In particular, given an affine projection $V \epi W$ there is no canonical choice of `affine kernel' $U \into V$ corresponding to the linear kernel $(\vec U \into \vec V) = \ker(\vec V \epi \vec W)$.
\end{obs}

\nid The observation leads to the following asymmetry in the translation of affine proframes and indframes.

\begin{rmk}[Affine data of affine proframes] \label{obs:affine-proframes} Given an affine space $V$ and a proframe $\vec V = \vec V^m \epi \vec V^{m-1} \epi ... \epi \vec V^0$ on $\vec V$, we obtain a sequence of affine surjective maps
\[
    V = V^m \epi V^{m-1} \epi V^{m-2} \epi ... \epi V^1 \epi V^0 = 0
\]
where the map $V \to V^m$ is the quotient of $V$ induced by the projection $\vec V \to \vec V^m$. Applying the affine vector functor $\hat -$, the above sequence becomes a sequence of affine spaces of affine vectors, which factors through the original proframe on $\vec V$ by the basepoint forgetting maps $\mathsf{unbase} : \hat V^{i} \to \vec V^{i}$.
\end{rmk}

\begin{rmk}[Affine data of affine indframes] Given an affine space $V$ and an indframe $0 = \vec V^0 \into \vec V^{1} \into ... \into \vec V^m = \vec V$ on $\vec V$, this sequence of linear subspaces does \emph{not} canonically induce a sequence of affine injective maps ending in $V$. Instead, we may think `affinely' of the indframe by pulling it back along $\textsf{unbase} : \hat V \to \vec V$ to a filtration of $\hat V$ (which simply `bases' a copy of the indframe at all points of $V$).
\end{rmk}

\nid In other words, affine proframes can be canonically expressed in terms of sequences of affine projections; however, affine indframes cannot be canonically expressed as sequence of affine subspaces. This observation has an analog in affine combinatorics of simplices; indeed, in \autoref{ssec:combinatorial-frames} we saw see that `simplicial subspaces' cannot be expressed canonically as simplicial maps while `simplicial projections' can (namely, in terms of degeneracies).

\subsecunnum{Simplicial trivializations and framed realizations}

Let us now relate the notion of framed realizations of framed simplices to our discussion of (linear and affine) trivializations here. Recall framed realizations (introduced for (partial) framed simplices in \autoref{defn:real-fr-simp} and \autoref{defn:real-fr-simp}, and generalized to the (partial) \emph{embedded} case by \autoref{defn:framed-real-emb} resp.\ \autoref{defn:framed-real-emb-part}) are affine maps from simplices to $\lR^n$ such that simplicial vectors with frame label $i$ are mapped to the $i$th positive component $\eps^+_i$ of the standard indframe of $\lR^n$. As we now explain, framed realizations can be understood as orthoequivalence classes of `simplicial' trivializations (which capture trivialization of affine simplicial space that, up to orthoequivalence, can be obtained from frames of simplicial vectors).

\begin{term}[Simplicial embedded trivialization] \label{term:simp-emb-triv} For an $m$-simplex $S$, an injective affine map $r : \abs{S} \into \lR^n$ (see \autoref{rmk:simplices-as-affine-sp}) is a `simplicial embedded trivialization' if the linear $n$-embedded trivialization $\vec r : \vec V(S) \into \lR^n$ is orthoequivalent to a trivialization of $\vec V(S)$ determined by some $n$-embedded frame $\{\vec v_i\}$ of vectors $v_i$ in $S$ (see \autoref{rmk:relation-emb}). We call the ordered list of $n$ vectors $\{v_i\}$ an `$r$-compatible frame'. Two such maps $r, r' : \abs{S} \into \lR^n$ are `orthoequivalent' if they are orthoequivalent as linear embedded trivializations $\vec r, \vec r' : \vec V(S) \into \lR^n$. If further $m = n$, we speak of `simplicial trivializations' instead.
\end{term}

\begin{obs}[Simplicial embedded trivialization determine embedded frames] \label{obs:simp-emb-triv-determine-emb-fr} Given a simplicial embedded trivialization $r : \abs{S} \into \lR^n$ and an $r$-compatible frame $\{v_i\}$, then exactly $m$ of the vectors $v_i$ in $S$ must be nonzero. The choice of nonzero vectors $v_i$ becomes unique if we require them to form a set of spine vectors under some identification $S \iso [m]$. Moreover, if $v_i$ is a zero vector then it is exchangeable for any other zero vector, and thus choices of zero vectors are redundant: we may replace zero vectors $v_i$ simply by the affine zero $\simpzero$. This results in a unique ordered list $\{v_i\}$ of $n$ vectors in $S$, containing $m$ nonzero vectors forming the spine of $S \iso [m]$, and $(n-m)$ zeros $\simpzero$. This defines a simplicial $n$-embedded framed simplex $(S \iso [m],\cF)$ with $v_i = \cF(j)$ iff $v_i$ is the $j$th nonzero vector in the list $\{v_i\}$.
\end{obs}

\begin{obs}[Framed realizations are simplicial trivialization classes] \label{rmk:emb-fr-real-orthoeq} Given an $n$-embedded framed $m$-simplex $(S \iso [m],\cF)$, note that any framed realization is a simplicial embedded trivialization. In fact, the set of framed realizations of $(S \iso [m],\cF)$ is exactly the orthoequivalence class of simplicial embedded trivialization of $S$ which determine the embedded framed simplex $(S \iso [m],\cF)$. This observation of course specializes to the non-embedded case $n = m$.
\end{obs}

Let us next consider the partial case: again, framed realizations of a partial embedded framed simplex can be understood as describing an orthoequivalence class of partial embedded trivializations of the underlying simplex, as follows.

\begin{term}[Simplicial embedded partial trivialization] \label{term:simp-emb-part-triv} Given an $m$-simplex $S$, we say that an affine map $r : \abs{S} \to \lR^n$ is `simplicial embedded partial trivialization' if it factors by affine maps $\abs{S} \epi \abs{T} \into \lR^n$ such that the first map realizes a degeneracy and the second is a simplicial embedded partial trivialization. Two such maps $r, r' : \abs{S} \to \lR^n$ are `orthoequivalent' if they are orthoequivalent as linear embedded partial trivializations $\vec r, \vec r' : \vec V(S) \to \lR^n$. If $m = n$, we speak of a `simplicial partial trivialization' instead.
\end{term}

\begin{obs}[Simplicial embedded partial trivializations determined frames] \label{prev:simp-emb-part-frame} Any simplicial embedded partial trivialization $r : \abs{S} \to \lR^n$ factors uniquely by maps $\abs{S} \epi \abs{T} \into \lR^n$; since the map $\abs{T} \into \lR^n$ determines an $n$-embedded frame $(T \iso [k], \cF)$, it follows that $r$ determines a $n$-embedded $k$-partial framed $(S \epi T \iso [k], \cF)$.
\end{obs}

\begin{obs}[Framed realizations are simplicial trivialization classes] \label{rmk:emb-part-fr-real-orthoeq} Given an $n$-embedded $k$-partial framed $m$-simplex $(S \epi T \iso [m],\cF)$, note that any framed realization is a simplicial embedded partial trivialization. In fact, the set of framed realizations of $(S \epi T \iso [m],\cF)$ is exactly the orthoequivalence class of simplicial embedded partial trivializations of $S$ which determine embedded partial framed simplex $(S \epi T \iso [m],\cF)$. The observation specializes to the non-embedded case $n = k$.
\end{obs}

\nid In summary, we may understand combinatorial frames in affine algebraic terms as follows.

\begin{rmk}[Combinatorial frames describe simplicial trivialization classes] Combinatorial (embedded partial) frames correspond exactly to orthoequivalence classes of (embedded partial) trivializations of simplices. In light of the role of orthoequivalence (see \autoref{punch:orthoequivalence} and \autoref{punch:orthoequivalence-emb}), they thus provide a combinatorial description of `generalized orthonormal frames' of simplices (realized as affine spaces).
\end{rmk}

\chapter{Stratified topology} \label{app:stratifications}

The notion of `stratified space' (or `singular space') refers to a decomposition of a space into `strata', usually ordered by dimension or `depth'. Frequently, such order is enforced by working with filtrations $X_0 \subset X_1 \subset ... \subset X_{k-1} \subset X_k$ of spaces $X = X_k$ where $X_{i-1}$ is required to be closed in $X_i$. Equivalently, and more concisely, such a filtration may be expressed by a continuous function $f : X \to [k]\op$ (where $[k]\op$ is the poset $(0 \ot 1 \ot ... \ot k)$ topologized such that downward closed subposets are open sets) which allows us to recover $X_i$ as the preimage $f\inv[i]\op$, $i \leq k$. This has been generalized by defining stratifications as continuous maps of spaces to any poset, yielding, for instance, definitions of `$\sS$-filtered spaces' in \cite[\S III.2.2.1]{goresky1988stratified} and of `$P$-stratifications' in \cite[Defn. A.5.1]{lurie2012higher}. Note, however, posets in the domain of such continuous maps may contain information that is unrelated to the decomposition of the underlying space, even when the map is surjective. In this appendix, we develop a notion of stratification which is similarly general, but in which the role of posets faithfully represents topological information about the stratification. The different definitions coincide in many cases, for instance for locally finite stratifications.
We use the following conventions.

\begin{conv}[Convenient spaces]
We assume all spaces are compactly generated. Denote the cartesian closed category of such spaces by $\Top$, and denote the internal hom in this category by $\Map(-,-)$.
\end{conv}

\begin{conv}[Specialization topology] \label{conv:spec-top}
Given a preorder $(X,\to)$ we regard it as a topological space with the \emph{specialization topology}, declaring the open subsets to be those that are downward closed; a subset $U$ is downward closed if $x \leq y$ and $y \in U$ implies that $x \in U$.\footnote{We frequently write the relation $x \leq y$ as $x \to y$, interpreting preorders and posets as categories.}
\end{conv}

\begin{conv}[Specialization order] \label{rmk:specialization-order} Given a topological space $X \in \Top$, we denote by $\Spec X$ its \textit{specialization order}: this is the preorder whose objects are the objects of the underlying set $X$, and whose morphisms $x \to y$ are given whenever $y$ is contained in the closure of $x$. Note that for a poset $P$, $\Spec P = P$.
\end{conv}

\begin{rmk}[Equivalence of finite preorders and finite spaces] The specialization topology provides a functor, which is an adjoint equivalence between finite preorders and finite topological spaces. The inverse functor is given by the specialization order functor $\Spec$.
\end{rmk}

\section{Stratified spaces}

\subsecunnum{Entrance paths, stratifications, and characteristic functions} In this section we introduce the basic notions of `stratified spaces', their `entrance path posets' and their `characteristic functions'.

A robust definition of stratified spaces is obtained by letting the topological decomposition of a space into strata determine the corresponding poset structure, in terms of the existence of so-called entrance paths between strata, as follows.

\begin{defn}[Entrance path] \label{def:ep}
Given a space $X$ and two subspaces $X_r$ and $X_s$, an \textbf{entrance path} from $X_r$ to $X_s$ is a path $p : [0,1] \to X$ with $p(x < 1) \subset X_r$ and $p(1) \in X_s$.
\end{defn}

\nid Here, the path is thought of as `entering' from the former subspace $X_r$ into the latter subspace $X_s$.

\begin{defn}[Formal entrance path] \label{defn:fep}
Given a space $X$ and two subspaces $X_r$ and $X_s$, we say there exists a \textbf{formal entrance path} from $X_r$ to $X_s$, denoted $r \to s$, when the closure of $X_r$ has nonempty intersection with $X_s$.
\end{defn}

\nid In contrast to entrance paths, note that the structure of formal entrance paths is boolean: either there exists a formal entrance path between subspaces or there doesn't. If there is an entrance path from a subspace $X_r$ to a subspace $X_s$ of a space $X$ this implies the existence of a formal entrance path, but the converse need not hold unless additional conditions are imposed.

\begin{term}[Formal entrance path relation of a decomposition] \label{defn:prestrat-and-strat}
Given a decomposition $X = \bigsqcup_{s \in \mathsf{Dec}} X_s$ of a space $X$ into a disjoint collection of subspaces, the `formal entrance path relation' on the set $\mathsf{Dec}$ of subspaces (the `decomposition set') is the relation that has an arrow $r \to s$ exactly when there is a formal entrance path from $X_r$ to $X_s$.
\end{term}

\nid Note that the formal entrance path relation of a decomposition is reflexive, but need not be antisymmetric or transitive.  Stratifications are exactly those decompositions for which this relation has no cycles, that is for which it is a directed acyclic graph.

\begin{defn}[Prestratifications and stratifications] \label{defn:stratification}
A \textbf{prestratification} $(X,f)$ of a topological space $X$ is a disjoint decomposition $f = \{X_s \subset X\}_{s \in \mathsf{Dec}(f)}$ of $X$ into nonempty connected subspaces indexed by a set $\mathsf{Dec}(f)$. The subspaces $X_s$ are called \textbf{strata} of $(X,f)$. A \textbf{stratification} $(X,f)$ is a prestratification such that the formal entrance path relation on the decomposition set $\mathsf{Dec}(f)$ has no cycles.
\end{defn}

\begin{notn}[Shorthand for (pre)stratifications]
\label{notn:abbreviating-stratifications}
We frequently abbreviate a (pre)stratification $(X,f)$ simply by $f$, referring to $f$ as a `(pre)stratification on $X$'. We will often abbreviate a stratum $X_s \subset X$ simply by its index $s \in \mathsf{Dec}(f)$.
\end{notn}

Observe that, given a stratification $(X,f)$, the transitive closure of the formal entrance path relation on the decomposition set $\mathsf{Dec}(f)$ is a partially ordered set, which has an arrow $r \to s$ exactly when there is a chain of formal entrance paths beginning at $r$ and ending at $s$. Note that the equivalence classes of the resulting transitive relation remain exactly the decomposition set $\mathsf{Dec}(f)$. This does not hold true if $(X,f)$ is merely a \textit{pre}stratification, in which case $\mathsf{Dec}(f)$ obtains the structure of a \textit{pre}ordered set.


\begin{defn}[Entrance path preorder and poset] \label{defn:entr}
For a prestratification $(X,f)$, the \textbf{entrance path preorder} $\Entr(f)$ is the decomposition set of the stratification together with the transitive closure of the formal entrance path relation. If $(X,f)$ is a stratification, then we refer to $\Entr(f)$ as the \textbf{entrance path poset} of $(X,f)$.
\end{defn}

\begin{rmk}[Exit paths and the exit path preorder] \label{rmk:entr-vs-exit}
Given a prestratification $(X,f)$, the opposite preorder of the entrance path preorder is called the exit path preorder, and is denoted $\mathsf{Exit}(f) = \Entr(f)\op$.  An exit path from $X_s$ to $X_r$ is a path $p: [0,1] \to X$ with $p(0) \in X_s$ and $p(x > 0) \subset X_r$; the path is `exiting' from the stratum $X_s$ into the stratum $X_r$.  Whether to focus on entrance or exit paths is a matter of convention and convenience; in this book, we find that entrance paths have more natural functoriality dependencies and so we work entirely with them.
\end{rmk}

\begin{eg}[Entrance path poset] \label{eg:epp}
\autoref{fig:a-stratification-and-its-entrance-path-poset} shows a stratification of the open 2-\disk{} into five strata, along with its entrance path poset (shown on the right, together with an indication of which poset elements correspond to which strata).
\begin{figure}[ht]
    \centering
    \def\svgwidth{1\columnwidth}
    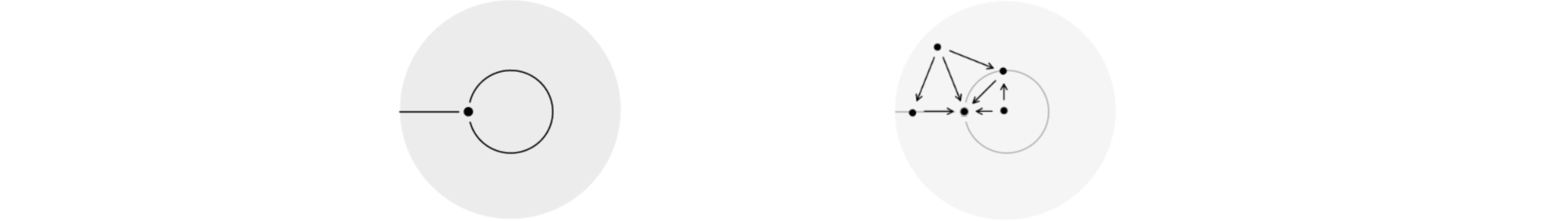

    \caption{A stratification and its entrance path poset.}
    \label{fig:a-stratification-and-its-entrance-path-poset}
\end{figure}
\end{eg}

\begin{eg}[Entrance path poset requiring transitive closure] \label{eg:eppt}
\autoref{fig:entrance-path-poset-as-the-transitive-closure} depicts a stratification of the open interval, into one open interval and two half-open interval strata, together with its entrance path poset.
\begin{figure}[ht]
    \centering
    \def\svgwidth{1\columnwidth}
    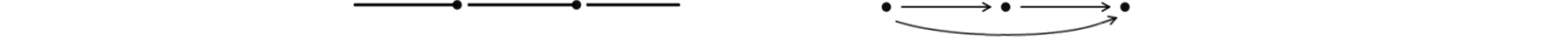

    \caption[Entrance path poset as the transitive closure.]{A stratification with entrance path poset as the transitive closure of the entrance path relation.}
    \label{fig:entrance-path-poset-as-the-transitive-closure}
\end{figure}
\end{eg}

\begin{eg}[Entrance path preorder]
In \autoref{fig:not-a-stratification-but-a-prestratification} we depict a decomposition of the circle that is not a stratification but a prestratification, because the formal entrance path relation has a cycle.
\begin{figure}[ht]
    \centering
    \def\svgwidth{1\columnwidth}
    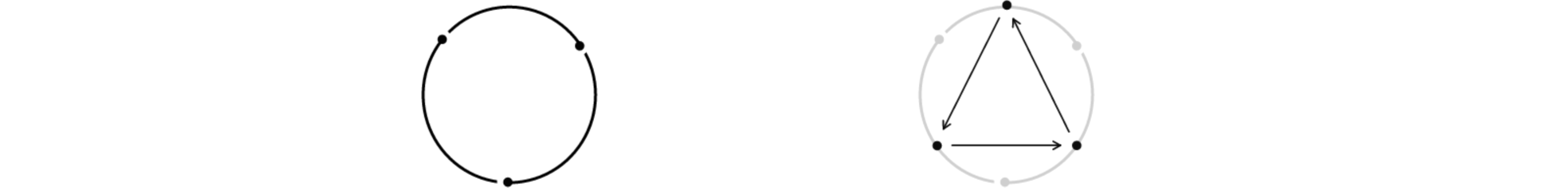

    \caption[Not a stratification but a prestratification.]{A decomposition that is not a stratification but a prestratification, and its formal entrance path relation.}
    \label{fig:not-a-stratification-but-a-prestratification}
\end{figure}
\end{eg}

\begin{term}[Discrete and indiscrete stratifications] \label{rmk:discrete-and-indiscrete-strat} Every space $X$ has an `indiscrete stratification' whose strata are the connected components of $X$. The entrance path preorder of the indiscrete stratification is the set of connected components of $X$. Conversely, every space also has a `discrete prestratification', such that each point becomes its own stratum. The entrance path poset of the discrete stratification of $X$ is the specialization order $\Spec X$ (in particular, the definition of specialization orders can be recovered from the definition of entrance path preorders).
\end{term}

\nid Unless indicated otherwise, a bare topological space is implicitly given the indiscrete stratification.

\begin{term}[Finite (pre)stratifications]
    We call a (pre)stratification $(X,f)$ `finite' if its entrance path preorder $\Entr(f)$ is finite, and call it `infinite' otherwise.
\end{term}

\begin{defn}[Characteristic function] \label{term:characteristic_function}
    Given a prestratification $(X,f)$, we refer to the function $X \to \Entr(f)$ sending each point $x \in X_r$ to its corresponding stratum $r \in \Entr(f)$, as the \textbf{characteristic function} of the prestratification; we denote the characteristic function of a prestratification $(X,f)$ by $f: X \to \Entr(f)$.
\end{defn}

\nid A fundamental property of characteristic functions is that they are `finitely continuous', as follows.

\begin{term}[Finitely continuous maps] A function of topological spaces $F : X \to Y$ is called `finitely continuous' for each finite subspace $Q \subset Y$ the function restricts to a continuous map $F : F\inv(Q) \to Q$.
\end{term}

\begin{lem}[Finite continuity in prestratifications]  \label{rmk:charact_fun_infinite} Characteristic functions are finitely continuous.
\end{lem}

\begin{proof} Consider a prestratification $(X,f)$ with characteristic function $f : X \to \Entr(f)$. Consider a finite subposet $Q \subset \Entr(f)$, and let $U \subset Q$ be a downward closed subposet. Arguing by contradiction, assume $f\inv(U) \subset f\inv(Q)$ is not open. Then there is a point $p \in f\inv(U)$ such that each neighborhood of $p$ intersects a preimage $f\inv(q)$ of some $q \in Q \setminus U$ not in $U$. Since $Q$ is finite, we can pick a $q \in Q \setminus U$ such that $f\inv(q)$ intersects \textit{all} neighborhoods of $p$. This means $p$ lies in the closure $f\inv(q)$ which entails there is an arrow from $q$ into some object of $U$, contradicting downward closure of the latter subposet. Thus, $f\inv(U) \subset f\inv(Q)$ is open, showing finite continuity of $f$.
\end{proof}

\nid In the case of \emph{finite} prestratifications, this of course implies that their characteristic functions are continuous in the usual sense.

\begin{term}[Locally finite prestratifications]
    A prestratification $(X,f)$ is `locally finite' if every stratum $s$ has an open neighborhood $s \subset N(s)$ which contains only finitely many strata.
\end{term}

\begin{term}[Covering relation] \label{rmk:covering-relation} Given a poset $(P,\leq)$ its covering relation is usually defined as follows: we say $x \in P$ `covers' $y \in P$, written $y <\cov x$, if $y < x$ is non-refinable (that is, for any $y < z < x$ we have either $y = z$ or $z = x$). Note that $(P,\leq)$ is the reflexive transitive closure of its covering relation.
\end{term}

\begin{lem}[Locally finite characteristic functions are continuous] \label{rmk:locally-finite}
     If $(X,f)$ is a locally finite prestratification, then its characteristic function $f: X \to \Entr(f)$ is continuous.
\end{lem}

\begin{proof} We need to show that for each $s \in \Entr(f)$ the downward closure $\Entr(f)^{\leq s}$ of $s$ has open preimage in $X$ under $f$. For $r \in \Entr(f)^{\leq s}$, let $Q_r \subset \Entr(f)^{\leq s}$ be the full subposet containing only $r$ and the elements $r' <\cov r$ that it covers. The assumption that $(X,f)$ is locally finite implies that $Q_r$ is finite, and that $f\inv(Q_r)$ contains an open neighborhood of $r$. It follows that $f\inv(\Entr(f)^{\leq s})$ is open as required.
\end{proof}

\begin{rmk}[Characteristic `maps']  Whenever a characteristic function is continuous we usually refer to it as a `characteristic map'. Note, characteristic functions of general infinite prestratifications need not be continuous, as the next example illustrates.
\end{rmk}

\begin{eg}[Infinite characteristic functions can be discontinuous] In \autoref{fig:a-stratification-with-non-continuous-characteristic-function} we depict a stratification of the closed interval with non-continuous characteristic function. In particular, the stratification is not locally finite.
\begin{figure}[ht]
    \centering
    \def\svgwidth{1\columnwidth}
    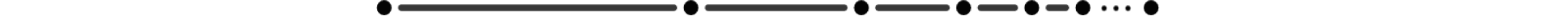

    \caption{A stratification with non-continuous characteristic function.}
    \label{fig:a-stratification-with-non-continuous-characteristic-function}
\end{figure}
\end{eg}

\nid As we will see in \autoref{lem:cquot_are_charact}, there is a precise characterization of those functions $f : X \to P$ from a space to a finite poset which are characteristic maps of stratifications.

From now on, we will focus most of our attention on stratifications instead of working in the more general context of prestratifications. In particular, most definitions will be given for stratifications only---however, the reader will notice that many immediately generalize if one replaces `stratifications' by `prestratifications', and `posets' by `preorders'.

\subsecunnum{Poset structures and quotient maps}

In this section we relate our definitions of stratifications with general `poset structures' on spaces.

\begin{defn}[Poset structures] Given a poset $P$, a \textbf{$P$-structured} space $(X,f)$ is a space $X$ together with a continuous map $f : X \to P$.
\end{defn}

\nid We will first show that characteristic maps of finite stratifications can be understood as a certain class of poset structures. Later we will show that, conversely, every poset structure can be universally `split' into a stratification.

Recall, a surjective continuous map $f : X \to Y$ of spaces is a quotient map if for each subset $U \subset Y$ we have that $U$ is open if and only if $f\inv(U)$ is open. If $Y$ is the specialization topology of a poset, we call $f$ a `poset quotient'. Poset quotients (to finite posets) admit the following useful characterization.

\begin{lem}[Quotient maps to finite posets] \label{lem:quotient_maps} For a space $X$, a finite poset $P$, and a surjective continuous map $f : X \to P$, the following are equivalent:
\begin{enumerate}
\item $f$ is a quotient map,
\item for all covers $p <\cov p'$ in $P$ there is a formal entrance path from $f\inv(p)$ to $f\inv(p')$.
\end{enumerate}
\end{lem}

\begin{rmk}[A quotient of posets is a map that is surjective on objects and on covers] \label{rmk:quotient_maps}
In the lemma, if $X$ is itself the specialization topology on a poset $Q$, then the lemma simplifies to saying that $f : Q \to P$ is a quotient map if and only if $f$ is surjective on objects and on covers.
\end{rmk}

\begin{proof}[Proof of \autoref{lem:quotient_maps}]  For $p \in P$, define $K^p_0$ to be the preimage $f\inv(p)$. Set $I^p_0 = \Set{p}$. Let $I^p_1$ be the set of $q \in P$ such that $f\inv(q)$ intersects the closure $\overline K^p_0$ of $K^p_0$, and define $K^p_1$ to be the union of preimages $f\inv(q)$ of $q \in I^p_1$. Set $I^p_2$ to be the set of $q \in P$ such that $f\inv(q)$ intersects the closure $\overline K^p_1$ of $K^p_1$, and define $K^p_2$ to be the union of preimages $f\inv(q)$ of $q \in I^p_2$. Repeating this process, since $P$ is finite, we find an index $j$ with $I^p_j = I^p_{j +1}$ and $K^p_j = K^p_{j+1} = \overline K^p_j$. Denote these sets by $I^p$ and $K^p$ respectively.

First, assume $f$ is a quotient map. Consider a cover $p < p'$. We claim it is impossible that $p' \notin I^p$: indeed, the complement $X \setminus K^p$ is the preimage of $P \setminus I^p$. Since $X \setminus K^p$ is open and since $f$ is a quotient map, it follows that $P \setminus I^p$ is open which contradicts the assumption that $p < p'$ and $p' \notin I^p$. Thus assume $p' \in I^p$. This implies $f\inv(p')$ intersects $K^p$ (and thus $f\inv(p') \subset K^p$). Then there is a sequence of arrows $p = p_0 < p_1 < ... < p_k = p'$ with $p_i \in K^p_i$. Since $p < p'$ is a cover we must have $k = 1$. Thus $f\inv(p')$ intersects the closure of $f\inv(p)$.

Next, assume $f$ satisfies that for any cover $p < p'$ in $P$, the preimage $f\inv(p')$ intersects the closure of the preimage $f\inv(p')$. Let $Q \subset P$ be a subposet. Let $I^{P \setminus Q}$ and $K^{P \setminus Q}$ be the respective unions of all $I^p$ and $K^p$ for each $p \in P \setminus Q$. If $Q$ is open then $f\inv(Q)$ is open by continuity of $f$. If $f\inv(Q)$ is open, then must be disjoint from $K^{P \setminus Q}$ (by construction of $K^{P \setminus Q}$). Thus $I^{P \setminus Q} = P\setminus Q$. Since $I^{P \setminus Q}$ is upward closed, it follows that $Q$ is downward closed, i.e.\ open as required.
\end{proof}

\nid A central role will be played by poset quotients whose `equivalence classes are connected' in the following sense.

\begin{defn}[Connected-quotient maps] \label{defn:cquot-map} For a space $X$ and a finite poset $P$, a continuous map $f : X \to P$ is called a \textbf{connected-quotient map} if it is a poset quotient map whose preimages of points $p \in P$ are connected. (Note, we take `connected' to also entail `non-empty'.)
\end{defn}

\begin{rmk}[Connected-quotient maps between posets] \label{defn:cquot-map-of-posets} A connected-quotient map $f : Q \to P$ where $Q$ is a poset (endowed with specialization topology) is a poset quotient whose preimages are connected subposets of $Q$.
\end{rmk}

\begin{figure}[ht]
    \centering
    \def\svgwidth{1\columnwidth}
    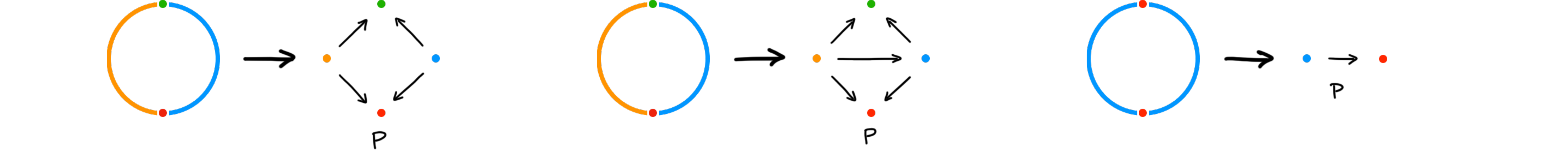

    \caption{A connected-quotient map and non-examples.}
    \label{fig:connected-quotient-maps}
\end{figure}

\begin{eg}[Connected-quotient map]
    In \autoref{fig:connected-quotient-maps} we depict three maps from the circle to three different posets (color-coding images and preimages in the same color). The first map is a connected quotient map; the second maps fails to be a quotient map despite have connected preimages, the third map is a quotient map but fails to have connected preimages.
\end{eg}

We now characterize stratifications among $P$-structures.

\begin{lem}[Characteristic maps are connected-quotient maps] \label{lem:cquot_are_charact} For a space $X$, a finite poset $P$, and a $P$-structure $f : X \to P$, the following are equivalent:
\begin{enumerate}
\item $f$ is the characteristic map of a stratification (that is, the decomposition of $X$ into preimages of $f$ is a stratification with characteristic map $f$ and entrance path poset $\Entr(f) = P$),
\item $f$ is a connected-quotient map.
\end{enumerate}
\end{lem}

\begin{proof} If $f$ is a characteristic map then, by definition, it has connected preimages and satisfies the second condition in \autoref{lem:quotient_maps}. Thus $f$ is a connected-quotient map.

Conversely, if $f$ is a connected-quotient map, then $f$ defines a stratification by decomposing $X$ into the preimages of $f$ (which are connected by \autoref{defn:cquot-map}). By \autoref{lem:quotient_maps} the map $f : X \to P$ is exactly the characteristic map of this stratification.
\end{proof}

The correspondence of characteristic maps and connected-quotient maps may further be generalized to the context of infinite stratifications, characterizing characteristic functions as `finitely connected-quotient' maps (analogous to the notion of `finite continuity' in \autoref{rmk:charact_fun_infinite}), but we forego a discussion of the infinite case here. We end with the following observation.

\begin{obs}[Connected-quotient maps compose] Using the definition of connected-quotient maps one verifies that, given a connected-quotient map $X \to P$ (of some stratification on $X$) and a connected-quotient map $P \to Q$ (of some stratification of $P$), their composite $X \to Q$ yields another connected-quotient map.
\end{obs}

\nid In particular, in the case of finite stratifications (and, under appropriate conditions, in the infinite case as well), we find that compositions of characteristic functions are again characteristic functions. As we will see, compositions of characteristic functions describe coarsenings (see \autoref{rmk:coarsenings_as_cquotients}).

\subsecunnum{Factoring poset structures into stratifications and labelings}

In this section we show that any poset structure factors into a stratification followed by a `labeling'. A labeling of a stratification in a category $\iC$ functorially associates data in $\iC$ to formal entrance paths in that stratification. We have seen conceptually similar notions of labelings in our discussions of labeled trusses (see \autoref{defn:n-trusses-labeled}), and in their relation to stratified meshes (see \autoref{defn:stratified-meshes}).

\begin{term}[Labelings] \label{defn:labelings} Let $\iC$ be a category, and $(X,f)$ a (pre)stratification. A `$\iC$-labeling' (or simply a `labeling') of $(X,f)$ in $\iC$ is a functor $L : \Entr(f) \to \iC$. If $\iC$ is a poset, we also call $L$ a `poset labeling'.
\end{term}

\nid There are other notions of entrance path categories (such as entrance path \textit{$\infty$-categories}, which we will meet later in \autoref{defn:tentr}) that can be considered in place of entrance path preorders. Working only with preorders provides the `$0$-categorical' base case of such fundamental category structures.

\begin{eg}[Specialization labelings] Let $(X,f)$ be a finite (pre)stratification. The `specialization labeling' of $X$ associated to $f$ is the labeling of the discrete prestratification of $X \to \Spec X$ given by the functor $\Spec f : \Spec X \to \Entr(f)$ (obtained by applying the specialization topology functor to the continuous map $f : X \to \Entr(f)$).
\end{eg}

We now show that any $P$-structure canonically factors as a stratification with a discrete labeling on that stratification. This factorization is referred to as the $P$-structure's `connected component splitting'. Discreteness of the labeling will mean the following.

\begin{term}[Discrete map] \label{defn:discretes} A map of posets $F : Q \to P$ is called a `discrete map' if its preimages are discrete, that is, for each $q \in Q$ the preimage $F\inv(q)$ contains no non-identity arrows.
\end{term}

\begin{constr}[Connected component splittings] \label{prop:conn_comp_split} For a  $P$-structure $f : X \to P$, the \textbf{connected component splitting} of $f$ is the factorization
\begin{equation}
f = (X \xto {\qqq f} \conn f \xto {\sss f} P)
\end{equation}
defined as follows. The map $\qqq f$ is the characteristic function of the stratification decomposing $X$ into the connected components of preimages of $f$; note that the formal entrance path graph cannot have cycles since $P$ is assumed to be a poset and $f$ to be continuous. The map $\sss f : \conn f \to P$ maps a given connected component $X$ of a preimage $f\inv(p)$ back to $p$.
\end{constr}

Note that even if $f : X \to P$ is continuous, the characteristic function $\qqq f$ need not be continuous (see \autoref{eg:p-stratification-translation}). We point out three universal properties of connected component splittings: universality among connected-quotient factorizations, universality among discrete map factorizations, and uniqueness among connected-quotient and discrete map factorizations.

\begin{lem}[Universality of connected component splitting] \label{prop:conn_comp_split_spaces_uniq} Let $f : X \to P$ be a $P$-structure. Assume $f$ factors into maps $g : X \to Q$ and $b : Q \to P$, where $g$ is continuous and $b$ a map of posets. Consider the following diagram:
\begin{equation}
\begin{tikzcd}[row sep=5pt, baseline=(W.base)]
 & Q \arrow[rd, "b"] \arrow[dd, dashed, no head] & \\
X \arrow[rr, "f", near start] \arrow[ru, "g"] \arrow[rd, "\qqq f"'] & & P \\
 & \conn f \arrow[ru, "\sss f"'] & |[alias=W]|
\end{tikzcd} .
\end{equation}
\begin{enumerate}
\item \emph{Characteristic map universality}: If $g$ is characteristic, then there is a unique poset map $Q \to \conn f$ making the above diagram commute,
\item \emph{Discrete map universality}: If $b$ is a discrete map, then there is a unique poset map $\conn f \to Q$ making the above diagram commute
\item \emph{Combined universality}: If $g$ is characteristic and $b$ is a discrete map, then there is a unique poset isomorphism $Q \iso \conn f$ making the above diagram commute.
\end{enumerate}
\end{lem}

\begin{proof} We first prove statement (1). Since $g$ is characteristic it has connected preimages. Thus its preimages must lie in the connected components of preimages of $f$. The map $\conn f \to Q$ is the inclusion of strata of $g$ into strata of $\qqq f$.

    We next prove statement (2). We first show that preimages of $g$ are unions of strata of $\qqq f$ (i.e.\ connected components of preimages of $f$). Let $Z$ be a connected component of a preimage of $f$. Let $\{q^Z_{i}\}_{i \in I}$ be the set of objects in $Q$ whose preimages $r^Z_{i} = g\inv(q^Z_{i})$ intersect $Z$. Note that, since $b$ is assumed to be a discrete map, there are no arrows between any $q^Z_{i}$ in $Q$. Let $Q^Z_{i}$ be the downward closure of $q^Z_{i}$ in $Q$. Since $g$ is assumed continuous, we have a disjoint open cover $\sqcup_i g\inv(Q^Z_{i}) \cap Z$ of $Z$. Since $Z$ is connected, the indexing set $I$ must be of cardinality $1$. This shows that preimages $g\inv(q)$ of $g$ are unions of connected components $Z$ of preimages of $f$. The map $\conn f \to Q$ can then be defined by mapping the strata $Z \subset g\inv(q)$ back to $q$.

The final statement (3) follows from combining statements (1) and (2).
\end{proof}

\begin{figure}[ht]
    \centering
    \def\svgwidth{1\columnwidth}
    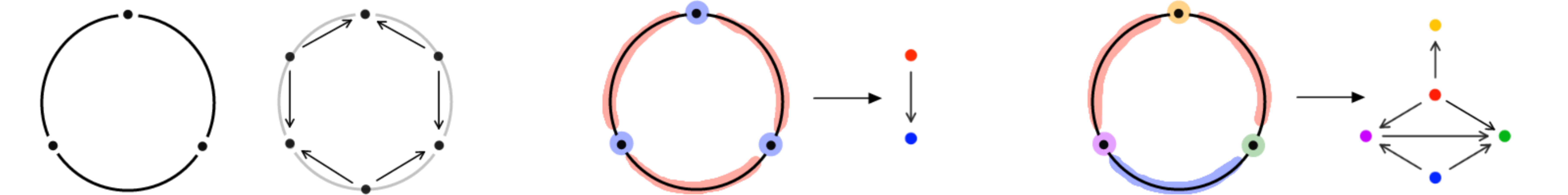

    \caption{Poset structures with the same underlying stratification.}
    \label{fig:translating-p-stratifications-into-stratifications}
\end{figure}

\begin{figure}[ht]
    \centering
    \def\svgwidth{1\columnwidth}
    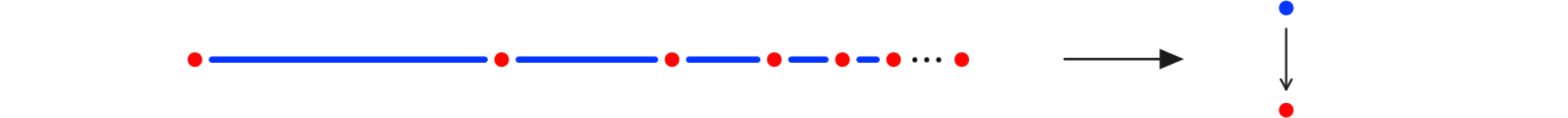

    \caption{A non-continuously splitting poset structure.}
    \label{fig:a-non-continuously-splitting-p-stratification}
\end{figure}

\begin{eg}[Translating $P$-structures into stratifications] \label{eg:p-stratification-translation} In \autoref{fig:translating-p-stratifications-into-stratifications} we depict a stratification of the circle on the left. To its right, we depict two $P$-structures (by color-coding images and preimages in the same color). Both $P$-structures recover the stratification on the left after connected component splitting. In particular, there are many $P$-structures with the same `underlying stratification'. In \autoref{fig:a-non-continuously-splitting-p-stratification} we depict another $P$-structure; it's connected component splitting recovers the stratification from \autoref{fig:a-stratification-with-non-continuous-characteristic-function} with non-continuous characteristic function.
\end{eg}

Using the above constructions, we describe relations of our notion of stratifications with two other frequently used definitions of stratifications, namely to `$P$-stratifications' and `$\sS$-filtered spaces'.

\begin{rmk}[Relation of stratifications and `$P$-stratifications'] \label{rmk:strat-vs-p-strat} Our notion of a `$P$-structure', given by a continuous maps from $X$ to $P$, is also known as a `$P$-stratification' (see \cite[Def. A.5.1]{lurie2012higher}). Indeed, by the above results, given a $P$-stratification $f : X \to P$ we can think of it as the stratification with characteristic map $\qqq f$ obtained by connected component splitting (note that there may be many different $P$-stratifications $f : X \to P$ that lead to the same stratification in this way). Conversely, every locally finite stratification $(X,f)$ arises as the connected component splitting of a $P$-stratification; indeed, by \autoref{rmk:locally-finite}, we can simply set $P = \Entr(f)$ and the characteristic map $f : X \to P$ will be continuous.
\end{rmk}

\begin{rmk}[Relation of stratifications and `$\sS$-filtered spaces'] \label{rmk:strat-vs-filtration} Given a poset $\sS$ with unique maximal element $\top$, a `$\sS$-filtration' of a space $X$ is a collection of closed subset $X_s$, $s \in \sS$, such that  $X_\top = X$ and $X_s \subset X_t$ whenever $t \leq s$ in $\sS$ (see \cite[\S III.2.2.1]{goresky1988stratified}). This defines a continuous map $f_\sS : X \to \sS$, mapping points in the subspace $X_t \setminus \bigcup_{s > t} X_s$ to $t \in \sS$. The characteristic function $\qqq {f_\sS}$ of the connected component splitting of $f_\sS$ yields a stratification in our sense. Conversely, every stratification $(X,f)$ with continuous characteristic map $f : X \to \Entr(f)$ yields an $\Entr(f)^\adjtop$-filtration of $X$ by setting $X_s = f\inv(\Entr(f)^{\geq s})$ (here, $\Entr(f)^\adjtop$ is the poset obtained by adjoining a new top element $\top$ to $\Entr(f)$, and $\Entr(f)^{\geq s}$ is the upper closure of an element $s$ in $\Entr(f)$).\footnote{In fact, in the case of finite stratifications, we can always recover $(X,f)$ from an $\lN$-filtration of $X$. Define a `depth map' $\mathrm{depth} : \Entr(f) \to \lN\op$, mapping $s \in \Entr(f)$ to $k$ if chains in $\Entr(f)$ starting at $s$ have maximal length $(k+1)$ (e.g. greatest elements have depth $0$). Define the filtration $X_0 \subset X_1 \subset ... \subset X_{k_{\mathrm{max}}} = X$ (where $k_{\mathrm{max}}$ is the maximal depth of elements in $\Entr(f)$) by setting $X_i$ to be the preimage of $[0,i]$ under the composite $\mathrm{depth}\circ f$.}
\end{rmk}

\subsecunnum{Classifying stratifications of posets}

In this section we discuss that the classifying space of any poset $P$ itself carries canonically the structures of a stratification, called the `classifying stratification' of $P$. While we introduced the concept already in \autoref{term:classifying-stratifications}, here we will revisit its construction in more concrete terms.

\begin{rmk}[Nerves of posets] \label{rmk:nerves} Recall the nerve $NP$ of a poset $(P,\leq)$ is the simplicial set whose $m$-simplices $S$ are the length-$m$ strings of composable arrows in $P$; in other words, an $m$-simplex is a map of posets $S: [m] \to P$.  The simplex $S: [m] \to P$ is called nondegenerate if it is injective.
\end{rmk}

\begin{rmk}[Classifying space of posets] \label{rmk:realizing-posets} Recall the `classifying space' $\abs{P}$ of a poset $P$ is the geometric realization of the nerve of $P$. (Abusing terminology, we sometimes refer to the classifying space of a poset itself as the `geometric realization' of the poset.) Explicitly, $\abs{P}$ is the space of functions $w: \obj(P) \to \lR_{\geq 0}$ whose support $\supp(w) \subset \obj(P)$ is the object set of a nondegenerate simplex in $P$, and whose total weight is fixed, i.e.\ $\sum_{p \in \obj(P)} w(p) = 1$.  We think of such a function $w$ as an `affine combination' of objects of the poset.
\end{rmk}

\begin{constr}[Classifying stratifications of posets] \label{constr:classstrat}
    The classifying space $\abs{P}$ of a poset $P$ has a stratification $\CStr {P}$, called the \textbf{classifying stratification}, with entrance path poset $P$ itself, constructed as follows.  The characteristic function of this stratification sends an affine combination $w$ of objects of the poset to the minimal object $\min(\supp(w))$ (in $P$) of the support of that affine combination:
\[
\CStr {P} : \abs{P} \to P \quad,\quad w \mapsto \min(\supp(w)) \in P.
\]
The stratum corresponding to the object $p \in P$ is denoted $\cstratum(p) \subset \CStr {P}$; it consists of all affine combinations $w$ of objects weakly greater than $p$, whose value at $p$ is nonzero.
\end{constr}

\nid In fact, the classifying stratification construction is functorial (as we will see in \autoref{constr:classstrat-functor}).


\section{Stratified maps}

\subsecunnum{Maps, coarsenings, and substratifications}

\begin{defn}[Map of stratifications] \label{defn:strat_maps}
A \textbf{map of stratifications} $F : (X, f: X \to \Entr(f)) \to (Y, g: Y \to \Entr(g))$, also called a `stratified map', is a continuous map $F : X \to Y$ for which there exists a (necessarily unique) map of posets $\Entr(F): \Entr(f) \to \Entr(g)$ such that $\Entr(F) \circ f = g \circ F$.
\end{defn}

\begin{notn}[Shorthand for stratified maps] \label{notn:abbreviating-strat-maps} In analogy to \autoref{notn:abbreviating-stratifications}, we often abbreviate stratified maps $F : (X,f) \to (Y,g)$ by $F : f \to g$.
\end{notn}

\begin{eg}[Map of stratifications]
    In \autoref{fig:a-stratified-map-and-a-non-stratified-map} we depict a stratified map on the left and a non-stratified map on the right. In both case, the underlying map of topological spaces is given by vertical projection.
\begin{figure}[ht]
    \centering
    \def\svgwidth{1\columnwidth}
    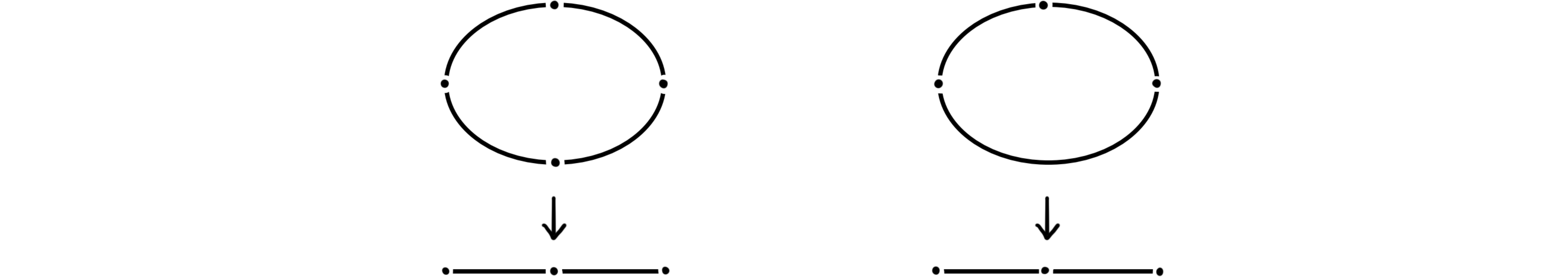

    \caption{A stratified map and a non-stratified map.}
    \label{fig:a-stratified-map-and-a-non-stratified-map}
\end{figure}
\end{eg}

\begin{defn}[Coarsenings and refinements of stratifications] \label{defn:coarsenings-and-refinements}
    A map of stratifications $F: (X,f) \to (Y,g)$ is a \textbf{coarsening} of $(X,f)$ to $(Y,g)$, or, synonymously, a \textbf{refinement} of $(Y,f)$ by $(X,f)$, if the underlying map of spaces $F: X \to Y$ is a homeomorphism.
\end{defn}

\nid Note that we use coarsening and refinement as synonyms of `dual flavor', i.e.\ describing dual processes: a coarsenings `coarsens' the domain, while a refinement, in opposite direction, `refines' the codomain.



\begin{eg}[Coarsening and refinement]
In \autoref{fig:a-coarsening-or-refinement-between-stratifications} we illustrate a coarsening of stratifications on the circle, along with the corresponding the refinement indicated by a dashed arrow in opposite direction.
\begin{figure}[ht]
    \centering
    \def\svgwidth{1\columnwidth}
    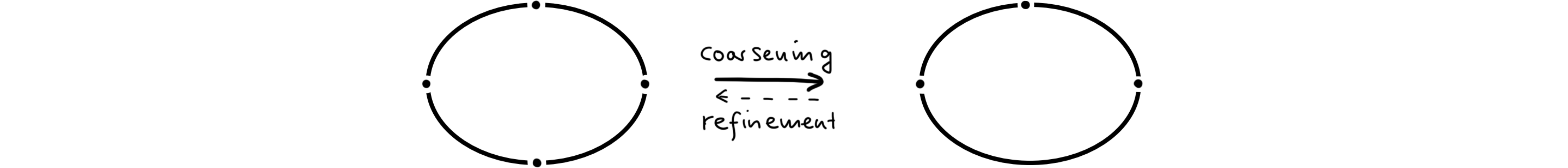

    \caption{A coarsening and its corresponding refinement of stratifications.}
    \label{fig:a-coarsening-or-refinement-between-stratifications}
\end{figure}
\end{eg}

\begin{defn}[Substratification]
A stratified map $(X,f) \to (Y,g)$ is a \textbf{substratification} if the underlying map $X \subset Y$ is an inclusion and if every stratum $s$ of $f$ is a connected component of $X \cap t$ for some stratum $t$ of $g$.
\end{defn}

\nid By extension we refer to stratified maps that are not literally inclusions, but whose underlying map is injective and a stratified homeomorphism onto a substratification, also as `substratifications'.

\begin{eg}[Substratification]
In \autoref{fig:two-stratified-maps-one-of-which-is-a-substratification} we depict two stratified maps: the first is a substratification, which though is not injective on entrance path posets; the second is a stratified map whose underlying map is injective, but which is not a substratification.
\begin{figure}[ht]
    \centering
    \def\svgwidth{1\columnwidth}
    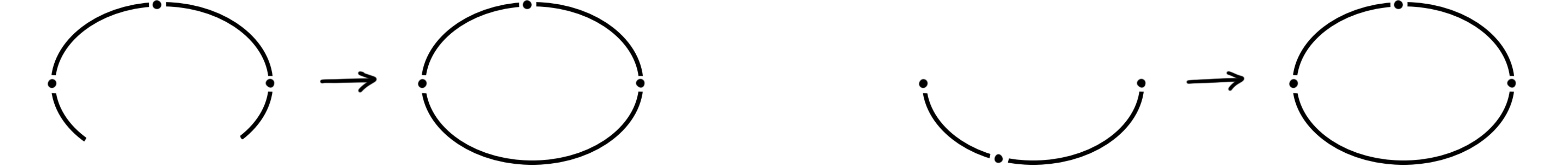

    \caption{Two stratified maps one of which is a substratification.}
    \label{fig:two-stratified-maps-one-of-which-is-a-substratification}
\end{figure}
\end{eg}

\begin{term}[Restricting stratifications] \label{term:strat-restr}
Given a stratification $(Y,g)$ and a subspace $X \subset Y$, the `restriction' $(X,\rest g X)$ is the substratification of $(Y,g)$ whose strata are the connected components of intersections $X \cap t$ for all strata $t$ of $g$.
\end{term}

\begin{defn}[Constructible substratifications] \label{defn:constr_substrat}
A substratification $(X,f) \to (Y,g)$ is \textbf{constructible} if every stratum of $(X,f)$ is exactly a stratum of $(Y,g)$.
\end{defn}


We now revisit the notion `boundary-constructible' stratifications (see \autoref{term:boundary-constructibility}): recall, a boundary-constructible stratification is a stratification $(Y,g)$ in which the topological closure $\overline s$ of each stratum $s$ yields a constructible substratification $(\overline s, \rest g {\overline s})$ of $g$ (by restricting $g$ to $\overline s$). Boundary-constructibility has an alternative description purely in terms of characteristic functions as follows.

\begin{lem}[Boundary-constructible stratifications are those with open characteristic function] \label{lem:boundary-constructibility-characteric-map}
A stratification $(X,f)$ is boundary-constructible if and only if the characteristic function $f : X \to \Entr(X)$ is an open map.
\end{lem}

\begin{proof}
Assume $f$ is boundary-constructible.  Let $U \subset X$ be an open subset.  We need to show that $f(U) \subset \Entr(X)$ is open, which in the specialization topology means that $f(U)$ is downward closed.  It suffices (because the entrance path poset is generated by formal entrance paths) to check that given an element $s \in f(U)$ and a formal entrance path $r \to s$, we have $r \in f(U)$.  The existence of the formal entrance path $r \to s$ means $s \cap \partial r \neq \emptyset$; boundary constructibility then implies that $s \subset \partial r$.  As $s \in f(U)$, there is some point of the stratum $s$ that is in $U$, and because $U$ is open, there must be a point of the stratum $r$ that is in $U$. Thus $r \in f(U)$ as required.

Conversely, assume $f: X \to \Entr(f)$ is open. It suffices to show that if there is a formal entrance path $r \to s$, i.e.\ $s \cap \partial r \neq \emptyset$, then $s \subset \partial r$.  Suppose there is such an entrance path but by contrast there is a point $p \in s \setminus \partial r = s \setminus \overline r$.  Then we can chose an open neighborhood $p \in U \subset X$ disjoint from the closure $\overline r$.  By assumption it follows that the image $f(U)$ is open, which is to say downward closed; thus $s \in f(U)$ and $r \to s$ implies that $r \in f(U)$, contradicting the fact that the neighborhood $U$ does not intersect even the closure of $r$.
\end{proof}

We can also characterize substratifications and coarsenings in terms of entrance path poset maps, as follows.

\begin{lem}[Substratification from discrete maps] \label{rmk:substratification_discrete} A map of finite stratified spaces $F : (X,g) \to (Y,f)$ is a substratification if and only if $F : X \to Y$ is a subspace inclusion and $\Entr(F) : \Entr(g) \to \Entr(f)$ is a discrete map.
\end{lem}

\begin{proof} By definition every substratification is a subspace inclusion of underlying spaces. The fact that $\Entr(F)$ is a discrete map follows since strata of substratifications are defined as connected components of the intersection of the subspace $X$ with strata of $f$, and since $\Entr(f)$ is a poset.

Conversely, assume the stratified map $F$ is a subspace inclusion and that $\Entr(F)$ is a discrete map. Note that the substratification $(X,\rest f X)$ of $f$ can be obtained by connected component splitting of the restriction of $f : Y \to \Entr(f)$ to $X \subset Y$. Since $g$ is a continuous characteristic map, and since $\Entr(F)$ is a discrete map, statement (3) of \autoref{prop:conn_comp_split_spaces_uniq} (applied to $f : Y \to \Entr(f)$ restricted to $X \subset Y$) shows that $g$ is a substratification of $f$ as claimed.
\end{proof}

\begin{lem}[Coarsenings from connected-quotient maps] \label{rmk:coarsenings_as_cquotients} Let $(X,f)$ be a finite stratification. Coarsenings of $f$ (up to isomorphism) are canonically in bijection with connected-quotients of $\Entr(f)$: namely, the bijection takes coarsenings $F$ to their entrance path poset maps $\Entr(F)$.
\end{lem}

\begin{proof} Let $F : (X,f) \to (X,g)$ be a coarsening. Then $\Entr(F)$ is a connected-quotient map since its preimages are connected and it satisfies condition (2) in \autoref{lem:quotient_maps}.

Conversely, let $H : \Entr(f) \to P$ be a connected-quotient map. Define a stratification $(X,g)$ whose strata are unions of those strata in $f$ that are mapped to the same object in $P$ under $H$. Since preimages of $H$ are connected, these unions are connected subspaces of $X$ and thus define a prestratification. Since $H$ is a connected-quotient map to a poset $P$, this prestratification is in fact a stratification with entrance path poset $\Entr(g) = P$. The resulting coarsening $(X,f) \to (X,g)$ is the identity on the underlying space $X$, and maps entrance path posets by $H$.
\end{proof}

\subsecunnum{The category of stratifications}

Having defined stratified spaces and maps, we now obtain the category of stratifications.

\begin{term}[The ordinary category of stratifications] Denote by $\Strat$ the category of stratifications and their stratified maps.
\end{term}

\nid Posets faithfully embed into stratifications by the classifying stratification functor as follows (recall the construction of classifying stratification from \autoref{constr:classstrat}).

\begin{constr}[Classifying stratification functor] \label{constr:classstrat-functor} Given a map of posets $F: P \to Q$, the induced map on classifying spaces, namely the map $\CStr F : \CStr P \to \CStr Q$ mapping $\CStr F(w)(q) = \sum_{p \in F\inv(q)} w(p)$, is a stratified map. This yields the `classifying stratification' functor
\begin{equation}
    \CStr - : \Pos \to \Strat
\end{equation}
from the category of posets to the category of stratifications.
\end{constr}

\nid Conversely, the entrance path poset construction previously described yields a functor from the category of stratifications to the category of posets.

\begin{constr}[Entrance path poset functor] \label{defn:entr-functor}
The association of the entrance path poset $\Entr(f)$ to the stratification $(X,f)$, and of the map of posets $\Entr(F)$ to the map of stratifications $F: (X,f) \to (Y,g)$ provides the `entrance path poset' functor
\begin{equation}
\Entr : \Strat \to \Pos
\end{equation}
from the category of stratifications to the category of posets.
\end{constr}

\begin{obs}[Entrance paths invert classifying stratifications] The preceding functors form a section-retraction pair: namely, $\Entr \circ \CStr{} = \id$.
\end{obs}

\nid In the case of sufficiently nice finite stratifications, we can further promote the  entrance path poset functor to a topological functor, for instance as follows. Here, we use the term `topological category' to refer to a category enriched in the category $\Top$ of topological spaces.

\begin{term}[The topological category of stratifications] Denote by $\TStrat$ the topological category of stratified spaces and their stratified maps (whose hom spaces are topologized as subspaces of the hom spaces in the topological category $\Top$ of topological spaces).
\end{term}

\begin{defn}[Topological entrance path poset functor]\label{rmk:entrz_cont}
Let $\TPos$ denote the $\Pos$-enriched category of posets; that is, arrows in the poset $\TPos(P,Q)$ are natural transformations.  Topologizing the hom posets in $\TPos$ using the specialization topology, we may consider $\TPos$ as a topological category. Let $\TStratfinlch$ denote the full subcategory of $\TStrat$ consisting of finite stratifications whose underlying spaces are locally compact Hausdorff spaces. The entrance path poset provides a functor of topological categories
\begin{equation}
\Entr : \TStratfinlch \to \TPos \quad . \qedhere
\end{equation}
\end{defn}

\nid The continuity of the functor $\Entr$ on hom spaces follows from standard arguments.\footnote{Recall, for compactly generated spaces $X,Y$ the internal hom $\Map(X,Y)$ is the $k$-ification of the natural topology on the set of continuous functions $C(X,Y)$ (see \cite[Prop. 5.11 and Thm. 5.15]{escardo2004comparing}). Denote the function set with natural topology by $\Top(X,Y)$. If $X$ is locally compact Hausdorff, then $\Top(X,Y)$ can be defined using the usual compact-open topology with subbase elements of the form $M(K,U) = \Set{F : X \to Y ~|~ F(K) \subset U}$ for compact $K \subset X$ and open $U \subset Y$ (see \cite[Theorem 5.6]{bradley2020topology}).
If $A,B$ are specialization topologies of finite posets then, firstly, $\Top(A,B)$ is finite too, and has subbase elements of the form $M(U,V) = \Set{F : A \to B ~|~ F(U) \subset V}$ for open sets $U,V$ (see \cite[Prop. 5.11]{escardo2004comparing}). Note also that $\Top(A,B) = \TPos(A,B)$ in this case (and further, being finite implies that $\Top(A,B)$ is compact and thus $\Top(A,B) = \Map(A,B)$). Now, let $X,Y$ be locally compact Hausdorff spaces, with finite stratifications given by characteristic maps $f : X \to \Entr(f)$ and $g : Y \to \Entr(g)$ respectively. Let $M(U,V) \subset \Map(\Entr(f),\Entr(g))$ be a subbase element.
We show $\Entr\inv M(U,V)$ is open in the subspace $\TStratfinlch(f,g)$ of $\Map(X,Y)$ (consisting of stratification maps $F : f \to g$) by showing it is open in the corresponding subspace of $\Top(X,Y)$. Let $F \in \Entr\inv M(U,V)$. Then $F$ maps strata in $U$ to strata in $V$. Pick a point $p_s$ of each stratum $s$ in $U$, and define $K = \Set{ p_s \in f\inv(s) ~|~ s \in U }$. Since $\Entr(f)$ is finite, $K$ is compact. Then $M(K,g\inv(V))$ is a subbase element of $\Top(X,Y)$ containing $F$ and contained in $\Entr\inv M(U,V)$ which shows the latter subspace is open.}

\begin{term}[Products of stratifications] \label{term:strat-products} Given two stratifications $(X,f)$ and $(Y,g)$, the product stratification is simply $(X \times Y, f \times g)$ where $f \times g$ is the characteristic function $X \times Y \to \Entr(f) \to \Entr(g)$ obtained by taking the product of characteristic functions $f : X \to \Entr(f)$ and $g : Y \to \Entr(g)$. One further defines products of stratified maps by taking products of their underlying continuous maps. This yields a topological `product' functor
\begin{equation}
    (- \times -) : \TStrat \times \TStrat \to \TStrat . \qedhere
\end{equation}
\end{term}

\begin{obs}[Fiberwise $\Top$-tensoredness of $\TStrat$] \label{rmk:strat_cc}
    Taking products with topological spaces provides a `fiberwise $\Top$-tensor' on the category of stratified spaces as follows. Let $(X,f)$ and $(Y,g)$ be finite locally compact Hausdorff stratifications and $F : \Entr(f) \to \Entr(g)$ a map of their entrance path posets. Denote by $\TStrat(f,g)_F$ the preimage of $F$ of the map $\Entr : \TStrat(f,g) \to \TPos(\Entr(f),\Entr(g))$. Using cartesian closedness of $\Top$, identify $\Map(Z,\Map(X,Y)) \iso \Map(Z \times X, Y)$ (for $Z \in \Top$); in particular, we obtain a homeomorphism
    \begin{equation}
        \Map(Z,\TStrat(f,g)_F) \iso \TStrat(Z \times f,g)_F
    \end{equation}
where the right hand side denotes the space of stratified maps $Z \times (X, f) \to (Y,g)$ whose underlying map of entrance path posets is $F$ (noting $\Entr(Z \times (X,f)) \iso \Entr(f))$).
\end{obs}

\subsecunnum{Stratified bundles}

A `stratified bundle' is a stratified map that is locally trivial along each stratum of the base. The notion generalizes the ordinary notion of `fiber bundles' of topological spaces.

\begin{defn}[Stratified bundles] \label{defn:stratified-bundle}
A stratified map $p : (Y,g) \to (X,f)$ is a \textbf{stratified bundle} if for each stratum $s$ of $f$ and each point $x \in s$, there is a neighborhood $U_x \subset s$ inside the stratum $s$, such that there is a stratification $(Z,h)$ together with an isomorphism of stratifications $T: U_x \times h \iso (p\inv(U_x),f)$ for which $p \circ T : U_x \times h \to U_x$ is the projection.  The stratification $(Z,h)$ is called the fiber of $F$ over the stratum $s$.
\end{defn}

\nid Note that every fiber bundle is naturally a stratified bundle with indiscrete stratifications on both base and total space. We will usually further assume that all the fibers of a stratified bundle are non-empty, in other words that the underlying map of spaces is surjective. In this case the stratification of the total space determines the stratification of the base space.

\begin{obs}[The base stratification is determined by the total stratification] \label{prop:base_unique}
Suppose $(Y,f) \to (X,g)$ and $(Y,f) \to (X,g')$ are stratified bundles with the same underlying surjective map $F: Y \to X$.  Then the stratifications $g$ and $g'$ are equal.
\end{obs}

Just as fiber bundles can be pulled back along continuous map, stratified bundles can be pulled back along stratified maps. Here, a `pullback' of stratified maps means the following.

\begin{term}[Pullbacks of stratifications] \label{term:strat-pullback} Given stratifications $(X,f)$, $(Y,g)$, and $(Z,h)$ and maps $F : f \to h$ and $G : g \to h$, the `pullback stratification' $(X,f) \times_{(Z,h)} (Y,g)$ is the stratification $(X \times_Z Y, f \times_h g)$, where $X \times_Z Y$ is the pullback of spaces and $f \times_h g$ is the restriction $\rest {f\times g} {X \times_Z Y}$ of the product stratification $f \times g$ to the pullback space $X \times_Z Y \subset X \times Y$.
\end{term}

\begin{eg}[Pullback stratification need not be finite or have continuous characteristic function]
In \autoref{fig:pullback-of-stratifications} we depict a pullback of finite stratifications that is not finite and does not have continuous characteristic function.
\begin{figure}[ht]
    \centering
    \def\svgwidth{1\columnwidth}
    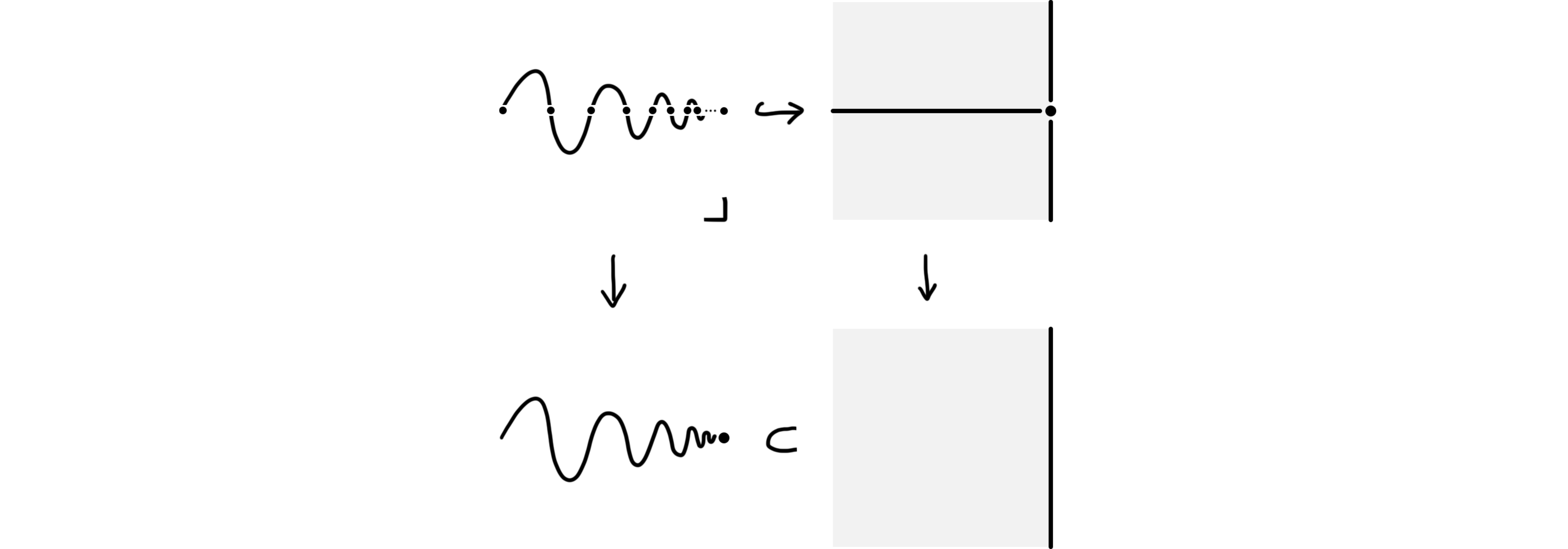

    \caption{Pullbacks of stratifications need not preserve finiteness.}
    \label{fig:pullback-of-stratifications}
\end{figure}
\end{eg}

\begin{obs}[Pullbacks of stratified bundles] Given a stratified bundle $p : (X,f) \to (Y,g)$ and a stratified map $F : (Y',g') \to (Y,g)$, then the pullback map $(X \times_Y Y', f \times_g g') \to (Y',g')$ is a stratified bundle itself, usually denoted by $F^*p : F^*f \to g'$. This follows since stratified maps map strata into strata, which allows us to `pull back' trivializations of $p$ over neighborhoods in strata of $g$ to trivialization of $F^*p$ over neighborhoods in strata of $g'$.
\end{obs}

\begin{rmk}[Constructible bundles] There is a natural strengthening of the notion of stratified bundles, namely to so-called `constructible bundles'.  As discussed in \autoref{ch:meshes}, `constructibility' requires that bundles can be reconstructed up to isomorphism from functorial data associated to fundamental categories (such as entrance path posets, or entrance path $\infty$-categories) of their base stratifications. There are several approaches to making this precise (see \cite[\S 6]{ayala2015stratified} \cite{curry2020classification}).
\end{rmk}

\section{Conical and cellular stratifications}

\subsecunnum{Conical stratifications} \label{sssec:conicality}

Many of the stratifications in this book satisfy an additional regularity condition called `conicality'. This condition requires neighborhoods of strata to locally look like a cone `normal' to the stratum and an open set `tangential' to the stratum. Let us first formalize the operation of taking cones on stratifications.

\begin{term}[Cones of stratifications] \label{term:strat-cone} Given a stratification $(X,f)$, we can define its `open cone' $f^\adjtop$ to be the stratification $(X^\adjtop, f^\adjtop : X^\adjtop \to \Entr(f)^\adjtop)$ as follows. The space $X^\adjtop$ is the usual open cone $X \times (0,1] \cup_{X \times \Set{1}} \top$, where $\top$ is the space with the single point $\top$. The poset $\Entr(f)^\adjtop$ is obtained from $\Entr(f)$ by adding a new top element $\top$.  The map $f^\adjtop$ sends the cone point $\top$ to $\top \in \Entr(f)^\adjtop$ and is otherwise given by $X \times (0,1) \xto{\pi_X} X \xto{f} \Entr(f)$.
\end{term}

\begin{defn}[Conical stratification] \label{defn:conical}
    A \textbf{tubular neighborhood} of a point $x$ of a stratification $(X,f)$ is a neighborhood $U_x$ of $x$, together with a stratified space $(Y_x,l_x)$, called a \textbf{link} at $x$, a connected topological space $Z_x$, called the \textbf{tangential neighborhood}, and a stratified homeomorphism
\[
Z_x \times l_x^\adjtop \iso (U_x,f |_{U_x})
\]
sending $z \times \top$ to $x$, for some $z \in Z_x$.  (Here $\top$ is the cone point in the cone $l_x^\adjtop$.)  A stratification is \textbf{conical} if it has a tubular neighborhood at every point.
\end{defn}

\begin{figure}[ht]
    \centering
    \def\svgwidth{1\columnwidth}
    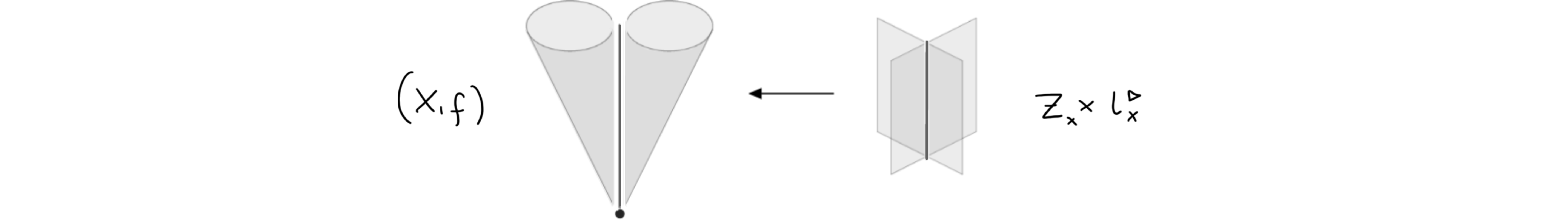

    \caption{A conical stratification with a tubular neighborhood.}
    \label{fig:a-conical-stratification-with-a-tubular-neighborhood}
\end{figure}

\begin{figure}[ht]
    \centering
    \def\svgwidth{1\columnwidth}
    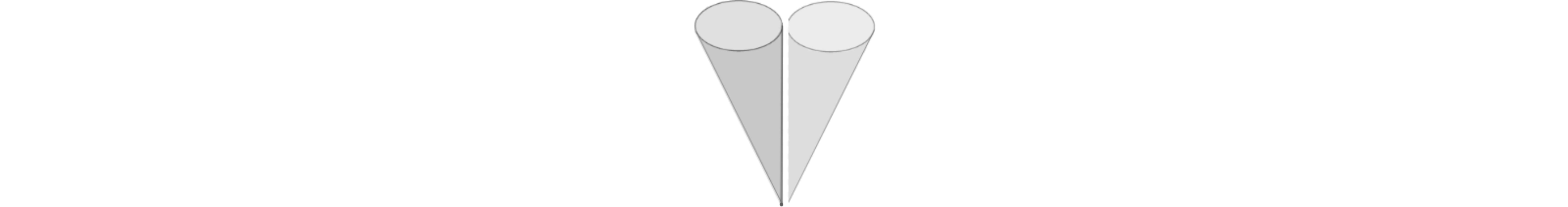

    \caption{A non-conical stratification.}
    \label{fig:a-non-conical-stratification}
\end{figure}

\begin{eg}[Conical and non-conical stratifications]
In \autoref{fig:a-conical-stratification-with-a-tubular-neighborhood} we illustrate an example of a conical stratification, together with an illustration of a tubular neighborhood as a product of a cone and a space. By contrast, in \autoref{fig:a-non-conical-stratification} we depict a stratification (of the same space, but now decomposed into only two strata) which is not a conical stratification.
\end{eg}

\begin{figure}[ht]
    \centering
    \def\svgwidth{1\columnwidth}
    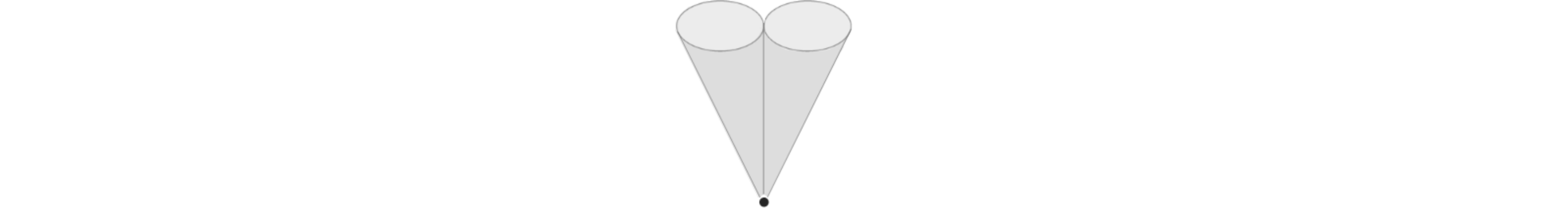

    \caption{A conical stratification that is not topological.}
    \label{fig:a-conical-stratification-that-is-not-topological}
\end{figure}

\begin{rmk}[Topological stratification]
    A conical stratification in which every stratum is a topological manifold is usually called a `topological stratification'.  Note that in that situation, the tangential spaces $Z_x$ can always be chosen to be euclidean spaces. The conical stratification shown in \autoref{fig:a-conical-stratification-with-a-tubular-neighborhood} is a topological stratification. An instance of a conical stratification (with two strata) that is not a topological stratification is depicted in \autoref{fig:a-conical-stratification-that-is-not-topological}.
\end{rmk}

\begin{rmk}[Constructible substratifications inherit conicality] \label{rmk:constr_substrat}
If the stratification $(X,f)$ is conical and $(Y,g) \into (X,f)$ is a constructible substratification (\autoref{defn:constr_substrat}), then the stratification $(Y,g)$ is also conical.
\end{rmk}

\begin{prop}[Classifying stratifications are conical] \label{rmk:nervestrat_conical}
The classifying stratification $\CStr P$ (see \autoref{constr:classstrat}) of any finite poset $P$ is conical.
\end{prop}

\begin{proof}
    For any poset element $p \in P$, and any point $x \in \cstratum(p)$ of the corresponding stratum, we construct a tubular neighborhood, which is in fact independent of the point $x$.  The points of the classifying stratification $\CStr P$, the stratum $\cstratum(p)$ and the link $\conelink(p)$ are given by certain affine combinations $w$, $w_{\cstratum(p)}$ and $w_{\conelink(p)}$, and $w$ of objects in $P$, $P^{\geq p}$ and $P^{<p}$ respectively. Since $P^{\geq p} \into P$ (respectively $P^{<p} \into P$) we may trivially extend $w_{\cstratum(p)}$ (respectively $w_{\conelink(p)}$) to affine combination of objects in $P$. The inclusion $\cstratum(p) \subset \CStr P$ now extends to a tubular neighborhood $\cstratum(p) \times \conestr(\conelink(p)) \to \CStr P$ by setting (cf. \cite[Prop. A.6.8]{lurie2012higher})
\begin{align*}
    \cstratum(p) \times \conelink(p) \times (0,1) & \to \CStr P\\
    (w_{\cstratum(p)},w_{\conelink(p)},t) & \mapsto  \left(w(q) := t \cdot w_{\cstratum(p)}(q) + (1-t)\cdot w_{\conelink^P_p}(q) \right)
\end{align*}
away from the conepoint $\conelink(p) \times (0,1) \into \conestr(\conelink(p))$. This verifies conicality of $\CStr P$ as claimed.
\end{proof}

\subsecunnum{Cellular stratifications}

In this final section we revisit the notion of regular cell complexes, and their properties as stratifications. Recall, a cell complex is called regular when the attaching maps are homeomorphisms onto their images, or equivalently, if the closure of each cell is a ball. We will show that regular cell complexes are both conical and `entrance path homotopically trivial' (more precisely, we will show that their `entrance path $\infty$-categories' are $0$-truncated).

We begin with the question of conicality. Note that the stratification of general (non-regular) cell complexes need not be conical in general, even when the complex is finite (i.e.\ has only finitely many cells). In contrast, in the regular case, we have the following.

\begin{prop}[Regularity implies conicality] \label{prop:regular-cell-cplx-conical} Finite regular cell complexes are conically stratified.
\end{prop}

\begin{proof} In \autoref{prop:cellular-posets-vs-reg-cell-cplx} we have seen that there is a stratified homeomorphism $X \iso \CStr \Entr X$ of $X$ with the classifying stratification of its entrance path poset. From \autoref{rmk:nervestrat_conical} it then follows that finite regular cell complexes are conically stratified.
\end{proof}

\nid Note that the result can be improved to hold for `locally finite' regular cell complexes without difficulties. Recall, a regular cell stratification is a constructible substratification of a regular cell complex (see \autoref{defn:regular-cell-strat}).


\begin{obs}[Regular cell stratifications are conical] \label{rmk:regular-cell-implies-conicality} Combining \autoref{prop:regular-cell-cplx-conical} with \autoref{rmk:constr_substrat} it follows that regular cell stratifications are conical.
\end{obs}

    Next let us explain how regular cell complexes are `entrance path homotopically trivial'. This requires us to first ask whether there is a good notion of `entrance path homotopies' between entrance paths. This leads to a definition of entrance path $(\infty,\infty)$-category which, in the context of conical stratifications, reduces to a notion of entrance path $(\infty,1)$-category (usually abbreviated simply as `entrance path $\infty$-category').  Modeled in quasicategories, this $\infty$-category can be expressed as follows.

\begin{defn}[{Entrance path $\infty$-category, \cite[Rmk. A.6.5]{lurie2012higher}}] \label{defn:tentr}
The \textbf{entrance path $\infty$-category} $\TEntr(f)$ of a conical stratification $(X,f)$ is the $\infty$-category with underlying simplicial set having $m$-simplices the maps from the stratified $m$-simplex to $f$; that is $\TEntr(f)_m := \Strat(\CStr{[m]},f)$.
\end{defn}

\begin{rmk}[The class of entrance path $\infty$-categories] The class of $\infty$-categories that are (up to equivalence) the entrance path $\infty$-categories of a conical stratifications have been described as $\infty$-categories with a conservative functor to a poset (see \cite[\S 2.1]{barwick2018exodromy}).
\end{rmk}

\nid Note that the conicality condition is essential to ensure that the given simplicial set is indeed an $\infty$-category.

\begin{lem}[Regular cell complexes are entrance path $0$-truncated] \label{lem:reg_use_Entrz}
For the stratification $(X,f)$ associated to a finite regular cell complex, the entrance path $\infty$-category is equivalent to the nerve of the entrance path poset:
\[
\TEntr(f) \eqv N\Entr(f)
\]
\end{lem}

\begin{proof}
Recall that a $\infty$-category is called $0$-truncated if its hom spaces are $(-1)$-types, meaning they are either empty or contractible.  We first check that $\TEntr(f)$ is $0$-truncated.  It suffices \cite[Prop. 3.12]{campbell2018truncated} to show that any sphere $\partial \spx {m} \to \TEntr(f)$, for $m > 1$, has a spacer.  Given such a map $s: \partial \spx {m} \to \TEntr(f)$, by the definition of $\TEntr(f)$, there is a lift $\tilde{s} : \abs{\partial \spx {m}} \to X$.  By construction, the map $\tilde{s}$ is built out of a collection of entrance paths in $X$; in particular, there is a sequence of entrance paths from the image $\tilde{s}(0)$ of the initial vertex $0 \in \partial \spx m$ to any other point $\tilde{s}(p)$, and so the image $\mathrm{im}(\tilde{s})$ is contained in the closure of the cell containing $\tilde{s}(0)$ in its interior.  By the regularity assumption, that closure is a closed $n$-\disk{} $D$.  Identify $\spx m \iso \big(\partial \spx m \times [0,1]\big)\big\slash_{\partial \spx m \times \Set{0}}$, and identify the \disk{} $D$ with a cone on its boundary.  Define the spacer $r : \spx m \to d \subset X$ by $(x,t) \mapsto t \cdot s(x)$.  By construction, the map $r$ sends the interior of $\spx m$ to the stratum containing $\tilde{s}(0)$, and so the map $r$ is stratified, as needed.

Because the entrance path $\infty$-category $\TEntr(f)$ is 0-truncated, it is equivalent to $N(\ho(\TEntr(X,f))$ \cite[Prop. 3.8]{campbell2018truncated}.
Furthermore, the homotopy category of any 0-truncated $\infty$-category has a skeleton that is a poset \cite[Prop. 3.10]{campbell2018truncated}.  Let $\overline{\ho(\TEntr(X,f))}$ denote the posetal skeleton of $\ho(\TEntr(X,f))$.  Observe that poset must be equivalent to $\Entr(f)$ as there are fully surjective maps of posets $\overline{\ho(\TEntr(X,f))} \to \Entr(f)$ and $\Entr(f) \to \overline{\ho(\TEntr(X,f))}$ sending a point to the stratum containing it, respectively sending a stratum to the equivalence class of any point in the stratum.  Altogether we have an equivalence as desired:
\[
\TEntr(f) \eqv N(\ho(\TEntr(f))) \eqv N(\overline{\ho(\TEntr(f)})  = N(\Entr(f)) \qedhere
\]
\end{proof}

\begin{cor}[Regular cell stratifications are entrance path $0$-truncated] Given a regular cell stratification $(X,f)$, we have a canonical equivalence $\TEntr(f) \eqv N\Entr(f)$, that is, the entrance path $\infty$-category of $f$ is $0$-truncated. \qed
\end{cor}

\backmatter

\listoffigures

\bibliographystyle{alpha}

\bibliography{cubical}

\end{document}

%% file: cubical-pre.tex
\usepackage{graphicx}
\usepackage{amssymb, amsmath, amsthm, stmaryrd, mathrsfs, mathpartir, color, cmll, tabulary, enumitem, mathtools}
\SetSymbolFont{stmry}{bold}{U}{stmry}{m}{n}

\usepackage{import}
\usepackage{pdfpages}
\usepackage{transparent}
\usepackage{xcolor}
\usepackage{epigraph}

\numberwithin{figure}{chapter}

\usepackage{setspace,wrapfig,verbatim,changepage,adjustbox, eso-pic}
\usepackage{braket,xspace,url,scalerel,microtype}
\usepackage{tikz, tikz-cd}
\usetikzlibrary{arrows,positioning,calc,decorations.markings,decorations.pathmorphing,shapes}
\usepackage{enumitem}
\usepackage[pdfencoding=auto,psdextra]{hyperref}
\hypersetup{bookmarksdepth=4,
            colorlinks,breaklinks,
            urlcolor=[RGB]{7,7,96},
            linkcolor=[RGB]{7,7,96},
            citecolor=[RGB]{7,7,96}}
\usepackage{float,bookmark,caption}
\usepackage[normalem]{ulem}
\usepackage[bbgreekl]{mathbbol}
\DeclareSymbolFontAlphabet{\mathbbm}{bbold}
\DeclareSymbolFontAlphabet{\mathbb}{AMSb}%

\newcommand\dout{\bgroup \markoverwith{\rule[0.2ex]{0.1pt}{0.4pt}\rule[0.8ex]{0.1pt}{0.4pt}}\ULon}

\usepackage{bclogo}
\usepackage{tabto}
\newcommand\marginsymbol[1][0pt]{%
  \tabto*{0cm}\makebox[-2cm][c]{\bcdz}\tabto*{\TabPrevPos}}

\graphicspath{{../geometry/}}

\usepackage{etoolbox} 

\makeatletter
\def\subsection{\@startsection{subsection}{2}%
    \z@{.5\linespacing\@plus.7\linespacing}{-.5em}%
    {\normalfont\bfseries}}
\makeatother

\DeclareRobustCommand{\SkipTocEntry}[5]{}

\makeatletter
\renewcommand\subsubsection{\@startsection{subsubsection}{3}%
  \z@{.5\baselineskip\@plus.7\baselineskip}{-.5em}%
  {\normalfont\bfseries}}
\renewcommand\paragraph{\@startsection{paragraph}{4}%
  \z@{.5\baselineskip\@plus.35\baselineskip}{-.5em}
  {\normalfont\itshape}}
\makeatother

\makeatletter
\newcommand{\firstToC}{{%
  \@fileswfalse
  \renewcommand{\contentsname}{Contents}%
  \@starttoc{toc}{\contentsname}%
}}
\newcommand{\anotherToC}{{%
  \@fileswfalse
  \renewcommand{\contentsname}{Contents (TEMPORARY)}%
  \@starttoc{toc}{\contentsname}%
}}
\newcommand{\lastToC}{{%
  \renewcommand{\@tocwrite}[2]{}
  \renewcommand{\contentsname}{Contents}%
  \@starttoc{toc}{\contentsname}%
}}
\makeatother

\makeatletter
\@ifclassloaded{amsbook}
  {
		\def\l@paragraph{\@tocline{4}{0pt}{1pc}{7pc}{}}
		\def\l@subparagraph{\@tocline{5}{0pt}{1pc}{7pc}{}}
		\@xp\gdef\csname r@tocindent4\endcsname{0pt}
	}%
  {}%
\makeatother

\newcommand{\subsecunnum}[1]{
\subsection*{%
    \for{toc}{\protect\hphantom{\thesubsection.\quad}#1}%
    \except{toc}{#1}%
}}

\newcommand{\subsubsecunnum}[1]{
\subsubsection*{%
    \for{toc}{\protect\hphantom{\thesubsubsection.\quad}#1}%
    \except{toc}{#1}%
}}

\makeatletter
\renewcommand*\l@figure{\@tocline{1}{3pt plus2pt}{0pt}{}{}} 
\let\l@table\l@figure
\makeatother

\newcommand{\changelocaltocdepth}[1]{%
  \addtocontents{toc}{\protect\setcounter{tocdepth}{#1}}%
  \setcounter{tocdepth}{#1}%
}

\makeatletter
    { \let\if@nobreak\iffalse
  \let\if@noskipsec\iffalse
  \let\par\@@par
  \let\-\@dischyph
  \let\'\@acci\let\`\@accii\let\=\@acciii
  \parindent\z@ \parskip\z@skip
  \everypar{}%
     \linewidth\hsize \@totalleftmargin\z@
     \leftskip\z@skip \rightskip\z@skip \@rightskip\z@skip
     \vspace{7pt}
     \begin{adjustwidth}{}{\leftmargin} \centering %
    }{\end{adjustwidth} \vspace{0pt}
     }
\makeatother

 \tikzset{double line with arrow/.style args={#1,#2}{decorate,decoration={markings,%
mark=at position 0 with {\coordinate (ta-base-1) at (0,1pt);
\coordinate (ta-base-2) at (0,-1pt);},
mark=at position 1 with {\draw[#1] (ta-base-1) -- (0,1pt);
\draw[#2] (ta-base-2) -- (0,-1pt);
}}}}
 \tikzset{Equal/.style={-,double line with arrow={-,-}}}

\tikzset{
    labls/.style={anchor=south, rotate=90, inner sep=.5mm}
}
\tikzset{
    labln/.style={anchor=north, rotate=90, inner sep=1mm}
}
\tikzset{
    lablh/.style={anchor=north, rotate=45, inner sep=1mm}
}

\tikzset{mid vert/.style={/utils/exec=\tikzset{every node/.append style={outer sep=0.8ex}},
postaction=decorate,decoration={markings,
mark=at position 0.5 with {\draw[-] (0,#1) -- (0,-#1);}}},
mid vert/.default=0.75ex}

\definecolor{brightgreen}{rgb}{0.4, 1.0, 0.0}
\definecolor{brightturquoise}{rgb}{0.03, 0.91, 0.87}
\definecolor{brightpink}{rgb}{1.0, 0.0, 0.5}
\definecolor{carrotorange}{rgb}{0.93, 0.57, 0.13}

\definecolor{NBcolor}{rgb}{0.6,0.6,0.0}		

\newcounter{cdcnt}
\newcounter{cldcnt}
\newcounter{tdcnt}

\definecolor{CDcolor}{rgb}{0.0,0.5,0.75}	
\definecolor{XCDcolor}{rgb}{0.8,0.8,1}	
\definecolor{CLDcolor}{RGB}{220,100,10}
\definecolor{CLDEcolor}{rgb}{.6,0,.15}                    
\definecolor{XCLDcolor}{rgb}{1,.8,.8}		
\definecolor{TDcolor}{rgb}{0.4,0.4,0.0}		
\definecolor{RVcolor}{rgb}{0.4,0.0,0.8}		

\makeatletter%
\@ifclassloaded{amsbook}%
  {
		
	}%
  {}%
\makeatother

\makeatletter%
\newif\ifhyperref
\@ifpackageloaded{hyperref}{\hyperreftrue}{\hyperreffalse}

\@ifclassloaded{amsbook}%
  {
    \ifhyperref
    \def\defthm#1#2#3{
      \newtheorem{#1}{#2}[chapter]%
      \numberwithin{#1}{chapter}%
      \expandafter\def\csname #1autorefname\endcsname{#2}%
      \expandafter\let\csname c@#1\endcsname\c@thm}
	}%
  {
      \ifhyperref
    \def\defthm#1#2#3{
      \newtheorem{#1}{#2}[section]%
      \numberwithin{#1}{section}%
      \expandafter\def\csname #1autorefname\endcsname{#2}%
      \expandafter\let\csname c@#1\endcsname\c@thm}
     }%

\theoremstyle{plain} 

\newtheorem{introthm}{Theorem}
\expandafter\def\csname introthmautorefname\endcsname{Theorem}

\newtheorem{introcor}{Corollary}
\expandafter\def\csname introcorautorefname\endcsname{Corollary}
\expandafter\let\csname c@introcor\endcsname\c@introthm

\newtheorem{introdef}{Definition}
\expandafter\def\csname introdefautorefname\endcsname{Definition}
\expandafter\let\csname c@introdef\endcsname\c@introthm

\expandafter\def\csname intropropautorefname\endcsname{Proposition}
\expandafter\let\csname c@introprop\endcsname\c@introthm

\newtheorem{introconj}{Conjecture}
\expandafter\def\csname introconjautorefname\endcsname{Conjecture}
\expandafter\let\csname c@introconj\endcsname\c@introthm

\@ifclassloaded{amsbook}%
  {
		\newtheorem{thm}{Theorem}[chapter]
		\numberwithin{thm}{chapter}
		
	}%
  {
		\newtheorem{thm}{Theorem}[section]
		\numberwithin{thm}{section}
		
	}%
\defthm{cor}{Corollary}{Corollaries}
\defthm{prop}{Proposition}{Propositions}
\defthm{lem}{Lemma}{Lemmas}
\defthm{klem}{Key Lemma}{Key Lemmas}
\defthm{conj}{Conjecture}{Conjectures}
\defthm{disconj}{Disproven Conjecture}{Disproven Conjectures}
\defthm{hyp}{Hypothesis}{Hypotheses}
\defthm{indhyp}{Inductive Hypothesis}{Inductive Hypotheses}

\theoremstyle{definition}

\defthm{defnx}{Definition}{Definitions}
\newenvironment{defn}
  {\pushQED{\qed}\renewcommand{\qedsymbol}{$\envmarksym$}\defnx}
  {\popQED\enddefnx}

\numberwithin{altdefnx}{chapter}%
\expandafter\def\csname altdefnxautorefname\endcsname{Definition}%
\expandafter\let\csname c@altdefnx\endcsname\c@thm%

\defthm{nonegx}{Non-Example}{Non-Examples}

\defthm{egx}{Example}{Examples}
\newenvironment{eg}
  {\pushQED{\qed}\renewcommand{\qedsymbol}{$\envmarksym$}\egx}
  {\popQED\endegx}

\defthm{constrx}{Construction}{Constructions}
\newenvironment{constr}
  {\pushQED{\qed}\renewcommand{\qedsymbol}{$\envmarksym$}\constrx}
  {\popQED\endconstrx}

\theoremstyle{remark}

\defthm{rmkx}{Remark}{Remarks}
\newenvironment{rmk}
  {\pushQED{\qed}\renewcommand{\qedsymbol}{$\envmarksym$}\rmkx}
  {\popQED\endrmkx}

\defthm{punchx}{Punchline}{Punchlines}

\defthm{previewx}{Preview}{Previews}
\newenvironment{preview}
  {\pushQED{\qed}\renewcommand{\qedsymbol}{$\envmarksym$}\previewx}
  {\popQED\endpreviewx}

\defthm{obsx}{Observation}{Observations}
\newenvironment{obs}
  {\pushQED{\qed}\renewcommand{\qedsymbol}{$\envmarksym$}\obsx}
  {\popQED\endobsx}

\defthm{probx}{Problem}{Problems}
\newenvironment{prob}
  {\pushQED{\qed}\renewcommand{\qedsymbol}{$\envmarksym$}\probx}
  {\popQED\endprobx}

\defthm{notnx}{Notation}{Notations}
\newenvironment{notn}
  {\pushQED{\qed}\renewcommand{\qedsymbol}{$\envmarksym$}\notnx}
  {\popQED\endnotnx}

\defthm{termx}{Terminology}{Terminology}
\newenvironment{term}
  {\pushQED{\qed}\renewcommand{\qedsymbol}{$\envmarksym$}\termx}
  {\popQED\endtermx}

\defthm{convx}{Convention}{Conventions}
\newenvironment{conv}
  {\pushQED{\qed}\renewcommand{\qedsymbol}{$\envmarksym$}\convx}
  {\popQED\endconvx}

\defthm{outx}{Outline}{Outlines}

\@ifclassloaded{amsbook}%
  {
		\let\c@equation\c@thm
		\numberwithin{equation}{chapter}
	}%
  {
		\let\c@equation\c@thm
		\numberwithin{equation}{section}
	}%

\def\thmqedhere{\expandafter\csname\csname @currenvir\endcsname @qed\endcsname}

\def\definecref{\newif\ifcref}
\ifx\creftrue\undefined
  \definecref
  \creffalse
\fi

\ifcref\usepackage{cleveref,aliascnt}\fi

\ifcref\else
  \@ifpackageloaded{mathtools}{\mathtoolsset{showonlyrefs,showmanualtags}}{}
\fi

\makeatother

\newcommand{\eleq}{\ensuremath{\unlhd}}
\newcommand{\egeq}{\ensuremath{\unrhd}}
\newcommand{\eles}{\ensuremath{\lhd}}
\newcommand{\egre}{\ensuremath{\rhd}}
\newcommand{\fleq}{\ensuremath{\preceq}}

\newcommand{\fles}{\ensuremath{\prec}}
\newcommand{\fgre}{\ensuremath{\succ}}
\newcommand{\cov}{\ensuremath{^{\mathrm{cov}}}}

\newcommand{\lbl}[1]{\mathrm{lbl}_{#1}}
\newcommand{\fp}[2]{\pi^{#1}_{#2}}
\newcommand{\fpu}[2]{\und\pi^{#1}_{#2}}
\newcommand{\fcl}[2]{\chi^{#1}_{#2}}
\newcommand{\flblo}[1]{\mathsf{lbl}_{#1}}
\newcommand{\flbl}[2]{\mathsf{lbl}^{#1}_{#2}}

\makeatletter
\newcommand{\oset}[3][0ex]{%
  \mathrel{\mathop{#3}\limits^{
    \vbox to#1{\kern-2\ex@
    \hbox{$\scriptstyle#2$}\vss}}}}
\makeatother
\newcommand{\TT}{\ensuremath{\mathbb{T}}}
\newcommand{\OTT}{\ensuremath{\mathring{\TT}}}
\newcommand{\CTT}{\ensuremath{\bar{\TT}}}
\newcommand{\OCTT}{\ensuremath{\oset{\multimapinv}{\TT}}}
\newcommand{\COTT}{\ensuremath{\oset{\multimap}{\TT}}}

\newcommand{\truss}[1]{\ensuremath{{\mathsf{T}}\!\mathsf{rs}_{#1}}}
\newcommand{\rstruss}[1]{\ensuremath{{\mathsf{T}}\!\mathsf{rs}^{\mathsf{rs}}_{#1}}}
\newcommand{\struss}[1]{\ensuremath{{\mathsf{T}}\!\mathsf{rs}^{\mathsf{s}}_{#1}}}
\newcommand{\rtruss}[1]{\ensuremath{{\mathsf{T}}\!\mathsf{rs}^{\mathsf{r}}_{#1}}}

\newcommand{\otruss}[1]{\ensuremath{\mathring{\mathsf{T}}\!\mathsf{rs}_{#1}}}
\newcommand{\ctruss}[1]{\ensuremath{\bar{\mathsf{T}}\!\mathsf{rs}_{#1}}}
\newcommand{\sotruss}[1]{\ensuremath{\mathring{\mathsf{T}}\!\mathsf{rs}^{\mathsf{s}}_{#1}}}
\newcommand{\sctruss}[1]{\ensuremath{\bar{\mathsf{T}}\!\mathsf{rs}^{\mathsf{s}}_{#1}}}
\newcommand{\rotruss}[1]{\ensuremath{\mathring{\mathsf{T}}\!\mathsf{rs}^{\mathsf{r}}_{#1}}}
\newcommand{\rctruss}[1]{\ensuremath{\bar{\mathsf{T}}\!\mathsf{rs}^{\mathsf{r}}_{#1}}}
\newcommand{\rsctruss}[1]{\ensuremath{\bar{\mathsf{T}}\!\mathsf{rs}^{\mathsf{rs}}_{#1}}}
\newcommand{\rsotruss}[1]{\ensuremath{\mathring{\mathsf{T}}\!\mathsf{rs}^{\mathsf{rs}}_{#1}}}
\newcommand{\crstruss}[1]{\ensuremath{{\mathsf{T}}\!\mathsf{rs}^{\mathrm{crs}}_{#1}}}
\newcommand{\degtruss}[1]{\ensuremath{{\mathsf{T}}\!\mathsf{rs}^{\mathrm{deg}}_{#1}}}
\newcommand{\ttruss}[1]{\ensuremath{{\mathcal{T}\mkern-2mu\mathit{rs}_{#1}}}}

\newcommand{\trussbun}[1]{\ensuremath{\mathsf{TrsBun}_{#1}}}

\newcommand{\otrussbun}[1]{\ensuremath{\mathring{\mathsf{T}}\!\mathsf{rsBun}_{#1}}}
\newcommand{\ctrussbun}[1]{\ensuremath{\bar{\mathsf{T}}\!\mathsf{rsBun}_{#1}}}

\newcommand{\trusslbl}[1]{\ensuremath{\mathsf{Lbl}{\mathsf{T}}\!\mathsf{rs}_{#1}}}

\newcommand{\trussbunlbl}[1]{\ensuremath{\mathsf{Lbl}{\mathsf{T}}\!\mathsf{rsBun}_{#1}}}

\newcommand{\ttr}[1]{\ensuremath{{\mathsf{T}}\mathsf{Bord}^{#1}}}
\newcommand{\ottr}[1]{\ensuremath{\mathring{\mathsf{T}}\mathsf{Bord}^{#1}}}
\newcommand{\cttr}[1]{\ensuremath{\bar{\mathsf{T}}\mathsf{Bord}^{#1}}}
\newcommand{\lttr}[2]{\ensuremath{{\mathsf{T}}\mathsf{Bord}^{#1}_{\vslash #2}}}

\newcommand{\nlttr}[2]{\ensuremath{{#1}{\mathsf{T}}\mathsf{Bord}_{\vslash #2}}}

\newcommand{\trussconc}[1]{\ensuremath{\mathsf{TrsConc}_{#1}}}
\newcommand{\trussconclbl}[1]{\ensuremath{\mathsf{TrsConc}_{#1}}}

\newcommand{\qlttr}[2]{\ensuremath{\mathcal{T}\mkern-2mu\mathit{Bord}^{#1}_{\vslash #2}}}

\newcommand{\blcat}[1]{\ensuremath{\mathsf{B}\mathsf{lk}_{#1}}}
\newcommand{\blset}[1]{\ensuremath{\mathsf{Blk}\mathsf{Set}_{#1}}}
\newcommand{\rblset}[1]{\ensuremath{\mathsf{Reg}\mathsf{Blk}\mathsf{Set}_{#1}}}

\newcommand{\brcat}[1]{\ensuremath{\mathsf{B}\mathsf{rc}_{#1}}}
\newcommand{\brset}[1]{\ensuremath{\mathsf{Brc}\mathsf{Set}_{#1}}}

\DeclareMathOperator\bnrv{N_{\mathsf{Blk}}}

\newcommand{\qqq}[1]{\ensuremath{\mathsf{char}(#1)}}
\newcommand{\conn}[1]{\ensuremath{\mathsf{cc}(#1)}}
\newcommand{\sss}[1]{\ensuremath{\mathsf{discr}(#1)}}

\DeclareMathOperator\spine{\ensuremath{\mathsf{spine}}}

\DeclareMathOperator\FrEntr{\mathsf{Entr^{fr}}}
\DeclareMathOperator\Entr{\mathsf{Entr}}
\DeclareMathOperator\stratify{\mathsf{CStr}}
\newcommand{\CStr}[1]{\ensuremath{\stratify {#1}}}
\DeclareMathOperator\cellify{\mathsf{cell}}
\newcommand{\cellstrat}[1]{\ensuremath{\cellify {#1}}}

\DeclareMathOperator\FG{\ensuremath{\Phi^{\Delta}_1}}
\DeclareMathOperator\FC{\ensuremath{\Phi_1}}
\newcommand{\ino}{\ensuremath{\bot}}

\newcommand{\proto}[1]{\ensuremath{\xslashedrightarrow{#1}}}
\DeclareMathOperator\URel{\mathsf{rel}}

\newcommand{\Pos}{\ensuremath{\mathsf{Pos}}}
\newcommand{\CPos}{\ensuremath{\mathsf{CellPos}}}
\newcommand{\CStrat}{\ensuremath{\mathsf{CellStrat}}}
\newcommand{\TCStrat}{\mathcal{C}\mkern-2mu\mathit{ellStrat}}
\newcommand{\SCplx}{\ensuremath{\mathsf{SimpCplx}}}
\newcommand{\OrdSCplx}{\ensuremath{\mathsf{SimpCplx}^{\mathrm{ord}}}}
\newcommand{\CCplx}{\ensuremath{\mathsf{CellCplx}}}
\newcommand{\CCplxPL}{\ensuremath{\mathsf{CellCplx}^{\mathrm{PL}}}}

\newcommand{\FrSCplx}[1]{\ensuremath{\mathsf{FrSimpCplx}_{#1}}}
\newcommand{\PFrSCplx}[1]{\ensuremath{\mathsf{ProFrSimpCplx}_{#1}}}
\newcommand{\FlFrSCplx}[1]{\ensuremath{\mathsf{FlatFrSimpCplx}_{#1}}}
\newcommand{\FlPFrSCplx}[1]{\ensuremath{\mathsf{FlatProFrSimpCplx}_{#1}}}
\newcommand{\FrDelta}[1]{\ensuremath{\mathsf{FrSimp}_{#1}}}
\newcommand{\PFrDelta}[1]{\ensuremath{\mathsf{ProFrSimp}_{#1}}}
\newcommand{\PaFrDelta}[1]{\ensuremath{\mathsf{PartFrSimp}_{#1}}}
\newcommand{\PaPFrDelta}[1]{\ensuremath{\mathsf{PartProFrSimp}_{#1}}}
\newcommand{\FrCCell}[1]{\ensuremath{\mathsf{FrCell}_{#1}}}
\newcommand{\PFrCCell}[1]{\ensuremath{\mathsf{ProFrCell}_{#1}}}
\newcommand{\FrCDiag}[1]{\ensuremath{\mathsf{FlatFrCellCplx}_{#1}}}
\newcommand{\PFrCDiag}[1]{\ensuremath{\mathsf{FlatProFrCellCplx}_{#1}}}
\newcommand{\FrCCplx}[1]{\ensuremath{\mathsf{FrCellCplx}_{#1}}}
\newcommand{\PFrCCplx}[1]{\ensuremath{\mathsf{ProFrCellCplx}_{#1}}}
\newcommand{\FrTCCplx}[1]{\ensuremath{\mathsf{FrCellCplx}^{\mathrm{TOP}}_{#1}}}

\newcommand{\tmesh}[1]{\ensuremath{{\mathcal{M}\mkern-2mu\mathit{esh}_{#1}}}}
\newcommand{\tmeshbun}[1]{\ensuremath{{\mathcal{M}\mkern-2mu\mathit{eshBun}_{#1}}}}
\newcommand{\otmesh}[1]{\ensuremath{\mathring{\mathcal{M}}\mkern-2mu\mathit{esh}_{#1}}}
\newcommand{\ctmesh}[1]{\ensuremath{\bar{\mathcal{M}}\mkern-2mu\mathit{esh}_{#1}}}
\newcommand{\sctmesh}[1]{\ensuremath{\bar{\mathcal{M}}\mkern-2mu\mathit{esh}^{\mathsf{s}}_{#1}}}
\newcommand{\rotmesh}[1]{\ensuremath{\mathring{\mathcal{M}}\mkern-2mu\mathit{esh}^{\mathsf{r}}_{#1}}}

\newcommand{\mesh}[1]{\ensuremath{{\mathsf{M}}\mathsf{esh}_{#1}}}
\newcommand{\omesh}[1]{\ensuremath{\mathring{\mathsf{M}}\mathsf{esh}_{#1}}}
\newcommand{\cmesh}[1]{\ensuremath{\bar{\mathsf{M}}\mathsf{esh}_{#1}}}
\newcommand{\scmesh}[1]{\ensuremath{\bar{\mathsf{M}}\mathsf{esh}^{\mathsf{s}}_{#1}}}
\newcommand{\romesh}[1]{\ensuremath{\mathring{\mathsf{M}}\mathsf{esh}^{\mathsf{r}}_{#1}}}
\newcommand{\crsmesh}[1]{\ensuremath{{\mathsf{M}}\mathsf{esh}^{\mathrm{crs}}_{#1}}}
\newcommand{\degmesh}[1]{\ensuremath{{\mathsf{M}}\mathsf{esh}^{\mathrm{deg}}_{#1}}}
\newcommand{\crstmesh}[1]{\ensuremath{{\mathcal{M}\mkern-2mu\mathit{esh}^{\mathrm{crs}}_{#1}}}}
\newcommand{\degtmesh}[1]{\ensuremath{{\mathcal{M}\mkern-2mu\mathit{esh}^{\mathrm{deg}}_{#1}}}}

\DeclareMathOperator\ETrs{\mathsf{ETrs}}
\DeclareMathOperator\CMsh{\mathsf{CMsh}}
\DeclareMathOperator\CrsMsh{\mathsf{CrsMsh}}
\DeclareMathOperator\CellMsh{\mathsf{CellMsh}}
\DeclareMathOperator\MshCplx{\mathsf{MshCplx}}

\newcommand{\bdreg}[1]{\ensuremath{\ic_{#1}}}
\newcommand{\bdregbd}[1]{\ensuremath{\partial\ic_{#1}}}
\newcommand{\bdregclos}[1]{\ensuremath{\overline\ic_{#1}}}
\newcommand{\crreg}[1]{\ensuremath{\mathsf{cr}_{#1}}}
\newcommand{\crregclos}[1]{\ensuremath{\overline{\mathsf{cr}}_{#1}}}

\newcommand{\Unord}[1]{\ensuremath{{#1}^{\mathsf{un}}}}
\newcommand{\Order}[1]{\ensuremath{{#1}^{\mathsf{ord}}}}
\DeclareMathOperator\Unframe{\ensuremath{\mathsf{Unframe}}}
\DeclareMathOperator\Gradfr{\ensuremath{\nabla}}
\DeclareMathOperator\Intfr{\ensuremath{\smallint}}

\newcommand{\bnum}[1]{\ensuremath{\und{#1}}}
\newcommand{\unsimp}[1]{\ensuremath{\Unord {[#1]}}}
\newcommand{\lincplx}[1]{\ensuremath{\mathsf{lin}[{#1}]}}

\newcommand{\sstar}[1]{\ensuremath{\mathsf{star}({#1})}}

\newcommand{\akin}{\!\nperp\!}
\DeclareMathOperator\svec{\ensuremath{\mathsf{vec}}}

\usepackage{calc}
\newcommand{\vardbtilde}[1]{\tilde{\raisebox{0pt}[0.85\height]{$\tilde{#1}$}}}

\DeclareFontFamily{U}  {MnSymbolB}{}

\DeclareSymbolFont{MnSyB}         {U}  {MnSymbolB}{m}{n}

\SetSymbolFont{MnSyB}       {bold}{U}  {MnSymbolB}{b}{n}

\DeclareFontShape{U}{MnSymbolB}{m}{n}{
    <-6>  MnSymbolB5
   <6-7>  MnSymbolB6
   <7-8>  MnSymbolB7
   <8-9>  MnSymbolB8
   <9-10> MnSymbolB9
  <10-12> MnSymbolB10
  <12->   MnSymbolB12}{}
\DeclareFontShape{U}{MnSymbolB}{b}{n}{
    <-6>  MnSymbolB-Bold5
   <6-7>  MnSymbolB-Bold6
   <7-8>  MnSymbolB-Bold7
   <8-9>  MnSymbolB-Bold8
   <9-10> MnSymbolB-Bold9
  <10-12> MnSymbolB-Bold10
  <12->   MnSymbolB-Bold12}{}
  
\DeclareMathSymbol{\nperp}{\mathrel}{MnSyB}{217}

\makeatletter

\newcommand\mathrule[3][0pt]{%
  \ifdim#2>#3\math@hrule[#1]{#2}{#3}\else\math@vrule[#1]{#2}{#3}\fi}

\newcommand\math@hrule[3][0pt]{%
  \gdef\mystery@factor{0.07}%
 \@tempdima=#3%
  \rule[#1]{0pt}{#3}
  \raisebox{.5\@tempdima+#1}{%
    \makebox[#2][l]{\kern-.5\@tempdima\@@mathrule{#2}{#3}}}%
}

\newcommand\math@vrule[3][0pt]{%
  \gdef\mystery@factor{0.0}%
 \@tempdima=#2%
  \rule[#1]{0pt}{#3}
  \raisebox{-.0\@tempdima+#1}{%
    \kern0.5\@tempdima%
    \rotatebox{90}{\kern-0.5\@tempdima\makebox[#3][l]{\@@mathrule{#3}{#2}}}%
    \kern0.5\@tempdima}%
}

\def\@@mathrule#1#2{%
  \@tempdimb=#2%
  \@tempdima=\dimexpr#1-\mystery@factor\@tempdimb
  \pdfliteral{%
    q []0 d %
    1 J 
    \strip@pt\@tempdimb\space w \strip@pt\@tempdimb\space 0 m %
    \strip@pt\@tempdima\space 0 l S Q }}
\makeatother

\newif\iffootnoterule

\makeatletter
\AtBeginDocument{%
\let\latex@@footnoterule\footnoterule

\renewcommand\footnoterule{%
  \iffootnoterule
  \latex@@footnoterule%
  \else
  \advance\skip\footins 4\p@\@plus2\p@\relax%
  \fi
  }
}
\makeatother

\newcommand{\disk}{disk}

\newcommand{\conestr}{\ensuremath{\mathsf{cone}}}

\newcommand{\cstratum}{\ensuremath{\mathsf{str}}}

\newcommand{\singbif}{\ensuremath{\mathsf{sing}}}
\newcommand{\regbif}{\ensuremath{\mathsf{reg}}}

\newcommand{\fib}[1]{\mathsf{fib}(#1)}

\newcommand{\sff}{\mathsf{f}}

\newcommand{\adjtop}{{\triangleright}}

\newcommand{\NF}[1]{\left\llbracket{#1}\right\rrbracket}
\newcommand{\TEntr}{\ensuremath{\mathcal{E}\mkern-2mu\mathit{ntr}}}

\newcommand{\sing}{\ensuremath{\mathsf{sing}}}

\newcommand{\reg}{\ensuremath{\mathsf{reg}}}

\newcommand{\ra}{\rightarrow}

\newcommand{\ept}{\mathsf{end}}
\newcommand{\gend}{\ensuremath{\hat\gamma}}
\newcommand{\scaff}{\mathsf{scaff}}

\newcommand{\successor}{\ensuremath{\mathtt{s}}}
\newcommand{\predecessor}{\ensuremath{\mathtt{p}}}

\newcommand{\monolabel}{\mathsf{M}}
\newcommand{\epilabel}{\mathsf{E}}

\newcommand{\allsig}[1]{\ensuremath{n \cdot \cP}}

\newcommand{\Strat}{\mathsf{S}\mathsf{trat}}
\newcommand{\TStrat}{\mathcal{S}\mkern-2mu\mathit{trat}}

\newcommand{\TStratfinlch}{\mathcal{S}\mathsf{trat}^{\mathrm{lch}}_{\mathrm{fin}}}

\newcommand{\lvl}[2]{#2_{(#1)}}

\newcommand{\spx}[1]{\Delta[{#1}]}

\newcommand{\Tot}[2]{\mathrm{Tot}^{#1} #2} 
\newcommand{\Totb}[2]{\mathrm{Tot}^{#1}\mathopen{}\left(#2\right)\mathclose{}} 
\newcommand{\ff}{f\kern-0.04em f}

\newcommand{\Ka}{K} 
\newcommand{\La}{L} 

\renewcommand{\lvl}[2]{\ensuremath{{#2}_{#1}}}

\newcommand{\cret}{\ensuremath{\mathsf{cr}}}
\newcommand{\cint}{\ensuremath{\mathsf{ci}}}

\newcommand{\II}{\ensuremath{\mathbb{I}}}

\DeclareMathAlphabet{\mathpzc}{OT1}{pzc}{m}{it}

\newcommand{\powerset}{\raisebox{.15\baselineskip}{\Large\ensuremath{\wp}}}

\newcommand{\vslash}{\mathbin{/\mkern-6mu/}}

\newcommand{\und}[1]{{\underline{#1}}}

\newcommand{\abs}[1]{\left|{#1}\right|}

\newcommand{\conelink}{\mathsf{link}}

\newcommand{\avg}[1]{\left\langle{#1}\right\rangle}
\newcommand{\rest}[2]{{\left.\kern-\nulldelimiterspace #1 \right|_{#2}}}

\newcommand{\eps}{\epsilon}

\newcommand{\obj}{\mathrm{obj}}
\newcommand{\Bool}{\ensuremath{\mathsf{Bool}}}
\newcommand{\BoolProf}{\ensuremath{\mathsf{BoolProf}}}

\newcommand{\Prof}{\ensuremath{\mathsf{Prof}}}

\newcommand{\Rel}{\ensuremath{\mathsf{Rel}}}
\newcommand{\TPos}{\ensuremath{\mathcal{P}\mkern-2mu\mathit{os}}}

\newcommand{\SetCat}{\ensuremath{\mathsf{Set}}}
\newcommand{\UnSimp}{\ensuremath{\und{\Delta}}}

\newcommand{\Cat}{\ensuremath{\mathsf{Cat}}}
\newcommand{\Top}{\ensuremath{\mathsf{Top}}}

\newcommand{\SSet}{\ensuremath{\mathsf{SSet}}}

\newcommand{\ho}{\mathrm{ho}}

\newcommand{\nid}{\noindent}

\newcommand{\sectionline}[4]{%
  \nointerlineskip 
  \vspace{#2\baselineskip}\hspace{\fill}
  {\centering\color{#1}
    \rule{.3\textwidth}{#4}
  }%
    \hspace{\fill}
    \par\nointerlineskip 
    \vspace{#3\baselineskip}
}

\newcommand{\pause}{\sectionline{black}{1}{1}{.5pt}}
\newcommand{\pauseae}{\sectionline{black}{.5}{1}{.5pt}}

\newcommand{\ignore}[1]{}

\DeclareRobustCommand{\SkipTocEntry}[5]{}

\makeatletter
\let\ea\expandafter

\def\mdef#1#2{\ea\ea\ea\gdef\ea\ea\noexpand#1\ea{\ea\ensuremath\ea{#2}\xspace}}
\def\alwaysmath#1{\ea\ea\ea\global\ea\ea\ea\let\ea\ea\csname your@#1\endcsname\csname #1\endcsname
  \ea\def\csname #1\endcsname{\ensuremath{\csname your@#1\endcsname}\xspace}}

\DeclareRobustCommand\widecheck[1]{{\mathpalette\@widecheck{#1}}}
\def\@widecheck#1#2{%
    \setbox\z@\hbox{\m@th$#1#2$}%
    \setbox\tw@\hbox{\m@th$#1%
       \widehat{%
          \vrule\@width\z@\@height\ht\z@
          \vrule\@height\z@\@width\wd\z@}$}%
    \dp\tw@-\ht\z@
    \@tempdima\ht\z@ \advance\@tempdima2\ht\tw@ \divide\@tempdima\thr@@
    \setbox\tw@\hbox{%
       \raise\@tempdima\hbox{\scalebox{1}[-1]{\lower\@tempdima\box
\tw@}}}%
    {\ooalign{\box\tw@ \cr \box\z@}}}

\newcount\foreachcount

\def\foreachletter#1#2#3{\foreachcount=#1
  \ea\loop\ea\ea\ea#3\@alph\foreachcount
  \advance\foreachcount by 1
  \ifnum\foreachcount<#2\repeat}

\def\foreachLetter#1#2#3{\foreachcount=#1
  \ea\loop\ea\ea\ea#3\@Alph\foreachcount
  \advance\foreachcount by 1
  \ifnum\foreachcount<#2\repeat}

\let\oldit\it

\def\definescr#1{\ea\gdef\csname s#1\endcsname{\ensuremath{\mathscr{#1}}\xspace}}
\foreachLetter{1}{27}{\definescr}
\def\definecal#1{\ea\gdef\csname c#1\endcsname{\ensuremath{\mathcal{#1}}\xspace}}
\foreachLetter{1}{27}{\definecal}
\def\definebold#1{\ea\gdef\csname b#1\endcsname{\ensuremath{\mathbf{#1}}\xspace}}
\foreachLetter{1}{27}{\definebold}
\def\definebb#1{\ea\gdef\csname l#1\endcsname{\ensuremath{\mathbb{#1}}\xspace}}
\foreachLetter{1}{27}{\definebb}
\def\definefrak#1{\ea\gdef\csname k#1\endcsname{\ensuremath{\mathfrak{#1}}\xspace}}
\foreachletter{1}{27}{\definefrak}
\foreachLetter{1}{27}{\definefrak}
\def\definesf#1{\ea\gdef\csname i#1\endcsname{\ensuremath{\mathsf{#1}}\xspace}}
\foreachletter{1}{6}{\definesf}
\foreachletter{7}{14}{\definesf}
\foreachletter{15}{27}{\definesf}
\foreachLetter{1}{27}{\definesf}
\def\definebar#1{\ea\gdef\csname #1bar\endcsname{\ensuremath{\overline{#1}}\xspace}}
\foreachLetter{1}{27}{\definebar}
\foreachletter{1}{8}{\definebar} 
\foreachletter{9}{15}{\definebar} 
\foreachletter{16}{27}{\definebar}
\def\definetil#1{\ea\gdef\csname #1til\endcsname{\ensuremath{\widetilde{#1}}\xspace}}
\foreachLetter{1}{27}{\definetil}
\foreachletter{1}{27}{\definetil}
\def\definehat#1{\ea\gdef\csname #1hat\endcsname{\ensuremath{\widehat{#1}}\xspace}}
\foreachLetter{1}{27}{\definehat}
\foreachletter{1}{27}{\definehat}
\def\definechk#1{\ea\gdef\csname #1chk\endcsname{\ensuremath{\widecheck{#1}}\xspace}}
\foreachLetter{1}{27}{\definechk}
\foreachletter{1}{27}{\definechk}
\def\defineul#1{\ea\gdef\csname u#1\endcsname{\ensuremath{\underline{#1}}\xspace}}
\foreachLetter{1}{27}{\defineul}
\foreachletter{1}{27}{\defineul}

\let\it\oldit

\def\autofmt@b#1\autofmt@end{\mathbf{#1}}
\def\autofmt@l#1#2\autofmt@end{\mathbb{#1}\mathsf{#2}}
\def\autofmt@c#1#2\autofmt@end{\mathcal{#1}\mathit{#2}}
\def\autofmt@s#1#2\autofmt@end{\mathscr{#1}\mathit{#2}}
\def\autofmt@f#1\autofmt@end{\mathsf{#1}}
\def\autofmt@k#1\autofmt@end{\mathfrak{#1}}
\def\autofmt@u#1\autofmt@end{\underline{\smash{\mathsf{#1}}}}
\def\autofmt@U#1\autofmt@end{\underline{\underline{\smash{\mathsf{#1}}}}}
\def\autofmt@h#1\autofmt@end{\widehat{#1}}
\def\autofmt@r#1\autofmt@end{\overline{#1}}
\def\autofmt@t#1\autofmt@end{\widetilde{#1}}
\def\autofmt@k#1\autofmt@end{\check{#1}}

\def\auto@drop#1{}
\def\autodef#1{\ea\ea\ea\@autodef\ea\ea\ea#1\ea\auto@drop\string#1\autodef@end}
\def\@autodef#1#2#3\autodef@end{%
  \ea\def\ea#1\ea{\ea\ensuremath\ea{\csname autofmt@#2\endcsname#3\autofmt@end}\xspace}}
\def\autodefs@end{blarg!}
\def\autodefs#1{\@autodefs#1\autodefs@end}
\def\@autodefs#1{\ifx#1\autodefs@end%
  \def\autodefs@next{}%
  \else%
  \def\autodefs@next{\autodef#1\@autodefs}%
  \fi\autodefs@next}

\DeclareSymbolFont{bbold}{U}{bbold}{m}{n}
\DeclareSymbolFontAlphabet{\mathbbb}{bbold}

\alwaysmath{alpha}
\alwaysmath{beta}
\alwaysmath{gamma}
\alwaysmath{Gamma}
\alwaysmath{delta}
\alwaysmath{Delta}
\alwaysmath{epsilon}
\mdef\ep{\varepsilon}
\alwaysmath{zeta}
\alwaysmath{eta}
\alwaysmath{theta}
\alwaysmath{Theta}
\alwaysmath{iota}
\alwaysmath{kappa}
\alwaysmath{lambda}
\alwaysmath{Lambda}
\alwaysmath{mu}
\alwaysmath{nu}
\alwaysmath{xi}
\alwaysmath{pi}
\alwaysmath{rho}
\alwaysmath{sigma}
\alwaysmath{Sigma}
\alwaysmath{tau}
\alwaysmath{upsilon}
\alwaysmath{Upsilon}
\alwaysmath{phi}
\alwaysmath{Pi}
\alwaysmath{Phi}
\mdef\ph{\varphi}
\alwaysmath{chi}
\alwaysmath{psi}
\alwaysmath{Psi}
\alwaysmath{omega}
\alwaysmath{Omega}

\definecolor{cblue}{RGB}{0,0,255}
\definecolor{cdarkblue}{RGB}{0, 38, 153}
\definecolor{clightblue}{RGB}{128, 170, 255}
\definecolor{cturquoise}{RGB}{0, 255, 255}
\definecolor{cgreen}{RGB}{0, 230, 0}
\definecolor{cdarkgreen}{RGB}{0, 128, 0}
\definecolor{corange}{RGB}{255, 153, 0}
\definecolor{cred}{RGB}{255, 51, 0}
\definecolor{cpurple}{RGB}{172, 0, 230}
\definecolor{cyellow}{RGB}{230, 230, 0}
\definecolor{cpink}{RGB}{255, 51, 204}

\DeclareSymbolFont{extraup}{U}{zavm}{m}{n}
\DeclareMathSymbol{\varheart}{\mathalpha}{extraup}{86}
\DeclareMathSymbol{\vardiamond}{\mathalpha}{extraup}{87}

\renewcommand{\Set}[1]{\{#1\}}

\mdef\delbar{\overline{\partial}}

\newcommand{\inv}{^{-1}}

\mdef\hf{\textstyle\frac12 }
\mdef\thrd{\textstyle\frac13 }
\mdef\qtr{\textstyle\frac14 }

\newcommand{\op}{^{\mathrm{op}}}

\let\iso\cong
\let\eqv\simeq

\mdef\Id{\mathrm{Id}}
\mdef\id{\mathrm{id}}
\mdef\uni{\mathrm{is\_uni}}
\mdef\sym{\mathrm{is\_sym}}
\mdef\rel{\mathrm{is\_rel}}
\mdef\Sym{\mathrm{Sym}}
\mdef\SYM{\textsc{sym}}
\mdef\REL{\textsc{rel}}
\alwaysmath{ell}
\alwaysmath{infty}
\alwaysmath{odot}
\def\frc#1/#2.{\frac{#1}{#2}}   
\mdef\ten{\mathrel{\otimes}}

\mdef\sqten{\mathrel{\boxtimes}}

\DeclareMathOperator\keraff{ker^{aff}}
\DeclareMathOperator\cokeraff{coker^{aff}}
\DeclareMathOperator\imaff{im^{aff}}
\newcommand{\intoaff}{\ensuremath{\;\dottedarrow{}\;}}
\newcommand{\simpzero}{\ensuremath{\mathbbm{0}}}

\DeclareMathOperator\supp{supp}

\DeclareMathOperator\coker{coker}

\DeclareMathOperator\dom{dom}
\DeclareMathOperator\cod{cod}
\DeclareMathOperator\mor{mor}
\DeclareMathOperator\PSh{PSh}

\DeclareMathOperator\Fun{Fun}

\mdef\Im{\mathrm{Im}}
\mdef\im{\mathrm{im}}
\let\lim\relax
\DeclareMathOperator\lim{lim}

\DeclareMathOperator\colim{colim}

\DeclareMathOperator\Spec{Spcl}

\DeclareMathOperator\Hom{Hom}
\DeclareMathOperator\Map{Map}
\DeclareRobustCommand{\ocirc}{%
  \mathbin{\mathpalette\on@ntimes\relax}%
}
\newcommand{\on@ntimes}[2]{%
  \vcenter{\hbox{%
    \sbox0{\m@th$#1\otimes$}%
    \setlength\unitlength{\wd0}%
    \begin{picture}(1,1)
    \linethickness{0.35pt}
    \put(.5,.5){\circle{.8}}
    \end{picture}%
  }}%
}

\newcommand{\ot}{\ensuremath{\leftarrow}}

\let\toot\rightleftarrows

\let\imp\Rightarrow

\let\into\hookrightarrow

\mdef\we{\overset{\sim}{\longrightarrow}}
\mdef\leftwe{\overset{\sim}{\longleftarrow}}

\let\epi\twoheadrightarrow

\newcounter{sarrow}

\let\xto\xrightarrow
\let\xot\xleftarrow
\def\rightarrowtailfill@{\arrowfill@{\Yright\joinrel\relbar}\relbar\rightarrow}
\newcommand\xrightarrowtail[2][]{\ext@arrow 0055{\rightarrowtailfill@}{#1}{#2}}

\def\twoheadrightarrowfill@{\arrowfill@{\relbar\joinrel\relbar}\relbar\twoheadrightarrow}
\newcommand\xtwoheadrightarrow[2][]{\ext@arrow 0055{\twoheadrightarrowfill@}{#1}{#2}}
\let\xepi\xtwoheadrightarrow

\def\slashedarrowfill@#1#2#3#4#5{%
  $\m@th\thickmuskip0mu\medmuskip\thickmuskip\thinmuskip\thickmuskip
   \relax#5#1\mkern-7mu%
   \cleaders\hbox{$#5\mkern-2mu#2\mkern-2mu$}\hfill
   \mathclap{#3}\mathclap{#2}%
   \cleaders\hbox{$#5\mkern-2mu#2\mkern-2mu$}\hfill
   \mkern-7mu#4$%
}
\def\rightslashedarrowfill@{%
  \slashedarrowfill@\relbar\relbar\mapstochar\rightarrow}
\newcommand\xslashedrightarrow[2][]{%
  \ext@arrow 0055{\rightslashedarrowfill@}{#1}{#2}}
\def\leftslashedarrowfill@{%
  \slashedarrowfill@\relbar\relbar\mapstochar\leftarrow}
\newcommand\xslashedleftarrow[2][]{%
  \ext@arrow 0055{\leftslashedarrowfill@}{#1}{#2}}
\def\hookrightslashedarrowfill@{%
  \slashedarrowfill@\relbar\relbar\mapstochar\hookrightarrow}
\newcommand\xslashedhookrightarrow[2][]{%
  \ext@arrow 0055{\hookrightslashedarrowfill@}{#1}{#2}}
\mdef\hto{\xslashedrightarrow{}}
\mdef\htoo{\xslashedrightarrow{\quad}}

\newbox\dottedarrow@box
\setbox\dottedarrow@box\hbox
  {%
    \begin{tikzpicture}
      \draw[dotted, semithick ,right hook->] (0,0) -- (1.5em,0);
    \end{tikzpicture}%
  }
\newcommand*\dottedarrow
  {\relax\ifmmode\expandafter\dottedarrow@m\else\expandafter\dottedarrow@t\fi}
\newcommand*\dottedarrow@t[1][1.3em]
  {\resizebox{#1}{!}{\raisebox{-.04ex}{\usebox\dottedarrow@box}}}
\newcommand*\dottedarrow@m[1][]
  {%
    \if\relax\detokenize{#1}\relax
      \mathchoice
        {\dottedarrow@t}
        {\dottedarrow@t}
        {\dottedarrow@t[1.1em]}
        {\dottedarrow@t[0.9em]}%
    \else
      \dottedarrow@t[#1]%
    \fi
  }

\newbox\envmarksym@box
\setbox\envmarksym@box\hbox
  {%
    \begin{tikzpicture}
      \draw (0,0) -- (1.5em,0) -- (1.5em,.5em);
    \end{tikzpicture}%
  }
\newcommand*\envmarksym[1][1.5em]
  {\resizebox{#1}{!}{\raisebox{.0ex}{\usebox\envmarksym@box}}}

\def\toiso{\xto{\smash{\raisebox{-.5mm}{$\scriptstyle\sim$}}}}



%% file: figuresused/1-truss-bordism.pdf_tex
\begingroup%
  \makeatletter%
  \providecommand\color[2][]{%
    \errmessage{(Inkscape) Color is used for the text in Inkscape, but the package 'color.sty' is not loaded}%
    \renewcommand\color[2][]{}%
  }%
  \providecommand\transparent[1]{%
    \errmessage{(Inkscape) Transparency is used (non-zero) for the text in Inkscape, but the package 'transparent.sty' is not loaded}%
    \renewcommand\transparent[1]{}%
  }%
  \providecommand\rotatebox[2]{#2}%
  \newcommand*\fsize{\dimexpr\f@size pt\relax}%
  \newcommand*\lineheight[1]{\fontsize{\fsize}{#1\fsize}\selectfont}%
  \ifx\svgwidth\undefined%
    \setlength{\unitlength}{680.31496063bp}%
    \ifx\svgscale\undefined%
      \relax%
    \else%
      \setlength{\unitlength}{\unitlength * \real{\svgscale}}%
    \fi%
  \else%
    \setlength{\unitlength}{\svgwidth}%
  \fi%
  \global\let\svgwidth\undefined%
  \global\let\svgscale\undefined%
  \makeatother%
  \begin{picture}(1,0.35416667)%
    \lineheight{1}%
    \setlength\tabcolsep{0pt}%
    \put(0,0){\includegraphics[width=\unitlength,page=1]{1-truss-bordism.pdf}}%
    \put(0.34328351,0.0077031){\makebox(0,0)[lt]{\lineheight{1.25}\smash{\begin{tabular}[t]{l}$T$\end{tabular}}}}%
    \put(0.64288583,0.0077031){\makebox(0,0)[lt]{\lineheight{1.25}\smash{\begin{tabular}[t]{l}$S$\end{tabular}}}}%
    \put(0,0){\includegraphics[width=\unitlength,page=2]{1-truss-bordism.pdf}}%
    \put(0.48780015,0.0077031){\makebox(0,0)[lt]{\lineheight{1.25}\smash{\begin{tabular}[t]{l}$R$\end{tabular}}}}%
    \put(0,0){\includegraphics[width=\unitlength,page=3]{1-truss-bordism.pdf}}%
  \end{picture}%
\endgroup%

%% file: figuresused/1-truss-bundle.pdf_tex
\begingroup%
  \makeatletter%
  \providecommand\color[2][]{%
    \errmessage{(Inkscape) Color is used for the text in Inkscape, but the package 'color.sty' is not loaded}%
    \renewcommand\color[2][]{}%
  }%
  \providecommand\transparent[1]{%
    \errmessage{(Inkscape) Transparency is used (non-zero) for the text in Inkscape, but the package 'transparent.sty' is not loaded}%
    \renewcommand\transparent[1]{}%
  }%
  \providecommand\rotatebox[2]{#2}%
  \newcommand*\fsize{\dimexpr\f@size pt\relax}%
  \newcommand*\lineheight[1]{\fontsize{\fsize}{#1\fsize}\selectfont}%
  \ifx\svgwidth\undefined%
    \setlength{\unitlength}{680.31496063bp}%
    \ifx\svgscale\undefined%
      \relax%
    \else%
      \setlength{\unitlength}{\unitlength * \real{\svgscale}}%
    \fi%
  \else%
    \setlength{\unitlength}{\svgwidth}%
  \fi%
  \global\let\svgwidth\undefined%
  \global\let\svgscale\undefined%
  \makeatother%
  \begin{picture}(1,0.375)%
    \lineheight{1}%
    \setlength\tabcolsep{0pt}%
    \put(0,0){\includegraphics[width=\unitlength,page=1]{1-truss-bundle.pdf}}%
    \put(0.56181571,0.0500117){\makebox(0,0)[lt]{\lineheight{1.25}\smash{\begin{tabular}[t]{l}$B$\end{tabular}}}}%
    \put(0.37137966,0.35536919){\makebox(0,0)[lt]{\lineheight{1.25}\smash{\begin{tabular}[t]{l}$T$\end{tabular}}}}%
    \put(0,0){\includegraphics[width=\unitlength,page=2]{1-truss-bundle.pdf}}%
    \put(0.50123606,0.14035454){\makebox(0,0)[lt]{\lineheight{1.25}\smash{\begin{tabular}[t]{l}$p$\end{tabular}}}}%
    \put(0,0){\includegraphics[width=\unitlength,page=3]{1-truss-bundle.pdf}}%
  \end{picture}%
\endgroup%

%% file: figuresused/1-truss-bordisms-composition-as-a-bundle.pdf_tex
\begingroup%
  \makeatletter%
  \providecommand\color[2][]{%
    \errmessage{(Inkscape) Color is used for the text in Inkscape, but the package 'color.sty' is not loaded}%
    \renewcommand\color[2][]{}%
  }%
  \providecommand\transparent[1]{%
    \errmessage{(Inkscape) Transparency is used (non-zero) for the text in Inkscape, but the package 'transparent.sty' is not loaded}%
    \renewcommand\transparent[1]{}%
  }%
  \providecommand\rotatebox[2]{#2}%
  \newcommand*\fsize{\dimexpr\f@size pt\relax}%
  \newcommand*\lineheight[1]{\fontsize{\fsize}{#1\fsize}\selectfont}%
  \ifx\svgwidth\undefined%
    \setlength{\unitlength}{680.31496063bp}%
    \ifx\svgscale\undefined%
      \relax%
    \else%
      \setlength{\unitlength}{\unitlength * \real{\svgscale}}%
    \fi%
  \else%
    \setlength{\unitlength}{\svgwidth}%
  \fi%
  \global\let\svgwidth\undefined%
  \global\let\svgscale\undefined%
  \makeatother%
  \begin{picture}(1,0.5)%
    \lineheight{1}%
    \setlength\tabcolsep{0pt}%
    \put(0,0){\includegraphics[width=\unitlength,page=1]{1-truss-bordisms-composition-as-a-bundle.pdf}}%
    \put(0.69390221,0.42867902){\color[rgb]{0,0,0}\makebox(0,0)[lt]{\lineheight{1.25}\smash{\begin{tabular}[t]{l}$T$\end{tabular}}}}%
    \put(0,0){\includegraphics[width=\unitlength,page=2]{1-truss-bordisms-composition-as-a-bundle.pdf}}%
    \put(0.823758,0.21366437){\makebox(0,0)[lt]{\lineheight{1.25}\smash{\begin{tabular}[t]{l}$p$\end{tabular}}}}%
    \put(0,0){\includegraphics[width=\unitlength,page=3]{1-truss-bordisms-composition-as-a-bundle.pdf}}%
    \put(0.12256761,0.34302699){\color[rgb]{0,0,0}\makebox(0,0)[lt]{\lineheight{1.25}\smash{\begin{tabular}[t]{l}$R_1$\end{tabular}}}}%
    \put(0,0){\includegraphics[width=\unitlength,page=4]{1-truss-bordisms-composition-as-a-bundle.pdf}}%
    \put(0.3958444,0.34152298){\makebox(0,0)[lt]{\lineheight{1.25}\smash{\begin{tabular}[t]{l}$R_2$\end{tabular}}}}%
    \put(0.6961935,0.0565328){\makebox(0,0)[lt]{\lineheight{1.25}\smash{\begin{tabular}[t]{l}$[2]$\end{tabular}}}}%
  \end{picture}%
\endgroup%

%% file: figuresused/suspension-bundle.pdf_tex
\begingroup%
  \makeatletter%
  \providecommand\color[2][]{%
    \errmessage{(Inkscape) Color is used for the text in Inkscape, but the package 'color.sty' is not loaded}%
    \renewcommand\color[2][]{}%
  }%
  \providecommand\transparent[1]{%
    \errmessage{(Inkscape) Transparency is used (non-zero) for the text in Inkscape, but the package 'transparent.sty' is not loaded}%
    \renewcommand\transparent[1]{}%
  }%
  \providecommand\rotatebox[2]{#2}%
  \newcommand*\fsize{\dimexpr\f@size pt\relax}%
  \newcommand*\lineheight[1]{\fontsize{\fsize}{#1\fsize}\selectfont}%
  \ifx\svgwidth\undefined%
    \setlength{\unitlength}{680.31496063bp}%
    \ifx\svgscale\undefined%
      \relax%
    \else%
      \setlength{\unitlength}{\unitlength * \real{\svgscale}}%
    \fi%
  \else%
    \setlength{\unitlength}{\svgwidth}%
  \fi%
  \global\let\svgwidth\undefined%
  \global\let\svgscale\undefined%
  \makeatother%
  \begin{picture}(1,0.41666667)%
    \lineheight{1}%
    \setlength\tabcolsep{0pt}%
    \put(0,0){\includegraphics[width=\unitlength,page=1]{suspension-bundle.pdf}}%
    \put(0.17769941,0.1244098){\makebox(0,0)[lt]{\lineheight{1.25}\smash{\begin{tabular}[t]{l}$p$\end{tabular}}}}%
    \put(0,0){\includegraphics[width=\unitlength,page=2]{suspension-bundle.pdf}}%
    \put(0.7581164,0.14117997){\makebox(0,0)[lt]{\lineheight{1.25}\smash{\begin{tabular}[t]{l}$\Sigma p$\end{tabular}}}}%
    \put(0,0){\includegraphics[width=\unitlength,page=3]{suspension-bundle.pdf}}%
    \put(0.14473515,0.35347818){\makebox(0,0)[lt]{\lineheight{1.25}\smash{\begin{tabular}[t]{l}$T$\end{tabular}}}}%
    \put(0.76248531,0.39267553){\makebox(0,0)[lt]{\lineheight{1.25}\smash{\begin{tabular}[t]{l}$\Sigma T$\end{tabular}}}}%
    \put(0.13375988,0.01481323){\makebox(0,0)[lt]{\lineheight{1.25}\smash{\begin{tabular}[t]{l}$X$\end{tabular}}}}%
    \put(0.50378272,0.00697378){\makebox(0,0)[lt]{\lineheight{1.25}\smash{\begin{tabular}[t]{l}$\Sigma X$\end{tabular}}}}%
  \end{picture}%
\endgroup%

%% file: figuresused/fillers-vs-fiber-morphisms.pdf_tex
\begingroup%
  \makeatletter%
  \providecommand\color[2][]{%
    \errmessage{(Inkscape) Color is used for the text in Inkscape, but the package 'color.sty' is not loaded}%
    \renewcommand\color[2][]{}%
  }%
  \providecommand\transparent[1]{%
    \errmessage{(Inkscape) Transparency is used (non-zero) for the text in Inkscape, but the package 'transparent.sty' is not loaded}%
    \renewcommand\transparent[1]{}%
  }%
  \providecommand\rotatebox[2]{#2}%
  \newcommand*\fsize{\dimexpr\f@size pt\relax}%
  \newcommand*\lineheight[1]{\fontsize{\fsize}{#1\fsize}\selectfont}%
  \ifx\svgwidth\undefined%
    \setlength{\unitlength}{680.31496063bp}%
    \ifx\svgscale\undefined%
      \relax%
    \else%
      \setlength{\unitlength}{\unitlength * \real{\svgscale}}%
    \fi%
  \else%
    \setlength{\unitlength}{\svgwidth}%
  \fi%
  \global\let\svgwidth\undefined%
  \global\let\svgscale\undefined%
  \makeatother%
  \begin{picture}(1,0.3125)%
    \lineheight{1}%
    \setlength\tabcolsep{0pt}%
    \put(0,0){\includegraphics[width=\unitlength,page=1]{fillers-vs-fiber-morphisms.pdf}}%
    \put(0.52224496,0.0855196){\makebox(0,0)[lt]{\lineheight{1.25}\smash{\begin{tabular}[t]{l}$p$\end{tabular}}}}%
    \put(0.53044287,0.00961583){\makebox(0,0)[lt]{\lineheight{1.25}\smash{\begin{tabular}[t]{l}$[2]$\end{tabular}}}}%
    \put(0.58023542,0.2915198){\color[rgb]{0,0,0}\makebox(0,0)[lt]{\lineheight{1.25}\smash{\begin{tabular}[t]{l}$T$\end{tabular}}}}%
    \put(0,0){\includegraphics[width=\unitlength,page=2]{fillers-vs-fiber-morphisms.pdf}}%
  \end{picture}%
\endgroup%

%% file: figuresused/scaffold-order-on-sections.pdf_tex
\begingroup%
  \makeatletter%
  \providecommand\color[2][]{%
    \errmessage{(Inkscape) Color is used for the text in Inkscape, but the package 'color.sty' is not loaded}%
    \renewcommand\color[2][]{}%
  }%
  \providecommand\transparent[1]{%
    \errmessage{(Inkscape) Transparency is used (non-zero) for the text in Inkscape, but the package 'transparent.sty' is not loaded}%
    \renewcommand\transparent[1]{}%
  }%
  \providecommand\rotatebox[2]{#2}%
  \newcommand*\fsize{\dimexpr\f@size pt\relax}%
  \newcommand*\lineheight[1]{\fontsize{\fsize}{#1\fsize}\selectfont}%
  \ifx\svgwidth\undefined%
    \setlength{\unitlength}{680.31496063bp}%
    \ifx\svgscale\undefined%
      \relax%
    \else%
      \setlength{\unitlength}{\unitlength * \real{\svgscale}}%
    \fi%
  \else%
    \setlength{\unitlength}{\svgwidth}%
  \fi%
  \global\let\svgwidth\undefined%
  \global\let\svgscale\undefined%
  \makeatother%
  \begin{picture}(1,0.22916667)%
    \lineheight{1}%
    \setlength\tabcolsep{0pt}%
    \put(0,0){\includegraphics[width=\unitlength,page=1]{scaffold-order-on-sections.pdf}}%
    \put(0.52836119,0.20437183){\color[rgb]{0,0,0}\makebox(0,0)[lt]{\lineheight{1.25}\smash{\begin{tabular}[t]{l}$T$\end{tabular}}}}%
    \put(0,0){\includegraphics[width=\unitlength,page=2]{scaffold-order-on-sections.pdf}}%
    \put(0.58834896,0.20390905){\color[rgb]{0.6,0.6,0.6}\makebox(0,0)[lt]{\lineheight{1.25}\smash{\begin{tabular}[t]{l}$\Sigma T$\end{tabular}}}}%
    \put(0,0){\includegraphics[width=\unitlength,page=3]{scaffold-order-on-sections.pdf}}%
  \end{picture}%
\endgroup%

%% file: figuresused/scaffold-norm-on-fillers.pdf_tex
\begingroup%
  \makeatletter%
  \providecommand\color[2][]{%
    \errmessage{(Inkscape) Color is used for the text in Inkscape, but the package 'color.sty' is not loaded}%
    \renewcommand\color[2][]{}%
  }%
  \providecommand\transparent[1]{%
    \errmessage{(Inkscape) Transparency is used (non-zero) for the text in Inkscape, but the package 'transparent.sty' is not loaded}%
    \renewcommand\transparent[1]{}%
  }%
  \providecommand\rotatebox[2]{#2}%
  \newcommand*\fsize{\dimexpr\f@size pt\relax}%
  \newcommand*\lineheight[1]{\fontsize{\fsize}{#1\fsize}\selectfont}%
  \ifx\svgwidth\undefined%
    \setlength{\unitlength}{680.31496063bp}%
    \ifx\svgscale\undefined%
      \relax%
    \else%
      \setlength{\unitlength}{\unitlength * \real{\svgscale}}%
    \fi%
  \else%
    \setlength{\unitlength}{\svgwidth}%
  \fi%
  \global\let\svgwidth\undefined%
  \global\let\svgscale\undefined%
  \makeatother%
  \begin{picture}(1,0.20833333)%
    \lineheight{1}%
    \setlength\tabcolsep{0pt}%
    \put(0,0){\includegraphics[width=\unitlength,page=1]{scaffold-norm-on-fillers.pdf}}%
    \put(0.6214543,0.18190015){\color[rgb]{0,0,0}\makebox(0,0)[lt]{\lineheight{1.25}\smash{\begin{tabular}[t]{l}$T$\end{tabular}}}}%
    \put(0,0){\includegraphics[width=\unitlength,page=2]{scaffold-norm-on-fillers.pdf}}%
    \put(0.68144207,0.18143737){\color[rgb]{0.6,0.6,0.6}\makebox(0,0)[lt]{\lineheight{1.25}\smash{\begin{tabular}[t]{l}$\Sigma T$\end{tabular}}}}%
    \put(0,0){\includegraphics[width=\unitlength,page=3]{scaffold-norm-on-fillers.pdf}}%
  \end{picture}%
\endgroup%

%% file: figuresused/boundaries-of-fillers.pdf_tex
\begingroup%
  \makeatletter%
  \providecommand\color[2][]{%
    \errmessage{(Inkscape) Color is used for the text in Inkscape, but the package 'color.sty' is not loaded}%
    \renewcommand\color[2][]{}%
  }%
  \providecommand\transparent[1]{%
    \errmessage{(Inkscape) Transparency is used (non-zero) for the text in Inkscape, but the package 'transparent.sty' is not loaded}%
    \renewcommand\transparent[1]{}%
  }%
  \providecommand\rotatebox[2]{#2}%
  \newcommand*\fsize{\dimexpr\f@size pt\relax}%
  \newcommand*\lineheight[1]{\fontsize{\fsize}{#1\fsize}\selectfont}%
  \ifx\svgwidth\undefined%
    \setlength{\unitlength}{680.31496063bp}%
    \ifx\svgscale\undefined%
      \relax%
    \else%
      \setlength{\unitlength}{\unitlength * \real{\svgscale}}%
    \fi%
  \else%
    \setlength{\unitlength}{\svgwidth}%
  \fi%
  \global\let\svgwidth\undefined%
  \global\let\svgscale\undefined%
  \makeatother%
  \begin{picture}(1,0.22916667)%
    \lineheight{1}%
    \setlength\tabcolsep{0pt}%
    \put(0,0){\includegraphics[width=\unitlength,page=1]{boundaries-of-fillers.pdf}}%
    \put(0.58572655,0.1974479){\color[rgb]{0,0,0}\makebox(0,0)[lt]{\lineheight{1.25}\smash{\begin{tabular}[t]{l}$T$\end{tabular}}}}%
    \put(0,0){\includegraphics[width=\unitlength,page=2]{boundaries-of-fillers.pdf}}%
    \put(0.72586924,0.07343362){\makebox(0,0)[lt]{\lineheight{1.25}\smash{\begin{tabular}[t]{l}$L'$\end{tabular}}}}%
    \put(0.72586924,0.10392195){\makebox(0,0)[lt]{\lineheight{1.25}\smash{\begin{tabular}[t]{l}$\partial_+ L'$\end{tabular}}}}%
    \put(0.72586924,0.04294523){\makebox(0,0)[lt]{\lineheight{1.25}\smash{\begin{tabular}[t]{l}$\partial_- L'$\end{tabular}}}}%
    \put(0.16683281,0.12498525){\makebox(0,0)[lt]{\lineheight{1.25}\smash{\begin{tabular}[t]{l}$L$\end{tabular}}}}%
    \put(0.16683281,0.15547358){\makebox(0,0)[lt]{\lineheight{1.25}\smash{\begin{tabular}[t]{l}$\partial_+ L$\end{tabular}}}}%
    \put(0.16683281,0.09449685){\makebox(0,0)[lt]{\lineheight{1.25}\smash{\begin{tabular}[t]{l}$\partial_- L$\end{tabular}}}}%
  \end{picture}%
\endgroup%